\documentclass[11pt,a4paper,oneside]{amsbook}

\usepackage[T1]{fontenc}
\usepackage[utf8]{inputenc}
\usepackage{lmodern}
\usepackage[english]{babel}
\usepackage[a4paper,margin=25mm,headheight=14pt]{geometry}
\usepackage{amsmath,amssymb,amsthm,amscd}
\allowdisplaybreaks
\usepackage{mathtools}
\usepackage[export]{adjustbox}
\usepackage[normalem]{ulem}
\usepackage{xcolor}
\usepackage{verbatim}
\usepackage{graphicx}
\usepackage{float}
\usepackage{placeins}
\usepackage[all]{xy}
\usepackage{wrapfig}
\usepackage{cases}
\usepackage{array}
\usepackage{longtable,booktabs}
\makeindex
\usepackage{mathrsfs}
\usepackage{tikz}
\usepackage{centernot}
\usepackage{multirow}
\usepackage{url}
\usepackage[final]{microtype}
\usepackage{fix-cm}
\usetikzlibrary{patterns}
\usetikzlibrary{intersections,calc,decorations.markings}
\usetikzlibrary{arrows.meta}

\definecolor{annred}{RGB}{255,0,0}

\theoremstyle{plain}
\newtheorem{theorem}{Theorem}[section]
\newtheorem{prop}[theorem]{Proposition}
\newtheorem{conj}[theorem]{Conjecture}
\newtheorem{prob}[theorem]{Problem}
\newtheorem{lemm}[theorem]{Lemma}
\newtheorem{coro}[theorem]{Corollary}

\theoremstyle{definition}
\newtheorem{defi}[theorem]{Definition}
\newtheorem{exam}[theorem]{Example}

\theoremstyle{remark}
\newtheorem{rema}[theorem]{Remark}

\numberwithin{equation}{section}

\makeatletter
\renewcommand{\p@section}{\thechapter.}
\makeatother

\def\TT{\mathbb{T}}

\def\ZZ{\mathbb{Z}}

\DeclareMathOperator{\g}{\gamma}

\title{An Introduction to the Lagrange and Markov Spectra\\
through the Lens of Generalized Markov Numbers}
\author{Yasuaki Gyoda}
\address{Institute for Advanced Research, Nagoya University, Furo-cho, Chikusa-ku, Nagoya-shi, 464-8601, Japan}
\email{ygyoda@math.nagoya-u.ac.jp}

\makeatletter
\@ifundefined{subjclassname@2020}{%
  \@namedef{subjclassname@2020}{\textup{2020} Mathematics Subject Classification}%
}{}
\makeatother
\usepackage[
  bookmarks=true,
  bookmarksnumbered=true,
  bookmarksopen=true,
  bookmarksopenlevel=1,
  colorlinks=true,
  linkcolor=blue!100!black,
  citecolor=blue!100!black,
  urlcolor=red!100!black,
  pdfauthor={Yasuaki Gyoda},
  pdftitle={An Introduction to the Lagrange and Markov Spectra through the Lens of Generalized Markov Numbers},
  pdfsubject={Lagrange and Markov spectra; generalized Markov numbers},
  pdfkeywords={Lagrange spectrum, Markov spectrum, generalized Markov number, cluster algebra}
]{hyperref}
\usepackage{bookmark}
\begin{document}

\frontmatter

\begin{abstract}
This text is a self-contained expository survey of the Lagrange and Markov spectra, centered on a comprehensive exposition of Markov's theorem and its generalizations. Its purpose is to provide a systematic text for learning the theory, with detailed proofs and explanations of the connections among its arithmetic, combinatorial, and geometric descriptions. The necessary background in continued fractions, quadratic irrationals, binary quadratic forms, and bi-infinite sequences is developed step by step, followed by an exposition of generalized Markov numbers, fence posets, curve lengths, and generalized Cohn matrices.

One goal of this exposition is to explain the formula connecting generalized Markov numbers with the two spectra. For nonnegative integer parameters $(k_1,k_2,k_3)$, a permutation $\sigma\in\mathfrak S_3$, and a fraction label $t\in\mathbb Q_{\geq0}\cup\{\infty\}$, let $m_t$ be the associated generalized Markov number and let $k_t=k_{i_t}$ be the parameter assigned to its position $i_t\in\{1,2,3\}$. The text explains the construction of an associated finite sequence $S(t)$ of positive integers and the identity
\[
\mathcal L(\alpha_{S(t)})
=\mathcal M(Q_{S(t)})
=\frac{\sqrt{((3+k_1+k_2+k_3)m_t-k_t)^2-4}}{m_t},
\]
where $\alpha_{S(t)}=[\overline{S(t)}]$ and
$Q_{S(t)}=(x-\alpha_{S(t)}y)(x-\alpha'_{S(t)}y)$,
with the prime denoting quadratic conjugation. Here $\mathcal L$ and $\mathcal M$ denote the Lagrange and Markov constants, respectively.

The survey explains how this identity relates generalized discrete Markov spectra to the classical theory and how Markov's theorem is recovered when the parameters vanish. The account also includes boundary values arising from irrational slopes and generalizations of Frobenius's uniqueness conjecture, providing a unified perspective on the classical theorem and its extensions.
\end{abstract}  

\maketitle

\setcounter{tocdepth}{1}
\tableofcontents

\mainmatter
\chapter{Background and Organization of the Text}\label{chap:background-organization}
The Lagrange and Markov spectra enter number theory through two elementary-looking questions.  One concerns the approximation of real numbers by rational numbers, and the other concerns the values of indefinite binary quadratic forms on integral points.  The definitions are simple, but the resulting spectra are far from elementary.  Already at the first stage, one is led to continued fractions, quadratic irrationals, binary quadratic forms, and two-sided infinite sequences.

This text studies these spectra through the lens of generalized Markov numbers.  The role of this opening chapter is to explain the background needed for that viewpoint and to indicate how the rest of the text is organized.  Rather than assuming that the reader is already familiar with the classical theory in detail, we use this chapter as a guide to the prerequisite material: which notions are needed, why they are introduced, and how they later become connected with generalized Markov numbers.

\section{History of the Lagrange and Markov Spectra and Generalized Markov Numbers}
We begin with a brief historical overview of the objects studied in this text and of the surrounding theory.

\subsection{Emergence of Continued Fractions and Diophantine Approximation}
Before discussing the Lagrange and Markov spectra, let us review the part of Diophantine approximation theory that underlies them. \index{Diophantine approximation}Roughly speaking, Diophantine approximation asks how well an irrational number $\alpha$ can be approximated by rational numbers $\frac pq$. There are many possible meanings of ``how well,'' but the most elementary question is the following.
\begin{prob}
Let $\alpha$ be an irrational number. For every $\varepsilon>0$, does there always exist a rational number $\frac pq$ satisfying
\[
\left|\alpha-\frac{p}{q}\right|<\varepsilon?
\]
\end{prob}
Yes. From the modern construction of the real numbers this is immediate, and in fact infinitely many such rational numbers exist. It is nevertheless natural to go one step further and ask the following question.
\begin{prob}
How can one construct a sequence of rational numbers with good approximation properties that converges to an irrational number $\alpha$?
\end{prob}
One answer is provided by the sequence of rational numbers obtained by truncating the continued-fraction expansion. This point of view goes back to Euler's paper \cite{euler1737}, written in 1737 and published in 1744. The relation between irrational numbers and infinite continued fractions discovered by Euler is now understood as the following correspondence. This does not mean, however, that Euler himself proved the theorem in this modern form.
\begin{theorem}
Let $\mathscr S$ be the set of all infinite sequences whose first entry is an integer and whose subsequent entries are positive integers. Then the map
\[
F:\mathscr S\to\mathbb R\setminus\mathbb Q,\quad (a_k)_{k=0}^{\infty}\mapsto [a_0;a_1,a_2,\dots]:=a_0+\frac{1}{a_1+\frac{1}{a_2+\frac{1}{\ddots}}}
\]
is a bijection.
\end{theorem}
Lagrange also used continued fractions in his 1770 paper \cite{lagrange1770} to characterize quadratic irrationals, that is, irrational numbers that occur as roots of quadratic equations with rational coefficients.
\begin{theorem}[Lagrange's Theorem]
The continued-fraction expansion of an irrational number $\alpha$ is eventually periodic if and only if $\alpha$ is a quadratic irrational.
\end{theorem}
The theory of continued fractions developed during this period later became a central tool in Diophantine approximation.

In the nineteenth century, the basic question of how well irrational numbers can be approximated by rational numbers came to be studied in the following quantitative form.
\begin{prob}\label{prob:approximation}
For an irrational number $\alpha$, how large can one take $L>0$ and $n>0$ so that there exist infinitely many rational numbers $\frac pq$ satisfying
\[
\left|\alpha-\frac{p}{q}\right|<\frac{1}{Lq^n}?
\]
\end{prob}
Let us examine this problem through a concrete example, namely the irrational number $\pi=3.141592\dots$. Rational numbers very close to $\pi$ exist no matter how small an error tolerance we impose. For instance, rational numbers satisfying $\left|\pi-\frac pq\right|<\frac1{1000}$ can be produced as
\begin{equation}\label{eq:p/q}
\frac pq=\frac{3141}{1000},\frac{6283}{2000}\left(=\frac{31415}{10000}\right),\frac{314159}{100000},\frac{392699}{125000}\left(=\frac{3141592}{1000000}\right),\frac{15707963}{5000000}\left(=\frac{31415926}{10000000}\right),\dots.
\end{equation}
However, all of these rational numbers have relatively large denominators. In general, the smaller the required error is, the larger the denominator of a rational number satisfying it must be. Fractions with small denominator are sparse on the number line, and hence are less likely to lie close to a specified irrational number.

Thus, in rational approximation, one must control not only the error but also the size of the denominator. The ``smallness of the denominator relative to the approximation error'' is measured by the parameters $n$ and $L$ when the error bound is written in the form $\frac{1}{Lq^n}$. Since the scale is difficult to interpret if both $L$ and $n$ are allowed to vary simultaneously, one usually fixes one of them and studies the other.

First fix $L=1$ and consider the supremum of the possible values of $n$ for the rational numbers appearing in \eqref{eq:p/q}. This value is computed as $n=-\frac{\log\left|\pi-\frac pq\right|}{\log q}$, and for the fractions displayed above it is approximately $1.076, 1.222, 1.115, 1.213, 1.085$, from left to right. On the other hand, the same computation for $\frac{22}{7}$ and $\frac{355}{113}$ gives approximately $3.429$ and $3.202$, respectively, which are much larger. This means that $\frac{22}{7}$ and $\frac{355}{113}$ give much better approximations than one would expect from the size of their denominators.

If an approximation for which $L$ or $n$ can be taken large is called a good approximation, then Problem~\ref{prob:approximation} asks the following: for a given irrational number $\alpha$, how far can we raise the parameters $L$ and $n$ measuring the quality of approximation before rational approximations of that quality cease to exist infinitely often?

The meaning of this question is not yet completely clear. What property of an irrational number is being measured? To clarify this, we first ask what it means for a number to have few good rational approximations. The following fact is fundamental.
\begin{theorem}\label{thm:rational}
For any rational number $\alpha$, the supremum of the real numbers $n$ for which there exist infinitely many reduced fractions $\frac pq$ satisfying
\[
\left|\alpha-\frac{p}{q}\right|<\frac{1}{q^n}
\]
is $1$.
\end{theorem}
Let us compare this with Dirichlet's theorem \cite{dirichlet1842}, one of the starting points of Diophantine approximation theory.
\begin{theorem}[Dirichlet's Theorem]\index{Dirichlet's theorem}
For any irrational number $\alpha$, the supremum of the real numbers $n$ for which there exist infinitely many rational numbers $\frac pq$ satisfying
\[
\left|\alpha-\frac{p}{q}\right|<\frac{1}{q^n}
\]
is at least $2$.
\end{theorem}
These two theorems show that rational numbers themselves are the real numbers with the fewest good rational approximations. Equivalently, when $L=1$ is fixed, the supremum of the possible exponents $n$ can be regarded as a measure of how far a real number is, arithmetically, from being rational.

\index{Irrationality exponent}This supremum is called the \emph{irrationality exponent} and is denoted by $\mu(\alpha)$. By Theorem~\ref{thm:rational}, the irrationality exponent of a rational number is $1$, while Dirichlet's theorem implies that $\mu(\alpha)\geq 2$ for every irrational number $\alpha$. The existence of irrational numbers with $\mu(\alpha)=2$ follows from the following theorem of Liouville \cite{Liouville}.
\begin{theorem}[Liouville's Theorem]\index{Liouville's theorem}
Let $\alpha$ be an algebraic irrational number of degree $d$. Then there exists a constant $L>0$ such that
\[
\left|\alpha-\frac{p}{q}\right|>\frac{1}{Lq^d}
\]
for every rational number $\frac pq$.
\end{theorem}
Combining this theorem with Dirichlet's theorem, it follows that every quadratic irrational has irrationality exponent $2$. Indeed, suppose that a quadratic irrational $\alpha$ satisfied $\mu(\alpha)>2$. Then, for some $\varepsilon>0$, there would be infinitely many reduced fractions $\frac pq$ satisfying $\left|\alpha-\frac pq\right|<q^{-2-\varepsilon/2}$.

On the other hand, Liouville's theorem gives a constant $L>0$ such that $\left|\alpha-\frac pq\right|>\frac{1}{Lq^2}$ for all $\frac pq$. For sufficiently large $q$ we have $\frac1{Lq^2}>q^{-2-\varepsilon/2}$, a contradiction.

In this way, the elementary problem of approximating an irrational number by rational numbers developed into a theory that measures arithmetic properties of numbers through the irrationality exponent. Pursuing this topic further would take us away from the main theme of the present text, so we close this discussion by recalling Roth's theorem \cite{roth1955}.
\begin{theorem}[Roth's Theorem]\index{Roth's theorem}
If $\alpha$ is an algebraic irrational number, then its irrationality exponent is $2$.
\end{theorem}
Roth's theorem is a decisive strengthening of Liouville's theorem. Not only quadratic irrationals but all algebraic irrational numbers have irrationality exponent at most $2$. Hence any irrational number whose irrationality exponent is larger than $2$ must be transcendental; Roth's theorem therefore also gives a powerful sufficient condition for transcendence.

\subsection{Minimization Problems for the Lagrange and Markov Constants}
In view of Dirichlet's and Liouville's theorems, for any irrational number $\alpha$ there are infinitely many rational numbers $\frac pq$ satisfying $\left|\alpha-\frac{p}{q}\right|<\frac1{q^n}$ up to the exponent $n=2$, while the exponent cannot be uniformly increased beyond this. The next natural problem is therefore to fix the exponent at $n=2$ and ask how large the constant $L$ can be.

For an irrational number $\alpha$, the supremum of the real numbers $L$ for which there exist infinitely many rational numbers $\frac pq$ satisfying
\[
\left|\alpha-\frac{p}{q}\right|<\frac{1}{Lq^2}
\]
is called the \emph{Lagrange constant} of $\alpha$ and is denoted by $\mathcal L(\alpha)$. The problem above asks for the smallest possible value of the Lagrange constant. Hurwitz gave the answer in 1891 \cite{hurwitz}.
\begin{theorem}[Hurwitz's Theorem]\label{thm:hurwitz}\index{Hurwitz's theorem}
For every irrational number $\alpha$ one has $\mathcal L(\alpha)\geq\sqrt{5}$, and for example $\mathcal L(\alpha)=\sqrt{5}$ when $\alpha=\frac{1+\sqrt{5}}{2}$.
\end{theorem}
Hurwitz also stated that the next smallest Lagrange constant after $\sqrt{5}$ is $2\sqrt{2}$, and wrote that this fact follows from Markov's work. The proof of this point, however, is not given in \cite{hurwitz}.

Markov's work approached these values through a related minimization problem for indefinite binary quadratic forms. For such a form $Q(x,y)=ax^2+bxy+cy^2$, put $D=b^2-4ac$ and define $\mathcal M(Q)$ by
\[
\inf_{(x,y)\in\mathbb Z^2\setminus\{(0,0)\}}|Q(x,y)|=\frac{\sqrt D}{\mathcal M(Q)}.
\]
This value is called the \emph{Markov constant}. The smallest possible Markov constant is $\sqrt{5}$, attained by $Q(x,y)=x^2-xy-y^2$, and the next smallest value is $2\sqrt{2}$, attained by $Q(x,y)=x^2-2xy-y^2$.

These two extremal results had already been announced by Korkin--Zolotarev in 1873 \cite{korkin-zolotarev1873}. Markov developed the problem into a systematic theory in his papers of 1879 and 1880 \cite{mar1,mar2}, using continued fractions to determine the smallest Markov constants. In his 1891 paper, Hurwitz explicitly pointed out the consequences of this theory for rational approximation. Continued fractions thus connect these two minimization problems.

The central result of Markov's theory is now called \emph{Markov's theorem}.
\begin{theorem}[Markov's Theorem]
\index{Markov equation}Let $M$ be the set of positive integers that occur in positive integer solutions of
\[
x^2+y^2+z^2=3xyz.
\]
If $\mathcal M_0$ denotes the set of Markov constants less than $3$, then
\[
\mathcal M_0=\left\{\frac{\sqrt{9m^2-4}}{m}\ \middle|\ m\in M\right\}.
\]
\end{theorem}
The elements of $M$ are called \emph{Markov numbers}. They should not be confused with Markov constants. The theorem says that Markov constants less than $3$ are completely described by Markov numbers. For example, the Markov number $1$ gives $\sqrt{5}$, and the Markov number $2$ gives $2\sqrt{2}$. These values also occur as Lagrange constants via continued fraction theory, but it took some time before this relation was organized in a clear form.

\subsection{Research on the Lagrange and Markov Spectra}
After the minimization problems for the Lagrange and Markov constants were solved, attention turned to the problem of understanding the sets formed by all such constants. Let $\mathcal L$ be the set of all Lagrange constants and let $\mathcal M$ be the set of all Markov constants. They are called the \emph{Lagrange spectrum} and the \emph{Markov spectrum}, respectively.\footnote{Markov's contribution to the Lagrange spectrum is also substantial, and some authors may prefer the name Markov--Lagrange spectrum. The author is sympathetic to this view. In this text, however, we use the standard name in order to avoid a cumbersome terminology and confusion with the Markov spectrum.}

The fact that these sets can be described by continued fractions is now formulated as follows.
\begin{theorem}\label{thm:perron-history}
For a bi-infinite sequence $\mathbf a=(a_n)_{n\in\mathbb Z}\in\mathbb Z_{\geq1}^{\mathbb Z}$, put
\[
\ell_n(\mathbf a):=[a_n;a_{n+1},a_{n+2},\dots]+[0;a_{n-1},a_{n-2},\dots].
\]
Then
\[
\mathcal L=\left\{\limsup_{n\to+\infty}\ell_n(\mathbf a)\ \middle|\ \mathbf a\in\mathbb Z_{\geq1}^{\mathbb Z}\right\},\quad
\mathcal M=\left\{\sup_{n\in\mathbb Z}\ell_n(\mathbf a)\ \middle|\ \mathbf a\in\mathbb Z_{\geq1}^{\mathbb Z}\right\}.
\]
\end{theorem}
Theorem~\ref{thm:perron-history} is a modern formulation of the continued-fraction descriptions discussed in Perron's papers of 1921. In Part~I, \S1, Perron derived the limsup formula for the approximation constant \cite{perron1}. In Part~II, \S1, he compared this formula with Markov's supremum problem for bi-infinite continued-fraction sequences \cite{perron2}. The displayed equalities are commonly referred to as \emph{Perron's identity}.

Perron's identity allows both the Lagrange spectrum and the Markov spectrum to be treated as the limit superior or the supremum of a function on bi-infinite continued fraction sequences. From this point of view, one also naturally obtains the inclusion $\mathcal L\subset\mathcal M$. Moreover, below $3$ the two spectra coincide completely, and the Lagrange constants are described by Markov numbers just as in Markov's theorem.
\begin{theorem}
For Lagrange constants below $3$,
\[
\mathcal L\cap(-\infty,3)=\left\{\frac{\sqrt{9m^2-4}}{m}\ \middle|\ m\in M\right\}.
\]
\end{theorem}

Since Perron's formulation, much work has been done on the parts of $\mathcal L$ and $\mathcal M$ above $3$. The following topics are somewhat outside the main line of this text, but they are important for understanding the global structure of the spectra.

\index{Markov spectrum!difference from the Lagrange spectrum}First, as mentioned above, one has $\mathcal L\subset\mathcal M$, and Freiman proved that this inclusion is strict \cite{freiman}.
\begin{theorem}
$\mathcal L\subsetneq\mathcal M$. In other words, $\mathcal M\setminus\mathcal L\ne\emptyset$.
\end{theorem}
Another historically important result was proved by Hall in 1947 \cite{hall1947}.
\begin{theorem}
The interval $[6,\infty)$ is contained in $\mathcal L$. Consequently, $[6,\infty)$ is also contained in $\mathcal M$.
\end{theorem}
This means that every sufficiently large real number belongs to the Lagrange spectrum, and hence also to the Markov spectrum. Such a half-line is called a \emph{Hall ray}. Freiman later determined, in 1975, the smallest possible initial point of such a ray \cite{freiman1975}.
\begin{theorem}
The largest half-line contained in $\mathcal L$ is $[c_F,\infty)$, where
\[
 c_F=\frac{2221564096+283748\sqrt{462}}{491993569}\approx4.5278295661\cdots.
\]
In particular, $c_F\in\mathcal L$.
\end{theorem}
The number $c_F$ is called the \emph{Freiman constant}. It follows that on $[c_F,\infty)$ both $\mathcal L$ and $\mathcal M$ contain an entire real half-line, and any set-theoretic difference between $\mathcal L$ and $\mathcal M$ is contained in the interval $[3,c_F)$.

\index{Hausdorff dimension}The structure of $\mathcal L$ and $\mathcal M$ in the remaining interval $[3,c_F)$ is still an active subject of research. For example, Moreira proved in 2018 the following result from the viewpoint of Hausdorff dimension \cite{moreira18}.
\begin{theorem}
For every $t\in\mathbb R$,
\[
\dim_H(\mathcal L\cap(-\infty,t))=\dim_H(\mathcal M\cap(-\infty,t)).
\]
If this common value is denoted by $d(t)$, then $d(t)$ is nondecreasing and
\[
\max\{t\in\mathbb R\mid d(t)=0\}=3.
\]
\end{theorem}
In 2024 Erazo--Lima--Matheus--Moreira--Vieira proved the following \cite{elmmv}.
\begin{theorem}
$\inf(\mathcal M\setminus\mathcal L)=3$.
\end{theorem}
These results show that $\mathcal L$ and $\mathcal M$ have closely related fractal structures inside $[3,c_F)$, while their set-theoretic difference already appears immediately after $3$. In this way, the study of the Lagrange and Markov spectra, although rooted in classical continued fraction theory, continues to develop today.

\subsection{Markov Numbers and Reduced Fractions}
The Markov numbers that describe the part of the Lagrange and Markov spectra below $3$ have been studied in many contexts beyond their original motivation in Diophantine approximation. A starting point for this development was Frobenius's 1913 paper \cite{frobenius}. In that paper, Frobenius related Markov numbers to reduced fractions. This correspondence shows that Markov numbers are deeply connected with rational numbers, lattice points, and line segments in the plane, and it still plays a fundamental role in modern work on Markov numbers.

To explain this relation, take a reduced fraction $t=p/q\geq1$. We state Frobenius's correspondence using the fraction-label convention adopted later in this text; see also \cite{frobenius,RabideauSchiffler2020}. Set $m_1=2$, and assume below that $p>q\geq1$.
\begin{itemize}
\item[(1)] For $i=0,1,\dots,p$, let $r_i$ be the remainder of $iq$ upon division by $p$. For each $i=1,2,\dots,p-2$, write $c$ if $r_i<r_{i+1}$ and $d$ if $r_i>r_{i+1}$, obtaining a word $s$.
\item[(2)] Replace each $c$ by $1,1$ and each $d$ by $2,2$, obtaining an integer sequence $S$. Define $m_t$ as the numerator of the reduced fraction $[2;S,2]$. For $p=2$, both $s$ and $S$ are empty and the continued fraction is $[2;2]$.
\end{itemize}
\begin{theorem}
The integer $m_t$ constructed above is a Markov number. Moreover, $t\mapsto m_t$ is a surjection from the reduced fractions at least $1$ onto the Markov numbers other than $1$.
\end{theorem}
For example, $m_2=5$. For $t=3$, the word is $s=c$ and $[2;1,1,2]=13/5$, so $m_3=13$. For $t=3/2$, the word is $s=d$ and $[2;2,2,2]=29/12$, so $m_{3/2}=29$. These agree with the fraction labels of the $(0,0,0)$-GM tree introduced later.

The theorem has two significant features. First, reduced fractions provide labels for Markov numbers. A Markov number is initially defined as an entry of a positive integer solution of the Markov equation. Fraction labels locate these entries in the tree and make individual numbers easier to describe. This viewpoint is used throughout the text.

Second, the construction has a geometric interpretation. For $p>q$, let $\ell$ join $(0,0)$ to $(p,q)$. Its intersection with $x=i$ has $y$-coordinate $iq/p$ and fractional part $r_i/p$. Thus comparing consecutive remainders describes its passage through the vertical strips
\[
i\leq x\leq i+1\qquad(1\leq i\leq p-2).
\]
Since $q/p<1$, the segment crosses at most one horizontal lattice line in each strip. It crosses none when $r_i<r_{i+1}$ and exactly one when $r_i>r_{i+1}$. The word $s$ records the crossings after omitting the first and last strips; the two endpoint entries $2$ are supplied separately in the continued fraction. The remainder construction therefore records how a rational-slope segment crosses the lattice, linking Markov numbers with lattice geometry and word combinatorics.

Frobenius also formulated a simple but very important conjecture about Markov numbers.
\begin{conj}[Frobenius's Uniqueness Conjecture]
For any Markov number $c$, a positive integer solution of the Markov equation whose largest component is $c$ is uniquely determined up to permutation of its components.
\end{conj}
In terms of fraction labels, this can be understood as the question of whether distinct reduced fractions give distinct Markov numbers.
\begin{conj}
The fraction-label map $t\mapsto m_t$ on reduced fractions $t\geq1$ is injective.
\end{conj}
Although the statement is concise, it is a difficult problem about the internal structure of Markov numbers, and it remains open in full generality. Frobenius's work not only gave a way to describe Markov numbers by reduced fractions, but also introduced a central problem that has continued to be studied ever since.

\subsection{Matrix Realizations of Markov Numbers and Their Hyperbolic-Geometric Interpretation}
Around the 1950s, a point of view developed in which Markov numbers are realized through elements of the modular group. Using Fricke's trace identity, Cohn \cite{cohn} related the traces of free generating pairs of the commutator subgroup $[SL(2,\mathbb Z),SL(2,\mathbb Z)]$ to the Markov equation in his 1955 paper. For any such pair $A,B$, the identity gives
\[
(\operatorname{tr}(A))^2+(\operatorname{tr}(B))^2+(\operatorname{tr}(AB))^2=\operatorname{tr}(A)\operatorname{tr}(B)\operatorname{tr}(AB).
\]
This gives the equation $x^2+y^2+z^2=xyz$, which contains the same information as the Markov equation $x^2+y^2+z^2=3xyz$ up to a scaling. Indeed, from a solution $(a,b,c)$ of the Markov equation one obtains a solution $(3a,3b,3c)$ of the former equation, and conversely positive integer solutions of the former equation give positive integer solutions of the Markov equation after division by $3$. Thus the Markov equation can also be viewed as a problem about traces of matrices.

The importance of Cohn's work is not merely that an equation resembling the Markov equation appears. To a matrix $M=\begin{bsmallmatrix}a&b\\ c&d\end{bsmallmatrix}\in SL(2,\mathbb Z)$ he associated the binary quadratic form $Q(x,y)=cx^2+(d-a)xy-by^2$, thereby providing a way to reinterpret the minimization problem for indefinite binary quadratic forms studied by Markov in matrix language. This made clear that the Markov numbers appearing in continued fractions and binary quadratic forms also arise naturally in the theory of discrete groups.

This line of thought led to a more hyperbolic-geometric interpretation in Cohn's 1971 paper \cite{cohn1971}. If the commutator subgroup acts on the upper half-plane $\mathbb H=\{x+iy\in\mathbb C\mid y>0\}$, the quotient is a once-punctured torus. Conjugacy classes of hyperbolic elements determine free homotopy classes of closed curves on the torus. In particular, conjugacy classes represented by primitive elements, meaning elements belonging to a free basis, correspond to essential simple closed curves, whose geodesic representatives are simple closed geodesics. This notion of primitivity is stronger than merely not being a proper power. If $A$ is the corresponding matrix and $\ell(A)$ is the length of the geodesic, then $|\operatorname{tr}(A)|=2\cosh\left(\frac{\ell(A)}{2}\right)$. Thus the fact that traces occur as three times Markov numbers means that Markov numbers are directly connected with lengths of simple closed geodesics on the once-punctured torus.

This geometric interpretation is closely related to fraction labels. An ordinary torus is obtained from the plane by identifying points differing by an integer vector. Primitive homotopy classes of closed curves are represented there by rational-slope lines in primitive lattice directions. This is a topological description of the slope; a curve on the punctured torus must avoid the puncture. Fraction labels record these directions. The resulting links among reduced fractions, lattice segments, continued fractions, quadratic forms, matrices, and hyperbolic surfaces also connect with Penner's decorated Teichm\"uller theory \cite{penner1987} and cluster algebras.

\subsection{Cluster Algebras and Generalized Markov Numbers}\index{Cluster algebra}
Cluster algebras were introduced by Fomin--Zelevinsky \cite{fzi,fziv}, and their connections with higher Teichm\"uller theory and related geometric structures were developed by Fock--Goncharov \cite{fg06,fg09}. This theory has had a major influence on the theory of Markov numbers. A cluster algebra is generated from collections of elements called clusters, whose entries are called cluster variables, together with exchange operations called mutations.

For cluster algebras associated with marked surfaces, Fomin--Shapiro--Thurston \cite{fst} developed a description using tagged arcs and tagged triangulations: cluster variables correspond to tagged arcs, clusters to tagged triangulations, and mutations to flips. For the once-punctured torus, one uses the component in which all tags are plain, so ordinary arcs and ideal triangulations suffice. Fomin--Thurston \cite{ft} realized cluster variables as suitably renormalized $\lambda$-lengths on decorated Teichm\"uller spaces, with coefficients encoded by laminations. Thus changes of triangulations correspond to transformations of variables in an algebra.

When the surface is the once-punctured torus, this framework is directly related to the classical theory of Markov numbers. In the Markov cluster algebra associated with this surface, suitable specializations of cluster variables give Markov numbers, and the three variables in one cluster give a solution of the Markov equation. Fraction labels, the combinatorics of lattice segments, Cohn matrices, and closed curves on the once-punctured torus are all organized under the common language of cluster algebras. This direction was broadened by the generalized cluster algebras introduced by Chekhov--Shapiro \cite{chsh14}. Generalized cluster algebras form a wider class of algebras containing ordinary cluster algebras, obtained by generalizing the exchange rules used in mutation. A natural question is then how much of the symmetry and good combinatorics of the classical Markov cluster algebra remains in this generalized setting.

The \emph{generalized Markov numbers}, introduced by Gyoda--Matsushita \cite{gyomatsu}, arose from this question. Since classical Markov numbers are connected, through the Markov cluster algebra, with reduced fractions, lattice segments, and curves on the torus, it is natural to ask whether Chekhov--Shapiro's generalized cluster algebras contain a well-behaved class with properties analogous to those of the classical Markov cluster algebra. The resulting equation is the following extension of the classical Markov equation:
\[
 x^2+y^2+z^2+k_1yz+k_2zx+k_3xy=(3+k_1+k_2+k_3)xyz.
\]
Here $k_1,k_2,k_3$ are nonnegative integers, and the integers appearing in positive integer solutions of this equation are called \emph{$(k_1,k_2,k_3)$-generalized Markov numbers}. This equation is not only a formal deformation of the classical Markov equation. It appears naturally when one tries, inside generalized cluster algebras, to preserve the symmetries and mutation-based generation mechanism familiar from the classical theory. Thus generalized Markov numbers give a way to reinterpret the structure behind the classical theory in a wider setting.

Subsequent work has reconstructed many aspects of the classical theory for generalized Markov numbers. In the equal-parameter case $k_1=k_2=k_3=k$, Gyoda--Maruyama \cite{gyo-maru} introduced generalized Cohn matrices, and Gyoda--Maruyama--Sato \cite{gyoda-maruyama-sato} developed related matrix, geometric, and combinatorial descriptions, including continued-fraction formulas indexed by reduced fractions. For arbitrary nonnegative parameters $k_1,k_2,k_3$, the author's paper \cite{gyoda-generalized} connected generalized Markov numbers with the Lagrange and Markov spectra. More precisely, it contains the following result.
\begin{theorem}
Let $m$ be a $(k_1,k_2,k_3)$-generalized Markov number, and suppose that it appears as the $i$-th component of a positive integer solution of the $(k_1,k_2,k_3)$-generalized Markov equation. We set 
\[\Delta(k_1,k_2,k_3,m,i):=((3+k_1+k_2+k_3)m-k_i)^2-4.\]
Then 
\[
\frac{\sqrt{\Delta(k_1,k_2,k_3,m,i)}}{m}\in \mathcal L.
\]
In particular, if
\[
\mathcal M_{k_1,k_2,k_3}:=\left\{\frac{\sqrt{\Delta(k_1,k_2,k_3,m,i)}}{m}\ \middle |\ \begin{aligned}
&\text{$m$ is a $(k_1,k_2,k_3)$-generalized Markov number}\\
&\text{appearing as the $i$-th component of a positive integer}\\
&\text{solution of the $(k_1,k_2,k_3)$-generalized Markov equation}
\end{aligned}\right\},
\]
then $\mathcal M_{k_1,k_2,k_3}\subset \mathcal L$.
\end{theorem}
Taking $(k_1,k_2,k_3)=(0,0,0)$ recovers the direction of Markov's theorem asserting that the classical discrete Markov values below $3$ occur as Lagrange constants. The irrational numbers realizing these values as Lagrange constants, and the binary quadratic forms realizing them as Markov constants, can also be given explicitly using simple closed curves on the once-punctured torus. The proof does not proceed by a direct generalization of the classical proof of Markov's theorem; rather, it uses a cluster-algebraic reinterpretation of the combinatorics of Markov numbers.

These facts indicate that generalized Markov numbers fit naturally into the arithmetic, geometric, and combinatorial structures already present in the classical theory. In this sense, they form a natural extension of classical Markov numbers.

\section{Organization of the Text}
The text is organized as follows. Part I develops the classical theory of the Lagrange and Markov spectra. Part II introduces generalized Markov numbers, uses them to construct generalized discrete Markov spectra, and relates the classical Markov theorem to this generalized framework.

\begin{figure}[tp]
\centering
\begin{tikzpicture}[
  route box/.style={draw=black!55,fill=black!3,rounded corners=2pt,
    line width=.55pt,text width=4.2cm,minimum height=1.35cm,
    align=center,inner sep=4pt,font=\small},
  route arrow/.style={-{Stealth[length=2mm,width=1.4mm]},
    draw=black!70,line width=.7pt},
  route goal/.style={anchor=west,font=\small\bfseries,inner sep=0pt},
  route note/.style={anchor=north west,text width=15.5cm,
    align=left,font=\small,inner sep=0pt}
]
\node[route goal] at (0,0) {1. Foundations of the Lagrange and Markov spectra};
\node[route box] (a1) at (2.3,-1)
  {\textbf{Chapter~\ref{chap:continued-fraction}}\\Continued fractions};
\node[route box] (a2) at (7.8,-1)
  {\textbf{Chapter~\ref{chap:lagrange-spectrum}}\\Lagrange spectrum};
\node[route box] (a3) at (13.3,-1)
  {\textbf{Chapter~\ref{chap:markov-spectrum}}\\Markov spectrum};
\draw[route arrow] (a1.east) -- (a2.west);
\draw[route arrow] (a2.east) -- (a3.west);
\node[route note] at (0,-1.9)
  {This route covers definitions, examples, bi-infinite sequence descriptions,
   and the relation between the two spectra.};

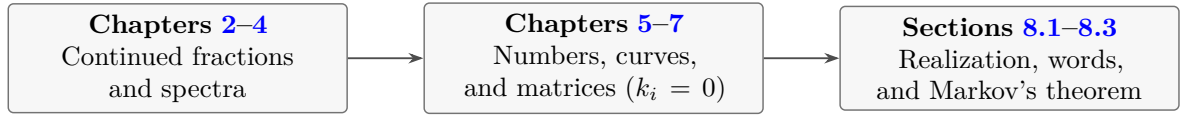
\node[route goal] at (0,-3.5) {2. The classical Markov theorem};
\node[route box] (b1) at (2.3,-4.5)
  {\textbf{Chapters~\ref{chap:continued-fraction}--\ref{chap:markov-spectrum}}\\
   Continued fractions\\and spectra};
\node[route box] (b2) at (7.8,-4.5)
  {\textbf{Chapters~\ref{chap:generalized-markov-equations-numbers}--\ref{chap:generalized-cohn-matrices}}\\
   Numbers, curves,\\and matrices ($k_i=0$)};
\node[route box] (b3) at (13.3,-4.5)
  {\textbf{Sections~\ref{sec:generalized-discrete-definition}--\ref{sec:markov-theorem}}\\
   Realization, words,\\and Markov's theorem};
\draw[route arrow] (b1.east) -- (b2.west);
\draw[route arrow] (b2.east) -- (b3.west);
\node[route note] at (0,-5.4)
  {Read Chapters~\ref{chap:generalized-markov-equations-numbers}--\ref{chap:generalized-cohn-matrices}
   and Section~\ref{sec:generalized-discrete-definition} with $(k_1,k_2,k_3)=(0,0,0)$.
   This route gives the full proof. For an overview, begin with the statement
   of Theorem~\ref{thm:markov}.};

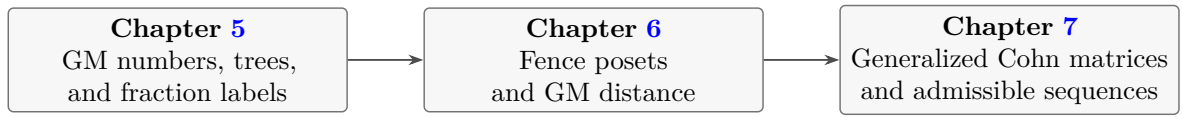
\node[route goal] at (0,-7) {3. Generalized Markov numbers and their combinatorics};
\node[route box] (c1) at (2.3,-8)
  {\textbf{Chapter~\ref{chap:generalized-markov-equations-numbers}}\\
   GM numbers, trees,\\and fraction labels};
\node[route box] (c2) at (7.8,-8)
  {\textbf{Chapter~\ref{chap:fence-posets-gm-distance}}\\
   Fence posets\\and GM distance};
\node[route box] (c3) at (13.3,-8)
  {\textbf{Chapter~\ref{chap:generalized-cohn-matrices}}\\
   Generalized Cohn matrices\\and admissible sequences};
\draw[route arrow] (c1.east) -- (c2.west);
\draw[route arrow] (c2.east) -- (c3.west);
\node[route note] at (0,-8.9)
  {Chapter~\ref{chap:continued-fraction} supplies the continued-fraction tools.
   For spectral applications, also read Chapters~\ref{chap:lagrange-spectrum}--\ref{chap:markov-spectrum}
   and continue to Chapter~\ref{chap:generalized-discrete-markov-spectra}.};

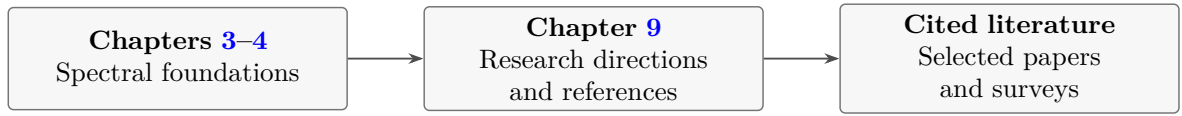
\node[route goal] at (0,-10.5) {4. Recent research on the spectra};
\node[route box] (d1) at (2.3,-11.5)
  {\textbf{Chapters~\ref{chap:lagrange-spectrum}--\ref{chap:markov-spectrum}}\\
   Spectral foundations};
\node[route box] (d2) at (7.8,-11.5)
  {\textbf{Chapter~\ref{chap:more-topics}}\\
   Research directions\\and references};
\node[route box] (d3) at (13.3,-11.5)
  {\textbf{Cited literature}\\Selected papers\\and surveys};
\draw[route arrow] (d1.east) -- (d2.west);
\draw[route arrow] (d2.east) -- (d3.west);
\node[route note] at (0,-12.4)
  {Chapter~\ref{chap:more-topics} provides an overview of selected research directions.
   For the generalized spectral results proved in this text, follow Routes~1 and~3,
   then read Chapter~\ref{chap:generalized-discrete-markov-spectra}.};
\end{tikzpicture}
\caption{Suggested reading routes.}
\label{fig:reading-routes}
\end{figure}

Figure~\ref{fig:reading-routes} suggests four reading routes, according to the reader's interests. Readers are, of course, also welcome to read the text from beginning to end. Follow the arrows from left to right, and read chapter or section ranges in numerical order. Chapters~\ref{chap:generalized-markov-equations-numbers}--\ref{chap:generalized-cohn-matrices} can also be read without first working through Chapters~\ref{chap:continued-fraction}--\ref{chap:markov-spectrum}; the continued-fraction background in Chapter~\ref{chap:continued-fraction} can be consulted as needed. A chapter-by-chapter description follows.

Chapter~\ref{chap:continued-fraction} summarizes the theory of continued fractions needed later. After reviewing reduced fractions and finite regular continued fractions, it treats infinite regular continued fractions, convergents, continued-fraction matrices, and the decomposition of irrational numbers into orbits under the unimodular group. The final section recalls Lagrange's characterization of quadratic irrationals by periodic continued fractions. This prepares the connection, used in Chapters~\ref{chap:lagrange-spectrum} and~\ref{chap:markov-spectrum}, between Lagrange and Markov constants and quadratic irrationals or binary quadratic forms.

Chapter~\ref{chap:lagrange-spectrum} deals with the Lagrange spectrum. We first define the Lagrange constant and give basic examples, and then interpret it as a limit superior of quantities obtained from convergents. By introducing the representation in terms of bi-infinite sequences, we formulate it in a way that can be compared with the Markov constant in the next chapter. For quadratic irrationals, we reduce to reduced quadratic irrationals up to $GL(2,\mathbb Z)$-equivalence and show that the Lagrange constant can be computed explicitly from the periodic part and the associated matrix. This establishes the method for computing values from periodic sequences used from Chapter~\ref{chap:generalized-markov-equations-numbers} onward.

Chapter~\ref{chap:markov-spectrum} turns to the Markov spectrum. We define the Markov constant for binary quadratic forms and organize representatives using canonical reduced binary quadratic forms and unimodular group orbits. We then express the Markov constant by bi-infinite sequences in a form parallel to the Lagrange spectrum. Finally, we show that the Markov constant of a binary quadratic form with rational coefficients coincides with the Lagrange constant of the corresponding quadratic irrational. This clarifies that quadratic irrationals, rational-coefficient binary quadratic forms, and periodic bi-infinite sequences give the same values. This will be the key point when the generalized theory is connected to spectra in the second half of the text.

Chapter~\ref{chap:generalized-markov-equations-numbers} begins Part II and introduces generalized Markov equations and generalized Markov numbers. We first give the definition and basic properties of the $(k_1,k_2,k_3)$-generalized Markov equation, and then construct the generalized Markov tree corresponding to the classical Markov tree. Through the correspondence with the Farey tree, we assign fraction labels to generalized Markov numbers and thereby organize the numbers appearing at vertices by reduced fractions. At the end of the chapter we introduce characteristic numbers, which later serve as auxiliary quantities for describing the components of generalized Cohn matrices. Thus the role of this chapter is to carry over the classical picture of Markov numbers and fraction labels to the generalized setting and to prepare the data needed for later computations.

Chapter~\ref{chap:fence-posets-gm-distance} introduces fence posets and generalized Markov distance. It relates order-ideal counts to continued fractions and transfer matrices, assigns generalized Markov lengths to curves, and defines distance by minimizing these lengths. The proof that line segments realize the distance uses sign-word reductions, minimal representatives, and local straightening of bends.

Chapter~\ref{chap:generalized-cohn-matrices} introduces generalized Cohn matrices and translates the numerical and curve-theoretic information from the preceding chapters into the language of $2\times2$ matrices. We first define the generalized Cohn tree and show that the entries of generalized Cohn matrices can be described explicitly using generalized Markov numbers and characteristic numbers. We then prove relations among characteristic numbers. Finally, by introducing generalized strongly admissible sequences, we show that generalized Cohn matrices can be expressed as products of elementary matrices. In this chapter, the arithmetic data of Chapter~\ref{chap:generalized-markov-equations-numbers} and the combinatorial-geometric data of Chapter~\ref{chap:fence-posets-gm-distance} are unified through matrix representations. In particular, it becomes clear that Cohn matrices, which play a central role in the classical theory, retain an essential role in the generalized setting.

Chapter~\ref{chap:generalized-discrete-markov-spectra} defines the generalized discrete Markov spectrum and presents the main spectral results discussed in this text. We first define a family of discrete values constructed from generalized Markov numbers, and then show that these values are realized as Lagrange constants of quadratic irrationals and as Markov constants of binary quadratic forms with rational coefficients. By specializing the general theory to $(k_1,k_2,k_3)=(0,0,0)$, we explain how the classical Markov theorem is embedded in the framework developed here. We then consider irrational-slope limits of the generalized strongly admissible sequences obtained from rational slopes, and show that a bi-infinite sequence obtained from a line of irrational slope avoiding the points of the lifted triangulation gives the boundary value $3+k_1+k_2+k_3$. We also discuss the correspondence between the $(0,0,0)$ type and the $(2,2,2)$ type, and then consider a generalization of Frobenius's uniqueness conjecture.

Chapter~\ref{chap:more-topics} collects several related directions for further reading and places the constructions of the text in a broader context.

\FloatBarrier
\section*{Acknowledgments}
The author thanks Esther Banaian for her advice during the preparation of this text. This work was supported by JSPS KAKENHI Grant Number JP25K17224.

\section*{Declaration of AI Use}
ChatGPT Pro 5.5 and 5.6, as well as GPT-6 Astra, were used to assist with checking and refining details of proofs, proofreading the text, and exploring related areas.
The author takes full responsibility for the content and accuracy of this text, including its mathematical arguments and references.

\part[Lagrange and Markov Spectra]{Lagrange and Markov Spectra}

\chapter{Continued Fractions}\label{chap:continued-fraction}
The main theme of this text is the approximation of irrational numbers by rational numbers. Continued fractions that converge to a given irrational number are indispensable for studying such approximations. This chapter collects the basic facts about continued fractions that will be used from Chapter~\ref{chap:lagrange-spectrum} onward.

Although we call them facts about continued fractions, a large part of the theory is, in effect, a theory of products of matrices in $GL(2,\mathbb Z)$. Continued-fraction calculations can be interpreted as products of such matrices. For this reason, matrix calculations in $GL(2,\mathbb Z)$ are an unavoidable tool in the modern treatment of continued fractions.

We first recall the elementary notions concerning reduced fractions. We then discuss finite regular continued-fraction expansions of rational numbers and infinite regular continued-fraction expansions of irrational numbers. In the final section we prove Lagrange's characterization of quadratic irrationals by periodic continued fractions.

The exposition and organization of this chapter are based largely on the corresponding chapters of \cite{kida}.

\section{Reduced Fractions}
We begin with the notion of a reduced fraction. Although this is familiar to many readers, we fix the precise convention used in this text.
\begin{defi}
Let $a,b\in \mathbb Z$. If there exists $k\in\mathbb Z$ such that $b=ak$, then $a$ is called a \emph{divisor} of $b$. We write this as $a\mid b$.
\end{defi}
This definition of divisibility also applies when $a$ or $b$ is zero or negative. For example, if $a\neq 0$, then $a\mid 0$ always holds, and hence every nonzero integer is a divisor of $0$. Conversely, if $b\neq 0$, then $0\mid b$ never holds, so $0$ is not a divisor of any nonzero integer. The relation $0\mid0$ also holds, since $0=0\cdot k$ for every integer $k$.
\begin{defi}
Let $a_1,\dots,a_n$ be integers that are not all zero. The \emph{greatest common divisor} $\gcd(a_1,\dots,a_n)$ is the positive integer $d$ satisfying the following two conditions:
\begin{enumerate}
  \item $d\mid a_i$ for every $i=1,\dots,n$.
  \item If an integer $c$ satisfies $c\mid a_i$ for every $i=1,\dots,n$, then $c\mid d$.
\end{enumerate}
\end{defi}
The greatest common divisor always exists.
\begin{defi}
If
\[
\gcd(a_1,a_2,\dots,a_n)=1,
\]
then $a_1,a_2,\dots,a_n$ are said to be \emph{relatively prime}.
\end{defi}
With this definition, it also makes sense to ask whether a pair involving $0$, or a pair involving both positive and negative numbers, is relatively prime.
\begin{exam}
Let us check from the definition whether $0$ is relatively prime to some small integers.
\begin{itemize}
\item The divisors of $1$ are $\pm 1$, whereas $0$ is divisible by every nonzero integer. Hence the greatest common divisor of $0$ and $1$ is $1$. Thus $0$ and $1$ are relatively prime.
\item The divisors of $2$ are $\pm 1,\pm 2$. Hence the greatest common divisor of $0$ and $2$ is $2$. Thus $0$ and $2$ are not relatively prime.
\item The divisors of $2$ are $\pm 1,\pm 2$, and the divisors of $-3$ are $\pm 1,\pm 3$. Hence the greatest common divisor of $2$ and $-3$ is $1$. Thus $2$ and $-3$ are relatively prime.
\end{itemize}
\end{exam}
We now define fractions and reducedness.
\begin{defi}
\index{Reduced fraction}Let $a,b\in\mathbb R$ and assume that $(a,b)\neq (0,0)$. The formal symbol $\frac{a}{b}$ is called a \emph{fraction}. If $a$ and $b$ are integers, if they are relatively prime, and if either $b>0$ or $(a,b)=(1,0)$, then the fraction $\frac{a}{b}$ is said to be \emph{reduced}.
\end{defi}
The symbol $a/b$ is a formal fraction. When $b\neq0$, it represents the real number $ab^{-1}$, which is rational if $a,b$ are integers. Formal fractions such as $1/0$ do not represent real numbers.

We identify a fraction with its real value whenever the denominator is nonzero. The reduced fraction representing an integer $n$ is $n/1$. The fractions $1/(-2)$ and $2/0$ are not reduced. Our convention admits $1/0$ as the unique reduced fraction with denominator zero; $-1/0$ is not reduced.

\section{Finite Regular Continued Fractions}
In this section and the next one, we review the basic properties of continued fractions. We first define finite regular continued fractions, which correspond to rational numbers, and study their elementary properties.
\begin{defi}\label{def:regular-continued-fraction}\index{Partial quotient}\index{Regular continued fraction}
Let $(a_i)_{i=0}^n=(a_0,a_1,\dots,a_n)$ be a finite sequence of real numbers. For every $1\leq k\leq n$, define recursively
\[
r_k^{(k)}:=a_k,
\qquad
r_j^{(k)}:=a_j+\frac{1}{r_{j+1}^{(k)}}
\quad (j=k-1,k-2,\dots,1),
\]
whenever these expressions are defined. We call the sequence \emph{continued-fraction admissible} if every $r_j^{(k)}$ with $1\leq j\leq k\leq n$ is defined and nonzero. Every one-term sequence $(a_0)$ is declared admissible. For an admissible sequence, set
\[
[a_0]:=a_0,
\qquad
[a_0;a_1,a_2,\dots,a_n]
:=a_0+\frac{1}{r_1^{(n)}}
\quad(n\geq1).
\]
The numbers $a_i$ are called the \emph{partial quotients}. In particular, a sequence with $a_0\in\mathbb Z$ and $a_k\in\mathbb Z_{\geq1}$ for every $1\leq k\leq n$ is automatically admissible. If, in addition, $a_n\neq1$ whenever $n\neq0$, the resulting continued fraction is called a \emph{finite regular continued fraction}.
\end{defi}
The last condition removes the ambiguity $[a_0;\dots,a_{n-1},1]=[a_0;\dots,a_{n-1}+1]$.
\begin{defi}\label{def:convergent}\index{Convergent}
Given a continued-fraction admissible sequence $(a_0,a_1,\dots,a_n)$ and an index $0\leq k\leq n$, we call $[a_0;a_1,a_2,\dots,a_k]$ its $k$-th \emph{convergent}.
\end{defi}
The following proposition computes the convergents.
\begin{prop}\label{prop:recursion}\index{Convergent}
Let $(a_0,a_1,\dots,a_n)$ be a continued-fraction admissible sequence of real numbers. Define two sequences by
\begin{alignat}{6}
  & p_0 &{}={}& a_0        &\qquad&
    p_1 &{}={}& a_0 a_1 + 1 &\qquad&
    p_k &{}={}& a_k p_{k-1}+ p_{k-2},\label{eq:p}\\
  & q_0 &{}={}& 1          &\qquad&
    q_1 &{}={}& a_1         &\qquad&
    q_k &{}={}& a_k q_{k-1}+ q_{k-2}.\label{eq:q}
\end{alignat}
Then $q_k\neq0$ and
\[
[a_0;a_1,\dots,a_k]=\frac{p_k}{q_k}
\]
for every $0\leq k\leq n$.
Moreover, for $2\leq k\leq n$, the right-hand side may be written as
\[
\frac{p_k}{q_k}=\frac{a_kp_{k-1}+p_{k-2}}{a_kq_{k-1}+q_{k-2}}.
\]
\end{prop}
\begin{proof}
We argue by induction on $k$. The cases $k=0$ and $k=1$ follow by direct calculation, and admissibility gives $q_1=a_1\neq0$. Let $k\geq2$, and assume the assertion for admissible sequences whose final index is less than $k$. Since the original sequence is admissible, $a_k\neq0$, and the shortened sequence
\[
(a_0,a_1,\dots,a_{k-2},a_{k-1}+\tfrac1{a_k})
\]
is also admissible. Let $p_i',q_i'$ be the quantities defined by the same recurrences for this shortened sequence. By the induction hypothesis, $q_{k-1}'\neq0$ and
\[
[a_0;a_1,\dots,a_{k-1}+\tfrac1{a_k}]
=\frac{p_{k-1}'}{q_{k-1}'}.
\]
The recurrences give
\[
p_{k-1}'=\frac{p_k}{a_k},
\qquad
q_{k-1}'=\frac{q_k}{a_k}.
\]
Therefore $q_k=a_kq_{k-1}'\neq0$, and
\[
[a_0;a_1,\dots,a_k]
=[a_0;a_1,\dots,a_{k-1}+\tfrac1{a_k}]
=\frac{p_{k-1}'}{q_{k-1}'}
=\frac{p_k}{q_k}.
\]
\end{proof}
Notice that this proposition does \emph{not} assume that $[a_0;a_1,a_2,\dots,a_n]$ is a finite regular continued fraction. Until Lemma~\ref{lem:det}, only continued-fraction admissibility will be assumed.

The convergents are conveniently computed using matrices.
\begin{theorem}\label{thm:continued-fraction-matrix}\index{Continued-fraction matrix}
Let $(a_0,a_1,\dots,a_n)$ be a continued-fraction admissible sequence and let $0\leq k\leq n$. Put $p_{-1}=1$ and $q_{-1}=0$, and define $p_k,q_k$ by \eqref{eq:p} and \eqref{eq:q}. Then
\begin{align}\label{eq:continued-fraction-matrix}
\begin{bmatrix}p_{k}&p_{k-1}\\q_k&q_{k-1}\end{bmatrix}
=
\begin{bmatrix}a_0&1\\1&0\end{bmatrix}
\begin{bmatrix}a_1&1\\1&0\end{bmatrix}
\cdots
\begin{bmatrix}a_k&1\\1&0\end{bmatrix}.
\end{align}
\end{theorem}
\begin{proof}
For $k=0$, the identity
\[
\begin{bmatrix}p_{0}&p_{-1}\\q_0&q_{-1}\end{bmatrix}
=
\begin{bmatrix}a_0&1\\1&0\end{bmatrix}
\]
is immediate from the definition. Let $k\geq 1$, and assume that the theorem has been proved up to $k-1$. Then
\[
\begin{bmatrix}a_0&1\\1&0\end{bmatrix}
\begin{bmatrix}a_1&1\\1&0\end{bmatrix}
\cdots
\begin{bmatrix}a_k&1\\1&0\end{bmatrix}
=
\begin{bmatrix}p_{k-1}&p_{k-2}\\q_{k-1}&q_{k-2}\end{bmatrix}
\begin{bmatrix}a_k&1\\1&0\end{bmatrix}
=
\begin{bmatrix}
a_kp_{k-1}+p_{k-2}&p_{k-1}\\
a_kq_{k-1}+q_{k-2}&q_{k-1}
\end{bmatrix}
=
\begin{bmatrix}p_k&p_{k-1}\\q_k&q_{k-1}\end{bmatrix}.
\]
Thus the formula holds for $k$ as well.
\end{proof}
The matrix \(\begin{bsmallmatrix}p_{k}&p_{k-1}\\ q_k&q_{k-1}\end{bsmallmatrix}\) is the \emph{continued-fraction matrix} of $[a_0;a_1,\dots,a_k]$.

The next lemma follows immediately from the matrix formula. In what follows, unless otherwise stated, $p_k$ and $q_k$ are used in the sense of Theorem~\ref{thm:continued-fraction-matrix}.
\begin{lemm}\label{lem:det}
For a continued-fraction admissible sequence $(a_0,a_1,\dots,a_n)$ and every $0\leq k\leq n$, one has
\begin{equation}\label{eq:det}
p_kq_{k-1}-q_kp_{k-1}=(-1)^{k+1}.
\end{equation}
\end{lemm}
\begin{proof}
Take determinants on both sides of \eqref{eq:continued-fraction-matrix}.
\end{proof}
We now derive several consequences for finite regular continued fractions.
\begin{coro}\label{cor:reduced-fraction}\index{Reduced fraction}
Let $[a_0;a_1,\dots,a_n]$ satisfy the conditions for a finite regular continued fraction, except that we also allow $a_n=1$. Then $\frac{p_k}{q_k}$ is a reduced fraction for every $0\leq k\leq n$.
\end{coro}
\begin{proof}
The case $n=0$ is clear, so assume $n\neq 0$. We have $a_0\in\mathbb Z$ and $a_k\in\mathbb Z_{\geq 1}$ for $1\leq k\leq n$. It is immediate from the recurrence that $q_k>0$. Let $d_k$ be the greatest common divisor of $p_k$ and $q_k$, and write $p_k=d_kp'_k$ and $q_k=d_kq'_k$. Then $p_kq_{k-1}-q_kp_{k-1}=d_k(p'_kq_{k-1}-q'_kp_{k-1})=(-1)^{k+1}$. Since $d_k\in\mathbb Z_{\geq 1}$ and $p'_kq_{k-1}-q'_kp_{k-1}\in\mathbb Z$, we must have $d_k=1$. Hence $\frac{p_k}{q_k}$ is reduced.
\end{proof}
\begin{coro}\label{cor:qk>k}
Let $[a_0;a_1,\dots,a_n]$ satisfy the conditions for a finite regular continued fraction, except that we also allow $a_n=1$. Then the sequence $(q_1,\dots,q_n)$ is strictly increasing, and $q_k\geq k$ for every $1\leq k\leq n$.
\end{coro}
\begin{proof}
We prove $q_k\geq k$ by induction. For $k=1$ this is clear from the definition. Let $k\geq 2$ and assume $q_{k-1}\geq k-1$. Since $a_k\geq 1$, we obtain $q_k=a_kq_{k-1}+q_{k-2}\geq q_{k-1}+1\geq k$. The same inequality also shows that the sequence is strictly increasing.
\end{proof}
Lemma~\ref{lem:det} can be rewritten as follows.
\begin{coro}\label{cor:pk/qk-pk-1/qk-1}
Let $[a_0;a_1,\dots,a_n]$ satisfy the conditions for a finite regular continued fraction, except that we also allow $a_n=1$. Then, for every $1\leq k\leq n$,
\begin{equation}\label{eq:pk/qk-pk-1/qk-1}
\frac{p_k}{q_{k}}-\frac{p_{k-1}}{q_{k-1}}
=
\frac{(-1)^{k+1}}{q_kq_{k-1}}.
\end{equation}
\end{coro}
\begin{proof}
Divide both sides of \eqref{eq:det} by $q_kq_{k-1}$.
\end{proof}
\begin{coro}\label{cor:size-pk/qk}
Let $[a_0;a_1,\dots,a_n]$ satisfy the conditions for a finite regular continued fraction, except that we also allow $a_n=1$. Then, for all $k\in\mathbb Z_{\geq 0}$ and $\ell\in\mathbb Z_{\geq 1}$ satisfying $1\leq 2k+2\ell+1\leq n$, one has
\[
\frac{p_{2k}}{q_{2k}}
<
\frac{p_{2k+2\ell}}{q_{2k+2\ell}}
<
\frac{p_{2k+2\ell+1}}{q_{2k+2\ell+1}}
<
\frac{p_{2k+1}}{q_{2k+1}}.
\]
\end{coro}
\begin{proof}
We first prove
\[
\frac{p_{2k}}{q_{2k}}<
\frac{p_{2k+2\ell}}{q_{2k+2\ell}}
\quad\text{and}\quad
\frac{p_{2k+2\ell+1}}{q_{2k+2\ell+1}}<
\frac{p_{2k+1}}{q_{2k+1}}.
\]
For $2\leq m\leq n$, we compute
\begin{align*}
\frac{p_m}{q_m}-\frac{p_{m-2}}{q_{m-2}}
&=\frac{p_mq_{m-2}-p_{m-2}q_m}{q_mq_{m-2}}\\
&=\frac{(a_mp_{m-1}+p_{m-2})q_{m-2}
-p_{m-2}(a_mq_{m-1}+q_{m-2})}{q_mq_{m-2}}\\
&=\frac{a_m(p_{m-1}q_{m-2}-q_{m-1}p_{m-2})}{q_mq_{m-2}}
=\frac{(-1)^m a_m}{q_mq_{m-2}},
\end{align*}
where the last equality follows from \eqref{eq:det}. Since $a_m>0$, taking $m=2k+2$ gives $\frac{p_{2k}}{q_{2k}}<\frac{p_{2k+2}}{q_{2k+2}}$, and taking $m=2k+3$ gives $\frac{p_{2k+3}}{q_{2k+3}}<\frac{p_{2k+1}}{q_{2k+1}}$. These are the desired inequalities for $\ell=1$. The inequalities for general $\ell$ follow by repeated application and transitivity.

It remains to prove
\[
\frac{p_{2k+2\ell}}{q_{2k+2\ell}}
<
\frac{p_{2k+2\ell+1}}{q_{2k+2\ell+1}}.
\]
Substituting $2k+2\ell+1$ for $k$ in \eqref{eq:pk/qk-pk-1/qk-1}, we obtain
\[
\frac{p_{2k+2\ell+1}}{q_{2k+2\ell+1}}
-\frac{p_{2k+2\ell}}{q_{2k+2\ell}}
=
\frac{(-1)^{2k+2\ell+2}}{q_{2k+2\ell+1}q_{2k+2\ell}}>0.
\]
This proves the assertion.
\end{proof}
The indexing in Corollary~\ref{cor:size-pk/qk} may obscure the simple meaning of the statement. Applying it for all possible $k$ gives, for example,
\[
\frac{p_0}{q_0}<\frac{p_2}{q_2}<\frac{p_4}{q_4}<\cdots<
\frac{p_5}{q_5}<\frac{p_3}{q_3}<\frac{p_1}{q_1},
\]
which may be easier to visualize.

We finish this section by proving the bijective correspondence between rational numbers and finite regular continued fractions. We first define the set of sequences that represent finite regular continued fractions. Since finite sequences have varying lengths, we regard them as infinite sequences that become zero from some point on. Define
\[
\mathscr Z:=
\{(a_k)_{k=0}^\infty\mid a_0\in \mathbb Z,\ a_i=0\text{ for every }i\in\mathbb Z_{\geq 1}\}.
\]
The map $z\colon\mathscr Z\to\mathbb Z$ defined by $z((a_k)_{k=0}^\infty)=a_0$ is clearly a bijection. Thus $z$ gives a bijective correspondence between integers and their finite regular continued-fraction expansions.

Next consider rational numbers that are not integers. Define
\[
\mathscr Q:=
\left\{(a_k)_{k=0}^\infty\ \middle|\
\begin{array}{l}
a_0\in \mathbb Z,\ \exists n\in\mathbb Z_{\geq 1}\text{ such that }
 a_1,\dots,a_{n-1}\in\mathbb Z_{\geq 1},\\
a_n\in\mathbb Z_{\geq 2},\text{ and }a_i=0\text{ for all }i\in\mathbb Z_{\geq n+1}
\end{array}
\right\}.
\]
Define $f\colon\mathscr Q\to\mathbb Q\setminus\mathbb Z$ by
\[
f((a_k)_{k=0}^\infty)=[a_0;a_1,\dots,a_n],
\]
where $n$ is the integer such that $a_n\in\mathbb Z_{\geq 2}$ and $a_i=0$ for all $i\geq n+1$. It is not immediate that this map is bijective, so we construct its inverse. We begin with the following lemma.
\begin{lemm}\label{lem:sequence-stop}
Let $\alpha\in\mathbb Q\setminus\mathbb Z$. Construct $\alpha_k$ and $a_k$ recursively by
\begin{equation}\label{eq:rational-to-sequence}
\alpha_0=\alpha,\qquad
a_k=\lfloor\alpha_k\rfloor,\qquad
\alpha_{k+1}=\frac{1}{\alpha_k-a_k}.
\end{equation}
Then there always exists $n\in\mathbb Z_{>0}$ such that $a_n=\alpha_n$. Thus the process stops at that point.
\end{lemm}
\begin{proof}
Write $\alpha=\alpha_0=\frac{r_0}{s_0}$ as a reduced fraction. Since $\alpha\in\mathbb Q\setminus\mathbb Z$, we have $s_0\geq 2$. As long as $a_{k-1}\neq\alpha_{k-1}$, the number $\alpha_k$ is rational; write it as a reduced fraction $\alpha_k=\frac{r_k}{s_k}$. If $s_k=1$, then $\alpha_k\in\mathbb Z$, hence $a_k=\alpha_k$, and we may take $n=k$.

Assume $s_k\geq 2$. Then
\[
\alpha_k-a_k=\frac{r_k}{s_k}-a_k=\frac{r_k-a_ks_k}{s_k}.
\]
Since $\alpha_k-a_k\neq 0$ and $0<\alpha_k-a_k<1$, we have $0<r_k-a_ks_k<s_k$. By definition,
\[
\alpha_{k+1}=\frac{s_k}{r_k-a_ks_k}.
\]
If this is written as the reduced fraction $\frac{r_{k+1}}{s_{k+1}}$, then $s_{k+1}$ divides $r_k-a_ks_k$. Hence
\[
s_{k+1}\leq r_k-a_ks_k<s_k.
\]
Thus, as long as $s_k\geq 2$, the denominators strictly decrease. Since $s_k\geq 1$, there must be some $n$ such that $s_n=1$. Then $\alpha_n\in\mathbb Z$, so $a_n=\alpha_n$.
\end{proof}
\begin{rema}
The procedure in Lemma~\ref{lem:sequence-stop} is exactly the Euclidean algorithm. In the next section we will carry out the analogous procedure for infinite continued fractions.
\end{rema}
\begin{theorem}\label{thm:well-def-inverse}
Let $\alpha\in\mathbb Q\setminus\mathbb Z$. Use \eqref{eq:rational-to-sequence} and let $n$ be the smallest index such that $a_n=\alpha_n$. Consider the finite sequence $(a_k)_{k=0}^n$, and extend it by putting $a_i=0$ for all $i\geq n+1$. Then the resulting infinite sequence $(a_k)_{k=0}^\infty$ belongs to $\mathscr Q$.
\end{theorem}
\begin{proof}
For every $k\geq 0$, the equality $a_k=\lfloor\alpha_k\rfloor$ gives $a_k\in\mathbb Z$. Since $\alpha\notin\mathbb Z$, we have $n\neq 0$. It remains to prove that $a_1,\dots,a_{n-1}\geq 1$ and $a_n\geq 2$. For every $0\leq k\leq n-1$, we have $\alpha_k-a_k\neq 0$, hence $0<\alpha_k-a_k<1$. Therefore
\[
\alpha_{k+1}=\frac{1}{\alpha_k-a_k}>1,
\]
and so $a_{k+1}\geq 1$. Thus $a_1,\dots,a_n\geq 1$.

It remains to show $a_n\geq 2$. Suppose $a_n=1$. Since $a_n=\alpha_n$, we would have $\alpha_n=1$. Then \eqref{eq:rational-to-sequence} gives
\[
\alpha_{n-1}=a_{n-1}+\frac{1}{\alpha_n}=a_{n-1}+1,
\]
so $\alpha_{n-1}\in\mathbb Z$. Hence $\alpha_{n-1}=a_{n-1}$, contradicting the minimality of $n$. Therefore $a_n\geq 2$.
\end{proof}
Theorem~\ref{thm:well-def-inverse} says that the correspondence $g\colon\mathbb Q\setminus\mathbb Z\to\mathscr Q$ defined by $g(\alpha)=(a_k)_{k=0}^\infty$ is well-defined. We now prove that $f$ and $g$ are inverse maps.
\begin{theorem}\label{thm:bijective-rational}
The maps $f$ and $g$ are inverse to each other. In particular, $f$ is a bijection.
\end{theorem}
Before proving this theorem, we record a lemma.
\begin{lemm}\label{lem:integer-part-rational}
For a finite regular continued fraction $[a_0;a_1,\dots,a_n]$, one has
\[
a_0=\lfloor [a_0;a_1,\dots,a_n]\rfloor.
\]
\end{lemm}
\begin{proof}
If $n=0$, then the assertion is clear. Assume $n\geq 1$. It suffices to prove
\[
a_0\leq [a_0;a_1,\dots,a_n]<a_0+1,
\quad\text{or equivalently}\quad
0<[0;a_1,\dots,a_n]<1.
\]
If $n=1$, then the last partial quotient satisfies $a_1\geq 2$, and hence
\[
[0;a_1]=\frac{1}{a_1}<1.
\]
If $n\geq 2$, then
\[
[0;a_1,\dots,a_n]
=\frac{1}{a_1+\frac{1}{[a_2;\dots,a_n]}}
<\frac{1}{a_1}\leq 1.
\]
Positivity is clear, and the claim follows.
\end{proof}
\begin{proof}[Proof of Theorem~\ref{thm:bijective-rational}]
We first show $f\circ g=\mathrm{id}_{\mathbb Q\setminus\mathbb Z}$. It suffices to prove the following: for $\alpha\in\mathbb Q\setminus\mathbb Z$, if $(\alpha_k)_{k=0}^n$ and $(a_k)_{k=0}^n$ are constructed by \eqref{eq:rational-to-sequence}, then for every $0\leq k\leq n-1$,
\[
\alpha=[a_0;a_1,\dots,a_k,\alpha_{k+1}]=[a_0;a_1,\dots,a_n].
\]
We prove the first equality. For $k=0$ it follows by direct computation. Suppose $k\geq 1$ and assume
\[
\alpha=[a_0;a_1,\dots,a_{k-1},\alpha_k].
\]
Then
\[
[a_0;a_1,\dots,a_{k-1},a_k,\alpha_{k+1}]
=
[a_0;a_1,\dots,a_{k-1},a_k+\alpha_k-a_k]
=
[a_0;a_1,\dots,a_{k-1},\alpha_k]
=
\alpha.
\]
Taking $k=n-1$ and using $a_n=\alpha_n$ gives the second equality. Hence $f\circ g=\mathrm{id}_{\mathbb Q\setminus\mathbb Z}$.

Next we show $g\circ f=\mathrm{id}_{\mathscr Q}$. Take $(b_k)_{k=0}^\infty\in\mathscr Q$, and write
\[
f((b_k)_{k=0}^\infty)=[b_0;b_1,\dots,b_n].
\]
Let $(a_k)_{k=0}^\infty$ be the sequence constructed from this rational number by \eqref{eq:rational-to-sequence}, with undefined terms completed by zero after the process stops. We must show that $(a_k)_{k=0}^\infty=(b_k)_{k=0}^\infty$. By Lemma~\ref{lem:integer-part-rational}, we first obtain $a_0=b_0$. Hence
\[
\alpha_1
=
\frac{1}{[b_0;b_1,\dots,b_n]-b_0}
=
\frac{1}{[0;b_1,\dots,b_n]}
=
[b_1;b_2,\dots,b_n].
\]
Applying Lemma~\ref{lem:integer-part-rational} again gives $a_1=b_1$. Repeating the same argument shows that $(a_k)_{k=0}^n=(b_k)_{k=0}^n$. Finally, since
\[
\alpha_n=[b_n]=b_n=a_n,
\]
the algorithm stops at this point, and hence $a_i=0$ for all $i\geq n+1$. Thus $(a_k)_{k=0}^\infty=(b_k)_{k=0}^\infty$.
\end{proof}
Since the sequences in $\mathscr Z\cup\mathscr Q$ correspond to finite regular continued fractions, and since $z\colon\mathscr Z\to\mathbb Z$ and $f\colon\mathscr Q\to\mathbb Q\setminus\mathbb Z$ are bijections, we obtain the following theorem.
\begin{theorem}
For every rational number $\alpha$, there exists a unique finite regular continued fraction whose value is $\alpha$. The sequence of partial quotients is obtained by the algorithm \eqref{eq:rational-to-sequence}.
\end{theorem}
Although Lemma~\ref{lem:sequence-stop} was stated for nonintegral rational numbers, the same procedure includes the integral case as a special case.
\begin{defi}
The finite regular continued fraction whose value is a rational number $\alpha$ is called the \emph{finite regular continued-fraction expansion} of $\alpha$.
\end{defi}
\begin{rema}\label{rem:two-finite-expansions}
Let
\[
\theta=[b_0;b_1,\dots,b_k]\qquad(k\geq 1)
\]
be the finite regular continued-fraction expansion of a rational number $\theta$. Suppose that a finite continued fraction $[c_0;c_1,\dots,c_j]$, where $c_0\in\mathbb Z$ and $c_1,\dots,c_j\in\mathbb Z_{\geq 1}$, also has value $\theta$. Then its sequence of partial quotients is either
\[
(b_0,b_1,\dots,b_k)
\quad\text{or}\quad
(b_0,b_1,\dots,b_{k-1},b_k-1,1).
\]
Indeed, if $c_j\geq 2$, this follows from the uniqueness of the finite regular continued-fraction expansion. If $c_j=1$, first absorb the terminal $1$ by
\[
[c_0;c_1,\dots,c_{j-1},1]
=[c_0;c_1,\dots,c_{j-1}+1],
\]
and then apply uniqueness.
\end{rema}

\section{Infinite Continued Fraction Expansions of Irrational Numbers}
In the preceding section, we saw that finite regular continued fractions are in bijection with rational numbers. We now consider infinite sequences of integers whose entries after the first are positive. Taking limits of their finite continued fractions gives a representation of irrational numbers. We justify this construction in this section. The strategy is almost the same as in the rational case, but limits make some parts slightly more delicate.

Consider the set
\[
\mathscr S:=
\{(a_k)_{k=0}^{\infty}\mid a_0\in \mathbb Z,\ a_k\in \mathbb Z_{\geq 1}\ (k\geq 1)\}.
\]
We define the infinite analogue of finite regular continued fractions as follows.
\begin{defi}\label{def:infinite-regular-continued-fraction}\index{Regular continued fraction}
For $(a_k)_{k=0}^{\infty}\in\mathscr S$, consider
\[
\lim_{n\to\infty}[a_0;a_1,a_2,\dots,a_n].
\]
This limit, when it exists, is called an \emph{infinite regular continued fraction}.
\end{defi}
We also write an infinite regular continued fraction as
\[
[a_0;a_1,a_2,\dots],
\qquad
a_0+\dfrac{1}{a_1+\dfrac{1}{a_2+\dfrac{1}{\ddots}}},
\]
in order to display the underlying infinite sequence of partial quotients.

Since the definition involves a limit, it is not a priori clear that an infinite regular continued fraction has a real value. The next theorem shows that it does, and that the value is irrational.
\begin{theorem}\label{thm:limit-existence}
For every $(a_k)_{k=0}^\infty\in\mathscr S$, the limit
\[
\lim_{n\to\infty}[a_0;a_1,a_2,\dots,a_n]
\]
exists and is irrational.
\end{theorem}
\begin{proof}
Put
\[
b_k:=[a_0;a_1,\dots,a_k]=\frac{p_k}{q_k}.
\]
Extending Corollary~\ref{cor:size-pk/qk} to the infinite sequence, we see that the sequence $\{b_{2k+1}\}_{k=0}^\infty$ is decreasing and bounded, and that $\{b_{2k}\}_{k=0}^\infty$ is increasing and bounded. Hence the limits
\[
\lim_{k\to\infty}b_{2k+1}
\quad\text{and}\quad
\lim_{k\to\infty}b_{2k}
\]
exist as real numbers. Denote them by $\alpha_1$ and $\alpha_2$, respectively. We prove that $\alpha_1=\alpha_2$. By Corollaries~\ref{cor:pk/qk-pk-1/qk-1} and~\ref{cor:qk>k}, for $k\geq 1$ we have
\[
0<b_{2k+1}-b_{2k}
=
\frac{(-1)^{2k+2}}{q_{2k+1}q_{2k}}
\leq
\frac{1}{2k(2k+1)}.
\]
Letting $k\to\infty$, we obtain $\lim_{k\to\infty}(b_{2k+1}-b_{2k})=0$. Hence
\[
\alpha_1-\alpha_2
=
\lim_{k\to\infty}b_{2k+1}-\lim_{k\to\infty}b_{2k}
=
\lim_{k\to\infty}(b_{2k+1}-b_{2k})
=
0.
\]
Thus the two limits coincide. Let $\alpha:=\alpha_1=\alpha_2$.

We next prove that $\alpha$ is irrational. For every $k\geq 0$ we have
\[
b_{2k}<\alpha<b_{2k+1}.
\]
Therefore
\[
0<\alpha-b_{2k}<b_{2k+1}-b_{2k}
=\frac{1}{q_{2k}q_{2k+1}}.
\]
Since $b_{2k}=\frac{p_{2k}}{q_{2k}}$, multiplying by $q_{2k}$ gives
\[
0<\alpha q_{2k}-p_{2k}<\frac{1}{q_{2k+1}}.
\]
Suppose, for a contradiction, that $\alpha=\frac{a}{b}$ is rational, where $\frac{a}{b}$ is reduced and $b>0$. Then
\[
0<aq_{2k}-bp_{2k}<\frac{b}{q_{2k+1}}.
\]
The middle term is an integer for every $k$, while $\frac{b}{q_{2k+1}}<1$ for all sufficiently large $k$ by Corollary~\ref{cor:qk>k}. This is impossible. Hence $\alpha$ is irrational.
\end{proof}
The theorem implies that the correspondence
\[
F\colon\mathscr S\to\mathbb R\setminus\mathbb Q,\qquad
F((a_k)_{k=0}^\infty)=\lim_{n\to\infty}[a_0;a_1,\dots,a_n],
\]
is well-defined. We now construct its inverse. First we prove the following theorem.
\begin{theorem}\label{thm:irrational-to-sequence}\index{Complete quotient}\index{Gauss map}
Let $\alpha\in\mathbb R\setminus\mathbb Q$. Construct $(\alpha_k)_{k=0}^\infty$ and $(a_k)_{k=0}^\infty$ by
\begin{equation}\label{eq:irrational-to-sequence}
\alpha_0=\alpha,\qquad
a_k=\lfloor\alpha_k\rfloor,\qquad
\alpha_{k+1}=\frac{1}{\alpha_k-a_k}.
\end{equation}
Then $(a_k)_{k=0}^\infty\in\mathscr S$.
\end{theorem}
\begin{proof}
For every $k\geq 0$, $a_k=\lfloor\alpha_k\rfloor$, so $a_k\in\mathbb Z$. Moreover, every $\alpha_k$ is irrational. Indeed, if some $\alpha_k$ were rational, then
\[
\alpha_{k-1}=a_{k-1}+\frac{1}{\alpha_k}
\]
would also be rational, and repeating this argument would imply that $\alpha_0$ is rational, a contradiction. Hence, for every $k\geq 1$, we have
\[
0<\alpha_{k-1}-a_{k-1}<1,
\]
and therefore $\alpha_k>1$. It follows that $a_k\geq 1$ for every $k\geq 1$. Thus $(a_k)_{k=0}^\infty\in\mathscr S$.
\end{proof}
Let $G\colon\mathbb R\setminus\mathbb Q\to\mathscr S$ be the correspondence in this theorem; that is, $G(\alpha)=(a_n)_{n=0}^\infty$, where the sequence is constructed by \eqref{eq:irrational-to-sequence}. Theorem~\ref{thm:irrational-to-sequence} says that $G$ is well-defined. The main result of this section is the following.
\begin{theorem}\label{thm:bijective}\index{Regular continued fraction}
The maps $F$ and $G$ are inverse to each other. In particular, $F$ is a bijection.
\end{theorem}
Before proving this theorem, we record a lemma.
\begin{lemm}\label{lem:integer-part}
For an infinite regular continued fraction, one has
\[
a_0=
\left\lfloor
\lim_{n\to\infty}[a_0;a_1,\dots,a_n]
\right\rfloor.
\]
\end{lemm}
\begin{proof}
It suffices to prove
\[
a_0<
\lim_{n\to\infty}[a_0;a_1,\dots,a_n]
<a_0+1,
\quad\text{or equivalently}\quad
0<
\lim_{n\to\infty}[0;a_1,\dots,a_n]
<1.
\]
Apply the ordering of even and odd convergents, proved in Theorem~\ref{thm:limit-existence}, to the infinite continued fraction $[0;a_1,a_2,\dots]$. Its value lies strictly between its first two convergents. Hence
\[
0
<
\lim_{n\to\infty}[0;a_1,\dots,a_n]
<
\frac{1}{a_1}
\leq 1,
\]
and the desired inequality follows.
\end{proof}
\begin{proof}[Proof of Theorem~\ref{thm:bijective}]
We first show $F\circ G=\mathrm{id}_{\mathbb R\setminus\mathbb Q}$. It suffices to prove the following: for $\alpha\in\mathbb R\setminus\mathbb Q$, construct $(\alpha_k)_{k=0}^\infty$ and $(a_k)_{k=0}^\infty$ by \eqref{eq:irrational-to-sequence}. Then, for every $k\in\mathbb Z_{\geq 0}$,
\[
\alpha=[a_0;a_1,\dots,a_k,\alpha_{k+1}]
=
\lim_{n\to\infty}[a_0;a_1,\dots,a_n].
\]
We first prove the first equality. For $k=0$ it follows by direct computation. Suppose $k\geq 1$ and assume
\[
\alpha=[a_0;a_1,\dots,a_{k-1},\alpha_k].
\]
Then
\[
[a_0;a_1,\dots,a_{k-1},a_k,\alpha_{k+1}]
=
[a_0;a_1,\dots,a_{k-1},a_k+\alpha_k-a_k]
=
[a_0;a_1,\dots,a_{k-1},\alpha_k]
=
\alpha.
\]
We now prove
\[
\alpha=\lim_{n\to\infty}[a_0;a_1,\dots,a_n].
\]
By the equality just proved and Proposition~\ref{prop:recursion}, for every $n\geq 1$ we have
\[
\alpha=[a_0;a_1,\dots,a_n,\alpha_{n+1}]
=
\frac{\alpha_{n+1}p_n+p_{n-1}}{\alpha_{n+1}q_n+q_{n-1}}.
\]
Therefore
\begin{align}\label{eq:difference-approximation}
\alpha-\frac{p_n}{q_n}
&=
\frac{\alpha_{n+1}p_n+p_{n-1}}{\alpha_{n+1}q_n+q_{n-1}}
-\frac{p_n}{q_n}\notag\\
&=
-\frac{p_nq_{n-1}-q_np_{n-1}}
{(\alpha_{n+1}q_n+q_{n-1})q_n}
=
-\frac{(-1)^{n+1}}
{(\alpha_{n+1}q_n+q_{n-1})q_n},
\end{align}
where the last equality follows from \eqref{eq:det}. Since $\alpha_{n+1}>a_{n+1}$, we have
\[
\alpha_{n+1}q_n+q_{n-1}
>
a_{n+1}q_n+q_{n-1}
=
q_{n+1}
\geq
n+1,
\]
where the final weak inequality follows from Corollary~\ref{cor:qk>k}. Hence
\begin{equation}\label{eq:valuation-approximation}
\left|\alpha-\frac{p_n}{q_n}\right|
=
\frac{1}{(\alpha_{n+1}q_n+q_{n-1})q_n}
<
\frac{1}{q_{n+1}q_n}
\leq
\frac{1}{n(n+1)}.
\end{equation}
Letting $n\to\infty$, we obtain
\[
\lim_{n\to\infty}
\left|\alpha-\frac{p_n}{q_n}\right|
=0.
\]
Since the limit $\lim_{n\to\infty}p_n/q_n$ exists by Theorem~\ref{thm:limit-existence}, this estimate implies
\[
\alpha=\lim_{n\to\infty}\frac{p_n}{q_n}
=\lim_{n\to\infty}[a_0;a_1,\dots,a_n].
\]
This proves $F\circ G=\mathrm{id}_{\mathbb R\setminus\mathbb Q}$.

Next we show $G\circ F=\mathrm{id}_{\mathscr S}$. Take $(b_n)_{n=0}^\infty\in\mathscr S$, and put
\[
\alpha=\lim_{n\to\infty}[b_0;b_1,\dots,b_n].
\]
Let $(a_n)_{n=0}^\infty$ be the sequence obtained from $\alpha$ by \eqref{eq:irrational-to-sequence}. We must prove $(a_n)_{n=0}^\infty=(b_n)_{n=0}^\infty$. By Lemma~\ref{lem:integer-part}, we first get $a_0=b_0$. Hence
\begin{align*}
\alpha_1
&=
\frac{1}{\lim_{n\to\infty}[b_0;b_1,\dots,b_n]-b_0}
=
\lim_{n\to\infty}
\frac{1}{[b_0;b_1,\dots,b_n]-b_0}\\
&=
\lim_{n\to\infty}
\frac{1}{[0;b_1,\dots,b_n]}
=
\lim_{n\to\infty}[b_1;b_2,\dots,b_n].
\end{align*}
Applying Lemma~\ref{lem:integer-part} again gives $a_1=b_1$. Repeating the same argument proves $(a_n)_{n=0}^\infty=(b_n)_{n=0}^\infty$.
\end{proof}
This bijection yields the following statement.
\begin{theorem}\label{thm:euler}
For every irrational number $\alpha$, there exists a unique infinite regular continued fraction whose limit is $\alpha$. The sequence of partial quotients is obtained by the algorithm \eqref{eq:irrational-to-sequence}.
\end{theorem}
\begin{defi}
The infinite regular continued fraction whose limit is an irrational number $\alpha$ is called the \emph{infinite regular continued-fraction expansion} of $\alpha$.
\end{defi}

We next use infinite continued fractions to examine how an irrational number is related to real numbers that approximate it.

\begin{lemm}\label{lem:stability-of-initial-partial-quotients}
Let $(\alpha^{(i)})_{i=1}^{\infty}$ be a sequence of real numbers converging to an irrational number $\alpha$. Fix $k\geq 0$, and assume that the regular continued-fraction expansion of each $\alpha^{(i)}$ has at least $k+1$ partial quotients. Then, for all sufficiently large $i$, the first $k+1$ partial quotients of $\alpha^{(i)}$ and $\alpha$ coincide.
\end{lemm}
\begin{proof}
Write
\[
\alpha=[a_0;a_1,a_2,\dots],
\quad\text{and let}\quad
\frac{p_j}{q_j}=[a_0;a_1,\dots,a_j]
\]
be its $j$-th convergent, with $p_{-1}=1$ and $q_{-1}=0$. \index{Complete quotient}Also write $\alpha_j=[a_j;a_{j+1},\dots]$ for the complete quotients of $\alpha$.

Fix $k\geq 0$. For every $\xi>1$, Proposition~\ref{prop:recursion} gives
\[
[a_0;a_1,\dots,a_k,\xi]
=\frac{\xi p_k+p_{k-1}}{\xi q_k+q_{k-1}}.
\]
Moreover, Lemma~\ref{lem:det} gives
\[
p_kq_{k-1}-p_{k-1}q_k=(-1)^{k+1}.
\]
Consequently, as $\xi$ ranges over $(1,\infty)$, the expression on the right varies monotonically over the open interval whose endpoints are
\[
\frac{p_k+p_{k-1}}{q_k+q_{k-1}}
\quad\text{and}\quad
\frac{p_k}{q_k}.
\]
Every real number in this interval has a regular continued-fraction expansion with at least $k+2$ partial quotients, and its first $k+1$ partial quotients are $a_0,a_1,\dots,a_k$. Conversely, every real number with this property belongs to the same interval.

Since
\[
\alpha=[a_0;a_1,\dots,a_k,\alpha_{k+1}]
\quad\text{and}\quad
\alpha_{k+1}>1,
\]
the number $\alpha$ lies in the interior of this interval. Its distances from the two endpoints are
\[
\left|\alpha-\frac{p_k}{q_k}\right|
=\frac{1}{q_k(\alpha_{k+1}q_k+q_{k-1})}
\quad\text{and}\quad
\left|\alpha-\frac{p_k+p_{k-1}}{q_k+q_{k-1}}\right|
=\frac{\alpha_{k+1}-1}
{(\alpha_{k+1}q_k+q_{k-1})(q_k+q_{k-1})}.
\]
Thus, if we set
\[
\varepsilon_k(\alpha):=
\frac{1}{\alpha_{k+1}q_k+q_{k-1}}
\min\left\{
\frac{1}{q_k},
\frac{\alpha_{k+1}-1}{q_k+q_{k-1}}
\right\}>0,
\]
then every real number $x$ satisfying
\[
|x-\alpha|<\varepsilon_k(\alpha)
\]
belongs to the interval above. Hence the first $k+1$ partial quotients of $x$ and $\alpha$ coincide.

Since $\alpha^{(i)}\to\alpha$, for all sufficiently large $i$ we have
\[
|\alpha^{(i)}-\alpha|<\varepsilon_k(\alpha).
\]
The asserted agreement of the first $k+1$ partial quotients follows.
\end{proof}

It is natural to expect that truncating the infinite continued-fraction expansion of an irrational number $\alpha$ gives good rational approximations to $\alpha$. Conversely, it is also known that all sufficiently good rational approximations arise as such truncations. We now explain this. First we record an inequality that follows from the preceding discussion.
\begin{prop}\label{prop:valuation-a-p/q}
Let $\alpha$ be irrational, and let $\frac{p_n}{q_n}$ be its $n$-th convergent. Then
\[
\frac{1}{q_{n}q_{n+2}}
<
\left|\alpha-\frac{p_n}{q_n}\right|
<
\frac{1}{q_nq_{n+1}}.
\]
\end{prop}
\begin{proof}
The second inequality follows from \eqref{eq:valuation-approximation}. For the first one, use \eqref{eq:difference-approximation}. Since $\alpha_{n+1}=a_{n+1}+1/\alpha_{n+2}<a_{n+1}+1$, we estimate the denominator as follows:
\[
(\alpha_{n+1}q_n+q_{n-1})q_n
<
((a_{n+1}+1)q_n+q_{n-1})q_n
=
(q_{n+1}+q_n)q_n
\leq
(a_{n+2}q_{n+1}+q_n)q_n
=
q_{n+2}q_n.
\]
This gives the desired lower bound.
\end{proof}
The next proposition says that convergents are especially good approximations.
\begin{prop}\label{prop:comparing-alpha-p/q}
Let $\alpha$ be irrational and let $n\geq 1$. For the $n$-th convergent $\frac{p_n}{q_n}$ of $\alpha$, and for any rational number $\frac pq$ satisfying
\[
\frac{p_n}{q_n}\neq \frac pq,
\qquad
0<q\leq q_n,
\]
one has
\[
|q\alpha-p|\geq |q_{n-1}\alpha-p_{n-1}|>|q_n\alpha-p_n|.
\]
In particular,
\[
\left|\alpha-\frac pq\right|
>
\left|\alpha-\frac{p_n}{q_n}\right|.
\]
\end{prop}
\begin{proof}
If $\frac pq$ is not reduced, reducing it preserves the assumptions of the proposition. Thus we may assume from the beginning that $\frac pq$ is reduced. By the lower bound in Proposition~\ref{prop:valuation-a-p/q}, we have
\[
\frac{1}{q_{n+1}}<|q_{n-1}\alpha-p_{n-1}|,
\]
and by the upper bound we have
\[
|q_n\alpha-p_n|<\frac{1}{q_{n+1}}.
\]
Thus
\[
|q_n\alpha-p_n|
<
\frac{1}{q_{n+1}}
<
|q_{n-1}\alpha-p_{n-1}|,
\]
which gives the second inequality in the first assertion.

It remains to prove the first inequality. Consider the linear system
\[
\begin{bmatrix}
 p_{n}&p_{n-1}\\
 q_n & q_{n-1}
\end{bmatrix}
\begin{bmatrix}
c\\
d
\end{bmatrix}
=
\begin{bmatrix}
p\\q
\end{bmatrix}.
\]
Since $\begin{bsmallmatrix}p_n&p_{n-1}\\ q_n&q_{n-1}\end{bsmallmatrix}\in GL(2,\mathbb Z)$, the system has an integral solution. Solving it gives
\[
c=(-1)^{n+1}(pq_{n-1}-p_{n-1}q),
\qquad
d=(-1)^{n+1}(p_nq-q_np).
\]
Here $c,d\in\mathbb Z$, and the assumption $\frac pq\neq\frac{p_n}{q_n}$ implies $d\neq 0$.

If $c=0$, then the system gives $dp_{n-1}=p$ and $dq_{n-1}=q$, hence
\[
|q\alpha-p|
=
|d|\,|q_{n-1}\alpha-p_{n-1}|
\geq
|q_{n-1}\alpha-p_{n-1}|,
\]
as desired.

Now assume $c\neq 0$. Since $q\leq q_n$, if $c$ and $d$ had the same sign, then they would both have to be positive, since $q=cq_n+dq_{n-1}>0$. Hence $q=cq_n+dq_{n-1}$ would imply $q\geq q_n$, with equality only when $c=1$ and $d=0$. In that case we would have $p=p_n$ and $q=q_n$, contradicting the assumption. Hence $c$ and $d$ have opposite signs.

On the other hand, the proof of Theorem~\ref{thm:limit-existence} shows that either
\[
\frac{p_{n-1}}{q_{n-1}}<\alpha<\frac{p_n}{q_n}
\quad\text{or}\quad
\frac{p_n}{q_n}<\alpha<\frac{p_{n-1}}{q_{n-1}}.
\]
Thus $q_{n-1}\alpha-p_{n-1}$ and $q_n\alpha-p_n$ have opposite signs. It follows that
\[
c(q_n\alpha-p_n)
\quad\text{and}\quad
d(q_{n-1}\alpha-p_{n-1})
\]
have the same sign. Therefore
\begin{align*}
 |q\alpha-p|
&=
|(cq_n+dq_{n-1})\alpha-(cp_n+dp_{n-1})|\\
&=
|c(q_n\alpha-p_n)+d(q_{n-1}\alpha-p_{n-1})|\\
&=
|c(q_n\alpha-p_n)|+|d(q_{n-1}\alpha-p_{n-1})|\\
&\geq
|d(q_{n-1}\alpha-p_{n-1})|
\geq
|q_{n-1}\alpha-p_{n-1}|.
\end{align*}
This proves the first assertion.

The final assertion follows from
\[
\left|\alpha-\frac pq\right|
=
\frac{1}{q}|q\alpha-p|
>
\frac{1}{q}|q_n\alpha-p_n|
\geq
\frac{1}{q_n}|q_n\alpha-p_n|
=
\left|\alpha-\frac{p_n}{q_n}\right|.\qedhere
\]
\end{proof}
Using this, we prove that every sufficiently good approximation comes from a convergent.
\begin{theorem}\label{thm:good-is-approximation}\index{Diophantine approximation}\index{Convergent}
Let $\alpha$ be irrational. If a reduced fraction $\frac pq$ with $q>0$ satisfies
\[
\left|\alpha-\frac pq\right|<\frac{1}{2q^2},
\]
then $\frac pq$ is a convergent of $\alpha$; equivalently, it is obtained by truncating the continued-fraction expansion of $\alpha$, allowing a final partial quotient equal to $1$.
\end{theorem}
\begin{proof}
Suppose, for a contradiction, that $\frac pq\neq\frac{p_k}{q_k}$ for every convergent $\frac{p_k}{q_k}$ of $\alpha$. Since $\alpha$ is irrational, Corollary~\ref{cor:qk>k} implies that there exists $N$ such that $q<q_N$. Choose the smallest $n$ such that $q_n>q$. Then
\[
q_{n-1}\leq q<q_n.
\]
By Proposition~\ref{prop:comparing-alpha-p/q} and the assumption,
\[
|q_{n-1}\alpha-p_{n-1}|
\leq
|q\alpha-p|
<
\frac{1}{2q}.
\]
Therefore
\[
\frac{1}{qq_{n-1}}
\leq
\frac{|qp_{n-1}-pq_{n-1}|}{qq_{n-1}}
=
\left|\frac{p_{n-1}}{q_{n-1}}-\frac pq\right|
\leq
\left|\frac{p_{n-1}}{q_{n-1}}-\alpha\right|
+
\left|\alpha-\frac pq\right|
<
\frac{1}{2qq_{n-1}}+\frac{1}{2q^2}.
\]
This inequality implies $q<q_{n-1}$, contradicting $q_{n-1}\leq q$. Hence $\frac pq$ is a convergent of $\alpha$.
\end{proof}
\begin{rema}\label{rem:finite-remark}
Let $\alpha=[a_0;a_1,\dots,a_N]=p_N/q_N$ be a finite regular continued fraction. If $N\geq1$, then for $0\leq n\leq N-2$,
\[
\frac1{q_nq_{n+2}}<\left|\alpha-\frac{p_n}{q_n}\right|<\frac1{q_nq_{n+1}}.
\]
For the last proper convergent, Lemma~\ref{lem:det} instead gives
\[
\left|\alpha-\frac{p_{N-1}}{q_{N-1}}\right|=\frac1{q_{N-1}q_N}.
\]
Thus the two strict bounds apply only through $n=N-2$; the upper bound becomes an equality at $n=N-1$.

Proposition~\ref{prop:comparing-alpha-p/q} remains valid for rational $\alpha$ when $1\leq n\leq N$. For $1\leq n\leq N-1$, its proof applies with the preceding bounds, using the equality for the last proper convergent when $n=N-1$. For $n=N$, if $p/q\neq p_N/q_N$ and $0<q\leq q_N$, then
\[
|q\alpha-p|=\frac{|qp_N-pq_N|}{q_N}\geq\frac1{q_N}
=|q_{N-1}\alpha-p_{N-1}|>|q_N\alpha-p_N|=0.
\]
The inequality $|\alpha-p/q|>|\alpha-p_N/q_N|$ is immediate as well.

Theorem~\ref{thm:good-is-approximation} also holds for rational $\alpha$. Suppose that a reduced fraction $p/q$, with $q>0$, satisfies
\[
\left|\alpha-\frac pq\right|<\frac1{2q^2}.
\]
The case $p/q=\alpha$ is immediate. Assume otherwise and suppose that $p/q$ is not a convergent of $\alpha$. If $N=0$, then $\alpha$ is an integer, and
\[
\left|\alpha-\frac pq\right|=\frac{|q\alpha-p|}{q}\geq\frac1q>\frac1{2q^2},
\]
a contradiction. If $N\geq1$, then
\[
\frac1{q_N}\leq\frac{|qp_N-pq_N|}{q_N}=|q\alpha-p|<\frac1{2q},
\]
so $q<q_N$. Choose the smallest $n$ such that $q_n>q$. Then $1\leq n\leq N$ and $q_{n-1}\leq q<q_n$. The rational version of Proposition~\ref{prop:comparing-alpha-p/q} gives
\[
|q_{n-1}\alpha-p_{n-1}|\leq|q\alpha-p|<\frac1{2q}.
\]
Since $p/q$ is not a convergent, $qp_{n-1}-pq_{n-1}\neq0$. Consequently,
\[
\frac1{qq_{n-1}}\leq\left|\frac{p_{n-1}}{q_{n-1}}-\frac pq\right| 
\leq\left|\frac{p_{n-1}}{q_{n-1}}-\alpha\right|+\left|\alpha-\frac pq\right|
<\frac1{2qq_{n-1}}+\frac1{2q^2}.
\]
This implies $q<q_{n-1}$, contradicting $q_{n-1}\leq q$.
\end{rema}

\section{Unimodular Group Orbits of Irrational Numbers}
In this section we decompose irrational numbers into orbits under a matrix group. The reason this decomposition is useful will become clear in the next chapter: irrational numbers lying in the same orbit have the same Lagrange constant. Hence, when studying possible values of Lagrange constants, this orbit decomposition removes redundant work. We first introduce the group that acts on irrational numbers.
\begin{defi}\label{def:unimodular-group}\index{Unimodular group}
Consider the set
\[
GL(2,\mathbb Z):=
\left\{
 \begin{bmatrix}a&b\\ c&d\end{bmatrix}
\ \middle |\
a,b,c,d\in\mathbb Z,\ |ad-bc|=1
\right\}.
\]
It is a group under matrix multiplication. We call it the \emph{unimodular group}.
\end{defi}
For \(A=\begin{bsmallmatrix}a&b\\c&d\end{bsmallmatrix}\in GL(2,\mathbb Z)\) and $\alpha\in\mathbb R\setminus\mathbb Q$, define the action of $A$ on $\alpha$ by the fractional linear transformation
\[
A\alpha:=\frac{a\alpha+b}{c\alpha+d}.
\]
We have the following.
\begin{prop}\label{prop:action-irrational}\index{Unimodular group}
For $A\in GL(2,\mathbb Z)$ and $\alpha\in\mathbb R\setminus\mathbb Q$, one has $A\alpha\in\mathbb R\setminus\mathbb Q$. Moreover, the fractional linear transformations of $GL(2,\mathbb Z)$ define a left action
\[
GL(2,\mathbb Z)\curvearrowright\mathbb R\setminus\mathbb Q.
\]
\end{prop}
\begin{proof}
Let \(A=\begin{bsmallmatrix}a&b\\c&d\end{bsmallmatrix}\) and \(\alpha\in\mathbb R\setminus\mathbb Q\).
Suppose $A\alpha$ is rational. Then
\[
\alpha=A^{-1}(A\alpha)
=
\frac{d(A\alpha)-b}{-c(A\alpha)+a}
\]
would also be rational, a contradiction. Hence $A\alpha$ is irrational.

Let \(E_2=\begin{bsmallmatrix}1&0\\0&1\end{bsmallmatrix}\). Then $E_2\alpha=\alpha$. If \(B=\begin{bsmallmatrix}e&f\\g&h\end{bsmallmatrix}\), then
\[
A(B\alpha)
=
A\frac{e\alpha+f}{g\alpha+h}
=
\frac{a\frac{e\alpha+f}{g\alpha+h}+b}
{c\frac{e\alpha+f}{g\alpha+h}+d}
=
\frac{(ae+bg)\alpha+(af+bh)}
{(ce+dg)\alpha+(cf+dh)}.
\]
On the other hand,
\[
AB=
\begin{bmatrix}
    a &b\\c&d
\end{bmatrix}
\begin{bmatrix}
    e &f\\g&h
\end{bmatrix}
=
\begin{bmatrix}
    ae+bg&af+bh\\ce+dg&cf+dh
\end{bmatrix},
\]
so
\[
(AB)\alpha=
\frac{(ae+bg)\alpha+(af+bh)}
{(ce+dg)\alpha+(cf+dh)}.
\]
Therefore $A(B\alpha)=(AB)\alpha$, and the assertion follows.
\end{proof}
\begin{rema}
By definition of the action, $A\alpha=(-A)\alpha$.
\end{rema}
We now introduce unimodular equivalence of irrational numbers.
\begin{defi}\label{def:unimodular-equivalence-irrational}\index{Unimodular orbit}\index{Unimodular equivalence!of irrational numbers}
Let $\alpha,\beta$ be irrational numbers. If there exists $A\in GL(2,\mathbb Z)$ such that $\beta=A\alpha$, then $\alpha$ and $\beta$ are said to be \emph{unimodularly equivalent}, or simply \emph{equivalent}. We write $\alpha\sim\beta$. The equivalence class
\[
O_\alpha=\{\beta\mid \beta\sim\alpha\}
\]
is called the \emph{unimodular orbit}, or simply the \emph{orbit}, of $\alpha$.
\end{defi}
The relation between the unimodular action on irrational numbers and continued fractions can be expressed succinctly using continued-fraction matrices.
\begin{prop}\label{prop:continued-fraction-matrix-representation}\index{Continued-fraction matrix}
Let $a_0,a_1,\dots,a_k,\alpha\in\mathbb R$. Whenever both sides are defined,
\[
[a_0;a_1,\dots,a_k,\alpha]
=
\begin{bmatrix}
       a_0&1\\1&0
\end{bmatrix}
\begin{bmatrix}
       a_1&1\\1&0
\end{bmatrix}
\cdots
\begin{bmatrix}
       a_k&1\\1&0
\end{bmatrix}\alpha.
\]
\end{prop}
\begin{proof}
By Proposition~\ref{prop:recursion} and Theorem~\ref{thm:continued-fraction-matrix},
\[
[a_0;a_1,\dots,a_k,\alpha]
=
\frac{\alpha p_k+p_{k-1}}{\alpha q_k+q_{k-1}}
=
\begin{bmatrix}
    p_k&p_{k-1}\\q_k&q_{k-1}
\end{bmatrix}\alpha
=
\begin{bmatrix}
    a_0&1\\1&0
\end{bmatrix}
\begin{bmatrix}
    a_1&1\\1&0
\end{bmatrix}
\cdots
\begin{bmatrix}
    a_k&1\\1&0
\end{bmatrix}\alpha.
\]
\end{proof}
The next theorem characterizes equivalence of two irrational numbers in terms of their infinite continued-fraction expansions.
\begin{theorem}[Serret's theorem]\label{thm:characterization-equivalent}\index{Serret's theorem}\index{Unimodular equivalence!of irrational numbers}
Let
\[
\alpha=[a_0;a_1,\dots],
\qquad
\beta=[b_0;b_1,\dots]
\]
be irrational numbers. Then $\alpha\sim\beta$ if and only if there exist $n,m\in\mathbb Z_{\geq 0}$ such that
\[
[a_n;a_{n+1},\dots]=[b_m;b_{m+1},\dots].
\]
In particular, by uniqueness of regular continued-fraction expansions, $a_{n+h}=b_{m+h}$ for every $h\in\mathbb Z_{\geq 0}$.
\end{theorem}
We first prove a lemma.
\begin{lemm}\label{lemm:generate}
Let \(A=\begin{bsmallmatrix}a&b\\c&d\end{bsmallmatrix}\in GL(2,\mathbb Z)\) and assume $c>d>0$. Then there exist $\ell\in\mathbb Z_{\geq 0}$ and integers $c_0\in\mathbb Z$, $c_1,\dots,c_\ell\in\mathbb Z_{\geq 1}$ such that
\[
A=
\begin{bmatrix}c_0&1\\1&0\end{bmatrix}
\begin{bmatrix}c_1&1\\1&0\end{bmatrix}
\cdots
\begin{bmatrix}c_\ell&1\\1&0\end{bmatrix}.
\]
\end{lemm}
\begin{proof}
Take the finite continued-fraction expansion
\[
\frac{a}{c}=[a_0;a_1,\dots,a_k].
\]
Let $\frac{p_i}{q_i}$ be its $i$-th convergent. Then
\[
\begin{bmatrix}
p_k&p_{k-1}\\q_k&q_{k-1}
\end{bmatrix}
=
\begin{bmatrix}a_0&1\\1&0\end{bmatrix}
\begin{bmatrix}a_1&1\\1&0\end{bmatrix}
\cdots
\begin{bmatrix}a_k&1\\1&0\end{bmatrix}.
\]
If necessary, choose the expansion so that $a_k\geq 2$. Then we may also view
\[
[a_0;a_1,\dots,a_k]
=
[a_0;a_1,\dots,a_k-1,1].
\]
If $\frac{p'_i}{q'_i}$ denotes the $i$-th convergent for this latter expression, then
\[
\begin{bmatrix}
p'_{k+1}&p'_k\\q'_{k+1}&q'_k
\end{bmatrix}
=
\begin{bmatrix}a_0&1\\1&0\end{bmatrix}
\begin{bmatrix}a_1&1\\1&0\end{bmatrix}
\cdots
\begin{bmatrix}a_k-1&1\\1&0\end{bmatrix}
\begin{bmatrix}1&1\\1&0\end{bmatrix}.
\]
Here $a_k-1\geq 1$.

Since $A\in GL(2,\mathbb Z)$ and
\[
\det
\begin{bmatrix}
p_k&p_{k-1}\\q_k&q_{k-1}
\end{bmatrix}
=(-1)^{k+1},
\qquad
\det
\begin{bmatrix}
p'_{k+1}&p'_k\\q'_{k+1}&q'_k
\end{bmatrix}
=(-1)^{k+2},
\]
exactly one of the following equalities holds:
\[
\det A=
\det
\begin{bmatrix}
p_k&p_{k-1}\\q_k&q_{k-1}
\end{bmatrix},
\qquad
\det A=
\det
\begin{bmatrix}
p'_{k+1}&p'_k\\q'_{k+1}&q'_k
\end{bmatrix}.
\]
We show that, in the first case,
\[
A=
\begin{bmatrix}
p_k&p_{k-1}\\q_k&q_{k-1}
\end{bmatrix},
\]
and in the second case the analogous primed equality holds. We prove the first case; the second is identical.

Since $A\in GL(2,\mathbb Z)$, the integers $a$ and $c$ are relatively prime. Since $c>0$, the fraction $\frac ac$ is reduced. By Corollary~\ref{cor:reduced-fraction}, $\frac{p_k}{q_k}$ is also reduced. Hence $a=p_k$ and $c=q_k$. Using the equality of determinants gives $p_kd-bq_k=p_kq_{k-1}-q_kp_{k-1}$. Equivalently, $p_k(d-q_{k-1})=q_k(b-p_{k-1})$.
Since $p_k$ and $q_k$ are relatively prime, we have $q_k\mid d-q_{k-1}$. If $d-q_{k-1}\geq 0$, then the assumption $d<c$ gives $d-q_{k-1}<c-q_{k-1}=q_k-q_{k-1}<q_k$. If $d-q_{k-1}<0$, then $q_{k-1}-d<q_{k-1}<q_k$.
Thus in either case $|d-q_{k-1}|<q_k$. Since $q_k$ divides $d-q_{k-1}$, we must have $d=q_{k-1}$. The equality $b=p_{k-1}$ then follows from the other three entries and the determinant equality. Hence
\[
A=
\begin{bmatrix}
p_k&p_{k-1}\\q_k&q_{k-1}
\end{bmatrix}
=
\begin{bmatrix}a_0&1\\1&0\end{bmatrix}
\begin{bmatrix}a_1&1\\1&0\end{bmatrix}
\cdots
\begin{bmatrix}a_k&1\\1&0\end{bmatrix}.
\]
In the second determinant case, the same argument gives
\[
A=
\begin{bmatrix}
p'_{k+1}&p'_k\\q'_{k+1}&q'_k
\end{bmatrix}
=
\begin{bmatrix}a_0&1\\1&0\end{bmatrix}
\begin{bmatrix}a_1&1\\1&0\end{bmatrix}
\cdots
\begin{bmatrix}a_k-1&1\\1&0\end{bmatrix}
\begin{bmatrix}1&1\\1&0\end{bmatrix}.
\]
This proves the lemma.
\end{proof}
\begin{proof}[Proof of Theorem~\ref{thm:characterization-equivalent}]
First assume that there exist $m,n\in\mathbb Z_{\geq 0}$ such that
\[
[a_n;a_{n+1},\dots]=[b_m;b_{m+1},\dots].
\]
Put
\[
\alpha_n:=[a_n;a_{n+1},\dots],
\qquad
\beta_m:=[b_m;b_{m+1},\dots].
\]
Then
\[
\alpha=[a_0;a_1,\dots,a_{n-1},\alpha_n],
\qquad
\beta=[b_0;b_1,\dots,b_{m-1},\beta_m],
\]
where, if $n=0$, the first expression simply means $\alpha=\alpha_0$. By Proposition~\ref{prop:continued-fraction-matrix-representation},
\[
\alpha=
\begin{bmatrix}a_0&1\\1&0\end{bmatrix}
\begin{bmatrix}a_1&1\\1&0\end{bmatrix}
\cdots
\begin{bmatrix}a_{n-1}&1\\1&0\end{bmatrix}\alpha_n,
\]
and
\[
\beta=
\begin{bmatrix}b_0&1\\1&0\end{bmatrix}
\begin{bmatrix}b_1&1\\1&0\end{bmatrix}
\cdots
\begin{bmatrix}b_{m-1}&1\\1&0\end{bmatrix}\beta_m.
\]
Since $\alpha_n=\beta_m$, we may solve the first equality for $\alpha_n$ and substitute into the second equality. We obtain
\[
\beta=
\begin{bmatrix}b_0&1\\1&0\end{bmatrix}
\begin{bmatrix}b_1&1\\1&0\end{bmatrix}
\cdots
\begin{bmatrix}b_{m-1}&1\\1&0\end{bmatrix}
\begin{bmatrix}a_{n-1}&1\\1&0\end{bmatrix}^{-1}
\cdots
\begin{bmatrix}a_1&1\\1&0\end{bmatrix}^{-1}
\begin{bmatrix}a_0&1\\1&0\end{bmatrix}^{-1}
\alpha.
\]
Each factor lies in $GL(2,\mathbb Z)$. Hence the product also lies in $GL(2,\mathbb Z)$, and therefore $\alpha$ and $\beta$ are equivalent.

Conversely, assume $\alpha\sim\beta$. Then there exists \(A=\begin{bsmallmatrix}a&b\\c&d\end{bsmallmatrix}\in GL(2,\mathbb Z)\) such that $\beta=A\alpha$. Since $A\alpha=(-A)\alpha$, we may replace $A$ by $-A$ if necessary and assume $c\alpha+d>0$. For any $n\geq 2$,
\[
\beta=A\alpha
=
\begin{bmatrix}a&b\\c&d\end{bmatrix}
\begin{bmatrix}
p_{n-1}&p_{n-2}\\
q_{n-1}&q_{n-2}
\end{bmatrix}
\alpha_n,
\]
where
\[
\begin{bmatrix}
p_{n-1}&p_{n-2}\\
q_{n-1}&q_{n-2}
\end{bmatrix}
:=
\begin{bmatrix}a_0&1\\1&0\end{bmatrix}
\begin{bmatrix}a_1&1\\1&0\end{bmatrix}
\cdots
\begin{bmatrix}a_{n-1}&1\\1&0\end{bmatrix}.
\]
Set
\[
\begin{bmatrix}a'&b'\\c'&d'\end{bmatrix}
:=
\begin{bmatrix}a&b\\c&d\end{bmatrix}
\begin{bmatrix}
p_{n-1}&p_{n-2}\\
q_{n-1}&q_{n-2}
\end{bmatrix}.
\]
We show that $n$ can be chosen so that $c'>d'>0$. Direct computation gives
\begin{align}
c'&=cp_{n-1}+dq_{n-1}
=q_{n-1}\left(\frac{cp_{n-1}}{q_{n-1}}+d\right),\label{eq:c-prime}\\
d'&=cp_{n-2}+dq_{n-2}
=q_{n-2}\left(\frac{cp_{n-2}}{q_{n-2}}+d\right).\label{eq:d-prime}
\end{align}
If $c=0$, then after replacing $A$ by $-A$ if necessary, we may assume $d>0$. Then $c'=dq_{n-1}$ and $d'=dq_{n-2}$, so for sufficiently large $n$ we have $q_{n-1}>q_{n-2}>0$, and hence $c'>d'>0$.

Thus assume $c\neq 0$. Since
\[
\lim_{n\to\infty}
\left(\frac{cp_{n-1}}{q_{n-1}}+d\right)
=
\lim_{n\to\infty}
\left(\frac{cp_{n-2}}{q_{n-2}}+d\right)
=
c\alpha+d>0,
\]
we may take $n$ sufficiently large so that both factors in parentheses in \eqref{eq:c-prime} and \eqref{eq:d-prime} are positive. Then $c',d'>0$. Moreover,
\begin{align*}
c'-d'
&=
q_{n-1}\left(\frac{cp_{n-1}}{q_{n-1}}+d\right)
-
q_{n-2}\left(\frac{cp_{n-2}}{q_{n-2}}+d\right)>
q_{n-1}\left(\frac{cp_{n-1}}{q_{n-1}}+d\right)
-
q_{n-1}\left(\frac{cp_{n-2}}{q_{n-2}}+d\right)\\
&=
\frac{c}{q_{n-2}}
(p_{n-1}q_{n-2}-q_{n-1}p_{n-2})
=
\frac{(-1)^n c}{q_{n-2}}.
\end{align*}
Since $n$ may still be chosen with either parity, we choose it so that $(-1)^n c>0$. Then $c'-d'>0$. Thus for such an $n$ we have $c'>d'>0$.

By Lemma~\ref{lemm:generate}, there exist $\ell\in\mathbb Z_{\geq 0}$ and $c_0\in\mathbb Z$, $c_1,\dots,c_\ell\in\mathbb Z_{\geq 1}$ such that
\[
\begin{bmatrix}a'&b'\\c'&d'\end{bmatrix}
=
\begin{bmatrix}c_0&1\\1&0\end{bmatrix}
\begin{bmatrix}c_1&1\\1&0\end{bmatrix}
\cdots
\begin{bmatrix}c_\ell&1\\1&0\end{bmatrix}.
\]
Hence
\[
\beta
=
\begin{bmatrix}c_0&1\\1&0\end{bmatrix}
\begin{bmatrix}c_1&1\\1&0\end{bmatrix}
\cdots
\begin{bmatrix}c_\ell&1\\1&0\end{bmatrix}\alpha_n
=
[c_0;c_1,\dots,c_\ell,\alpha_n].
\]
Since $\alpha_n=[a_n;a_{n+1},\dots]$, this is
\[
\beta=[c_0;c_1,\dots,c_\ell,a_n,a_{n+1},\dots].
\]
Because $n\geq 1$, we have $a_n\geq 1$, so this is the infinite regular continued-fraction expansion of $\beta$. Taking $m=\ell+1$, we obtain
\[
[a_n;a_{n+1},\dots]=[b_m;b_{m+1},\dots].
\]
This completes the proof.
\end{proof}

\section{Periodic Continued Fractions and Quadratic Irrationals}\label{sec:quadratic-irrational}
In this section we study the case where the partial quotients in the regular continued-fraction expansion of an irrational number eventually become periodic. We characterize this phenomenon in terms of a property of the irrational number itself. We begin by defining periodic continued fractions.
\begin{defi}\label{def:periodic-continued-fraction}\index{Purely periodic continued fraction}\index{Periodic continued fraction}
An infinite regular continued fraction whose partial quotients are eventually periodic, as in
\[
[a_0;a_1,a_2,\dots,a_k,b_0,b_1,\dots,b_\ell,b_0,b_1,\dots,b_\ell,b_0,b_1,\dots],
\]
is called a \emph{periodic continued fraction}. We write it as
\[
[a_0;a_1,a_2,\dots,a_k,\overline{b_0,b_1,\dots,b_\ell}].
\]
The finite sequence $(b_0,\dots,b_\ell)$ is called the \emph{period}. A periodic continued fraction whose periodic part starts at the beginning is called a \emph{purely periodic continued fraction}.
\end{defi}
Next we define quadratic irrationals.
\begin{defi}\label{def:quadratic-irrational}\index{Discriminant!of a quadratic irrational}\index{Quadratic conjugate}\index{Quadratic irrational}\index{Reduced quadratic irrational}
An irrational number $\alpha$ is called a \emph{quadratic irrational} if it is a root of a quadratic equation with integer coefficients. Equivalently, there exist $a\in\mathbb Z_{\geq 1}$ and $b,c\in\mathbb Z$ with $\gcd(a,b,c)=1$ such that
\[
\alpha=\frac{-b+\varepsilon\sqrt{D}}{2a},
\]
where $D=b^2-4ac$, $\varepsilon\in\{1,-1\}$, $D>0$, and $D$ is not a square. In this case $D$ is called the \emph{discriminant} of $\alpha$. For such a quadratic irrational $\alpha$, its quadratic conjugate is denoted by $\alpha'$, namely
\[
\alpha'=\frac{-b-\varepsilon\sqrt{D}}{2a}.
\]
If $\alpha>1$ and $-1<\alpha'<0$, then $\alpha$ is called a \emph{reduced quadratic irrational}.
\end{defi}
If $\alpha$ is reduced, then $\alpha>\alpha'$, and hence necessarily $\varepsilon=1$. We also impose $a\geq 1$ and $\gcd(a,b,c)=1$ in order to make the discriminant $D$ a well-defined quantity associated with $\alpha$, rather than depending on a nonprimitive multiple of a quadratic equation.

Let $I_2$ be the set of all quadratic irrationals, and let $R_2$ be the set of all reduced quadratic irrationals. For a positive nonsquare integer $d$, let $I_2(d)$ be the set of quadratic irrationals with discriminant $d$, and let $R_2(d)$ be the set of reduced quadratic irrationals with discriminant $d$.

The main theorem of this section is the following.
\begin{theorem}[Lagrange's theorem]\label{thm:characterization-quadratic}\index{Purely periodic continued fraction}\index{Lagrange's theorem}\index{Periodic continued fraction}\index{Quadratic irrational}\index{Reduced quadratic irrational}
The following statements hold.
\begin{itemize}
    \item[(1)] The infinite regular continued-fraction expansion of an irrational number $\alpha$ is periodic if and only if $\alpha$ is a quadratic irrational.
    \item[(2)] The infinite regular continued-fraction expansion of an irrational number $\alpha$ is purely periodic if and only if $\alpha$ is a reduced quadratic irrational.
\end{itemize}
\end{theorem}
We prepare for the proof.
\begin{theorem}\label{thm:action-on-I2}\index{Discriminant!of a quadratic irrational}
The action of the unimodular group on $\mathbb R\setminus\mathbb Q$ restricts to an action on $I_2$, and for every positive nonsquare integer $d\in\mathbb Z_{>0}$ it restricts to an action on $I_2(d)$.
\end{theorem}
\begin{proof}
It suffices to verify that $I_2(d)$ is preserved. Let $\alpha\in I_2(d)$, and let $ax^2+bx+c$ be a quadratic polynomial having $\alpha$ as a root, with $\gcd(a,b,c)=1$. Then the discriminant $b^2-4ac$ is the discriminant of $\alpha$. Let \(M=\begin{bsmallmatrix}s&t\\u&v\end{bsmallmatrix}\in GL(2,\mathbb Z)\) and put $\beta=M\alpha$. By Proposition~\ref{prop:action-irrational}, $\beta$ is irrational. Solving $\beta=M\alpha$ for $\alpha$, we obtain
\[
\alpha=M^{-1}\beta=\frac{v\beta-t}{-u\beta+s}.
\]
Substituting
\[
x=\alpha=\frac{v\beta-t}{-u\beta+s}
\]
into $ax^2+bx+c=0$ and clearing denominators gives
\[
(av^2+cu^2-buv)\beta^2
+
(-2atv+bsv+but-2csu)\beta
+
at^2+cs^2-bst
=
0.
\]
Thus $\beta$ is a root of the quadratic polynomial
\[
(av^2+cu^2-buv)x^2
+
(-2atv+bsv+but-2csu)x
+
at^2+cs^2-bst.
\]
Put
\[
A=av^2+cu^2-buv,\qquad
B=-2atv+bsv+but-2csu,\qquad
C=at^2+cs^2-bst.
\]
To see that $B^2-4AC$ is the discriminant of $\beta$, we check that $\gcd(A,B,C)=1$ and $A\neq 0$. If $A<0$, we multiply the polynomial by $-1$; this makes the leading coefficient positive and does not change the root or the discriminant. The condition $A\neq 0$ follows because $\beta$ is irrational. Now suppose $e\mid A$, $e\mid B$, and $e\mid C$. Then
\begin{align*}
    s^2A+suB+u^2C&=a(sv-tu)^2=a,\\
    2stA+(sv+tu)B+2uvC&=b(sv-tu)^2=b,\\
    t^2A+tvB+v^2C&=c(sv-tu)^2=c.
\end{align*}
Since $\gcd(a,b,c)=1$, any common divisor of $A,B,C$ must divide $1$. Thus $\gcd(A,B,C)=1$. Finally,
\[
B^2-4AC=(b^2-4ac)(sv-tu)^2=b^2-4ac=d,
\]
because $M\in GL(2,\mathbb Z)$. Hence $\beta\in I_2(d)$.
\end{proof}
\begin{theorem}\label{thm:reduced-equiv}\index{Reduced quadratic irrational}
Fix a positive nonsquare integer $d$. For every
\[
\alpha=[a_0;a_1,\dots]\in I_2(d),
\]
the complete quotient $\alpha_n=[a_n;a_{n+1},\dots]$ belongs to $R_2(d)$ for all sufficiently large $n$. Moreover, once this holds for some $n$, it holds for every $m\geq n$.
\end{theorem}
\begin{proof}
Let
\[
\alpha=[a_0;a_1,\dots]\in I_2(d),
\qquad
\alpha_k:=[a_k;a_{k+1},\dots].
\]
For $k\geq 1$,
\[
\alpha
=
[a_0;a_1,\dots,a_{k-1},\alpha_k]
=
\begin{bmatrix}a_0&1\\1&0\end{bmatrix}
\begin{bmatrix}a_1&1\\1&0\end{bmatrix}
\cdots
\begin{bmatrix}a_{k-1}&1\\1&0\end{bmatrix}\alpha_k,
\]
where the last equality follows from Proposition~\ref{prop:continued-fraction-matrix-representation}. Since each matrix \(\begin{bsmallmatrix}a_i&1\\1&0\end{bsmallmatrix}\) lies in $GL(2,\mathbb Z)$, the numbers $\alpha$ and $\alpha_k$ are equivalent. By Theorem~\ref{thm:action-on-I2}, we obtain $\alpha_k\in I_2(d)$ for every $k\geq 0$.

Moreover, the construction algorithm \eqref{eq:irrational-to-sequence} gives
\[
\alpha_{k+1}=\frac{1}{\alpha_k-\lfloor\alpha_k\rfloor}.
\]
Since $0<\alpha_k-\lfloor\alpha_k\rfloor<1$, we have $\alpha_{k+1}>1$. Thus $\alpha_k>1$ for all $k\geq 1$. It remains to show that, for all sufficiently large indices, the conjugate lies in the interval $(-1,0)$.

From
\[
\alpha=
\begin{bmatrix}
p_k&p_{k-1}\\q_k&q_{k-1}
\end{bmatrix}\alpha_{k+1},
\]
we solve for $\alpha_{k+1}$ and obtain
\[
\alpha_{k+1}
=
\begin{bmatrix}
p_k&p_{k-1}\\q_k&q_{k-1}
\end{bmatrix}^{-1}\alpha
=
\begin{bmatrix}
q_{k-1}&-p_{k-1}\\-q_k&p_k
\end{bmatrix}\alpha.
\]
The inverse matrix may differ from the displayed matrix by an overall sign, but multiplying the matrix by $-1$ does not change the fractional linear transformation. Taking quadratic conjugates of both sides gives
\[
\alpha'_{k+1}
=
-\frac{q_{k-1}\alpha'-p_{k-1}}{q_k\alpha'-p_k}
=
-\frac{q_{k-1}}{q_k}
\frac{\alpha'-\frac{p_{k-1}}{q_{k-1}}}{\alpha'-\frac{p_k}{q_k}}.
\]
Since $\lim_{k\to\infty}\frac{p_k}{q_k}=\alpha$, we have
\[
\lim_{k\to\infty}
\frac{\alpha'-\frac{p_{k-1}}{q_{k-1}}}
{\alpha'-\frac{p_k}{q_k}}
=1.
\]
Therefore, there exists $N$ such that, for every $k\geq N$,
\[
\frac{\alpha'-\frac{p_{k-1}}{q_{k-1}}}
{\alpha'-\frac{p_k}{q_k}}
>0.
\]
Since $q_{k-1}/q_k>0$, this implies $\alpha'_{k+1}<0$ for every $k\geq N$. Now take conjugates in the defining relation
\[
\alpha_{k+2}=\frac{1}{\alpha_{k+1}-a_{k+1}}.
\]
We obtain
\[
\alpha'_{k+2}
=
-\frac{1}{a_{k+1}-\alpha'_{k+1}}.
\]
Because $\alpha'_{k+1}<0$ and $a_{k+1}\geq 1$, this gives
\[
-1<\alpha'_{k+2}<0
\]
for every $k\geq N$. Thus $\alpha_{k+2}\in R_2(d)$ for every $k\geq N$. Moreover, if $\alpha_m$ is reduced, then $\alpha_m>1$ and $-1<\alpha_m'<0$ imply
\[
\alpha_{m+1}=\frac{1}{\alpha_m-a_m}>1,
\qquad
\alpha_{m+1}'=\frac{1}{\alpha_m'-a_m}\in(-1,0),
\]
so the reduced condition persists for all later indices.
\end{proof}
The preceding theorem immediately gives the following corollary.
\begin{coro}\label{cor:reduced-equiv}
For every positive nonsquare integer $d$ and every $\alpha\in I_2(d)$, there exists $\beta\in R_2(d)$ such that $\alpha\sim\beta$. Thus every orbit in $I_2(d)$ has a representative in $R_2(d)$.
\end{coro}
\begin{lemm}\label{lem:finiteness-R2d}
Let $\alpha\in R_2(d)$. Suppose that $\alpha$ is a root of $ax^2+bx+c$, where $a\in\mathbb Z_{\geq 1}$ and $\gcd(a,b,c)=1$. Then $0<-b<\sqrt d$.
In particular, $R_2(d)$ is finite for every positive nonsquare integer $d$.
\end{lemm}
\begin{proof}
Since $\alpha$ is reduced,
\[
\alpha=\frac{-b+\sqrt d}{2a}>1,
\qquad
-1<\alpha'=\frac{-b-\sqrt d}{2a}<0.
\]
Multiplying by $2a$ gives
\[
-b+\sqrt d>2a>b+\sqrt d>0.
\]
In particular, $-b+\sqrt d>b+\sqrt d$ gives $0<-b$, and $b+\sqrt d>0$ gives $-b<\sqrt d$. Thus $0<-b<\sqrt d$. For fixed $d$, there are only finitely many possible values of $b$. Since $d-b^2=-4ac$ and $a,c\in\mathbb Z$, there are also only finitely many possible pairs $(a,c)$. Hence $R_2(d)$ is finite.
\end{proof}
\begin{lemm}\label{lem:count-down}
Let $\alpha$ be a reduced quadratic irrational. If the continued fraction algorithm gives
\[
\alpha=[a_0;\alpha_1],
\quad\text{then}\quad
\lfloor\alpha\rfloor=a_0=\left\lfloor-\frac{1}{\alpha_1'}\right\rfloor.
\]
\end{lemm}
\begin{proof}
From
\[
\alpha=[a_0;\alpha_1]=a_0+\frac{1}{\alpha_1},
\]
taking conjugates gives
\[
\alpha'=a_0+\frac{1}{\alpha'_1}.
\]
Rearranging, we obtain
\[
-\frac{1}{\alpha'_1}=a_0+(-\alpha').
\]
Since $-1<\alpha'<0$, the integer part of $-\frac{1}{\alpha'_1}$ is $a_0$.
\end{proof}
\begin{proof}[Proof of Theorem~\ref{thm:characterization-quadratic}]
We first prove the implication in (2) that a purely periodic continued fraction represents an element of $R_2$. Let
\[
\alpha=[\overline{a_0;a_1,\dots,a_{n-1}}]
\]
with $n\geq1$.
Then
\[
\alpha=[a_0;a_1,\dots,a_{n-1},\alpha]
=
\frac{\alpha p_{n-1}+p_{n-2}}{\alpha q_{n-1}+q_{n-2}}.
\]
Rearranging gives
\[
q_{n-1}\alpha^2+(q_{n-2}-p_{n-1})\alpha-p_{n-2}=0.
\]
Thus $\alpha$ is a root of a quadratic equation. Since $\alpha\notin\mathbb Q$, we have $\alpha\in I_2$. Since the continued fraction is purely periodic,
\[
\alpha
=
[a_0;a_1,\dots,a_{n-1},\alpha]
=
[a_0;a_1,\dots,a_{n-1},a_0,a_1,\dots,a_{n-1},\alpha]
=
\cdots.
\]
Thus $\alpha=\alpha_{kn}$ for every $k\in\mathbb Z_{\geq 0}$. By Theorem~\ref{thm:reduced-equiv}, this implies $\alpha\in R_2$.

Next we prove the implication in (1) that a periodic continued fraction represents an element of $I_2$. Suppose
\[
\alpha=[a_0;a_1,a_2,\dots,a_{n-1},
\overline{a_n,a_{n+1},\dots,a_{n+k-1}}]
\]
for some $n\geq1$ and $k\geq1$.
Put
\[
\alpha_n=[\overline{a_n;a_{n+1},\dots,a_{n+k-1}}].
\]
By the previous paragraph, $\alpha_n\in R_2\subset I_2$. Since
\[
\alpha=[a_0;a_1,a_2,\dots,a_{n-1},\alpha_n]
=
\begin{bmatrix}
p_{n-1}&p_{n-2}\\q_{n-1}&q_{n-2}
\end{bmatrix}\alpha_n
\]
and the matrix lies in $GL(2,\mathbb Z)$, Theorem~\ref{thm:action-on-I2} shows that $\alpha\in I_2$.

We now prove the converse implication in (1): if $\alpha\in I_2$, then the continued-fraction expansion of $\alpha$ is periodic. Let $d$ be the discriminant of $\alpha$. Then $d$ is a positive nonsquare integer and $\alpha\in I_2(d)$. By Theorem~\ref{thm:reduced-equiv}, for a sufficiently large $n$, $m\geq n$ implies $\alpha_m\in R_2(d)$. Since $R_2(d)$ is finite by Lemma~\ref{lem:finiteness-R2d}, there exist indices
\[
n\leq \ell<\ell'
\]
such that $\alpha_\ell=\alpha_{\ell'}$. Then
\[
\alpha_\ell
=
[a_\ell;a_{\ell+1},\dots,a_{\ell'-1},\alpha_{\ell'}]
=
[a_\ell;a_{\ell+1},\dots,a_{\ell'-1},\alpha_\ell].
\]
Therefore
\[
\alpha_\ell=[\overline{a_\ell;a_{\ell+1},\dots,a_{\ell'-1}}],
\]
and the continued-fraction expansion of $\alpha$ is periodic.

\index{Gauss map}Finally, we prove the converse implication in (2): if $\alpha\in R_2$, then the continued-fraction expansion of $\alpha$ is purely periodic. Let $d$ be the discriminant of $\alpha$, so $\alpha\in R_2(d)$. The Gauss map sends a reduced quadratic irrational to another reduced quadratic irrational with the same discriminant. Indeed, if $a_0=\lfloor\alpha\rfloor$, then $\alpha>1$ and $-1<\alpha'<0$ imply
\[
\alpha_1=\frac{1}{\alpha-a_0}>1,
\qquad
\alpha'_1=\frac{1}{\alpha'-a_0}\in(-1,0),
\]
and the discriminant is preserved by the $GL(2,\mathbb Z)$-action. Hence $\alpha_n\in R_2(d)$ for every $n\geq 0$. Since $R_2(d)$ is finite, there exist $0\leq \ell<\ell'$ such that $\alpha_\ell=\alpha_{\ell'}$. If $\ell=0$, then $\alpha_0=\alpha_{\ell'}$, and the same argument as above gives pure periodicity.

Assume $\ell\geq 1$. Then
\[
\alpha_{\ell-1}=[a_{\ell-1};\alpha_\ell],
\qquad
\alpha_{\ell'-1}=[a_{\ell'-1};\alpha_{\ell'}].
\]
By Lemma~\ref{lem:count-down},
\[
a_{\ell-1}
=
\left\lfloor-\frac{1}{\alpha'_\ell}\right\rfloor
=
\left\lfloor-\frac{1}{\alpha'_{\ell'}}\right\rfloor
=
a_{\ell'-1}.
\]
Thus $\alpha_{\ell-1}=\alpha_{\ell'-1}$. Repeating this step, we obtain
\[
\alpha_0=\alpha_{\ell'-\ell}.
\]
The same argument as before then gives
\[
\alpha=[\overline{a_0;a_1,\dots,a_{\ell'-\ell-1}}],
\]
so the expansion is purely periodic.
\end{proof}
\begin{rema}
In the proof that $\alpha\in R_2$ implies pure periodicity, one cannot simply assert at the outset that there exists $\ell$ with $\alpha_0=\alpha_\ell$. The finiteness of $R_2(d)$ alone does not force $\alpha_0$ to repeat. A priori, it could happen that $\alpha_0$ is distinct from every later $\alpha_i$, while the later $\alpha_i$ take only finitely many values. Lemma~\ref{lem:count-down} is used precisely to rule out this possibility.
\end{rema}
We end by recording a relation between a reduced quadratic irrational and the quadratic irrational obtained by reversing the period. This fact will be needed in later chapters.
\begin{prop}\label{prop:quadratic-conjugate}\index{Quadratic conjugate}
If, for some $n\geq1$,
\[
\alpha=[\overline{a_0;a_1,\dots,a_{n-1}}]\in R_2,
\quad\text{then}\quad
-\frac{1}{\alpha'}=[\overline{a_{n-1};a_{n-2},\dots,a_0}],
\]
where $\alpha'$ is the quadratic conjugate of $\alpha$.
\end{prop}
\begin{proof}
From $\alpha=[\overline{a_0;a_1,\dots,a_{n-1}}]$, we have
\[
\alpha
=
[a_0;a_1,\dots,a_{n-1},\alpha]
=
\begin{bmatrix}
p_{n-1}&p_{n-2}\\q_{n-1}&q_{n-2}
\end{bmatrix}\alpha
=
\frac{p_{n-1}\alpha+p_{n-2}}{q_{n-1}\alpha+q_{n-2}}.
\]
Thus $\alpha$ is a root of
\begin{equation}\label{eq:equation-solution-alpha}
q_{n-1}x^2+(q_{n-2}-p_{n-1})x-p_{n-2}=0.
\end{equation}
On the other hand,
\[
\begin{bmatrix}
p_{n-1}&p_{n-2}\\q_{n-1}&q_{n-2}
\end{bmatrix}
=
\begin{bmatrix}a_0&1\\1&0\end{bmatrix}
\begin{bmatrix}a_1&1\\1&0\end{bmatrix}
\cdots
\begin{bmatrix}a_{n-1}&1\\1&0\end{bmatrix}.
\]
Taking transposes gives
\[
\begin{bmatrix}
p_{n-1}&q_{n-1}\\p_{n-2}&q_{n-2}
\end{bmatrix}
=
\begin{bmatrix}a_{n-1}&1\\1&0\end{bmatrix}
\begin{bmatrix}a_{n-2}&1\\1&0\end{bmatrix}
\cdots
\begin{bmatrix}a_0&1\\1&0\end{bmatrix}.
\]
Put
\[
\beta:=[\overline{a_{n-1};a_{n-2},\dots,a_0}].
\]
Then
\[
\beta
=
[a_{n-1};a_{n-2},\dots,a_0,\beta]
=
\begin{bmatrix}
p_{n-1}&q_{n-1}\\p_{n-2}&q_{n-2}
\end{bmatrix}\beta
=
\frac{p_{n-1}\beta+q_{n-1}}{p_{n-2}\beta+q_{n-2}}.
\]
Rearranging, we obtain
\[
p_{n-2}\beta^2+(q_{n-2}-p_{n-1})\beta-q_{n-1}=0.
\]
Dividing by $-\beta^2$ gives
\[
-p_{n-2}
+
(q_{n-2}-p_{n-1})\left(-\frac1\beta\right)
+
q_{n-1}\left(-\frac1\beta\right)^2
=
0.
\]
Thus $-\frac1\beta$ is a root of \eqref{eq:equation-solution-alpha}. Since $a_{n-1}\geq 1$, Lemma~\ref{lem:integer-part} gives
\[
\beta>\lfloor\beta\rfloor=a_{n-1}\geq 1,
\quad\text{and hence}\quad
-1<-\frac1\beta<0.
\]
Therefore the conjugate root is $\alpha'=-\frac1\beta$. Equivalently,
\[
\beta=-\frac{1}{\alpha'}.
\]
This proves the claim.
\end{proof}

\chapter[Lagrange Spectrum]{Lagrange Spectrum}\label{chap:lagrange-spectrum}
Since Chapter~\ref{chap:continued-fraction} prepared the basic facts on continued fractions and quadratic irrationals, we now use them in this chapter to study the fundamental properties of the Lagrange spectrum. The Lagrange spectrum is the set of all Lagrange constants, which measure how well irrational numbers can be approximated by rational numbers, and it is one of the most fundamental objects in Diophantine approximation theory. Looking only at the definition, one might get no more than the impression that it is a multiplicative analogue of the irrationality exponent, which measures the quality of approximation in terms of powers of the denominator. However, when viewed through continued fraction expansions, one sees that its values are deeply connected with infinite sequences and periodicity.

We first define the Lagrange constant, give a computable expression for it, and compute basic examples. We then introduce a description in terms of bi-infinite sequences, which makes it easier to compare the Lagrange spectrum with the Markov spectrum in the next chapter. Finally, for quadratic irrationals, we show that the action of $GL(2,\mathbb Z)$ and the theory of periodic continued fractions developed in the previous chapter make it possible to compute the Lagrange constant explicitly.

Standard texts that are written with Markov's theorem in mind, such as \cite{aig,bombieri2,reutenauer}, often impose the assumption that the Lagrange constant is at most $3$ from the beginning. In this chapter, however, we do not impose such an assumption and work in the general setting.

The discussion in this chapter follows \cite{aig}.

\section{Definitions and First Examples}
In this section we introduce the Lagrange spectrum and check the simplest examples.
\begin{defi}\label{def:lagrange1770-spectrum}\index{Lagrange constant}\index{Lagrange spectrum}
Let $\alpha\in \mathbb R\setminus\mathbb Q$. We define $\mathcal L(\alpha)$ to be the supremum of all positive real numbers $L$ satisfying the following condition:
\begin{itemize}
    \item there exist infinitely many reduced fractions $p/q$ with $q>0$ such that
    \[
    \left| \alpha - \frac{p}{q}\right| < \frac{1}{L q^2}.
    \]
\end{itemize}
If this set is unbounded above, we set $\mathcal L(\alpha)=\infty$. We call $\mathcal L(\alpha)$ the \emph{Lagrange constant} of $\alpha$. The set
\[
\mathcal L := \{\mathcal L(\alpha) \mid \alpha \in \mathbb R\setminus \mathbb Q\}
\]
is called the \emph{Lagrange spectrum}. With this convention, $\mathcal L\subset\mathbb R\cup\{\infty\}$, and its finite part is $\mathcal L\cap\mathbb R$.
\end{defi}
Although rational numbers are not included in the definition, it is useful to first see what would happen if the same condition were applied to a rational number.
\begin{prop}\label{prop:rational-equal-0}
Let $\alpha$ be a rational number, and let $L$ be any positive real number. Then there are only finitely many reduced fractions $p/q$ with $q>0$ satisfying
\[
\left|\alpha - \frac{p}{q}\right|< \frac{1}{Lq^2}.
\]
\end{prop}
\begin{proof}
There is at most one reduced fraction $p/q$ equal to $\alpha$, so we assume below that $p/q\neq\alpha$. Write $\alpha=a/b$, where $a/b$ is reduced. Then
\[
\left|\alpha-\frac{p}{q}\right|=
\frac{|aq-bp|}{bq}\geq\frac{1}{bq}.
\]
Thus, if $\left|\alpha-p/q\right|<1/(Lq^2)$, then $1/(Lq^2)>1/(bq)$, and hence $q<b/L$. Since $q$ is a positive integer, there are only finitely many possible values of $q$. For each fixed such $q$, solving the same inequality for $p$ gives
\[
\frac{aq}{b}-\frac{1}{Lq}<p<\frac{aq}{b}+\frac{1}{Lq}.
\]
Since $p$ is an integer, there are only finitely many such $p$ for each fixed $q$. Hence only finitely many reduced fractions satisfy the inequality.
\end{proof}
The proposition shows that, for a rational number $\alpha$, the set of positive real numbers $L$ satisfying the condition in Definition~\ref{def:lagrange1770-spectrum} is empty. Thus rational numbers do not lead to a meaningful Lagrange constant in this sense. In the rest of this section we consider irrational numbers. We first give a characterization of the Lagrange constant in terms of continued fractions.
\begin{theorem}\label{thm:characterization-lagrange1770}\index{Complete quotient}\index{Lagrange constant}
Let $\alpha$ be an irrational number with infinite continued-fraction expansion
$\alpha=[a_0;a_1,\dots]$. Put
\[
\alpha_n:=[a_n;a_{n+1},\dots],\qquad
\beta_n:=[a_n;a_{n-1},\dots,a_1].
\]
Then
\begin{equation}\label{eq:alpha+beta}
\mathcal L(\alpha)=\limsup_{n\to\infty}
\left(\alpha_{n+1}+\frac{1}{\beta_n}\right).
\end{equation}
In what follows we write
\[
\lambda_n(\alpha):=\alpha_{n+1}+\frac{1}{\beta_n}.
\]
\end{theorem}
\begin{lemm}\label{lem:beta-n}
Let $\alpha=[a_0;a_1,\dots]$ be the infinite continued-fraction expansion of an irrational number. Then, for every $k\geq 1$,
\[
\frac{q_k}{q_{k-1}}=[a_k;a_{k-1},\dots,a_1](=\beta_k).
\]
\end{lemm}
\begin{proof}
We prove this by induction. For $k=1$, we have $q_1/q_0=a_1$. Assume the claim holds for $k-1$. Then
\[
\frac{q_k}{q_{k-1}}
=\frac{a_kq_{k-1}+q_{k-2}}{q_{k-1}}
=a_k+[0;a_{k-1},\dots,a_1]
=[a_k;a_{k-1},\dots,a_1].
\]
This proves the claim for $k$.
\end{proof}
\begin{proof}[Proof of Theorem~\ref{thm:characterization-lagrange1770}]
For the convergents of $\alpha$, we have
\[
\left|\alpha-\frac{p_n}{q_n}\right|
=\left|\frac{\alpha_{n+1}p_n+p_{n-1}}{\alpha_{n+1}q_n+q_{n-1}}-\frac{p_n}{q_n}\right|
=\frac{1}{\alpha_{n+1}q_n^2+q_{n-1}q_n}
=\frac{1}{\lambda_n(\alpha)q_n^2},
\]
where the last equality uses Lemma~\ref{lem:beta-n}.

We first show that $\mathcal L(\alpha)\le\limsup_{n\to\infty}\lambda_n(\alpha)$. Let $L>0$ be such that infinitely many reduced fractions $p/q$ satisfy
\[
\left|\alpha-\frac pq\right|<\frac{1}{Lq^2}.
\]
If $L>2$, Theorem~\ref{thm:good-is-approximation} shows that all these fractions are convergents. Hence, for infinitely many $n$,
\[
\frac{1}{\lambda_n(\alpha)q_n^2}<\frac{1}{Lq_n^2},
\]
so $\lambda_n(\alpha)>L$ infinitely often. Therefore $L\le\limsup_{n\to\infty}\lambda_n(\alpha)$.

If $L\le2$, the same conclusion follows from $\limsup_{n\to\infty}\lambda_n(\alpha)\ge2$. Indeed, if partial quotients $a_j\ge2$ occur infinitely often, then $\lambda_{j-1}(\alpha)>2$ infinitely often. Otherwise, $a_j=1$ for all sufficiently large $j$, and in this case
\[
\limsup_{j\to\infty}\lambda_j(\alpha)
=\lim_{j\to\infty}\lambda_j(\alpha)=\sqrt5>2.
\]
Thus $L\le\limsup_{n\to\infty}\lambda_n(\alpha)$ in every case. Taking the supremum over all such $L$ gives
\[
\mathcal L(\alpha)\le\limsup_{n\to\infty}\lambda_n(\alpha).
\]

For the opposite inequality, put
\[
R:=\limsup_{n\to\infty}\lambda_n(\alpha)\in[2,\infty].
\]
First suppose that $R<\infty$. For every $0<\varepsilon<R$, the definition of the limit superior gives infinitely many $n$ satisfying
\[
\lambda_n(\alpha)>R-\varepsilon.
\]
For these $n$,
\[
\left|\alpha-\frac{p_n}{q_n}\right|
<\frac{1}{(R-\varepsilon)q_n^2}.
\]
Hence $R-\varepsilon\le\mathcal L(\alpha)$. Letting $\varepsilon\downarrow0$ gives $R\le\mathcal L(\alpha)$.

Now suppose that $R=\infty$. For every $M>0$, infinitely many $n$ satisfy $\lambda_n(\alpha)>M$, and their convergents satisfy
\[
\left|\alpha-\frac{p_n}{q_n}\right|<\frac{1}{Mq_n^2}.
\]
Thus every $M>0$ satisfies the condition in Definition~\ref{def:lagrange1770-spectrum}, and $\mathcal L(\alpha)=\infty=R$. This proves the theorem.
\end{proof}

\begin{exam}\label{ex:concrete-example}
Let us compute some examples using Theorem~\ref{thm:characterization-lagrange1770}.
\begin{itemize}
    \item[(1)] Let $\alpha=(1+\sqrt{5})/2$. Since
    $\frac{1+\sqrt{5}}{2}=[1;\frac{1+\sqrt{5}}{2}]$, we have
    $\frac{1+\sqrt{5}}{2}=[\overline{1}]$. Therefore
    \[
    \mathcal L\left(\frac{1+\sqrt{5}}{2}\right)
    =\limsup_{n\to\infty}([\overline 1]+[0;1,\dots,1])
    =[\overline{1}]+\lim_{n\to\infty}[0;1,\dots,1]
    =[\overline 1]+[0;\overline{1}]=\sqrt{5}.
    \]
    \item[(2)] Let $\alpha=1+\sqrt{2}$. Since $1+\sqrt{2}=[2;1+\sqrt{2}]$, we have
    $1+\sqrt{2}=[\overline{2}]$. Therefore
    \[
    \mathcal L(1+\sqrt{2})
    =\limsup_{n\to\infty}([\overline 2]+[0;2,\dots,2])
    =[\overline{2}]+\lim_{n\to\infty}[0;2,\dots,2]
    =[\overline 2]+[0;\overline{2}]=2\sqrt{2}.
    \]
\end{itemize}
\end{exam}
We finish this section by observing that the Lagrange spectrum can be described using bi-infinite sequences and a limit superior. Let
$\mathbf a=(\dots,a_{-1},a_0,a_1,\dots)$ be a bi-infinite sequence with $a_i\in\mathbb Z_{\geq1}$ for every $i\in\mathbb Z$. Define
\[
\ell_n(\mathbf a):=[a_n;a_{n+1},\dots]+[0;a_{n-1},a_{n-2},\dots].
\]
Then the following result holds.
\begin{coro}[Perron's Identity]\label{cor:perron-formula}\index{Bi-infinite sequence}\index{Lagrange spectrum}\index{Perron's identity}
The Lagrange spectrum is characterized as
\[
\mathcal L=
\left\{
\limsup_{n\to+\infty}\ell_n(\mathbf a)
\ \middle|\
\mathbf a\in\mathbb Z_{\geq1}^{\mathbb Z}
\right\}.
\]
\end{coro}
We begin with the following lemma.

\begin{lemm}\label{lem:finite-infinite-comparison}
Let $a_0\in\mathbb Z$, let $a_1,\dots,a_m\in\mathbb Z_{\geq1}$ with $m\geq1$, and let $\xi>1$. Put
\[
x=[a_0;a_1,\dots,a_m],\qquad
y=[a_0;a_1,\dots,a_m,\xi].
\]
Write $x=p_m/q_m$, and let $q_{m-1}$ be the denominator of the preceding convergent. Then
\[
|x-y|
=\frac{1}{q_m(\xi q_m+q_{m-1})}
<\frac{1}{q_m(q_m+q_{m-1})}.
\]
\end{lemm}
\begin{proof}
By Proposition~\ref{prop:continued-fraction-matrix-representation},
\[
y=\frac{\xi p_m+p_{m-1}}{\xi q_m+q_{m-1}}.
\]
Therefore, using Lemma~\ref{lem:det}, we obtain
\[
\left|y-\frac{p_m}{q_m}\right|
=\frac{|p_{m-1}q_m-p_mq_{m-1}|}
 {q_m(\xi q_m+q_{m-1})}
=\frac{1}{q_m(\xi q_m+q_{m-1})}.
\]
The final inequality follows from $\xi>1$.
\end{proof}
\begin{proof}[Proof of Corollary~\ref{cor:perron-formula}]
Let the set on the right-hand side be
\[
\mathcal R:=
\left\{
\limsup_{n\to+\infty}\ell_n(\mathbf a)
\ \middle|\
\mathbf a\in\mathbb Z_{\geq1}^{\mathbb Z}
\right\}.
\]
We prove $\mathcal R=\mathcal L$.

First we show $\mathcal R\subset\mathcal L$. Take an arbitrary
\[
\mathbf a=(\ldots,a_{-1},a_0,a_1,\ldots)
\in \mathbb Z_{\geq1}^{\mathbb Z}
\quad\text{and put}\quad
\alpha:=[0;a_1,a_2,\ldots].
\]
Set
\[
u_n:=[0;a_{n-1},a_{n-2},\ldots,a_1],\qquad
v_n:=[0;a_{n-1},a_{n-2},\ldots,a_1,a_0,a_{-1},\ldots].
\]
By Theorem~\ref{thm:characterization-lagrange1770},
\[
\mathcal L(\alpha)=
\limsup_{n\to\infty}
\left([a_n;a_{n+1},a_{n+2},\ldots]+u_n\right),
\]
whereas
\[
\ell_n(\mathbf a)=
[a_n;a_{n+1},a_{n+2},\ldots]+v_n.
\]
Let $Q_{n-1}$ and $Q_{n-2}$ be the denominators of the last and the preceding convergents of $u_n$, respectively. Since $v_n$ is obtained by appending to $u_n$ an infinite continued fraction greater than $1$, Lemma~\ref{lem:finite-infinite-comparison} gives
\[
|u_n-v_n|
<\frac{1}{Q_{n-1}(Q_{n-1}+Q_{n-2})}
\leq\frac{1}{n(n-1)}
\qquad(n\geq2).
\]
For the last inequality, we used Corollary~\ref{cor:qk>k}, which gives $Q_{n-1}\geq n-1$, together with $Q_{n-2}\geq1$.
Therefore
\[
\left|
\ell_n(\mathbf a)-
\left([a_n;a_{n+1},a_{n+2},\ldots]+u_n\right)
\right|
\leq\frac{1}{n(n-1)}\to0
\qquad(n\to\infty),
\]
and so
\[
\limsup_{n\to\infty}\ell_n(\mathbf a)=\mathcal L(\alpha)\in\mathcal L.
\]
Thus $\mathcal R\subset\mathcal L$.

Next we show $\mathcal L\subset\mathcal R$. Take $r\in\mathcal L$. Then there is an irrational number
\[
\alpha=[a_0;a_1,a_2,\ldots]
\]
such that $r=\mathcal L(\alpha)$. Put
\[
\widetilde{\alpha}:=[0;a_1,a_2,\ldots].
\]
Since the expression in \eqref{eq:alpha+beta} does not depend on $a_0$, we have
\[
\mathcal L(\widetilde{\alpha})=\mathcal L(\alpha)=r.
\]
Define a bi-infinite sequence $\mathbf b=(b_n)_{n\in\mathbb Z}$ by
\[
b_n=
\begin{cases}
a_n & (n\geq1),\\
1   & (n\leq0).
\end{cases}
\]
Set
\[
u_n:=[0;a_{n-1},a_{n-2},\ldots,a_1],\qquad
w_n:=[0;a_{n-1},a_{n-2},\ldots,a_1,1,1,\ldots].
\]
By Theorem~\ref{thm:characterization-lagrange1770},
\[
\mathcal L(\widetilde{\alpha})=
\limsup_{n\to\infty}
\left([a_n;a_{n+1},a_{n+2},\ldots]+u_n\right),
\]
whereas
\[
\ell_n(\mathbf b)=
[a_n;a_{n+1},a_{n+2},\ldots]+w_n.
\]
Let $Q_{n-1}$ and $Q_{n-2}$ be the denominators of the last and the preceding convergents of $u_n$, respectively. Since $w_n$ is obtained by appending $[1;1,\dots]>1$ to $u_n$, Lemma~\ref{lem:finite-infinite-comparison} and Corollary~\ref{cor:qk>k} give
\[
|u_n-w_n|
<\frac{1}{Q_{n-1}(Q_{n-1}+Q_{n-2})}
\leq\frac{1}{n(n-1)}
\qquad(n\geq2).
\]
Hence
\[
\limsup_{n\to\infty}\ell_n(\mathbf b)
=\mathcal L(\widetilde{\alpha})=r.
\]
Thus $r\in\mathcal R$, and $\mathcal L\subset\mathcal R$. This proves the desired equality.
\end{proof}
Corollary~\ref{cor:perron-formula} was presented as Perron's characterization of $\mathcal L$, but for the purposes of this text it is not the most useful form. In its proof, the left-hand side of the bi-infinite sequence attached to an irrational number $\alpha$ was filled with $1$'s. Nothing essential depends on this choice: any sequence could have been placed on the left, and the same argument would still work. Thus Corollary~\ref{cor:perron-formula} mainly reformulates Theorem~\ref{thm:characterization-lagrange1770}. What will be more useful later is the formula \eqref{eq:alpha+beta} itself and the construction, from a Lagrange constant, of a bi-infinite sequence for which the relevant value is realized as a \emph{supremum}. This is the subject of the next section.

\section{A Supremum Construction from Bi-infinite Sequences}
In this section we associate to finite Lagrange constants certain bi-infinite sequences for which the relevant value is realized as a supremum. Notice that here we use a supremum, not a limit superior.  We begin with the following proposition.

\begin{prop}\label{prop:fin-or-infin}\index{Bounded partial quotients}\index{Partial quotient}
Let $\alpha$ be an irrational number with infinite continued-fraction expansion $\alpha=[a_0;a_1,\dots]$.
\begin{itemize}
    \item[(1)] If $(a_n)_{n=0}^\infty$ is bounded, then $\mathcal L(\alpha)<\infty$.
    \item[(2)] If $(a_n)_{n=0}^\infty$ is unbounded, then $\mathcal L(\alpha)=\infty$.
\end{itemize}
\end{prop}
\begin{proof}
For each $n\geq1$, put
\[
x_n:=\alpha_{n+1}+\frac{1}{\beta_n}.
\]
Lemmas~\ref{lem:integer-part-rational} and~\ref{lem:integer-part} imply that, for every $n\geq1$,
\[
a_{n+1}<x_n<a_{n+1}+2.
\]
Consequently,
\[
\limsup_{n\to\infty}a_{n+1}
\leq
\limsup_{n\to\infty}x_n
\leq
\limsup_{n\to\infty}a_{n+1}+2.
\]

The sequence $(a_n)_{n=0}^{\infty}$ is bounded if and only if $(a_{n+1})_{n=1}^{\infty}$ is bounded. Hence $(x_n)_{n=1}^{\infty}$ is bounded when $(a_n)_{n=0}^{\infty}$ is bounded, whereas $(x_n)_{n=1}^{\infty}$ is unbounded when $(a_n)_{n=0}^{\infty}$ is unbounded. Thus
\[
\limsup_{n\to\infty}x_n
\]
is finite in the bounded case and is $\infty$ in the unbounded case.

By Theorem~\ref{thm:characterization-lagrange1770},
\[
\mathcal L(\alpha)
=
\limsup_{n\to\infty}
\left(\alpha_{n+1}+\frac{1}{\beta_n}\right)
=
\limsup_{n\to\infty}x_n.
\]
This proves both assertions.
\end{proof}
It follows that $\mathcal L(\alpha)$ is a meaningful finite number only when the partial quotients of $\alpha$ are bounded. Thus, when studying finite elements of the Lagrange spectrum, we may assume that $(a_n)_{n=0}^\infty$ is bounded.
\begin{prop}
Let $\alpha$ be an irrational number with infinite continued-fraction expansion $\alpha=[a_0;a_1,\dots]$, and assume that $(a_n)_{n=0}^\infty$ is bounded. Put
\[
\alpha_n=[a_n;a_{n+1},\dots],\qquad
\beta_n=[a_n;a_{n-1},\dots,a_1],\qquad
\lambda_n(\alpha)=\alpha_{n+1}+\frac1{\beta_n}.
\]
Then the following hold:
\begin{itemize}
\item[(1)] $(\alpha_n)_{n=0}^\infty$ is bounded.
\item[(2)] $(\beta_n)_{n=1}^\infty$ is bounded.
\item[(3)] $(\lambda_n(\alpha))_{n=1}^\infty$ is bounded.
\end{itemize}
\end{prop}
\begin{proof}
Let $a$ be the maximum of $(a_n)_{n=1}^\infty$. For (1), Lemma~\ref{lem:integer-part} gives $a_k<\alpha_k<a_k+1$ for every $k\geq1$, and hence $1<\alpha_n<a+1$ for every $n\geq1$. Thus $(\alpha_n)_{n=0}^\infty$ is bounded. For (2), $\beta_1=a_1$, and for $n\geq2$,
\[
\beta_n=a_n+\frac1{[a_{n-1};a_{n-2},\dots,a_1]}.
\]
The denominator is at least $1$, so $a_n\leq\beta_n\leq a_n+1$ and $1\leq\beta_n\leq a+1$ for every $n\geq1$. This proves boundedness. Equality in the upper bound can occur, for example when $\beta_2=[a_2;1]=a_2+1$. Finally, (3) follows from (1) and (2), since
\[
1<\alpha_{n+1}+\frac1{\beta_n}<a+2
\]
for every $n\geq1$.
\end{proof}
Since $(\alpha_{n+1}+1/\beta_n)_{n=1}^\infty$ is bounded, the Bolzano--Weierstrass theorem gives a convergent subsequence
$(\alpha_{n_i+1}+1/\beta_{n_i})_{i=1}^\infty$. Let its limit be $r$. Since $(\beta_{n_i})_{i=1}^\infty$ is also bounded, it has a further convergent subsequence; write its limit as $\eta$. Along the same indices, the corresponding subsequence of $(\alpha_{n+1})_{n=1}^\infty$ converges to $r-1/\eta$. Put
\[
\theta:=r-\frac1\eta.
\]
The pair $(\theta,\eta)$ obtained in this way from the accumulation point $r$ will be called a \emph{pair associated with the accumulation point $r$}. By the definition of the limit superior, $\mathcal L(\alpha)$ is the supremum of the accumulation points $r$ obtained from such subsequences. To avoid excessive subscripts, we shall henceforth denote the chosen subsequence simply by the indices $n_i$.

\begin{prop}\label{prop:theta-eta}\index{Accumulation point}
Let $\alpha$ be an irrational number, and let $(a_n)_{n=0}^\infty$ be the sequence giving its infinite continued-fraction expansion. Assume that $(a_n)_{n=1}^\infty$ is bounded, and let $a$ be the maximum of this sequence, excluding $a_0$. Let $r$ be an accumulation point of $(\lambda_n(\alpha))_{n=1}^\infty$, and let $(\theta,\eta)$ be a pair associated with $r$. Then $\theta$ and $\eta$ are irrational numbers. If
\[
\theta=[b_0;b_1,\dots],\qquad
\eta=[b_{-1};b_{-2},\dots],
\quad\text{then}\quad
1\leq b_i\leq a
\]
for every $i\in\mathbb Z$. In particular, $1<\theta,\eta<a+1$.
\end{prop}
\begin{proof}
First, every complete quotient with $m\geq1$ satisfies
\begin{equation}\label{eq:complete-quotient-bounds}
1+\frac1{a+1}<\alpha_m<a+1.
\end{equation}
This follows from $\alpha_m=a_m+1/\alpha_{m+1}$, $1\leq a_m\leq a$, and $1<\alpha_{m+1}<a+1$.

We first prove that $\theta$ is irrational. Suppose otherwise, and write
\[
\theta=[b_0;b_1,\dots,b_k]=p_k/q_k.
\]
If $k=0$, then $\theta$ is an integer. However, for every $i$,
\[
\frac1{a+1}<\alpha_{n_i+1}-\lfloor\alpha_{n_i+1}\rfloor
=\frac1{\alpha_{n_i+2}}<\frac{a+1}{a+2}.
\]
Thus $\alpha_{n_i+1}$ stays at distance greater than $1/(a+2)$ from every integer, contradicting convergence to $\theta$.

Assume $k\geq1$, and choose $i$ large enough that $|\alpha_{n_i+1}-\theta|<1/(2q_k^2)$. Theorem~\ref{thm:good-is-approximation} implies that $\theta$ is a convergent of $\alpha_{n_i+1}$. Passing to a further subsequence if necessary, the same one of the two representations in Remark~\ref{rem:two-finite-expansions} occurs each time. Its initial partial quotients are
\[
(b_0,b_1,\dots,b_k)\quad\text{or}\quad(b_0,b_1,\dots,b_{k-1},b_k-1,1).
\]
Denote this block by $(c_0,\dots,c_j)$ and the denominator of its preceding convergent by $q'$. The current denominator is $q_k$, and $0<q'<q_k$. Hence
\[
\alpha_{n_i+1}=[c_0;c_1,\dots,c_j,\alpha_{n_i+j+2}].
\]
The calculation used in \eqref{eq:difference-approximation}, together with \eqref{eq:complete-quotient-bounds}, gives
\[
|\alpha_{n_i+1}-\theta|=\frac1{q_k(\alpha_{n_i+j+2}q_k+q')}
>\frac1{(a+2)q_k^2}.
\]
This uniform positive lower bound contradicts convergence to $\theta$.

Now write $\theta=[b_0;b_1,\dots]$. For any fixed $r\geq0$, Lemma~\ref{lem:stability-of-initial-partial-quotients} implies that its first $r+1$ partial quotients agree with those of $\alpha_{n_i+1}$ for large $i$. Thus $b_r=a_{n_i+r+1}$ and $1\leq b_r\leq a$.

Next suppose that $\eta$ is rational, and write
\[
\eta=[b_{-1};b_{-2},\dots,b_{-k}]=p_{-k}/q_{-k}.
\]
If $k=1$, then $\eta$ is an integer. For $n\geq3$, the bounds $1\leq[a_{n-2};\dots,a_1]\leq a+1$ give
\[
1+\frac1{a+1}\leq\beta_{n-1}\leq a+1,
\quad\text{and therefore}\quad
\frac1{a+1}\leq\beta_n-\lfloor\beta_n\rfloor
=\frac1{\beta_{n-1}}\leq\frac{a+1}{a+2}.
\]
Thus $\beta_n$ stays at distance at least $1/(a+2)$ from every integer, contradicting $\beta_{n_i}\to\eta$.

Assume $k\geq2$, and choose $i$ large enough that
\[
|\beta_{n_i}-\eta|<\frac1{2q_{-k}^2},\qquad n_i\geq k+3.
\]
If $a_1\geq2$, the finite regular expansion is $[a_{n_i};a_{n_i-1},\dots,a_1]$. If $a_1=1$, it is instead
\[
\beta_{n_i}=[a_{n_i};a_{n_i-1},\dots,a_3,a_2+1].
\]
In either case its first $k+1$ partial quotients agree with the original reversed expression, and it has at least $k+2$ partial quotients. In particular, $\beta_{n_i}\neq\eta$. Theorem~\ref{thm:good-is-approximation} and Remark~\ref{rem:finite-remark} show that $\eta$ is a proper convergent of $\beta_{n_i}$. By Remark~\ref{rem:two-finite-expansions}, the initial block representing it is
\[
(b_{-1},b_{-2},\dots,b_{-k})
\quad\text{or}\quad
(b_{-1},b_{-2},\dots,b_{-(k-1)},b_{-k}-1,1).
\]
Pass to a further subsequence on which the same representation occurs, denote it by $(c_0,\dots,c_j)$, and let $q'$ be the preceding denominator. At least one partial quotient remains after this block. Writing the remaining finite continued fraction as $\zeta_i$, we have
\[
\beta_{n_i}=[c_0;c_1,\dots,c_j,\zeta_i],\qquad1<\zeta_i\leq a+1.
\]
Even if only one partial quotient remains, it is at least $2$: it is $a_1\geq2$, or $a_2+1\geq2$ when $a_1=1$. Since the current denominator is $q_{-k}$ and $0<q'<q_{-k}$, the same calculation gives
\[
|\beta_{n_i}-\eta|=\frac1{q_{-k}(\zeta_iq_{-k}+q')}
>\frac1{(a+2)q_{-k}^2}.
\]
This contradicts convergence to $\eta$, so $\eta$ is irrational.

Write $\eta=[b_{-1};b_{-2},\dots]$. For any fixed $r\geq1$ and sufficiently large $i$, the finite regular expansion of $\beta_{n_i}$ has at least $r$ partial quotients, whose first $r$ entries are $a_{n_i},a_{n_i-1},\dots,a_{n_i-r+1}$. By Lemma~\ref{lem:stability-of-initial-partial-quotients}, they agree with those of $\eta$. Thus $b_{-r}=a_{n_i-r+1}$ and $1\leq b_{-r}\leq a$.

We have proved $1\leq b_i\leq a$ for all $i\in\mathbb Z$. Lemma~\ref{lem:integer-part} then gives $1<\theta,\eta<a+1$.
\end{proof}

Let $r$ be an accumulation point of $(\lambda_n(\alpha))_{n=1}^\infty$, and let $(\theta,\eta)$ be a pair associated with it. Write
\[
\theta=[b_0;b_1,\dots],\qquad
\eta=[b_{-1};b_{-2},\dots],
\]
and consider the bi-infinite sequence
\[
\mathbf b=(\dots,b_{-1},b_0,b_1,\dots).
\]
We call this the \emph{bi-infinite sequence determined by $(\theta,\eta)$}. Then
\[
\ell_0(\mathbf b)
=[b_0;b_1,\dots]+[0;b_{-1},b_{-2},\dots]
=\theta+\frac1\eta=r.
\]
\begin{theorem}\label{thm:limit-point}\index{Accumulation point}
Let $\alpha$ be an irrational number whose partial quotients are bounded, and consider the sequence $(\lambda_n(\alpha))_{n=1}^\infty$. Let $r$ be any accumulation point of this sequence, and let
\[
\mathbf b=(\dots,b_{-2},b_{-1},b_0,b_1,\dots)
\]
be the bi-infinite sequence determined by a pair $(\theta,\eta)$ associated with $r$. For any $h\in\mathbb Z$, put
\[
\theta':=[b_h;b_{h+1},\dots],\qquad
\eta':=[b_{h-1};b_{h-2},\dots]
\quad\text{and define}\quad
r':=\ell_h(\mathbf b)=\theta'+\frac1{\eta'}.
\]
Then $r'$ is also an accumulation point of $(\lambda_n(\alpha))_{n=1}^\infty$.
\end{theorem}
\begin{proof}
We prove the case $h>0$. Since
\[
\theta=[b_0;b_1,\dots,b_{h-1},\theta'],
\]
if $p_k/q_k$ denotes the convergents of $\theta$, then
\begin{equation}\label{eq:theta}
\theta=\frac{\theta'p_{h-1}+p_{h-2}}{\theta'q_{h-1}+q_{h-2}}.
\end{equation}
Choose indices $n_i$ such that $\alpha_{n_i+1}\to\theta$ and $\beta_{n_i}\to\eta$. By the argument in the proof of Proposition~\ref{prop:theta-eta}, for all sufficiently large $i$ the first $h$ partial quotients of $\alpha_{n_i+1}$ agree with those of $\theta$. Thus
\[
\alpha_{n_i+1}=[b_0;b_1,\dots,b_{h-1},\alpha_{n_i+h+1}],
\]
and hence
\begin{equation}\label{eq:alphani+1}
\alpha_{n_i+1}
=\frac{\alpha_{n_i+h+1}p_{h-1}+p_{h-2}}{\alpha_{n_i+h+1}q_{h-1}+q_{h-2}}.
\end{equation}
Solving \eqref{eq:theta} for $\theta'$ and \eqref{eq:alphani+1} for $\alpha_{n_i+h+1}$, we obtain integers $A,B,C,D$ such that
\[
\theta'=\frac{\theta A+B}{\theta C+D},\qquad
\alpha_{n_i+h+1}=\frac{\alpha_{n_i+1}A+B}{\alpha_{n_i+1}C+D}.
\]
Since $\alpha_{n_i+1}\to\theta$, it follows that $\alpha_{n_i+h+1}\to\theta'$.

For $\eta$ we use
\[
\eta'=[b_{h-1};b_{h-2},\dots,b_0,\eta],\qquad
\beta_{n_i+h}=[b_{h-1};b_{h-2},\dots,b_0,\beta_{n_i}].
\]
The same argument gives $\beta_{n_i+h}\to\eta'$. Therefore
\[
\lambda_{n_i+h}(\alpha)\to \theta'+\frac1{\eta'}=r'.
\]
Thus $r'$ is an accumulation point. The case $h<0$ is proved in the same way, with the roles of $\theta$ and $\eta$ interchanged.
\end{proof}
Let $\mathscr A$ be the set of all bi-infinite sequences $(a_i)_{i=-\infty}^\infty$ with $a_i\in\mathbb Z_{\geq1}$ for every $i\in\mathbb Z$. For $\mathbf a\in\mathscr A$, define
\[
\mathcal S(\mathbf a):=\sup_{h\in\mathbb Z}\ell_h(\mathbf a),
\]
and set $\mathcal S(\mathbf a)=\infty$ when the set is unbounded above. The preceding discussion gives the following theorem.
\begin{theorem}\label{thm:LsubsetS}\index{Bi-infinite sequence}\index{Lagrange spectrum}
Let
\[
\mathcal S:=\{\mathcal S(\mathbf a)\mid \mathbf a\in\mathscr A\}.
\]
Then $\mathcal L\subset\mathcal S$.
\end{theorem}
\begin{proof}
For $\infty\in\mathcal L$, take a bi-infinite sequence $\mathbf a$ whose right-hand partial quotients $a_h$ are unbounded. Since $\alpha_h=[a_h;a_{h+1},\dots]$ satisfies $\alpha_h>a_h$, we have
\[
\sup_{h\in\mathbb Z}\ell_h(\mathbf a)\ge
\sup_{h\in\mathbb Z}\alpha_h=\infty.
\]
Thus $\infty\in\mathcal S$.

Now suppose $\mathcal L(\alpha)=r<\infty$. By the proposition above, the partial quotients of $\alpha$ are bounded. The number $r$ can be taken as an accumulation point of $(\lambda_n(\alpha))_{n=1}^\infty$. Let $\mathbf b=(\dots,b_{-2},b_{-1},b_0,b_1,b_2,\dots)$ be the bi-infinite sequence determined by a pair associated with $r$. Then $\ell_0(\mathbf b)=r$. By Theorem~\ref{thm:limit-point}, for every $h\in\mathbb Z$, the number $\ell_h(\mathbf b)$ is an accumulation point of $(\lambda_n(\alpha))_{n=1}^\infty$. By Theorem~\ref{thm:characterization-lagrange1770}, $r$ is the supremum of all accumulation points of this sequence. Hence
\[
r=\mathcal S(\mathbf b)\in\mathcal S.
\]
This proves the assertion.
\end{proof}
We end this section with several cautions. In general, it is difficult to compute an accumulation point $r$ of $(\lambda_n(\alpha))_{n=1}^\infty$ and a pair $(\theta,\eta)$ associated with it directly from the infinite continued-fraction expansion of an irrational number $\alpha$. Moreover, the pair $(\theta,\eta)$ associated with $r$ need not be unique, because there may be choices in taking $\eta$ (or, if the order of the two limits is reversed, in taking $\theta$). Thus a single accumulation point $r$ may give rise to more than one bi-infinite sequence.

Theorem~\ref{thm:limit-point} says that once one fixes a bi-infinite sequence associated with a single accumulation point $r$, one can produce further accumulation points $r'$ from it. However, it does not assert that all accumulation points of $(\lambda_n(\alpha))_{n=1}^\infty$ are obtained from that one bi-infinite sequence.

It might seem that Theorem~\ref{thm:LsubsetS} supplies this missing assertion, but it does not. That theorem only says that if one takes a bi-infinite sequence associated with the accumulation point which realizes the Lagrange constant, then the Lagrange constant is the supremum of the accumulation points obtained from that particular sequence. Therefore, even if one finds a bi-infinite sequence $\mathbf b$ arising from an accumulation point $r$ and an associated pair, it need not be true that $\mathcal S(\mathbf b)=\mathcal L(\alpha)$.

For these reasons, the material in this section alone does not provide a practical method for computing $\mathcal L(\alpha)$ for a general irrational number $\alpha$.

\section{Lagrange Constants of Quadratic Irrationals}
At the end of the previous section we explained that several difficulties make it hard to compute the Lagrange constant of a general irrational number by means of bi-infinite sequences. However, when the infinite continued-fraction expansion of $\alpha$ has a particularly simple form, these difficulties can be overcome. The simple form in question is pure periodicity, or equivalently, the case where $\alpha$ is a reduced quadratic irrational. The same computation also gives the Lagrange constant of a non-reduced quadratic irrational. We explain this in this section.

We first show that the study of quadratic irrationals can be reduced to the study of reduced quadratic irrationals. The following proposition actually holds for arbitrary irrational numbers.
\begin{prop}\label{prop:equivalent-lagrange1770}\index{Unimodular equivalence!of irrational numbers}
Let $\alpha$ and $\beta$ be $GL(2,\mathbb Z)$-equivalent irrational numbers. Then
\[
\mathcal L(\alpha)=\mathcal L(\beta).
\]
\end{prop}
\begin{proof}
Write
\[
\alpha=[a_0;a_1,\ldots],\qquad
\beta=[b_0;b_1,\ldots].
\]
By Theorem~\ref{thm:characterization-equivalent}, there exist $n_0,m_0\in\mathbb Z_{\ge0}$ such that $a_{n_0+h}=b_{m_0+h}$ for every $h\ge0$. Choose $r\in\mathbb Z_{\ge0}$ so that $n:=n_0+r\ge1$ and $m:=m_0+r\ge1$. Then $a_{n+h}=b_{m+h}$ for every $h\ge0$. Write the common tail as $[c_0;c_1,c_2,\ldots]$. Thus the two expansions have the same tail after finitely many initial terms have been deleted.

For $j\ge1$, put
\begin{align*}
\rho_j&:=[c_j;c_{j+1},\ldots]+[0;c_{j-1},\ldots,c_0,a_{n-1},\ldots,a_1],\\
\sigma_j&:=[c_j;c_{j+1},\ldots]+[0;c_{j-1},\ldots,c_0,b_{m-1},\ldots,b_1].
\end{align*}
Since $n,m\ge1$, these backward tails use the same index convention as the sequences computing the Lagrange constants. Apart from finitely many initial terms, the sequences computing $\mathcal L(\alpha)$ and $\mathcal L(\beta)$ are $(\rho_j)_{j\ge1}$ and $(\sigma_j)_{j\ge1}$, respectively. Hence
\[
\mathcal L(\alpha)=\limsup_{j\to\infty}\rho_j,
\qquad
\mathcal L(\beta)=\limsup_{j\to\infty}\sigma_j.
\]

The first terms of $\rho_j$ and $\sigma_j$ are equal, so their difference comes only from the second terms. For $j\ge3$, put
\[
u_j:=[0;c_{j-1},c_{j-2},\ldots,c_2],
\]
and let $Q_{j-2}$ and $Q_{j-3}$ be the denominators of its last and preceding convergents. Also put
\[
\xi_\alpha:=[c_1;c_0,a_{n-1},\ldots,a_1],\qquad
\xi_\beta:=[c_1;c_0,b_{m-1},\ldots,b_1],
\]
omitting either finite list when it is empty. Since $c_0,c_1\ge1$, we have $\xi_\alpha,\xi_\beta>1$, and
\begin{align*}
[0;c_{j-1},\ldots,c_0,a_{n-1},\ldots,a_1]
&=[0;c_{j-1},\ldots,c_2,\xi_\alpha],\\
[0;c_{j-1},\ldots,c_0,b_{m-1},\ldots,b_1]
&=[0;c_{j-1},\ldots,c_2,\xi_\beta].
\end{align*}
Applying Lemma~\ref{lem:finite-infinite-comparison} to compare each expression with $u_j$, and then using the triangle inequality, gives
\[
|\rho_j-\sigma_j|
<\frac{2}{Q_{j-2}(Q_{j-2}+Q_{j-3})}
\le\frac{2}{(j-2)(j-1)}.
\]
For the last inequality, Corollary~\ref{cor:qk>k} gives $Q_{j-2}\ge j-2$, and $Q_{j-3}\ge1$. Thus $|\rho_j-\sigma_j|\to0$. Two real sequences whose difference tends to zero have the same limit superior, so $\mathcal L(\alpha)=\mathcal L(\beta)$.
\end{proof}
\begin{coro}
Let $\alpha$ be a quadratic irrational with continued-fraction expansion
\[
\alpha=[a_0;a_1,\dots,a_k,\overline{c_1,\dots,c_n}].
\]
Define reduced quadratic irrationals
\[
\gamma_1:=[\overline{c_1,c_2,\dots,c_n}],\quad
\gamma_2:=[\overline{c_2,c_3,\dots,c_n,c_1}],\quad\dots,
\gamma_n:=[\overline{c_n,c_1,\dots,c_{n-1}}].
\]
Then
\[
\mathcal L(\alpha)=\mathcal L(\gamma_1)=\mathcal L(\gamma_2)=\cdots=\mathcal L(\gamma_n).
\]
\end{coro}
\begin{proof}
This follows from Theorem~\ref{thm:characterization-equivalent} and Proposition~\ref{prop:equivalent-lagrange1770}.
\end{proof}
The preceding argument shows that, in order to study the Lagrange constant of a quadratic irrational, it suffices to study an equivalent reduced quadratic irrational. The main theorem of this section is the following.
\begin{theorem}\label{thm:L=S}\index{Purely periodic continued fraction}\index{Lagrange constant}\index{Periodic continued fraction}
Let $k\geq1$ and let $b_0,\dots,b_{k-1}\in\mathbb Z_{\geq1}$. Suppose that a reduced quadratic irrational $\alpha$ has infinite continued-fraction expansion
\[
\alpha=[\overline{b_0,b_1,\dots,b_{k-1}}].
\]
Let $\mathbf b$ be the bi-infinite sequence obtained by repeating the period $b_0,b_1,\dots,b_{k-1}$ indefinitely in both directions:
\[
\mathbf b=(\dots,b_0,b_1,\dots,b_{k-1},b_0,b_1,\dots,b_{k-1},b_0,b_1,\dots,b_{k-1},\dots).
\]
Then
\[
\mathcal L(\alpha)=\mathcal S(\mathbf b).
\]
\end{theorem}
\begin{proof}
Extend $b_j$ periodically with period $k$ to all $j\in\mathbb Z$, and write $a_n=b_n$ for the partial quotients of $\alpha$. For $0\leq i<k$, put
\[
\theta_i=[b_i;b_{i+1},b_{i+2},\dots],\qquad
\eta_i=[b_{i-1};b_{i-2},b_{i-3},\dots].
\]
Restrict $n\geq1$ to the residue class $n+1\equiv i\pmod{k}$. Then
\[
\alpha_{n+1}=\theta_i,\qquad
\beta_n=[a_n;a_{n-1},\dots,a_1]=[b_{i-1};b_{i-2},\dots,b_{i-n}].
\]
The last expression is the truncation after the first $n$ terms of the infinite continued fraction for $\eta_i$. By Theorem~\ref{thm:limit-existence}, $\beta_n\to\eta_i$ along this residue class. Hence
\[
\lambda_n(\alpha)=\alpha_{n+1}+\frac1{\beta_n}
\longrightarrow r_i:=\theta_i+\frac1{\eta_i}=\ell_i(\mathbf b)
\qquad(n+1\equiv i\pmod{k}).
\]
These finitely many residue classes partition the indices, and each contains infinitely many indices. Therefore
\[
\mathcal L(\alpha)=\limsup_{n\to\infty}\lambda_n(\alpha)=\max_{0\leq i<k}r_i.
\]
Since $\ell_i(\mathbf b)$ is periodic in $i$ with period $k$, this maximum equals $\mathcal S(\mathbf b)$. Distinct residue classes may have the same limit; no assumption that $k$ is the least period is needed.
\end{proof}
As a consequence we obtain the following statement.
\begin{coro}\label{cor:quadratic-lagrange1770}
Let $\alpha$ be a quadratic irrational with infinite continued-fraction expansion
\[
\alpha=[a_0;a_1,\dots,a_n,\overline{b_0,b_1,\dots,b_{k-1}}]
\]
for some $n\geq0$ and $k\geq1$.
Let $\mathbf b$ be the bi-infinite sequence obtained by repeating the period $b_0,b_1,\dots,b_{k-1}$ indefinitely in both directions:
\[
\mathbf b=(\dots,b_0,b_1,\dots,b_{k-1},b_0,b_1,\dots,b_{k-1},b_0,b_1,\dots,b_{k-1},\dots).
\]
Then
\[
\mathcal L(\alpha)=\mathcal S(\mathbf b).
\]
\end{coro}
In this corollary, the candidates for the supremum value of $\mathcal S(\mathbf b)$ are finite in number, so $\mathcal L(\alpha)$ can in principle be computed by hand.

We end this chapter by expressing the Lagrange constant of $\alpha$ in terms of continued-fraction matrices associated with the continued-fraction expansion. This makes the computation of $\mathcal L(\alpha)$ still easier. For a finite sequence $(a_0,\dots,a_k)$, define
\[
F_{(a_0,a_1,\dots,a_k)}:=
\begin{bmatrix} a_0&1\\1&0 \end{bmatrix}
\begin{bmatrix} a_1&1\\1&0 \end{bmatrix}
\cdots
\begin{bmatrix} a_k&1\\1&0 \end{bmatrix}.
\]
\begin{theorem}\label{thm:lagrange1770-continued-fraction-matrix}\index{Continued-fraction matrix}
Let $\alpha$ be a quadratic irrational, and suppose that a nonempty period of its infinite continued-fraction expansion is $(b_0,\dots,b_{k-1})$, so $k\geq1$. For $0\leq i\leq k-1$, put
\[
S_i:=(b_i,b_{i+1},\dots,b_{k-1},b_0,\dots,b_{i-1}).
\]
Then
\[
\mathcal L(\alpha)=
\max\left\{
\frac{\sqrt{(\mathrm{tr}(F_{S_i}))^2-(-1)^k\cdot4}}{(F_{S_i})_{21}}
\ \middle|\
0\leq i\leq k-1
\right\}.
\]
Here $(F_{S_i})_{21}$ denotes the $(2,1)$-entry of $F_{S_i}$.
\end{theorem}
\begin{proof}
Read the indices modulo $k$. Put
\[
\theta_i=[\overline{b_i,b_{i+1},\dots,b_{k-1},b_0,\dots,b_{i-1}}],
\quad\text{and}\quad
\eta_i=[\overline{b_{i-1},b_{i-2},\dots,b_0,b_{k-1},\dots,b_i}].
\]
By Corollary~\ref{cor:quadratic-lagrange1770}, it suffices to show
\[
\theta_i+\frac1{\eta_i}
=\frac{\sqrt{(\mathrm{tr}(F_{S_i}))^2-(-1)^k\cdot4}}{(F_{S_i})_{21}}.
\]
Since
\[
\theta_i=[b_i;b_{i+1},\dots,b_{i-1},\theta_i],
\]
Proposition~\ref{prop:continued-fraction-matrix-representation} gives, where $F_{S_i}$ acts by fractional linear transformations,
\[
\theta_i=F_{S_i}\theta_i.
\]
Write
\[
F_{S_i}=\begin{bmatrix}
    p_{k-1}&p_{k-2}\\q_{k-1}&q_{k-2}
\end{bmatrix}.
\]
From the definition of the action,
\[
\theta_i=\frac{\theta_i p_{k-1}+p_{k-2}}{\theta_i q_{k-1}+q_{k-2}}.
\]
Solving this quadratic equation gives
\begin{align*}
\theta_i&=\frac{p_{k-1}-q_{k-2}+\sqrt{(p_{k-1}-q_{k-2})^2+4p_{k-2}q_{k-1}}}{2q_{k-1}}\\
&=\frac{p_{k-1}-q_{k-2}+\sqrt{(p_{k-1}+q_{k-2})^2-4(p_{k-1}q_{k-2}-p_{k-2}q_{k-1})}}{2q_{k-1}}\\
&=\frac{p_{k-1}-q_{k-2}+\sqrt{(\mathrm{tr}(F_{S_i}))^2-(-1)^k\cdot4}}{2(F_{S_i})_{21}}.
\end{align*}
We take the plus sign in front of the square root because $\theta_i$ is larger than its quadratic conjugate $\theta_i'$. On the other hand, by Proposition~\ref{prop:quadratic-conjugate},
\[
\eta_i=-\frac1{\theta_i'}.
\]
Therefore
\[
\frac1{\eta_i}
=-\frac{p_{k-1}-q_{k-2}-\sqrt{(\mathrm{tr}(F_{S_i}))^2-(-1)^k\cdot4}}{2(F_{S_i})_{21}}.
\]
Adding the two expressions gives the desired formula.
\end{proof}
In the proof of Theorem~\ref{thm:lagrange1770-continued-fraction-matrix}, the numerator in the expression for $\theta_i+1/\eta_i$ depends on the trace of $F_{S_i}$. In fact this trace is independent of $i$. Indeed,
\[
F_{S_{i+1}}=
\begin{bmatrix} b_i&1\\1&0 \end{bmatrix}^{-1}
F_{S_i}
\begin{bmatrix} b_i&1\\1&0 \end{bmatrix},
\]
where $S_k:=S_0$. Since $\mathrm{tr}(AB)=\mathrm{tr}(BA)$, we have
\begin{align*}
\mathrm{tr}(F_{S_{i+1}})
=\mathrm{tr}\left(
\begin{bmatrix} b_i&1\\1&0 \end{bmatrix}^{-1}
F_{S_i}
\begin{bmatrix} b_i&1\\1&0 \end{bmatrix}
\right)=\mathrm{tr}\left(
F_{S_i}
\begin{bmatrix} b_i&1\\1&0 \end{bmatrix}
\begin{bmatrix} b_i&1\\1&0 \end{bmatrix}^{-1}
\right)
=\mathrm{tr}(F_{S_i}).
\end{align*}
Thus, when computing the candidates for $\mathcal L(\alpha)$, the numerator only has to be computed once. Since this common numerator is positive, the maximum is attained precisely when the denominator $(F_{S_i})_{21}$ is minimal. We record this as a proposition.
\begin{prop}
In the setting of Theorem~\ref{thm:lagrange1770-continued-fraction-matrix}, let $j\in\{0,\dots,k-1\}$ satisfy
\[
\min_{0\leq i\leq k-1}\{(F_{S_i})_{21}\}=(F_{S_j})_{21}.
\]
Then
\[
\mathcal L(\alpha)=
\frac{\sqrt{(\mathrm{tr}(F_{S_j}))^2-(-1)^k\cdot4}}{(F_{S_j})_{21}}.
\]
\end{prop}
Let us compute the Lagrange constant for a concrete quadratic irrational.
\begin{exam}\label{ex:1+sqrt3-lagrange1770}
We compute the Lagrange constant of $\alpha=1+\sqrt{3}$. First we find its infinite continued-fraction expansion. Since
\[
\alpha-2=\sqrt{3}-1=\frac{2}{\sqrt{3}+1},
\]
we have
\[
\frac1{\alpha-2}=\frac{\sqrt{3}+1}{2}=1+\frac{\sqrt{3}-1}{2}=1+\frac1\alpha.
\]
Therefore
\[
\alpha=2+\frac1{1+\frac1\alpha}.
\]
Hence $\alpha=[\overline{2,1}]$, the period is $(2,1)$, and its length is $k=2$. We now compute $\mathcal L(\alpha)$ using Theorem~\ref{thm:lagrange1770-continued-fraction-matrix}.

The rotations of the period are
\[
S_0=(2,1),\qquad S_1=(1,2).
\]
The corresponding continued-fraction matrices are
\[
F_{(2,1)}=
\begin{bmatrix}2&1\\1&0\end{bmatrix}
\begin{bmatrix}1&1\\1&0\end{bmatrix}
=
\begin{bmatrix}3&2\\1&1\end{bmatrix},
\qquad
F_{(1,2)}=
\begin{bmatrix}1&1\\1&0\end{bmatrix}
\begin{bmatrix}2&1\\1&0\end{bmatrix}
=
\begin{bmatrix}3&1\\2&1\end{bmatrix}.
\]
In both cases $\mathrm{tr}(F_{S_i})=4$. On the other hand,
\[
(F_{(2,1)})_{21}=1,
\qquad
(F_{(1,2)})_{21}=2.
\]
Since $k=2$, we have $(-1)^k=1$. Hence Theorem~\ref{thm:lagrange1770-continued-fraction-matrix} gives
\[
\mathcal L(\alpha)
=\max_{i=0,1}
\left\{
\frac{\sqrt{(\mathrm{tr}(F_{S_i}))^2-4}}{(F_{S_i})_{21}}
\right\}
=\max\left\{
\frac{\sqrt{12}}{1},\frac{\sqrt{12}}{2}
\right\}
=2\sqrt{3}.
\]
Thus $\mathcal L(1+\sqrt{3})=2\sqrt{3}$.
\end{exam}
The preceding discussion shows that, for the Lagrange spectrum of quadratic irrationals, once we know which cyclic cut of the period makes $(F_{S_j})_{21}$ minimal, the value can be computed. At present, however, the most direct way to determine this cut is simply to compute all the candidates and compare them. For certain classes of quadratic irrationals there are methods that avoid this brute-force comparison, and these classes will be studied more deeply in Part~II.

\chapter{Markov Spectrum}\label{chap:markov-spectrum}
In this chapter we discuss the Markov spectrum. The Markov spectrum arises from a minimization problem for indefinite binary quadratic forms, and at first sight it may look quite different from the Lagrange spectrum studied in the preceding chapter. Once both spectra are rewritten in terms of continued fractions and bi-infinite sequences, however, their structures become very similar.

We first define the Markov constant attached to a binary quadratic form and examine its meaning through concrete examples. Next, using canonical reduced binary quadratic forms and the action of the unimodular group, we choose representatives in each class of quadratic forms for which the Markov constant is easier to compute. Then, in analogy with the Lagrange constant from the preceding chapter, we express the Markov constant by means of a bi-infinite sequence, so that the two spectra can be compared within a common framework. Finally, we show that the Markov constant of a binary quadratic form with rational coefficients agrees with the Lagrange constant of the corresponding quadratic irrational. Thus three objects give the same value: the Lagrange constant of a quadratic irrational, the Markov constant of a rational binary quadratic form, and the value of $\mathcal S$ obtained from a periodic bi-infinite sequence.

This chapter follows mainly \cite{reutenauer,classical-dynamical}.

\section{Definitions and First Examples}
We begin by introducing the set called the Markov spectrum.
\begin{defi}\label{def:markov-spectrum}\index{Discriminant!of a binary quadratic form}\index{Binary quadratic form}\index{Markov constant}\index{Markov spectrum}
Let $Q$ be a real binary quadratic form. We assume that $Q$ is indefinite; namely, if $Q(x,y)=ax^2+bxy+cy^2$, then $D(Q):=b^2-4ac>0$.
We also assume that $Q(x,y)\neq 0$ for every lattice point $(x,y)\in\mathbb Z^2\setminus\{(0,0)\}$. Then
\[
\mathcal M(Q):=\frac{\sqrt{D(Q)}}{\inf_{(x,y)\in \mathbb{Z}^2 \setminus \{(0,0)\}}|Q(x,y)|}
\]
is called the \emph{Markov constant} attached to $Q$. If the infimum in the denominator is $0$, we put $\mathcal M(Q)=\infty$. The set of all Markov constants
\[
\mathcal{M} := \left\{ \mathcal M(Q)\;\middle|\; \begin{aligned}
&Q(x,y) = ax^2 + bxy + cy^2, \; a,b,c\in \mathbb R,\; D = b^2 - 4ac > 0,\\
&\text{$Q(x,y)\neq 0$ for every $(x,y)\in\mathbb Z^2\setminus\{(0,0)\}$}
\end{aligned}\right\}
\]
is called the \emph{Markov spectrum}. With this convention, $\mathcal M\subset\mathbb R\cup\{\infty\}$, and its finite part is $\mathcal M\cap\mathbb R$.
\end{defi}
In what follows, all quadratic forms are real binary quadratic forms, so we will simply call them quadratic forms. Let us examine more carefully the conditions under which the Markov constant of $Q(x,y)=ax^2+bxy+cy^2$ is defined in the sense used in this text. The definition requires both $D(Q)>0$ and the absence of nonzero lattice zeros. The infimum may nevertheless be zero, in which case $\mathcal M(Q)=\infty$.\footnote{In this text, the assertion that $\mathcal M(Q)=\infty$ is distinguished from the assertion that $\mathcal M(Q)$ is not well defined.} We first record the following proposition.
\begin{prop}\label{prop:characterization-Q=0}
Let $Q(x,y)=ax^2+bxy+cy^2$, where $a,b,c\in \mathbb R$, $D>0$, and $a\neq0$. Then the following two conditions are equivalent:
\begin{enumerate}
\item there is no lattice point $(\alpha,\beta)\in\mathbb Z^2\setminus\{(0,0)\}$ such that $Q(\alpha,\beta)=0$;
\item the polynomial $Q(x,1)$ has two distinct irrational roots.
\end{enumerate}
\end{prop}
\begin{proof}
Since $D(Q)>0$, the polynomial $Q(x,1)$ has no multiple root. Suppose that a nonzero lattice point $(\alpha,\beta)$ satisfies $Q(\alpha,\beta)=0$. If $\beta=0$, then $a\alpha^2=0$, contradicting $a\neq0$ and $\alpha\neq0$. Hence $\beta\neq0$. Dividing $Q(\alpha,\beta)=0$ by $\beta^2$, we obtain
\[
a\left(\frac{\alpha}{\beta}\right)^2
+b\left(\frac{\alpha}{\beta}\right)
+c=0.
\]
Thus $x=\alpha/\beta\in\mathbb Q$ is a root of $Q(x,1)=0$. Therefore at least one root of $Q(x,1)$ is rational.

Conversely, suppose that one of the two roots of $Q(x,1)$ is rational. Then $Q(x,1)=0$ has a rational solution $x=\alpha/\beta$, written in lowest terms. Multiplying by $\beta^2$, we obtain $Q(\alpha,\beta)=0$.
\end{proof}
We will also use the following elementary observation.
\begin{prop}\label{prop:necessary-condition}
If $\mathcal M(Q)$ is well defined, then $a\neq 0$ and $c\neq 0$.
\end{prop}
\begin{proof}
If $a=0$, then $Q(1,0)=0$, so $\mathcal M(Q)$ is not defined. If $c=0$, then $Q(0,1)=0$, and again $\mathcal M(Q)$ is not defined. Hence both $a$ and $c$ must be nonzero.
\end{proof}
The preceding two propositions give the following corollary.
\begin{coro}\label{cor:characterization-Q=0}
The constant $\mathcal M(Q)$ is well defined if and only if $Q(x,1)$ is a quadratic polynomial with two distinct irrational roots.
\end{coro}
\begin{proof}
If $\mathcal M(Q)$ is well defined, then $a\neq0$ by Proposition~\ref{prop:necessary-condition}, and $Q$ has no nonzero lattice zero. Hence Proposition~\ref{prop:characterization-Q=0} shows that $Q(x,1)$ has two distinct irrational roots. Conversely, if $Q(x,1)$ is a quadratic polynomial with two distinct irrational roots, then $a\neq0$ and $D(Q)>0$. Proposition~\ref{prop:characterization-Q=0} shows that $Q$ has no nonzero lattice zero. Thus $\mathcal M(Q)$ is well defined.
\end{proof}
\begin{exam}\label{ex:concrete-example2}
As in the case of the Lagrange spectrum, let us look at concrete examples.
\begin{itemize}
\item[(1)] Let $Q(x,y)=x^2-xy-y^2$. Then $D(Q)=5$, and the roots of $Q(x,1)$ are $(1\pm \sqrt{5})/2$. Hence $Q(x,y)\neq0$ for every nonzero lattice point. Since $Q$ has integer coefficients, $|Q(\alpha,\beta)|\in\mathbb Z_{\geq 1}$ for every lattice point $(\alpha,\beta)$ with $Q(\alpha,\beta)\neq0$. Moreover $Q(1,0)=1$, so
\[
\inf_{(x,y)\in \mathbb{Z}^2 \setminus \{(0,0)\}}|Q(x,y)|=1.
\]
Therefore
\[
\mathcal M(Q)=\frac{\sqrt{5}}{1}=\sqrt{5}.
\]
\item[(2)] Let $Q(x,y)=x^2-2xy-y^2$. Then $D(Q)=8$, and the roots of $Q(x,1)$ are $1\pm\sqrt{2}$. Hence $Q(x,y)\neq0$ for every nonzero lattice point. Again $Q$ has integer coefficients, so $|Q(\alpha,\beta)|\in\mathbb Z_{\geq 1}$ for every lattice point $(\alpha,\beta)$ with $Q(\alpha,\beta)\neq0$. Since $Q(1,0)=1$, we have
\[
\inf_{(x,y)\in \mathbb{Z}^2 \setminus \{(0,0)\}}|Q(x,y)|=1.
\]
Thus
\[
\mathcal M(Q)=\frac{\sqrt{8}}{1}=2\sqrt{2}.
\]
\end{itemize}
\end{exam}

\section{Unimodular Group Orbits of Binary Quadratic Forms}
In this section, as preparation for computing Markov constants, we decompose quadratic forms into orbits under the unimodular group. We will also see that in each orbit one may choose a representative with good properties, called a canonical reduced quadratic form. In Chapter~\ref{chap:continued-fraction}, Section~\ref{sec:quadratic-irrational}, we carried out an analogous discussion for quadratic irrationals; the present discussion can be viewed as the counterpart for quadratic forms.

First we introduce the unimodular group action on the set of quadratic forms for which $\mathcal M(Q)$ is well defined. Put
\[
\mathcal Q
:= \left\{Q\colon\mathbb R^2\to \mathbb R \;\middle|\; \begin{aligned}
&Q(x,y) = ax^2 + bxy + cy^2, \; a,b,c\in \mathbb R,\; D(Q) = b^2 - 4ac > 0,\\
&\text{$Q(x,y)\neq 0$ for every $(x,y)\in\mathbb Z^2\setminus\{(0,0)\}$}
\end{aligned}\right\}.
\]
By Corollary~\ref{cor:characterization-Q=0}, $Q(x,1)$ always has two irrational roots for $Q\in\mathcal Q$. For \(A=\begin{bsmallmatrix}p&q\\r&s\end{bsmallmatrix}\in GL(2,\mathbb Z)\) and $Q(x,y)=ax^2+bxy+cy^2\in\mathcal Q$, define $QA$ by
\begin{equation}\label{eq:action-form}
QA(x,y):=Q(px+qy,rx+sy).
\end{equation}
We have the following.
\begin{theorem}\label{thm:action-form}\index{Unimodular group}\index{Discriminant!of a binary quadratic form}\index{Binary quadratic form}
For $A\in GL(2,\mathbb Z)$ and $Q\in\mathcal Q$, one has $QA\in\mathcal Q$. Moreover, this operation gives a right action $\mathcal Q\curvearrowleft GL(2,\mathbb Z)$.
\end{theorem}
\begin{proof}
We first show that $QA\in\mathcal Q$ for $Q\in\mathcal Q$. We have
\begin{align*}
QA(x,y)&=Q(px+qy,rx+sy)\\
&=(ap^2+bpr+cr^2)x^2
+(2apq+b(ps+qr)+2crs)xy+(aq^2+bqs+cs^2)y^2.
\end{align*}
Thus
\[
D(QA)
=(2apq+b(ps+qr)+2crs)^2
-4(ap^2+bpr+cr^2)(aq^2+bqs+cs^2)
=(b^2-4ac)(ps-qr)^2.
\]
Since $ps-qr=\pm1$, we obtain $D(QA)=D(Q)$, and in particular $D(QA)>0$. Next regard $A$ as the linear transformation of $\mathbb R^2$ given by
\[
A\begin{bmatrix}
    x\\y
\end{bmatrix}=\begin{bmatrix}
    p&q\\r&s
\end{bmatrix}\begin{bmatrix}
    x\\y
\end{bmatrix}=\begin{bmatrix}
    px+qy\\rx+sy
\end{bmatrix}.
\]
Since $A\in GL(2,\mathbb Z)$, this restricts to a bijection
\[
\mathbb Z^2\setminus\{(0,0)\}\longrightarrow \mathbb Z^2\setminus\{(0,0)\}.
\]
Hence, because $Q$ has no nonzero lattice zero, the same is true of $QA$. Therefore $QA\in\mathcal Q$. Finally, since $E_2$ denotes the $2\times2$ identity matrix, we have
$QE_2=Q$, and $(QB)A=Q(BA)$ follows immediately from viewing
$A$ and $B$ as linear transformations of $\mathbb R^2$.
\end{proof}
We now introduce unimodular equivalence on $\mathcal Q$.
\begin{defi}\label{def:unimodular-equivalence-forms}\index{Unimodular orbit}\index{Unimodular equivalence!of quadratic forms}
Let $Q,R\in\mathcal Q$. If there exists $A\in GL(2,\mathbb Z)$ such that $R=QA$, then $Q$ and $R$ are called \emph{unimodularly equivalent}; in what follows we simply say \emph{equivalent}. We write $Q\sim R$. The equivalence class
\[
O_Q=\{R\mid R\sim Q\}
\]
is called the \emph{unimodular orbit} of $Q$, or simply the \emph{orbit} of $Q$.
\end{defi}
Let $\mathcal Q(d)$ denote the set of all elements of $\mathcal Q$ with $D(Q)=d$. In the proof of Theorem~\ref{thm:action-form}, we saw that the unimodular group action preserves the discriminant. Hence we have the following corollary.
\begin{coro}
The action \eqref{eq:action-form} of the unimodular group on $\mathcal Q$ restricts to a right action on $\mathcal Q(d)$.
\end{coro}
In Corollary~\ref{cor:reduced-equiv} of Chapter~\ref{chap:continued-fraction}, Section~\ref{sec:quadratic-irrational}, we saw that a quadratic irrational $\alpha$ is unimodularly equivalent to a reduced quadratic irrational $\beta$. A similar statement holds for elements of $\mathcal Q$. To state it, we first define reduced quadratic forms.
\begin{defi}\label{def:canonical-reduced-form}\index{Reduced quadratic form}\index{Canonical reduced quadratic form}
Let $Q(x,y)=ax^2+bxy+cy^2$ be an indefinite quadratic form. Suppose that $Q(x,1)$ has two distinct roots $\alpha,\beta$ satisfying
\[
|\alpha|>1,\qquad |\beta|<1,\qquad \alpha\beta<0.
\]
Then $Q$ is called a \emph{reduced quadratic form}. If, in addition, the root with absolute value greater than $1$ satisfies $\alpha>1$, then $Q$ is called a \emph{canonical reduced quadratic form}.\footnote{This terminology is not universal; in many texts such a form is simply called a reduced quadratic form.}
\end{defi}
It may seem asymmetric to define reducedness by looking at the roots of $Q(x,1)$, since this appears to distinguish $x$ and $y$. In fact no symmetry is lost. The roots of $Q(x,1)$ and those of $Q(1,y)$ are reciprocal to each other, as long as the roots are nonzero. Therefore one obtains an equivalent definition by imposing the corresponding condition on the roots of $Q(1,y)$.

When the coefficients $a,b,c$ of $Q(x,y)=ax^2+bxy+cy^2$ are rational, if one root of $Q(x,1)$ is irrational, then the other root is also irrational, and the two roots are quadratic conjugates. If $Q$ is a canonical reduced quadratic form, then its roots satisfy $\alpha>1$ and $-1<\beta<0$, so $\alpha$ is a reduced quadratic irrational. Conversely, if $\alpha$ is a reduced quadratic irrational, then any rational-coefficient quadratic form having $\alpha$ as a root is a canonical reduced quadratic form. From this viewpoint, canonical reduced quadratic forms may be regarded as a generalization of reduced quadratic irrationals to the setting of arbitrary irrational roots.

Let $\mathcal R$ be the set of all canonical reduced quadratic forms in $\mathcal Q$, and let $\mathcal R(d)$ be the subset consisting of those with discriminant $d$. We prove the following theorem.
\begin{theorem}\label{thm:equivalent-reduced}\index{Unimodular orbit}\index{Canonical reduced quadratic form}\index{Unimodular equivalence!of quadratic forms}
For every $Q\in\mathcal Q(d)$, there exists $R\in\mathcal R(d)$ such that $Q\sim R$. Hence every orbit in $\mathcal Q(d)$ has a representative in $\mathcal R(d)$. In particular, every orbit in $\mathcal Q$ has a representative in $\mathcal R$.
\end{theorem}
We first prove a lemma.
\begin{lemm}\label{lem:action-quadratic-irrational}
Let $A=\begin{bsmallmatrix}p&q\\r&s\end{bsmallmatrix}\in GL(2,\mathbb Z)$ and let $Q\in\mathcal Q$. If the roots of $Q(x,1)$ are $\alpha$ and $\beta$, then the roots of $QA^{-1}(x,1)$ are
\[
\frac{p\alpha+q}{r\alpha+s},\qquad \frac{p\beta+q}{r\beta+s}.
\]
Similarly, if the roots of $Q(1,y)$ are $\alpha$ and $\beta$, then the roots of $QA^{-1}(1,y)$ are
\[
\frac{s\alpha+r}{q\alpha+p},\qquad \frac{s\beta+r}{q\beta+p}.
\]
\end{lemm}
\begin{proof}
We prove the first assertion. Since the roots $\alpha,\beta$ are irrational and $p,q,r,s$ are integers, the denominators $r\alpha+s$ and $r\beta+s$ do not vanish. Since $\det A=\pm1$, the inverse of $A$ is either $\begin{bsmallmatrix}s&-q\\-r&p\end{bsmallmatrix}$ or its negative. Because a quadratic form is homogeneous of degree two, this overall sign does not affect the value of the form. Hence, for the purpose of finding the roots, we may compute
\begin{align*}
QA^{-1}(x,1)
&=Q(sx-q,-rx+p)=(-rx+p)^2 Q\left(\frac{sx-q}{-rx+p},1\right).
\end{align*}
Since $Q(\alpha,1)=0$, the equality
\[
\alpha=\frac{sx-q}{-rx+p}
\]
implies that $QA^{-1}(x,1)=0$. Solving this equality for $x$ gives
\[
x=\frac{p\alpha+q}{r\alpha+s}.
\]
Thus this is a root of $QA^{-1}(x,1)$. The argument for $\beta$ is identical. The second assertion is proved in the same way, using the roots of $Q(1,y)$ instead of those of $Q(x,1)$.
\end{proof}
\begin{proof}[Proof of Theorem~\ref{thm:equivalent-reduced}]
Let $\alpha,\beta$ be the roots of $Q(x,1)$, ordered so that $\alpha>\beta$. By Proposition~\ref{prop:characterization-Q=0}, both roots are irrational. If the roots already satisfy $\alpha>0>\beta$, we do nothing at this stage. If $0>\alpha>\beta$, choose a sufficiently large integer $h$ and replace $Q$ by \(Q\begin{bsmallmatrix}1&h\\0&1\end{bsmallmatrix}^{-1}\). By Lemma~\ref{lem:action-quadratic-irrational}, the new roots are $\alpha+h$ and $\beta+h$, so we are reduced to the case where both roots are positive.

It remains to handle the case $\alpha>\beta>0$. During the process of deleting common initial partial quotients, the order of the two roots may be reversed. Whenever this happens, we rename the larger root $\alpha$ and the smaller root $\beta$. Write the infinite regular continued-fraction expansions as
\[
\alpha=[a_0;a_1,\dots],\qquad \beta=[b_0;b_1,\dots].
\]
Let $m$ be the smallest index such that $a_m\neq b_m$. We construct a quadratic form $Q'$ equivalent to $Q$ whose roots $\alpha',\beta'$ satisfy $\alpha'>0>\beta'$. If $m=0$, then after interchanging the names of the roots if necessary, we may assume $a_0>b_0$. Then $\alpha-a_0>0>\beta-a_0$. By Lemma~\ref{lem:action-quadratic-irrational}, a form whose roots are $\alpha-a_0$ and $\beta-a_0$ is obtained by taking \(Q'=Q\begin{bsmallmatrix}1&-a_0\\0&1\end{bsmallmatrix}^{-1}\). This gives the desired $Q'$. If $m\neq0$, put \(Q_1:=Q\begin{bsmallmatrix}0&1\\1&-a_0\end{bsmallmatrix}^{-1}\). Then the roots of $Q_1(x,1)$ are
\[
\alpha_1=[a_1;a_2,\dots],\qquad \beta_1=[b_1;b_2,\dots].
\]
Thus replacing $Q$ by the equivalent form $Q_1$ decreases the value of $m$ by at least one. Repeating this operation until $m=0$, and then applying the argument above, we obtain a form $Q'$ whose two roots have opposite signs.

Therefore, after replacing $Q$ by an equivalent form if necessary, we may assume that the roots of $Q(x,1)$ are $\alpha,\beta$ with $\alpha>0>\beta$. If $|\alpha|>1$ and $|\beta|<1$, then $R=Q$ is already the desired form. If $|\alpha|>1$ and $|\beta|>1$, choose an integer $h$ such that $-1<\beta+h<0$ and take \(R=Q\begin{bsmallmatrix}1&h\\0&1\end{bsmallmatrix}^{-1}\). If $|\alpha|<1$ and $|\beta|<1$, take \(\widetilde Q:=Q\begin{bsmallmatrix}0&1\\1&0\end{bsmallmatrix}^{-1}\). The roots of $\widetilde Q(x,1)$ are $1/\alpha$ and $1/\beta$, so this case is reduced to one of the preceding cases. Finally, if $|\alpha|<1$ and $|\beta|>1$, then $R=\widetilde Q$ is canonical reduced. This proves the theorem.
\end{proof}
Lemma~\ref{lem:action-quadratic-irrational} shows that the unimodular action on quadratic forms induces the unimodular action on their roots. In particular, continued-fraction reduction gives another proof of Corollary~\ref{cor:reduced-equiv}.

The eventual reduction assertion of Theorem~\ref{thm:reduced-equiv} also extends to arbitrary $Q\in\mathcal Q$. Let $\alpha,\beta$ be the distinct irrational roots of $Q(x,1)$, and let $p_n/q_n$ be the convergents of $\alpha=[a_0;a_1,\dots]$. Put
\[
M_n=\begin{bmatrix}p_n&p_{n-1}\\q_n&q_{n-1}\end{bmatrix}.
\]
The roots of $QM_n$ are $\alpha_{n+1}=[a_{n+1};a_{n+2},\dots]$ and
\[
\beta_{n+1}=-\frac{q_{n-1}}{q_n}
\frac{\beta-p_{n-1}/q_{n-1}}{\beta-p_n/q_n}\qquad(n\geq1).
\]
Since $p_n/q_n\to\alpha\neq\beta$, the last ratio is positive for all sufficiently large $n$, so $\beta_{n+1}<0$. It follows that
\[
\beta_{n+2}=\frac1{\beta_{n+1}-a_{n+1}}\in(-1,0),\qquad \alpha_{n+2}>1.
\]
Thus all sufficiently late transformed forms are canonical reduced. Each transformation is unimodular and preserves the discriminant. If $\alpha$ is quadratic and $Q$ is the homogenization of its primitive integral minimal polynomial, then $\beta_{n+2}$ is the conjugate of $\alpha_{n+2}$, yielding Corollary~\ref{cor:reduced-equiv}.

The property specific to quadratic irrationals is the finiteness of reduced roots of a fixed discriminant when primitive integral minimal polynomials are used (Lemma~\ref{lem:finiteness-R2d}). For arbitrary real coefficients, canonical reduced forms of a fixed discriminant need not form a finite set, so reduction alone does not imply periodicity. For example, for fixed $d>0$, the forms
\[
Q_\theta(x,y)=\frac{\sqrt d}{\theta+1/\sqrt2}(x-\theta y)(x+y/\sqrt2),
\qquad \theta>1,\quad\theta\notin\mathbb Q,
\]
are infinitely many distinct elements of $\mathcal R(d)$.

\section{A Bi-infinite Sequence Formula for the Markov Constant}
We now explain how to compute the Markov constant by means of a bi-infinite sequence, in parallel with the Lagrange constant. The next corollary is a refinement of Theorem~\ref{thm:action-form} and follows directly from it.
\begin{coro}\label{cor:MQ=MR}\index{Unimodular equivalence!of quadratic forms}
If $Q,R\in\mathcal Q$ and $Q\sim R$, then $\mathcal M(Q)=\mathcal M(R)$.
\end{coro}
\begin{proof}
By definition,
\[
\mathcal M(Q)=\frac{\sqrt{D(Q)}}{\inf_{(x,y)\in \mathbb{Z}^2 \setminus \{(0,0)\}}|Q(x,y)|}.
\]
It suffices to show that $D(Q)=D(R)$ and that the two infima in the denominators are equal. The equality $D(Q)=D(R)$ was proved in the proof of Theorem~\ref{thm:action-form}. Moreover, an element $A\in GL(2,\mathbb Z)$ acts on $\mathbb R^2$ as a bijection that sends lattice points to lattice points and the origin to the origin. Hence
\[
\inf_{(x,y)\in \mathbb{Z}^2 \setminus \{(0,0)\}}|R(x,y)|
=\inf_{(x,y)\in \mathbb{Z}^2 \setminus \{(0,0)\}}|Q(px+qy,rx+sy)|
=\inf_{(x,y)\in \mathbb{Z}^2 \setminus \{(0,0)\}}|Q(x,y)|.
\]
The claim follows.
\end{proof}

\begin{coro}\label{cor:M=MQR}
The following equality holds:
\[
\mathcal M=\left\{\mathcal M(Q)\  \middle|\  Q\in \mathcal R\right\}.
\]
\end{coro}
\begin{proof}
This follows immediately from Theorem~\ref{thm:equivalent-reduced} and Corollary~\ref{cor:MQ=MR}.
\end{proof}

The next theorem is the main result of this section.
\begin{theorem}\label{thm:markov-infinity-sequence}\index{Bi-infinite sequence}\index{Canonical reduced quadratic form}\index{Markov constant}
Let $Q\in\mathcal R$, and let the two irrational roots of $Q(x,1)$ be $
\theta$ and $-\frac{1}{\eta}$, 
where $\theta,\eta>1$. From the continued-fraction expansions
\[
\theta=[a_0;a_1,\dots],\qquad \eta=[a_{-1};a_{-2},\dots],
\]
form the bi-infinite sequence
\[
\mathbf a=(\dots,a_{-2},a_{-1},a_0,a_1,a_2,\dots).
\]
Then
\[
\mathcal M(Q)=\mathcal S(\mathbf a).
\]
\end{theorem}
Before proving the theorem, let us spell out what it means. Corollary~\ref{cor:perron-formula} represents the Lagrange constant of any irrational number $\alpha$ as $\limsup_{n\to+\infty}\ell_n(\mathbf a)$ for a suitable bi-infinite sequence $\mathbf a$. For a quadratic irrational, repeating its continued-fraction period gives an explicit sequence satisfying $\mathcal L(\alpha)=\mathcal S(\mathbf a)$. For a general irrational number with bounded partial quotients, the proof of Theorem~\ref{thm:LsubsetS} obtains this supremum representation by taking accumulation points. In the present theorem, the two irrational roots of $Q(x,1)$ directly supply the two continued-fraction tails whenever $Q\in\mathcal R$. Conversely, by choosing arbitrary irrational numbers $\theta,
\eta>1$ and setting
\[
Q=(x-\theta y)\left(x+\frac{1}{\eta}y\right),
\]
one obtains an element $Q\in\mathcal R$ corresponding to any given bi-infinite sequence $\mathbf a$ of positive integers. We therefore get the following corollary.
\begin{coro}\label{cor:M=S}\index{Lagrange spectrum!inclusion in the Markov spectrum}\index{Markov spectrum}
One has $\mathcal M=\mathcal S$. In particular, Theorem~\ref{thm:LsubsetS} implies $\mathcal L\subset \mathcal M$.
\end{coro}
We now prepare the proof of Theorem~\ref{thm:markov-infinity-sequence}. For $\delta\in\mathbb Z$, put \(F_\delta:=\begin{bsmallmatrix}\delta&1\\1&0\end{bsmallmatrix}\).
\begin{lemm}\label{lem:uniqueness-neighbor}
Let $Q\in\mathcal R$, and suppose that the roots of $Q(x,1)$ are $\theta$ and $-1/\eta$, where $\theta,\eta>1$. Then $QF_\delta^{-1}$ and $QF_\varepsilon$ belong to $\mathcal R$ if and only if, respectively,
\[
\delta=\delta_0:=\left\lfloor\eta\right\rfloor,
\qquad
\varepsilon=\varepsilon_0:=\left\lfloor\theta\right\rfloor.
\]
\end{lemm}
\index{Neighbor of a quadratic form}We call $QF_{\delta_0}^{-1}$ and $QF_{\varepsilon_0}$ the \emph{left neighbor} and the \emph{right neighbor} of $Q$, respectively.
\begin{proof}
By Lemma~\ref{lem:action-quadratic-irrational}, the roots of $QF_\delta^{-1}(x,1)$ are
\[
\delta+\frac1\theta,\qquad\delta-\eta.
\]
If $\delta\leq0$, the first is less than $1$ and the second less than $-1$, so the form is not canonical reduced. If $\delta\geq1$, the first root is greater than $1$, and the necessary and sufficient condition is
\[
-1<\delta-\eta<0,\qquad\text{equivalently}\qquad\eta-1<\delta<\eta.
\]
Since $\eta$ is irrational, the unique integer satisfying this is $\delta=\lfloor\eta\rfloor$.

Next, $F_\varepsilon^{-1}=\begin{bsmallmatrix}0&1\\1&-\varepsilon\end{bsmallmatrix}$ and $QF_\varepsilon=Q(F_\varepsilon^{-1})^{-1}$, so the roots of $QF_\varepsilon(x,1)$ are
\[
\frac1{\theta-\varepsilon},\qquad
\frac1{-1/\eta-\varepsilon}=-\frac1{\varepsilon+1/\eta}.
\]
If $\varepsilon=0$, these roots lie in $(0,1)$ and $(-\infty,-1)$, respectively. If $\varepsilon\leq-1$, the first lies in $(0,1)$ and the second is positive. Neither case gives a canonical reduced form.

If $\varepsilon\geq1$, the second root lies in $(-1,0)$. The form is then canonical reduced if and only if the first root is greater than $1$, namely if and only if
\[
0<\theta-\varepsilon<1.
\]
Since $\theta$ is irrational, the unique integer satisfying this is $\varepsilon=\lfloor\theta\rfloor$. Both assertions follow.
\end{proof}
There is one point about the relation between right and left neighbors that must be checked. Namely, the right neighbor of the left neighbor of $Q$ should return to $Q$, and conversely. This sounds obvious, since the matrices used to move to the right and to the left are inverse to each other. However, the integers $\delta$ and $\varepsilon$ defining the left and right neighbors are determined separately for each quadratic form. Thus we must check that the value of $\delta$ for $Q$ agrees with the value of $\varepsilon$ for the left neighbor of $Q$.
\begin{prop}\label{prop:neighbor-well-defined}\index{Neighbor of a quadratic form}
Let $Q_0\in\mathcal R$. If $Q_1$ is the right neighbor of $Q_0$, then the left neighbor of $Q_1$ is $Q_0$. If $Q_{-1}$ is the left neighbor of $Q_0$, then the right neighbor of $Q_{-1}$ is $Q_0$.
\end{prop}
\begin{proof}
Let the two roots of $Q_0(x,1)$ be $\theta_0$ and $-1/\eta_0$, where $\theta_0,\eta_0>1$, and write
\[
\theta_0=[a_0;a_1,\dots],\qquad \eta_0=[a_{-1};a_{-2},\dots].
\]
For $i=\pm1$, let the two roots of $Q_i(x,1)$ be $\theta_i$ and $-1/\eta_i$, where $\theta_i,\eta_i>1$. We prove
\[
Q_0=Q_0F_{\delta_0}^{-1}F_{\varepsilon_{-1}},
\qquad
Q_0=Q_0F_{\varepsilon_0}F_{\delta_1}^{-1},
\quad\text{where}\quad
\delta_i:=\left\lfloor\eta_i\right\rfloor,
\qquad
\varepsilon_i:=\left\lfloor\theta_i\right\rfloor.
\]
By definition, $\delta_0=a_{-1}$ and $\varepsilon_0=a_0$. If $Q_{-1}=Q_0F_{\delta_0}^{-1}$, then Lemma~\ref{lem:uniqueness-neighbor} gives $\delta_0=\lfloor\eta_0\rfloor$, and hence
\[
\theta_{-1}=\lfloor\eta_0\rfloor+\frac{1}{\theta_0}
=[a_{-1};a_0,\dots],
\qquad
-\frac{1}{\eta_{-1}}=\lfloor\eta_0\rfloor-\eta_0
=-[0;a_{-2},a_{-3},\dots].
\]
Thus $\varepsilon_{-1}=a_{-1}=\delta_0$, and therefore
\[
Q_0=Q_0F_{\delta_0}^{-1}F_{\varepsilon_{-1}}.
\]
Next, if $Q_1=Q_0F_{\varepsilon_0}$, then Lemma~\ref{lem:uniqueness-neighbor} gives $\varepsilon_0=\lfloor\theta_0\rfloor$, and so
\[
\theta_1=\frac{1}{\theta_0-\left\lfloor\theta_0\right\rfloor}
=\frac{1}{[0;a_1,a_2,\dots]}=[a_1;a_2,\dots],
\qquad
-\frac{1}{\eta_1}=\frac{1}{-\frac{1}{\eta_0}-\lfloor\theta_0\rfloor}
=-[0;a_0,a_{-1},\dots].
\]
Therefore $\eta_1=[a_0;a_{-1},\dots]$, and hence $\delta_1=a_0=\varepsilon_0$. It follows that
\[
Q_1F_{\delta_1}^{-1}=Q_0F_{\varepsilon_0}F_{\delta_1}^{-1}=Q_0.
\]
\end{proof}
The preceding proposition makes the following definition well defined.
\begin{defi}
Let $Q_0\in\mathcal R$. Define a bi-infinite sequence of quadratic forms
\[
(\dots,Q_{-2},Q_{-1},Q_0,Q_1,Q_2,\dots)
\]
inductively by requiring that $Q_{i+1}$ be the unique right neighbor of $Q_i$ and that $Q_{i-1}$ be the unique left neighbor of $Q_i$. We call this sequence the \emph{chain} of $Q_0$.
\end{defi}
The next corollary describes the roots of the forms appearing in the chain of $Q_0$.
\begin{coro}\label{cor:root-description}
Let $Q_0\in\mathcal R$, and consider the chain
\[
(\dots,Q_{-2},Q_{-1},Q_0,Q_1,Q_2,\dots)
\]
of $Q_0$. Suppose that the two roots of $Q_0(x,1)$ are $\theta_0$ and $-1/\eta_0$, where $\theta_0,\eta_0>1$, and write
\[
\theta_0=[a_0;a_1,\dots],\qquad \eta_0=[a_{-1};a_{-2},\dots].
\]
For $i\in\mathbb Z$, let the two roots of $Q_i(x,1)$ be $\theta_i$ and $-1/\eta_i$, where $\theta_i,\eta_i>1$. Then
\[
\theta_i=[a_i;a_{i+1},\dots], \qquad \eta_i=[a_{i-1};a_{i-2},\dots].
\]
\end{coro}
\begin{proof}
The case $i=\pm1$ is contained in the proof of Proposition~\ref{prop:neighbor-well-defined}. The general case follows by induction using the same argument.
\end{proof}

In the next proposition, the defining expression of the Markov constant appears naturally inside continued fraction theory. This makes the connection with continued fractions, and the strategy of the proof, more transparent.
\begin{prop}\label{prop:D/Q}
Let $Q_0\in\mathcal R$, let
\[
(\dots,Q_{-2},Q_{-1},Q_0,Q_1,Q_2,\dots)
\]
be the chain of $Q_0$, and let $\theta_i,-1/\eta_i$ be the roots of $Q_i(x,1)$, with $\theta_i,\eta_i>1$. If
\[
\theta_0=[a_0;a_1,\dots],\qquad \eta_0=[a_{-1};a_{-2},\dots],
\]
then
\[
[a_i;a_{i+1},\dots]+[0;a_{i-1},a_{i-2},\dots]
=\theta_i+\frac{1}{\eta_i}
=\frac{\sqrt{D(Q_0)}}{|Q_i(1,0)|}.
\]
\end{prop}
\begin{proof}
The first equality follows immediately from Corollary~\ref{cor:root-description}. We prove the second one. Write
\[
Q_i(x,y)=ax^2+bxy+cy^2.
\]
The numbers $\theta_i$ and $-1/\eta_i$ are the roots of $Q_i(x,1)=ax^2+bx+c$, and $\theta_i>-1/\eta_i$. Thus, if $a>0$, then
\[
\theta_i=\frac{-b+\sqrt{D(Q_i)}}{2a},\qquad
-\frac{1}{\eta_i}=\frac{-b-\sqrt{D(Q_i)}}{2a},
\]
whereas if $a<0$, then
\[
\theta_i=\frac{-b-\sqrt{D(Q_i)}}{2a},\qquad
-\frac{1}{\eta_i}=\frac{-b+\sqrt{D(Q_i)}}{2a}.
\]
In the first case,
\[
\theta_i+\frac{1}{\eta_i}=\frac{\sqrt{D(Q_i)}}{a}
=\frac{\sqrt{D(Q_i)}}{|a|},
\]
and in the second case,
\[
\theta_i+\frac{1}{\eta_i}=-\frac{\sqrt{D(Q_i)}}{a}
=\frac{\sqrt{D(Q_i)}}{|a|}.
\]
Finally, $a=Q_i(1,0)$, and $Q_i$ is equivalent to $Q_0$, so $D(Q_i)=D(Q_0)$. The result follows.
\end{proof}
The next proposition is the final preparation for Theorem~\ref{thm:markov-infinity-sequence}. Its proof is the most delicate part of this section.
\begin{prop}\label{prop:inf=inf}
Let $Q_0\in\mathcal R$, and let
\[
(\dots,Q_{-2},Q_{-1},Q_0,Q_1,Q_2,\dots)
\]
be the chain of $Q_0$. Then
\begin{equation}
\inf_{h\in\mathbb Z}{|Q_h(1,0)|}
=
\inf_{(x,y)\in \mathbb{Z}^2 \setminus \{(0,0)\}}|Q_0(x,y)|.
\end{equation}
\end{prop}
\begin{proof}
For every $h\in\mathbb Z$, the forms $Q_h$ and $Q_0$ are equivalent. Hence
\[
\inf_{(x,y)\in \mathbb{Z}^2 \setminus \{(0,0)\}}|Q_h(x,y)|
=
\inf_{(x,y)\in \mathbb{Z}^2 \setminus \{(0,0)\}}|Q_0(x,y)|.
\]
By the definition of an infimum,
\[
|Q_h(1,0)|\geq
\inf_{(x,y)\in \mathbb{Z}^2 \setminus \{(0,0)\}}|Q_h(x,y)|.
\]
Therefore
\[
\inf_{h\in \mathbb{Z}}|Q_h(1,0)|\geq
\inf_{(x,y)\in \mathbb{Z}^2 \setminus \{(0,0)\}}|Q_0(x,y)|.
\]
It remains to prove
\begin{equation}\label{eq:inf<inf}
\inf_{h\in \mathbb{Z}}|Q_h(1,0)|\leq
\inf_{(x,y)\in \mathbb{Z}^2 \setminus \{(0,0)\}}|Q_0(x,y)|.
\end{equation}

First suppose that the bi-infinite sequence
\[
(\dots,a_{-2},a_{-1},a_0,a_1,a_2,\dots)
\]
determined by the two irrational roots $\theta_0$ and $-1/\eta_0$ of $Q_0(x,1)$ satisfies $a_i=1$ for all $i\in\mathbb Z$. Then, by Example~\ref{ex:concrete-example} (1), we have $\theta_0=\eta_0=(1+\sqrt{5})/2$. Every quadratic form having $\theta_0$ and $-1/\eta_0$ as roots is of the form
\[
Q(x,y)=\lambda(x-\theta_0 y)\left(x+\frac{1}{\eta_0}y\right)
\]
for some nonzero real number $\lambda$.
Replacing $Q$ by $\lambda^{-1}Q$ multiplies both sides of the desired equality of infima by $|\lambda|^{-1}$, and hence does not affect whether the equality holds. Thus we may assume $\lambda=1$ and $Q_0(x,y)=x^2-xy-y^2$. By Example~\ref{ex:concrete-example2} (1),
\[
\inf_{h\in \mathbb{Z}}|Q_h(1,0)|
=|Q_0(1,0)|=1=
\inf_{(x,y)\in \mathbb{Z}^2 \setminus \{(0,0)\}}|Q_0(x,y)|.
\]
In particular, \eqref{eq:inf<inf} holds in this case.

Now suppose that $a_i\geq2$ for some $i\in\mathbb Z$. By moving from $Q_0$ to a suitable left or right neighbor and then renaming it $Q_0$, we may assume that $a_0\geq2$. Then $\theta_0>2$. It suffices to show that for every $(x,y)\in\mathbb Z^2\setminus\{(0,0)\}$ there exists an index $i$ such that
\[
|Q_0(x,y)|\geq |Q_i(1,0)|.
\]
If $|Q_0(x,y)|\geq |Q_0(1,0)|$, then we take $i=0$. Hence assume $|Q_0(x,y)|< |Q_0(1,0)|$. Under this assumption, $y=0$ is impossible, since
\[
|Q_0(x,0)|=x^2|Q_0(1,0)|\geq |Q_0(1,0)|.
\]
If $x=0$, then, because $Q_0$ is the right neighbor of $Q_{-1}$,
\begin{align*}
|Q_0(0,y)|
&=|Q_{-1}F_{\varepsilon_{-1}}(0,y)|\\
&=\left|Q_{-1}\begin{bmatrix}
  a_{-1}&1\\1&0
\end{bmatrix}(0,y)\right|\\
&=\left|Q_{-1}(y,0)\right|
=y^2\left|Q_{-1}(1,0)\right|
\geq \left|Q_{-1}(1,0)\right|.
\end{align*}
Thus the case $x=0$ is settled. Henceforth assume $x\neq0$ and $y\neq0$. Write
\begin{equation}\label{eq:Q=a0}
Q_0(x,y)=ay^2\left(\frac{x}{y}-\theta_0\right)\left(\frac{x}{y}+\frac{1}{\eta_0}\right).
\end{equation}
First consider the case $x/y>0$. Since $Q(x,y)=Q(-x,-y)$, we may assume $x,y>0$. We claim that in fact $x/y\geq2$. Suppose, to the contrary, that $x/y<2$. Since $\theta_0>2$, we have $x/y<\theta_0$. Hence
\[
|Q_0(x,y)|
=|a|y^2\left|\frac{x}{y}-\theta_0\right|\left|\frac{x}{y}+\frac{1}{\eta_0}\right|
>|a|y^2\left(2-\frac{x}{y}\right)\left|\frac{x}{y}\right|
= |a|(2y-x)x\geq |a|.
\]
This contradicts $|Q_0(x,y)|<|Q_0(1,0)|=|a|$. Thus $x/y\geq2$. Comparing the assumption $|Q_0(x,y)|<|Q_0(1,0)|=|a|$ with \eqref{eq:Q=a0}, we obtain
\[
2y^2\left|\frac{x}{y}-\theta_0\right|
<y^2\left|\frac{x}{y}-\theta_0\right|\left(\frac{x}{y}+\frac{1}{\eta_0}\right)<1.
\]
Therefore
\[
\left|\theta_0-\frac{x}{y}\right|<\frac{1}{2y^2}.
\]
Write $x/y=p/q$ in lowest terms. Then $q\leq y$, and the above inequality gives
\[
\left|\theta_0-\frac{p}{q}\right|<\frac{1}{2y^2}\leq\frac{1}{2q^2}.
\]
By Theorem~\ref{thm:good-is-approximation}, $p/q$ is a convergent of $\theta_0$. Thus there exists $n\geq0$ such that
\[
\frac{p}{q}=[a_0;a_1,\dots,a_n]=:\frac{p_n}{q_n}.
\]
Then
\[
\begin{bmatrix}
    p_n\\q_n
\end{bmatrix}=F_{a_0}F_{a_1}\cdots F_{a_n}\begin{bmatrix}
    1\\0
\end{bmatrix}.
\]
Applying $Q_0$ to both sides gives
\[
Q_0(p_n,q_n)=Q_{n+1}(1,0).
\]
Since $p_n/q_n$ is in lowest terms, there exists $d\in\mathbb Z_{\geq1}$ such that $(x,y)=(dp_n,dq_n)$. Therefore
\[
|Q_0(x,y)|=d^2|Q_0(p_n,q_n)|\geq |Q_0(p_n,q_n)|=|Q_{n+1}(1,0)|.
\]
This proves the case $x/y>0$.

Now assume $x/y<0$. Replacing $(x,y)$ by $(-x,-y)$ if necessary, we may assume $y>0$. Since
\[
\left|\frac{x}{y}-\theta_0\right|>2,
\]
the assumption $|Q_0(x,y)|<|Q_0(1,0)|=|a|$ and \eqref{eq:Q=a0} imply
\[
2y^2\left|\frac{x}{y}+\frac{1}{\eta_0}\right|
<y^2\left|\frac{x}{y}-\theta_0\right|\left|\frac{x}{y}+\frac{1}{\eta_0}\right|<1.
\]
Hence
\[
\left|\frac{x}{y}+\frac{1}{\eta_0}\right|<\frac{1}{2y^2}.
\]
As in the argument for $\theta_0$, write the reduced form of $-x/y$ as $p'/q'$. Then $q'\leq y$, and $p'/q'$ is a convergent of $1/\eta_0$. Since $p'/q'>0$, there exists $m\geq1$ such that
\[
\frac{p'}{q'}=[0;a_{-1},a_{-2},\dots,a_{-m}].
\]
Write this convergent as $p'_m/q'_m$, and put
\[
J:=\begin{bmatrix}0&-1\\1&0\end{bmatrix}.
\]
Direct calculation gives $F_a^{-1}=-JF_aJ^{-1}$, so for $m\geq1$,
\[
F_{a_{-1}}^{-1}\cdots F_{a_{-(m+1)}}^{-1}
=(-1)^{m+1}J F_{a_{-1}}\cdots F_{a_{-(m+1)}}J^{-1}.
\]
On the other hand,
\[
F_0F_{a_{-1}}\cdots F_{a_{-m}}
=\begin{bmatrix}p'_m&p'_{m-1}\\q'_m&q'_{m-1}\end{bmatrix}.
\]
Thus the second column of $F_{a_{-1}}\cdots F_{a_{-(m+1)}}$ is $\begin{bsmallmatrix}q'_m\\p'_m\end{bsmallmatrix}$. Since $J^{-1}\begin{bsmallmatrix}1\\0\end{bsmallmatrix}=-\begin{bsmallmatrix}0\\1\end{bsmallmatrix}$, we obtain
\[
F_{a_{-1}}^{-1}\cdots F_{a_{-(m+1)}}^{-1}\begin{bmatrix}1\\0\end{bmatrix}
=(-1)^{m+2}\begin{bmatrix}-p'_m\\q'_m\end{bmatrix}.
\]
The definition of the chain and homogeneity of the quadratic form therefore give
\[
Q_0(-p'_m,q'_m)=Q_{-(m+1)}(1,0).
\]
Since $x/y=-p'_m/q'_m$, $y>0$, and the fraction is reduced, $(x,y)=(-dp'_m,dq'_m)$ for some positive integer $d$. Consequently,
\[
|Q_0(x,y)|=d^2|Q_0(-p'_m,q'_m)|\geq|Q_0(-p'_m,q'_m)|=|Q_{-(m+1)}(1,0)|.
\]
This completes the proof.
\end{proof}
We are now ready to prove Theorem~\ref{thm:markov-infinity-sequence}. With the preparations in place, the proof is short.
\begin{proof}[Proof of Theorem~\ref{thm:markov-infinity-sequence}]
For $Q=Q_0$, we must show
\[
\sup_{h\in \mathbb Z}([a_h;a_{h+1},\dots]+[0;a_{h-1},a_{h-2},\dots])
=\frac{\sqrt{D(Q)}}{\inf_{(x,y)\in \mathbb{Z}^2 \setminus \{(0,0)\}}|Q(x,y)|}.
\]
By Proposition~\ref{prop:D/Q}, it suffices to prove the equality of extended real numbers
\[
\frac{\sqrt{D(Q_0)}}{\inf_{h\in\mathbb Z}|Q_h(1,0)|}
=\frac{\sqrt{D(Q_0)}}{\inf_{(x,y)\in \mathbb{Z}^2 \setminus \{(0,0)\}}|Q_0(x,y)|}.
\]
If both denominators are $0$, both sides are read as $\infty$. In all cases, the equality follows from Proposition~\ref{prop:inf=inf}.
\end{proof}
We close this section by commenting on the inclusion relation between $\mathcal L$ and $\mathcal M$. Corollary~\ref{cor:M=S} shows that $\mathcal L\subset\mathcal M$. This inclusion is known to be proper, as was proved in various works including \cite{freiman}. Concrete values that belong to $\mathcal M$ but not to $\mathcal L$ are still actively studied; see, for example, \cite{classical-dynamical} for a detailed account of known constructions in $\mathcal M\setminus\mathcal L$.

\section{Markov Constants of Quadratic Forms with Rational Coefficients}
We finish this chapter by considering $\mathcal M(Q)$ in the case where the coefficients $a,b,c$ of $Q(x,y)=ax^2+bxy+cy^2$ are rational. From the previous section, we know how to compute the Markov constant when $Q$ is canonical reduced. Thus we would like to replace $Q$ by a canonical reduced form $Q'$ satisfying $\mathcal M(Q)=\mathcal M(Q')$. The existence of such a form is guaranteed by Theorem~\ref{thm:equivalent-reduced}, and the proof also gives an algorithm to find it. When $a,b,c$ are rational, however, the theory of quadratic irrationals tells us more directly which form $Q'$ should be used.

Let
\[
\mathcal Q_\mathbb Q:=\{Q(x,y)=ax^2+bxy+cy^2\in \mathcal Q\mid a,b,c\in \mathbb Q\}
\]
be the set of quadratic forms in $\mathcal Q$ with rational coefficients.
\begin{prop}\label{prop:equivalent-calculatable}
Let $Q(x,y)=ax^2+bxy+cy^2\in\mathcal Q_\mathbb Q$. Then $Q(x,1)$ has two quadratic irrational roots, which are conjugate to each other. Choose one of them and call it $\alpha$. Suppose that the infinite continued-fraction expansion of $\alpha$ is
\[
\alpha=[a_0;a_1,\dots,a_n,\overline{b_0,\dots,b_{k-1}}]
\]
for some $n\geq0$ and $k\geq1$.
Put
\[
\beta:=[\overline{b_0,b_1,\dots,b_{k-1}}],
\]
and let $\beta'$ be the quadratic conjugate of $\beta$. Define
\[
Q'(x,y):=(x-\beta y)(x-\beta' y).
\]
Then
\[
\mathcal M(Q)=\mathcal M(Q'),
\]
and $Q'$ is canonical reduced.
\end{prop}
\begin{proof}
By Theorem~\ref{thm:reduced-equiv}, the numbers $\alpha$ and $\beta$ are equivalent. Thus there exists $A\in GL(2,\mathbb Z)$ such that $\beta=A\alpha$. Consider the quadratic form $QA^{-1}$. By Lemma~\ref{lem:action-quadratic-irrational}, the polynomial $QA^{-1}(x,1)$ has $\beta$ as a root. Since $QA^{-1}(x,y)$ has rational coefficients, the roots of $QA^{-1}(x,1)$ are $\beta$ and $\beta'$. Therefore there exists $a'\in\mathbb Q\setminus\{0\}$ such that
\[
QA^{-1}(x,y)=a'(x-\beta y)(x-\beta'y).
\]
Since $\beta$ is reduced, we have $\beta>1$ and $-1<\beta'<0$. Hence $QA^{-1}$ is canonical reduced. By Corollary~\ref{cor:MQ=MR},
\[
\mathcal M(Q)=\mathcal M(QA^{-1}).
\]
Moreover $Q'=(1/a')QA^{-1}$. Canonical reducedness is preserved under multiplication by a nonzero scalar, and the Markov constant is unchanged by such a scalar multiple. Hence $Q'$ is canonical reduced and $\mathcal M(Q')=\mathcal M(Q)$.
\end{proof}
Proposition~\ref{prop:equivalent-calculatable} determines $\mathcal M(Q)$ from one of the two roots of $Q(x,1)$ or $Q(1,y)$. One might therefore ask whether determining the Markov constant requires information about both irrational roots. In the present case, however, the form has rational coefficients: if $\alpha$ is an irrational root of $Q(x,1)$, then its quadratic conjugate $\alpha'$ is the other root. Thus, in this setting, the Markov constant is essentially determined by a single irrational number. In this case there is a clean correspondence between the Lagrange constant and the Markov constant.
\begin{theorem}\label{thm:quadratic-markov-lagrange1770}\index{Binary quadratic form}\index{Markov constant}\index{Quadratic irrational}
Let $Q\in\mathcal Q_\mathbb Q$, and let $\alpha$ be one of the quadratic irrational roots of $Q(x,1)=0$. Then
\[
\mathcal M(Q)=\mathcal L(\alpha).
\]
\end{theorem}
\begin{proof}
Suppose that the infinite continued-fraction expansion of $\alpha$ is
\[
\alpha=[a_0;a_1,\dots,a_n,\overline{b_0,\dots,b_{k-1}}]
\]
for some $n\geq0$ and $k\geq1$.
Put $\beta=[\overline{b_0,b_1,\dots,b_{k-1}}]$ and $Q'(x,y)=(x-\beta y)(x-\beta'y)$. By Proposition~\ref{prop:equivalent-calculatable} and Theorem~\ref{thm:markov-infinity-sequence}, if
\[
\mathbf b=(\dots,b_0,b_1,\dots,b_{k-1},b_0,b_1,\dots,b_{k-1},b_0,b_1,\dots,b_{k-1},\dots)
\]
is the bi-infinite sequence obtained by repeating the period in both directions, then
\[
\mathcal M(Q)=\mathcal M(Q')=\mathcal S(\mathbf b).
\]
On the other hand, Proposition~\ref{prop:equivalent-lagrange1770} and Theorem~\ref{thm:L=S} give
\[
\mathcal L(\alpha)=\mathcal L(\beta)=\mathcal S(\mathbf b).
\]
Therefore $\mathcal M(Q)=\mathcal L(\alpha)$.
\end{proof}
Let us compute the Markov constant of a quadratic form analogous to Example~\ref{ex:1+sqrt3-lagrange1770}.
\begin{exam}
We compute the Markov constant $\mathcal M(Q)$ of the rational-coefficient quadratic form $Q(x,y)=x^2-2xy-2y^2$. The roots of $Q(x,1)=x^2-2x-2$ are $x=1\pm\sqrt{3}$. Thus we may take $\alpha=1+\sqrt{3}$ as one quadratic irrational root of $Q(x,1)$. Example~\ref{ex:1+sqrt3-lagrange1770} gives $\alpha=[\overline{2,1}]$; hence the period is $(2,1)$ and its length is $k=2$. The purely periodic continued fraction determined by this period is $\beta=[\overline{2,1}]$. In this case $\beta=\alpha$, and its quadratic conjugate is $\beta'=1-\sqrt{3}$. Since $\beta>1$ and $-1<\beta'<0$, the form $Q$ is canonical reduced. By Proposition~\ref{prop:equivalent-calculatable} and Theorem~\ref{thm:markov-infinity-sequence}, for the bi-infinite sequence $\mathbf b=(\dots,2,1,2,1,2,1,2,1,\dots)$ obtained by repeating the period $(2,1)$ in both directions, we have $\mathcal M(Q)=\mathcal S(\mathbf b)$. On the other hand, the previous example showed that $\mathcal L(1+\sqrt{3})=2\sqrt{3}$. Therefore Theorem~\ref{thm:quadratic-markov-lagrange1770} gives $\mathcal M(Q)=\mathcal L(\alpha)=2\sqrt{3}$.
\end{exam}
In fact, Example~\ref{ex:concrete-example} and Example~\ref{ex:concrete-example2} also correspond to each other through Theorem~\ref{thm:quadratic-markov-lagrange1770}. The reader may check this directly.

We conclude the chapter by summarizing the relation between the Markov spectrum and the Lagrange spectrum obtained so far. Corollary~\ref{cor:M=S} shows that in general $\mathcal L\subset\mathcal M$, and it is known that the reverse inclusion does not hold. However, appropriate restrictions of the two sets are exactly the same.
\begin{coro}\label{cor:coincide}
Let
\[
\mathcal L_2:=\{\mathcal L(\alpha)\mid \alpha\in I_2\},
\qquad
\mathcal M_\mathbb Q:=\{\mathcal M(Q)\mid Q\in \mathcal Q_{\mathbb Q}\},
\]
and
\[
\mathcal S_{\mathrm{period}}:=\{\mathcal S(\mathbf b)\mid \text{$\mathbf b$ is a periodic bi-infinite sequence}\}.
\]
Then
\[
\mathcal L_2=\mathcal M_\mathbb Q=\mathcal S_{\mathrm{period}}.
\]
\end{coro}
\begin{proof}
This follows from Corollary~\ref{cor:quadratic-lagrange1770} and Theorem~\ref{thm:quadratic-markov-lagrange1770}.
\end{proof}
From the next chapter onward, we study values in $\mathcal L_2$ and $\mathcal M_\mathbb Q$ that admit special descriptions in terms of generalized Markov numbers.

\part[Generalized Markov Numbers]{Generalized Markov Numbers}
\chapter[Generalized Markov Equations and Generalized Markov Numbers]{Generalized Markov Equations and Generalized Markov Numbers}\label{chap:generalized-markov-equations-numbers}

In this chapter we introduce generalized Markov equations and generalized Markov numbers, which form the starting point of the second part of the text. Up to Chapter~\ref{chap:markov-spectrum}, we studied the Lagrange spectrum and the Markov spectrum. The purpose of the present chapter is to prepare the arithmetic objects that will later be used to describe a discrete family of values in that theory. We first define the generalized Markov equation and record its basic properties. We then construct generalized Markov trees and organize all positive integer solutions. After that, we introduce fraction labels through the correspondence with the Farey tree. Finally, we define characteristic numbers, which will later be used to describe generalized Cohn matrices and generalized discrete Markov spectra.

The material in this chapter is based on papers by the author and collaborators \cite{gyomatsu,gyo-maru,gyoda-maruyama-sato,bana-gyo}. The papers \cite{gyo-maru,gyoda-maruyama-sato} give proofs only in the symmetric case $k_1=k_2=k_3$; here we rewrite the arguments in full generality. The paper \cite{bana-gyo} treats a more general framework coming from the theory of cluster algebras; here we specialize those arguments so that the discussion remains within elementary number theory.

\section{Definitions and Basic Properties}
\begin{defi}\label{def:gm-equation}\index{Markov equation}\index{Generalized Markov equation}\index{Generalized Markov number}\index{Generalized Markov triple}\index{Markov number}\index{Markov triple}\index{GM number|seeonly{Generalized Markov number}}
For $k_1,k_2,k_3\in \mathbb Z_{\geq 0}$, the \emph{$(k_1,k_2,k_3)$-generalized Markov equation} is
\begin{equation}\label{Diophantine}
x^2 + y^2 + z^2 + k_1 yz + k_2 zx + k_3 xy = (3 + k_1 + k_2 + k_3) xyz.
\end{equation}
A permutation of a positive integer solution $(x,y,z)$ of this equation is called a \emph{$(k_1,k_2,k_3)$-generalized Markov triple}. A positive integer that occurs in a $(k_1,k_2,k_3)$-generalized Markov triple is called a \emph{$(k_1,k_2,k_3)$-generalized Markov number}. When $k_1=k_2=k_3=0$, we simply call them Markov triples and Markov numbers.
\end{defi}
Since these names are long, we henceforth abbreviate ``generalized Markov'' to ``GM''. One point should be kept in mind. The $(k_1,k_2,k_3)$-GM triples are the positive integer solutions of the $(k_1,k_2,k_3)$-GM equation together with all permutations of those solutions. If $k_1=k_2=k_3$, then the equation is symmetric in the three variables $x,y,z$, and hence every permutation of a positive integer solution is again a positive integer solution. In that case there is no difference between the two notions. If the parameters are not all equal, however, a permutation of a positive integer solution need not be a positive integer solution.

Let us first discuss basic properties of the positive integer solutions of the GM equation. At this stage the order of the three components is taken into account, so we are not yet speaking about GM triples. We begin with an algorithm that enumerates all positive integer solutions. Define a tree $\mathbb T(k_1,k_2,k_3)$ whose vertices are triples of positive integers as follows.
\begin{itemize}
    \item[(1)] The initial vertex is $(1,1,1)$.
    \item[(2)] The triple $(1,1,1)$ has the following three children:
    $(k_1+2,1,1)$, $(1,k_2+2,1)$, and $(1,1,k_3+2)$.
    \item[(3)] At every vertex other than the initial vertex, the generation rule is as follows.
    \begin{itemize}
        \item[(i)] If $a$ is the largest component of $(a,b,c)$, then $(a,b,c)$ has the following two children:
        \[
        \left(a,\dfrac{a^2+k_2ac+c^2}{b},c\right),
        \quad
        \left(a,b,\dfrac{a^2+k_3ab+b^2}{c}\right).
        \]
        \item[(ii)] If $b$ is the largest component of $(a,b,c)$, then $(a,b,c)$ has the following two children:
        \[
        \left(\dfrac{b^2+k_1bc+c^2}{a},b,c\right),
        \quad
        \left(a,b,\dfrac{a^2+k_3ab+b^2}{c}\right).
        \]
        \item[(iii)] If $c$ is the largest component of $(a,b,c)$, then $(a,b,c)$ has the following two children:
        \[
        \left(\dfrac{b^2+k_1bc+c^2}{a},b,c\right),
        \quad
        \left(a,\dfrac{a^2+k_2ac+c^2}{b},c\right).
        \]
    \end{itemize}
\end{itemize}

\begin{exam}
For $k_1=1,k_2=2,k_3=0$, the first few vertices of $\mathbb T(1,2,0)$ are as follows:
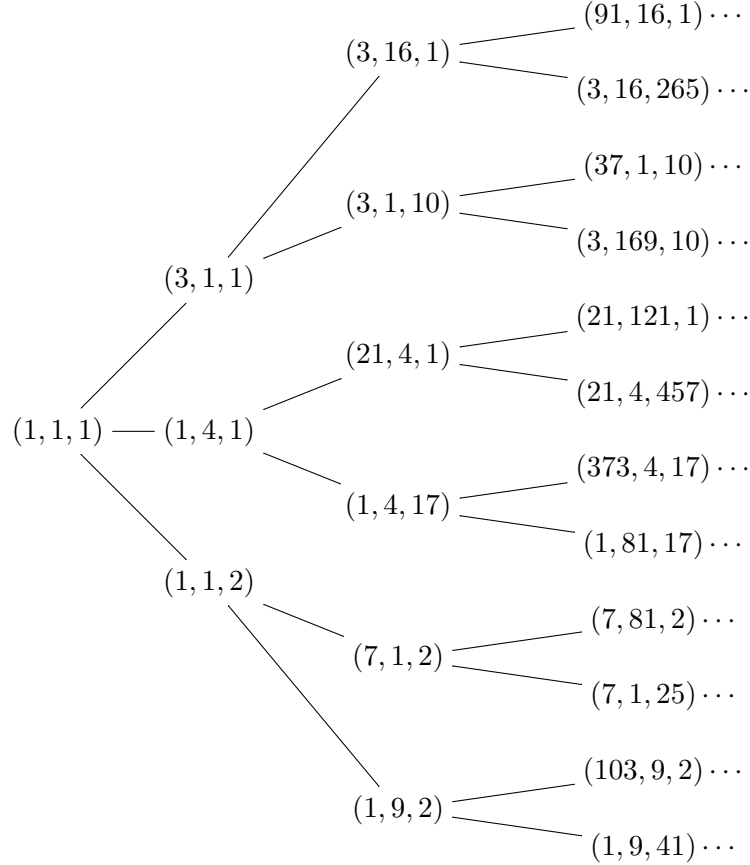
\[
\begin{xy}
(0,0)*+{(1,1,1)}="0",
(20,20)*+{(3,1,1)}="1",
(20,0)*+{(1,4,1)}="1'",
(20,-20)*+{(1,1,2)}="1''",
(45,50)*+{(3,16,1)}="20",
(45,30)*+{(3,1,10)}="21",
(45,10)*+{(21,4,1)}="22",
(45,-10)*+{(1,4,17)}="23",
(45,-30)*+{(7,1,2)}="24",
(45,-50)*+{(1,9,2)}="25",
(80,55)*+{(91,16,1)\cdots}="40",
(80,45)*+{(3,16,265)\cdots}="41",
(80,35)*+{(37,1,10)\cdots}="42",
(80,25)*+{(3,169,10)\cdots}="43",
(80,15)*+{(21,121,1)\cdots}="44",
(80,5)*+{(21,4,457)\cdots}="45",
(80,-5)*+{(373,4,17)\cdots}="46",
(80,-15)*+{(1,81,17)\cdots}="47",
(80,-25)*+{(7,81,2)\cdots}="48",
(80,-35)*+{(7,1,25)\cdots}="49",
(80,-45)*+{(103,9,2)\cdots}="410",
(80,-55)*+{(1,9,41)\cdots}="411",
\ar@{-}"0";"1"\ar@{-}"0";"1'"\ar@{-}"0";"1''"
\ar@{-}"1";"20"\ar@{-}"1";"21"
\ar@{-}"1'";"22"\ar@{-}"1'";"23"
\ar@{-}"1''";"24"\ar@{-}"1''";"25"
\ar@{-}"20";"40"\ar@{-}"20";"41"
\ar@{-}"21";"42"\ar@{-}"21";"43"
\ar@{-}"22";"44"\ar@{-}"22";"45"
\ar@{-}"23";"46"\ar@{-}"23";"47"
\ar@{-}"24";"48"\ar@{-}"24";"49"
\ar@{-}"25";"410"\ar@{-}"25";"411"
\end{xy}
\]
\end{exam}

We have the following theorem.

\begin{theorem}\label{Diophantinetheorem}\index{Generalized Markov equation!descent algorithm}\index{Generalized Markov equation}
Every positive integer solution of the $(k_1,k_2,k_3)$-GM equation appears exactly once in $\mathbb T(k_1,k_2,k_3)$.
\end{theorem}
We prepare the proof with the following proposition.

\begin{prop}\label{inductive-solution}
Suppose that $(x,y,z)=(a,b,c)$ is a positive integer solution of \eqref{Diophantine}. Then
\[
\left(\dfrac{b^2+k_1bc+c^2}{a},b,c\right),\quad
\left(a,\dfrac{a^2+k_2ac+c^2}{b},c\right),\quad
\left(a,b,\dfrac{a^2+k_3ab+b^2}{c}\right)
\]
are also positive integer solutions of \eqref{Diophantine}.
\end{prop}

\begin{proof}
It suffices to prove the assertion for
$\left(\dfrac{b^2+k_1bc+c^2}{a},b,c\right)$.
Positivity is clear, so it remains to prove that this is a triple of integers and satisfies \eqref{Diophantine}. Since $(a,b,c)$ is a solution of \eqref{Diophantine}, we have
\[
\dfrac{b^2+k_1bc+c^2}{a}
=(3+k_1+k_2+k_3)bc-a-k_3b-k_2c.
\]
Thus $\left(\dfrac{b^2+k_1bc+c^2}{a},b,c\right)$ is a triple of integers. To make the following computation easier to read, put $A:=(3+k_1+k_2+k_3)bc-a-k_3b-k_2c$. We show that $(A,b,c)$ is a solution of \eqref{Diophantine}.

The sum and product of $a$ and $A$ are
\[
a+A=(3+k_1+k_2+k_3)bc-k_3b-k_2c,\qquad a\cdot A=b^2+k_1bc+c^2.
\]
By the relation between roots and coefficients, $a$ and $A$ are the two roots of
\[
X^2-\{(3+k_1+k_2+k_3)bc-k_3b-k_2c\}X+b^2+k_1bc+c^2=0.
\]
Substituting $X=A$ into this quadratic equation and rearranging, we obtain
\[
A^2+b^2+c^2+k_3Ab+k_1bc+k_2cA=(3+k_1+k_2+k_3)Abc.
\]
This is exactly \eqref{Diophantine} with $(x,y,z)=(A,b,c)$.
\end{proof}

\index{Vieta jump}We call the three operations
\begin{align*}
(a,b,c)&\mapsto\left(\dfrac{b^2+k_1bc+c^2}{a},b,c\right),
\\
(a,b,c)&\mapsto\left(a,\dfrac{a^2+k_2ac+c^2}{b},c\right),\\
(a,b,c)&\mapsto\left(a,b,\dfrac{a^2+k_3ab+b^2}{c}\right)
\end{align*}
the \emph{first, second, and third Vieta jumps}, respectively. Each Vieta jump is an involution, namely applying the same operation again returns the original triple. We next determine the solutions that contain two equal components.

\begin{lemm}\label{singular}\index{Singular solution}
The positive integer solutions of \eqref{Diophantine} that contain equal components are only
\[
(1,1,1),\ (k_1+2,1,1),\ (1,k_2+2,1),\ (1,1,k_3+2).
\]
\end{lemm}

\begin{proof}
Let $(a,b,c)$ be a positive integer solution of \eqref{Diophantine} that contains equal components. We prove only the case $a=b$. Substituting $(x,y,z)=(a,a,c)$ into \eqref{Diophantine}, we obtain
\[
(2+k_3)a^2+c^2+(k_1+k_2)ac=(3+k_1+k_2+k_3)a^2c.
\]
Therefore
\[
c=\frac{1}{2}\Bigl(
 a^2k_3+a^2k_1+a^2k_2+3a^2-ak_1-ak_2
\pm a\sqrt{(ak_3+(a-1)k_1+(a-1)k_2+3a)^2-4(k_3+2)}
\Bigr).
\]
Put $k:=ak_3+(a-1)k_1+(a-1)k_2+3a>0$. Since $c$ is an integer, the expression under the square root must be a square. Hence there exists a positive integer $l$ such that
$l^2=k^2-4(k_3+2)$.
Because $a\ge1$, we have $k\ge k_3+3$, and hence $k+l>k_3+2$. From $(k+l)(k-l)=4(k_3+2)$, it follows that
$1\le k-l\le3$,
and so
\[
(k-l,k+l)=(1,4(k_3+2)),\ (2,2(k_3+2)),\ \left(3,\frac{4(k_3+2)}{3}\right).
\]
Since $k=((k+l)+(k-l))/2$ must be an integer, the first and the third possibilities are impossible. In the case $(k-l,k+l)=(2,2(k_3+2))$, we obtain
$k=k_3+3$ and $l=k_3+1$, and hence obtain
$(a,a,c)=(1,1,1)$ or $(1,1,k_3+2)$.
The cases $a=c$ and $b=c$ are proved in the same way, with the corresponding parameter $k_2$ or $k_1$ in place of $k_3$.
\end{proof}

\index{Singular solution}We call the four triples below
\[
(1,1,1),\ (k_1+2,1,1),\ (1,k_2+2,1),\ (1,1,k_3+2)
\]
\emph{singular}. \index{Nonsingular solution}We call every other positive integer solution of \eqref{Diophantine} \emph{nonsingular}.

\begin{prop}\label{nonsingular-induction}
Let $(a,b,c)$ be a nonsingular positive integer solution of \eqref{Diophantine}, and denote its three Vieta jumps by $(a',b,c)$, $(a,b',c)$, and $(a,b,c')$. If $a>\max\{b,c\}$, then
\[
a'<\max\{b,c\}<a,\qquad b'>a,\qquad c'>a.
\]
Thus the jump at the largest component strictly decreases the maximum, whereas each of the other two jumps creates a new largest component.
\end{prop}
\begin{proof}
The two increasing components satisfy
\[
b'=\frac{a^2+k_2ac+c^2}{b}>\frac{a^2}{b}>a,\qquad
c'=\frac{a^2+k_3ab+b^2}{c}>\frac{a^2}{c}>a.
\]
To estimate $a'$, note that $a,a'$ are the roots of
\[
f(X)=X^2-\bigl((3+k_1+k_2+k_3)bc-k_3b-k_2c\bigr)X+(b^2+k_1bc+c^2).
\]
For nonnegative integers $\rho,\mu,\nu$, put
\[
K:=3+\rho+\mu+\nu,\qquad
G_{\rho;\mu,\nu}(u,v):=(2+\rho)u^2-Ku^2v+(\mu+\nu)uv+v^2.
\]
Rearranging gives
\[
G_{\rho;\mu,\nu}(u,v)
=v^2-u^2-(3+\rho)u^2(v-1)-(\mu+\nu)uv(u-1).
\]
If $u>v\geq1$, then $v^2-u^2<0$, while the other two terms are nonpositive because $\rho,\mu,\nu\geq0$. Hence
\begin{equation}\label{equation}
G_{\rho;\mu,\nu}(u,v)<0\qquad(u>v\geq1).
\end{equation}
By Lemma~\ref{singular}, the three components of a nonsingular solution are distinct. If $b>c$, then $f(b)=G_{k_3;k_1,k_2}(b,c)<0$. Since $f$ is monic and $a>b$, the number $b$ lies between its roots, so $a'<b$. If $c>b$, use $f(c)=G_{k_2;k_1,k_3}(c,b)<0$ to obtain $a'<c$. Thus $a'<\max\{b,c\}$ in either case.
\end{proof}

\begin{coro}\label{nonsingular-induction-general}\index{Generalized Markov equation!descent algorithm}\index{Nonsingular solution}\index{Vieta jump}
Let $(a,b,c)$ be a nonsingular positive integer solution of \eqref{Diophantine}, and denote its three Vieta jumps by $(a',b,c)$, $(a,b',c)$, and $(a,b,c')$.
\begin{itemize}
\item[(1)] If $a$ is largest, then $a'<\max\{b,c\}<a$, $b'>a$, and $c'>a$.
\item[(2)] If $b$ is largest, then $b'<\max\{a,c\}<b$, $a'>b$, and $c'>b$.
\item[(3)] If $c$ is largest, then $c'<\max\{a,b\}<c$, $a'>c$, and $b'>c$.
\end{itemize}
In particular, the jump at the largest component is the unique Vieta jump that decreases the maximum.
\end{coro}
\begin{proof}
Assertion (1) is Proposition~\ref{nonsingular-induction}. For (2), observe that $(b,c,a)$ solves the $(k_2,k_3,k_1)$-GM equation. Apply the proposition after simultaneously rotating the variables and their opposite coefficients. Similarly, $(c,a,b)$ solves the $(k_3,k_1,k_2)$-GM equation, which proves (3).
\end{proof}

\begin{rema}\label{remark-adjacent}
Corollary~\ref{nonsingular-induction-general} gives a canonical descent for every nonsingular solution. Replacing its unique largest component by the other root makes the new maximum equal to the larger of the two unchanged components, strictly below the old maximum. Each other Vieta jump produces a component larger than all three old components.

For the three singular solutions other than the root, direct calculation gives
\[
(k_1+2,1,1)\longmapsto(1,1,1),\quad
(1,k_2+2,1)\longmapsto(1,1,1),\quad
(1,1,k_3+2)\longmapsto(1,1,1)
\]
under the jump at the largest component. Thus every nonroot vertex of $\mathbb T(k_1,k_2,k_3)$ has a unique parent, obtained by that jump. The other two jumps give its two children, with strictly larger maxima. At the root the three jumps give the three specified singular children. Consequently, the three Vieta jumps at a vertex give exactly its parent and children, or its three children at the root. A nonroot singular solution occurs only as the corresponding child of the root.
\end{rema}

\begin{proof}[Proof of Theorem~\ref{Diophantinetheorem}]
Proposition~\ref{inductive-solution} and the fact that $(1,1,1)$ solves \eqref{Diophantine} show that every vertex is a positive integer solution. Conversely, take any positive integer solution $(a,b,c)$. If it is singular, Lemma~\ref{singular} and the definition of the first generation place it in the tree. Otherwise its components are distinct, so it has a unique largest component. Repeatedly jump at that component. At every nonsingular stage, Corollary~\ref{nonsingular-induction-general} strictly decreases the positive integer maximum. The process therefore reaches a singular solution after finitely many steps. Reversing these involutive jumps gives a path in the tree by Remark~\ref{remark-adjacent}, so $(a,b,c)$ occurs in the tree.

The same descent proves uniqueness. Each nonroot singular solution has the root as its unique parent. For a nonsingular solution, the only adjacent solution with smaller maximum is obtained by jumping at the largest component. Thus the path back to the root is determined by the solution itself. Two occurrences of the same triple have the same parent at each stage and hence occupy the same position. Every positive integer solution therefore occurs exactly once.
\end{proof}

\begin{coro}\label{relatively-prime}\index{Generalized Markov number!coprimality}
For every positive integer solution $(a,b,c)$ of \eqref{Diophantine}, any two of $a,b,c$ are relatively prime.
\end{coro}

\begin{proof}
The assertion is clear for $(a,b,c)=(1,1,1)$. We prove only that $a$ and $b$ are relatively prime. Rewrite \eqref{Diophantine} as
\[
z^2=(3+k_1+k_2+k_3)xyz-x^2-y^2-k_1yz-k_2zx-k_3xy,
\]
and substitute $(x,y,z)=(a,b,c)$. Suppose that $a$ and $b$ have a positive common divisor $d>1$, and let $d'$ be a prime divisor of $d$. Reducing the displayed equation modulo $d'$ gives $c^2\equiv0\pmod{d'}$, because $d'$ divides both $a$ and $b$. Hence $d'$ divides $c$, and so $d'$ is a common divisor of $a,b,c$.

By Proposition~\ref{inductive-solution}, the adjacent triple $(a',b',c')$ in $\mathbb T(k_1,k_2,k_3)$ whose largest component is smaller than $\max\{a,b,c\}$ is obtained by replacing the largest component by the other root of the corresponding quadratic equation. If the first component is replaced, then the two roots have sum $(3+k_1+k_2+k_3)bc-k_2c-k_3b$, which is divisible by $d'$; since the old first component is also divisible by $d'$, the new one is divisible by $d'$ as well. The arguments for the second and third components are identical. Hence common divisibility by $d'$ is preserved along the descent. Repeating the operation, we conclude that $d'$ is eventually a common divisor of $(1,1,1)$, which is impossible. Thus no positive common divisor $d>1$ exists, and $a$ and $b$ are relatively prime.
\end{proof}

\section{Generalized Markov Trees}
The tree $\mathbb T(k_1,k_2,k_3)$ is economical from the point of view of enumerating all positive integer solutions of the $(k_1,k_2,k_3)$-GM equation. However, when we later compare it with the matrix theory, it will often be more convenient to decompose this tree into several binary trees. We therefore introduce new binary trees.
\begin{defi}
The full subtrees of $\mathbb{T}(k_1,k_2,k_3)$ whose initial vertices are respectively
$\left(k_1+2,1,1\right)$, $\left(1,k_2+2,1\right)$, and $(1,1,k_3+2)$
are called the \emph{first, second, and third branches} of $\mathbb{T}(k_1,k_2,k_3)$ and are denoted by
$\mathbb{T}_1(k_1,k_2,k_3)$,
$\mathbb{T}_2(k_1,k_2,k_3)$, and
$\mathbb{T}_3(k_1,k_2,k_3)$.
\end{defi}
\begin{exam}\label{ex:k2T}
The first few vertices of $\mathbb T_2(1,2,0)$ are as follows:
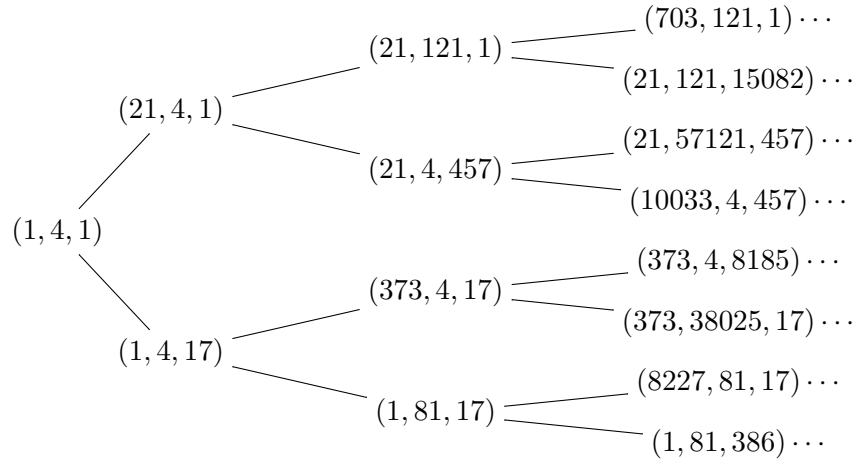
\[
\begin{xy}
(10,0)*+{(1,4,1)}="1",
(25,16)*+{(21,4,1)}="2",
(25,-16)*+{(1,4,17)}="3",
(60,24)*+{(21,121,1)}="4",
(60,8)*+{(21,4,457)}="5",
(60,-8)*+{(373,4,17)}="6",
(60,-24)*+{(1,81,17)}="7",
(100,28)*+{(703,121,1)\cdots}="8",
(100,20)*+{(21,121,15082)\cdots}="9",
(100,12)*+{(21,57121,457)\cdots}="10",
(100,4)*+{(10033,4,457)\cdots}="11",
(100,-4)*+{(373,4,8185)\cdots}="12",
(100,-12)*+{(373,38025,17)\cdots}="13",
(100,-20)*+{(8227,81,17)\cdots}="14",
(100,-28)*+{(1,81,386)\cdots}="15",
\ar@{-}"1";"2"\ar@{-}"1";"3"
\ar@{-}"2";"4"\ar@{-}"2";"5"
\ar@{-}"3";"6"\ar@{-}"3";"7"
\ar@{-}"4";"8"\ar@{-}"4";"9"
\ar@{-}"5";"10"\ar@{-}"5";"11"
\ar@{-}"6";"12"\ar@{-}"6";"13"
\ar@{-}"7";"14"\ar@{-}"7";"15"
\end{xy}
\]
\end{exam}
By definition, each $\mathbb{T}_i(k_1,k_2,k_3)$ is a complete binary tree. In the theory below, however, we will actually use the following complete binary trees, obtained by rearranging the components at each vertex so that the newly produced component is written in the middle. The branches above will be used to prove that the vertices of the binary trees defined below enumerate all positive integer solutions.

\begin{defi}\label{def:gm-tree}\index{Generalized Markov tree}
Let $\mathfrak S_3$ be the symmetric group of degree $3$, acting on the left on $\{1,2,3\}$. For $\sigma\in \mathfrak S_3$, define the \emph{$(k_1,k_2,k_3,\sigma)$-generalized Markov tree} (or simply the \emph{GM tree}) $\mathrm{M}\mathbb T(k_1,k_2,k_3,\sigma)$ as follows.
\begin{itemize}
\item[(1)] The initial vertex is
\[
\big((1,\sigma(1)), (k_{\sigma(2)}+2, \sigma(2)), (1,\sigma(3))\big).
\]
\item[(2)] For every vertex $((a,\alpha),(b,\beta),(c,\gamma))$, define the following two children, distinguishing the left child from the right child:
\[
\begin{xy}
(0,0)*+{((a,\alpha),(b,\beta),(c,\gamma))}="1",
(-40,-15)*+{\left((a,\alpha),\left(\frac{a^2+k_\gamma ab+b^2}{c},\gamma\right),(b,\beta)\right)}="2",
(40,-15)*+{\left((b,\beta),\left(\frac{b^2+k_\alpha bc+c^2}{a},\alpha\right),(c,\gamma)\right)}="3",
\ar@{-}"1";"2"\ar@{-}"1";"3"
\end{xy}
\]
\end{itemize}
\end{defi}
\begin{exam}\label{ex:CMT(0,1,2,id)}
The first few vertices of $\mathrm{M}\mathbb T(1,2,0,\mathrm{id})$ are as follows:
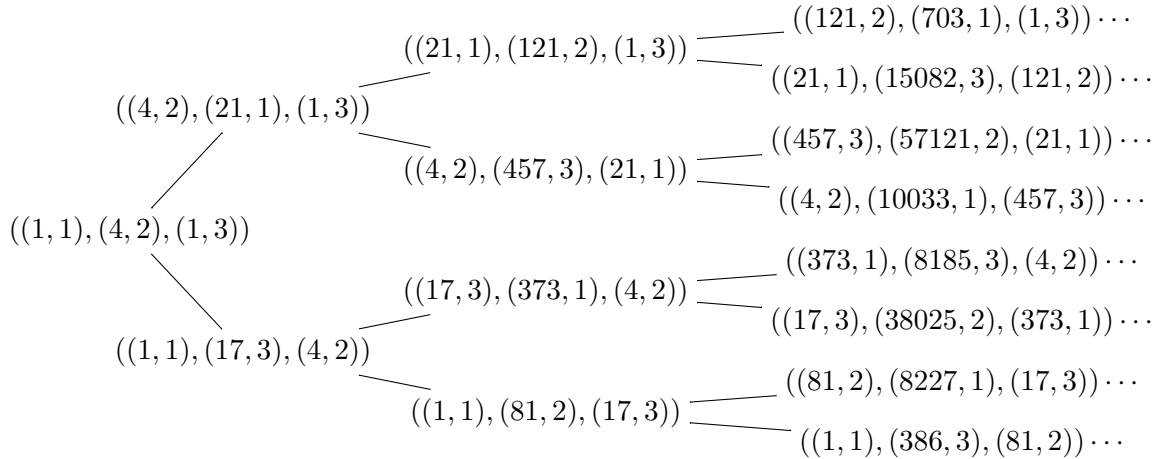
\begin{center}
\adjustbox{max width=.98\linewidth}{$
\begin{xy}
(10,0)*+{((1,1),(4,2),(1,3))}="1",
(25,16)*+{((4,2),(21,1),(1,3))}="2",
(25,-16)*+{((1,1),(17,3),(4,2))}="3",
(65,24)*+{((21,1),(121,2),(1,3))}="4",
(65,8)*+{((4,2),(457,3),(21,1))}="5",
(65,-8)*+{((17,3),(373,1),(4,2))}="6",
(65,-24)*+{((1,1),(81,2),(17,3))}="7",
(120,28)*+{((121,2),(703,1),(1,3))\cdots}="8",
(120,20)*+{((21,1),(15082,3),(121,2))\cdots}="9",
(120,12)*+{((457,3),(57121,2),(21,1))\cdots}="10",
(120,4)*+{((4,2),(10033,1),(457,3))\cdots}="11",
(120,-4)*+{((373,1),(8185,3),(4,2))\cdots}="12",
(120,-12)*+{((17,3),(38025,2),(373,1))\cdots}="13",
(120,-20)*+{((81,2),(8227,1),(17,3))\cdots}="14",
(120,-28)*+{((1,1),(386,3),(81,2))\cdots}="15",
\ar@{-}"1";"2"\ar@{-}"1";"3"
\ar@{-}"2";"4"\ar@{-}"2";"5"
\ar@{-}"3";"6"\ar@{-}"3";"7"
\ar@{-}"4";"8"\ar@{-}"4";"9"
\ar@{-}"5";"10"\ar@{-}"5";"11"
\ar@{-}"6";"12"\ar@{-}"6";"13"
\ar@{-}"7";"14"\ar@{-}"7";"15"
\end{xy}
$}
\end{center}
\end{exam}
\index{Position label}Each vertex of this tree consists of three pairs, hence of six entries in total. The first entry of each pair is a GM number, while the second records its original coordinate in the positive integer solution. When the parameters are not all equal, permuting the coordinates need not preserve the GM equation. Since the components are rearranged at each generation, we retain each number's original coordinate as its position label. \index{Generalized Markov pair}We call a GM number together with its position label a \emph{$(k_1,k_2,k_3)$-GM pair}.

For $((a,h),(b,i),(c,j))\in \mathrm{M}\mathbb T(k_1,k_2,k_3,\sigma)$, we call the $m$-th entry of the $n$-th pair the $(n,m)$-entry.

A GM pair at a specified position (first, second, or third) of a specified vertex is called an \emph{occurrence} of that pair. Equal numerical GM pairs may have distinct occurrences. Occurrences inherited from a parent to its children will be assigned the same fraction label. Whether different fraction labels can have equal numerical values is a separate question addressed later.

We now show that $\mathrm{M}\mathbb T(k_1,k_2,k_3,\sigma)$ is the same tree as a branch of $\mathbb T(k_1,k_2,k_3)$, up to the difference caused by rearrangements. To do this, we first introduce isomorphisms of complete binary trees.

\begin{defi}
Let $A$ and $B$ be regarded as complete binary trees, that is, assume that their elements correspond to the vertices of complete binary trees. A bijection $f:A\to B$ is called a \emph{complete binary tree isomorphism} if, for vertices $u$ and $v$ of $A$, the vertex $v$ is a child of $u$ if and only if $f(v)$ is a child of $f(u)$. If, in addition, $A$ and $B$ are ordered complete binary trees, meaning that the two children of each vertex are distinguished as the left and right children, and if $f$ also preserves the distinction between left and right children, then $f$ is called an \emph{ordered complete binary tree isomorphism}.
\end{defi}

Note that $\mathbb T_i(k_1,k_2,k_3)$ is an unordered complete binary tree, whereas $\mathrm{M}\mathbb T(k_1,k_2,k_3,\sigma)$ is an ordered complete binary tree. We next construct a bijection that gives an unordered complete binary tree isomorphism between $\mathbb T_i(k_1,k_2,k_3)$ and $\mathrm{M}\mathbb T(k_1,k_2,k_3,\sigma)$.

For $\tau\in\mathfrak S_3$, let $^\tau(a,b,c)$ denote the triple obtained by permuting $(a,b,c)$ so that $a$, $b$, and $c$ become the $\tau(1)$-st, $\tau(2)$-nd, and $\tau(3)$-rd components, respectively. For a tree $\mathbb T$, let $V(\mathbb T)$ denote its vertex set. Define a map
\[
\pi_{\sigma}:V(\mathrm{M}\mathbb T(k_1,k_2,k_3,\sigma))
\to\mathbb Z_{> 0}^3
\]
as follows. For $v=((a,\alpha),(b,\beta),(c,\gamma))$, define $\tau_v\in\mathfrak S_3$ by
\[
\tau_v(1)=\alpha,\qquad\tau_v(2)=\beta,\qquad\tau_v(3)=\gamma,
\]
and put $\pi_\sigma(v):=\prescript{\tau_v}{}{(a,b,c)}$. Thus the $\alpha$-, $\beta$-, and $\gamma$-components of the image are $a,b,c$, respectively. This permutation is determined by position labels, even at the root where $a=c=1$. At the initial vertex,
\[
\pi_\sigma((1,\sigma(1)),(k_{\sigma(2)}+2,\sigma(2)),(1,\sigma(3)))
=\prescript{\sigma}{}{(1,k_{\sigma(2)}+2,1)}.
\]
If $\tau=\tau_v$, inspection of the child labels gives
\begin{align*}
\pi_\sigma\left((a,\alpha),\left(\frac{a^2+k_\gamma ab+b^2}{c},\gamma\right),(b,\beta)\right)
&=\prescript{\tau\circ(2\ 3)}{}{\left(a,\frac{a^2+k_\gamma ab+b^2}{c},b\right)},\\
\pi_\sigma\left((b,\beta),\left(\frac{b^2+k_\alpha bc+c^2}{a},\alpha\right),(c,\gamma)\right)
&=\prescript{\tau\circ(1\ 2)}{}{\left(b,\frac{b^2+k_\alpha bc+c^2}{a},c\right)}.
\end{align*}
\begin{prop}\label{prop:maximal-middle-component}\index{Generalized Markov tree!maximal middle component}
For every vertex $((a,\alpha),(b,\beta),(c,\gamma))$ of $\mathrm M\mathbb T(k_1,k_2,k_3,\sigma)$, the middle GM number satisfies
\[
b>\max\{a,c\}.
\]
\end{prop}
\begin{proof}
At the initial vertex, $b=k_{\sigma(2)}+2>1=a=c$. If $b>\max\{a,c\}$ holds at a vertex, the new middle components of its left and right children satisfy, respectively,
\[
\frac{a^2+k_\gamma ab+b^2}{c}>\frac{b^2}{c}>b,
\qquad
\frac{b^2+k_\alpha bc+c^2}{a}>\frac{b^2}{a}>b.
\]
Each exceeds both inherited components, so the assertion follows by induction on the distance from the initial vertex.
\end{proof}

\begin{prop}\label{prop:isomorphism-CMT-kT}\index{Generalized Markov tree}
The map $\pi_\sigma$ induces the following complete binary tree isomorphisms:
\begin{align*}
\pi_{(1\ 2)}&\colon\mathrm{M}\mathbb T(k_1,k_2,k_3,(1\ 2))\simeq \TT_1(k_1,k_2,k_3),\\
\pi_{(1\ 3\ 2)}&\colon\mathrm{M}\mathbb T(k_1,k_2,k_3,(1\ 3\ 2))\simeq \TT_1(k_1,k_2,k_3),\\
\pi_{\mathrm{id}}&\colon\mathrm{M}\mathbb T(k_1,k_2,k_3,\mathrm{id})\simeq\TT_2(k_1,k_2,k_3),\\
\pi_{(1\ 3)}&\colon\mathrm{M}\mathbb T(k_1,k_2,k_3,(1\ 3))\simeq \TT_2(k_1,k_2,k_3),\\
\pi_{(2\ 3)}&\colon\mathrm{M}\mathbb T(k_1,k_2,k_3,(2\ 3))\simeq \TT_3(k_1,k_2,k_3),\\
\pi_{(1\ 2\ 3)}&\colon\mathrm{M}\mathbb T(k_1,k_2,k_3,(1\ 2\ 3))\simeq \TT_3(k_1,k_2,k_3).
\end{align*}
\end{prop}

\begin{proof}
By Proposition~\ref{prop:maximal-middle-component}, the middle GM number at every vertex is strictly larger than the other two. We prove the isomorphism assertion for $\sigma=(1\ 2)$; the other cases are identical. We use the following six states. A vertex
$((a,\alpha),(b,\beta),(c,\gamma))$ is said to be in the state indexed by $\tau\in\mathfrak S_3$ if
\[
\pi_{(1\ 2)}((a,\alpha),(b,\beta),(c,\gamma))=\prescript{\tau}{}{(a,b,c)}
\quad\text{and}\quad
(\alpha,\beta,\gamma)=(\tau(1),\tau(2),\tau(3)).
\]
For the six possible values of $\tau$, the corresponding index triple and the largest component of $\prescript{\tau}{}{(a,b,c)}$ are as follows:
\[
\begin{array}{|c|c|c|c|}
\hline
\text{State} & \tau & (\alpha,\beta,\gamma) &
\begin{array}{c}
\text{Largest component of }\prescript{\tau}{}{(a,b,c)}
\end{array}\\
\hline
(1)&\mathrm{id}&(1,2,3)&\text{second component}\\
(2)&(2\ 3)&(1,3,2)&\text{third component}\\
(3)&(1\ 2)&(2,1,3)&\text{first component}\\
(4)&(1\ 3)&(3,2,1)&\text{second component}\\
(5)&(1\ 2\ 3)&(2,3,1)&\text{third component}\\
(6)&(1\ 3\ 2)&(3,1,2)&\text{first component}\\
\hline
\end{array}
\]
In each row this means precisely that the middle entry $b$ of the displayed vertex is the largest of $a,b,c$.

Suppose that a vertex is in the state indexed by $\tau$. By the definition of $\pi_{(1\ 2)}$, its left child is in the state indexed by $\tau\circ(2\ 3)$, and its right child is in the state indexed by $\tau\circ(1\ 2)$. Thus the transitions are
\[
\begin{array}{|c|c|c|}
\hline
\text{State} & \text{left child} & \text{right child}\\
\hline
(1)&(2)&(3)\\
(2)&(1)&(6)\\
(3)&(5)&(1)\\
(4)&(6)&(5)\\
(5)&(3)&(4)\\
(6)&(4)&(2)\\
\hline
\end{array}
\]
The initial vertex of $\mathrm{M}\mathbb T(k_1,k_2,k_3,(1\ 2))$ is
$((1,2),(k_1+2,1),(1,3))$, which is in state $(3)$ and is mapped by $\pi_{(1\ 2)}$ to $(k_1+2,1,1)$, the initial vertex of $\TT_1(k_1,k_2,k_3)$.

Now assume that a noninitial vertex is in one of the six states. Since its middle entry is the largest of $a,b,c$, Corollary~\ref{nonsingular-induction-general} and Remark~\ref{remark-adjacent} show that replacing either the first or the third entry produces the two children in the corresponding branch of $\mathbb T(k_1,k_2,k_3)$. The same conclusion for the initial vertex is checked directly. The formulas defining the two children of $\mathrm{M}\mathbb T(k_1,k_2,k_3,(1\ 2))$ are exactly these two Vieta jumps, written after the rearrangement that places the newly produced entry in the middle. Hence $\pi_{(1\ 2)}$ sends the two children of each vertex of $\mathrm{M}\mathbb T(k_1,k_2,k_3,(1\ 2))$ to the two children of its image in $\TT_1(k_1,k_2,k_3)$.

It follows by induction on the distance from the initial vertex that $\pi_{(1\ 2)}$ gives a complete binary tree isomorphism from $\mathrm{M}\mathbb T(k_1,k_2,k_3,(1\ 2))$ onto $\TT_1(k_1,k_2,k_3)$. The proof for the other five values of $\sigma$ is the same.
\end{proof}
\begin{rema}\label{rem:mirror}
The isomorphism $\pi$ does not distinguish left and right children. Hence, if $\sigma^{\ast}:=\sigma\circ(1\ 3)$, then for every $\sigma$ the two trees $\mathrm{M}\mathbb T(k_1,k_2,k_3,\sigma)$ and $\mathrm{M}\mathbb T(k_1,k_2,k_3,\sigma^{\ast})$ are sent by $\pi$ to the same $\TT_i(k_1,k_2,k_3)$.
\end{rema}
\begin{coro}\label{cor:uniqueness-vertex}
Let $V_{(k_1,k_2,k_3)}$ be the set of all vertices of the six GM trees $\mathrm{M}\mathbb T(k_1,k_2,k_3,\sigma)$ with $\sigma\in\mathfrak S_3$. For every $v\in V_{(k_1,k_2,k_3)}$, the position at which $v$ appears is unique in the union of these six trees.
\end{coro}
\begin{proof}
First, if $v=((a,\alpha),(b,\beta),(c,\gamma))\in V_{(k_1,k_2,k_3)}$, then the same displayed vertex cannot appear more than once in the tree $\mathrm{M}\mathbb T(k_1,k_2,k_3,\sigma)$ to which $v$ belongs; otherwise the bijectivity of Proposition~\ref{prop:isomorphism-CMT-kT} would be contradicted. We prove that $v$ cannot appear in two different trees $\mathrm{M}\mathbb T(k_1,k_2,k_3,\sigma)$ and $\mathrm{M}\mathbb T(k_1,k_2,k_3,\sigma')$.
Whether $v$ is the initial vertex, and if it is not the initial vertex whether it is the left or right child of its parent, can be determined from the data of $v$ itself. Indeed, if $a=c$, then the image of $v$ under the corresponding map $\pi_\sigma$ is a singular solution. Since singular solutions occur only at the initial vertices of the branches by Remark~\ref{remark-adjacent}, the vertex $v$ itself is the initial vertex. Here Proposition~\ref{prop:maximal-middle-component} shows that neither $a=b$ nor $b=c$ can occur. If $a\neq c$, then $v$ is not the initial vertex and has a parent. In this case, according as $v$ is the left child or the right child of its parent, the parent is one of
\[\left((a,\alpha),(c,\gamma),\left(\frac{a^2+k_\beta ac+c^2}{b},\beta\right)\right),\quad \left(\left(\frac{a^2+k_{\beta}ac+c^2}{b},\beta\right),(a,\alpha),(c,\gamma)\right).\]
By Proposition~\ref{prop:maximal-middle-component}, the $(2,1)$-entry is strictly larger than the $(1,1)$- and $(3,1)$-entries. Therefore, if $a<c$ then the parent is the former, and if $a>c$ then the parent is the latter. Thus the parent is uniquely determined. It follows that the path from $v$ to the initial vertex is uniquely determined by $v$, and the displayed form of the initial vertex is also uniquely determined. Hence $v$ cannot belong to two different trees $\mathrm{M}\mathbb T(k_1,k_2,k_3,\sigma)$ and $\mathrm{M}\mathbb T(k_1,k_2,k_3,\sigma')$.
\end{proof}

The next theorem says that $\mathrm{M}\mathbb T(k_1,k_2,k_3,\sigma)$ enumerates GM triples whose second component is the largest one.

\begin{theorem}\label{thm:all-solution}\index{Generalized Markov tree}\index{Generalized Markov number}\index{Generalized Markov triple}
Let $(a,b,c)$ be a $(k_1,k_2,k_3)$-GM triple satisfying $b>\max\{a,c\}$, and suppose that for $\tau\in \mathfrak S_3$ the triple $\prescript{\tau}{}{(a,b,c)}$ is a solution of the $(k_1,k_2,k_3)$-GM equation. Then there exists a unique $\sigma\in \mathfrak S_3$ and unique vertices $v$ and $v^{\ast}$, where $\sigma^{\ast}:=\sigma\circ(1\ 3)$, such that
$v\in \mathrm{M}\mathbb T(k_1,k_2,k_3,\sigma)$, $v^{\ast}\in \mathrm{M}\mathbb T(k_1,k_2,k_3,\sigma^{\ast})$, $v=((a,\tau(1)),(b,\tau(2)),(c,\tau(3)))$, and $v^{\ast}=((c,\tau(3)),(b,\tau(2)),(a,\tau(1)))$. Moreover, the position of $v^{\ast}$ in $\mathrm{M}\mathbb T(k_1,k_2,k_3,\sigma^{\ast})$ is the mirror image of the position of $v$ in $\mathrm{M}\mathbb T(k_1,k_2,k_3,\sigma)$, obtained by interchanging left and right at every level.
\end{theorem}
\begin{proof}
Since $\prescript{\tau}{}{(a,b,c)}\neq (1,1,1)$, Theorem~\ref{Diophantinetheorem} implies that $\prescript{\tau}{}{(a,b,c)}$ belongs to one of $\mathbb T_1(k_1,k_2,k_3)$, $\mathbb T_2(k_1,k_2,k_3)$, and $\mathbb T_3(k_1,k_2,k_3)$. Suppose that $\prescript{\tau}{}{(a,b,c)}$ belongs to $\mathbb T_{1}(k_1,k_2,k_3)$; the other cases are handled in the same way. By Proposition~\ref{prop:isomorphism-CMT-kT}, there exist unique vertices $v_1$ of $\mathrm{M}\mathbb T(k_1,k_2,k_3,(1\ 2))$ and $v_2$ of $\mathrm{M}\mathbb T(k_1,k_2,k_3,(1\ 3\ 2))$ such that
$\pi_{(1\ 2)}(v_1)=\prescript{\tau}{}{(a,b,c)}$ and
$\pi_{(1\ 3\ 2)}(v_2)=\prescript{\tau}{}{(a,b,c)}$.
The triple $\prescript{\tau}{}{(a,b,c)}$ can be displayed in the following six ways:
\[
\prescript{\tau}{}{(a,b,c)},\ \prescript{\tau\circ(1\ 3)}{}{(c,b,a)},\ \prescript{\tau\circ(1\ 2)}{}{(b,a,c)},\ \prescript{\tau\circ(2\ 3)}{}{(a,c,b)},\ \prescript{\tau\circ(1\ 2\ 3)}{}{(b,c,a)},\ \prescript{\tau\circ(1\ 3\ 2)}{}{(c,a,b)}.
\]
Since $b>\max\{a,c\}$ by assumption, Proposition~\ref{prop:maximal-middle-component} leaves the following two cases:
\begin{itemize}
    \item[(1)] $v_1=((a,\tau(1)),(b,\tau(2)),(c,\tau(3)))$, $\pi_{(1\ 2)}(v_1) =\prescript{\tau}{}{(a,b,c)}$, and\\
    $v_2=((c,\tau\circ(1\ 3)(1)),(b,\tau\circ(1\ 3)(2)),(a,\tau\circ(1\ 3)(3)))$, $\pi_{(1\ 3\ 2)}(v_2) =\prescript{\tau\circ(1\ 3)}{}{(c,b,a)}$.
    \item[(2)] $v_1=((c,\tau\circ(1\ 3)(1)),(b,\tau\circ(1\ 3)(2)),(a,\tau\circ(1\ 3)(3)))$, $\pi_{(1\ 2)}(v_1) =\prescript{\tau\circ(1\ 3)}{}{(c,b,a)}$, and\\
    $v_2=((a,\tau(1)),(b,\tau(2)),(c,\tau(3)))$, $\pi_{(1\ 3\ 2)}(v_2) =\prescript{\tau}{}{(a,b,c)}$.
\end{itemize}
By Corollary~\ref{cor:uniqueness-vertex}, the same displayed vertex cannot occur in two different GM trees. Hence $v_1$ and $v_2$ cannot have the same displayed form, and exactly one of (1) and (2) occurs. In case (1), take $\sigma=(1\ 2)$, $\sigma^{\ast}=(1\ 3\ 2)$, $v=v_1$, and $v^{\ast}=v_2$. In case (2), take $\sigma=(1\ 3\ 2)$, $\sigma^{\ast}=(1\ 2)$, $v=v_2$, and $v^{\ast}=v_1$. Corollary~\ref{cor:uniqueness-vertex} also shows that no vertex with the same displayed form as $v$ or $v^{\ast}$ appears in any other tree. It remains to show that the position of $v^{\ast}$ in $\mathrm{M}\mathbb T(k_1,k_2,k_3,\sigma^{\ast})$ is obtained from the position of $v$ in $\mathrm{M}\mathbb T(k_1,k_2,k_3,\sigma)$ by interchanging left and right at every level. This follows by induction on the distance from the initial vertex. Indeed, the initial vertices of the two trees are obtained from one another by interchanging the first and third pairs. If two vertices are related in this way, then the left child of one is related to the right child of the other, and the right child of one is related to the left child of the other, by the defining formulas for the two children.
\end{proof}

\section{Farey Trees and Fraction Labels}\label{sec:farey-label}
We next introduce fraction labels for GM pairs. For this purpose, we first define Farey triples and the Farey tree.
\begin{defi}\label{def:farey-triple}\index{Farey triple}
For two fractions $\frac{a}{b}$ and $\frac{c}{d}$, write $\det\!\left(\frac{a}{b},\frac{c}{d}\right)$ for $ad-bc$.
A triple $\left(\frac{a}{b},\frac{c}{d},\frac{e}{f}\right)$ is called a \emph{Farey triple} if it satisfies the following conditions:
\begin{itemize}
\item[(1)] Each of $\frac{a}{b},\frac{c}{d},\frac{e}{f}$ is a reduced fraction.
\item[(2)]
\[
\left|\det\!\left(\frac{a}{b},\frac{c}{d}\right)\right|
=
\left|\det\!\left(\frac{c}{d},\frac{e}{f}\right)\right|
=
\left|\det\!\left(\frac{e}{f},\frac{a}{b}\right)\right|
=1.
\]
\end{itemize}
\end{defi}

\index{Farey tree}Define the \emph{Farey tree} $\mathrm{F}\mathbb T$ as follows.
\begin{itemize}
\item[(1)] The root vertex is $\left(\frac{0}{1},\frac{1}{1},\frac{1}{0}\right)$.
\item[(2)] Each vertex $\left(\frac{a}{b},\frac{c}{d},\frac{e}{f}\right)$ has the following two children:
\[
\begin{xy}(0,0)*+{\left(\dfrac{a}{b},\dfrac{c}{d},\dfrac{e}{f}\right)}="1",(-30,-15)*+{\left(\dfrac{a}{b},\dfrac{a}{b}\oplus\dfrac{c}{d},\dfrac{c}{d}\right)}="2",(30,-15)*+{\left(\dfrac{c}{d},\dfrac{c}{d}\oplus\dfrac{e}{f},\dfrac{e}{f}\right)}="3",\ar@{-}"1";"2"\ar@{-}"1";"3"
\end{xy}
\]
\index{Mediant}where $\frac{a}{b}\oplus\frac{c}{d}=\frac{a+c}{b+d}$. Here $\frac{1}{0}$ is treated as the endpoint $\infty$, and it is regarded as larger than every finite nonnegative fraction.
\end{itemize}
The first few vertices of $\mathrm{F}\mathbb T$ are as follows:
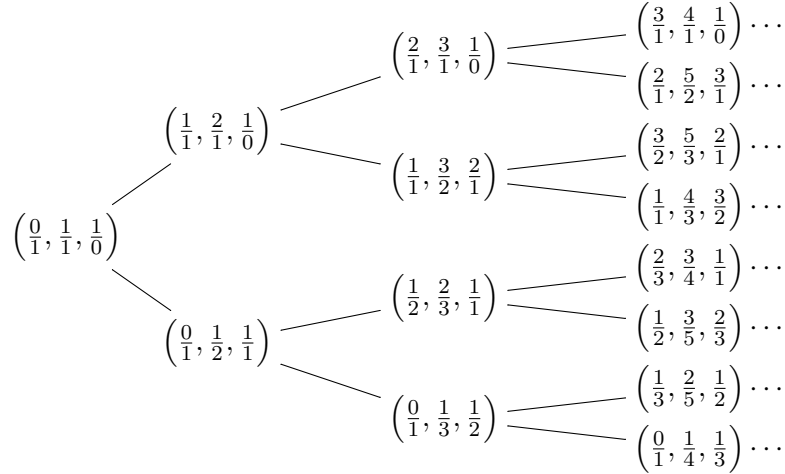
\[
\begin{xy}(0,0)*+{\left(\frac{0}{1},\frac{1}{1},\frac{1}{0}\right)}="1",(20,-14)*+{\left(\frac{0}{1},\frac{1}{2},\frac{1}{1}\right)}="2",(20,14)*+{\left(\frac{1}{1},\frac{2}{1},\frac{1}{0}\right)}="3",(50,-24)*+{\left(\frac{0}{1},\frac{1}{3},\frac{1}{2}\right)}="4",(50,-8)*+{\left(\frac{1}{2},\frac{2}{3},\frac{1}{1}\right)}="5",(50,8)*+{\left(\frac{1}{1},\frac{3}{2},\frac{2}{1}\right)}="6",(50,24)*+{\left(\frac{2}{1},\frac{3}{1},\frac{1}{0}\right)}="7",(85,-28)*+{\left(\frac{0}{1},\frac{1}{4},\frac{1}{3}\right)\cdots}="8",(85,-20)*+{\left(\frac{1}{3},\frac{2}{5},\frac{1}{2}\right)\cdots}="9",(85,-12)*+{\left(\frac{1}{2},\frac{3}{5},\frac{2}{3}\right)\cdots}="10",(85,-4)*+{\left(\frac{2}{3},\frac{3}{4},\frac{1}{1}\right)\cdots}="11",(85,4)*+{\left(\frac{1}{1},\frac{4}{3},\frac{3}{2}\right)\cdots}="12",(85,12)*+{\left(\frac{3}{2},\frac{5}{3},\frac{2}{1}\right)\cdots}="13",(85,20)*+{\left(\frac{2}{1},\frac{5}{2},\frac{3}{1}\right)\cdots}="14",(85,28)*+{\left(\frac{3}{1},\frac{4}{1},\frac{1}{0}\right)\cdots}="15",\ar@{-}"1";"2"\ar@{-}"1";"3"\ar@{-}"2";"4"\ar@{-}"2";"5"\ar@{-}"3";"6"\ar@{-}"3";"7"\ar@{-}"4";"8"\ar@{-}"4";"9"\ar@{-}"5";"10"\ar@{-}"5";"11"\ar@{-}"6";"12"\ar@{-}"6";"13"\ar@{-}"7";"14"\ar@{-}"7";"15"
\end{xy}
\]
We first prove the basic properties of the Farey tree that will be used later.
\begin{prop}\label{prop:property-farey}\index{Farey tree}\index{Farey triple}
The following hold.
\begin{itemize}
\item[(1)] If $(r,t,s)$ is a Farey triple, then $(r,r\oplus t,t)$ and $(t,t\oplus s,s)$ are also Farey triples. In particular, every vertex of $\mathrm{F}\mathbb T$ is a Farey triple.
\item[(2)] For every reduced fraction $t\in(0,\infty)$, there exists a unique Farey triple $F$ in $\mathrm{F}\mathbb T$ whose second component is $t$.
\item[(3)] For every $(r,t,s)$ in $\mathrm{F}\mathbb T$, the inequalities $r<t<s$ hold.
\end{itemize}
\end{prop}
We first record several lemmas.

\begin{lemm}\label{lem:mediant-det}
Let $x=\frac{a}{b}$ and $y=\frac{c}{d}$ be reduced fractions, and put $x\oplus y=\frac{a+c}{b+d}$. Then
\[
\det(x,x\oplus y)=\det(x,y),\qquad
\det(x\oplus y,y)=\det(x,y).
\]
In particular, if $|\det(x,y)|=1$, then $x\oplus y$ is a reduced fraction.
\end{lemm}

\begin{proof}
A direct computation gives
\[
\det\!\left(\frac{a}{b},\frac{a+c}{b+d}\right)=a(b+d)-b(a+c)=ad-bc=\det\!\left(\frac{a}{b},\frac{c}{d}\right),
\]
and
\[
\det\!\left(\frac{a+c}{b+d},\frac{c}{d}\right)=(a+c)d-(b+d)c=ad-bc=\det\!\left(\frac{a}{b},\frac{c}{d}\right).
\]
If $g=\gcd(a+c,b+d)$, then $g\mid\bigl((a+c)d-(b+d)c\bigr)=ad-bc$. Thus, if $|\det(x,y)|=1$, then $g\mid1$, so $g=1$. Hence $x\oplus y$ is reduced.
\end{proof}

\begin{lemm}\label{lem:mediant-between}\index{Mediant}
Consider the extended nonnegative rationals, consisting of nonnegative reduced fractions together with $\frac{1}{0}$. Suppose that $\frac{a}{b}<\frac{c}{d}$ and $ad-bc=-1$. Then the mediant $\frac{a+c}{b+d}$ satisfies $\frac{a}{b}<\frac{a+c}{b+d}<\frac{c}{d}$. Here $\frac{1}{0}$ is understood to be larger than every nonnegative rational number.
\end{lemm}
\begin{proof}
If $b,d>0$, then
\[
\frac{a+c}{b+d}-\frac{a}{b}=\frac{bc-ad}{b(b+d)}=\frac1{b(b+d)}>0,\qquad\text{and}\qquad
\frac{c}{d}-\frac{a+c}{b+d}=\frac{bc-ad}{d(b+d)}=\frac1{d(b+d)}>0.
\]
If $b=0$ or $d=0$, then one endpoint is $\frac{1}{0}=\infty$, and the remaining inequality follows immediately from the definition.
\end{proof}

\begin{lemm}\label{lem:middle-is-mediant}\index{Mediant}
Let $(r,t,s)=\left(\frac{a}{b},\frac{c}{d},\frac{e}{f}\right)$ be a Farey triple satisfying $r<t<s$. Then $t=r\oplus s=\frac{a+e}{b+f}$.
\end{lemm}

\begin{proof}
Since $r<t<s$ and $|\det(r,s)|=1$, we have $\det(r,s)=af-be=-1$. Put $\mathbf u=(a,b)$ and $\mathbf w=(e,f)$. Then $\det(\mathbf u^T,\mathbf w^T)=-1$, where $\det$ is used in the ordinary matrix sense. Hence $\mathbf u$ and $\mathbf w$ form a basis of $\mathbb Z^2$. Therefore $\mathbf v:=(c,d)$ can be written uniquely as $\mathbf{v}=\alpha \mathbf{u}+\beta \mathbf{w}$. On the other hand, $r<t$ and $|\det(r,t)|=1$ imply $\det(r,t)=ad-bc=-1$, and therefore
\[
-1=\det(\mathbf u^T,\mathbf v^T)=\det\bigl(\mathbf u^T,\alpha \mathbf u^T+\beta \mathbf w^T\bigr)=\beta\,\det(\mathbf u^T,\mathbf w^T)=\beta(-1).
\]
Thus $\beta=1$. Similarly, $t<s$ and $|\det(t,s)|=1$ imply $\det(t,s)=cf-de=-1$, and hence
\[
-1=\det(\mathbf v^T,\mathbf w^T)=\det(\alpha\mathbf u^T+\mathbf w^T,\mathbf w^T)=\alpha\,\det(\mathbf u^T,\mathbf w^T)=\alpha(-1).
\]
Thus $\alpha=1$. Hence $\mathbf v=\mathbf u+\mathbf w$, that is, $(c,d)=(a+e,b+f)$. Therefore $t=\frac{c}{d}=\frac{a+e}{b+f}=r\oplus s$.
\end{proof}

\begin{proof}
We prove (1). Let $(r,t,s)$ be a Farey triple. By Lemma~\ref{lem:mediant-det},
\[
|\det(r,r\oplus t)|=|\det(r,t)|=1,\qquad |\det(r\oplus t,t)|=|\det(r,t)|=1,
\]
and $r\oplus t$ is reduced. Since $|\det(t,r)|=|\det(r,t)|=1$, it follows that $(r,r\oplus t,t)$ is a Farey triple. The same argument shows that $(t,t\oplus s,s)$ is a Farey triple. Thus the property of being a Farey triple is preserved when we take children. The initial vertex $(0/1,1/1,1/0)$ is a Farey triple, and therefore every vertex of $\mathrm{F}\mathbb T$ is a Farey triple.

We next prove (3). At the initial vertex, we have $\frac{0}{1}<\frac{1}{1}<\frac{1}{0}(=\infty)$. Suppose that $r<t<s$ at a vertex $(r,t,s)$. By Lemma~\ref{lem:mediant-between}, we have $r<r\oplus t<t$ and $t<t\oplus s<s$. Thus the same inequalities hold for the left child $(r,r\oplus t,t)$ and the right child $(t,t\oplus s,s)$. Hence (3) follows by induction on the depth.

We prove (2). Take a reduced fraction $t\in(0,\infty)$. Put $L_0=\frac{0}{1}$ and $R_0=\frac{1}{0}$. For $i\in \mathbb Z_{\geq 0}$, compute inductively $M_i:=L_i\oplus R_i$. If $t<M_i$, put $L_{i+1}:=L_i$ and $R_{i+1}:=M_i$; if $t>M_i$, put $L_{i+1}:=M_i$ and $R_{i+1}:=R_i$. These two cases correspond to taking the left child and the right child, respectively. If $t=M_i$, stop the procedure. Writing $L_i=\frac{a_i}{b_i}$ and $R_i=\frac{c_i}{d_i}$ at each step, we have
\[
M_i=\frac{a_i+c_i}{b_i+d_i},
\]
and Lemma~\ref{lem:mediant-det} shows that $|\det(L_i,R_i)|=1$ is always preserved. In particular, $L_i$ and $R_i$ are always reduced. Lemma~\ref{lem:mediant-between} also gives $L_i<t<R_i$ for all $i$.
We first show that the procedure always stops after finitely many steps. Once both endpoints have positive denominators and the procedure continues, the denominator of the next mediant $M_{i+1}$ is strictly larger than the denominator $b_i+d_i$ of $M_i$. The only time denominators may fail to increase monotonically is while one endpoint is $\frac{1}{0}$. In that case the procedure simply moves through the integer part: if $t$ is an integer, it stops there, and otherwise after finitely many steps it enters an interval between two consecutive integers containing $t$. From then on, both endpoint denominators are positive.
Now write $t=\frac pq$ in lowest terms. If both endpoints have positive denominators and
$L_i=\frac{a_i}{b_i}<\frac pq<R_i=\frac{c_i}{d_i}$, then
\[
 b_ip-a_iq\ge1,
 \qquad
 c_iq-pd_i\ge1,
\]
and from $b_ic_i-a_id_i=1$ we obtain
\[
q=d_i(b_ip-a_iq)+b_i(c_iq-pd_i)\ge b_i+d_i.
\]
Therefore, after both endpoint denominators become positive, the desired mediant must be reached before a mediant denominator larger than $q$ would be required. Hence the procedure cannot continue indefinitely. Thus $t=M_n$ for some finite $n$. The vertex $(L_n,t,R_n)$ is obtained by repeatedly taking children, so it is a vertex of $\mathrm{F}\mathbb T$. This proves existence of a vertex whose second component is $t$.
Finally, we prove uniqueness. If a vertex $(r,t,s)$ with second component $t$ is given, Lemma~\ref{lem:middle-is-mediant} implies that $t=r\oplus s$, and (3) gives $r<t<s$. Thus the sequence of left and right choices from the initial vertex to this vertex is determined uniquely by $t$. If two distinct vertices had the same second component $t$, then there would be two different paths from the initial vertex, a contradiction. Hence such a vertex is unique.
\end{proof}
The Farey tree enumerates reduced fractions and is closely related to Farey sequences and the Stern--Brocot tree.

Comparing the generation rules of the Farey tree and the GM tree, one sees that the rearrangements of components are the same. This suggests a correspondence from positive reduced fractions to GM numbers.
\begin{defi}\label{def:gm-fraction-label}\index{Farey tree}\index{Fraction label}\index{Generalized Markov pair}\index{Position label}
Let
\[
f\colon\mathrm F\mathbb T\longrightarrow\mathrm M\mathbb T(k_1,k_2,k_3,\sigma)
\]
be the unique ordered complete binary tree isomorphism sending the root to the initial vertex. If $v=(r,t,s)$ and $f(v)=((a,h),(b,i),(c,j))$, label the first, second, and third \emph{occurrences} of GM pairs at $f(v)$ by $r,t,s$, respectively.

This labeling is consistent between generations. A child in the Farey tree inherits two entries of its parent in the first and third positions and has their mediant in the middle. The corresponding GM-tree child inherits the same two GM pairs in its first and third positions and has a new GM pair in the middle. By induction on depth, all occurrences with the same label $t$ have the same numerical value $(m_t,i_t)$. We call $t$ the \emph{fraction label} of the occurrence.
\end{defi}
Strictly speaking, $m_t,i_t$ also depend on $k_1,k_2,k_3,\sigma$; we suppress these parameters when a single tree is fixed. Occurrences obtained by inheriting a GM pair from a parent to its children share its fraction label. The injectivity of the numerical maps $t\mapsto(m_t,i_t)$ and $t\mapsto m_t$ will be examined later.

We close this section with a simple but important corollary.
\begin{coro}\label{cor:dual-remark}\index{Fraction label!reciprocal labels}\index{Fraction label}
Here we agree that $1/0=\infty$ and $1/\infty=0$. Let $t$ be a reduced fraction in $[0,\infty]$. Let $(m_t,i_t)$ be the $(k_1,k_2,k_3)$-GM pair with fraction label $t$ in $\mathrm M\mathbb T(k_1,k_2,k_3,\sigma)$, and let $(m^{\ast}_{\frac{1}{t}},i^{\ast}_{\frac{1}{t}})$ be the $(k_1,k_2,k_3)$-GM pair with fraction label $\frac{1}{t}$ in $\mathrm M\mathbb T(k_1,k_2,k_3,\sigma^{\ast})$, where $\sigma^{\ast}$ is the one appearing in Theorem~\ref{thm:all-solution}. Then
$(m_t,i_t)=(m^{\ast}_{\frac{1}{t}},i^{\ast}_{\frac{1}{t}})$.
\end{coro}
\begin{proof}
By Theorem~\ref{thm:all-solution}, it suffices to show that the fraction located at the mirror-symmetric position to $t$ in the Farey tree is $1/t$. The endpoint cases $t=0$ and $t=\infty$ are immediate from the root. For positive finite $t$, this follows by induction on the distance from the root to the unique Farey triple whose middle component is $t$. The root case is $t=1$. Suppose that a vertex $(r,t,s)$ has mirror-symmetric vertex $(s^{-1},t^{-1},r^{-1})$, with the convention $0^{-1}=\infty$ and $\infty^{-1}=0$. The left child of $(r,t,s)$ has middle label $r\oplus t$, while the right child of the mirror-symmetric vertex has middle label $t^{-1}\oplus r^{-1}=(r\oplus t)^{-1}$. The right-child case is the same.
\end{proof}

\section{Characteristic Numbers}\label{sec:characteristic-number}
We finish this chapter by defining numbers called characteristic numbers. We first prove the following theorem.
\begin{theorem}\label{thm:exist-unique-ut}
Fix one tree $\mathrm{M}\mathbb T(k_1,k_2,k_3,\sigma)$, and let $(m_r,m_t,m_s)$ be the GM triple corresponding, via fraction labels, to a vertex $(r,t,s)\in\mathrm F\mathbb T$ of the Farey tree. In particular, $0\leq r<t<s\leq\infty$ and $0<t<\infty$. Then there exists a unique integer $u$ satisfying
\begin{align}\label{eq:characteristic-def}
\begin{cases}
m_r u\equiv m_s \pmod{m_t},\\
0<u<m_t.
\end{cases}
\end{align}
\end{theorem}
To prove this theorem, we use the following consequence of the Euclidean algorithm. The proof here uses facts about finite regular continued fractions.
\begin{lemm}\label{lem:ax+by}
Let $x,y\in \mathbb Z$, and suppose that at least one of them is nonzero. Then there exist $a,b\in \mathbb Z$ such that $ax+by=\gcd(x,y)$.
\end{lemm}
\begin{proof}
First assume that $x,y\geq0$. If $y=0$, then $x>0$ by assumption, and we may take $a=1,b=0$. If $x=0$, then $y>0$, and we may take $a=0,b=1$. We therefore assume $x,y>0$. Put $d=\gcd(x,y)$ and write $x=dx'$, $y=dy'$. Let
$\frac{x'}{y'}=[a_0;a_1,\dots,a_n]$ be the finite regular continued-fraction expansion. With the convention $p_{-1}=1,q_{-1}=0$ when $n=0$, put $\frac{p_i}{q_i}=[a_0;a_1,\dots,a_i]$. By Lemma~\ref{lem:det}, $x'q_{n-1}-y'p_{n-1}=(-1)^{n+1}$. If $(-1)^{n+1}=1$, take $a=q_{n-1}$ and $b=-p_{n-1}$; if $(-1)^{n+1}=-1$, take $a=-q_{n-1}$ and $b=p_{n-1}$. Then $ax'+by'=1$. Multiplying both sides by $d$ gives $ax+by=d$, proving the assertion in this case. If $x<0$ or $y<0$, first apply the preceding argument to $|x|$ and $|y|$, and then replace $a$ by $-a$ if $|x|=-x$, and replace $b$ by $-b$ if $|y|=-y$. The same conclusion follows.
\end{proof}
\begin{proof}
We first prove existence. By Corollary~\ref{relatively-prime}, we have $\gcd(m_r,m_t)=1$. Hence Lemma~\ref{lem:ax+by} gives integers $a,b$ such that $m_r a+m_t b=1$.

Thus $a$ is an inverse of $m_r$ modulo $m_t$. Put $u_0:=am_s$. Then
\[
m_r u_0=m_r a m_s\equiv 1\cdot m_s\equiv m_s\pmod{m_t},
\]
so $u_0$ is a solution of the congruence $m_r x\equiv m_s\pmod{m_t}$. Choose the integer $u$ satisfying $0\le u<m_t$ and
\[
u\equiv u_0\pmod{m_t}.
\]
This $u$ also solves the same congruence. If $u=0$, then $m_s\equiv0\pmod{m_t}$, that is, $m_t\mid m_s$. Since $\gcd(m_t,m_s)=1$, this implies $m_t=1$. However, for a label $t\in(0,\infty)$, the corresponding entry is a middle entry of a GM-tree vertex; at the initial vertex it is $k_{\sigma(2)}+2\ge2$, and the same lower bound is preserved along the tree. Hence $m_t\ge2$, a contradiction. Therefore $u\neq0$, and we obtain a solution satisfying $0<u<m_t$.

We next prove uniqueness. Suppose that $u$ and $u'$ both satisfy
\[
m_r u\equiv m_s\pmod{m_t},\qquad 0<u<m_t,\qquad\text{and}\qquad
m_r u'\equiv m_s\pmod{m_t},\qquad 0<u'<m_t.
\]
Subtracting the two congruences gives
\[
m_r(u-u')\equiv0\pmod{m_t},
\]
that is, $m_t\mid m_r(u-u')$. Since $\gcd(m_r,m_t)=1$, we obtain $m_t\mid(u-u')$. Moreover, $0<u,u'<m_t$ implies
\[
-(m_t-1)\le u-u'\le m_t-1.
\]
The only multiple of $m_t$ in this interval is $0$, so $u-u'=0$. Hence $u=u'$. Therefore the required integer $u_t$ exists and is unique.
\end{proof}
\begin{defi}\label{def:characteristic-number}\index{Characteristic number}
Let $t\in\mathbb Q\cap(0,\infty)$ be an interior fraction label, let $(r,t,s)$ be the unique Farey-tree vertex with middle component $t$, and let $(m_r,m_t,m_s)$ be the corresponding GM triple. The unique integer satisfying \eqref{eq:characteristic-def} is denoted by $u_t$ and called the \emph{characteristic number} of the label $t$.

To state later formulas uniformly, introduce the auxiliary endpoint values
\[
u_{\frac01}:=-k_{\sigma(1)},\qquad u_{\frac10}:=1.
\]
These are auxiliary values outside the domain $\mathbb Q\cap(0,\infty)$ of characteristic numbers. The inequalities $0<u_t<m_t$ are asserted only for interior labels.
\end{defi}
It may seem strange to use notation depending only on $t$ for a characteristic number defined from a GM triple. However, by Proposition~\ref{prop:property-farey} (2), the Farey triple $(r,t,s)$ is uniquely determined by $t$. Thus $u_t$ depends only on $t$, and the notation is justified. We call this the \emph{fraction labeling of characteristic numbers}. We finish the section by recording an important property of characteristic numbers.
\begin{prop}\label{prop:t-1/t-relation-gen}\index{Characteristic number}
For any reduced fraction $t\in[0,1]\cap \mathbb Q$, let $u_t$ be the characteristic number (or auxiliary endpoint value when $t=0$) with fraction label $t$ in $\mathrm M\mathbb T(k_1,k_2,k_3,\sigma)$, and put $k_t:=k_{i_t}$. Let $u_{\frac{1}{t}}^{\ast}$ be the characteristic number (or auxiliary endpoint value when $1/t=\infty$) with fraction label $\frac{1}{t}$ in $\mathrm M\mathbb T(k_1,k_2,k_3,\sigma^{\ast})$, where $\sigma^{\ast}$ is the one appearing in Theorem~\ref{thm:all-solution}. Then
\[
u^{\ast}_{\frac{1}{t}}=m_t-u_t-k_t.
\]
\end{prop}
The proof will be given in Section~\ref{section:relation-of-characteristic}, using results established there.

\chapter[Fence Posets and Generalized Markov Distance]{Fence Posets and Generalized Markov Distance}\label{chap:fence-posets-gm-distance}
\tikzset{
  chsix triangulation/.style={draw=black!30,line width=.36pt,
    line cap=butt,line join=miter},
  chsix distinguished edge/.style={draw=black!70,line width=.78pt},
  chsix red curve/.style={draw=annred,line width=1.10pt,
    line cap=round,line join=round},
  chsix blue curve/.style={draw=blue!72!black,line width=1.10pt,
    line cap=round,line join=round},
  chsix reference segment/.style={draw=black!30,densely dotted,
    line width=.42pt},
  chsix transition arrow/.style={draw=black!65,line width=.55pt,
    -{Stealth[length=1.50mm,width=1.05mm]}},
  chsix point/.style={circle,fill=black,inner sep=1.05pt},
  chsix red point/.style={circle,fill=annred,inner sep=1.10pt},
  chsix fence edge/.style={draw=black!80,line width=.85pt},
  chsix fence vertex/.style={font=\small,inner sep=1.6pt},
  chsix curve label/.style={font=\scriptsize,fill=white,inner sep=.7pt},
  chsix forward at/.style={postaction={decorate},decoration={markings,
    mark=at position #1 with {\arrow{Stealth[length=1.65mm,width=1.10mm]}}}},
  chsix forward/.style={chsix forward at=.55},
  chsix backward/.style={postaction={decorate},decoration={markings,
    mark=at position .55 with {\arrowreversed{Stealth[length=1.65mm,width=1.10mm]}}}},
  chsix panel title/.style={font=\small}
}
In \hyperref[chap:generalized-markov-equations-numbers]{Chapter~\ref*{chap:generalized-markov-equations-numbers}}, we introduced generalized Markov equations, generalized Markov trees, fraction labels, and characteristic numbers. To describe the relations among generalized Markov numbers, fraction labels, and characteristic numbers in combinatorial and geometric terms, we now introduce fence posets and generalized Markov distance. Our aim is to make this theory accessible through the structure of posets and the intersections of curves.

We first associate a fence poset with a finite integer sequence and show that its number of order ideals is closely related to continued fractions and continued-fraction matrices. This reinterprets the computations in \hyperref[chap:continued-fraction]{Chapter~\ref*{chap:continued-fraction}} as the combinatorics of posets. We then assign a generalized Markov length to curves and generalized arcs in the plane and define generalized Markov distance by minimizing this length. Thus the generalized Markov numbers introduced in the preceding chapter can also be understood as invariants of geometric objects.

The background of this chapter lies in cluster algebra theory. Snake graphs, introduced in \cite{msw,msw2} to describe cluster algebra generators combinatorially, are also used in combinatorial descriptions of Markov numbers. We use the equivalent formulation in terms of fence posets.

Our treatment of generalized Markov distance is based on \cite{llrs,banaian}. These works treat, respectively, the cases $k_1=k_2=k_3=0$ and $k_1=k_2=k_3$; here we extend the framework to arbitrary $(k_1,k_2,k_3)$.
\section{Order Ideals of Fence Posets and Continued Fractions}
\begin{defi}
Let $(P,\preceq)$ be a poset. For $x,y\in P$, write $x\lessdot y$ if
\[
x\prec y \quad\text{and there is no } z\in P \text{ with } x\prec z\prec y.
\]
\index{Cover relation}The relation $\lessdot$ is called the \textbf{cover relation} of $(P,\preceq)$. \index{Hasse diagram}The \textbf{Hasse diagram} of $(P,\preceq)$ is the graph satisfying the following conditions.
\begin{itemize}
  \item Its vertex set is $P$.
  \item There is an edge, without multiplicity, between distinct vertices $x,y\in P$ if and only if $x\lessdot y$ or $y\lessdot x$.
\end{itemize}
Each edge is drawn with the larger element above the smaller one in the cover relation, and its orientation is omitted.
\end{defi}
We consider posets whose undirected Hasse diagrams are finite paths, as specified in the following definition.
\begin{defi}\label{def:fence-poset}\index{Path-order labeling}\index{Hasse diagram}\index{Fence poset}
A finite poset $(P,\preceq)$ is a \textbf{fence poset} if its undirected Hasse diagram is a finite path. A poset with one vertex is allowed, and the empty poset is regarded as a degenerate fence poset. In the nonempty case, there are an integer $m\geq1$ and a labeling
\[
P=\{P(1),P(2),\dots,P(m)\}
\]
such that the edges of its Hasse diagram are exactly
\[
\bigl\{\{P(i),P(i+1)\}\bigm|1\leq i<m\bigr\}
\]
This labeling is unique up to reversal. We fix one of the two choices and call it a \textbf{path-order labeling} of $P$.
\end{defi}
Henceforth every nonempty fence poset is equipped with a path-order labeling $P(1),\dots,P(m)$, with these vertices placed from left to right in diagrams. For simplicity, we may write $i$ for $P(i)$ and represent $P$ by an order on $\{1,\dots,m\}$. Given a finite sequence $S=(a_0,\dots,a_n)$ of positive integers, put $s_k:=\sum_{i=0}^{k}a_i\ (k=0,\dots,n)$ and define the poset $P_S:=\bigl(\{1,2,\dots,s_n-1\},\preceq\bigr)$ as follows. For each \(x\in P_S\), if $s_{k-1}\le x<s_k\ (s_{-1}:=0)$, set
\[
\varepsilon(x):=(-1)^{k}
\]
For adjacent vertices $x,x+1\ (1\leq x<s_n-1)$, define the order $\prec$ on $P_S$ by the following cover relations:
\[
\begin{cases}
x\lessdot x+1 & \text{if } \varepsilon(x)=1,\\
x+1\lessdot x & \text{if } \varepsilon(x)=-1.
\end{cases}
\]The last vertex is labeled $s_n-1$, not $s_n$. If $n\geq1$, the numbers of consecutive cover edges with the same orientation are $a_0-1,a_1,\ldots,a_{n-1},a_n-1$. If $n=0$, then $P_{(a_0)}$ is a chain with $a_0-1$ vertices and $\max\{a_0-2,0\}$ edges.
\begin{exam}
For $S=(3,2,1,2)$, the Hasse diagram of $P_S$ is as follows.
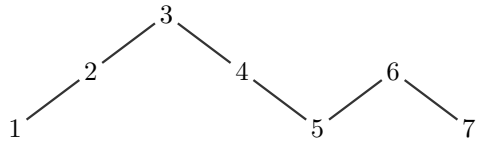
\begin{center}
\begin{tikzpicture}[scale=1.0,baseline=(current bounding box.center),
  every node/.style={chsix fence vertex}]
  \tikzset{
    e/.style={chsix fence edge}
  }

  \node[] (v0) at (0,0) {$1$};
  \node[] (v2) at (4,0) {$5$};
  \node[] (v4) at (6,0) {$7$};

  \node[] (u1) at (1,0.75) {$2$};
  \node[] (w2) at (3,0.75) {$4$};
  \node[] (t1) at (5,0.75) {$6$};

  \node[] (v1) at (2,1.5) {$3$};

  \draw[e] (v0)--(u1);
  \draw[e] (u1)--(v1);

  \draw[e] (v1)--(w2);
  \draw[e] (w2)--(v2);

  \draw[e] (t1)--(v4);
  \draw[e] (v2)--(t1);
\end{tikzpicture}
\end{center}
\end{exam}
\begin{defi}\label{def:order-ideal}\index{Order ideal}
Let $(P,\le)$ be a poset. A subset $I\subset P$ is an \emph{order ideal} if, for every $x\in I$ and $y\in P$, the inequality $y\le x$ implies $y\in I$. Thus $I$ is downward closed. Write $\mathcal J(P)$ for the set of all order ideals of $P$.
\end{defi}
Unlike ideals in a ring, the empty set is also regarded as an order ideal.
\begin{exam}
The following are all $27$ order ideals of $P_{(3,2,1,2)}$.
\begin{align*}
&\emptyset, \{1\},\{5\},\{7\},\{1,2\},\{1,5\},\{1,7\},\{4,5\},\{5,7\},\{1,2,5\},\{1,2,7\},\{1,4,5\},\{1,5,7\},\{4,5,7\},\\
&\{5,6,7\},\{1,2,4,5\},\{1,2,5,7\},\{1,4,5,7\},\{1,5,6,7\},\{4,5,6,7\},\{1,2,3,4,5\},\{1,2,4,5,7\},\\&\{1,2,5,6,7\},\{1,4,5,6,7\},\{1,2,4,5,6,7\},\{1,2,3,4,5,7\},\{1,2,3,4,5,6,7\}.
\end{align*}
\end{exam}
The key quantity for us is the number of order ideals of a fence poset:
\[
N(P):=\#\mathcal J(P)
\]
For a sequence of positive integers, put $N(a_0,\ldots,a_n):=N(P_{(a_0,\ldots,a_n)})$ and $N():=1$. If $P^\ast$ is the dual of a finite poset $P$, the complement map $I\mapsto P\setminus I$ is a bijection from $\mathcal J(P)$ to $\mathcal J(P^\ast)$, so
\[
N(P^\ast)=N(P)
\]
Although listing the order ideals by hand is cumbersome, their number is easily computed by the following theorem.
\begin{theorem}\label{thm:continued-fraction-order-ideal}\index{Fence poset}\index{Order ideal}
For a finite sequence $(a_0,a_1,\dots,a_n)$ of positive integers, write $N(a_0,\dots,a_n)$ for $N(P_{(a_0,\dots,a_n)})$. If
\[\frac{p_n}{q_n}:=[a_0;a_1,\dots,a_n]\]is in lowest terms, then
\[N(a_0,\dotsm,a_n)=p_n,\quad N(a_1,\dotsm,a_n)=q_n\]
Here $N()=1$.
\end{theorem}
\begin{proof}
It suffices to show that the numerator and denominator of $[a_0;a_1,\dots,a_n]$ in lowest terms are $N(a_0,a_1,\dots,a_n)$ and $N(a_1,\dots,a_n)$. For $n=0$, we have $p_0=a_0,q_0=1$, whereas $P_{(a_0)}$ is the chain $1\prec2\prec\cdots\prec a_0-1$, whose order ideals are
\[\{\emptyset,\{1\},\{1,2\},\dots,\{1,2,\dots,a_0-1\}\}.\]
There are $a_0$ of them, so $N(a_0)=a_0$; also $N()=1$ by definition. This proves the claim for $n=0$. For $n=1$, we have $p_1=a_0a_1+1,q_1=a_1$. The Hasse diagram of $P_{(a_0,a_1)}$ is
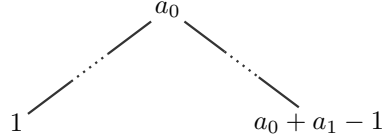
\begin{center}
\begin{tikzpicture}[scale=1.0,baseline=(current bounding box.center),
  every node/.style={chsix fence vertex}]
  \tikzset{e/.style={chsix fence edge}}

  \node (v0) at (0,0) {$1$};
  \node (v2) at (4,0) {$a_0+a_1-1$};
  \node (v1) at (2,1.5) {$a_0$};

  \path (v0) -- (v1)
    coordinate[pos=0.35] (m01a)
    coordinate[pos=0.65] (m01b);
  \draw[e] (v0)--(m01a);
  \draw[e,dotted] (m01a)--(m01b);
  \draw[e] (m01b)--(v1);

  \path (v1) -- (v2)
    coordinate[pos=0.35] (m12a)
    coordinate[pos=0.65] (m12b);
  \draw[e] (v1)--(m12a);
  \draw[e,dotted] (m12a)--(m12b);
  \draw[e] (m12b)--(v2);
\end{tikzpicture}
\end{center}
Every order ideal other than the whole poset $P_{(a_0,a_1)}$ is therefore the disjoint union of an order ideal $I_0$ contained in $\{1,2,\dots,a_0-1\}$ and an order ideal $I_1$ contained in $\{a_0+1,a_0+2,\dots,a_0+a_1-1\}$. There are $a_0$ choices for $I_0$ and $a_1$ for $I_1$, giving $a_0a_1+1$ order ideals in total. Thus $N(a_0,a_1)=a_0a_1+1$, and $N(a_1)=a_1$, proving the claim for $n=1$.
For $n>1$, proceed by induction, assuming the assertion through $n-1$. The recurrences \eqref{eq:p}, \eqref{eq:q} for $p_n,q_n$ show that it suffices to prove
\begin{align}
N(a_0,a_1,\dots,a_n)&=a_nN(a_0,\dots,a_{n-1})+N(a_0,\dots,a_{n-2})\label{eq:N=aN+N}\\
N(a_1,\dots,a_n)&=a_nN(a_1,\dots,a_{n-1})+N(a_1,\dots,a_{n-2})\nonumber
\end{align}
The second identity follows by the same argument as the first, so we prove only the first.

Suppose first that $n$ is odd. Then $s_{n-1}$ is maximal. If $a_n\geq2$, the vertex $s_n-1$ is also minimal, whereas if $a_n=1$, then $s_n-1=s_{n-1}$. Divide the order ideals according to whether they contain $s_{n-1}$.

An order ideal containing $s_{n-1}$ must contain $\{s_{n-2},s_{n-2}+1,\dots,s_n-1\}$ by downward closedness. It is therefore the disjoint union of this set and an order ideal contained in $\{1,\dots,s_{n-2}-1\}$. By induction, there are $N(a_0,\dots,a_{n-2})$ choices for the latter.

On the other hand, the order ideals not containing $s_{n-1}$ are naturally in bijection with the order ideals of the poset obtained by deleting $s_{n-1}$ from $P_{(a_0,\dots,a_n)}$. This poset is
\[P_{(a_0,a_1,\dots,a_{n-1})}\sqcup P^{\ast}_{(a_n)}\]
where $P^\ast$ denotes the poset with the order reversed. Its number of order ideals is
\[N(a_0,\dots,a_{n-1})N(a_n)=a_nN(a_0,\dots,a_{n-1}).\]

We have thus proved, for odd $n$, that
\[N(a_0,a_1,\dots,a_n)=a_nN(a_0,\dots,a_{n-1})+N(a_0,\dots,a_{n-2})\]
If $n$ is even, then $s_{n-1}$ is minimal. An order ideal not containing it contains none of the vertices $s_{n-2},\dots,s_{n-1}-1,s_{n-1}+1,\dots,s_n-1$ above it, and is naturally identified with an order ideal of $P_{(a_0,\dots,a_{n-2})}$. There are $N(a_0,\dots,a_{n-2})$ such ideals. An order ideal containing $s_{n-1}$, on the other hand, is uniquely determined by independent choices of an order ideal in the left poset $P_{(a_0,\dots,a_{n-1})}$ and an initial segment of the right chain $s_{n-1}+1\prec\cdots\prec s_n-1$. There are $a_n$ choices for the latter, including the empty initial segment, and thus $a_nN(a_0,\dots,a_{n-1})$ such ideals. The same recurrence follows for even $n$.
\end{proof}
Thus the number of order ideals of $P_{(a_0,\dots,a_n)}$ is the numerator of $[a_0;a_1,\dots,a_n]$ in lowest terms. The following corollary is immediate.
\begin{coro}\label{cor:matrix-combinatorial}\index{Order ideal}
Let $n\geq0$ and let $(a_0,\dots,a_n)$ be a finite sequence of positive integers. For $n\geq1$, we have
\[F_{(a_0,a_1,\dots,a_n)}=\begin{bmatrix}
    a_0&1\\1&0
\end{bmatrix}\begin{bmatrix}
    a_1&1\\1&0
\end{bmatrix}\cdots\begin{bmatrix}
    a_n&1\\1&0
\end{bmatrix}=\begin{bmatrix}
    N(a_0,\dots,a_n)&N(a_0,\dots,a_{n-1})\\N(a_1,\dots,a_n)&N(a_1,\dots,a_{n-1})
\end{bmatrix}\]
For $n=0$, interpret the right-hand side as
\[
\begin{bmatrix}N(a_0)&N()\\N()&0\end{bmatrix}.
\]
\end{coro}
\begin{proof}
This follows from Theorems~\ref{thm:continued-fraction-order-ideal} and~\ref{thm:continued-fraction-matrix}.
\end{proof}

This correspondence computes the number of order ideals of a fence poset as a continued-fraction numerator. It also yields the following proposition.
\begin{prop}\label{prop:reverse-property}
Let $(a_0,a_1,\dots,a_n)$ be a finite sequence of positive integers. Then
\[
N(a_0,a_1,\dots,a_{n-1},a_n)=N(a_n,a_{n-1},\dots,a_1,a_0).
\]
\end{prop}
\begin{proof}
For $n=0$, the assertion is immediate; assume $n\geq1$.
Each matrix $F_a=\begin{bsmallmatrix}a&1\\1&0\end{bsmallmatrix}$ is symmetric, so
\[
F_{a_n}F_{a_{n-1}}\cdots F_{a_0}=(F_{a_0}F_{a_1}\cdots F_{a_n})^T
\]
Applying Corollary~\ref{cor:matrix-combinatorial} to the right-hand side gives
\[
F_{a_n}F_{a_{n-1}}\cdots F_{a_0}
=
\begin{bmatrix}
    N(a_0,\dots,a_n)&N(a_1,\dots,a_n)\\
    N(a_0,\dots,a_{n-1})&N(a_1,\dots,a_{n-1})
\end{bmatrix}
\]
Applying the same corollary to the left-hand side shows that its $(1,1)$ entry is $N(a_n,a_{n-1},\dots,a_0)$. Comparing these entries gives $N(a_0,a_1,\dots,a_n)=N(a_n,a_{n-1},\dots,a_0)$.
\end{proof}
\section{Generalized Markov Length of Curves and Generalized Markov Distance}\label{sec:GM-distance}

Fix $(k_1,k_2,k_3)\in\mathbb Z_{\geq0}^3$ and $\sigma\in\mathfrak S_3$. Let $\mathcal V$ consist of the following points of $\mathbb R^2$.
\begin{enumerate}
\item All lattice points $(a,b)\in\mathbb Z^2$.
\item If $k_{\sigma(1)}\neq0$, all horizontal-edge midpoints
$(a+\frac12,b)$($(a,b)\in\mathbb Z^2$).
\item If $k_{\sigma(2)}\neq0$, all diagonal-edge midpoints
$(a+\frac12,b+\frac12)$($(a,b)\in\mathbb Z^2$).
\item If $k_{\sigma(3)}\neq0$, all vertical-edge midpoints
$(a,b+\frac12)$($(a,b)\in\mathbb Z^2$).
\end{enumerate}
In this section, unless otherwise specified, a ``point'' means an element of $\mathcal V$.

Call all lines of slopes $0,-1,\infty$ through lattice points \emph{triangulation lines}, and write $\widetilde{\mathbb R}^{\,2}$ for the plane equipped with the triangulation they determine and the point set $\mathcal V$. \index{Triangulation edge}We distinguish a \emph{triangulation edge}, meaning a whole side of a triangle, from an \emph{edge} of $\widetilde{\mathbb R}^{\,2}$, meaning the segment joining two consecutive points on a triangulation line. If the midpoint of a triangulation edge belongs to $\mathcal V$, that whole edge is divided into two edges. For an edge $e$ of $\widetilde{\mathbb R}^{\,2}$, write $\operatorname{side}(e)$ for the unique triangulation edge containing it and $\mu(e)$ for the midpoint of $\operatorname{side}(e)$. Thus $\mu(e)$ is the midpoint of the whole triangulation edge and, when subdivision occurs, is an endpoint of the subdivided edge $e$.

A \emph{curve segment} in $\widetilde{\mathbb R}^{\,2}$ is a curve whose endpoints belong to $\mathcal V$. Throughout this section, every curve segment is assumed to satisfy the following conditions.
\begin{itemize}
\item Its interior avoids $\mathcal V$.
\item Whenever it meets an edge of $\widetilde{\mathbb R}^{\,2}$ away from its own endpoints, the intersection avoids the endpoints of the edge, and the curve passes from one side of the edge to the other without tangency. We call this a \emph{transverse intersection}.
\item It passes through triangles only finitely many times.
\item It has finitely many self-intersections.
\end{itemize}
Parametrize each curve segment by a continuous map $[0,1]\to\mathbb R^2$.

A \emph{lattice point} means an element of $\mathbb Z^2$; the nonlattice points of $\mathcal V$ are midpoints of triangulation edges.

From now on, the endpoints of curve segments are assumed to be lattice points.

\begin{defi}\label{def:passage-occurrence}\index{Triangle-passage occurrence}\index{Edge-crossing occurrence}
Let $\gamma$ be a curve segment. For a triangle $\Delta$ of $\widetilde{\mathbb R}^{\,2}$, if
\[
I\subset\gamma^{-1}(\operatorname{int}\Delta)
\]
is a connected component, call $(\Delta,I)$ a \emph{triangle-passage occurrence} of $\gamma$. For an edge $e$ of $\widetilde{\mathbb R}^{\,2}$ and $t\in(0,1)$ with $\gamma(t)\in e^\circ$, call $(e,t)$ an \emph{edge-crossing occurrence} of $\gamma$. Here $e^\circ$ is the edge with both endpoints removed. Order occurrences by increasing parameter of $\gamma$.

Different components $I$ for the same triangle $\Delta$, and different crossing times $t$ for the same edge $e$, are distinct occurrences. Unless geometric triangles or edges themselves are explicitly meant, passages of $\gamma$ and sign assignments are henceforth interpreted occurrence by occurrence.

For an oriented infinite line $l$ that is not a triangulation line, use the same terminology after choosing an orientation-preserving parametrization $\lambda\colon\mathbb R\to l$ and replacing both $[0,1]$ and $(0,1)$ by $\mathbb R$. There are finitely many occurrences in every bounded region, and their order and neighboring curve portions do not depend on the orientation-preserving parametrization.
\end{defi}

\begin{defi}\label{def:general-arc}\index{Triangulation edge}\index{Generalized arc}
A curve segment $\gamma$ in $\widetilde{\mathbb R}^{\,2}$ with lattice endpoints is a \emph{generalized arc} if it satisfies all of the following conditions.
\begin{itemize}
\item For each triangle-passage occurrence $(\Delta,I)$, let $\overline I$ denote the closure of $I$ in the parameter interval $[0,1]$. The restriction $\gamma|_{\overline I}$ is injective. If $\overline I$ contains neither endpoint $0,1$, its image joins two distinct points of $\partial\Delta$. If $\overline I$ contains exactly one of $0,1$, its image joins the corresponding lattice endpoint to an intersection on the side opposite that vertex. If $\overline I$ contains both $0,1$, its image joins the two lattice endpoints of $\gamma$.
\item Every self-intersection of $\gamma$ is a transverse intersection at one point of the curve portions belonging to two distinct triangle-passage occurrences.
\item For any two consecutive edge-crossing occurrences $(e,t),(e',t')$ in parameter order,
\[
\operatorname{side}(e)\neq\operatorname{side}(e')
\]
holds. Thus the same whole triangulation edge cannot be crossed twice consecutively.
\end{itemize}
\end{defi}

To associate a finite sign sequence with an oriented generalized arc, assign signs in $\{+,-\}$ to triangle-passage and edge-crossing occurrences by the following rules. As in Definition~\ref{def:passage-occurrence}, the assignment is made independently for each occurrence, even for the same triangle or edge.

For $q\in\{+,-\}$ and an integer $a\geq0$, let $q^a$ denote $a$ repetitions of $q$, with $q^0$ the empty word.

\begin{defi}\label{def:sign-assignment}
Let $\gamma$ be an oriented generalized arc in $\widetilde{\mathbb R}^{\,2}$.
\begin{enumerate}
\item[(1)] \index{Sign assignment rules}For each triangle-passage occurrence $(\Delta,I)$ such that $\overline I$ contains an endpoint, choose one sign, either $-$ or $+$. This is the \emph{endpoint rule} for $\gamma$ (Figure~\ref{fig:minus-righttriangles-endpoint}).
\begin{figure}[ht]
\centering
\begin{tikzpicture}[x=.90cm,y=.90cm]
\foreach \offset/\start/\finish/\endpoint/\direction/\otherdirection in {
  0/{(0,0)}/{(0.58,0.42)}/{(0,0)}/forward/backward,
  2.7/{(0,1)}/{(0.38,0)}/{(0,1)}/forward/backward,
  5.4/{(0,0.38)}/{(1,0)}/{(1,0)}/backward/forward,
  8.1/{(0,0)}/{(0.58,0.42)}/{(0,0)}/backward/forward,
  10.8/{(0,1)}/{(0.38,0)}/{(0,1)}/backward/forward,
  13.5/{(0,0.38)}/{(1,0)}/{(1,0)}/forward/backward%
}{
  \begin{scope}[shift={(\offset,0)}]
    \draw[chsix triangulation] (0,0)--(1,0)--(0,1)--cycle;
    \draw[chsix red curve,chsix \direction] \start--\finish;
    \node[chsix red point] at \endpoint {};
    \begin{scope}[shift={(1.22,0)}]
      \begin{scope}[rotate around={180:(.5,.5)}]
        \draw[chsix triangulation] (0,0)--(1,0)--(0,1)--cycle;
        \draw[chsix red curve,chsix \otherdirection] \start--\finish;
        \node[chsix red point] at \endpoint {};
      \end{scope}
    \end{scope}
  \end{scope}
}
\end{tikzpicture}
\caption{Triangle passages containing an endpoint}
\label{fig:minus-righttriangles-endpoint}
\end{figure}
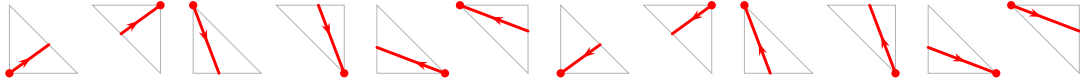

\item[(2)] \index{Sign assignment rules}For each triangle-passage occurrence $(\Delta,I)$ such that $\overline I$ contains neither endpoint, use $\gamma(\overline I)$ to assign one sign in $\{+,-\}$ as follows.
\begin{enumerate}
\item[(i)] Cut $\Delta$ along the image of $\gamma|_{\overline I}$. If the region on the left is a quadrilateral, assign $-$ (Figure~\ref{fig:minus-righttriangles}).
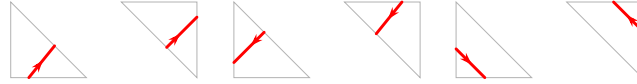
\begin{figure}[ht]
\centering
\begin{tikzpicture}[baseline=0mm]
\draw[chsix triangulation] (0,0) -- (1,0) -- (0,1) -- cycle;
\draw[chsix red curve,chsix forward] (0.24,0) -- (0.58,0.42);

\end{tikzpicture}\hspace{0.45cm}%
\begin{tikzpicture}[baseline=0mm]
\draw[chsix triangulation] (1,0) -- (1,1) -- (0,1) -- cycle;
\draw[chsix red curve,chsix backward] (1,0.8) -- (0.6,0.4);

\end{tikzpicture}\hspace{0.45cm}%
\rotatebox{180}{\begin{tikzpicture}[baseline=10mm]
\draw[chsix triangulation] (1,0) -- (1,1) -- (0,1) -- cycle;
\draw[chsix red curve,chsix backward] (1,0.8) -- (0.6,0.4);

\end{tikzpicture}}\hspace{0.45cm}%
\rotatebox{180}{\begin{tikzpicture}[baseline=10mm]
\draw[chsix triangulation] (0,0) -- (1,0) -- (0,1) -- cycle;
\draw[chsix red curve,chsix forward] (0.24,0) -- (0.58,0.42);

\end{tikzpicture}}\hspace{0.45cm}%
\begin{tikzpicture}[baseline=0mm]
\draw[chsix triangulation] (0,0) -- (1,0) -- (0,1) -- cycle;
\draw[chsix red curve,chsix forward] (0,0.38) -- (0.38,0);

\end{tikzpicture}\hspace{0.45cm}%
\begin{tikzpicture}[baseline=0mm]
\draw[chsix triangulation] (1,0) -- (1,1) -- (0,1) -- cycle;
\draw[chsix red curve,chsix forward] (1,0.62) -- (0.62,1);

\end{tikzpicture}
\caption{Negative signs for triangle passages}
\label{fig:minus-righttriangles}
\end{figure}

\item[(ii)] Assign $+$ to every remaining triangle-passage occurrence (Figure~\ref{fig:plus-righttriangles}).
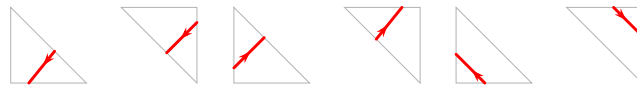
\begin{figure}[ht]
\centering
\begin{tikzpicture}[baseline=0mm]
\draw[chsix triangulation] (0,0) -- (1,0) -- (0,1) -- cycle;
\draw[chsix red curve,chsix backward] (0.24,0) -- (0.58,0.42);

\end{tikzpicture}\hspace{0.45cm}%
\begin{tikzpicture}[baseline=0mm]
\draw[chsix triangulation] (1,0) -- (1,1) -- (0,1) -- cycle;
\draw[chsix red curve,chsix forward] (1,0.8) -- (0.6,0.4);

\end{tikzpicture}\hspace{0.45cm}%
\rotatebox{180}{\begin{tikzpicture}[baseline=10mm]
\draw[chsix triangulation] (1,0) -- (1,1) -- (0,1) -- cycle;
\draw[chsix red curve,chsix forward] (1,0.8) -- (0.6,0.4);

\end{tikzpicture}}\hspace{0.45cm}%
\rotatebox{180}{\begin{tikzpicture}[baseline=10mm]
\draw[chsix triangulation] (0,0) -- (1,0) -- (0,1) -- cycle;
\draw[chsix red curve,chsix backward] (0.24,0) -- (0.58,0.42);

\end{tikzpicture}}\hspace{0.45cm}%
\begin{tikzpicture}[baseline=0mm]
\draw[chsix triangulation] (0,0) -- (1,0) -- (0,1) -- cycle;
\draw[chsix red curve,chsix backward] (0,0.38) -- (0.38,0);

\end{tikzpicture}\hspace{0.45cm}%
\begin{tikzpicture}[baseline=0mm]
\draw[chsix triangulation] (1,0) -- (1,1) -- (0,1) -- cycle;
\draw[chsix red curve,chsix backward] (1,0.62) -- (0.62,1);

\end{tikzpicture}
\caption{Positive signs for triangle passages}
\label{fig:plus-righttriangles}
\end{figure}
\end{enumerate}
This is the \emph{triangle-passage rule} for $\gamma$.

\item[(3)] \index{Edge type}\index{Sign assignment rules}For an edge-crossing occurrence $(e,t)$, let $i\in\{1,2,3\}$ be its edge type, corresponding to horizontal, diagonal, and vertical edges in that order. If $k_{\sigma(i)}=0$, assign no sign; this also applies when $\mu(e)\notin\mathcal V$ and the curve passes through $\mu(e)$. Now suppose $k_{\sigma(i)}>0$. Since $\mu(e)\in\mathcal V$ is an endpoint of the subdivided edge $e$, the curve avoids $\mu(e)$. Thus, near the crossing $\gamma(t)$, the midpoint $\mu(e)$ lies to the left or right of the oriented curve, and the sign is defined as follows.
\begin{enumerate}
\item[(i)] If $\mu(e)$ lies to the left of the oriented $\gamma$ near $\gamma(t)$, put $q=-$ and assign $k_{\sigma(i)}$ copies of $q$ (Figure~\ref{fig:minus-edge}).
\begin{figure}[ht]
\centering
\begin{tikzpicture}[baseline=0mm]
\draw[chsix distinguished edge] (-0.5,0)--(0.5,0); \node[chsix point] at (0,0) {};
\draw[chsix red curve,chsix forward] (-0.5,-0.3)--(0.5,0.1);

\end{tikzpicture}\hspace{0.5cm}
\begin{tikzpicture}[baseline=0mm]
\draw[chsix distinguished edge] (-0.5,0)--(0.5,0); \node[chsix point] at (0,0) {};
\draw[chsix red curve,chsix backward] (-0.5,-0.1)--(0.5,0.3);

\end{tikzpicture}\hspace{0.5cm}
\begin{tikzpicture}[baseline=0mm]
\draw[chsix distinguished edge] (-0.5,0.5)--(0.5,-0.5); \node[chsix point] at (0,0) {};
\draw[chsix red curve,chsix forward] (-0.5,-0.3)--(0.5,0);

\end{tikzpicture}\hspace{0.6cm}
\begin{tikzpicture}[baseline=0mm]
\draw[chsix distinguished edge] (-0.5,0.5)--(0.5,-0.5); \node[chsix point] at (0,0) {};
\draw[chsix red curve,chsix backward] (-0.5,0)--(0.5,0.3);

\end{tikzpicture}\hspace{0.7cm}
\begin{tikzpicture}[baseline=0mm]
\draw[chsix distinguished edge] (0,-0.5)--(0,0.5); \node[chsix point] at (0,0) {};
\draw[chsix red curve,chsix forward] (-0.5,-0.3)--(0.5,0);

\end{tikzpicture}\hspace{0.7cm}
\begin{tikzpicture}[baseline=0mm]
\draw[chsix distinguished edge] (0,-0.5)--(0,0.5); \node[chsix point] at (0,0) {};
\draw[chsix red curve,chsix backward] (-0.5,0)--(0.5,0.3);

\end{tikzpicture}
\caption{Negative signs for edge crossings}
\label{fig:minus-edge}
\end{figure}
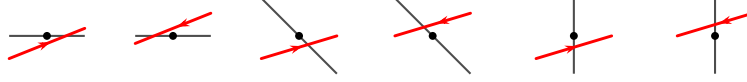

\item[(ii)] If $\mu(e)$ lies to the right of the oriented $\gamma$ near $\gamma(t)$, put $q=+$ and assign $k_{\sigma(i)}$ copies of $q$ (Figure~\ref{fig:plus-edge}).
\begin{figure}[ht]
\centering
\begin{tikzpicture}[baseline=0mm]
\draw[chsix distinguished edge] (-0.5,0)--(0.5,0); \node[chsix point] at (0,0) {};
\draw[chsix red curve,chsix backward] (-0.5,-0.3)--(0.5,0.1);

\end{tikzpicture}\hspace{0.5cm}
\begin{tikzpicture}[baseline=0mm]
\draw[chsix distinguished edge] (-0.5,0)--(0.5,0); \node[chsix point] at (0,0) {};
\draw[chsix red curve,chsix forward] (-0.5,-0.1)--(0.5,0.3);

\end{tikzpicture}\hspace{0.5cm}
\begin{tikzpicture}[baseline=0mm]
\draw[chsix distinguished edge] (-0.5,0.5)--(0.5,-0.5); \node[chsix point] at (0,0) {};
\draw[chsix red curve,chsix backward] (-0.5,-0.3)--(0.5,0);

\end{tikzpicture}\hspace{0.6cm}
\begin{tikzpicture}[baseline=0mm]
\draw[chsix distinguished edge] (-0.5,0.5)--(0.5,-0.5); \node[chsix point] at (0,0) {};
\draw[chsix red curve,chsix forward] (-0.5,0)--(0.5,0.3);

\end{tikzpicture}\hspace{0.7cm}
\begin{tikzpicture}[baseline=0mm]
\draw[chsix distinguished edge] (0,-0.5)--(0,0.5); \node[chsix point] at (0,0) {};
\draw[chsix red curve,chsix backward] (-0.5,-0.3)--(0.5,0);

\end{tikzpicture}\hspace{0.7cm}
\begin{tikzpicture}[baseline=0mm]
\draw[chsix distinguished edge] (0,-0.5)--(0,0.5); \node[chsix point] at (0,0) {};
\draw[chsix red curve,chsix forward] (-0.5,0)--(0.5,0.3);

\end{tikzpicture}
\caption{Positive signs for edge crossings}
\label{fig:plus-edge}
\end{figure}
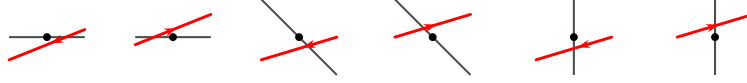

\end{enumerate}
This is the \emph{edge-crossing rule} for $\gamma$.
\end{enumerate}
\end{defi}

Fix $(k_1,k_2,k_3)$ and $\sigma$, and let $\gamma$ be an oriented generalized arc. Choose its endpoint signs and define $\varepsilon(\gamma)$ and $S(\gamma)$ as follows.
\begin{enumerate}
\item \index{Sign sequence}List the signs from rules (1)--(3) in occurrence order along $\gamma$. Denote this finite sequence by
\[
\varepsilon(\gamma)=(\varepsilon_1,\ldots,\varepsilon_m).
\]
\item \index{Run-length sequence}Let $a_0,\ldots,a_n$ be the lengths of the maximal constant-sign runs in $\varepsilon(\gamma)$, from left to right. Denote this positive integer sequence by
\[
S(\gamma)=(a_0,\ldots,a_n).
\]
\end{enumerate}

For $\varepsilon(\gamma)=(\varepsilon_1,\ldots,\varepsilon_m)$, define the fence poset $P_\gamma$ to have underlying set $\{1,2,\ldots,m-1\}$. For $1\leq i\leq m-2$, impose $i\lessdot i+1$ if $\varepsilon_{i+1}=+$, and $i\gtrdot i+1$ if $\varepsilon_{i+1}=-$.

The sign $\varepsilon_{i+1}$ records the cover relation between vertices $i$ and $i+1$, rather than a vertex itself.

If the first and last triangle-passage occurrences coincide, there are no crossings and $m=1$, so $P_\gamma$ is empty. Otherwise the endpoint choices affect only $\varepsilon_1$ or $\varepsilon_m$, neither of which defines a cover relation. Thus $P_\gamma$ is well-defined independently of the endpoint rule.

\index{Generalized Markov length}The number of order ideals of $P_\gamma$ is the \textbf{$(k_1,k_2,k_3,\sigma)$-generalized Markov length}, or \textbf{GM length}, of $\gamma$. We write it as $|\gamma|$, suppressing its dependence on $(k_1,k_2,k_3)$ and $\sigma$.

If the orientation is reversed and the chosen endpoint signs are also reversed, rules (1)--(3) give the new sign sequence
\[
(\bar\varepsilon_m,\ldots,\bar\varepsilon_1)
\]
where $\bar{+}=-$ and $\bar{-}=+$. Its run-length sequence is $(a_n,\ldots,a_0)$. Since $P_\gamma$ does not depend on the endpoint choices, Proposition~\ref{prop:reverse-property} shows that $|\gamma|$ does not depend on orientation.

\begin{prop}\label{prop:length}\index{Fence poset}\index{GM length|seeonly{Generalized Markov length}}
If $S(\gamma)=(a_0,a_1,\ldots,a_n)$, then $P_\gamma$ is isomorphic either to $P_{(a_0,\ldots,a_n)}$ or to its dual $P_{(a_0,\ldots,a_n)}^\ast$. In particular,
\[
|\gamma|=N(a_0,\ldots,a_n).
\]
\end{prop}
\begin{proof}
If $n=0$, all signs in $\varepsilon(\gamma)$ agree, and $P_\gamma$ is a chain with $a_0-1$ vertices, empty when $a_0=1$. Hence
\[
|\gamma|=a_0=N(a_0).
\]
Now suppose $n\geq1$. The sequence $S(\gamma)=(a_0,\ldots,a_n)$ lists the constant-sign run lengths from left to right, and the cover orientation of $P_\gamma$ changes precisely when the sign changes. Thus the numbers of consecutive cover relations of the same orientation are
\[
a_0-1,\ a_1,\ a_2,\ldots,\ a_{n-1},\ a_n-1
\]
as in the definition of $P_{(a_0,\ldots,a_n)}$. If the corresponding cover orientations agree in the two posets, then
$P_\gamma\cong P_{(a_0,\ldots,a_n)}$; if every corresponding cover orientation is reversed, then
$P_\gamma\cong P_{(a_0,\ldots,a_n)}^\ast$. The formula for $|\gamma|$ follows from the definition of $N(a_0,\ldots,a_n)$ and the identity $N(P^\ast)=N(P)$.
\end{proof}

\[
U_+:=\begin{bmatrix}1&0\\1&1\end{bmatrix},\qquad
U_-:=\begin{bmatrix}1&1\\0&1\end{bmatrix},\qquad
\mathbf1:=\begin{bmatrix}1\\1\end{bmatrix}
\]
Write $E_2$ for the $2\times2$ identity matrix. Let $\varepsilon(\eta)=\varepsilon_1\cdots\varepsilon_m$ be the sign sequence of an oriented generalized arc $\eta$, and suppose $m\geq2$. The vertex set of $P_\eta$ is
\[
V_\eta:=\{1,\ldots,m-1\}
\]
For $J\subseteq V_\eta$, record membership by
\[
x_i=\begin{cases}
0&(i\notin J),\\
1&(i\in J)
\end{cases}
\qquad(1\leq i\leq m-1)
\]
If $\varepsilon_{i+1}=+$, then $i\lessdot i+1$, so downward closedness at this cover is equivalent to $x_{i+1}\leq x_i$. The allowed adjacent states are
\[
(x_i,x_{i+1})=(0,0),(1,0),(1,1)
\]
Similarly, if $\varepsilon_{i+1}=-$, then $i\gtrdot i+1$, so $x_i\leq x_{i+1}$, and the allowed states are
\[
(x_i,x_{i+1})=(0,0),(0,1),(1,1)
\]
Let rows represent the current state $x_i$ and columns the next state $x_{i+1}$, with the first row and column corresponding to state $0$ and the second to state $1$. The matrices $U_+,U_-$ then give the following tables.
\[
\begin{array}{c|cc}
U_+&x_{i+1}=0&x_{i+1}=1\\ \hline
x_i=0&1&0\\
x_i=1&1&1
\end{array}
\qquad
\begin{array}{c|cc}
U_-&x_{i+1}=0&x_{i+1}=1\\ \hline
x_i=0&1&1\\
x_i=1&0&1
\end{array}
\]
An entry $1$ means the state pair is allowed, and $0$ means it is forbidden. For example, the upper-right entry of $U_+$ is $0$ because $(x_i,x_{i+1})=(0,1)$ violates downward closedness. Since matrix indices start at $1$,
\[
\bigl(U_{\varepsilon_{i+1}}\bigr)_{x_i+1,x_{i+1}+1}
=\begin{cases}
1&\text{if the state pair $(x_i,x_{i+1})$ is allowed},\\
0&\text{otherwise}.
\end{cases}
\]

Matrix multiplication sums over the states shared by consecutive factors. For example,
\[
\bigl(U_{\varepsilon_2}U_{\varepsilon_3}\bigr)_{x_1+1,x_3+1}
=
\sum_{x_2\in\{0,1\}}
\bigl(U_{\varepsilon_2}\bigr)_{x_1+1,x_2+1}
\bigl(U_{\varepsilon_3}\bigr)_{x_2+1,x_3+1}
\]
sums over the intermediate states $x_2=0,1$ with initial state $x_1$ and final state $x_3$ fixed. Repeating this and multiplying by $\mathbf1^{\mathsf T}$ and $\mathbf1$ on the left and right also sums over the endpoint states, giving
\[
\mathbf1^{\mathsf T}U_{\varepsilon_2}\cdots U_{\varepsilon_{m-1}}\mathbf1
=
\sum_{(x_1,\ldots,x_{m-1})\in\{0,1\}^{m-1}}
\prod_{i=1}^{m-2}
\bigl(U_{\varepsilon_{i+1}}\bigr)_{x_i+1,x_{i+1}+1}
\]
Each product on the right is $1$ exactly when all adjacent states are allowed, equivalently when $J_x:=\{i\in V_\eta\mid x_i=1\}$ is an order ideal of $P_\eta$, and is $0$ otherwise. Thus the sum counts the order ideals of $P_\eta$, and we obtain
\begin{equation}\label{eq:gm-transfer}\index{Transfer matrix}\index{Generalized Markov length}\index{Order ideal}
|\eta|=\mathbf1^{\mathsf T}U_{\varepsilon_2}\cdots U_{\varepsilon_{m-1}}\mathbf1
\end{equation}
The cover sign between vertices $i$ and $i+1$ is $\varepsilon_{i+1}$ ($1\leq i\leq m-2$), so $\varepsilon_1$ and $\varepsilon_m$ do not occur in the product. When $m=2$, the empty product is $E_2$, and $\mathbf1^{\mathsf T}E_2\mathbf1=2$ counts the two order ideals of a one-vertex poset. When $m=1$, the fence is empty and $|\eta|=1$.
\index{Transfer matrix}The matrices $U_+,U_-$ are called \emph{transfer matrices}: their entries record which pairs of adjacent membership states are allowed, and their products count allowed state sequences. We call \eqref{eq:gm-transfer}, which expresses GM length through their product, the \emph{transfer formula}.

\begin{exam}\label{exam:gm-length-calculation}
Take $\sigma=\operatorname{id}$ and $(k_1,k_2,k_3)=(1,2,1)$. Set $A=(0,1)$, $B=(2,2)$, and
\[
 C_1=\left(1,\frac14\right),\quad
 C_2=\left(\frac74,\frac14\right),\quad
 C_3=\left(\frac74,1\right),\quad
 C_4=\left(\frac54,\frac74\right)
\]
Let $\gamma$ be the polygonal path joining $A,C_1,C_2,C_3,C_4,B$ in this order (Figure~\ref{fig:gm-length-example}). It passes once through each of the five lightly shaded triangles, crossing a vertical edge, a diagonal edge, a horizontal edge, and another diagonal edge in order.

The black dots mark the midpoints of the crossed edges. At $C_1,C_2,C_3$, the midpoint lies to the left of the arc, so the assigned signs are $-$, $--$, and $-$, respectively. At $C_4$, the midpoint lies to the right, giving $++$. For the passages from $C_1$ to $C_2$ and from $C_2$ to $C_3$, the region on the left is a triangle, so both signs are $+$. From $C_3$ to $C_4$, it is a quadrilateral, so the sign is $-$. Choosing $+$ for both endpoint passages gives
\[
 \varepsilon(\gamma)=(+,-,+,-,-,+,-,-,+,+,+),
 \qquad S(\gamma)=(1,1,1,2,1,2,3)
\]
The poset $P_\gamma$ has ten vertices, with cover relations
\[
 1\gtrdot2\lessdot3\gtrdot4\gtrdot5\lessdot6\gtrdot7\gtrdot8\lessdot9\lessdot10.
\]
\begin{figure}[htbp]
\centering
\begin{tikzpicture}[line cap=round,line join=round]
  \begin{scope}[x=2.35cm,y=2.35cm]
    \node[chsix panel title] at (1,2.23) {$\gamma$};
    \fill[black!6] (0,1)--(1,0)--(2,0)--(2,2)--(1,2)--(1,1)--cycle;
    \foreach \i in {0,1,2}{
      \draw[chsix triangulation] (\i,0)--(\i,2);
      \draw[chsix triangulation] (0,\i)--(2,\i);
    }
    \draw[chsix triangulation] (0,1)--(1,0) (0,2)--(2,0) (1,2)--(2,1);
    \coordinate (gmA) at (0,1);
    \coordinate (gmB) at (2,2);
    \coordinate (gmC1) at (1,{1/4});
    \coordinate (gmC2) at ({7/4},{1/4});
    \coordinate (gmC3) at ({7/4},1);
    \coordinate (gmC4) at ({5/4},{7/4});
    \draw[chsix red curve,
      postaction={decorate},decoration={markings,
        mark=at position .13 with {\arrow{Stealth[length=1.65mm,width=1.10mm]}},
        mark=at position .56 with {\arrow{Stealth[length=1.65mm,width=1.10mm]}},
        mark=at position .88 with {\arrow{Stealth[length=1.65mm,width=1.10mm]}}}]
      (gmA)--(gmC1)--(gmC2)--(gmC3)--(gmC4)--(gmB);
    \foreach \P in {gmA,gmB,gmC1,gmC2,gmC3,gmC4}
      \node[chsix red point] at (\P) {};
    \foreach \P in {(1,.5),(1.5,.5),(1.5,1),(1.5,1.5)}
      \node[chsix point] at \P {};
    \node[below left=1pt] at (gmA) {$A$};
    \node[above right=1pt] at (gmB) {$B$};
    \node[below left=2pt,font=\scriptsize] at (gmC1) {$C_1$};
    \node[below right=2pt,font=\scriptsize] at (gmC2) {$C_2$};
    \node[right=2pt,font=\scriptsize] at (gmC3) {$C_3$};
    \node[above left=2pt,font=\scriptsize] at (gmC4) {$C_4$};
  \end{scope}
  \draw[chsix transition arrow] (5.20,2.35)--(5.85,2.35);
  \begin{scope}[shift={(6.25cm,2.64cm)},x=.58cm,y=.58cm]
    \node[chsix panel title] at (4.5,1.85) {$P_\gamma$};
    \foreach \i/\x/\y in {1/0/1,2/1/0,3/2/1,4/3/0,5/4/-1,6/5/0,7/6/-1,8/7/-2,9/8/-1,10/9/0}
      \node[chsix fence vertex] (gmp\i) at (\x,\y) {$\i$};
    \foreach \i/\j in {1/2,2/3,3/4,4/5,5/6,6/7,7/8,8/9,9/10}
      \draw[chsix fence edge] (gmp\i)--(gmp\j);
  \end{scope}
\end{tikzpicture}
\caption{An example of computing GM length}
\label{fig:gm-length-example}
\end{figure}
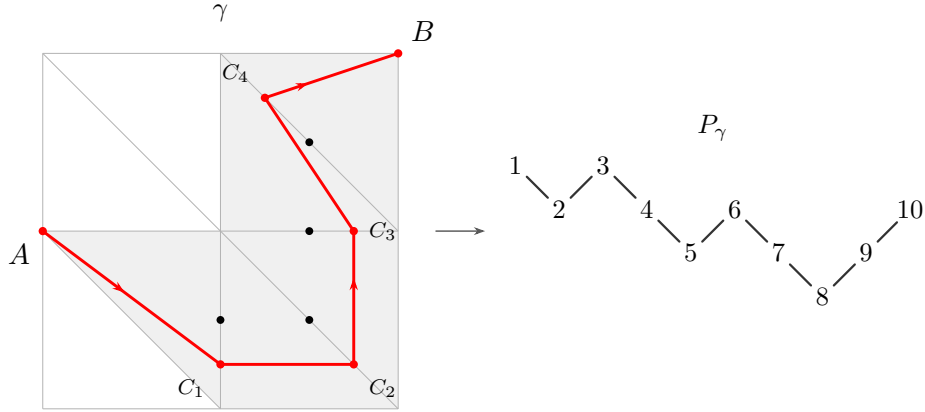
\FloatBarrier

Count the order ideals according to whether they contain vertex $5$. If they do not, they cannot contain $3,4,6$, and there are three choices on $\{1,2\}$. On $\{7,8,9,10\}$, there is one choice without $8$ and $2\cdot3$ choices with $8$, giving $3(1+2\cdot3)=21$ ideals. If $5$ is included, the choices on $\{1,2,3,4\}$ and $\{6,7,8,9,10\}$ are independent. For the former, there are $3\cdot2$ choices without $3$ and two with $3$. For the latter, there is one choice without $8$ and $3\cdot3$ with $8$, giving $(3\cdot2+2)(1+3\cdot3)=80$ ideals. Hence $|\gamma|=21+80=101$. The transfer formula, omitting the two endpoint signs, gives the same result:
\[
 |\gamma|
 =\mathbf1^{\mathsf T}U_-U_+U_-^2U_+U_-^2U_+^2\mathbf1
 =\mathbf1^{\mathsf T}\begin{bmatrix}45&19\\26&11\end{bmatrix}\mathbf1
 =101.
\]
\end{exam}

\begin{defi}\label{def:gm-distance}\index{Generalized Markov distance}\index{GM distance|seeonly{Generalized Markov distance}}
For distinct lattice points $A,B$, set
\[
\mathcal A_0(A,B)
:=\{\gamma\mid\gamma\text{ is a generalized arc joining }A,B\text{ with no self-intersections}\}
\]
and define
\[
d(A,B):=\inf_{\gamma\in\mathcal A_0(A,B)}|\gamma|
\]
We call $d(A,B)$ the \textbf{$(k_1,k_2,k_3,\sigma)$-generalized Markov distance}, or \textbf{GM distance}, between $A$ and $B$.
Here ``distance'' denotes this combinatorial quantity, not a metric in the usual sense; in particular, the triangle inequality need not hold.
\end{defi}

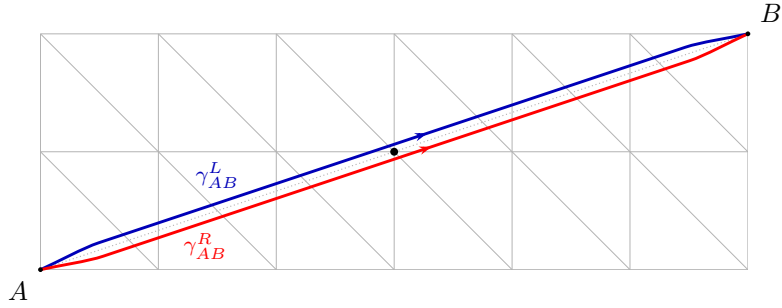
\begin{figure}[ht]
\centering
\begin{tikzpicture}[x=0.01cm,y=0.01cm,scale=0.6]
  \tikzset{
    gridline/.style={chsix triangulation},
    gammaL/.style={chsix blue curve,chsix forward},
    gammaR/.style={chsix red curve,chsix forward},
    enddot/.style={circle,fill=black,inner sep=0pt,minimum size=1.8pt}
  }

  \begin{scope}[shift={(180,170)}, x={(260,0)}, y={(0,260)}]

    \begin{scope}
      \clip (0,0) rectangle (6,2);

      \foreach \i in {0,...,5}{
        \foreach \j in {0,1}{
          \draw[gridline] (\i,{\j+1}) -- ({\i+1},\j);
        }
      }

      \foreach \i in {0,...,6}{
        \draw[gridline] (\i,0) -- (\i,2);
      }
      \foreach \j in {0,1,2}{
        \draw[gridline] (0,\j) -- (6,\j);
      }
    \end{scope}

    \coordinate (A) at (0,0);
    \coordinate (B) at (6,2);

    \def\tx{0.9486832981}
    \def\ty{0.3162277660}
    \def\nLx{-0.3162277660}
    \def\nLy{0.9486832981}
    \def\nRx{0.3162277660}
    \def\nRy{-0.9486832981}

    \def\eps{0.060}

    \def\sep{0.016}

    \def\endstraight{0.32}
    \def\endblend{0.22}

    \def\ctrl{0.055}

    \coordinate (AsepL) at ({0.070*\tx+\sep*\nLx},
                            {0.070*\ty+\sep*\nLy});
    \coordinate (AsepR) at ({0.070*\tx+\sep*\nRx},
                            {0.070*\ty+\sep*\nRy});

    \coordinate (ApreL) at ({\endstraight*\tx+\sep*\nLx},
                            {\endstraight*\ty+\sep*\nLy});
    \coordinate (ApreR) at ({\endstraight*\tx+\sep*\nRx},
                            {\endstraight*\ty+\sep*\nRy});

    \coordinate (ALpost) at ({(\endstraight+\endblend)*\tx+\eps*\nLx},
                             {(\endstraight+\endblend)*\ty+\eps*\nLy});
    \coordinate (ARpost) at ({(\endstraight+\endblend)*\tx+\eps*\nRx},
                             {(\endstraight+\endblend)*\ty+\eps*\nRy});

    \coordinate (BsepL) at ({6-0.070*\tx+\sep*\nLx},
                            {2-0.070*\ty+\sep*\nLy});
    \coordinate (BsepR) at ({6-0.070*\tx+\sep*\nRx},
                            {2-0.070*\ty+\sep*\nRy});

    \coordinate (BpreL) at ({6-\endstraight*\tx+\sep*\nLx},
                            {2-\endstraight*\ty+\sep*\nLy});
    \coordinate (BpreR) at ({6-\endstraight*\tx+\sep*\nRx},
                            {2-\endstraight*\ty+\sep*\nRy});

    \coordinate (BLpost) at ({6-(\endstraight+\endblend)*\tx+\eps*\nLx},
                             {2-(\endstraight+\endblend)*\ty+\eps*\nLy});
    \coordinate (BRpost) at ({6-(\endstraight+\endblend)*\tx+\eps*\nRx},
                             {2-(\endstraight+\endblend)*\ty+\eps*\nRy});

    \coordinate (CmidL) at ({3+\eps*\nLx},{1+\eps*\nLy});
    \coordinate (CmidR) at ({3+\eps*\nRx},{1+\eps*\nRy});

    \draw[chsix reference segment] (A)--(B);
    \node[chsix point] at (3,1) {};

    \draw[gammaL]
      (A)
      .. controls
        ({.17*\tx+.018*\nLx},
         {.17*\ty+.018*\nLy})
        and
        ({.40*\tx+\eps*\nLx},
         {.40*\ty+\eps*\nLy})
      .. (ALpost)
      -- (CmidL)
      -- (BLpost)
      .. controls
        ({6-.40*\tx+\eps*\nLx},
         {2-.40*\ty+\eps*\nLy})
        and
        ({6-.17*\tx+.018*\nLx},
         {2-.17*\ty+.018*\nLy})
      .. (B);

    \draw[gammaR]
      (A)
      .. controls
        ({.17*\tx+.018*\nRx},
         {.17*\ty+.018*\nRy})
        and
        ({.40*\tx+\eps*\nRx},
         {.40*\ty+\eps*\nRy})
      .. (ARpost)
      -- (CmidR)
      -- (BRpost)
      .. controls
        ({6-.40*\tx+\eps*\nRx},
         {2-.40*\ty+\eps*\nRy})
        and
        ({6-.17*\tx+.018*\nRx},
         {2-.17*\ty+.018*\nRy})
      .. (B);

    \node[enddot,label={[font=\small]below left:$A$}] at (A) {};
    \node[enddot,label={[font=\small]above right:$B$}] at (B) {};

    \node[text=blue!72!black,font=\scriptsize] at (1.5,0.8) {$\gamma_{AB}^{L}$};
    \node[text=annred,font=\scriptsize] at (1.4,0.2) {$\gamma_{AB}^{R}$};
  \end{scope}
\end{tikzpicture}
\caption{The pure left and right push-offs}
\label{fig:gamma-L-R}
\end{figure}

\index{Push-off!left}\index{Push-off!right}Orient the segment $AB$ from $A$ to $B$ and take a sufficiently narrow strip around it. A \emph{push-off} is a curve obtained by fixing the endpoints and moving the interior slightly to a specified side. Let $\gamma^R_{AB}$ be a simple, sufficiently small push-off that meets $AB$ only at its endpoints and detours to the right of every point on the open segment $AB$. If $AB$ lies on a triangulation line, push it into the triangles on its right; define $\gamma^L_{AB}$ similarly on the left. Between consecutive points on $AB$, keep the curve sufficiently close to $AB$ to avoid crossing additional triangulation edges. If $AB$ is not on a triangulation line, a detour around a midpoint in $\mathcal V$ crosses its whole triangulation edge exactly once. If it follows a triangulation line, cross edges in the other two directions alternately. Both push-offs are generalized arcs without self-intersections and thus belong to $\mathcal A_0(A,B)$ (Figure~\ref{fig:gamma-L-R}).

Two push-offs to the same side can be continuously deformed into one another while preserving the order of triangulation-edge crossings. Thus $P_{\gamma^R_{AB}}$, $P_{\gamma^L_{AB}}$, and their GM lengths are well-defined by the specified side of $AB$.

In Figure~\ref{fig:gamma-L-R}, the dotted line is $AB$ and the central black dot is a lattice point on its interior. The red curve $\gamma^R_{AB}$ passes on the right, and the blue curve $\gamma^L_{AB}$ on the left.

Since
\[
\gamma^R_{AB},\gamma^L_{AB}
\in\mathcal A_0(A,B)
\]
the set $\mathcal A_0(A,B)$ is nonempty. GM lengths are positive integers, so the infimum in Definition~\ref{def:gm-distance} is attained.

\begin{theorem}\label{thm:gm-distance}\index{Generalized Markov distance}\index{Push-off!left}\index{Push-off!right}
Let $p,q$ be relatively prime positive integers, and put $A=(0,0)$, $B=(q,p)$. Then
\[
d(A,B)=|\gamma^R_{AB}|=|\gamma^L_{AB}|.
\]
\end{theorem}

We prepare for the proof of Theorem~\ref{thm:gm-distance} by introducing terminology and proving several lemmas.

First we show that some arc of minimum GM length never passes through the same triangle twice. We then compare left and right midpoint choices using $U_+,U_-$. Next we glue the passage triangles in order to form a closed triangle strip, defined in Lemma~\ref{lem:fresh-canonical-geodesic}, and construct a polygonal path minimizing ordinary Euclidean length within it. We explicitly construct generalized arcs following this path and compare their GM lengths when the side on which a lattice point is avoided is changed.

We use the following terminology. A generalized arc without self-intersections is also called a \emph{simple arc}. Its \emph{triangle-passage sequence} lists, in occurrence order, the triangles of its triangle-passage occurrences from Definition~\ref{def:passage-occurrence}. Repeated passages through the same planar triangle are distinct terms. For a passage joining two different sides of a triangle, the common vertex of those sides is the \emph{cut-off vertex}.

Write $n(\gamma)$ for the number of edge-crossing occurrences, counting repeated crossings of the same edge separately. A \emph{crossed edge} will mean the whole triangulation edge crossed at an occurrence; subdivided edges will be specified explicitly. Horizontal, diagonal, and vertical edges have types $1,2,3$, respectively, and the integer $k_{\sigma(i)}$ assigned to a type-$i$ edge is its \emph{weight}.

A finite sign sequence is also called a \emph{word}, and juxtaposition denotes concatenation. A consecutive portion is a \emph{subword}. Write $w^a$ for the concatenation of $a$ copies of $w$, with $w^0$ the empty word. \emph{Sign reversal} replaces each sign by $\bar+=-$ or $\bar-=+$, and is distinct from reversing the order of the signs. For $w=w_1\cdots w_s$, put
\[
 M(w):=U_{w_1}\cdots U_{w_s},\qquad M(\varnothing):=E_2
\]
Inequalities between matrices, rows, and columns are entrywise: $X\geq Y$ means $X_{ij}\geq Y_{ij}$ for all corresponding entries. The notation $X>0$ means every entry is positive; a ``positive row'' or ``positive column'' has the same meaning. If there is a crossing, equivalently if the sign sequence has length $m\geq2$, the transfer formula \eqref{eq:gm-transfer} expresses GM length as the product of the word with its first and last triangle signs omitted, multiplied by $\mathbf1^{\mathsf T}$ on the left and $\mathbf1$ on the right. Call this word the \emph{interior word}. The constant-sign subword $q^{k_{\sigma(i)}}$ assigned to a type-$i$ crossing is an \emph{edge block}, with product $U_q^{k_{\sigma(i)}}$. At a crossing of weight $0$, the edge block is empty and its product is $E_2$. We repeatedly use
\begin{equation}\label{eq:fresh-word-domination}
 M(w)\geq E_2,\qquad
 w\text{ contains both signs }\Longrightarrow M(w)\geq U_+,\ U_-.
\end{equation}
Indeed, each factor is at least $E_2$, and products of nonnegative matrices preserve entrywise inequalities. An arc with no crossings has GM length $1$ by definition.

\begin{lemm}\label{lem:same-length-lemma0}
For any distinct lattice points $A,B$,
\[
 |\gamma^R_{AB}|=|\gamma^L_{AB}|.
\]
\end{lemm}
\begin{proof}
The half-turn $x\mapsto A+B-x$ preserves the triangulation, edge types, and point set $\mathcal V$, and exchanges $A,B$. Traversing the image of $\gamma^R_{AB}$ in reverse gives a left push-off from $A$ to $B$. The half-turn preserves signs and edge weights, and orientation reversal preserves GM length, proving the assertion.
\end{proof}

\begin{lemm}\label{lem:fresh-minimizer}\index{Generalized Markov distance!minimizing arcs}\index{Generalized arc}
Fix distinct lattice points $A,B$. Among all generalized arcs from $A$ to $B$, allowing self-intersections, choose one for which $(|\gamma|,n(\gamma))$ is lexicographically minimal. Thus we first minimize GM length and then minimize the number of crossings among arcs of that length. Such an arc exists and satisfies the following properties.
\begin{enumerate}
\item It never passes through the same geometric triangle twice. In particular, it has no self-intersections.
\item No passage triangle except the first has $A$ as a vertex, and none except the last has $B$ as a vertex.
\item $|\gamma|=d(A,B)$.
\end{enumerate}
\end{lemm}
\begin{proof}
The candidate set contains the pure push-offs, GM lengths are positive integers, and crossing counts are nonnegative integers. Hence a minimum pair exists.

First prove (2). Suppose a passage triangle other than the first has $A$ as a vertex, and let $\Delta$ be the last such triangle. If it is the last passage triangle, join $A,B$ directly inside $\Delta$. Otherwise the edge through which $\gamma$ leaves $\Delta$ is opposite $A$: leaving through an edge incident to $A$ would make the next triangle another one containing $A$, contrary to the choice of $\Delta$. Join $A$ directly to the crossing point on this opposite edge, deleting the initial portion. A small perturbation inside the triangle makes intersections with other passages transverse and gives a generalized arc satisfying the endpoint condition (Figure~\ref{fig:initial-part-shortcut}). If crossings remain, this merely replaces the deleted initial product by $E_2$ in the GM-length formula, so \eqref{eq:fresh-word-domination} shows that GM length does not increase. If all crossings disappear, the length is $1$. In either case the crossing count decreases, contradicting lexicographic minimality. The same argument at $B$ proves (2).

\begin{figure}[htbp]
\centering
\begin{tikzpicture}[
  x=2.00cm,y=2.00cm,
  line cap=round,line join=round,
  tri/.style={chsix triangulation},
  kept/.style={chsix red curve},
  joined/.style={chsix red curve},
  erased/.style={draw=annred!55,line width=.8pt,dash pattern=on 2.3pt off 1.9pt},
  every node/.style={font=\small,inner sep=1.4pt}
]
\foreach \panel in {0,1}{
  \begin{scope}[xshift={\panel*6.80cm}]
    \path[fill=black!5] (0,0)--(1,0)--(0,1)--cycle;
    \begin{scope}
      \clip (-1.55,-.72) rectangle (1.24,1.20);
      \foreach \i in {-1,0,1}{
        \draw[tri] (\i,-.72)--(\i,1.20);
      }
      \foreach \i in {0,1}{
        \draw[tri] (-1.55,\i)--(1.24,\i);
      }
      \foreach \i in {-2,-1,0,1,2}{
        \draw[tri] (-1.55,{\i+1.55})--(1.24,{\i-1.24});
      }
    \end{scope}
    \draw[chsix distinguished edge] (0,0)--(1,0)--(0,1)--cycle;
    \node at (.18,.45) {$\Delta$};
    \draw[kept,chsix forward at=.56]
      (.55,.45) .. controls (.74,.68) and (.87,.90) .. (1,1);
    \ifnum\panel=0
      \draw[kept,chsix forward at=.44]
        (0,0) .. controls (-.30,-.10) and (-.44,-.60) .. (-.90,-.60)
        .. controls (-1.55,-.60) and (-1.55,.42) .. (-.90,.55)
        .. controls (-.62,.64) and (-.27,.70) .. (0,.65);
      \draw[kept]
        (0,.65) .. controls (.16,.65) and (.40,.56) .. (.55,.45);
      \node[chsix point] at (0,.65) {};
    \else
      \draw[erased]
        (0,0) .. controls (-.30,-.10) and (-.44,-.60) .. (-.90,-.60)
        .. controls (-1.55,-.60) and (-1.55,.42) .. (-.90,.55)
        .. controls (-.62,.64) and (-.27,.70) .. (0,.65)
        .. controls (.16,.65) and (.40,.56) .. (.55,.45);
      \draw[joined,chsix forward]
        (0,0) .. controls (.20,.16) and (.43,.32) .. (.55,.45);
    \fi
    \node[chsix point] at (0,0) {};
    \node[chsix point] at (.55,.45) {};
    \node[chsix point] at (1,1) {};
    \node at (.08,-.12) {$A$};
    \node at (1.11,1.11) {$B$};
  \end{scope}
}
\draw[chsix transition arrow] (1.40,.38)--(1.70,.38);
\end{tikzpicture}
\caption{Deleting an initial portion}
\label{fig:initial-part-shortcut}
\end{figure}
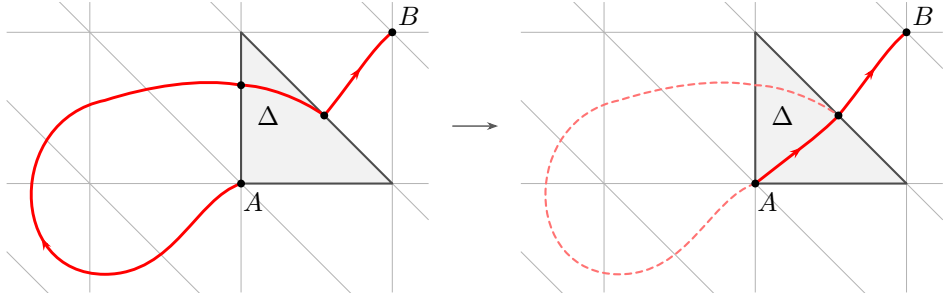

Next prove (1). Suppose $\gamma$ passes through the same triangle $\Delta$ twice. Consider deleting the portion from its entry point in the first chosen passage to its exit point in the later passage, and joining these points directly inside $\Delta$ (Figure~\ref{fig:repeated-triangle-shortcut}). If both points lie on the same triangulation edge, as on the right of the figure, the displayed connection alone is not a generalized arc. Include these two crossings in the deleted interval. Immediately before and after this enlarged interval, the curve lies in the same triangle on the opposite side of that edge, so reconsider the connection there. If the procedure reaches an endpoint, the same triangle containing that endpoint occurs at both ends of the interval, contradicting (2). Each step adds two crossings to the deleted interval. Hence after finitely many steps its ends lie on different sides of one triangle, while the original first and last passages remain.

After choosing the deletion interval in this way, join its ends inside the final triangle. Keep the positions and directions of the retained crossings and reconnect them by simple curves inside each triangle. If necessary, perturb within the triangle interiors so that distinct passages meet only in finitely many transverse intersections. No consecutive crossings of the same triangulation edge arise at the connection, and the first and last passages are unchanged. Thus the new curve $\eta$ is a generalized arc with fewer crossings. We compare its GM length with that of $\gamma$.

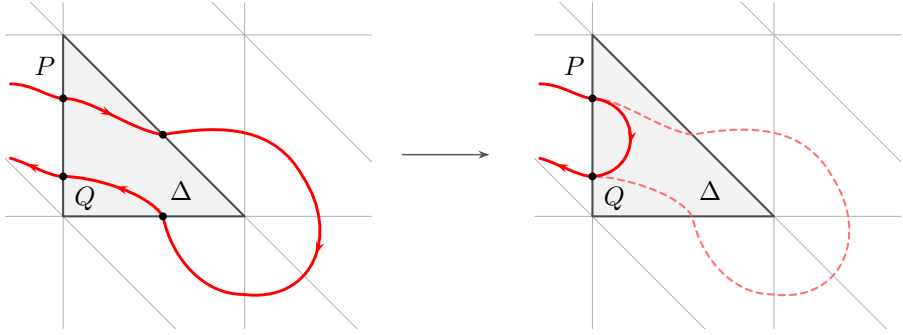
\begin{figure}[htbp]
\centering
\begin{tikzpicture}[
  x=2.40cm,y=2.40cm,
  line cap=round,line join=round,
  tri/.style={chsix triangulation},
  kept/.style={chsix red curve},
  joined/.style={chsix red curve},
  erased/.style={draw=annred!55,line width=.8pt,dash pattern=on 2.3pt off 1.9pt},
  every node/.style={font=\small,inner sep=1.4pt}
]
\foreach \panel in {0,1}{
  \begin{scope}[xshift={\panel*7.00cm}]
    \path[fill=black!5] (0,0)--(1,0)--(0,1)--cycle;
    \begin{scope}
      \clip (-.32,-.62) rectangle (1.70,1.18);
      \foreach \i in {0,1,2}{
        \draw[tri] (\i,-.62)--(\i,1.18);
      }
      \foreach \i in {0,1}{
        \draw[tri] (-.32,\i)--(1.70,\i);
      }
      \foreach \i in {-1,0,1,2}{
        \draw[tri] (-.32,{\i+.32})--(1.70,{\i-1.70});
      }
    \end{scope}
    \draw[chsix distinguished edge] (0,0)--(1,0)--(0,1)--cycle;
    \node at (.65,.14) {$\Delta$};
    \draw[kept]
      (-.29,.73) .. controls (-.17,.73) and (-.07,.65) .. (0,.65);
    \draw[kept,chsix forward at=.70]
      (0,.22) .. controls (-.08,.22) and (-.18,.29) .. (-.29,.32);
    \ifnum\panel=0
      \draw[kept,chsix forward at=.50]
        (0,.65) .. controls (.16,.65) and (.36,.48) .. (.55,.45);
      \draw[kept,chsix forward at=.53]
        (.55,.45) .. controls (.89,.51) and (1.17,.49) .. (1.32,.22)
        .. controls (1.53,-.10) and (1.38,-.48) .. (.99,-.43)
        .. controls (.78,-.43) and (.60,-.24) .. (.55,0);
      \draw[kept,chsix forward at=.52]
        (.55,0) .. controls (.47,.13) and (.16,.22) .. (0,.22);
      \node[chsix point] at (.55,.45) {};
      \node[chsix point] at (.55,0) {};
    \else
      \draw[erased]
        (0,.65) .. controls (.16,.65) and (.36,.48) .. (.55,.45)
        .. controls (.89,.51) and (1.17,.49) .. (1.32,.22)
        .. controls (1.53,-.10) and (1.38,-.48) .. (.99,-.43)
        .. controls (.78,-.43) and (.60,-.24) .. (.55,0)
        .. controls (.47,.13) and (.16,.22) .. (0,.22);
      \draw[joined,chsix forward]
        (0,.65) .. controls (.28,.62) and (.28,.24) .. (0,.22);
    \fi
    \node[chsix point] at (0,.65) {};
    \node[chsix point] at (0,.22) {};
    \node at (-.10,.83) {$P$};
    \node at (.12,.10) {$Q$};
  \end{scope}
}
\draw[chsix transition arrow] (1.87,.34)--(2.35,.34);
\end{tikzpicture}
\caption{Connecting within the same triangle}
\label{fig:repeated-triangle-shortcut}
\end{figure}

There are retained crossings on both sides of the deleted interval. Let $w$ be the original subword between the two initially chosen passages in $\Delta$, including their triangle signs. If $w$ contains both signs, so does the deleted subword even when the interval is enlarged. By \eqref{eq:fresh-word-domination}, its product is at least both $U_+$ and $U_-$, hence at least the matrix for the single triangle sign after reconnection. Now suppose $w$ consists of a single sign $q$. Consecutive passages with the same sign cut off the same endpoint of their shared edge, so all passages in this interval cut off one lattice point $V$. Figure~\ref{fig:same-sign-passages} shows the case $q=+$. Let $P$ be the entry point in the first passage through $\Delta$ and $Q$ the exit point in the later passage; the deleted interval is red and dashed. The curve proceeds counterclockwise around $V$ and revisits $\Delta$. Each passage cuts off $V$ on the left, while the edge midpoints, shown as black dots, are on the right. Thus both triangle signs and the signs contributed by edge crossings are all $+$.

\begin{figure}[htbp]
\centering
\begin{tikzpicture}[
  x=3.1cm,y=3.1cm,
  line cap=round,line join=round,
  every node/.style={font=\small,inner sep=1.2pt}
]
  \draw[chsix triangulation]
    (1,0)--(0,1)--(-1,1)--(-1,0)--(0,-1)--(1,-1)--cycle;
  \foreach \X/\Y in {1/0,0/1,-1/1,-1/0,0/-1,1/-1}
    \draw[chsix triangulation] (0,0)--(\X,\Y);
  \foreach \X/\Y in {.5/0,0/.5,-.5/.5,-.5/0,0/-.5,.5/-.5}
    \node[chsix point] at (\X,\Y) {};
  \node[chsix point] at (0,0) {};
  \node[below left=2pt] at (0,0) {$V$};
  \node at (.60,.28) {$\Delta$};
  \coordinate (P) at (.20,0);
  \coordinate (Q) at (0,.44);
  \draw[chsix red curve]
    plot[domain=-22:0,samples=13,variable=\t]
      ({(.20+.24*\t/450)*cos(\t)},{(.20+.24*\t/450)*sin(\t)});
  \draw[chsix red curve]
    plot[domain=450:474,samples=13,variable=\t]
      ({(.20+.24*\t/450)*cos(\t)},{(.20+.24*\t/450)*sin(\t)});
  \draw[draw=annred!80,line width=.9pt,dash pattern=on 2pt off 1.4pt,
    postaction={decorate},decoration={markings,
      mark=at position .10 with {\arrow{Stealth[length=1.5mm,width=1.05mm]}},
      mark=at position .35 with {\arrow{Stealth[length=1.5mm,width=1.05mm]}},
      mark=at position .62 with {\arrow{Stealth[length=1.5mm,width=1.05mm]}},
      mark=at position .89 with {\arrow{Stealth[length=1.5mm,width=1.05mm]}}}]
    plot[domain=0:450,samples=181,variable=\t]
      ({(.20+.24*\t/450)*cos(\t)},{(.20+.24*\t/450)*sin(\t)});
  \draw[chsix red curve,chsix forward at=.50]
    plot[domain=0:90,samples=46,variable=\t]
      ({(.20+.24*\t/90)*cos(\t)},{(.20+.24*\t/90)*sin(\t)});
  \node[chsix red point] at (P) {};
  \node[chsix red point] at (Q) {};
  \node[anchor=west] at (.29,.035) {$P$};
  \node[above left=2pt] at (Q) {$Q$};
  \foreach \angle/\radius in {45/.12,112.5/.17,157.5/.19,225/.22,292.5/.25,337.5/.28,405/.37}
    \node[text=annred,font=\scriptsize,inner sep=.2pt]
      at ({\radius*cos(\angle)},{\radius*sin(\angle)}) {$+$};
\end{tikzpicture}
\caption{Deleting a constant-sign interval}
\label{fig:same-sign-passages}
\end{figure}

Both passages through $\Delta$ enter through the bottom side and leave through the left side. Joining $P$ to $Q$ directly, as the solid red curve does, still gives sign $+$. Reverse orientation for $q=-$. In general, two passages with the same cut-off vertex and sign join the same two sides of $\Delta$ in the same direction.

In this case no enlargement is needed, and the direct connection also has sign $q$, so $M(w)\geq U_q$ applies. In either case the retained edge blocks are unchanged and $|\eta|\leq|\gamma|$. Together with the smaller crossing count, this contradicts lexicographic minimality and proves (1). Since self-intersections occur only between distinct passages, an arc passing through each triangle at most once has none.

Since $\mathcal A_0(A,B)$ is a subset of the generalized arcs over which we minimized, $|\gamma|\leq d(A,B)$. By (1), $\gamma\in\mathcal A_0(A,B)$, giving the reverse inequality and proving (3). Thus allowing self-intersections does not change the minimum GM length.
\end{proof}

We next minimize the choice of side at edge midpoints for a fixed triangle-passage sequence. Take a generalized arc with $n:=n(\gamma)\geq1$ edge crossings. List triangle passages and edge crossings in curve order, indexing them by $0,1,\ldots,2n$ starting with the triangle passage containing the initial endpoint:
\[
 \underset{0}{\text{triangle passage}},\quad
 \underset{1}{\text{edge crossing}},\quad
 \underset{2}{\text{triangle passage}},\quad\ldots,\quad
 \underset{2n-1}{\text{edge crossing}},\quad
 \underset{2n}{\text{triangle passage}}.
\]
Even positions are triangle passages, and odd positions are edge crossings. For odd $j=1,3,\ldots,2n-1$, let $i_j$ be the type of crossing $j$. If $k_{\sigma(i_j)}>0$, call the common sign assigned by rule (3)(i) or (ii) the \emph{midpoint sign} $b_j$. The edge block is $b_j^{k_{\sigma(i_j)}}$, consisting of $k_{\sigma(i_j)}$ copies of $b_j$. At a position with $k_{\sigma(i_j)}=0$, choose $b_j\in\{+,-\}$ arbitrarily for convenience, with empty edge block. Put $b=(b_1,b_3,\ldots,b_{2n-1})\in\{+,-\}^n$. Write $t_2,t_4,\ldots,t_{2n-2}$ for the interior triangle signs, indexed by their positions. The endpoint triangle signs are omitted from the transfer formula, so
\begin{equation}\label{eq:fresh-fixed-passage-length}
 F(b):=\mathbf1^{\mathsf T}
 U_{b_1}^{k_{\sigma(i_1)}}U_{t_2}U_{b_3}^{k_{\sigma(i_3)}}\cdots
 U_{t_{2n-2}}U_{b_{2n-1}}^{k_{\sigma(i_{2n-1})}}\mathbf1
\end{equation}
When $n=1$, this is $F(b)=\mathbf1^{\mathsf T}U_{b_1}^{k_{\sigma(i_1)}}\mathbf1$.

\begin{defi}\label{def:fresh-canonical-midpoints}\index{Canonical choice of midpoint signs}
For comparison, put $t_i:=0$ for even $i\leq0$ or $i\geq2n$. In particular, $t_0=t_{2n}=0$; these are markers for reaching an endpoint passage, not the endpoint triangle signs. Centered at the crossing in odd position $j$, call the pair of triangle signs
\[
 P_j(r):=(t_{j-(2r+1)},t_{j+(2r+1)})\qquad(r=0,1,\ldots)
\]
a \emph{comparison pair}. At $r=0$ we compare the immediately preceding and following signs $(t_{j-1},t_{j+1})$; each increase of $r$ moves one triangle outward on both sides. Define the \emph{comparison radius} by
\[
 R_j:=\min\{r\geq0\mid P_j(r)\notin\{(+,-),(-,+)\}\}
\]
This minimum exists because the sequence is finite. Call $P_j(R_j)$ the \emph{terminal pair}. If it is $(q,q),(0,q)$, or $(q,0)$, put $c_j:=q$; if it is $(0,0)$, either choice of $c_j$ is allowed. Call the resulting sequence $c=(c_1,c_3,\ldots,c_{2n-1})$ a \emph{canonical choice of midpoint signs}, and each $c_j$ a \emph{canonical sign}. A position whose terminal pair is not $(0,0)$ has a unique canonical sign and is called a \emph{forced position}.

The triangle signs are \emph{antisymmetric about $j$ up to radius $R$} if $P_j(r)\in\{(+,-),(-,+)\}$ for $0\leq r<R$. Thus the first $R$ signs on each side, taken outward from the center, are opposite at corresponding positions. The terminal pair $P_j(R_j)$ is not included in this condition.
\end{defi}

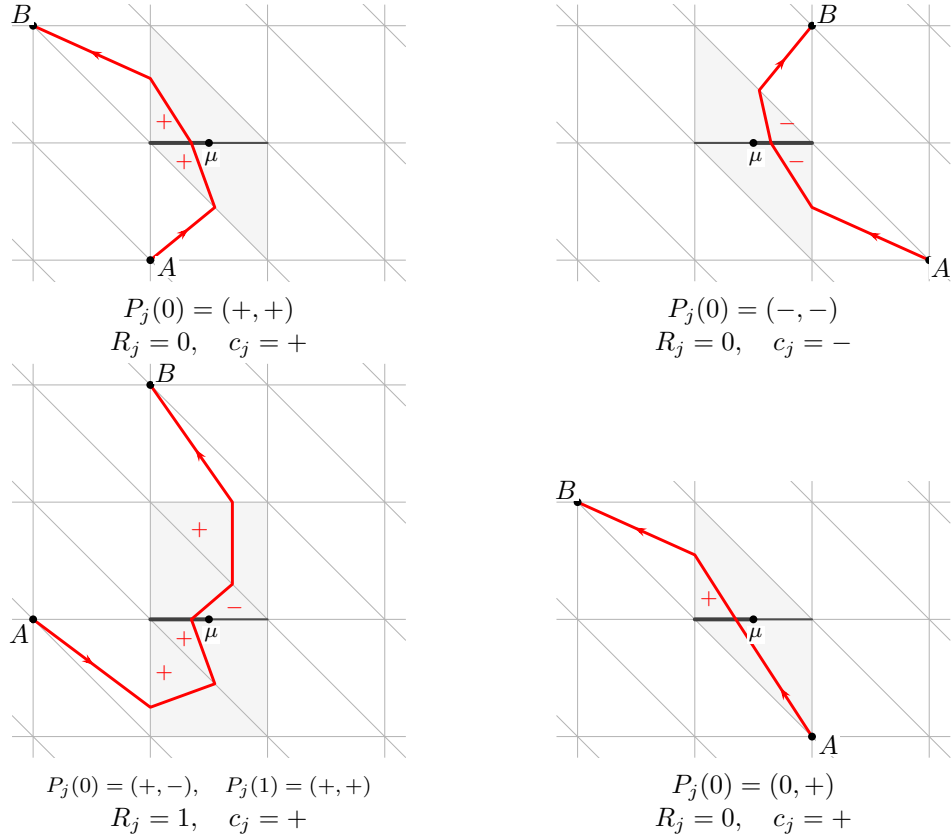
\begin{figure}[htbp]
\centering
\begin{tikzpicture}[
  x=1.55cm,y=1.55cm,
  line cap=round,line join=round,
  tri/.style={chsix triangulation},
  curve/.style={chsix red curve},
  standard/.style={draw=black!75,line width=1.5pt},
  sign/.style={text=annred,font=\scriptsize\bfseries,inner sep=.25pt},
  midpoint label/.style={font=\scriptsize,fill=white,inner sep=.7pt},
  endpoint label/.style={font=\small,fill=white,inner sep=.7pt},
  every node/.style={font=\small,inner sep=1.2pt}
]
\foreach \panel/\xoffset/\yoffset in {0/0/0,1/7.20/0,2/0/-6.30,3/7.20/-6.30}{
  \begin{scope}[xshift={\xoffset cm},yshift={\yoffset cm}]
    \path[fill=black!4] (0,0)--(1,-1)--(1,0)--(0,1)--cycle;
    \ifnum\panel=2
      \path[fill=black!4] (0,0)--(0,-1)--(1,-1)--cycle;
      \path[fill=black!4] (0,1)--(1,0)--(1,1)--cycle;
      \def\gridtop{2.18}
    \else
      \def\gridtop{1.18}
    \fi
    \begin{scope}
      \clip (-1.18,-1.18) rectangle (2.18,\gridtop);
      \foreach \i in {-1,0,1,2}{
        \draw[tri] (\i,-1.18)--(\i,\gridtop);
        \draw[tri] (-1.18,\i)--(2.18,\i);
      }
      \foreach \i in {-2,-1,0,1,2,3,4}{
        \draw[tri] (-1.18,{\i+1.18})--(2.18,{\i-2.18});
      }
    \end{scope}
    \draw[chsix distinguished edge] (0,0)--(1,0);
    \ifnum\panel=1
      \draw[standard] (.50,0)--(1,0);
    \else
      \draw[standard] (0,0)--(.50,0);
    \fi
    \ifcase\panel
      \draw[curve,chsix forward at=.33] (0,-1)--(.55,-.55)--(.35,0);
      \draw[curve,chsix forward at=.70] (.35,0)--(0,.55)--(-1,1);
      \node[sign] at (.29,-.16) {$+$};
      \node[sign] at (.12,.18) {$+$};
      \node[chsix point] at (0,-1) {};
      \node[chsix point] at (-1,1) {};
      \node[endpoint label] at (.14,-1.06) {$A$};
      \node[endpoint label] at (-1.10,1.10) {$B$};
      \node at (.50,-1.43) {$P_j(0)=(+,+)$};
      \node at (.50,-1.72) {$R_j=0,\quad c_j=+$};
    \or
      \draw[curve,chsix forward at=.33] (2,-1)--(1,-.55)--(.65,0);
      \draw[curve,chsix forward at=.70] (.65,0)--(.55,.45)--(1,1);
      \node[sign] at (.87,-.16) {$-$};
      \node[sign] at (.79,.16) {$-$};
      \node[chsix point] at (2,-1) {};
      \node[chsix point] at (1,1) {};
      \node[endpoint label] at (2.10,-1.07) {$A$};
      \node[endpoint label] at (1.13,1.10) {$B$};
      \node at (.50,-1.43) {$P_j(0)=(-,-)$};
      \node at (.50,-1.72) {$R_j=0,\quad c_j=-$};
    \or
      \draw[curve,chsix forward at=.27]
        (-1,0)--(0,-.75)--(.55,-.55)--(.35,0);
      \draw[curve,chsix forward at=.72]
        (.35,0)--(.70,.30)--(.70,1)--(0,2);
      \node[sign] at (.12,-.45) {$+$};
      \node[sign] at (.29,-.16) {$+$};
      \node[sign] at (.72,.10) {$-$};
      \node[sign] at (.42,.77) {$+$};
      \node[chsix point] at (-1,0) {};
      \node[chsix point] at (0,2) {};
      \node[endpoint label] at (-1.12,-.12) {$A$};
      \node[endpoint label] at (.13,2.10) {$B$};
      \node[font=\scriptsize] at (.50,-1.43)
        {$P_j(0)=(+,-),\quad P_j(1)=(+,+)$};
      \node at (.50,-1.72) {$R_j=1,\quad c_j=+$};
    \or
      \draw[curve,chsix forward at=.42] (1,-1)--(.35,0);
      \draw[curve,chsix forward at=.70] (.35,0)--(0,.55)--(-1,1);
      \node[sign] at (.12,.18) {$+$};
      \node[chsix point] at (1,-1) {};
      \node[chsix point] at (-1,1) {};
      \node[endpoint label] at (1.14,-1.07) {$A$};
      \node[endpoint label] at (-1.10,1.10) {$B$};
      \node at (.50,-1.43) {$P_j(0)=(0,+)$};
      \node at (.50,-1.72) {$R_j=0,\quad c_j=+$};
    \fi
    \node[chsix point] at (.50,0) {};
    \node[midpoint label] at (.52,-.13) {$\mu$};
  \end{scope}
}
\end{tikzpicture}
\caption{Examples of forced positions}
\label{fig:forced-midpoint-positions}
\end{figure}

Forced positions can also be detected geometrically. The vertex cut off by an interior passage is on the left for sign $+$ and on the right for sign $-$. Compare passages one at a time before and after crossing $j$, proceeding outward while their cut-off vertices are on opposite sides. The first pair with both on the left gives $c_j=+$, and the first with both on the right gives $c_j=-$. Even when examining the preceding portion, left and right are measured using the original curve orientation. If just one side reaches an endpoint passage, use the sign of the interior passage on the other side; if both sides reach endpoint passages simultaneously, the position is not forced. In Figure~\ref{fig:forced-midpoint-positions}, $\mu$ is the midpoint of the central crossed edge, and the thick black segment is the half-edge corresponding to the canonical sign. The red curve illustrates a crossing $j$ realizing that sign.

\begin{lemm}\label{lem:fresh-canonical-minimum}\index{Canonical choice of midpoint signs}
For \eqref{eq:fresh-fixed-passage-length}, every midpoint-sign sequence $b$ and canonical choice $c$ satisfy $F(c)\leq F(b)$.
\end{lemm}
\begin{proof}
We replace the entries of $b$ by their canonical signs $c_j$ one at a time, showing that $F$ never increases. To choose which entry to change, first establish two properties of canonical signs. The first is that if the triangle signs around crossing $j$ satisfy $t_{j-(2r+1)}=\bar t_{j+(2r+1)}$ for $0\leq r<R$, then crossed edges at equal distances from the center have equal weights. Explicitly,
\begin{equation}\label{eq:fresh-reflected-weights}
 k_{\sigma(i_{j-2d})}=k_{\sigma(i_{j+2d})}\qquad(1\leq d\leq R)
\end{equation}
For $R\geq1$, the two triangles adjacent to crossing $j$ correspond under the half-turn about the midpoint of their shared edge. Since $t_{j-1}=\bar t_{j+1}$, the edge used at crossing $j-2$ maps to the one used at crossing $j+2$. The half-turn preserves edge type, so their weights agree. Extending this correspondence outward one triangle at a time gives the equality for every $1\leq d\leq R$.

The second property is that, if $c_j=q$ is forced, then for each $1\leq d\leq R_j$, at least one of positions $j-2d,j+2d$ also has canonical sign $q$ and has comparison radius smaller than $R_j$. Put $R:=R_j$ and examine the two comparison procedures simultaneously. For $0\leq r<R-d$, antisymmetry about $j$ makes the compared pairs take the form $(\varepsilon,\eta)$ and $(\bar\eta,\bar\varepsilon)$, with $\varepsilon,\eta\in\{+,-\}$; no endpoint marker $0$ occurs in this range. If a first radius $r$ with $\varepsilon=\eta$ occurs here, it is the comparison radius at both positions. Their canonical signs are then $\varepsilon,\bar\varepsilon$, one of which is $q$. If opposite-sign pairs persist for all $0\leq r<R-d$, at $r=R-d$ the outer entries reach the first and second entries of $P_j(R)$. Put $(x,y):=P_j(R)$ and $z:=t_{j+2R-4d+1}$. The two comparison pairs are then
\[
 P_{j-2d}(R-d)=(x,z),\qquad P_{j+2d}(R-d)=(\bar z,y)
\]
Here $(x,y)$ is $(q,q),(0,q)$, or $(q,0)$, and $z\in\{+,-\}$. If $z=q$, the first pair $(x,z)$ forces $q$; if $z=\bar q$, the second pair $(\bar z,y)$ forces $q$. Its comparison radius is $R-d<R$, proving the second property.

Now compute the change in $F$ when only $b_j$ is varied. In \eqref{eq:fresh-fixed-passage-length}, let $X$ be the product before $U_{b_j}^{k_{\sigma(i_j)}}$ and $Y$ the product after it, using $E_2$ for an empty product. Thus $F(b)=\mathbf1^{\mathsf T}XU_{b_j}^{k_{\sigma(i_j)}}Y\mathbf1$. Write
$u:=\mathbf1^{\mathsf T}X=(u_1,u_2)$,
$v:=Y\mathbf1=(v_1,v_2)^{\mathsf T}$. With all other midpoint signs fixed, abbreviate the two values obtained by setting $b_j=+$ or $-$ as $F(b_j=+)$ and $F(b_j=-)$. Then
\begin{equation}\label{eq:fresh-one-midpoint-difference}
 F(b_j=+)-F(b_j=-)=k_{\sigma(i_j)}(u_2v_1-u_1v_2)
\end{equation}
Since $u_1,v_1>0$, when $k_{\sigma(i_j)}>0$ the right-hand side has the sign of $u_2/u_1-v_2/v_1$. Thus $b_j=+$ gives the smaller value when $u_2/u_1<v_2/v_1$, and $b_j=-$ when $u_2/u_1>v_2/v_1$. Compare these ratios through the products $u^{\mathsf T}=X^{\mathsf T}\mathbf1$ and $v=Y\mathbf1$. Since $U_q^{\mathsf T}=U_{\bar q}$, the factors of $X^{\mathsf T}$ correspond to the sign-reversed sequence taken from the center toward the left, while those of $Y$ correspond to the sequence taken toward the right. For a positive column $w=(w_1,w_2)^{\mathsf T}$,
\[
 \frac{(U_+w)_2}{(U_+w)_1}=1+\frac{w_2}{w_1},\qquad
 \frac{(U_-w)_2}{(U_-w)_1}=\frac{w_2}{w_1+w_2}
\]
If two positive columns have ratios $\alpha,\beta>0$, multiplication by $U_+$ changes them to $1+\alpha,1+\beta$, and multiplication by $U_-$ changes them to $\alpha/(1+\alpha),\beta/(1+\beta)$. The differences are
\[
 (1+\alpha)-(1+\beta)=\alpha-\beta,\qquad
 \frac{\alpha}{1+\alpha}-\frac{\beta}{1+\beta}
 =\frac{\alpha-\beta}{(1+\alpha)(1+\beta)}
\]
The denominators are positive, so in either case the sign of the difference is unchanged. Consequently, multiplying two positive columns on the left by the same matrix product preserves the order of their component ratios. We may therefore compare $X^{\mathsf T}$ and $Y$ from the beginning and remove matching factors successively. Comparing the ratios after applying the remaining products to $\mathbf1$ determines the order of the original ratios $u_2/u_1,v_2/v_1$. Here expand $U_q^a$ into $a$ factors $U_q$ and omit factors with exponent $0$.

Let $C$ be the full common initial product of $X^{\mathsf T}$ and $Y$, with $C=E_2$ if there are no matching initial factors. If both products still have factors immediately after $C$, these factors differ, so write
\[
 X^{\mathsf T}=C U_{\bar q}X_0,\qquad Y=C U_qY_0
\]
Here $q\in\{+,-\}$ is the sign of the first differing factor on the $Y$ side, and $X_0,Y_0$ are the remaining products, interpreted as $E_2$ if empty. The factor $U_{\bar q}$ on the $X^{\mathsf T}$ side was $U_q$ before transposition, so both pre-transposition signs agree and equal $q$. Let $\alpha,\beta>0$ be the second-to-first component ratios of the positive columns $X_0\mathbf1,Y_0\mathbf1$. Assume $k_{\sigma(i_j)}>0$.

If $q=+$, after removing $C$ the columns are $U_-X_0\mathbf1$ and $U_+Y_0\mathbf1$, whose ratios satisfy
\[
 \frac{(U_-X_0\mathbf1)_2}{(U_-X_0\mathbf1)_1}
 =\frac{\alpha}{1+\alpha}<1<1+\beta
 =\frac{(U_+Y_0\mathbf1)_2}{(U_+Y_0\mathbf1)_1}
\]
Multiplication by $C$ preserves this order, so $u_2/u_1<v_2/v_1$ also holds for the original columns. Thus the right-hand side of \eqref{eq:fresh-one-midpoint-difference} is negative, giving $F(b_j=+)<F(b_j=-)$. Choosing $b_j=+$ makes $F$ smaller.

If $q=-$, the columns after removing $C$ are $U_+X_0\mathbf1$ and $U_-Y_0\mathbf1$, so
\[
 \frac{(U_+X_0\mathbf1)_2}{(U_+X_0\mathbf1)_1}
 =1+\alpha>1>\frac{\beta}{1+\beta}
 =\frac{(U_-Y_0\mathbf1)_2}{(U_-Y_0\mathbf1)_1}
\]
Now $u_2/u_1>v_2/v_1$, so \eqref{eq:fresh-one-midpoint-difference} is positive and $F(b_j=-)<F(b_j=+)$. Thus choosing $b_j=-$ makes $F$ smaller. These calculations show that if the first differing factors have the same pre-transposition sign $q$, choosing $b_j=q$ is better. If removing the common factors leaves just one product empty, its column is $\mathbf1$, with ratio $1$. Let $q$ be the first pre-transposition sign on the other side; the same comparison shows that $b_j=q$ makes $F$ smaller. If both products are empty, the original ratios are equal and both choices of $b_j$ give the same $F$.

Use these comparisons to turn the whole sequence into its canonical choice. Starting from arbitrary $b$, we aim to obtain $c$ without increasing $F$. First, at positions with $k_{\sigma(i_j)}=0$, the factor is $E_2$, so set $b_j=c_j$ without changing $F$. Call a forced position with $b_j\neq c_j$ a \emph{disagreement position}. If any exist, choose one with smallest comparison radius $R_j$, and put $R:=R_j$, $q=c_j$. Its weight is positive because all zero-weight positions have already been corrected. Compare the associated $X^{\mathsf T}$ and $Y$ from the beginning. The left triangle sign $t_{j-(2r+1)}$ contributes $U_{t_{j-(2r+1)}}^{\mathsf T}$ to $X^{\mathsf T}$, to be compared with $U_{t_{j+(2r+1)}}$ in $Y$. For $0\leq r<R$, we have $t_{j-(2r+1)}=\bar t_{j+(2r+1)}$, so
\[
 U_{t_{j-(2r+1)}}^{\mathsf T}
 =U_{\overline{t_{j-(2r+1)}}}
 =U_{t_{j+(2r+1)}}.
\]
Thus, provided the closer portions agree, these two matrices can also be removed as common factors. Next compare the edge-crossing factors at $j-2d,j+2d$, where $d=r+1$. Their exponents agree by \eqref{eq:fresh-reflected-weights}. If the exponent is $0$, omit both blocks. If it is positive, the second property proved above says that at least one position has canonical sign $q$ and comparison radius smaller than $R$. It cannot be a disagreement position by the choice of $j$, so its actual midpoint sign is also $q$. The two midpoint signs are therefore either opposite or both equal to $q$. If opposite, the transposed factors agree; if both equal to $q$, the preceding calculation shows that $b_j=q$ makes $F$ smaller. If no strict comparison has occurred by radius $R$, the terminal pair $P_j(R)$ is $(q,q),(0,q)$, or $(q,0)$. The calculation for equal next signs $q$, or for one empty product, again shows that $b_j=q$ is better. Replace $b_j$ by $c_j$. Each such change decreases $F$ and removes one disagreement, so finitely many changes correct every forced position.

Finally consider a position with terminal pair $(0,0)$. If one exists, $n$ is odd and this position is necessarily the center $j=n$. The whole triangle-sign sequence is antisymmetric about $j$. The canonical-sign procedures at symmetric positions correspond under sign reversal, so $c_{j-2d}=\bar c_{j+2d}$ for $1\leq d\leq(j-1)/2$. These positions have already been corrected, and the corresponding edge weights agree. Hence the products on the two sides of the center satisfy $X^{\mathsf T}=Y$. The difference in \eqref{eq:fresh-one-midpoint-difference} is $0$, so the central sign may also be replaced by the specified $c_j$ without changing $F$. The final sequence is $c$, proving $F(c)\leq F(b)$.
\end{proof}

The preceding lemma determines midpoint signs minimizing GM length for a fixed triangle-passage sequence. We now show that a generalized arc with these signs can be chosen arbitrarily close to a shortest polygonal path for ordinary Euclidean length. The shortest path may pass through lattice points, whereas the generalized arc whose GM length we measure must avoid them, so the two are constructed separately.

We first arrange the passage triangles in order.

Let $\Delta_0,\Delta_2,\ldots,\Delta_{2n}$, with $n\geq1$, be the triangle-passage sequence of a generalized arc $\gamma$, and let $E_j$, for $j=1,3,\ldots,2n-1$, be the whole crossed triangulation edges, including their endpoints. Take a separate copy of each triangle and, for each odd $j$, glue $\Delta_{j-1}$ to $\Delta_{j+1}$ along $E_j$, identifying points with the same original planar position on that edge. Call the resulting space a \emph{closed triangle strip}, denoted by $S$. Here ``closed'' means that triangle edges and vertices are included. Even if the same planar triangle occurs repeatedly, use separate copies and glue only along the specified edges. Let $\pi_S\colon S\to\mathbb R^2$ return each copy to its original triangle. Distinct points of $S$ may therefore have the same image under $\pi_S$.

A closed triangle strip is also called a \emph{triangle strip}, or simply a \emph{strip}; a consecutive collection of its triangle copies is a \emph{substrip}. Vertices of $S$ mean the vertices corresponding to lattice points in its triangle copies, not edge midpoints. Edges not used for gluing are \emph{boundary edges}, and their endpoints are \emph{boundary vertices}. A \emph{closed disk} is a space homeomorphic to $\{(x,y)\mid x^2+y^2\leq1\}$, and will also be called a \emph{disk}.

For a curve in $S$, write $\ell_{\mathrm E}$ for its length obtained by summing ordinary Euclidean lengths within the triangles. Define the distance $d_S$ between two points as the infimum of $\ell_{\mathrm E}$ over curves joining them. For these lengths and distances, curves may run along the boundary and meet lattice vertices or edge midpoints. A curve formed by joining finitely many line segments end to end is a \emph{polygonal path}.

\begin{lemm}\label{lem:fresh-canonical-geodesic}\index{Generalized Markov distance!minimizing arcs}\index{Generalized Markov distance}
Construct the triangle strip $S$ from the passage sequence of a generalized arc $\gamma$ as above. Let $A,B$ be the vertices in the first and last triangle copies corresponding to its initial and final endpoints. Then the following hold.
\begin{enumerate}
\item There exists a polygonal path $g$ minimizing $\ell_{\mathrm E}$ among all curves in $S$ joining $A$ and $B$. Its image and orientation from $A$ to $B$ are unique, and it can bend only at boundary vertices.
\item For every canonical midpoint-sign sequence $c$ and every $\delta>0$, there exists an arc $\widehat\gamma_\delta$ in $S$ joining $A,B$ and following $g$ with the same passage sequence. Its planar image $\gamma_\delta:=\pi_S\circ\widehat\gamma_\delta$ is a generalized arc with midpoint signs $c$, and suitable parametrizations satisfy
\[
\sup_{t\in[0,1]}d_S\bigl(\widehat\gamma_\delta(t),g(t)\bigr)<\delta
\]
If $b$ is the midpoint-sign sequence of $\gamma$, then
\[
|\gamma_\delta|=F(c)\leq F(b)=|\gamma|.
\]
\end{enumerate}
\end{lemm}
\begin{proof}
First construct the shortest polygonal path in (1). By the generalized-arc condition, every nonendpoint passage enters and leaves its triangle through different sides. Adding triangle copies in passage order therefore replaces one boundary edge by two new edges. It follows that $S$ is a disk and all its lattice vertices lie on its boundary. Choose $x_j\in E_j$ for $j=1,3,\ldots,2n-1$, and consider
\[
\begin{aligned}
L(x_1,x_3,\ldots,x_{2n-1})
:={}&\bigl\|\pi_S(x_1)-\pi_S(A)\bigr\|_2\\
&+\sum_{a=1}^{n-1}\bigl\|\pi_S(x_{2a+1})-\pi_S(x_{2a-1})\bigr\|_2
+\bigl\|\pi_S(B)-\pi_S(x_{2n-1})\bigr\|_2
\end{aligned}
\]
where $\|(x,y)\|_2:=\sqrt{x^2+y^2}$. The points $A,x_1$ belong to $\Delta_0$, the points $x_{2a-1},x_{2a+1}$ to $\Delta_{2a}$ for $1\leq a<n$, and $x_{2n-1},B$ to $\Delta_{2n}$. Each pair can therefore be joined by a segment in that triangle, and $L$ is the Euclidean length of the polygonal path through the chosen points in order. Identifying each edge with $[0,1]$ makes $L$ continuous on $[0,1]^n$. Take a sequence approaching its infimum. Applying Theorem~\ref{thm:bolzano} successively to each coordinate yields a convergent subsequence; continuity shows that $L$ attains its minimum at the limit.

This minimum is a lower bound for the length of every $A$--$B$ curve in $S$, not just polygonal paths. Indeed, let $\eta$ be any such curve. Cutting the strip along $E_j$ separates an initial and a final portion whose intersection is exactly $E_j$. The final substrips $\Delta_{j+1}\cup\Delta_{j+3}\cup\cdots\cup\Delta_{2n}$ shrink as $j$ increases. Taking the first point at which $\eta$ reaches each gives points $x_1,x_3,\ldots,x_{2n-1}$ on $E_1,E_3,\ldots,E_{2n-1}$ in that order. Arrival times may coincide when several edges share a vertex. The lengths of the portions from $A$ to $x_1$, from $x_1$ to $x_3$, and so on, ending from $x_{2n-1}$ to $B$, are at least the straight-line distances between their endpoint images under $\pi_S$. These distances sum to $L(x_1,x_3,\ldots,x_{2n-1})$, giving
\[
\ell_{\mathrm E}(\eta)\geq L(x_1,x_3,\ldots,x_{2n-1})\geq\min L
\]
Since $\eta$ was arbitrary, $\min L$ is a lower bound for every $A$--$B$ curve. The polygonal path minimizing $L$ is therefore the required shortest path $g$. The same construction works for any two points of $S$, using the substrip between triangles containing them.

Every subpath of a shortest path is itself shortest between its endpoints, since shortening it would shorten the whole path. Suppose $g$ bends at a nonvertex point $p$ in the interior of a triangle or on an edge. A sufficiently small neighborhood of $p$ maps injectively and isometrically under $\pi_S$ to a planar disk, or a half-disk if $p$ is on a boundary edge. Choose points $x,y$ of $g$ on either side of $p$ in this neighborhood so that the portions from $x$ to $p$ and from $p$ to $y$ are segments. Both a disk and a half-disk contain the segment between any two of their points, so the segment $\pi_S(x)\pi_S(y)$ lifts to a segment in $S$. Since the direction changes at $p$,
\[
\bigl\|\pi_S(x)-\pi_S(y)\bigr\|_2
<\bigl\|\pi_S(x)-\pi_S(p)\bigr\|_2
+\bigl\|\pi_S(p)-\pi_S(y)\bigr\|_2
\]
Replacing the two segments by this segment shortens $g$, a contradiction. Thus bends occur only at boundary vertices. A positive-length portion returning to the same point could also be deleted. Since there are finitely many vertices, a shortest path is a finite polygonal path without repeated vertices.

To prove uniqueness, suppose two different shortest paths join the same points. Between a separation and a subsequent meeting, choose two subpaths with disjoint interiors that bound a closed region, called a \emph{bigon}. Its interior contains no lattice vertex of $S$. Regard the bigon as a polygon, including all bends of the subpaths among its vertices, and triangulate it, subdividing the triangles of $S$ as necessary. At each polygon vertex, let the interior angle $\alpha$ be the sum of the incident triangle angles inside the region, and call $\pi-\alpha$ the exterior angle. At the separation and reunion points the interior angles are positive, so their exterior angles sum to less than $2\pi$. At any other vertex, an interior angle less than $\pi$ would allow a shortening by a chord inside the bigon. Its exterior angle is therefore nonpositive. The total exterior angle is consequently less than $2\pi$, whereas summing the angles in the triangulation gives $2\pi$, a contradiction. All angles and notions of inside and outside are taken in $S$; portions overlapping under $\pi_S$ are counted separately. This proves (1).

For (2), we must choose crossing points on each $E_j$ arbitrarily close to $g$ and realizing $c_j$. First note that $g\cap E_j$ is a nonempty point or interval. It is nonempty because $E_j$ separates $A$ from $B$. The segment along $E_j$ between any two of its points realizes their planar distance and is therefore shortest in $S$. If $g$ passes through both points, uniqueness forces it to coincide with that segment between them. Thus $g\cap E_j$ is connected. For the same reason, after entering a substrip of consecutive triangles, $g$ cannot leave and later return: such a portion would cross an end edge of the substrip twice and would have to stay on the intervening edge segment.

Let $\mu$ be the midpoint of $E_j$, and let $V_j^{\mathrm L},V_j^{\mathrm R}$ be its left and right endpoints relative to the crossing direction. We want the crossing in $[V_j^{\mathrm L},\mu]$ if $c_j=+$ and in $[\mu,V_j^{\mathrm R}]$ if $c_j=-$. Call the corresponding interval the \emph{closed canonical half-edge}, and its interior the \emph{open canonical half-edge}. Crossing in the open half-edge gives sign $c_j$, so first show that $g\cap E_j$ meets the closed canonical half-edge.

As in the proof of Lemma~\ref{lem:fresh-canonical-minimum}, use the comparison radius $R_j$ from Definition~\ref{def:fresh-canonical-midpoints}. Put $R:=R_j$ and take the substrip formed by $R+1$ triangles on each side of the central edge:
\[
C:=\bigcup_{r=0}^{R}\bigl(\Delta_{j-(2r+1)}\cup\Delta_{j+(2r+1)}\bigr)
\]
For $r<R$, the signs $t_{j-(2r+1)}$ and $t_{j+(2r+1)}$ are opposite. Matching triangles outward from the center therefore identifies $\Delta_i$ with $\Delta_{2j-i}$ by the half-turn, or $180^\circ$ rotation, about $\mu$. Write $\rho\colon C\to C$ for this correspondence on $C$. The maps agree on glued edges and preserve curve lengths, so $\rho$ is an isometry. Again we match specified triangle copies rather than identifying points merely by their planar positions. The map $\rho$ exchanges the initial and final portions on the two sides of $E_j$, and on $E_j$ fixes $\mu$ and exchanges the endpoints.

Suppose $j$ is forced and $c_j=+$. In the initial half of $C$, compare the subpath of $g$ with the half-turn image of the final subpath, traversed in reverse. The outermost triangle on the initial side is $\Delta_{j-2R-1}$; let $P$ be its vertex opposite the side toward the center. The boundary of this half-strip consists of $E_j$, a boundary chain from $V_j^{\mathrm L}$ to $P$, and a boundary chain from $P$ to $V_j^{\mathrm R}$. Let $a,b$ be the points where the two subpaths begin traversing this half-strip. A half-turn preserves left and right relative to a curve, but reversing orientation exchanges the triangle signs $+$ and $-$. The sign definition therefore gives the following positions. If the first pair that is not opposite is $(+,+)$, then $a$ lies on the edge incident to $P$ in the first chain, and $b$ on the edge incident to $P$ in the second. For $(0,+)$, we have $a=P$; for $(+,0)$, we have $b=P$; the other point remains on its corresponding edge. Adding opposite-sign triangle pairs outward does not change this relative position: each pair corresponds under the half-turn, so the two boundary chains do not exchange sides. Figure~\ref{fig:canonical-entrance-positions} shows $R=0$. The solid red curve is the initial subpath, and the dashed red curve is the half-turn image of the final subpath traversed in reverse. Each panel is labeled by the original terminal pair.

\begin{figure}[H]
\centering
\begin{tikzpicture}[
  x=2.60cm,y=2.60cm,
  line cap=round,line join=round,
  tri/.style={chsix triangulation},
  curve/.style={chsix red curve},
  reflected/.style={chsix red curve,dash pattern=on 3.2pt off 2.2pt},
  every node/.style={font=\small,inner sep=1.2pt}
]
\foreach \panel in {0,1,2}{
  \begin{scope}[xshift={\panel*4.00cm}]
    \coordinate (P) at (1,-1);
    \ifcase\panel
      \coordinate (a) at (.50,-.50);
      \coordinate (b) at (1,{-1/3});
      \coordinate (x) at (.20,0);
      \coordinate (y) at (.80,0);
      \node at (.50,.39) {$(+,+)$};
    \or
      \coordinate (a) at (P);
      \coordinate (b) at (1,-.50);
      \coordinate (x) at ({1/3},0);
      \coordinate (y) at ({2/3},0);
      \node at (.50,.39) {$(0,+)$};
    \or
      \coordinate (a) at (.50,-.50);
      \coordinate (b) at (P);
      \coordinate (x) at ({1/3},0);
      \coordinate (y) at ({2/3},0);
      \node at (.50,.39) {$(+,0)$};
    \fi
    \draw[tri] (0,0)--(P)--(1,0);
    \draw[chsix distinguished edge] (0,0)--(1,0);
    \draw[curve,chsix forward at=.55] (a)--(x);
    \draw[reflected,chsix forward at=.58] (b)--(y);
    \foreach \point in {(0,0),(1,0),(.50,0),(P),(a),(b)}{
      \node[chsix point] at \point {};
    }
    \node[above=3pt] at (0,0) {$V_j^{\mathrm L}$};
    \node[above=3pt] at (1,0) {$V_j^{\mathrm R}$};
    \node[above=3pt] at (.50,0) {$\mu$};
    \node at (.50,-.10) {$E_j$};
    \ifcase\panel
      \node[below left=2pt] at (a) {$a$};
      \node[right=3pt] at (b) {$b$};
      \node[below=3pt] at (P) {$P$};
    \or
      \node[below=3pt] at (P) {$a=P$};
      \node[right=3pt] at (b) {$b$};
    \or
      \node[below left=2pt] at (a) {$a$};
      \node[below=3pt] at (P) {$b=P$};
    \fi
  \end{scope}
}
\end{tikzpicture}
\caption{The positions of $a$ and $b$}
\label{fig:canonical-entrance-positions}
\end{figure}
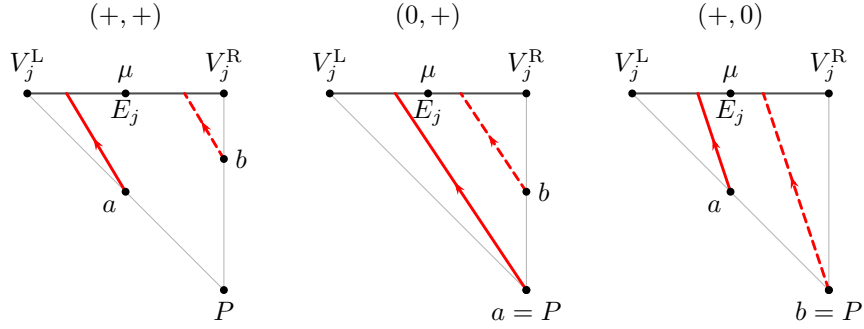

Suppose, for contradiction, that $g\cap E_j$ misses the closed canonical half-edge. Then $g$ reaches $E_j$ to the right of $\mu$, while the rotated subpath reaches it to the left. Let $x,y$ be their respective first points on $E_j$. If $a\in E_j$, then $a=V_j^{\mathrm L}$, contradicting the assumption that $g$ meets only the right half-edge. If $b\in E_j$, then $b=V_j^{\mathrm R}$, contradicting the analogous statement for the rotated subpath. Thus $a$ belongs to the boundary chain from $V_j^{\mathrm L}$ to $P$ excluding $V_j^{\mathrm L}$, and $b$ to the chain from $V_j^{\mathrm R}$ to $P$ excluding $V_j^{\mathrm R}$. If $a=b=P$, the subpaths meet there. Otherwise, along the boundary through $V_j^{\mathrm L},P,V_j^{\mathrm R}$, the four points occur in the order $a,b,x,y$, so the subpaths joining $a$ to $x$ and $b$ to $y$ intersect inside the disk. In either case, call an intersection $Z$. Figure~\ref{fig:canonical-crossing-contradiction} illustrates these arrangements schematically. The red curves represent the same subpaths as in Figure~\ref{fig:canonical-entrance-positions}, and the thick segment is the closed canonical half-edge.

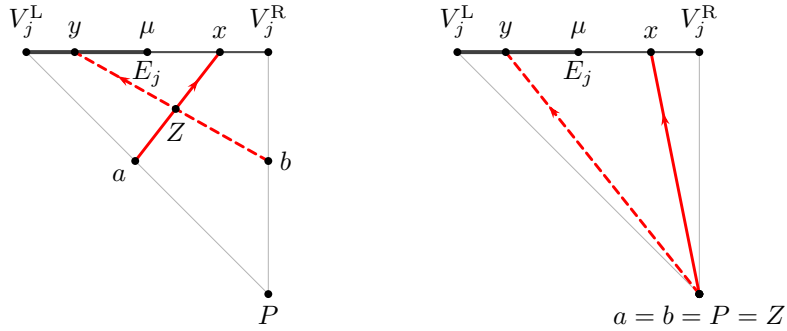
\begin{figure}[H]
\centering
\begin{tikzpicture}[
  x=3.20cm,y=3.20cm,
  line cap=round,line join=round,
  tri/.style={chsix triangulation},
  curve/.style={chsix red curve},
  reflected/.style={chsix red curve,dash pattern=on 3.2pt off 2.2pt},
  standard/.style={draw=black!75,line width=1.5pt},
  every node/.style={font=\small,inner sep=1.2pt}
]
\foreach \panel in {0,1}{
  \begin{scope}[xshift={\panel*5.70cm}]
    \coordinate (P) at (1,-1);
    \coordinate (x) at (.80,0);
    \coordinate (y) at (.20,0);
    \ifnum\panel=0
      \coordinate (a) at (.45,-.45);
      \coordinate (b) at (1,-.45);
      \coordinate (Z) at ({71/115},{-27/115});
    \else
      \coordinate (a) at (P);
      \coordinate (b) at (P);
      \coordinate (Z) at (P);
    \fi
    \draw[tri] (0,0)--(P)--(1,0);
    \draw[chsix distinguished edge] (0,0)--(1,0);
    \draw[standard] (0,0)--(.50,0);
    \draw[curve,chsix forward at=.75] (a)--(x);
    \draw[reflected,chsix forward at=.78] (b)--(y);
    \foreach \point in {(0,0),(1,0),(.50,0),(P),(a),(b),(x),(y),(Z)}{
      \node[chsix point] at \point {};
    }
    \node[above=4pt] at (0,0) {$V_j^{\mathrm L}$};
    \node[above=4pt] at (1,0) {$V_j^{\mathrm R}$};
    \node[above=4pt] at (.50,0) {$\mu$};
    \node[above=4pt] at (x) {$x$};
    \node[above=4pt] at (y) {$y$};
    \node at (.50,-.09) {$E_j$};
    \ifnum\panel=0
      \node[below left=2pt] at (a) {$a$};
      \node[right=3pt] at (b) {$b$};
      \node[below=4pt] at (Z) {$Z$};
      \node[below=3pt] at (P) {$P$};
    \else
      \node[below=3pt] at (P) {$a=b=P=Z$};
    \fi
  \end{scope}
}
\end{tikzpicture}
\caption{Intersection of the subpaths}
\label{fig:canonical-crossing-contradiction}
\end{figure}

The point $Z$ is not on $E_j$, and its half-turn image $\rho Z$ also lies on $g$. These are distinct points on opposite sides of $E_j$. The portion of $g$ from $Z$ to $\rho Z$ lies in $C$, and its half-turn image is another shortest path between the same points. By uniqueness the two coincide. Their intersection with $E_j$ is therefore a half-turn-invariant point or interval, so it contains $\mu$, contradicting the assumption. The case $c_j=-$ follows by exchanging left and right. If the terminal pair is $(0,0)$, then $C=S$ and $\rho$ exchanges $A,B$. Uniqueness makes $g$ itself half-turn-invariant, so again $g\cap E_j$ contains $\mu$. Thus either canonical sign gives a closed half-edge meeting $g$.

For each $j=1,3,\ldots,2n-1$, choose $x_j\in g\cap E_j$ in the closed canonical half-edge. These points occur along $g$ in passage order: a reversal would require leaving and returning to an edge separating the initial and final portions, contrary to its connected intersection with $g$. Several $x_j$ can coincide only at a common vertex of their edges. The points $A,x_1$ lie in $\Delta_0$, the points $x_{2a-1},x_{2a+1}$ in $\Delta_{2a}$ for $1\leq a<n$, and $x_{2n-1},B$ in $\Delta_{2n}$. The segment joining each pair in its triangle is shortest, so uniqueness identifies it with the corresponding portion of $g$.

Move each $x_j$ to $y_j$ in the open canonical half-edge, with $d_S(x_j,y_j)<\delta/2$, and join consecutive $y_j$ by segments within their triangles while keeping $A,B$ fixed. For every nonendpoint passage the two points lie on distinct sides. The first and last passages join a lattice endpoint to its opposite side. Thus each segment has interior in its triangle interior and avoids lattice vertices and edge midpoints. Perturb the $y_j$ slightly further within their allowed open half-edges so that distinct crossings of the same planar edge do not coincide. Even for multiple passages through the same planar triangle, their boundary endpoints can be moved independently to avoid tangencies, overlapping segments, and triple intersections. There are only finitely many forbidden conditions, and each point varies in an open interval, so all can be avoided simultaneously. The planar image then has only transverse intersections between distinct passages and is a generalized arc $\gamma_\delta$ with the original passage sequence and midpoint signs $c$.

Uniform closeness to $g$ follows from the displacements of the chosen points. Parametrize the segments of $g$ and $\widehat\gamma_\delta$ in each triangle linearly over the same parameter interval. If a segment of $g$ collapses to a point, keep it stationary on that interval. The distance between corresponding points is at most the larger endpoint displacement. Taking the additional perturbations sufficiently small therefore gives $d_S(\widehat\gamma_\delta(t),g(t))<\delta$ for every $t$. Finally, since $\gamma_\delta$ has the original passage sequence and canonical signs $c$, equation~\eqref{eq:fresh-fixed-passage-length} and Lemma~\ref{lem:fresh-canonical-minimum} give $|\gamma_\delta|=F(c)\leq F(b)=|\gamma|$. Euclidean minimization thus determines the geometric placement, while minimality of the canonical signs guarantees that GM length does not increase.
\end{proof}

\begin{exam}\label{exam:gm-length-canonical}
Consider the strip $S$ consisting of the five passage triangles in Example~\ref{exam:gm-length-calculation}, and put $D=(1,1)$. It is identified with the lightly shaded region in Figure~\ref{fig:gm-length-canonical}. Any curve in $S$ from $A=(0,1)$ to $B=(2,2)$ meets the edge from $(1,0)$ to $D$, say at $(1,y)$ with $0\leq y\leq1$. The distance from this point to $A$ is at least $1$, and its distance to $B$ is $\sqrt{1+(2-y)^2}\geq\sqrt2$. The curve therefore has length at least $1+\sqrt2$. The path $A,D,B$ attains this length, so it is the shortest polygonal path $g$ of Lemma~\ref{lem:fresh-canonical-geodesic}(1), bending at $D$.

The interior triangle signs are $(t_2,t_4,t_6)=(+,+,-)$. The comparison pairs at crossings $1,3,7$ are $(0,+),(+,+),(-,0)$, respectively, so $c_1=c_3=+$ and $c_7=-$. At crossing $5$, the pair $(+,-)$ is followed by $(+,0)$, giving $c_5=+$. Thus the unique canonical choice is
\[
 (c_1,c_3,c_5,c_7)=(+,+,+,-)
\]
For $0<\epsilon<\tfrac12$, realize it by the polygonal path $\gamma_\epsilon$ joining, in order,
\[
 A,\quad (1,1-\epsilon),\quad
 \left(\frac{3+\epsilon}{2},\frac{3-\epsilon}{2}\right),\quad B
\]
Its middle segment lies on $y=x-\epsilon$, and it has the same passage sequence as $\gamma$ and the stated canonical midpoint signs. As $\epsilon$ decreases, it approaches $g$ arbitrarily closely while avoiding $D$ and the edge midpoint $(\tfrac32,\tfrac32)$. The left panel of Figure~\ref{fig:gm-length-canonical} shows the shortest path $g$ in (1); the right panel shows the generalized arc in (2) in red for $\epsilon=\tfrac14$. The dashed curve on the right is $g$, black dots mark the crossed-edge midpoints, and the signs beside the red curve are the midpoint signs.

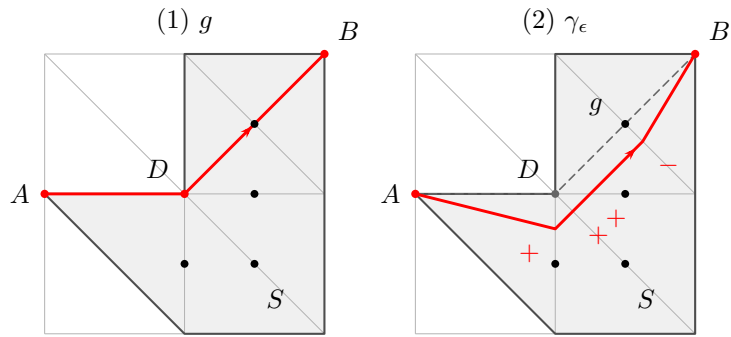
\begin{figure}[!htbp]
\centering
\begin{tikzpicture}[x=1.85cm,y=1.85cm,
  line cap=round,line join=round,every node/.style={font=\small},
  shortest line/.style={draw=black!60,densely dashed,line width=.6pt}]
\foreach \panel in {0,1}{
  \begin{scope}[shift={({2.65*\panel},0)}]
    \fill[black!6] (0,1)--(1,0)--(2,0)--(2,2)--(1,2)--(1,1)--cycle;
    \foreach \i in {0,1,2}{
      \draw[chsix triangulation] (\i,0)--(\i,2);
      \draw[chsix triangulation] (0,\i)--(2,\i);
    }
    \draw[chsix triangulation] (0,1)--(1,0) (0,2)--(2,0) (1,2)--(2,1);
    \draw[chsix distinguished edge] (0,1)--(1,0)--(2,0)--(2,2)--(1,2)--(1,1)--cycle;
    \node at (1.65,.25) {$S$};
    \ifnum\panel=0
      \node[chsix panel title] at (1,2.23) {(1) $g$};
      \draw[chsix red curve,chsix forward at=.70] (0,1)--(1,1)--(2,2);
      \node[chsix red point] at (1,1) {};
    \else
      \node[chsix panel title] at (1,2.23) {(2) $\gamma_\epsilon$};
      \draw[shortest line] (0,1)--(1,1)--(2,2);
      \node[circle,fill=black!60,inner sep=1.05pt] at (1,1) {};
      \draw[chsix red curve,chsix forward at=.70]
        (0,1)--(1,.75)--(1.625,1.375)--(2,2);
      \node[text=annred,below left=2pt] at (1,.75) {$+$};
      \node[text=annred,below right=2pt] at (1.125,.875) {$+$};
      \node[text=annred,below right=2pt] at (1.25,1) {$+$};
      \node[text=annred,below right=2pt] at (1.625,1.375) {$-$};
      \node[above left=2pt] at (1.45,1.45) {$g$};
    \fi
    \foreach \P in {(1,.5),(1.5,.5),(1.5,1),(1.5,1.5)}
      \node[chsix point] at \P {};
    \node[chsix red point,label={left:$A$}] at (0,1) {};
    \node[chsix red point,label={above right:$B$}] at (2,2) {};
    \node[above left=2pt] at (1,1) {$D$};
  \end{scope}
}
\end{tikzpicture}
\caption{The shortest polygonal path and an arc with canonical midpoint signs}
\label{fig:gm-length-canonical}
\end{figure}

The edge weights are $1,2,1,2$ in order, so the transfer formula gives
\[
 |\gamma_\epsilon|=\mathbf1^{\mathsf T}U_+^6U_-^3\mathbf1
 =\mathbf1^{\mathsf T}\begin{bmatrix}1&3\\6&19\end{bmatrix}\mathbf1=29
\]
which is smaller than the original length $|\gamma|=101$.
\end{exam}

Henceforth apply Lemma~\ref{lem:fresh-canonical-geodesic} to a minimizing arc with crossings obtained from Lemma~\ref{lem:fresh-minimizer}. Since each triangle occurs only once, reconnecting within triangles by segments creates no self-intersections. We may therefore replace the arc by $\gamma_\delta$ from Lemma~\ref{lem:fresh-canonical-geodesic}(2). Its crossing count is unchanged, and minimality makes its GM length unchanged as well. When $g$ passes through a boundary vertex, the generalized arc avoids that vertex and follows the corresponding consecutive triangles once each, in passage order. Statements about planar positions and directions concern $\pi_S\circ g$, while passages through different copies of the same planar vertex remain distinct in $S$.

The next lemma shows that if the polygonal path bends at an intermediate lattice point, or continues straight while the generalized arc changes its detour side, changing the preceding detour side produces a smaller GM length.

\begin{samepage}
\begin{lemm}\label{lem:fresh-first-bend}
Let $A,B$ be distinct lattice points with no lattice point on the open segment $AB$. Let $\gamma$ be a generalized arc joining them and satisfying (1),(2) of Lemma~\ref{lem:fresh-minimizer}. Assume also that $\gamma$ follows a shortest polygonal path by the construction of Lemma~\ref{lem:fresh-canonical-geodesic}(2). Suppose this path proceeds straight from $A$ to an intermediate lattice point $D$. List all lattice points on $AD$ in order as
\[
 A=Q_0,Q_1,\ldots,Q_n=D,\qquad Q_j=A+jv,\qquad n\geq1
\]
Here $v=Q_1-Q_0$, and each open segment $Q_{j-1}Q_j$ has no lattice point. Assume $\gamma$ makes a half-turn on the same side at every $Q_1,\ldots,Q_{n-1}$. If necessary, exchange left and right by the reflection $(x,y)\mapsto(y,x)$ so that this side is the right. If $n=1$, there is no intermediate lattice point, and we designate the detour side at $D$ as the right.

If the detour side changes to the left at $D$, or if it remains on the right and the angle traversed around $D$ is greater than $\pi$ but less than $2\pi$, then there exists a generalized arc from $A$ to $B$ with GM length strictly smaller than that of $\gamma$. The first case includes a straight polygonal path whose detour side alone changes. The angle is measured from the ray pointing back toward the incoming segment to the outgoing ray, in the order of triangles traversed by $\gamma$ around $D$.
\end{lemm}
\end{samepage}
\begin{proof}
We divide the proof into four cases.
\begin{enumerate}
\item\label{case:fresh-bend-left}
$AD$ does not follow triangulation edges, and the detour at $D$ is to the left.
\item\label{case:fresh-bend-right}
$AD$ does not follow triangulation edges, and the detour to the right of $D$ exceeds a half-turn, excluding case~(\ref{case:fresh-bend-exit-edge}).
\item\label{case:fresh-bend-along-edge}
$AD$ follows triangulation edges, excluding case~(\ref{case:fresh-bend-exit-edge}).
\item\label{case:fresh-bend-exit-edge}
The detour to the right of $D$ exceeds a half-turn, the outgoing polygonal path follows a triangulation edge, and during the detour $\gamma$ crosses that edge incident to $D$.
\end{enumerate}
Figure~\ref{fig:fresh-bend-cases} illustrates the four cases. Dashed lines indicate the shortest polygonal path and red curves the generalized arc; the red point in (\ref{case:fresh-bend-exit-edge}) marks the crossing of the edge in the outgoing direction.
\begin{figure}[!htbp]
\centering
\begin{tikzpicture}[x=1.45cm,y=1.45cm,
  every node/.style={font=\small},
  bend polyline/.style={draw=black!65,densely dashed,line width=.65pt},
  bend arrow/.style={postaction={decorate},decoration={markings,
    mark=at position .55 with {\arrow{Stealth[length=1.3mm,width=.9mm]}}}}]
\def\bendradius{.22}
\foreach \sx/\sy/\caselabel in {
  0/0/case:fresh-bend-left,
  4.4/0/case:fresh-bend-right,
  0/-4.2/case:fresh-bend-along-edge,
  4.4/-4.2/case:fresh-bend-exit-edge}{
  \begin{scope}[shift={(\sx,\sy)}]
    \begin{scope}
      \clip (-2.16,-2.16) rectangle (.96,1.34);
      \foreach \i in {-2,-1,0,1}{
        \draw[chsix triangulation] (\i,-3)--(\i,2);
        \draw[chsix triangulation] (-3,\i)--(2,\i);
      }
      \foreach \i in {-4,-3,...,3}{
        \draw[chsix triangulation] (-3,{\i+3})--(2,{\i-2});
      }
    \end{scope}
    \node at (-.6,1.62) {(\ref{\caselabel})};
  \end{scope}
}
\begin{scope}
  \draw[bend polyline] (-2,-2)--(.88,.88);
  \draw[chsix red curve] (-2,-2)--({-1+\bendradius*cos(235)},{-1+\bendradius*sin(235)});
  \draw[chsix red curve,bend arrow]
    ({-1+\bendradius*cos(235)},{-1+\bendradius*sin(235)})
    arc[start angle=235,end angle=395,radius=\bendradius];
  \draw[chsix red curve]
    ({-1+\bendradius*cos(395)},{-1+\bendradius*sin(395)})
    --({\bendradius*cos(205)},{\bendradius*sin(205)});
  \draw[chsix red curve,bend arrow]
    ({\bendradius*cos(205)},{\bendradius*sin(205)})
    arc[start angle=205,end angle=55,radius=\bendradius];
  \draw[chsix red curve,-{Stealth[length=1.65mm,width=1.10mm]}]
    ({\bendradius*cos(55)},{\bendradius*sin(55)})--(.78,.98);
  \node[chsix point,label={left:$A$}] at (-2,-2) {};
  \node[chsix point,label={above left:$Q_1$}] at (-1,-1) {};
  \node[chsix point,label={below right:$D$}] at (0,0) {};
\end{scope}
\begin{scope}[shift={(4.4,0)}]
  \draw[bend polyline] (-2,-2)--(0,0)--(.335,1.25);
  \draw[chsix red curve] (-2,-2)--({-1+\bendradius*cos(235)},{-1+\bendradius*sin(235)});
  \draw[chsix red curve,bend arrow]
    ({-1+\bendradius*cos(235)},{-1+\bendradius*sin(235)})
    arc[start angle=235,end angle=395,radius=\bendradius];
  \draw[chsix red curve]
    ({-1+\bendradius*cos(395)},{-1+\bendradius*sin(395)})
    --({\bendradius*cos(235)},{\bendradius*sin(235)});
  \draw[chsix red curve,bend arrow]
    ({\bendradius*cos(235)},{\bendradius*sin(235)})
    arc[start angle=235,end angle=425,radius=\bendradius];
  \draw[chsix red curve,-{Stealth[length=1.65mm,width=1.10mm]}]
    ({\bendradius*cos(425)},{\bendradius*sin(425)})--(.47,1.25);
  \node[chsix point,label={left:$A$}] at (-2,-2) {};
  \node[chsix point,label={above left:$Q_1$}] at (-1,-1) {};
  \node[chsix point,label={above left:$D$}] at (0,0) {};
\end{scope}
\begin{scope}[shift={(0,-4.2)}]
  \draw[bend polyline] (-2,-1)--(0,-1)--(.59,1.2);
  \draw[chsix red curve] (-2,-1)--({-1+\bendradius*cos(210)},{-1+\bendradius*sin(210)});
  \draw[chsix red curve,bend arrow]
    ({-1+\bendradius*cos(210)},{-1+\bendradius*sin(210)})
    arc[start angle=210,end angle=330,radius=\bendradius];
  \draw[chsix red curve]
    ({-1+\bendradius*cos(330)},{-1+\bendradius*sin(330)})
    --({\bendradius*cos(210)},{-1+\bendradius*sin(210)});
  \draw[chsix red curve,bend arrow]
    ({\bendradius*cos(210)},{-1+\bendradius*sin(210)})
    arc[start angle=210,end angle=425,radius=\bendradius];
  \draw[chsix red curve,-{Stealth[length=1.65mm,width=1.10mm]}]
    ({\bendradius*cos(425)},{-1+\bendradius*sin(425)})--(.73,1.2);
  \node[chsix point,label={above left:$A$}] at (-2,-1) {};
  \node[chsix point,label={above:$Q_1$}] at (-1,-1) {};
  \node[chsix point,label={above left:$D$}] at (0,-1) {};
\end{scope}
\begin{scope}[shift={(4.4,-4.2)}]
  \draw[chsix distinguished edge] (0,0)--(0,1);
  \draw[bend polyline] (-2,-2)--(0,0)--(0,1.25);
  \draw[chsix red curve] (-2,-2)--({-1+\bendradius*cos(235)},{-1+\bendradius*sin(235)});
  \draw[chsix red curve,bend arrow]
    ({-1+\bendradius*cos(235)},{-1+\bendradius*sin(235)})
    arc[start angle=235,end angle=395,radius=\bendradius];
  \draw[chsix red curve]
    ({-1+\bendradius*cos(395)},{-1+\bendradius*sin(395)})
    --({\bendradius*cos(235)},{\bendradius*sin(235)});
  \draw[chsix red curve,bend arrow]
    ({\bendradius*cos(235)},{\bendradius*sin(235)})
    arc[start angle=235,end angle=465,radius=\bendradius];
  \draw[chsix red curve,-{Stealth[length=1.65mm,width=1.10mm]}]
    ({\bendradius*cos(465)},{\bendradius*sin(465)})--(-.11,1.25);
  \node[chsix red point] at (0,\bendradius) {};
  \node[chsix point,label={left:$A$}] at (-2,-2) {};
  \node[chsix point,label={above left:$Q_1$}] at (-1,-1) {};
  \node[chsix point,label={above left:$D$}] at (0,0) {};
\end{scope}
\end{tikzpicture}
\caption{The four cases at the lattice point $D$}
\label{fig:fresh-bend-cases}
\end{figure}
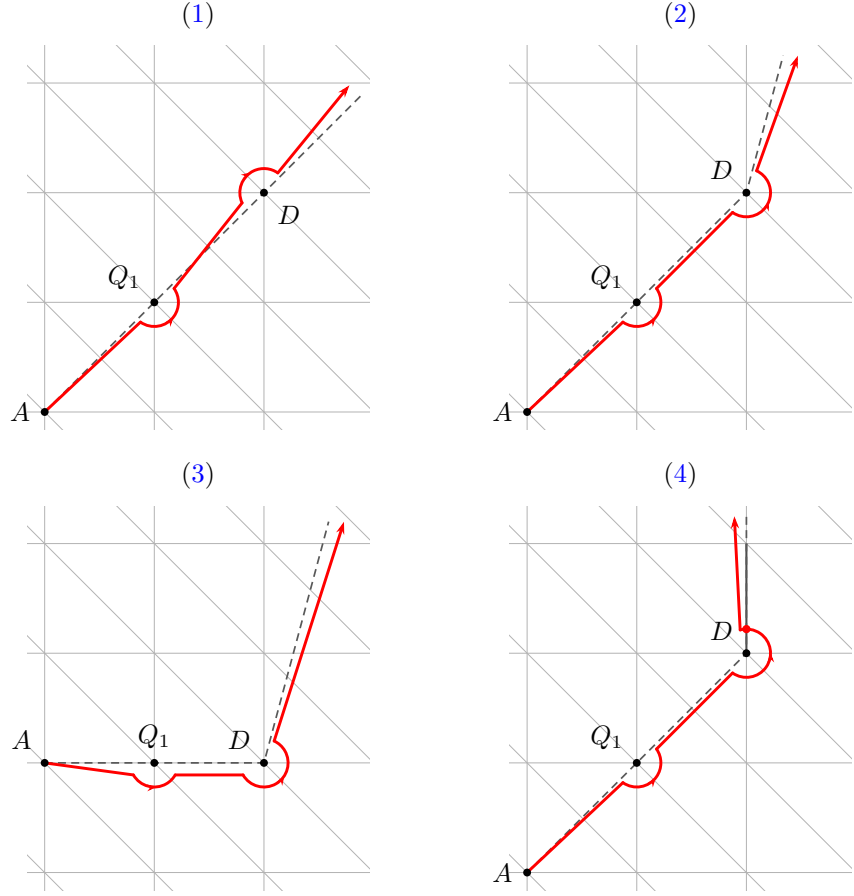
We first describe the strategy and calculations common to all cases. Replace the right detours up to $D$ by left detours and join to the original tail near $D$, showing that GM length decreases. If the connection creates consecutive crossings of the same triangulation edge or violates the condition that an endpoint joins its opposite side, remove the offending portions before forming a generalized arc. Perturb the retained crossings and the portions inside triangles slightly, staying on the same side of every positive-weight midpoint. Each passage can be made injective, with only finitely many transverse intersections between distinct passages, without changing the sign sequence. Figure~\ref{fig:fresh-left-detour-splice} illustrates case~(\ref{case:fresh-bend-left}). The original arc $\gamma$ is on the left and the new arc $\gamma'$ on the right; it joins the original tail at the red point in the shaded triangle. The portion after that point and the dashed shortest polygonal path are identical in both panels.
\begin{figure}[htbp]
\centering
\begin{tikzpicture}[x=1.65cm,y=1.65cm,
  every node/.style={font=\small},
  splice polyline/.style={draw=black!65,densely dashed,line width=.65pt},
  splice arrow/.style={postaction={decorate},decoration={markings,
    mark=at position .55 with {\arrow{Stealth[length=1.3mm,width=.9mm]}}}}]
\def\spliceradius{.22}
\foreach \panel in {0,1}{
  \begin{scope}[shift={({4.4*\panel},0)}]
    \begin{scope}
      \clip (-2.16,-2.16) rectangle (.96,1.12);
      \fill[black!6] (-1,0)--(0,-1)--(0,0)--cycle;
      \foreach \i in {-2,-1,0,1}{
        \draw[chsix triangulation] (\i,-3)--(\i,2);
        \draw[chsix triangulation] (-3,\i)--(2,\i);
      }
      \foreach \i in {-4,-3,...,3}{
        \draw[chsix triangulation] (-3,{\i+3})--(2,{\i-2});
      }
    \end{scope}
    \coordinate (A) at (-2,-2);
    \coordinate (Q) at (-1,-1);
    \coordinate (D) at (0,0);
    \coordinate (splice) at ({\spliceradius*cos(205)},{\spliceradius*sin(205)});
    \draw[splice polyline] (A)--(.88,.88);
    \ifnum\panel=0
      \def\startangle{235}
      \def\endangle{395}
    \else
      \def\startangle{215}
      \def\endangle{55}
    \fi
    \draw[chsix red curve]
      (A)--({-1+\spliceradius*cos(\startangle)},{-1+\spliceradius*sin(\startangle)});
    \draw[chsix red curve,splice arrow]
      ({-1+\spliceradius*cos(\startangle)},{-1+\spliceradius*sin(\startangle)})
      arc[start angle=\startangle,end angle=\endangle,radius=\spliceradius];
    \draw[chsix red curve]
      ({-1+\spliceradius*cos(\endangle)},{-1+\spliceradius*sin(\endangle)})--(splice);
    \draw[chsix red curve,splice arrow]
      (splice) arc[start angle=205,end angle=55,radius=\spliceradius];
    \draw[chsix red curve,-{Stealth[length=1.65mm,width=1.10mm]}]
      ({\spliceradius*cos(55)},{\spliceradius*sin(55)})--(.78,.98);
    \node[chsix point,label={left:$A$}] at (A) {};
    \node[chsix point,label={below right:$D$}] at (D) {};
    \ifnum\panel=0
      \node[chsix point,label={above left:$Q_1$}] at (Q) {};
      \node[text=annred,fill=white,inner sep=1pt] at (-.65,-1.22) {$\gamma$};
    \else
      \node[chsix point,label={below right:$Q_1$}] at (Q) {};
      \node[text=annred,fill=white,inner sep=1pt] at (-1.38,-.78) {$\gamma'$};
      \node[chsix red point] at (splice) {};
    \fi
  \end{scope}
}
\draw[chsix transition arrow] (1.26,-.50)--(1.92,-.50);
\end{tikzpicture}
\caption{A left detour and its connection}
\label{fig:fresh-left-detour-splice}
\end{figure}
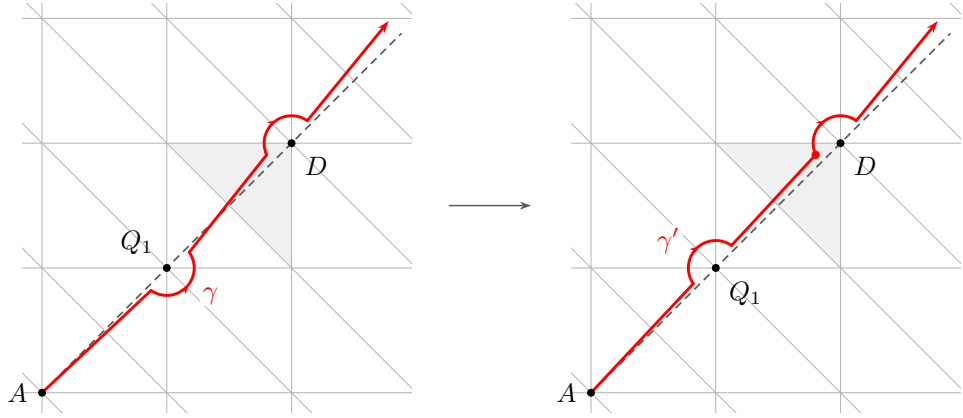

Splitting the interior word at the connection as $w=w'w''$ gives
\[
|\gamma|=\bigl(\mathbf1^{\mathsf T}M(w')\bigr)\bigl(M(w'')\mathbf1\bigr)
\]
Thus GM length is the product of the row for the initial portion and the column for the final portion. It suffices to compare the rows for two initial portions joining the same tail.

A reflection exchanging all signs preserves GM length: with $J=\begin{bsmallmatrix}0&1\\1&0\end{bsmallmatrix}$, we have $JU_+J=U_-$ and $J\mathbf1=\mathbf1$. When edge types are exchanged, exchange their weights as well. By Lemma~\ref{lem:same-length-lemma0}, let their common length be
\[
 N:=|\gamma^R_{AD}|=|\gamma^L_{AD}|
\]
If these arcs have crossings, let $r_0,t_0$ be $\mathbf1^{\mathsf T}$ multiplied by their respective interior-word products. Thus $r_0\mathbf1=t_0\mathbf1=N$. If there are no crossings, put $N=1$ and compute directly from the endpoint rule without using $r_0,t_0$.

\smallskip\noindent
\textbf{(\ref{case:fresh-bend-left}) Nonedge direction along $AD$, followed by a left detour at $D$.}
This case has $n\geq2$, since for $n=1$ the detour at $D$ was designated right. The integer linear transformation $T(x,y)=(-y,x+y)$ has determinant $1$ and cycles the six triangulation directions. Figure~\ref{fig:fresh-six-directions} shows their vectors $v_1,\ldots,v_6$ on the left and their images $Tv_1,\ldots,Tv_6$ on the right. Repeatedly applying $T$ and permuting the weights accordingly, we may assume $v=(q,p)$ with $p,q>0$. Since $T$ preserves lines and half-planes, it preserves the local shortest-path condition that the angle around a lattice point is at least $\pi$. It does not preserve Euclidean lengths themselves, so shortest-path comparisons use the lengths before transformation.

\begin{figure}[htbp]
\centering
\begin{tikzpicture}[x=1.55cm,y=1.55cm,
  every node/.style={font=\small},
  direction vector/.style={chsix red curve,-{Stealth[length=1.7mm,width=1.2mm]}}]
\foreach \panel in {0,1}{
  \begin{scope}[shift={({4.4*\panel},0)}]
    \begin{scope}
      \clip (-1.12,-1.12) rectangle (1.12,1.12);
      \foreach \i in {-1,0,1}{
        \draw[chsix triangulation] (\i,-2)--(\i,2);
        \draw[chsix triangulation] (-2,\i)--(2,\i);
      }
      \foreach \s in {-2,-1,...,2}{
        \draw[chsix triangulation] (-2,{\s+2})--(2,{\s-2});
      }
    \end{scope}
    \coordinate (O) at (0,0);
    \coordinate (v1) at (1,0);
    \coordinate (v2) at (0,1);
    \coordinate (v3) at (-1,1);
    \coordinate (v4) at (-1,0);
    \coordinate (v5) at (0,-1);
    \coordinate (v6) at (1,-1);
    \foreach \i in {1,...,6}{
      \draw[direction vector] (O)--(v\i);
    }
    \node[chsix point] at (O) {};
    \foreach \i/\preimage/\position in
      {1/6/right,2/1/above,3/2/above left,4/3/left,5/4/below,6/5/below right}{
      \ifnum\panel=0
        \node[\position=3pt] at (v\i) {$v_{\i}$};
      \else
        \node[\position=3pt] at (v\i) {$Tv_{\preimage}$};
      \fi
    }
  \end{scope}
}
\draw[chsix transition arrow] (1.55,0)--node[above=3pt] {$T$} (2.70,0);
\end{tikzpicture}
\caption{Cycling the six directions by $T$}
\label{fig:fresh-six-directions}
\end{figure}
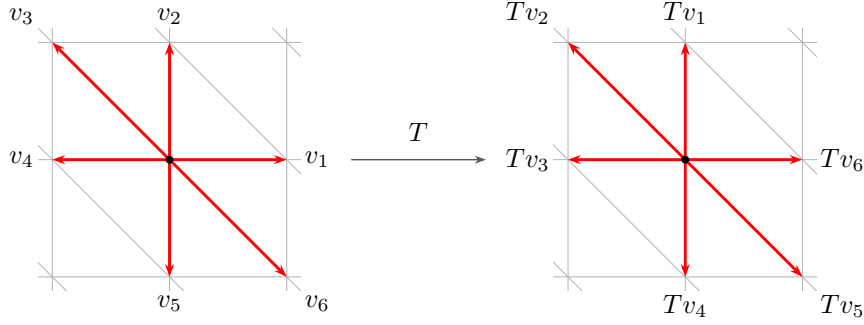

Let the comparison arc $\delta$ be the pure left push-off from $A$ to $D$. In the common triangle immediately before $D$, join the initial portion of $\delta$ to the tail of $\gamma$ to obtain an arc $\gamma'$ from $A$ to $B$. Let $\gamma''$ be the other reconnection, from the initial portion of $\gamma$ to $D$. Both initial portions enter this triangle through the same edge; $\gamma$ leaves through a different edge, while $\delta$ ends at the opposite vertex $D$. Thus reconnection does not create a return crossing of the entry edge. Figure~\ref{fig:fresh-four-arcs} shows $\gamma,\delta,\gamma',\gamma''$ in red, blue, green, and yellow, slightly separating overlapping portions. The figure depicts the case $A=Q_0,Q_1,D=Q_2$, with dashed lines for the corresponding polygonal path. The thick side of the shaded triangle is the common entry edge, and the black point marks the reconnection.
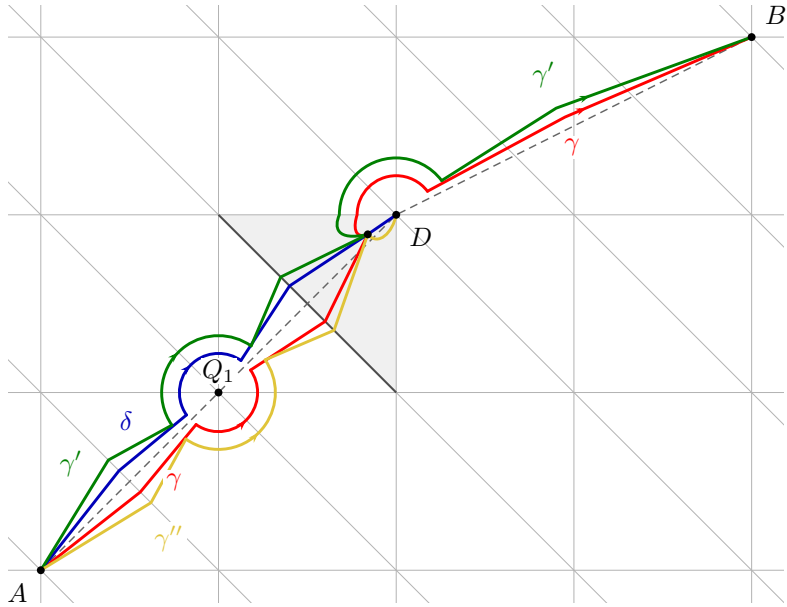
\begin{figure}[htbp]
\centering
\begin{tikzpicture}[x=2.35cm,y=2.35cm,
  every node/.style={font=\small},
  overview arrow/.style={postaction={decorate},decoration={markings,
    mark=at position .50 with {\arrow{Stealth[length=1.4mm,width=.95mm]}}}},
  arc label/.style={fill=white,inner sep=1.2pt}]
\begin{scope}
  \clip (-2.18,-2.18) rectangle (2.18,1.18);
  \fill[black!6] (-1,0)--(0,-1)--(0,0)--cycle;
  \foreach \i in {-2,-1,...,2}{
    \draw[chsix triangulation] (\i,-3)--(\i,2);
  }
  \foreach \i in {-2,-1,0,1}{
    \draw[chsix triangulation] (-3,\i)--(3,\i);
  }
  \foreach \s in {-4,-3,...,3}{
    \draw[chsix triangulation] (-3,{\s+3})--(3,{\s-3});
  }
\end{scope}
\draw[chsix distinguished edge] (-1,0)--(0,-1);
\coordinate (A) at (-2,-2);
\coordinate (Q) at (-1,-1);
\coordinate (D) at (0,0);
\coordinate (B) at (2,1);
\coordinate (join) at (-.16,-.11);
\draw[draw=black!60,densely dashed,line width=.55pt] (A)--(D)--(B);
\draw[chsix red curve]
  (A)--(-1.44,-1.56)--({-1+.22*cos(235)},{-1+.22*sin(235)});
\draw[chsix red curve,overview arrow]
  ({-1+.22*cos(235)},{-1+.22*sin(235)})
  arc[start angle=235,end angle=395,radius=.22];
\draw[chsix red curve]
  ({-1+.22*cos(395)},{-1+.22*sin(395)})--(-.4,-.6)--(join);
\draw[chsix red curve,draw=blue!72!black]
  (A)--(-1.56,-1.44)--({-1+.22*cos(215)},{-1+.22*sin(215)});
\draw[chsix red curve,draw=blue!72!black,overview arrow]
  ({-1+.22*cos(215)},{-1+.22*sin(215)})
  arc[start angle=215,end angle=55,radius=.22];
\draw[chsix red curve,draw=blue!72!black]
  ({-1+.22*cos(55)},{-1+.22*sin(55)})--(-.6,-.4)--(join)--(D);
\draw[chsix red curve,draw=green!50!black]
  (A)--(-1.62,-1.38)--({-1+.32*cos(215)},{-1+.32*sin(215)});
\draw[chsix red curve,draw=green!50!black,overview arrow]
  ({-1+.32*cos(215)},{-1+.32*sin(215)})
  arc[start angle=215,end angle=55,radius=.32];
\draw[chsix red curve,draw=green!50!black]
  ({-1+.32*cos(55)},{-1+.32*sin(55)})--(-.65,-.35)--(join);
\draw[chsix red curve,draw=yellow!80!orange!85!black]
  (A)--(-1.38,-1.62)--({-1+.32*cos(235)},{-1+.32*sin(235)});
\draw[chsix red curve,draw=yellow!80!orange!85!black,overview arrow]
  ({-1+.32*cos(235)},{-1+.32*sin(235)})
  arc[start angle=235,end angle=395,radius=.32];
\draw[chsix red curve,draw=yellow!80!orange!85!black]
  ({-1+.32*cos(395)},{-1+.32*sin(395)})--(-.35,-.65)--(join)
  .. controls (-.10,-.18) and (-.02,-.10) .. (D);
\draw[chsix red curve]
  (join) .. controls (-.24,-.14) and (-.24,-.06) .. (-.22,0)
  arc[start angle=180,end angle=37,radius=.22];
\draw[chsix red curve,overview arrow]
  ({.22*cos(37)},{.22*sin(37)})--(.95,.55)--(B);
\draw[chsix red curve,draw=green!50!black]
  (join) .. controls (-.34,-.15) and (-.35,-.07) .. (-.32,0)
  arc[start angle=180,end angle=37,radius=.32];
\draw[chsix red curve,draw=green!50!black,overview arrow]
  ({.32*cos(37)},{.32*sin(37)})--(.90,.60)--(B);
\node[chsix point] at (join) {};
\node[chsix point,label={below left:$A$}] at (A) {};
\node[chsix point] at (Q) {};
\node[arc label,inner sep=.5pt] at (-1,-.88) {$Q_1$};
\node[chsix point,label={below right:$D$}] at (D) {};
\node[chsix point,label={above right:$B$}] at (B) {};
\node[arc label,text=annred] at (-1.25,-1.50) {$\gamma$};
\node[arc label,text=blue!72!black] at (-1.52,-1.16) {$\delta$};
\node[arc label,text=green!50!black] at (-1.83,-1.40) {$\gamma'$};
\node[arc label,text=yellow!80!orange!85!black] at (-1.28,-1.82) {$\gamma''$};
\node[arc label,text=annred] at (.99,.38) {$\gamma$};
\node[arc label,text=green!50!black] at (.83,.79) {$\gamma'$};
\end{tikzpicture}
\caption{Reconnection in a common triangle}
\label{fig:fresh-four-arcs}
\end{figure}

To prove $|\gamma'|<|\gamma|$, compare the lengths of the four arcs $\gamma,\delta,\gamma',\gamma''$. First compute the product for a push-off of the last segment $Q_{n-1}D$. This segment has no interior lattice point, so the coordinates of $v=(q,p)$ are relatively prime. With its initial endpoint as the origin, an edge midpoint $tv$ satisfies $2tp,2tq\in\mathbb Z$. Since $\gcd(p,q)=1$, this implies $2t\in\mathbb Z$, so $0<t<1$ forces $t=1/2$. Hence the left and right push-offs can differ in sign only at the central edge midpoint. Let $k$ be the weight of its triangulation edge, and let $P$ be the product of the subword before crossing that edge, omitting the initial triangle sign. The two halves correspond under the half-turn about the midpoint. Reversing the order of the first half's signs and reversing each sign produces the second half, whose product is therefore $P^{\mathsf T}$. If $W_-,W_+$ are the interior-word products of the right and left push-offs of $Q_{n-1}D$, respectively, then
\[
 W_-=PU_-^kP^{\mathsf T},\qquad W_+=PU_+^kP^{\mathsf T}
\]
Choose where to cut the interior words of $\gamma,\delta$ according to the central sign of the last segment and its weight $k$. Let $H$ be the common product from the cut to just before $D$, distinguishing two cases.
\begin{enumerate}
\item[(i)] The arc $\gamma$ avoids the center of $Q_{n-1}D$ on the right and $k>0$. Cut immediately after the factor for crossing the edge containing this midpoint. That factor is $U_-^k$ for $\gamma$ and $U_+^k$ for $\delta$; the common product thereafter is $H=P^{\mathsf T}$.
\item[(ii)] The arc $\gamma$ avoids the center on the left, or $k=0$. Cut immediately before the product for the last segment. This product is $W_+$ for both arcs, so $H=W_+$.
\end{enumerate}
In either case, let $r=(r_1,r_2)$ and $t=(t_1,t_2)$ be $\mathbf1^{\mathsf T}$ multiplied by the interior-word products up to the cut for $\gamma$ and $\delta$, respectively. The triangle sign where $\gamma$ begins its left detour at $D$ is $+$. If $v_0$ is the positive column for the remaining interior word, the column after $H$ is $c=U_+v_0=(c_1,c_2)^{\mathsf T}$. Since $\delta$ ends at $D$ and reconnection exchanges the two tails, the four GM lengths are
\[
 |\gamma|=rHc,\qquad |\delta|=tH\mathbf1,\qquad
 |\gamma'|=tHc,\qquad |\gamma''|=rH\mathbf1.
\]
First compare the entries of $r,t$. In (i), their final factors are $U_-^k,U_+^k$. In (ii), the triangle ending the right detour at $Q_{n-1}$ has sign $-$, and the one ending the left detour has sign $+$, so the final factors are $U_-,U_+$. Put $m=k$ in (i) and $m=1$ in (ii). Let $(x,y),(x',y')$ be the rows immediately before these final factors. Their entries are positive because they are obtained from $\mathbf1^{\mathsf T}$ by multiplying $U_\pm$. Thus
\[
 r=(x,y)U_-^m=(x,mx+y),\qquad t=(x',y')U_+^m=(x'+my',y')
\]
Since $m\geq1$, we have $mx+y>x$ and $x'+my'>y'$, giving $r_1<r_2$ and $t_1>t_2$. Hence $r_1t_2<r_2t_2<r_2t_1$, so $\det\left(\begin{smallmatrix}r\\t\end{smallmatrix}\right)=r_1t_2-r_2t_1<0$. Every $U_\pm$ has determinant $1$, so $\det H=1$ for either $H=P^{\mathsf T}$ or $H=W_+$.

Write $(c,\mathbf1)$ for the matrix with these two columns in order. Since $c=U_+v_0$, we have $c_1<c_2$, so $\det(c,\mathbf1)=c_1-c_2<0$. Expressing the difference of products of the four lengths as a determinant and using multiplicativity gives
\begin{equation}\label{eq:fresh-cross-det}
 |\gamma||\delta|-|\gamma'||\gamma''|
 =\det\begin{pmatrix}rHc&rH\mathbf1\\tHc&tH\mathbf1\end{pmatrix}
 =\det\left(\begin{pmatrix}r\\t\end{pmatrix}H(c,\mathbf1)\right)
 =\det\begin{pmatrix}r\\t\end{pmatrix}\det(c,\mathbf1)>0
\end{equation}
The last equality uses $\det H=1$, and the strict inequality uses the negativity of both determinants. To deduce $|\gamma'|<|\gamma|$, we now prove $|\gamma''|\geq|\delta|$ by comparing the sign sequence of $\gamma''$ with that of the pure right push-off from $A$ to $D$.

By hypothesis, $\gamma$ comes from Lemma~\ref{lem:fresh-canonical-geodesic}(2) and already has canonical midpoint signs. For $j=1,\ldots,n-1$, the segment $Q_{j-1}Q_j$ is a translate of $Q_{n-1}D$, so its only interior edge midpoint is its center. Starting at the crossing of its midpoint edge, compare triangle signs outward toward $Q_{j-1}$ and $Q_j$. The half-turn about the midpoint identifies the two sides, so corresponding triangle signs remain opposite until reaching the lattice-point detours. For $j\geq2$, both endpoints are bypassed on the right, and the first equal pair consists of two $-$ signs. The canonical sign at the center is therefore $-$. For $j=1$, the endpoint marker $0$ at $Q_0=A$ is compared with the initial detour sign $-$, giving the same result. Thus, except on the last segment $Q_{n-1}D$, the signs of $\gamma$ agree with those of the pure right push-off.

For the last segment, comparison outward from its center reaches opposite signs $-$ on the $Q_{n-1}$ side and $+$ on the $D$ side, so it must continue further. The canonical central sign can be either $+$ or $-$, depending on the full passage sequence and its ends. If $\gamma$ avoids the central midpoint on the right, then $\gamma''$ is the pure right push-off $\gamma^R_{AD}$, and left--right equality gives $|\gamma''|=|\gamma^R_{AD}|=|\delta|$. If it avoids it on the left, then $\gamma''$ differs from $\gamma^R_{AD}$ only in this last midpoint detour. Split the interior product of $\gamma^R_{AD}$ immediately before the last segment product $W_-$, and let $z=(z_1,z_2)$ be the row for the preceding portion. Then
\begin{equation}\label{eq:fresh-last-centre}
 |\gamma''|-|\gamma^R_{AD}|
 =zP(U_+^k-U_-^k)P^{\mathsf T}\mathbf1
 =k\det\begin{pmatrix}\mathbf1^{\mathsf T}P\\zP\end{pmatrix}
 =k(z_2-z_1)\geq0
\end{equation}
The last equality uses $\det P=1$. The last factor producing $z$ is $U_-$, from the triangle ending the right detour at $Q_{n-1}$, so $z_2>z_1$. Hence again $|\gamma''|\geq|\gamma^R_{AD}|=|\delta|$, with equality if $k=0$. Therefore
\[
 |\gamma||\delta|>|\gamma'||\gamma''|
 \geq|\gamma'||\delta|
\]
Dividing by the positive number $|\delta|$ gives $|\gamma|>|\gamma'|$. The calculation remains valid for $P=E_2$, when the subwords on either side of the midpoint edge are empty.

\smallskip\noindent
\textbf{(\ref{case:fresh-bend-right}) Nonedge direction along $AD$, followed by more than a right half-turn at $D$.}
Now suppose the detour at $D$ also passes on the right, through an angle greater than $\pi$. If the outgoing polygonal path follows a triangulation edge, this case covers only arcs that do not cross that edge incident to $D$; crossings of it are treated in (\ref{case:fresh-bend-exit-edge}). Apply the same coordinate change as in (\ref{case:fresh-bend-left}) so that $v=(q,p)$ with $p,q>0$, and use the products $W_\pm$ found there for each segment. Since the detour at $D$ is also on the right, the same canonical-midpoint argument gives central sign $-$ on every segment, including the last. Thus the arc up to just before $D$ agrees with the pure right push-off.

Let $\Delta$ be the triangle where $\gamma$ finishes its detour around $D$, and distinguish three cases.
\begin{enumerate}
\item[(i)] The point $B$ is a vertex of $\Delta$, and $\gamma$ ends in this triangle.
\item[(ii)] The arc continues beyond $\Delta$, and $\Delta$ is the first triangle met by the ray from $D$ in direction $v$.
\item[(iii)] The arc continues beyond $\Delta$, and $\Delta$ differs from that first triangle.
\end{enumerate}
If $B$ is a vertex of $\Delta$, Lemma~\ref{lem:fresh-minimizer}(2) makes it the last passage triangle. Otherwise $\gamma$ leaves $\Delta$ through the side opposite $D$.

\smallskip\noindent
\textbf{(i) The endpoint $B$ is a vertex of $\Delta$.}
Replace the right detour at $D$ by a left detour and connect directly to $B$, obtaining $\gamma'$. Let $j$ be the number of edges crossed by the original detour. Since it exceeds a half-turn without revisiting a triangle, $j=3,4,5$. If $j=5$, the vertex $B$ opposite the entry edge of $\Delta$ is also a vertex of the triangle where the detour started, contrary to Lemma~\ref{lem:fresh-minimizer}(2). Thus $j=3,4$. A left route into $\Delta$ crosses $6-j$ edges, ending with $DB$. Connecting to $B$ after this last crossing violates the endpoint condition. Instead, connect directly to $B$ within the triangle just before $DB$. Its entry edge is opposite $B$, so this gives a generalized arc $\gamma'$. Figure~\ref{fig:fresh-terminal-shortcut} depicts $j=3$. The left panel shows the intermediate route passing to the left of $D$, crossing $DB$, and ending at $B$; the right panel connects directly before $DB$, with the removed portion dashed in gray.

\begin{figure}[!htbp]
\centering
\begin{tikzpicture}[x=2.20cm,y=2.20cm,
  line cap=round,line join=round,every node/.style={font=\small},
  removed path/.style={draw=black!45,dash pattern=on 2.8pt off 2pt,line width=.7pt}]
\def\terminalradius{.32}
\foreach \panel in {0,1}{
  \begin{scope}[shift={({3*\panel},0)}]
    \begin{scope}
      \clip (-1.16,-1.16) rectangle (1.08,1.18);
      \fill[black!6] (0,0)--(-1,1)--(0,1)--cycle;
      \foreach \i in {-1,0,1}{
        \draw[chsix triangulation] (\i,-1.2)--(\i,1.2);
        \draw[chsix triangulation] (-1.2,\i)--(1.2,\i);
      }
      \foreach \i in {-2,-1,0,1,2}{
        \draw[chsix triangulation] (-1.2,{\i+1.2})--(1.2,{\i-1.2});
      }
    \end{scope}
    \coordinate (A) at (-1,-1);
    \coordinate (entry) at ({\terminalradius*cos(135)},{\terminalradius*sin(135)});
    \draw[chsix distinguished edge] (0,0)--(0,1);
    \draw[chsix red curve,chsix forward at=.55]
      (A) .. controls (-.62,-.54) and (-\terminalradius,-.27) .. (-\terminalradius,0)
      arc[start angle=180,end angle=135,radius=\terminalradius];
    \ifnum\panel=0
      \draw[chsix red curve,chsix forward at=.69]
        (entry) arc[start angle=135,end angle=90,radius=\terminalradius]
        .. controls (.27,\terminalradius) and (.30,.64) .. (0,1);
    \else
      \draw[removed path]
        (entry) arc[start angle=135,end angle=90,radius=\terminalradius]
        .. controls (.27,\terminalradius) and (.30,.64) .. (0,1);
      \draw[chsix red curve,chsix forward at=.58] (entry)--(0,1);
    \fi
    \node[chsix point,label={[fill=white,inner sep=.6pt]below left:$A$}] at (A) {};
    \node[chsix point,label={[fill=white,inner sep=.6pt]below right:$D$}] at (0,0) {};
    \node[chsix red point,label={[fill=white,inner sep=.6pt]above:$B$}] at (0,1) {};
  \end{scope}
}
\draw[chsix transition arrow] (1.23,.53)--(1.77,.53);
\end{tikzpicture}
\caption{Connecting directly to the endpoint}
\label{fig:fresh-terminal-shortcut}
\end{figure}
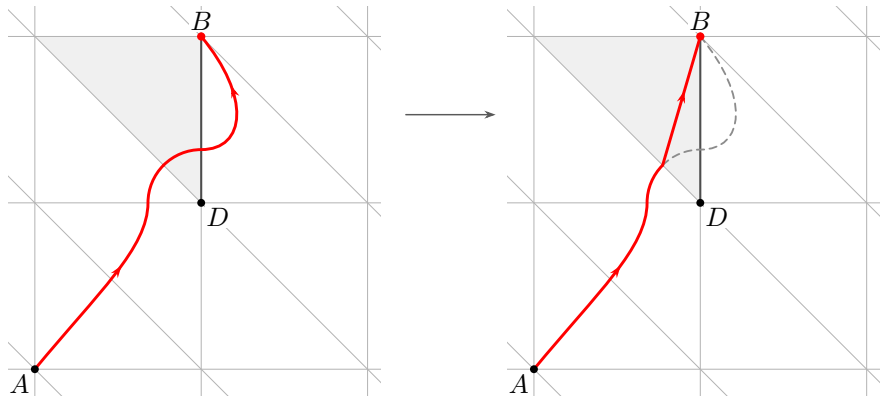

Count the detour signs to show that this direct connection decreases length. Let $h$ be the weight of $DB$. For the right and left routes into $\Delta$, let $a,b$ be the sum of crossed-edge weights plus the number of intervening triangles. These triangle counts are $j-1$ and $5-j$. Opposite edges have equal weights, and the right route's total edge weight is at least the left route's, so $a-b\geq2j-6\geq0$. The original detour begins with triangle sign $-$, all intermediate signs are $+$, and the last triangle sign is excluded from the interior word. For the left route the signs are reversed, but connecting before $DB$ also removes its last edge crossing and the preceding triangle sign. Thus, writing $r_0=(x,N-x)$ and $t_0=(u,N-u)$, we have
\[
 |\gamma|=r_0U_-U_+^a\mathbf1=x+(a+1)N,\qquad
 |\gamma'|=t_0U_+U_-^{b-h-1}\mathbf1=(b-h+1)N-u
\]
Here $0<x,u<N$, and $b-h-1$ is the nonnegative sum of the weights of the remaining $5-j$ edges and $4-j$. Hence $|\gamma|-|\gamma'|=(a-b+h)N+x+u>0$.

Figure~\ref{fig:fresh-exit-triangle-cases} illustrates the remaining cases (ii),(iii). The gray triangle is the first met by the ray from $D$ in direction $v$; the pale orange triangle in the right panel is $\Delta$. The red curve is $\gamma$, and the dashed lines indicate the shortest polygonal path and the ray in direction $v$.
\begin{figure}[htbp]
\centering
\begin{tikzpicture}[x=2.1cm,y=2.1cm,
  every node/.style={font=\small},
  folded path/.style={draw=black!55,densely dashed,line width=.6pt},
  reference ray/.style={draw=black!75,densely dashed,line width=.7pt,
    -{Stealth[length=1.6mm,width=1.05mm]}},
  detour arrow/.style={postaction={decorate},decoration={markings,
    mark=at position .55 with {\arrow{Stealth[length=1.5mm,width=1.0mm]}}}},
  curve label/.style={fill=white,inner sep=1pt}]
\def\exitradius{.24}
\foreach \panel/\outangle/\endangle/\exitx in {0/75/425/.47,1/120/470/-.55}{
  \begin{scope}[shift={({3.25*\panel},0)}]
    \begin{scope}
      \clip (-1.18,-1.18) rectangle (1.30,1.44);
      \fill[black!9] (0,0)--(1,0)--(0,1)--cycle;
      \ifnum\panel=1
        \fill[orange!15] (0,0)--(0,1)--(-1,1)--cycle;
      \fi
      \foreach \i in {-1,0,1}{
        \draw[chsix triangulation] (\i,-2)--(\i,2);
        \draw[chsix triangulation] (-2,\i)--(2,\i);
      }
      \foreach \s in {-2,-1,...,3}{
        \draw[chsix triangulation] (-2,{\s+2})--(2,{\s-2});
      }
    \end{scope}
    \draw[chsix distinguished edge] (0,0)--(1,0)--(0,1)--cycle;
    \ifnum\panel=1
      \draw[draw=orange!70!black,line width=.75pt] (0,0)--(0,1)--(-1,1)--cycle;
    \fi
    \draw[folded path] (-1,-1)--(0,0)--({1.31*cos(\outangle)/sin(\outangle)},1.31);
    \draw[reference ray] (0,0)--(1.18,1.18);
    \node[curve label,below right=2pt] at (.89,.89) {$v$};
    \draw[chsix red curve]
      (-1,-1)--({\exitradius*cos(235)},{\exitradius*sin(235)});
    \draw[chsix red curve,detour arrow]
      ({\exitradius*cos(235)},{\exitradius*sin(235)})
      arc[start angle=235,end angle=\endangle,radius=\exitradius];
    \draw[chsix red curve,-{Stealth[length=1.6mm,width=1.05mm]}]
      ({\exitradius*cos(\endangle)},{\exitradius*sin(\endangle)})--(\exitx,1.30);
    \node[chsix point,label={below left:$A$}] at (-1,-1) {};
    \node[chsix point,label={[fill=white,inner sep=.5pt]above left:$D$}] at (0,0) {};
    \ifnum\panel=0
      \node[chsix panel title] at (.05,1.72) {(ii) Same triangle};
      \node[curve label,text=annred] at (.61,1.01) {$\gamma$};
    \else
      \node[chsix panel title] at (.05,1.72) {(iii) Different triangles};
      \node[curve label,text=annred] at (-.27,1.12) {$\gamma$};
    \fi
  \end{scope}
}
\end{tikzpicture}
\caption{Comparing the triangles where the detour ends}
\label{fig:fresh-exit-triangle-cases}
\end{figure}
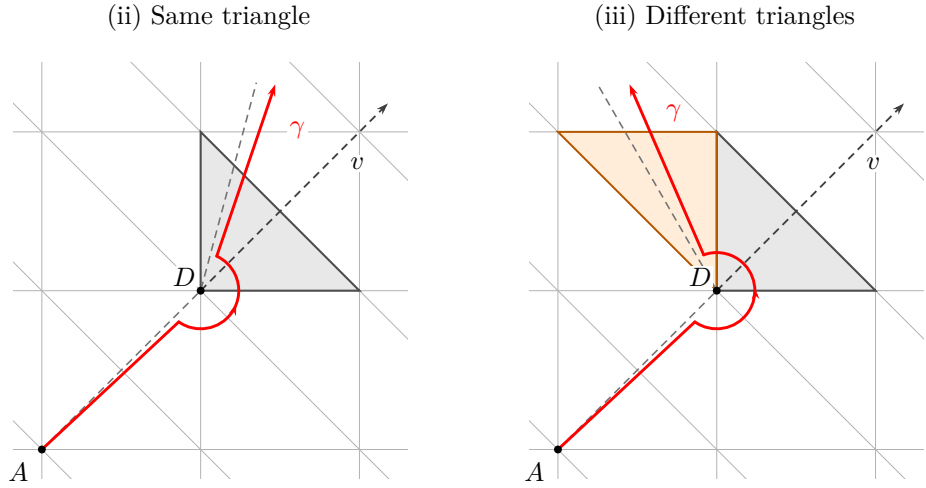

\smallskip\noindent
\textbf{(ii) The arc continues, and the two triangles coincide.}
Consider the left panel of Figure~\ref{fig:fresh-exit-triangle-cases}. Put $E=D+v$ and let $\delta$ be the pure left push-off from $A$ to $E$. Reconnect within the common triangle $\Delta$: the initial portion of $\delta$ joined to the tail of $\gamma$ gives $\gamma'$ from $A$ to $B$, and the remaining portions give $\gamma''$ from $A$ to $E$ (Figure~\ref{fig:fresh-after-d-splice}). Each connection enters through an edge incident to $D$ and leaves through the opposite side, so neither creates a backtracking crossing. Endpoint passages are unchanged, and both curves are generalized arcs. The figure shows $\gamma,\delta,\gamma',\gamma''$ in red, blue, green, and yellow, slightly separating overlaps. The shaded triangle is $\Delta$, with the exchange point marked in black. Dashed lines show the line through $A,D,E$ and direction $v$; the thick edge is opposite $D$.

\begin{figure}[htbp]
\centering
\begin{tikzpicture}[x=2.9cm,y=2.9cm,
  every node/.style={font=\small},
  arc label/.style={fill=white,inner sep=1.2pt}]
\begin{scope}
  \clip (-1.12,-1.12) rectangle (1.22,1.42);
  \fill[black!6] (0,0)--(1,0)--(0,1)--cycle;
  \foreach \i in {-1,0,1}{
    \draw[chsix triangulation] (\i,-2)--(\i,2);
    \draw[chsix triangulation] (-2,\i)--(2,\i);
  }
  \foreach \s in {-2,-1,0,1,2}{
    \draw[chsix triangulation] (-2,{\s+2})--(2,{\s-2});
  }
\end{scope}
\draw[chsix distinguished edge] (1,0)--(0,1);
\draw[draw=black!55,densely dashed,line width=.6pt,
  -{Stealth[length=1.4mm,width=.95mm]}] (-1,-1)--(1.13,1.13);
\node[below right=2pt,text=black!70] at (.70,.70) {$v$};
\coordinate (join) at (.225,.375);
\draw[chsix red curve,-{Stealth[length=1.5mm,width=1.0mm]}]
  (-1,-1)--({.20*cos(250)},{.20*sin(250)})
  arc[start angle=250,end angle=360,radius=.20]
  --(join)--(.28,.72)--(.42,1.25);
\draw[chsix blue curve,chsix forward at=.78]
  (-1,-1)--({.20*cos(200)},{.20*sin(200)})
  arc[start angle=200,end angle=90,radius=.20]
  --(join)--(.40,.60)--(1,1);
\draw[chsix red curve,draw=green!50!black,-{Stealth[length=1.5mm,width=1.0mm]}]
  (-1,-1)--({.32*cos(200)},{.32*sin(200)})
  arc[start angle=200,end angle=90,radius=.32]
  --(join)--(.17,.83)--(.20,1.30);
\draw[chsix red curve,draw=yellow!80!orange!85!black,chsix forward at=.80]
  (-1,-1)--({.32*cos(250)},{.32*sin(250)})
  arc[start angle=250,end angle=360,radius=.32]
  --(join)--(.48,.52)--(1,1);
\node[chsix point] at (join) {};
\node[chsix point,label={below left:$A$}] at (-1,-1) {};
\node[chsix point,label={[fill=white,inner sep=.7pt]below left:$D$}] at (0,0) {};
\node[chsix point,label={above right:$E$}] at (1,1) {};
\node[above=2pt] at (.30,1.34) {to $B$};
\node[arc label,text=annred] at (.48,1.03) {$\gamma$};
\node[arc label,text=blue!72!black] at (.50,.77) {$\delta$};
\node[arc label,text=green!50!black] at (.06,1.03) {$\gamma'$};
\node[arc label,text=yellow!80!orange!85!black] at (.69,.53) {$\gamma''$};
\end{tikzpicture}
\caption{Reconnection just after $D$}
\label{fig:fresh-after-d-splice}
\end{figure}
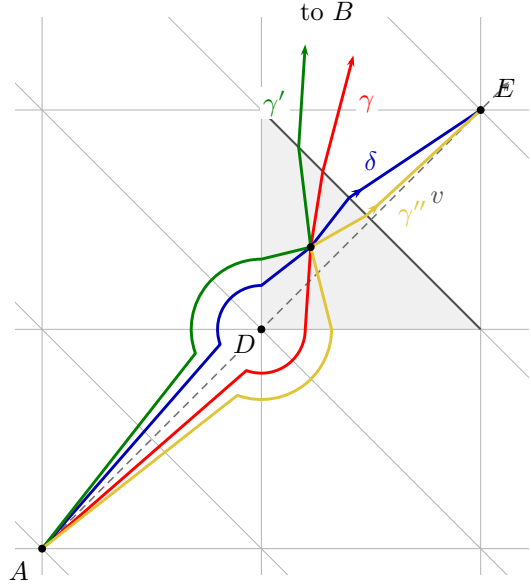

Compute the connecting products to compare the right and left routes around $D$. Let $C_-,C_+$ be the products for connecting the push-offs of $AD$ and $DE$ by right and left half-turns at $D$, including both end triangle signs. The right connection crosses one edge of each type and passes through four triangles. Its first and last triangle signs are $-$ and the two intermediate ones are $+$. With $K_0=k_1+k_2+k_3$, its word is $-+^{K_0+2}-$. The left connection reverses all signs, so
\begin{equation}\label{eq:fresh-primitive-fan}
 C_-=U_-U_+^{K_0+2}U_-,\qquad C_+=U_+U_-^{K_0+2}U_+
\end{equation}
Recall that $N=|\gamma^R_{AD}|=|\gamma^L_{AD}|$ is the common GM length of the pure push-offs of $AD$. Their interior rows $r_0,t_0$ satisfy $r_0\mathbf1=t_0\mathbf1=N$. Include the sign of the common triangle and put
\[
 r=r_0C_-=(r_1,r_2),\qquad t=t_0C_+=(t_1,t_2)
\]
Let $w=w_1\cdots w_s$ be the remaining interior word of $\gamma$, beginning with the crossing of the opposite side, and define its column by
\[
 c=M(w)\mathbf1=U_{w_1}\cdots U_{w_s}\begin{pmatrix}1\\1\end{pmatrix}
 =\begin{pmatrix}c_1\\c_2\end{pmatrix}
\]
The triangle sign containing $B$ is excluded from $w$. The tail of $\delta$ is the pure left push-off of $DE$, so its column is $d=W_+\mathbf1$. Write $(c,d)$ for the matrix with columns $c,d$ in that order. The four lengths are
\[
 |\gamma|=rc,\qquad |\delta|=td,\qquad
 |\gamma'|=tc,\qquad |\gamma''|=rd.
\]
The final factor producing $r$ is $U_-$, and that producing $t$ is $U_+$. Applying the component comparison in (\ref{case:fresh-bend-left}) with $m=1$ gives $0<r_1<r_2$ and $0<t_2<t_1$. Hence
\[
 \det\begin{pmatrix}r\\t\end{pmatrix}
 =\det\begin{pmatrix}r_1&r_2\\t_1&t_2\end{pmatrix}
 =r_1t_2-r_2t_1<0
\]
Next determine the sign of the determinant of the tail columns $c,d$. Let $g$ be the shortest polygonal path associated with $\gamma$. Looking from $D$ toward $E$, immediately after $D$ the path $g$ lies to the left of $DE$. While the two routes cross the same edges in the same order, glue their triangles into a common strip, ending where they first leave a triangle through different sides or where one route terminates. Suppose the portion of $g$ from $D$ toward $B$ meets $DE$ again in this strip, first at $Z$. The subpath of $g$ from $D$ to $Z$ and the segment $DZ$ on $DE$ are both shortest between $D,Z$. Uniqueness in the strip contradicts their differing immediately after $D$. Therefore $g$ stays to the left of $DE$ as long as they follow the same triangle sequence, up to the first distinct exit sides or the first endpoint.

In the first triangle where the two arcs enter through the same edge and leave through different edges, the left route $\gamma$ has sign $+$ and the push-off of $DE$ has sign $-$. If their midpoint detours first differ before this, the midpoint signs have the same order. For edge weight $h>0$, the difference is $U_+^h$ versus $U_-^h$; if $h=0$, it contributes no factor, and comparison continues. Let $G$ be the common initial product and $c',d'$ the remaining columns. Then $c=Gc'$, $d=Gd'$, and $\det G=1$, so $\det(c,d)=\det(c',d')$. At the first effective difference, $c'$ starts with $U_+^m$ and $d'$ with $U_-^\ell$, for $m,\ell\geq1$. Thus $c'_1<c'_2$ and $d'_1>d'_2$, giving $\det(c',d')=c'_1d'_2-c'_2d'_1<0$.

If one route ends first, compare the endpoints as well. The generalized-arc endpoint condition places that endpoint at the vertex opposite the common entry edge, so its remaining column is $\mathbf1$. If $\gamma$ ends first, the other column starts with $U_-$, giving $\det(\mathbf1,d')=d'_2-d'_1<0$. If the $DE$ route ends first, then $\det(c',\mathbf1)=c'_1-c'_2<0$. Both routes cannot end simultaneously at $E$, for that would give $B-A=(n+1)v$, placing the lattice point $D$ inside $AB$, contrary to hypothesis. Thus in all cases, including an earlier endpoint,
\begin{equation}\label{eq:fresh-column-order}
 \det(c,d)<0.
\end{equation}
As in (\ref{case:fresh-bend-left}), use these rows and columns to compute
\[
 (rc)(td)-(tc)(rd)
 =\det\begin{pmatrix}rc&rd\\tc&td\end{pmatrix}
 =\det\begin{pmatrix}r\\t\end{pmatrix}\det(c,d)>0
\]
Moreover, $rd$ is the length of the pure right push-off from $A$ to $E$ with only its last central midpoint changed to the left. Apply the calculation in \eqref{eq:fresh-last-centre} to the final segment $DE$ with $z=r$. Together with equality of the left and right push-off lengths, it gives $rd\geq|\gamma^R_{AE}|=|\delta|=td$. Hence
\[
 (rc)(td)>(tc)(rd)\geq(tc)(td)
\]
Dividing by the positive number $td$ yields $|\gamma|=rc>tc=|\gamma'|$.

\smallskip\noindent
\textbf{(iii) The arc continues, and the two triangles differ.}
Consider the right panel of Figure~\ref{fig:fresh-exit-triangle-cases}. Replace the long right detour around $D$ by a left detour through a smaller angle, joining the same tail (Figure~\ref{fig:fresh-long-detour-shortcut}). If $j$ is the original number of crossings around $D$, then $j=4,5$, since $j=6$ would revisit the starting triangle. The left detour crosses $6-j$ edges. Thus the new arc still has at least one crossing; it enters $\Delta$ through an edge incident to $D$ and leaves through the opposite side to join the old tail. It does not backtrack across the entry edge of $\Delta$, and its endpoint passages are unchanged, so it is a generalized arc $\gamma'$. The figure shows $j=4$, with the old arc $\gamma$ on the left and $\gamma'$ on the right. Both share the portion beyond the black point in the shaded triangle $\Delta$; the gray dashed curve on the right is the replaced portion.

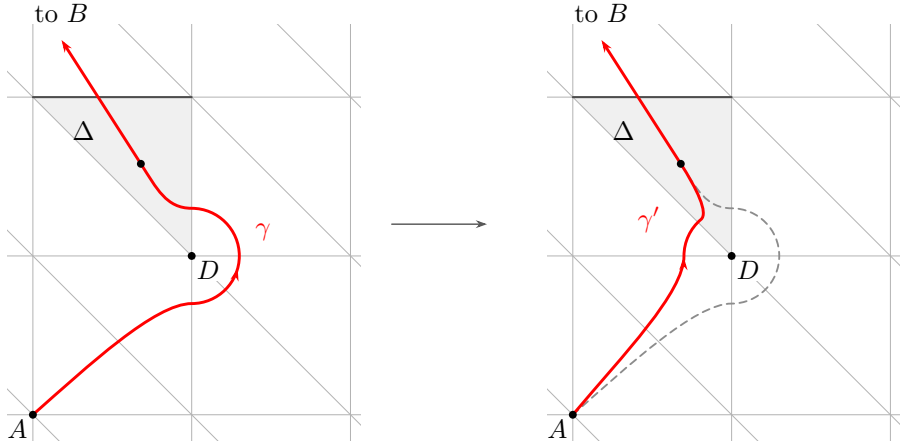
\begin{figure}[!htbp]
\centering
\begin{tikzpicture}[x=2.10cm,y=2.10cm,
  line cap=round,line join=round,every node/.style={font=\small},
  replaced path/.style={draw=black!45,dash pattern=on 2.8pt off 2pt,line width=.7pt}]
\def\longdetourradius{.30}
\foreach \panel in {0,1}{
  \begin{scope}[shift={({3.4*\panel},0)}]
    \begin{scope}
      \clip (-1.16,-1.16) rectangle (1.06,1.46);
      \fill[black!6] (0,0)--(0,1)--(-1,1)--cycle;
      \foreach \i in {-1,0,1}{
        \draw[chsix triangulation] (\i,-1.2)--(\i,1.5);
        \draw[chsix triangulation] (-1.2,\i)--(1.1,\i);
      }
      \foreach \s in {-2,-1,0,1,2}{
        \draw[chsix triangulation] (-1.2,{\s+1.2})--(1.1,{\s-1.1});
      }
    \end{scope}
    \coordinate (A) at (-1,-1);
    \coordinate (splice) at (-.32,.58);
    \draw[chsix distinguished edge] (-1,1)--(0,1);
    \ifnum\panel=0
      \draw[chsix red curve,chsix forward at=.62]
        (A) .. controls (-.65,-.69) and (-.23,-\longdetourradius) .. (0,-\longdetourradius)
        arc[start angle=270,end angle=450,radius=\longdetourradius]
        .. controls (-.16,\longdetourradius) and (-.22,.424) .. (splice);
      \node[text=annred,fill=white,inner sep=1pt] at (.45,.15) {$\gamma$};
    \else
      \draw[replaced path]
        (A) .. controls (-.65,-.69) and (-.23,-\longdetourradius) .. (0,-\longdetourradius)
        arc[start angle=270,end angle=450,radius=\longdetourradius]
        .. controls (-.16,\longdetourradius) and (-.22,.424) .. (splice);
      \draw[chsix red curve,chsix forward at=.65]
        (A) .. controls (-.68,-.62) and (-\longdetourradius,-.23) .. (-\longdetourradius,0)
        arc[start angle=180,end angle=130,radius=\longdetourradius]
        .. controls (-.145,.27) and (-.22,.424) .. (splice);
      \node[text=annred,fill=white,inner sep=1pt] at (-.52,.23) {$\gamma'$};
    \fi
    \draw[chsix red curve,-{Stealth[length=1.65mm,width=1.10mm]}]
      (splice)--(-.82,1.36);
    \node[chsix point] at (splice) {};
    \node[chsix point,label={[fill=white,inner sep=.6pt]below left:$A$}] at (A) {};
    \node[chsix point,label={[fill=white,inner sep=.6pt]below right:$D$}] at (0,0) {};
    \node at (-.68,.79) {$\Delta$};
    \node[above=3pt] at (-.82,1.36) {to $B$};
  \end{scope}
}
\draw[chsix transition arrow] (1.26,.20)--(1.86,.20);
\end{tikzpicture}
\caption{Changing the detour at $D$}
\label{fig:fresh-long-detour-shortcut}
\end{figure}

Write the products for the right and left detour words as $U_-U_+^aU_-$ and $U_+U_-^bU_+$, respectively. Here $a,b$ count all signs except the two end triangle signs: the sums of crossed-edge weights plus $j-1$ and $5-j$ intervening triangles, respectively. Opposite edge directions have equal weights, so the total edge weight on the right is at least that on the left. The difference in triangle counts is $2j-6$, hence $a-b\geq2j-6\geq2$.

Immediately before the shared tail, the old and new rows are $r_0U_-U_+^aU_-$ and $t_0U_+U_-^bU_+$. Both $r_0,t_0$ have component sum $N$, the GM length of the pure push-offs of $AD$. Set $x:=(r_0)_1$, $u:=(t_0)_1$, so $r_0=(x,N-x)$, $t_0=(u,N-u)$, with $0<x,u<N$ because all entries are positive. Their difference is
\begin{equation}\label{eq:fresh-long-fan}
 r_0U_-U_+^aU_--t_0U_+U_-^bU_+
 =\bigl((a-b-2)N+x+u,\ (a-b)N+x+u\bigr)>0
\end{equation}
The strict inequality means that both entries are positive. Multiplying both rows by the positive column for the remaining common tail yields $|\gamma|>|\gamma'|$.

\smallskip\noindent
\textbf{(\ref{case:fresh-bend-along-edge}) The path from $A$ to $D$ follows triangulation edges.}
Now consider a horizontal, vertical, or slope-$-1$ direction from $A$ to $D$, excluding (\ref{case:fresh-bend-exit-edge}). Use the coordinate change from (\ref{case:fresh-bend-left}) to make the direction horizontal. Let $h$ be the horizontal-edge weight, and let $a_*$ be one plus the sum of the other two weights. For $n\geq2$, the interior word of the pure right push-off from $A$ to $D$ is $+^{a_*}(-+^{a_*})^{n-2}$; the left push-off reverses every sign. Replacing signs in order by matrices and multiplying on the left by $\mathbf1^{\mathsf T}$ gives rows of the form
\[
 r_0=(\alpha,\beta),\qquad t_0=(\beta,\alpha),
 \qquad \alpha+\beta=N
\]
In \eqref{eq:fresh-long-fan}, this corresponds to $x=\alpha$ and $u=\beta=N-x$. For $n=1$, the points $A,D$ are adjacent lattice points, so both push-offs have no crossings and $N=1$. In that case compute each comparison directly, omitting the triangle sign containing the initial endpoint.

First suppose the detour side changes at $D$. This requires $n\geq2$, since the side at $D$ was designated right when $n=1$. Join the pure left push-off to the portion of $\gamma$ beyond $D$. This removes the horizontal crossing immediately before $D$ (Figure~\ref{fig:fresh-horizontal-crossing-removal}). The triangle signs before and after this crossing in the old curve are $+$ and $-$; the new connection retains just one triangle sign, $+$. The figure shows $n=2$, with the removed horizontal crossing on the thick edge $Q_1D$. The original arc is on the left and the new arc on the right; the portions beyond the black point agree, and the replaced portion is gray and dashed.

\begin{figure}[!htbp]
\centering
\begin{tikzpicture}[x=1.82cm,y=1.82cm,
  line cap=round,line join=round,every node/.style={font=\small},
  crossing sign/.style={text=annred,fill=white,inner sep=.8pt},
  replaced path/.style={draw=black!45,dash pattern=on 2.8pt off 2pt,line width=.7pt}]
\def\horizontalradius{.24}
\def\oldhorizontalpath{
  (A) .. controls (-1.70,-.04) and (-1.25,-.0473) ..
  ({-1+\horizontalradius*cos(210)},{\horizontalradius*sin(210)})
  arc[start angle=210,end angle=315,radius=\horizontalradius]
  .. controls (-.68,-.02) and (-\horizontalradius,-.18) .. (-\horizontalradius,0)
  arc[start angle=180,end angle=150,radius=\horizontalradius]
}
\foreach \panel in {0,1}{
  \begin{scope}[shift={({4.1*\panel},0)}]
    \begin{scope}
      \clip (-2.16,-.74) rectangle (1.04,1.06);
      \fill[black!6] (-1,0)--(0,0)--(-1,1)--cycle;
      \foreach \i in {-2,-1,0,1}{
        \draw[chsix triangulation] (\i,-1)--(\i,1.2);
      }
      \foreach \i in {0,1}{
        \draw[chsix triangulation] (-2.3,\i)--(1.2,\i);
      }
      \foreach \s in {-3,-2,-1,0,1,2}{
        \draw[chsix triangulation] (-2.3,{\s+2.3})--(1.2,{\s-1.2});
      }
    \end{scope}
    \coordinate (A) at (-2,0);
    \coordinate (splice) at ({\horizontalradius*cos(150)},{\horizontalradius*sin(150)});
    \draw[chsix distinguished edge] (-1,0)--(0,0);
    \ifnum\panel=0
      \draw[chsix red curve,chsix forward at=.47] \oldhorizontalpath;
      \node[crossing sign] at (-.54,-.30) {$+$};
      \node[crossing sign] at (-.47,.20) {$-$};
      \node[text=annred,fill=white,inner sep=1pt] at (-1.48,-.35) {$\gamma$};
    \else
      \draw[replaced path] \oldhorizontalpath;
      \draw[chsix red curve,chsix forward at=.46]
        (A) .. controls (-1.70,.04) and (-1.25,.0473) ..
        ({-1+\horizontalradius*cos(150)},{\horizontalradius*sin(150)})
        arc[start angle=150,end angle=45,radius=\horizontalradius]
        .. controls (-.68,.02) and (-.26,.03) .. (splice);
      \node[crossing sign] at (-.47,.20) {$+$};
      \node[text=annred,fill=white,inner sep=1pt] at (-1.48,.36) {$\gamma'$};
    \fi
    \draw[chsix red curve,-{Stealth[length=1.65mm,width=1.10mm]}]
      (splice) arc[start angle=150,end angle=45,radius=\horizontalradius]
      --(.96,{sqrt(2)*\horizontalradius-.96});
    \node[chsix point] at (splice) {};
    \node[chsix point,label={below left:$A$}] at (A) {};
    \node[chsix point,label={[fill=white,inner sep=.6pt]below right:$D$}] at (0,0) {};
    \node[chsix point] at (-1,0) {};
    \node[below] at (-1,-.31) {$Q_1$};
    \node[below right=1pt] at (.96,{sqrt(2)*\horizontalradius-.96}) {to $B$};
  \end{scope}
}
\draw[chsix transition arrow] (1.23,.20)--(1.69,.20);
\end{tikzpicture}
\caption{Removing a horizontal-edge crossing}
\label{fig:fresh-horizontal-crossing-removal}
\end{figure}
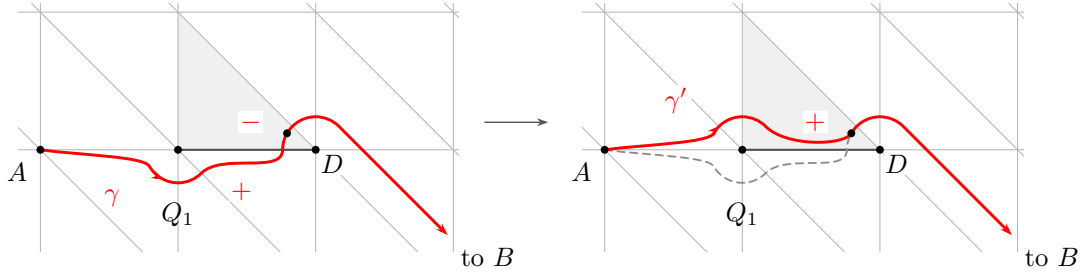

Let $\eta\in\{+,-\}$ be the midpoint sign on the horizontal edge. The rows immediately before rejoining the original route are $(\alpha,\beta)U_+U_\eta^hU_-$ for the old curve and $(\beta,\alpha)U_+$ for the new curve. Since $U_\eta^h\geq E_2$, their difference satisfies
\[
 (\alpha,\beta)U_+U_\eta^hU_--(\beta,\alpha)U_+
 \geq(\alpha,\beta)U_+U_--(\beta,\alpha)U_+=(0,2\beta).
\]
If the route continues beyond the connecting triangle, multiply both rows by the positive column obtained by applying the remaining interior product to $\mathbf1$. The second entry of the row difference is positive, so the old GM length is larger. If the route ends at a vertex $B$ of that triangle, follow the pure left push-off until it first enters a triangle having $B$ as a vertex, and connect directly to $B$. If crossings remain, the last entry edge is opposite $B$, and the interior word is an initial subword of the pure left push-off's word, so its GM length is at most $N$. If there are no crossings, its length is $1$, which satisfies the same bound. The old length, omitting the final triangle sign, is $r_0U_+U_\eta^h\mathbf1\geq r_0U_+\mathbf1=N+\beta>N$. Thus $|\gamma'|<|\gamma|$ in either case.

Next suppose the detour remains on the same side and exceeds a half-turn at $D$. First treat three or four crossings around $D$. Replace the initial portion from $A$ to $D$ by the pure left push-off and detour to the left of $D$ to join the old tail. This changes three crossings to two, or four to one. Figure~\ref{fig:fresh-edge-long-detours} shows $n=2$, with the three-to-two case above and the four-to-one case below. The old arc is on the left and the new arc on the right; red points mark crossings of edges incident to $D$. In each row the portions after the black point agree, and the old portion is gray and dashed in the right panel.

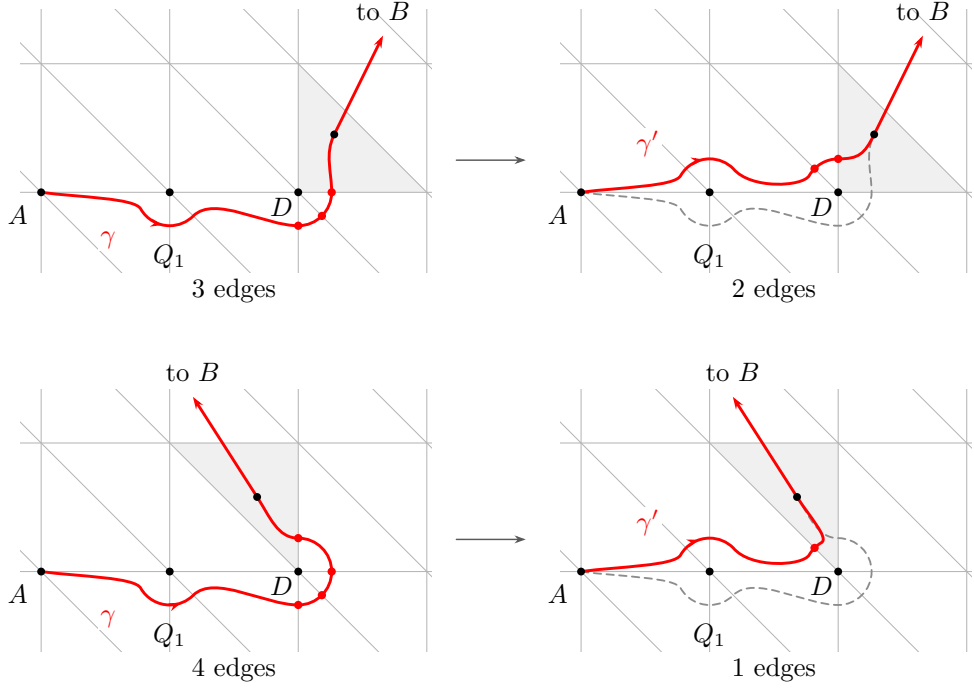
\begin{figure}[!htbp]
\centering
\begin{tikzpicture}[x=1.70cm,y=1.70cm,
  line cap=round,line join=round,every node/.style={font=\small},
  replaced path/.style={draw=black!45,dash pattern=on 2.8pt off 2pt,line width=.7pt}]
\def\edgefanradius{.26}
\def\oldedgefanprefix{
  (A) .. controls (-1.70,-.04) and (-1.27,-.0434) ..
  ({-1+\edgefanradius*cos(210)},{\edgefanradius*sin(210)})
  arc[start angle=210,end angle=315,radius=\edgefanradius]
  .. controls (-.65,-.02) and (-.26,-\edgefanradius) .. (0,-\edgefanradius)
}
\def\newedgefanprefix{
  (A) .. controls (-1.70,.04) and (-1.27,.0434) ..
  ({-1+\edgefanradius*cos(150)},{\edgefanradius*sin(150)})
  arc[start angle=150,end angle=45,radius=\edgefanradius]
  .. controls (-.68,.02) and (-.275,.0434) ..
  ({\edgefanradius*cos(150)},{\edgefanradius*sin(150)})
}
\foreach \row in {0,1}{
  \foreach \panel in {0,1}{
    \begin{scope}[shift={({4.2*\panel},{-2.95*\row})}]
      \begin{scope}
        \clip (-2.16,-.62) rectangle (1.04,1.42);
        \ifnum\row=0
          \fill[black!6] (0,0)--(1,0)--(0,1)--cycle;
        \else
          \fill[black!6] (0,0)--(0,1)--(-1,1)--cycle;
        \fi
        \foreach \i in {-2,-1,0,1}{
          \draw[chsix triangulation] (\i,-.8)--(\i,1.6);
        }
        \foreach \i in {0,1}{
          \draw[chsix triangulation] (-2.3,\i)--(1.2,\i);
        }
        \foreach \s in {-3,-2,-1,0,1,2}{
          \draw[chsix triangulation] (-2.3,{\s+2.3})--(1.2,{\s-1.2});
        }
      \end{scope}
      \coordinate (A) at (-2,0);
      \ifnum\row=0
        \coordinate (splice) at (.28,.45);
        \coordinate (tail) at (.66,1.22);
        \def\oldedgefanpath{
          \oldedgefanprefix
          arc[start angle=270,end angle=360,radius=\edgefanradius]
          .. controls (\edgefanradius,.16) and (.204,.296) .. (splice)
        }
        \def\newedgefanpath{
          \newedgefanprefix
          arc[start angle=150,end angle=90,radius=\edgefanradius]
          .. controls (.16,\edgefanradius) and (.204,.296) .. (splice)
        }
      \else
        \coordinate (splice) at (-.32,.58);
        \coordinate (tail) at (-.82,1.36);
        \def\oldedgefanpath{
          \oldedgefanprefix
          arc[start angle=270,end angle=450,radius=\edgefanradius]
          .. controls (-.16,\edgefanradius) and (-.22,.424) .. (splice)
        }
        \def\newedgefanpath{
          \newedgefanprefix
          arc[start angle=150,end angle=120,radius=\edgefanradius]
          .. controls (-.07,.26) and (-.22,.424) .. (splice)
        }
      \fi
      \ifnum\panel=0
        \draw[chsix red curve,chsix forward at=.34] \oldedgefanpath;
        \foreach \angle in {270,315,360}{
          \node[chsix red point] at
            ({\edgefanradius*cos(\angle)},{\edgefanradius*sin(\angle)}) {};
        }
        \ifnum\row=1
          \node[chsix red point] at (0,\edgefanradius) {};
        \fi
        \node[text=annred,fill=white,inner sep=1pt] at (-1.48,-.37) {$\gamma$};
        \pgfmathtruncatemacro{\crossingcount}{3+\row}
      \else
        \draw[replaced path] \oldedgefanpath;
        \draw[chsix red curve,chsix forward at=.40] \newedgefanpath;
        \node[chsix red point] at
          ({\edgefanradius*cos(135)},{\edgefanradius*sin(135)}) {};
        \ifnum\row=0
          \node[chsix red point] at (0,\edgefanradius) {};
        \fi
        \node[text=annred,fill=white,inner sep=1pt] at (-1.48,.38) {$\gamma'$};
        \pgfmathtruncatemacro{\crossingcount}{2-\row}
      \fi
      \draw[chsix red curve,-{Stealth[length=1.65mm,width=1.10mm]}]
        (splice)--(tail);
      \node[chsix point] at (splice) {};
      \node[chsix point,label={below left:$A$}] at (A) {};
      \node[chsix point,label={[fill=white,inner sep=.6pt]below left:$D$}] at (0,0) {};
      \node[chsix point] at (-1,0) {};
      \node[below] at (-1,-.33) {$Q_1$};
      \node[above=2pt] at (tail) {to $B$};
      \node at (-.50,-.76) {$\crossingcount$ edges};
    \end{scope}
  }
  \draw[chsix transition arrow] (1.23,{.25-2.95*\row})--(1.77,{.25-2.95*\row});
}
\end{tikzpicture}
\caption{Shortening detours when the incoming path follows an edge}
\label{fig:fresh-edge-long-detours}
\end{figure}

As before, write the connecting products as $U_-U_+^aU_-$ and $U_+U_-^bU_+$. The sign counts excluding the two ends satisfy $a-b\geq1$. For $n\geq2$, substituting $u=N-x$ in \eqref{eq:fresh-long-fan} gives
\[
 r_0U_-U_+^aU_--t_0U_+U_-^bU_+
 =\bigl((a-b-1)N,(a-b+1)N\bigr)\geq0
\]
The second entry is positive. For $n=1$, omitting the initial triangle sign gives the row difference
\[
 \mathbf1^{\mathsf T}U_+^aU_--\mathbf1^{\mathsf T}U_-^bU_+
 =(a-b-1,a-b+1)
\]
which agrees with the preceding difference for $N=1$. If the route continues beyond the connecting triangle, multiplication by the positive column for the remaining interior word proves that the new GM length is smaller.

If the curve ends at a vertex $B$ of that triangle, connect directly so that the last entry edge is opposite $B$. When the old curve has three crossings, let $DB$ be the last edge of the new detour and $h_0$ its weight. Connecting to $B$ just before it removes the last crossing and its adjacent triangle signs. Thus for $n\geq2$, with the terminal triangle sign omitted from the interior word,
\[
 |\gamma|=r_0U_-U_+^a\mathbf1=(a+1)N+\alpha,\qquad
 |\gamma'|=t_0U_+U_-^{b-h_0-1}\mathbf1=(b-h_0+1)N-\beta
\]
Here $b-h_0-1\geq0$ is the weight of the single remaining edge of the new detour. Since $\alpha+\beta=N$, the difference is $(a-b+h_0+1)N>0$. For $n=1$, also omit the initial triangle sign: the old length is $a+2$, the new length is $b-h_0+1$, and their difference is again $a-b+h_0+1>0$. When the old curve has four crossings, $B$ is a vertex of the last triangle of the pure left push-off. Connect to $B$ as soon as that push-off first enters a triangle containing $B$. The new GM length is at most $N$, while the old length is $(a+1)N+\alpha>N$ for $n\geq2$ and $a+2>1=N$ for $n=1$. Both cases give a strict decrease.

If the old curve crosses five edges around $D$, the connecting portion on the other side has no edge crossing. Directly joining the pure left push-off to the old tail creates backtracking, so the preceding row difference cannot be used as it stands. Whenever the connection crosses the same triangulation edge twice consecutively, delete both crossings and the intervening curve portion, and join the retained portions in the adjacent triangle (Figure~\ref{fig:fresh-five-edge-cancellation}). The figure shows the original arc on the left, the initial reconnection in the middle, and the curve after removing the backtracking across the thick edge on the right. The two red points in the middle mark the deleted crossings. In the right panel, the black points are joined in the triangle to the left of the thick edge, and the removed portion is gray and dashed.

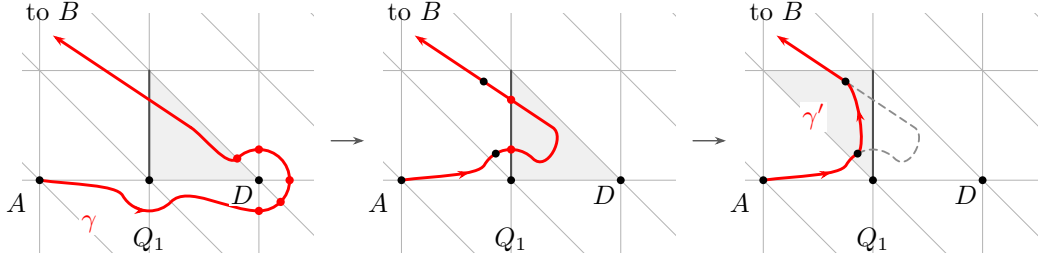
\begin{figure}[!htbp]
\centering
\begin{tikzpicture}[x=1.45cm,y=1.45cm,
  line cap=round,line join=round,every node/.style={font=\small},
  removed path/.style={draw=black!45,dash pattern=on 2.8pt off 2pt,line width=.7pt}]
\def\cancelradius{.28}
\foreach \panel in {0,1,2}{
  \begin{scope}[shift={({3.3*\panel},0)}]
    \begin{scope}
      \clip (-2.16,-.66) rectangle (.50,1.52);
      \ifnum\panel=2
        \fill[black!6] (-1,0)--(-1,1)--(-2,1)--cycle;
      \else
        \fill[black!6] (-1,0)--(0,0)--(-1,1)--cycle;
      \fi
      \foreach \i in {-2,-1,0}{
        \draw[chsix triangulation] (\i,-.8)--(\i,1.7);
      }
      \foreach \i in {0,1}{
        \draw[chsix triangulation] (-2.3,\i)--(.7,\i);
      }
      \foreach \s in {-3,-2,-1,0,1,2}{
        \draw[chsix triangulation] (-2.3,{\s+2.3})--(.7,{\s-.7});
      }
    \end{scope}
    \coordinate (A) at (-2,0);
    \coordinate (splice) at (-.62,.48);
    \coordinate (tail) at (-1.88,1.32);
    \coordinate (cut-in) at ({-1+\cancelradius*cos(120)},{\cancelradius*sin(120)});
    \coordinate (cut-out) at (-1.25,.90);
    \draw[chsix distinguished edge] (-1,0)--(-1,1);
    \ifnum\panel=0
      \draw[chsix red curve,chsix forward at=.28]
        (A) .. controls (-1.70,-.04) and (-1.30,-.04) ..
        ({-1+\cancelradius*cos(210)},{\cancelradius*sin(210)})
        arc[start angle=210,end angle=315,radius=\cancelradius]
        .. controls (-.65,-.04) and (-.28,-\cancelradius) .. (0,-\cancelradius)
        arc[start angle=270,end angle=495,radius=\cancelradius]
        .. controls (-.288,.108) and (-.50,.40) .. (splice);
      \foreach \angle in {270,315,360,450,495}{
        \node[chsix red point] at
          ({\cancelradius*cos(\angle)},{\cancelradius*sin(\angle)}) {};
      }
      \draw[chsix red curve,-{Stealth[length=1.65mm,width=1.10mm]}]
        (splice)--(tail);
      \node[text=annred,fill=white,inner sep=.8pt] at (-1.55,-.40) {$\gamma$};
    \else
      \draw[chsix red curve,chsix forward at=.68]
        (A) .. controls (-1.70,.04) and (-1.30,.04) ..
        ({-1+\cancelradius*cos(150)},{\cancelradius*sin(150)})
        arc[start angle=150,end angle=120,radius=\cancelradius];
      \ifnum\panel=1
        \def\cancellationstyle{chsix red curve}
      \else
        \def\cancellationstyle{removed path}
      \fi
      \draw[\cancellationstyle]
        (cut-in) arc[start angle=120,end angle=45,radius=\cancelradius]
        .. controls (-.67,.066) and (-.50,.40) .. (splice)--(cut-out);
      \ifnum\panel=1
        \node[chsix red point] at (-1,\cancelradius) {};
        \node[chsix red point] at (-1,{.48+2*.38/3}) {};
      \else
        \draw[chsix red curve,chsix forward at=.60]
          (cut-in) .. controls (-1.06,.29) and (-1.13,.82) .. (cut-out);
        \node[text=annred,fill=white,inner sep=.8pt] at (-1.54,.58) {$\gamma'$};
      \fi
      \draw[chsix red curve,-{Stealth[length=1.65mm,width=1.10mm]}]
        (cut-out)--(tail);
      \node[chsix point] at (cut-in) {};
      \node[chsix point] at (cut-out) {};
    \fi
    \node[chsix point,label={below left:$A$}] at (A) {};
    \node[chsix point,label={[fill=white,inner sep=.6pt]below left:$D$}] at (0,0) {};
    \node[chsix point] at (-1,0) {};
    \node[below] at (-1,-.34) {$Q_1$};
    \node[above=2pt] at (tail) {to $B$};
  \end{scope}
}
\foreach \offset in {0,3.3}{
  \draw[chsix transition arrow] ({.66+\offset},.36)--({.96+\offset},.36);
}
\end{tikzpicture}
\caption{A five-edge detour and cancellation of crossings}
\label{fig:fresh-five-edge-cancellation}
\end{figure}

Repeat if the new connection again has consecutive crossings of the same edge. Each step removes two crossings, so the procedure terminates. It stops when the connecting triangle has distinct entry and exit edges, when it reaches $B$, or when it returns to $A$. At the final end, connect to $B$ from the first remaining triangle having $B$ as a vertex. At the initial end, remove the initial sequence of triangles containing $A$ and connect so that the curve leaves through the side opposite $A$. Reconnecting the retained portions within their triangles thus gives a generalized arc satisfying the endpoint conditions. We now estimate the GM length of this generalized arc.

First suppose crossings remain and a tail can be shared with the old curve. Consider the first edge crossed after the new curve joins the old tail. Cross it at the original point in the original direction, then follow the old route to $B$. All signs from this crossing onward agree. Split both matrix products immediately before this crossing factor, and multiply each initial product on the left by $\mathbf1^{\mathsf T}=(1,1)$ to obtain its prefix row.

We will show that every entry of the new prefix row is at most $2N$, whereas every entry of the old prefix row is greater than $4N$. Multiplying these rows by the same positive tail column will then prove that GM length decreases.

For the new curve, first suppose $n\geq2$ and cancellation does not return to $A$. Before entering the connecting triangle it follows an initial portion of the pure left push-off, so its row is obtained from an initial portion of that push-off's interior product. Right multiplication by $U_+$ or $U_-$ replaces one component by the sum of the two and leaves the other unchanged, so neither component decreases. The row for the full pure left push-off has component sum $N$, hence each component of every intermediate row is at most $N$. Multiplying by the one additional factor for the connecting triangle makes both components at most $2N$. If $n=1$ or cancellation reaches $A$, the triangle containing the initial endpoint contributes no interior sign, so the comparison row remains $(1,1)$; again both components are at most $2N$.

The old curve cannot end in the triangle immediately after its five crossings. The endpoint condition would give $B=Q_{n-1}$, but it has already passed through a triangle containing $Q_{n-1}$, contradicting Lemma~\ref{lem:fresh-minimizer}(2). Thus its tail continues. The five crossings have at least four intervening triangle signs, so $a\geq4$. For $n\geq2$, the row after this portion is
\[
 (x,N-x)U_-U_+^aU_-
 =\bigl(x+aN,x+(a+1)N\bigr)
\]
Since $x>0$, both components exceed $4N$. Every subsequent factor before reaching the common tail is at least $E_2$, so this lower bound persists. For $n=1$, the old row is $(a+1,a+2)$; since $N=1$ and $a\geq4$, both components again exceed $4N$. Thus immediately before the shared tail, each old row component is strictly larger than its new counterpart. If $w$ is the remaining common interior word, its column $M(w)\mathbf1$ is positive. Multiplying both rows by it proves that the old GM length is larger.

If cancellation reaches $B$, end the new curve in a triangle entered through the side opposite $B$. Its GM length is the remaining row multiplied by $\mathbf1$, hence at most $4N$ because each entry is at most $2N$. If all crossings disappear, its length is $1$, satisfying the same bound. The old length is the row immediately after its five crossings multiplied by the positive integer column $M(w)\mathbf1$ for the remaining interior word. Both row entries exceed $4N$, so this product also exceeds $4N$. Thus the new length is smaller even when cancellation reaches the endpoint.

\smallskip\noindent
\textbf{(\ref{case:fresh-bend-exit-edge}) Crossing the edge followed by the outgoing segment from $D$.}
Finally suppose $\gamma$ detours to the right of $D$ through more than a half-turn, and the shortest polygonal path follows a triangulation edge immediately after $D$. Assume $\gamma$ crosses this edge incident to $D$ during the detour. It crosses at most once because triangles are not revisited. Let $\Delta$ be the triangle entered at that crossing, shaded in Figure~\ref{fig:fresh-exit-edge-connection}. The arc cannot end at $B$ within $\Delta$: the shortest path would then proceed straight from $D$ to $B$ in $\Delta$, making the crossed edge $DB$, whereas crossing $DB$ and then ending at $B$ violates the generalized-arc endpoint condition. Construct $\gamma'$ by passing to the left of $D$ into $\Delta$ and joining the original crossing point on the side opposite $D$. From there to $B$, retain the old route. The thick vertical edge in the figure is the one followed by the outgoing shortest path. The old and new arcs are on the left and right and share the tail beyond the black point on the opposite side.

\begin{figure}[!htbp]
\centering
\begin{tikzpicture}[x=1.90cm,y=1.90cm,
  line cap=round,line join=round,every node/.style={font=\small},
  replaced path/.style={draw=black!45,dash pattern=on 2.8pt off 2pt,line width=.7pt}]
\def\exitedgeradius{.28}
\def\oldexitedgepath{
  (A) .. controls (-.65,-.69) and (-.23,-\exitedgeradius) .. (0,-\exitedgeradius)
  arc[start angle=270,end angle=450,radius=\exitedgeradius]
  .. controls (-.16,\exitedgeradius) and (-.24,.60) .. (shared)
}
\foreach \panel in {0,1}{
  \begin{scope}[shift={({2.85*\panel},0)}]
    \begin{scope}
      \clip (-1.16,-1.16) rectangle (.78,1.43);
      \fill[black!6] (0,0)--(0,1)--(-1,1)--cycle;
      \foreach \i in {-1,0,1}{
        \draw[chsix triangulation] (\i,-1.2)--(\i,1.5);
        \draw[chsix triangulation] (-1.2,\i)--(.85,\i);
      }
      \foreach \s in {-2,-1,0,1,2}{
        \draw[chsix triangulation] (-1.2,{\s+1.2})--(.85,{\s-.85});
      }
    \end{scope}
    \coordinate (A) at (-1,-1);
    \coordinate (shared) at (-.24,1);
    \draw[chsix reference segment] (A)--(0,0)--(0,1.36);
    \draw[chsix distinguished edge] (0,0)--(0,1)--(-1,1);
    \ifnum\panel=0
      \draw[chsix red curve,chsix forward at=.56] \oldexitedgepath;
      \node[chsix red point] at (0,\exitedgeradius) {};
      \node[text=annred,fill=white,inner sep=1pt] at (.48,.10) {$\gamma$};
    \else
      \draw[replaced path] \oldexitedgepath;
      \draw[chsix red curve,chsix forward at=.51]
        (A) .. controls (-.24,-.24) and (-.24,.38) .. (shared);
      \node[text=annred,fill=white,inner sep=1pt] at (-.57,.27) {$\gamma'$};
    \fi
    \draw[chsix red curve,-{Stealth[length=1.65mm,width=1.10mm]}]
      (shared)--(-.24,1.37);
    \node[chsix point] at (shared) {};
    \node[chsix point,label={[fill=white,inner sep=.6pt]below left:$A$}] at (A) {};
    \node[chsix point,label={[fill=white,inner sep=.6pt]below right:$D$}] at (0,0) {};
    \node at (-.68,.77) {$\Delta$};
    \node[above=3pt] at (-.24,1.37) {to $B$};
    \node[right,fill=white,inner sep=.6pt] at (.04,1.22) {$g$};
  \end{scope}
}
\draw[chsix transition arrow] (1.00,.18)--(1.40,.18);
\end{tikzpicture}
\caption{Joining the same point on the opposite side}
\label{fig:fresh-exit-edge-connection}
\end{figure}
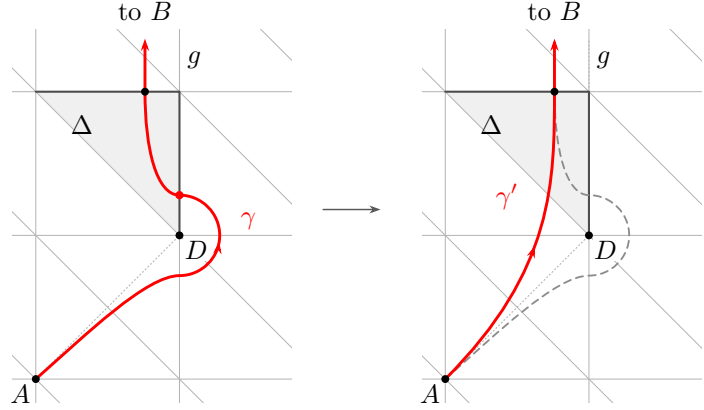

If $AD$ does not follow triangulation edges, the canonical-midpoint argument from (\ref{case:fresh-bend-right}) again shows that the arc up to just before $D$ agrees with the pure right push-off.

Write $\eta\in\{+,-\}$ for the midpoint sign where the old curve crosses the edge followed by the shortest path to enter $\Delta$. The comparison below applies to either choice of this sign. If its weight is $h$, the triangle immediately before the crossing has sign $+$ because it cuts off $D$ on the left, while the passage in $\Delta$ toward the side opposite $D$ has sign $-$. Split the products immediately before crossing that opposite side. The old row has the entrywise lower bound
\[
 r_0U_-U_+^aU_\eta^hU_-\ \geq\ r_0U_-U_+^aU_-
\]
Compare this lower bound with the new row. Here $a$ counts the signs around $D$ after removing the two end triangle signs and the factor $U_\eta^h$. We used $U_\eta^h\geq E_2$, valid for both $\eta=+$ and $\eta=-$. If $AD$ also follows triangulation edges and $n=1$, omit the initial $U_-$, which is the initial endpoint triangle sign. The lower bound is then $\mathbf1^{\mathsf T}U_+^aU_\eta^hU_-\geq\mathbf1^{\mathsf T}U_+^aU_-$.

Apply $T(x,y)=(-y,x+y)$ from (\ref{case:fresh-bend-left}) as many times as needed, permuting edge weights accordingly. Denote the transformed horizontal, slope-$-1$, and vertical weights by $k_H,k_D,k_V$, so $K_0=k_H+k_D+k_V$. Let $j$ count all old crossings around $D$, including the outgoing edge under consideration. The new arc enters $\Delta$ without crossing that edge and joins $\gamma$ on the side opposite $D$. Write its connecting product as $U_+U_-^bU_+$. Counting the edge weights and triangle signs around the lattice point gives the following values. ``Horizontal'' means the transformed case where $AD$ follows triangulation edges.
\[
\begin{array}{c|c|c|c|c}
\text{Incoming direction}&\text{Outgoing direction}&j&a&b\\ \hline
\text{Nonedge}&(0,1)&4&K_0+3&k_H+k_D+1\\
\text{Nonedge}&(-1,1)&5&K_0+k_V+4&k_H\\
\text{Horizontal}&(0,1)&4&K_0+3&k_D\\
\text{Horizontal}&(-1,1)&5&K_0+k_V+4&\text{No crossing}
\end{array}
\]
For a nonedge incoming direction, an outgoing direction $(-1,0)$ is excluded: crossing that edge would return to the triangle from which the arc approached $D$. No other outgoing triangulation-edge directions are possible when the angle around $D$ is greater than $\pi$ and less than $2\pi$.

For the first three rows, $a-b$ is respectively $k_V+2$, $k_D+2k_V+4$, and $k_H+k_V+3$. The first two satisfy $a-b\geq2$, so \eqref{eq:fresh-long-fan} applies. The third satisfies $a-b\geq1$, so the row comparison for a horizontal incoming path applies. In the last row the new connection has no crossing; remove consecutive crossings of the same edge near the connection as in the five-edge case above. Since $a\geq4$, each old row entry exceeds $4N$, whereas each new entry is at most $2N$. These comparisons use a lower bound for the old row, so they remain valid after restoring the omitted factor $U_\eta^h$. If a common tail remains, multiply both rows by its positive interior-word column to obtain a smaller new length. If cancellation in the last row reaches the endpoint, the new length is at most $4N$ and the old length exceeds $4N$, exactly as before. Thus all four rows give a strict decrease.

In every case, the construction gives a generalized arc of strictly smaller GM length. This proves the lemma.
\end{proof}

\begin{figure}[htbp]
\centering
\begin{tikzpicture}[x=1.1cm,y=1.1cm,line cap=round,line join=round,
  every node/.style={font=\small}]
  \coordinate (A) at (0,0);
  \coordinate (Q) at (2,0);
  \coordinate (D) at (4,0);
  \coordinate (B) at (5.1,2.2);
  \draw[gray,dashed] (A)--(D)--(B);
  \draw[annred,thick,-{Stealth[length=1.8mm]}]
    (A)--(1.8,0)
    arc[start angle=180,end angle=360,radius=.2]
    --(3.8,0)
    arc[start angle=180,end angle=423.435,radius=.2]
    --(B);
  \draw[gray!75,-{Stealth[length=1.3mm]}]
    (1.62,0) arc[start angle=180,end angle=360,radius=.38];
  \draw[gray!75,-{Stealth[length=1.3mm]}]
    (3.52,0) arc[start angle=180,end angle=423.435,radius=.48];
  \node[below] at (2,-.42) {$\theta_1=\pi$};
  \node[below] at (4,-.53) {$\pi<\theta_2<2\pi$};
  \foreach \P/\lab in {A/A,Q/Q_1,D/D,B/B}{
    \fill (\P) circle[radius=1.3pt];
    \node[above left,inner sep=3pt] at (\P) {$\lab$};
  }
  \node[above] at (1,0.08) {$g$};
\end{tikzpicture}
\caption{A shortest polygonal path and its angles around lattice points}
\label{fig:line}
\end{figure}
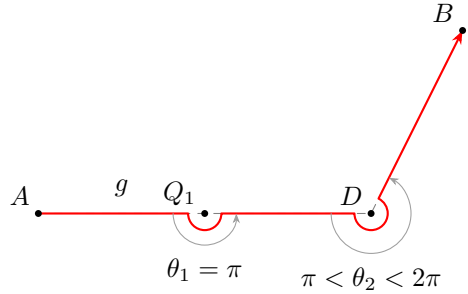

\begin{proof}[Proof of Theorem~\ref{thm:gm-distance}]
Using the preceding lemmas, we show that a minimizing arc is obtained as a push-off of $AB$. Choose a generalized arc $\gamma$ of minimum GM length, and among these one with the fewest crossings. By Lemma~\ref{lem:fresh-minimizer}, $|\gamma|=d(A,B)$ and no triangle is visited twice. Only the first passage triangle has $A$ as a vertex, and only the last has $B$ as a vertex. Since $p,q>0$, the endpoints do not belong to a single triangulation triangle, so there is at least one crossing. Form the closed triangle strip $S$ from its passage sequence and let $g$ be its shortest $A$--$B$ polygonal path. By Lemma~\ref{lem:fresh-canonical-geodesic}(2), replace $\gamma$ by a generalized arc following $g$ with canonical midpoint signs. Its crossing count is unchanged and its GM length does not increase, so it remains $d(A,B)$. Denote the replacement again by $\gamma$.

We show that the planar image of $g$ is the segment $AB$. By Lemma~\ref{lem:fresh-canonical-geodesic}(1), it can bend only at boundary vertices; if there are no intermediate boundary vertices, the conclusion is immediate. Otherwise examine the boundary vertices it meets in order, distinguishing different copies in $S$ even if they map to the same lattice point. At the $i$th vertex, rotate the ray pointing back along the incoming segment to the ray along the outgoing segment, following the side on which $\gamma$ avoids the vertex. Let $\theta_i$ be this angle, positive counterclockwise. For straight continuation, a right detour has $\theta_i=\pi$ and a left detour has $\theta_i=-\pi$. Figure~\ref{fig:line} shows the planar image of $g$ dashed and the corresponding generalized arc solid. If $|\theta_i|<\pi$, a chord across the triangles traversed around that vertex would shorten $g$, contradicting Euclidean minimality. Thus $|\theta_i|\geq\pi$. If $|\theta_i|>2\pi$, a triangle would be revisited, so $|\theta_i|\leq2\pi$.

If the first detour is to the left, use the reflection from Lemma~\ref{lem:fresh-first-bend} to exchange left and right. This preserves GM length, positivity of $p,q$, and their relative primality, so assume the first detour is to the right. If every $\theta_i=\pi$, the path continues straight and is $AB$. Otherwise let $D$ be the first vertex with $\theta_i\neq\pi$. Up to $D$, the path is straight and every intermediate lattice point is bypassed by a right half-turn. If $\theta_i<0$, the side changes to the left; if $\pi<\theta_i<2\pi$, the right detour exceeds a half-turn. In either case Lemma~\ref{lem:fresh-first-bend} produces an $A$--$B$ arc of smaller GM length, contradicting minimality.

The remaining possibility is $\theta_i=2\pi$, which also contradicts the minimizing-arc properties. If $AD$ is not along triangulation edges, a full turn around $D$ revisits the triangle from which it approached. If it is along edges, after the full turn $g$ returns along that line toward $A$, and the generalized arc passes on the side opposite the outward route. At every lattice point between $A$ and $D$, the outward route already used the three triangles in one half-plane. To avoid them while turning through at least $\pi$ on the return route, the arc must use the three triangles in the opposite half-plane. Thus the return path also goes straight through each such point and continues along the same line toward $A$. Since $B=(q,p)$ lies on none of the horizontal, vertical, or slope-$-1$ lines through $A$, it cannot reach $B$ along the way. It must then pass through a second triangle having $A$ as a vertex, contradicting Lemma~\ref{lem:fresh-minimizer}(2). No vertex with $\theta_i\neq\pi$ exists, and the planar image of $g$ is $AB$.

Finally determine the remaining choices for a generalized arc following this segment. Since $p,q$ are relatively prime, the open segment $AB$ contains no lattice point. If $tB$ is a triangulation-edge midpoint, then $2tp,2tq\in\mathbb Z$. Expressing $1$ as an integer linear combination of $p,q$ gives $2t\in\mathbb Z$, so $0<t<1$ implies $t=1/2$. At least one of $p,q$ is odd, so $B/2$ is indeed a triangulation-edge midpoint, and it is the only possible point requiring a choice of detour side. If that edge has weight $0$, avoiding the midpoint is unnecessary and either side contributes no crossing sign. Choosing the right or left side gives the same sign sequence as the pure right or left push-off, respectively. Lemma~\ref{lem:same-length-lemma0} makes their lengths equal, so
\[
d(A,B)=|\gamma|=|\gamma^R_{AB}|=|\gamma^L_{AB}|.
\]
\end{proof}

Theorem~\ref{thm:gm-distance} computes GM distance as the GM length of a push-off of a segment. The next theorem identifies this value with the generalized Markov number labeled by $p/q$, including horizontal and vertical directions. Its proof uses the generalized Cohn matrices introduced in the next chapter and is therefore deferred to Section~\ref{sec:description-cohn}.

\begin{theorem}\label{thm:length-is-number}\index{Generalized Markov distance}\index{Generalized Markov length}
Fix $(k_1,k_2,k_3)\in\mathbb Z_{\geq0}^3$ and $\sigma\in\mathfrak S_3$. Let $p,q$ be relatively prime nonnegative integers, not both zero, and put $A=(0,0)$, $B=(q,p)$. Set $t:=p/q$, with $t=\infty$ if $q=0$. If $m_t$ is the generalized Markov number defined in Section~\ref{sec:farey-label} for these parameters and fraction label $t$, then
\[
|\gamma^R_{AB}|=|\gamma^L_{AB}|=m_t.
\]
In particular,
\[
d(A,B)=m_t.
\]
\end{theorem}

\chapter{Generalized Cohn Matrices}\label{chap:generalized-cohn-matrices}
In Chapter~\ref{chap:fence-posets-gm-distance}, we described generalized Markov numbers using combinatorial and geometric objects: fence posets, GM lengths of curves, and GM distances. In this chapter we express these objects by $2\times2$ matrices and study generalized Markov numbers from a matrix-theoretic viewpoint. The matrices introduced for this purpose are generalized Cohn matrices. They generalize the Cohn matrices of the classical theory and express, within the same matrices, the generalized Markov numbers and the quantities associated with curves introduced in the preceding chapters.

We first define generalized Cohn trees and systematically construct the generalized Cohn matrix associated with each reduced fraction. We then describe their entries explicitly in terms of generalized Markov numbers and characteristic numbers. After proving relations among characteristic numbers, we introduce generalized strongly admissible sequences and show that generalized Cohn matrices are products of elementary matrices. Thus the generalized Markov equation and its solutions from Chapter~\ref{chap:generalized-markov-equations-numbers}, and the fence posets, GM lengths, and GM distances from Chapter~\ref{chap:fence-posets-gm-distance}, are expressed by the same matrices.

Generalized Markov numbers, characteristic numbers, and continued fractions appear in both the entries and the elementary matrix decompositions of generalized Cohn matrices. The matrix representations obtained here will also be used in the next chapter to define generalized discrete Markov spectra and realize their values explicitly.

This chapter is based on \cite{gyo-maru,gyoda-maruyama-sato,gyoda-generalized}. We use an equivalent normalization of generalized Cohn matrices compatible with our continued-fraction matrices.
\section{Definitions and Examples}
Fix $k_1,k_2,k_3\in\mathbb Z_{\geq0}$ and $\sigma\in\mathfrak S_3$. Define $C_{\frac{0}{1}},C_{\frac{1}{1}},C_{\frac{1}{0}}$ by
\begin{align*}
C_{\frac{0}{1}} &=
\begin{bmatrix}
3+k_1+k_2+k_3 & -(3+k_1+k_2+k_3)k_{\sigma(1)}-1\\
1 & -k_{\sigma(1)}
\end{bmatrix},\\
C_{\frac{1}{1}} &=
\begin{bmatrix}
(3+k_1+k_2+k_3)(k_{\sigma(2)}+2)-k_{\sigma(2)}-1 & 2+k_1+k_2+k_3\\
k_{\sigma(2)}+2 & 1
\end{bmatrix},\\
C_{\frac{1}{0}} &=
\begin{bmatrix}
2+k_1+k_2+k_3-k_{\sigma(3)} & 1+k_1+k_2+k_3-k_{\sigma(3)}\\
1 & 1
\end{bmatrix}.
\end{align*}
For every $(r,t,s)\in\mathrm{F}\mathbb T$, define recursively
\[
C_{r\oplus t}:=C_rC_t-D_s,\quad
C_{t\oplus s}:=C_tC_s-D_r
\]
where
\[
D_r=
\begin{bmatrix}
k_{i_r}& k_{i_r}(3+k_1+k_2+k_3)\\
0 & k_{i_r}
\end{bmatrix}
\]
\index{Generalized Cohn matrix}The matrix $C_t$ is called the \emph{$(k_1,k_2,k_3,\sigma)$-generalized Cohn matrix}, or \textbf{GC matrix}. \index{Generalized Cohn tree}We also define the \emph{$(k_1,k_2,k_3,\sigma)$-generalized Cohn tree}, or \textbf{GC tree}, by
\[
\mathrm{Co}\mathbb T(k_1,k_2,k_3,\sigma):=\mathrm{F}\mathbb T\big|_{(r,t,s)\mapsto(C_r,C_t,C_s)}.
\]
Each vertex of this GC tree is called a \emph{$(k_1,k_2,k_3,\sigma)$-generalized Cohn triple}, or \textbf{GC triple}.
\begin{exam}
The first few vertices of $\mathrm{Co}\mathbb T(1,2,0,\mathrm{id})$ are as follows:
\begin{center}
\adjustbox{max width=.98\linewidth}{$
\begin{xy}
(0,0)*+{\left(\begin{bmatrix}6&-7\\1&-1\end{bmatrix},\begin{bmatrix}21&5\\4&1\end{bmatrix},\begin{bmatrix}5&4\\1&1\end{bmatrix}\right)}="1",
(30,-14)*+{\left(\begin{bmatrix}6&-7\\1&-1\end{bmatrix},\begin{bmatrix}98&23\\17&4\end{bmatrix},\begin{bmatrix}21&5\\4&1\end{bmatrix}\right)}="2",
(30,14)*+{\left(\begin{bmatrix}21&5\\4&1\end{bmatrix},\begin{bmatrix}109&83\\21&16\end{bmatrix},\begin{bmatrix}5&4\\1&1\end{bmatrix}\right)}="3",
(95,-7)*+{\left(\begin{bmatrix}98&23\\17&4\end{bmatrix},\begin{bmatrix}2149&507\\373&88\end{bmatrix},\begin{bmatrix}21&5\\4&1\end{bmatrix}\right)}="4",
(95,-21)*+{\left(\begin{bmatrix}6&-7\\1&-1\end{bmatrix},\begin{bmatrix}467&98\\81&17\end{bmatrix},\begin{bmatrix}98&23\\17&4\end{bmatrix}\right)}="5",
(95,21)*+{\left(\begin{bmatrix}109&83\\21&16\end{bmatrix},\begin{bmatrix}626&507\\121&98\end{bmatrix},\begin{bmatrix}5&4\\1&1\end{bmatrix}\right)}="6",
(95,7)*+{\left(\begin{bmatrix}21&5\\4&1\end{bmatrix},\begin{bmatrix}2394&1823\\457&348\end{bmatrix},\begin{bmatrix}109&83\\21&16\end{bmatrix}\right)}="7",
\ar@{-}"1";"2"\ar@{-}"1";"3"\ar@{-}"2";"4"\ar@{-}"2";"5"\ar@{-}"3";"6"\ar@{-}"3";"7"
\end{xy}
$}
\end{center}
\end{exam}
\section{Entries in Terms of Generalized Markov Numbers and Characteristic Numbers}
Write $k_t$ for $k_{i_t}$. The aim of this section is to prove the following explicit description of the entries of generalized Cohn matrices.
\begin{theorem}\label{thm:Mt-description}\index{Characteristic number}\index{Generalized Cohn matrix}\index{GC matrix|seeonly{Generalized Cohn matrix}}
For every reduced fraction $t\in [0,\infty]$,
\[
C_t = \begin{bmatrix}
(3+k_1+k_2+k_3)m_t - k_t - u_t & \dfrac{(3+k_1+k_2+k_3)m_t u_t - k_t u_t - u_t^2 - 1}{m_t} \\
m_t & u_t
\end{bmatrix}.
\]
\end{theorem}
For brevity, put $K:=3+k_1+k_2+k_3$. We begin with the trace and determinant.
\begin{prop}\label{thm:Cohn-well-defined}
For every reduced fraction $t\in [0,\infty]$, the following hold.
\begin{itemize}
\item[(1)]$\mathrm{tr}(C_t)=K(C_t)_{21}-k_t$,
\item[(2)]$C_{t}\in SL(2,\mathbb Z)$.
\end{itemize}
\end{prop}

We prepare two lemmas for the proof of Proposition~\ref{thm:Cohn-well-defined}.

\begin{lemm}\label{lem:basic-property-trace}
Let $A,B\in SL(2,\mathbb Z)$. Then:
\begin{itemize}
\item[(1)]$\mathrm{tr}(A)=\mathrm{tr}(A^{-1})$
\item[(2)]$\mathrm{tr}(AB)=\mathrm{tr}(A)\mathrm{tr}(B)-\mathrm{tr}(AB^{-1})$
\item[(3)]$A^2=\mathrm{tr}(A)A-E_2$, where $E_2$ is the $2\times2$ identity matrix.
\end{itemize}
\end{lemm}
\begin{proof}
Assertion (3) follows from Corollary~\ref{cor:CH} in the appendix. Multiplying it by $A^{-1}$ gives
\[
A+A^{-1}=\mathrm{tr}(A)E_2
\]
Taking traces proves (1). Next replace $A$ by $B$ in (3),
multiply by $B^{-1}$ and then by $A$ on the left, and take traces to obtain
\[
\mathrm{tr}(AB)+\mathrm{tr}(AB^{-1})
=\mathrm{tr}(A)\mathrm{tr}(B)
\]
Rearranging gives (2).
\end{proof}
\begin{lemm}\label{lem:MtM}
For $M\in SL(2,\mathbb Z)$ satisfying $\mathrm{tr}(M)=Km_{21}-k_t$, the following identities hold:
\begin{align*}
M\begin{bmatrix}0&K\\ 0&0\end{bmatrix}M&=(\mathrm{tr}(M)+k_t)M+\begin{bmatrix}0&K\\ 0&0\end{bmatrix},\\
M^{-1}\begin{bmatrix}0&K\\ 0&0\end{bmatrix}M^{-1}&=-(\mathrm{tr}(M^{-1})+k_t)M^{-1}+\begin{bmatrix}0&K\\ 0&0\end{bmatrix}.
\end{align*}
\end{lemm}

\begin{proof}
We first prove the first identity. Write $M=\begin{bsmallmatrix}m_{11}&m_{12}\\ m_{21}&m_{22}\end{bsmallmatrix}$. Then
\begin{align*}
M\begin{bmatrix}0&K\\ 0&0\end{bmatrix}M
&=M\begin{bmatrix}1\\0\end{bmatrix}\begin{bmatrix}0&K\end{bmatrix}M
=\begin{bmatrix}m_{11}\\ m_{21}\end{bmatrix}\begin{bmatrix}Km_{21}&Km_{22}\end{bmatrix}\\
&=K\begin{bmatrix}m_{11}m_{21}&m_{11}m_{22}\\ m_{21}^2&m_{21}m_{22}\end{bmatrix}
=Km_{21}\begin{bmatrix}m_{11}&m_{12}\\ m_{21}&m_{22}\end{bmatrix}
+\begin{bmatrix}0&K\\ 0&0\end{bmatrix}\\
&=(\mathrm{tr}(M)+k_t)M+\begin{bmatrix}0&K\\ 0&0\end{bmatrix}.
\end{align*}
Here we used $m_{11}m_{22}-m_{21}m_{12}=1$.
For the second identity, we have

\[
M^{-1}=\begin{bmatrix}m_{22}&-m_{12}\\-m_{21}&m_{11}\end{bmatrix}
\]

and hence

\begin{align*}
M^{-1}\begin{bmatrix}0&K\\0&0\end{bmatrix}M^{-1}
&=K\begin{bmatrix}-m_{21}m_{22}&m_{11}m_{22}\\m_{21}^2&-m_{11}m_{21}\end{bmatrix}.
\end{align*}

On the other hand, $\mathrm{tr}(M^{-1})=\mathrm{tr}(M)=Km_{21}-k_t$ gives

\[
-(\mathrm{tr}(M^{-1})+k_t)M^{-1}
+\begin{bmatrix}0&K\\0&0\end{bmatrix}
=K\begin{bmatrix}
-m_{21}m_{22}&m_{12}m_{21}+1\\
m_{21}^2&-m_{11}m_{21}
\end{bmatrix}.
\]

These expressions also agree, since $m_{11}m_{22}-m_{12}m_{21}=1$.
\end{proof}
\begin{proof}[Proof of Proposition~\ref{thm:Cohn-well-defined}]
We prove the trace formula and the determinant assertion simultaneously by induction on the distance from the initial vertex of the GC tree. Both assertions for $C_{\frac{0}{1}},C_{\frac{1}{1}},C_{\frac{1}{0}}$ follow by direct calculation. Suppose that the three matrices in a GC triple $(C_r,C_t,C_s)$ belong to $SL(2,\mathbb Z)$ and satisfy the trace formula. We prove both assertions for the left child; the right child is treated similarly. First we prove the trace formulas
\[\mathrm{tr}(C_{r\oplus t})=K(C_{r\oplus t})_{21}-k_{r\oplus t},\quad \mathrm{tr}(C_{t\oplus s})=K(C_{t\oplus s})_{21}-k_{t\oplus s}.\]
Note that $k_{r\oplus t}=k_s,k_{t\oplus s}=k_r$. We prove the first equality; the second is similar. Since $C_{t}=C_{r}C_s-D_t$,
\begin{align*}
\mathrm{tr}(C_{r\oplus t})
&=\mathrm{tr}(C_rC_t-D_s)=\mathrm{tr}(C_r(C_rC_s-D_t))-2k_s
=\mathrm{tr}(C_r^2C_s)-\mathrm{tr}(C_rD_t)-2k_s\\
&\overset{\text{Lemma~\ref{lem:basic-property-trace}(2)}}{=}
\mathrm{tr}(C_r)\mathrm{tr}(C_rC_s)-\mathrm{tr}(C_rC_s^{-1}C_r^{-1})-\mathrm{tr}(C_rD_t)-2k_s\\
&=\mathrm{tr}(C_r)\mathrm{tr}(C_rC_s)-\mathrm{tr}(C_s)-\mathrm{tr}(C_rD_t)-2k_s.
\end{align*}
Now write $C_r=\begin{bsmallmatrix}
    r_{11}&r_{12}\\r_{21}&r_{22}
\end{bsmallmatrix}$. We have
\begin{align*}
\mathrm{tr}(D_tC_r^{-1})
&=\mathrm{tr}\!\left(\begin{bmatrix}r_{22}k_t-r_{21}Kk_t&\ast\\ \ast&r_{11}k_t\end{bmatrix}\right)
=k_t\mathrm{tr}(C_r)-Kk_tr_{21}
=-k_rk_t
\end{align*}
by the induction hypothesis for $C_r$. Consequently,
\begin{align*}
&\mathrm{tr}(C_r)\mathrm{tr}(C_rC_s)-\mathrm{tr}(C_rD_t)-\mathrm{tr}(C_s)-2k_s\\
&=\mathrm{tr}(C_r)\mathrm{tr}(C_rC_s)-\mathrm{tr}(C_rD_t)-\mathrm{tr}(C_s)-2k_s+k_rk_t-k_rk_t\\
&=\mathrm{tr}(C_r)\mathrm{tr}(C_rC_s)-\mathrm{tr}(C_rD_t)-\mathrm{tr}(C_r^{-1}D_t)-\mathrm{tr}(C_s)-k_rk_t-2k_s\\
&=\mathrm{tr}(C_r)\mathrm{tr}(C_rC_s)-\mathrm{tr}((C_r+C_r^{-1})D_t)-\mathrm{tr}(C_s)-k_rk_t-2k_s\\
&=\mathrm{tr}(C_r)\mathrm{tr}(C_rC_s)-\mathrm{tr}(C_r)\mathrm{tr}(D_t)-\mathrm{tr}(C_s)-k_rk_t-2k_s\\
&=\mathrm{tr}(C_r)\mathrm{tr}(C_rC_s-D_t)-\mathrm{tr}(C_s)-k_rk_t-2k_s\\
&=\mathrm{tr}(C_r)\mathrm{tr}(C_t)-\mathrm{tr}(C_s)-k_rk_t-2k_s\\
&=\left(\begin{bmatrix}
    0& K
\end{bmatrix}C_r\begin{bmatrix}
    1 \\ 0
\end{bmatrix}-k_r\right)\left(\begin{bmatrix}
    0& K
\end{bmatrix}C_t\begin{bmatrix}
    1 \\ 0
\end{bmatrix}-k_t\right)\\
&-\left(\begin{bmatrix}
    0& K
\end{bmatrix}C_s\begin{bmatrix}
    1\\ 0
\end{bmatrix}-k_s\right)-k_rk_t-2k_s\\
&=\begin{bmatrix}
    0& K
\end{bmatrix}C_r\begin{bmatrix}
    0&K\\0& 0
\end{bmatrix}C_rC_s\begin{bmatrix}
    1 \\ 0
\end{bmatrix}-k_r\begin{bmatrix}
    0& K
\end{bmatrix}C_rC_s\begin{bmatrix}
    1 \\ 0
\end{bmatrix}\\
&-k_t\begin{bmatrix}
    0& K
\end{bmatrix}C_r\begin{bmatrix}
    1 \\ 0
\end{bmatrix}-\begin{bmatrix}
    0& K
\end{bmatrix}C_s\begin{bmatrix}
    1 \\ 0
\end{bmatrix}-k_s\\
&\overset{\text{Lemma~\ref{lem:MtM}}}{=}\begin{bmatrix}
    0& K
\end{bmatrix}\left((\mathrm{tr}C_r+k_r)C_r+\begin{bmatrix}
    0&K\\0& 0\end{bmatrix}\right)C_s\begin{bmatrix}
    1 \\0 
\end{bmatrix}-k_r\begin{bmatrix}
    0& K
\end{bmatrix}C_rC_s\begin{bmatrix}
    1 \\ 0
\end{bmatrix}\\
&-k_t\begin{bmatrix}
    0& K
\end{bmatrix}C_r\begin{bmatrix}
    1 \\ 0
\end{bmatrix}-\begin{bmatrix}
    0& K
\end{bmatrix}C_s\begin{bmatrix}
    1 \\ 0
\end{bmatrix}-k_s\\
&\overset{\text{Lemma~\ref{lem:basic-property-trace} (3)}}{=}\begin{bmatrix}
    0& K
\end{bmatrix}C_r^2C_s\begin{bmatrix}
    1 \\ 0
\end{bmatrix}-k_t\begin{bmatrix}
    0& K
\end{bmatrix}C_r\begin{bmatrix}
    1 \\ 0
\end{bmatrix}-k_s\\
&=\begin{bmatrix}
    0& K
\end{bmatrix}(C_r^2C_s-C_rD_t-D_s)\begin{bmatrix}
    1 \\ 0
\end{bmatrix}-k_s\\
&=\begin{bmatrix}
    0& K
\end{bmatrix}(C_rC_t-D_s)\begin{bmatrix}
    1 \\ 0
\end{bmatrix}-k_s=K(C_{r\oplus t})_{21}-k_s.
\end{align*}
This proves the trace formula.
Next, in the same induction step, we prove that the determinant is $1$. Again, it suffices to prove $\det(C_{r\oplus t})=1$. Put $X:=C_rC_t=\begin{bsmallmatrix}x_{11}&x_{12}\\ x_{21}&x_{22}\end{bsmallmatrix}$. The trace formula just proved for $C_{r\oplus t}=X-D_s$ gives
\[
\mathrm{tr}(X)-2k_s=Kx_{21}-k_s,
\qquad\text{that is,}\qquad
\mathrm{tr}(X)=Kx_{21}+k_s
\]
Thus
\begin{align*}
\det(C_{r\oplus t})
&=\det(C_rC_t-D_s)\\
&=(x_{11}-k_s)(x_{22}-k_s)-x_{21}(x_{12}-Kk_s)\\
&=\det(X)-k_s\mathrm{tr}(X)+k_s^2+Kk_sx_{21}\\
&=1-k_s(Kx_{21}+k_s)+k_s^2+Kk_sx_{21}=1.
\end{align*}
Both assertions therefore pass to the child vertices, completing the simultaneous induction.
\end{proof}
We use this to prove the following proposition.
\begin{prop}\label{prop:C21=m}
For every reduced fraction $t\in[0,\infty]$, we have $(C_t)_{21}=m_t$.
\end{prop}
\begin{proof}
The cases $t=\frac{0}{1},\frac{1}{1},\frac{1}{0}$ follow by direct calculation. For the other cases, we induct on the distance from the initial vertex of the GC tree. Suppose that a GC triple $(C_r,C_t,C_s)$ satisfies the proposition. We must prove
 \[(C_{r\oplus t})_{21}=m_{r\oplus t}=\frac{m_{r}^2+k_sm_rm_t+m_t^2}{m_s},\quad (C_{t\oplus s})_{21}=m_{t\oplus s}=\frac{m_{t}^2+k_rm_tm_s+m_s^2}{m_r}\]
We prove the first equality; the second is similar. In the proof of Proposition~\ref{thm:Cohn-well-defined}, we obtained \[\mathrm{tr}(C_{r\oplus t})=\mathrm{tr}(C_r)\mathrm{tr}(C_t)-\mathrm{tr}(C_s)-k_rk_t-2k_s\]
(the expression immediately before the matrix calculations in that proof). Applying Proposition~\ref{thm:Cohn-well-defined} (1) and rearranging gives
 \begin{align}\label{eq:C21}(C_{r\oplus t})_{21}=K(m_rm_t)-m_rk_t-m_tk_r-m_s\end{align}
Since $(m_r,m_t,m_s)$ satisfies the GM equation, namely,
 \[m_r^2+m_t^2+m_s^2+k_rm_tm_s+k_tm_{s}m_{r}+k_{s}m_{r}m_{t}=Km_rm_tm_s\]
we can rewrite \eqref{eq:C21} as
 \[(C_{r\oplus t})_{21}=\frac{m_{r}^2+k_sm_rm_t+m_t^2}{m_s}\]
This proves the assertion.
\end{proof}
Finally, we prove the following proposition.
\begin{prop}\label{prop:C22=u}\index{Characteristic number}
For every reduced fraction $t\in[0,\infty]$, we have $(C_t)_{22}=u_t$.
\end{prop}
For this purpose, we introduce the index.
\begin{defi}\label{def:gc-index}\index{Index of a GC matrix}
For every $t\in\mathbb Q_{\geq0}\cup\{\infty\}$, the quantity
 \[I_{t}:=\frac{(C_{t})_{22}}{(C_{t})_{21}}\]
is called the \textbf{index} of $C_t$.
\end{defi}
\begin{prop}\label{fraction-labeling-inj}
The index is strictly increasing: if $s<t$, then $I_s<I_t$. Here $\frac{1}{0}$ is regarded as larger than every rational number.
\end{prop}
\begin{proof}
It suffices to prove $I_{r}<I_{t}<I_{s}$ for every Farey triple $(r,t,s)\in \mathrm{F}\mathbb T$. First we show $I_{t}<I_{s}$. Since $C_t=C_rC_s-D_t$, we have $C_r=(C_t+D_t)C_s^{-1}$. Comparing the $(2,1)$ entries gives
\[
r_{21}=s_{22}t_{21}-(t_{22}+k_t)s_{21}\le s_{22}t_{21}-t_{22}s_{21}.
\]
Hence
\[
0<\dfrac{r_{21}}{t_{21}s_{21}}\le\dfrac{s_{22}}{s_{21}}-\dfrac{t_{22}}{t_{21}}=I_{s}-I_{t}.
\]
This proves $I_{t}<I_{s}$. Next we prove $I_{r}<I_{t}$. The identity $C_t=C_rC_s-D_t$ gives $C_s=C_r^{-1}(C_t+D_t)$. Comparing the $(2,1)$ entries gives
\begin{align*}
s_{21}&=r_{11}t_{21}-r_{21}(t_{11}+k_t)\\
&=(Kr_{21}-k_r-r_{22})t_{21}-r_{21}(Kt_{21}-t_{22})\\
&=-k_rt_{21}+t_{22}r_{21}-t_{21}r_{22}\le t_{22}r_{21}-t_{21}r_{22}.
\end{align*}
Hence
\[
0<\dfrac{s_{21}}{r_{21}t_{21}}\le\dfrac{t_{22}}{t_{21}}-\dfrac{r_{22}}{r_{21}}=I_{t}-I_{r}.
\]
This proves $I_{r}<I_{t}$.
\end{proof}
\begin{lemm}\label{lem:1,1-positivity}\index{Generalized Cohn matrix!positivity of the lower-right entry}
For every reduced fraction $t\in(0,\infty]$, we have $(C_t)_{22}>0$.
\end{lemm}
\begin{proof}
By the strict monotonicity of the index (Proposition~\ref{fraction-labeling-inj}), it suffices to consider $t=\frac{1}{n}$. The $(2,2)$ entry of $C_{\frac{1}{1}}$ is $1$. The $(2,2)$ entry of $C_{\frac{1}{2}}$ is $k_{\sigma(2)}+2$, which is larger. Suppose that the $(2,2)$ entry of $C_{\frac{1}{i}}$ is positive and that the $(2,2)$ entry of $C_{\frac{1}{i+1}}$ is larger. We show that the $(2,2)$ entry of $C_{\frac{1}{i+2}}$ is larger still, and in particular positive. Write
\[
C_{\frac{1}{i}}=\begin{bmatrix}a&b\\c&d\end{bmatrix},\quad
C_{\frac{1}{i+1}}=\begin{bmatrix}a'&b'\\c'&d'\end{bmatrix},\quad
C_{\frac{1}{i+2}}=\begin{bmatrix}a''&b''\\c''&d''\end{bmatrix}.
\]
Then
\begin{align*}
C_{\frac{1}{i+1}}
&=C_{\frac{0}{1}}C_{\frac{1}{i}}-D_{\frac{1}{i-1}}\\
&=\begin{bmatrix}
K& -Kk_{\frac{0}{1}}-1\\1& -k_{\frac{0}{1}}
\end{bmatrix}
\begin{bmatrix}a&b\\c&d\end{bmatrix}
-\begin{bmatrix}k_{\frac{1}{i-1}}&Kk_{\frac{1}{i-1}}\\0&k_{\frac{1}{i-1}}\end{bmatrix}\\
&=\begin{bmatrix}\ast&Kb-Kk_{\frac{0}{1}}d-d-Kk_{\frac{1}{i-1}}\\\ast&b-k_{\frac{0}{1}}d-k_{\frac{1}{i-1}}\end{bmatrix}.
\end{align*}
Thus
\[
b'=Kb-Kk_{\frac{0}{1}}d-d-Kk_{\frac{1}{i-1}},\quad
d'=b-k_{\frac{0}{1}}d-k_{\frac{1}{i-1}}.
\]
By assumption, $b-k_{\frac{0}{1}}d-k_{\frac{1}{i-1}}\ge d>0$. Moreover,
\[
C_{\frac{1}{i+2}}=C_{\frac{0}{1}}C_{\frac{1}{i+1}}-D_{\frac{1}{i}}
\]
gives
\begin{align*}
d''-d'
&=b'-(k_{\frac{0}{1}}+1)d'-k_{\frac{1}{i}}\\
&=(K-k_{\frac{0}{1}}-1)b
+\bigl(k_{\frac{0}{1}}^2+(1-K)k_{\frac{0}{1}}-1\bigr)d
+(k_{\frac{0}{1}}+1-K)k_{\frac{1}{i-1}}
-k_{\frac{1}{i}}\\
&\ge(K-k_{\frac{0}{1}}-1)((k_{\frac{0}{1}}+1)d+k_{\frac{1}{i-1}})
+\bigl(k_{\frac{0}{1}}^2+(1-K)k_{\frac{0}{1}}-1\bigr)d
+(k_{\frac{0}{1}}+1-K)k_{\frac{1}{i-1}}
-k_{\frac{1}{i}}\\
&=\bigl(K-2-k_{\frac{0}{1}}\bigr)d-k_{\frac{1}{i}}\geq k_{\frac{1}{i+1}}+1>0
\end{align*}
For the final inequality, we used the fact that along the branch $1,\frac12,\frac13,\dots$, the position labels other than that of $\frac01$ alternate between $\sigma(2)$ and $\sigma(3)$, together with $K-2-k_{\frac01}=1+k_{\sigma(2)}+k_{\sigma(3)}$. Thus the $(2,2)$ entries remain positive and strictly increase along this branch. The strict monotonicity of the index proves the assertion for every positive reduced fraction.
\end{proof}
\begin{proof}[Proof of Proposition~\ref{prop:C22=u}]
The cases $t=\frac{0}{1},\frac{1}{0}$ follow by direct calculation. For $t\in (0,\infty)\cap \mathbb Q$, take the Farey triple $(r,t,s)$ whose middle entry is $t$. Then
\[C_s=C_r^{-1}(C_t+D_t)\]
By Proposition~\ref{thm:Cohn-well-defined} (1) and Proposition~\ref{prop:C21=m},
\[C_r=\begin{bmatrix}
    Km_r-k_r-(C_r)_{22}&\ast\\m_r&(C_r)_{22}\\
\end{bmatrix},C_t=\begin{bmatrix}
    Km_t-k_t-(C_t)_{22}&\ast\\m_t&(C_t)_{22}\\
\end{bmatrix}\]
Using these identities to calculate the $(2,1)$ entries on both sides gives
\[m_s=m_r(C_t)_{22}-m_t(k_r+(C_r)_{22})\]
Therefore
\[m_{r}(C_t)_{22}\equiv m_s\pmod{m_t}\]
Moreover, Proposition~\ref{fraction-labeling-inj} and Lemma~\ref{lem:1,1-positivity} give
\[0<I_t<I_{\frac{1}{0}}=1\]
so $0<(C_t)_{22}<m_t$. The uniqueness of $u_t$ now gives $u_t=(C_t)_{22}$.
\end{proof}
We can now finish the proof of Theorem~\ref{thm:Mt-description}.
\begin{proof}[Proof of Theorem~\ref{thm:Mt-description}]
Combine Propositions~\ref{thm:Cohn-well-defined}, \ref{prop:C21=m}, and~\ref{prop:C22=u}.
\end{proof}
\section{Relations among Characteristic Numbers}\label{section:relation-of-characteristic}
In this section we use GC matrices to prove the following proposition, postponed in Section~\ref{sec:characteristic-number}.
\par\medskip
\noindent\textbf{Proposition~\ref{prop:t-1/t-relation-gen}} (restated)\textbf{.}\,
For every reduced fraction $t\in[0,1]\cap \mathbb Q$, let $u_t$ be the characteristic number (or auxiliary endpoint value when $t=0$) with fraction label $t$ in $\mathrm M\mathbb T(k_1,k_2,k_3,\sigma)$, and put $k_t:=k_{i_t}$. Let $u_{\frac{1}{t}}^\ast$ be the characteristic number (or auxiliary endpoint value when $1/t=\infty$) with fraction label $\frac{1}{t}$ in $\mathrm M\mathbb T(k_1,k_2,k_3,\sigma^\ast)$, where $\sigma^\ast=\sigma\circ(1\ 3)$. Then
\[
u^\ast_{\frac{1}{t}}=m_t-u_t-k_t
\]
holds.
\par\medskip
\begin{lemm}\label{lem:minimal-value}
For the $(k_1,k_2,k_3)$-GM numbers with labels $t\in(0,1]\cap \mathbb Q$ in the $(k_1,k_2,k_3,\sigma)$-GM tree, the following hold.
\begin{itemize}
    \item[(1)] If $i_t=\sigma(1)$, then $m_t\geq m_{\frac{2}{3}}$.
    \item[(2)] If $i_t=\sigma(2)$, then $m_t\geq m_{\frac{1}{1}}$.
    \item[(3)] If $i_t=\sigma(3)$, then $m_t\geq m_{\frac{1}{2}}$.
\end{itemize} 
All three bounds are sharp.
\end{lemm}
\begin{proof}
As shown in the proof of Proposition~\ref{prop:maximal-middle-component}, the new middle GM number at a child vertex is larger than every GM number at its parent. The middle label at the root is $1$, proving (2). Every vertex with middle label $0<t<1$ descends from the root's left child, whose middle fraction label is $1/2$ and whose middle position label is $\sigma(3)$. Since new middle GM numbers increase along every branch, this proves (3).

To prove (1), consider the first vertex along the path from the root to the vertex with middle label $t$ whose middle position label is $\sigma(1)$. Before that step, the entry with position label $\sigma(1)$ at every vertex has fractional label $0/1$, inherited from the root. Indeed, when passing to a child, the two entries whose position labels are not replaced by the new middle entry are inherited unchanged.
Since the determinants of neighboring fractions are $\pm1$, the vertices of the Farey tree containing $0/1$ have the form
  \[
  \left(\frac01,\frac1{i+1},\frac1i\right)\qquad(i\geq0)
  \]
where $i=0$ denotes the root. The right child that replaces $0/1$ has middle label $2/(2i+1)$, and along a path with $t\leq1$ we have $i\geq1$.
Thus $m_t\geq m_{2/(2i+1)}$ for some $i\geq1$. It remains to compare these numbers. The parent of the vertex with middle label $\frac{2}{2i+1}$ is $\left(\frac{0}{1},\frac{1}{i+1},\frac{1}{i}\right)$, so
 \[m_{\frac{2}{2i+1}}=\frac{m^2_{\frac{1}{i+1}}+k_{\sigma(1)}m_{\frac{1}{i+1}}m_{\frac{1}{i}}+m^2_{\frac{1}{i}}}{m_{\frac{0}{1}}}=m^2_{\frac{1}{i+1}}+k_{\sigma(1)}m_{\frac{1}{i+1}}m_{\frac{1}{i}}+m^2_{\frac{1}{i}}\]
For every $i$, a triple of the form $((m_{\frac{0}{1}},\sigma(1)),(m_{\frac{1}{i+1}},\alpha),(m_{\frac{1}{i}},\beta))$ is a vertex of the $(k_1,k_2,k_3,\sigma)$-GM tree, and along the boundary branch $m_{\frac{1}{i}}<m_{\frac{1}{i+1}}$. Since all terms are nonnegative, the expression above is strictly increasing in each variable. Hence for every $i\geq1$,
 \[
 m_{\frac{2}{2i+1}}<m_{\frac{2}{2i+3}}
 \]
Thus $m_{2/3}$ is the minimum, proving (1). Each of the three lower bounds is attained at the displayed fraction, so the bounds are sharp.
\end{proof}
\begin{proof}[Proof of Proposition~\ref{prop:t-1/t-relation-gen}]
For $t=0$, both sides equal $1$ by the endpoint conventions. Suppose $0<t\leq1$, and let $(r,t,s)$ be the Farey triple with middle entry $t$. Put $x:=m_t-u_t-k_t$. By Corollary~\ref{cor:dual-remark} and the uniqueness of characteristic numbers, it suffices to show that $0<x<m_t$ and $m_sx\equiv m_r\pmod{m_t}$.

The GC triple satisfies
\[
C_r=(C_t+D_t)C_s^{-1}.
\]
Using Theorem~\ref{thm:Mt-description} to compare the $(2,1)$ entries and rearranging gives the identity
\[
m_sx=m_r+m_t(m_s-u_s).
\]
Since $s>0$, we have $0<u_s\leq m_s$: the inequality is strict for an interior label, and $u_\infty=m_\infty=1$. Thus the right-hand side is positive, so $x>0$. Also, $u_t>0$ and $k_t\geq0$ give $x<m_t$. Reducing the same identity modulo $m_t$ yields the required congruence. Hence $u^\ast_{1/t}=x=m_t-u_t-k_t$.
\end{proof}
\section{Generalized Cohn Matrices from Generalized Strongly Admissible Sequences}\label{sec:description-cohn}
In this section we show that generalized Cohn matrices admit factorizations into elementary matrices.

\begin{figure}[ht]
    \centering
   \begin{tikzpicture}[x=0.01cm,y=0.01cm,scale=0.5]
  \definecolor{annblue}{RGB}{255,0,0}

  \tikzset{
    gridline/.style={gray,line width=0.15pt,line cap=butt,line join=miter},
    gammaL/.style={draw=annblue,line width=1pt,line cap=round,line join=round},
    enddot/.style={circle,fill=black,inner sep=0pt,minimum size=1.8pt}
  }

  \begin{scope}[shift={(300,170)}, x={(260,0)}, y={(0,260)}]

    \begin{scope}
      \clip (-1,0) rectangle (5,2);

      \foreach \i in {-1,...,4}{
        \foreach \j in {0,1}{
          \draw[gridline] (\i,{\j+1}) -- ({\i+1},\j);
        }
      }

      \foreach \i in {-1,...,5}{
        \draw[gridline] (\i,0) -- (\i,2);
      }
      \foreach \j in {0,1,2}{
        \draw[gridline] (-1,\j) -- (5,\j);
      }
    \end{scope}

    \coordinate (A) at (0,0);
    \coordinate (B) at (5,2);

    \def\xleftshift{0.18}

    \coordinate (ALshift) at ({-\xleftshift},0);
    \coordinate (BLshift) at ({5-\xleftshift},2);

    \draw[gammaL] (ALshift) -- (BLshift);

    \node[enddot,label={[font=\small]below left:$A$}] at (A) {};
    \node[enddot,label={[font=\small]above right:$B$}] at (B) {};

    \node[annblue,font=\normalsize] at (1.45,0.85) {$\overline{L_t}$};

  \end{scope}
\end{tikzpicture}
    \caption{An example of $\overline{L_t}$}
    \label{fig:ex-extendedpresnakegraph}
\end{figure}
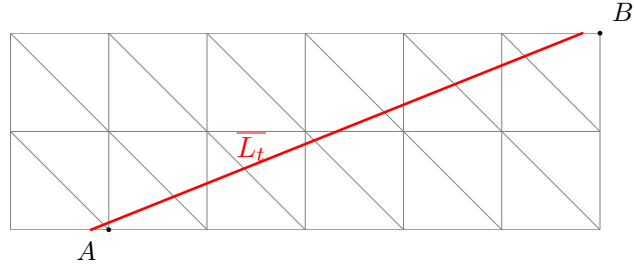

Let $t=\frac{p}{q}$ be a positive reduced fraction, and write $L_t$ for the pure left push-off $\gamma_{AB}^L$ joining $A=(0,0)$ to $B=(q,p)$. Write the corresponding segment as $\ell_t(u)=u(q,p)$, $0\leq u\leq1$.
For $0<\varepsilon<\varepsilon_t:=1/(2(p+q))$, put
\[
 L_{t,\varepsilon}(u)=u(q,p)-(\varepsilon,0)\qquad(0\leq u\leq1)
\]
\index{Admissible perturbation}This translated segment is called an \emph{admissible perturbation} in this section. Its intersection with the horizontal edge at its initial point $(-\varepsilon,0)$ is counted as an edge-crossing occurrence, whereas its contact at the terminal point $(q-\varepsilon,p)$ is not counted. For a triangle passage, apply the triangle-crossing rule to the two edges joined by the segment. A chosen admissible perturbation is denoted by $\overline{L_t}$. Figure~\ref{fig:ex-extendedpresnakegraph} shows the case $t=\frac25$. In Figures~\ref{fig:ex-extendedpresnakegraph} and~\ref{fig:ex-extendedpresnakegraph-signed}, the displacement is enlarged within a range that preserves the passage order and signs.

\begin{lemm}\label{lem:admissible-perturbation-Lt}\index{Admissible perturbation}
The admissible perturbation above avoids $\mathcal V$ and every edge midpoint, and meets each edge transversely. For every choice of $0<\varepsilon<\varepsilon_t$, the order of the triangle-passage and edge-crossing occurrences, and the signs assigned to them, are the same.
\end{lemm}
\begin{proof}
If the image passes through a half-lattice point $(i/2,j/2)$, substituting $u=j/(2p)$ gives
\[
2p\varepsilon=qj-pi\in\mathbb Z
\]
This is impossible, since $0<2p\varepsilon<p/(p+q)<1$. Thus no lattice point or edge midpoint is met. The direction $(q,p)$ is parallel to neither a horizontal edge, a vertical edge, nor an edge of slope $-1$, so all intersections are transverse.
The crossing times are $j/p$ for the horizontal line $y=j$, $(i+\varepsilon)/q$ for the vertical line $x=i$, and $(h+\varepsilon)/(p+q)$ for the diagonal line $x+y=h$. Equality of crossing times of different types would give a lattice point, and a change of the subdivided edge containing an intersection would require passing through an edge midpoint. Neither is possible. Each endpoint remains in the interior of the same horizontal edge. Thus, throughout the connected interval $0<\varepsilon<\varepsilon_t$, neither the crossing order nor the side of a midpoint changes. These data determine the sign rules, proving the claim.
\end{proof}

\begin{rema}\label{rem:one-punctured-torus}
We can also interpret $\overline{L_t}$ on a once-punctured torus. Let $\pi$ be the quotient map from $\mathbb R^2\setminus\mathbb Z^2$ to the once-punctured torus obtained by identifying integer translates. The interior of the segment $AB$ gives an arc with both ends at the puncture, and its left translate gives a closed curve. Lifting one period of this closed curve from the horizontal edge $y=0$ to $y=p$, counting the initial edge crossing but not the terminal one, gives the curve $\overline{L_t}$ and its sign sequence defined above. Figure~\ref{fig:interpretation-of-Lt} shows $t=\frac{1}{1}$.
Opposite sides with the same label $a$ or $b$ are identified in the directions indicated by the arrows; their images on the torus have the same labels. The upper row shows the arc $\pi(u,u)$ for $0<u<1$, and the lower row shows the closed curve $\pi(u-\varepsilon,u)$ for $0\leq u\leq1$, with $\varepsilon=0.15$ in the figure. The left side is defined relative to the direction of travel in the plane. Red arrows indicate the direction of travel, and blue arrows indicate the translation $(-\varepsilon,0)$. Open circles represent the puncture. The dotted line in the lower planar diagram is the original segment, and dashed parts on the torus lie on the back side.
\begin{figure}[ht]
\centering
\begin{tikzpicture}[
  font=\small,
  grid edge/.style={gray!75,line width=.4pt},
  diagonal edge/.style={gray!45,line width=.35pt},
  red curve/.style={red!85!black,line width=.95pt,line cap=round,line join=round},
  basis curve/.style={gray!70,line width=.4pt},
  hidden part/.style={dash pattern=on 1.5pt off 1.4pt},
  visible part/.style={solid},
  curve arrow/.style={-{Stealth[length=2mm,width=1.5mm]}},
  basis arrow/.style={-{Stealth[length=1.4mm,width=1mm]}},
  shift arrow/.style={blue!70!black,line width=.6pt,-{Stealth[length=1.5mm,width=1.1mm]}},
  puncture/.style={circle,draw=black,fill=white,line width=.4pt,inner sep=0pt,minimum size=3.2pt},
  declare function={
    torusX(\x,\y)=(2+.62*cos(360*\y+25))*cos(360*\x-110);
    torusY(\x,\y)=sin(35)*(2+.62*cos(360*\y+25))*sin(360*\x-110)+cos(35)*.62*sin(360*\y+25);
    contourPhi(\t)=atan2(cos(35)*sin(\t),sin(35));
    contourX(\t,\k)=(2+.62*cos(contourPhi(\t)+\k))*cos(\t);
    contourY(\t,\k)=sin(35)*(2+.62*cos(contourPhi(\t)+\k))*sin(\t)+cos(35)*.62*sin(contourPhi(\t)+\k);
  }
]
\newcommand{\torusSurface}{
  \path[fill=gray!7,draw=gray!65,line width=.45pt]
    plot[domain=0:360,samples=145,variable=\t]
      ({contourX(\t,0)},{contourY(\t,0)}) -- cycle;
  \path[fill=white,draw=gray!65,line width=.45pt]
    plot[domain=0:360,samples=145,variable=\t]
      ({contourX(\t,180)},{contourY(\t,180)}) -- cycle;
}
\newcommand{\unitFrame}{
  \draw[grid edge] (0,0) rectangle (1,1);
  \draw[diagonal edge] (0,1)--(1,0);
  \foreach \y in {0,1}{
    \draw[grid edge,basis arrow] (.32,\y)--(.58,\y);
  }
  \foreach \x in {0,1}{
    \draw[grid edge,basis arrow] (\x,.32)--(\x,.58);
  }
  \node[below] at (.70,0) {$a$};
  \node[above] at (.70,1) {$a$};
  \node[left] at (0,.70) {$b$};
  \node[right] at (1,.70) {$b$};
}
\begin{scope}[shift={(0,3.6)},x=2.4cm,y=2.4cm]
  \unitFrame
  \draw[red curve] (0,0)--(1,1);
  \draw[red curve,curve arrow] (.38,.38)--(.59,.59);
  \foreach \x/\y in {0/0,1/0,0/1,1/1}{\node[puncture] at (\x,\y) {};}
  \node[below left] at (0,0) {$A$};
  \node[above right] at (1,1) {$B$};
\end{scope}
\node at (1.2,3.06) {$(u,u)\quad(0<u<1)$};
\draw[-{Stealth[length=2mm]}] (3.05,4.8)--node[above] {$\pi$}(4.32,4.8);
\begin{scope}[x=2.4cm,y=2.4cm]
  \unitFrame
  \draw[diagonal edge] (-.30,0)--(0,0) (-.30,1)--(0,1) (-.30,.30)--(0,0);
  \draw[gray!60,densely dotted] (0,0)--(1,1);
  \draw[red curve] (-.15,0)--(.85,1);
  \draw[red curve,curve arrow] (.23,.38)--(.44,.59);
  \draw[shift arrow] (.30,.30)--(.15,.30);
  \draw[shift arrow] (.76,.76)--(.61,.76);
  \foreach \x/\y in {0/0,1/0,0/1,1/1}{\node[puncture] at (\x,\y) {};}
  \fill[red!85!black] (-.15,0) circle[radius=.7pt];
  \fill[red!85!black] (.85,1) circle[radius=.7pt];
\end{scope}
\node at (1.05,-.52) {$(u-\varepsilon,u)\quad(0\leq u\leq1)$};
\draw[-{Stealth[length=2mm]}] (3.05,1.2)--node[above] {$\pi$}(4.32,1.2);
\begin{scope}[shift={(7,4.8)},scale=.85]
  \torusSurface
  \draw[basis curve,hidden part] plot[domain=0.35849226:0.75261885,samples=71,variable=\t]
    ({torusX(\t,0)},{torusY(\t,0)});
  \draw[basis curve,visible part] plot[domain=0.00000000:0.35849226,samples=65,variable=\t]
    ({torusX(\t,0)},{torusY(\t,0)});
  \draw[basis curve,visible part] plot[domain=0.75261885:1.00000000,samples=45,variable=\t]
    ({torusX(\t,0)},{torusY(\t,0)});
  \draw[basis curve,basis arrow] plot[domain=0.16424613:0.19424613,samples=12,variable=\t]
    ({torusX(\t,0)},{torusY(\t,0)});
  \draw[basis curve,hidden part] plot[domain=0.28247629:0.78247629,samples=90,variable=\t]
    ({torusX(0,\t)},{torusY(0,\t)});
  \draw[basis curve,visible part] plot[domain=0.00000000:0.28247629,samples=51,variable=\t]
    ({torusX(0,\t)},{torusY(0,\t)});
  \draw[basis curve,visible part] plot[domain=0.78247629:1.00000000,samples=40,variable=\t]
    ({torusX(0,\t)},{torusY(0,\t)});
  \draw[basis curve,basis arrow] plot[domain=0.12623814:0.15623814,samples=12,variable=\t]
    ({torusX(0,\t)},{torusY(0,\t)});
  \draw[diagonal edge,hidden part] plot[domain=0.19013595:0.43945691,samples=45,variable=\t]
    ({torusX(\t,1-\t)},{torusY(\t,1-\t)});
  \draw[diagonal edge,visible part] plot[domain=0.00000000:0.19013595,samples=35,variable=\t]
    ({torusX(\t,1-\t)},{torusY(\t,1-\t)});
  \draw[diagonal edge,visible part] plot[domain=0.43945691:1.00000000,samples=101,variable=\t]
    ({torusX(\t,1-\t)},{torusY(\t,1-\t)});
  \draw[red curve,hidden part] plot[domain=0.58227223:0.85973768,samples=50,variable=\t]
    ({torusX(\t,\t)},{torusY(\t,\t)});
  \draw[red curve,visible part] plot[domain=0.00000000:0.58227223,samples=105,variable=\t]
    ({torusX(\t,\t)},{torusY(\t,\t)});
  \draw[red curve,visible part] plot[domain=0.85973768:1.00000000,samples=26,variable=\t]
    ({torusX(\t,\t)},{torusY(\t,\t)});
  \draw[red curve,curve arrow] plot[domain=0.27913611:0.30313611,samples=12,variable=\t]
    ({torusX(\t,\t)},{torusY(\t,\t)});
  \draw[red curve,curve arrow] plot[domain=0.91786884:0.94186884,samples=12,variable=\t]
    ({torusX(\t,\t)},{torusY(\t,\t)});
  \node[puncture] (puncturePoint) at ({torusX(0,0)},{torusY(0,0)}) {};
  \node[below left=2pt] at (puncturePoint) {puncture};
  \node[below=1pt] at ({torusX(.15,0)},{torusY(.15,0)}) {$a$};
  \node[above left=3pt] at ({torusX(0,.23)},{torusY(0,.23)}) {$b$};
\end{scope}
\node at (7,3.06) {$\pi(u,u)$};
\begin{scope}[shift={(7,1.2)},scale=.85]
  \torusSurface
  \draw[basis curve,hidden part] plot[domain=0.35849226:0.75261885,samples=71,variable=\t]
    ({torusX(\t,0)},{torusY(\t,0)});
  \draw[basis curve,visible part] plot[domain=0.00000000:0.35849226,samples=65,variable=\t]
    ({torusX(\t,0)},{torusY(\t,0)});
  \draw[basis curve,visible part] plot[domain=0.75261885:1.00000000,samples=45,variable=\t]
    ({torusX(\t,0)},{torusY(\t,0)});
  \draw[basis curve,basis arrow] plot[domain=0.16424613:0.19424613,samples=12,variable=\t]
    ({torusX(\t,0)},{torusY(\t,0)});
  \draw[basis curve,hidden part] plot[domain=0.28247629:0.78247629,samples=90,variable=\t]
    ({torusX(0,\t)},{torusY(0,\t)});
  \draw[basis curve,visible part] plot[domain=0.00000000:0.28247629,samples=51,variable=\t]
    ({torusX(0,\t)},{torusY(0,\t)});
  \draw[basis curve,visible part] plot[domain=0.78247629:1.00000000,samples=40,variable=\t]
    ({torusX(0,\t)},{torusY(0,\t)});
  \draw[basis curve,basis arrow] plot[domain=0.12623814:0.15623814,samples=12,variable=\t]
    ({torusX(0,\t)},{torusY(0,\t)});
  \draw[diagonal edge,hidden part] plot[domain=0.19013595:0.43945691,samples=45,variable=\t]
    ({torusX(\t,1-\t)},{torusY(\t,1-\t)});
  \draw[diagonal edge,visible part] plot[domain=0.00000000:0.19013595,samples=35,variable=\t]
    ({torusX(\t,1-\t)},{torusY(\t,1-\t)});
  \draw[diagonal edge,visible part] plot[domain=0.43945691:1.00000000,samples=101,variable=\t]
    ({torusX(\t,1-\t)},{torusY(\t,1-\t)});
  \draw[red curve,hidden part] plot[domain=0.28872092:0.94523818,samples=119,variable=\t]
    ({torusX(\t-.15,\t)},{torusY(\t-.15,\t)});
  \draw[red curve,visible part] plot[domain=0.00000000:0.28872092,samples=52,variable=\t]
    ({torusX(\t-.15,\t)},{torusY(\t-.15,\t)});
  \draw[red curve,visible part] plot[domain=0.94523818:1.00000000,samples=12,variable=\t]
    ({torusX(\t-.15,\t)},{torusY(\t-.15,\t)});
  \draw[red curve,curve arrow] plot[domain=0.08424031:0.10824031,samples=12,variable=\t]
    ({torusX(\t-.15,\t)},{torusY(\t-.15,\t)});
  \draw[red curve,curve arrow] plot[domain=0.18048062:0.20448062,samples=12,variable=\t]
    ({torusX(\t-.15,\t)},{torusY(\t-.15,\t)});
  \node[puncture] (puncturePoint) at ({torusX(0,0)},{torusY(0,0)}) {};
  \node[below left=2pt] at (puncturePoint) {puncture};
  \node[below=1pt] at ({torusX(.15,0)},{torusY(.15,0)}) {$a$};
  \node[above left=3pt] at ({torusX(0,.23)},{torusY(0,.23)}) {$b$};
\end{scope}
\node at (7,-.52) {$\pi(u-\varepsilon,u)$};
\end{tikzpicture}

\caption{The torus interpretation of $\overline{L_t}$}\label{fig:interpretation-of-Lt}
\end{figure}
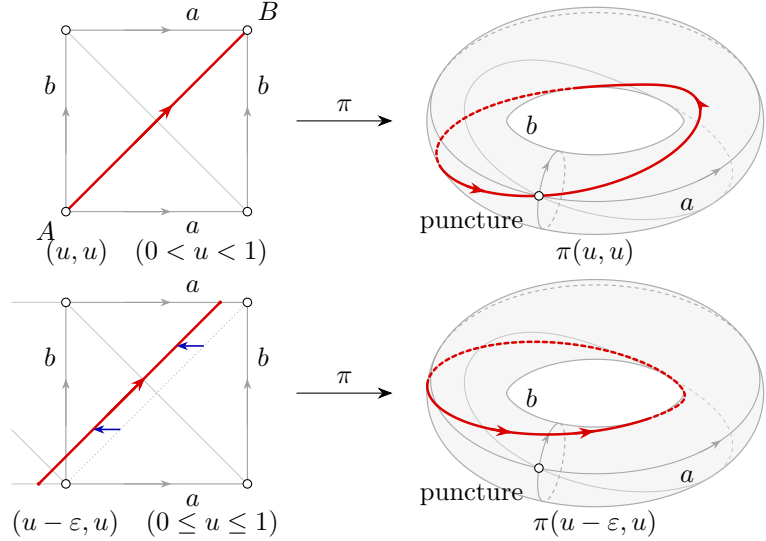
\end{rema}
Fix $(k_1,k_2,k_3)\in \mathbb Z^3_{\geq 0}$ and $\sigma\in \mathfrak S_3$. \index{Generalized strongly admissible sequence}We now define the \textbf{$(k_1,k_2,k_3,\sigma)$-generalized strongly admissible sequence} $S(t)$.
First put
\[
S(\tfrac{0}{1})=(1+k_{\sigma(2)}+k_{\sigma(3)},1),\quad
S(\tfrac{1}{0})=(1+k_{\sigma(1)}+k_{\sigma(2)},1)
\]
For every reduced fraction $t\in (0,\infty)$, define $S(t)$ as follows.
\begin{itemize}   
 \item[(1)] Orient $\overline{L_t}$ from the lower left to the upper right. Apply the triangle-crossing and edge-crossing rules to every triangle-passage and edge-crossing occurrence of $\overline L_t$, and list the resulting signs in their order of occurrence.
 \item[(2)] \index{Run-length sequence}Form the integer sequence $S(t)=(a_0,\dots,a_n)$ from the lengths of the consecutive runs of equal signs in (1).
\end{itemize}

\begin{lemm}\label{lem:endpoint-rule-Lt}
Let $t\in(0,\infty)$ and $S(t)=(a_0,\dots,a_n)$. Choose $+$ as both endpoint triangle signs of $L_t=\gamma^L_{(0,0),(q,p)}$. Its consecutive run lengths then satisfy
\[
S(L_t)=(a_1,\dots,a_n)
\]
In particular,
\[
|L_t|=N(a_1,\dots,a_n)
\]
holds.
\end{lemm}
\begin{proof}
The first occurrences along $\overline L_t$ are, in order, the horizontal edge $y=0$, a triangle, the diagonal edge $x+y=0$, a triangle, and the vertical edge $x=0$. Indeed, the last two crossings occur at times $\varepsilon/(p+q)<\varepsilon/q$, and for the chosen range of $\varepsilon$ no other subdivision line is crossed earlier. All five occurrences have negative signs, with total multiplicity
\[
k_{\sigma(1)}+1+k_{\sigma(2)}+1+k_{\sigma(3)}=K-1
\]
The next triangle is crossed from $x=0$ to $x+y=1$. Its cut-off vertex $(0,1)$ lies on the left, so its sign is $+$. Hence the initial maximal run is exactly $-^{K-1}$.

Remove this initial part and reconnect the initial point in the first remaining triangle to $A$, and the terminal point in the last triangle to $B$. Preserving the internal crossings and their sides relative to edge midpoints, and making the connections inside the endpoint triangles, gives a generalized arc with the same passages and signs as the pure left push-off of $AB$. The original sign of the last triangle is also $+$, since the common vertex of its two crossed edges lies on the left. Choosing both endpoint signs to be $+$ therefore gives the word of $L_t$ without changing any other sign. Crossings of weight zero remain empty words, so the same argument applies to them.
Thus $S(L_t)=(a_1,\ldots,a_n)$, and Proposition~\ref{prop:length} gives $|L_t|=N(a_1,\ldots,a_n)$.
\end{proof}

\begin{exam}\label{ex:t=2/5}
Let $(k_1,k_2,k_3,\sigma)=(1,2,0,\textrm{id})$ and $t=\frac25$. Figure~\ref{fig:ex-extendedpresnakegraph-signed} shows the signs along $\overline L_t$, giving
\[
S(\tfrac{2}{5})=(5,1,3,3,1,5,4,1,3,4)
\]
as the resulting sequence.

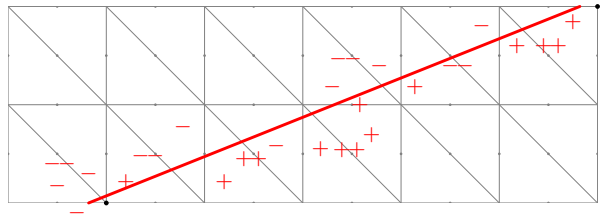
\begin{figure}[ht]
\centering
\begin{tikzpicture}[x=0.01cm,y=0.01cm,scale=0.5]
\definecolor{annblue}{RGB}{255,0,0}

\tikzset{
  gridline/.style={gray,line width=0.15pt,line cap=butt,line join=miter},
  gammaL/.style={draw=annblue,line width=1.1pt,line cap=round,line join=round},
   midpoint/.style={circle,fill=gray,inner sep=0pt,minimum size=1pt},
  enddot/.style={circle,fill=black,inner sep=0pt,minimum size=1.8pt},
  trisign/.style={text=annblue,font=\scriptsize\bfseries,inner sep=0.25pt},
  edgesign/.style={text=annblue,font=\scriptsize\bfseries,inner sep=0.15pt}
}

\begin{scope}[shift={(300,170)}, x={(260,0)}, y={(0,260)}]

\begin{scope}
  \clip (-1,0) rectangle (5,2);

  \foreach \i in {-1,...,4}{
    \foreach \j in {0,1}{
      \draw[gridline] (\i,{\j+1}) -- ({\i+1},\j);
    }
  }

  \foreach \i in {-1,...,5}{
    \draw[gridline] (\i,0) -- (\i,2);
  }
  \foreach \j in {0,1,2}{
    \draw[gridline] (-1,\j) -- (5,\j);
  }

  \foreach \i in {-1,...,4}{
    \foreach \j in {0,1,2}{
      \node[midpoint] at ({\i+0.5},\j) {};
    }
  }

  \foreach \i in {-1,...,5}{
    \foreach \j in {0,1}{
      \node[midpoint] at (\i,{\j+0.5}) {};
    }
  }

  \foreach \i in {-1,...,4}{
    \foreach \j in {0,1}{
      \node[midpoint] at ({\i+0.5},{\j+0.5}) {};
    }
  }
\end{scope}

\coordinate (A) at (0,0);
\coordinate (B) at (5,2);

\def\xleftshift{0.18}
\coordinate (ALshift) at ({-\xleftshift},0);
\coordinate (BLshift) at ({5-\xleftshift},2);
\draw[gammaL] (ALshift) -- (BLshift);

\node[trisign] at (-0.5,0.18) {$-$};
\node[trisign] at (-0.18,0.3) {$-$};
\node[trisign] at (0.20,0.22) {$+$};
\node[trisign] at (0.78,0.78) {$-$};
\node[trisign] at (1.20,0.22) {$+$};
\node[trisign] at (1.73,0.58) {$-$};
\node[trisign] at (2.18,0.55) {$+$};
\node[trisign] at (2.70,0.70) {$+$};
\node[trisign] at (2.30,1.18) {$-$};
\node[trisign] at (2.78,1.4) {$-$};
\node[trisign] at (3.14,1.18) {$+$};
\node[trisign] at (3.82,1.80) {$-$};
\node[trisign] at (4.18,1.6) {$+$};
\node[trisign] at (4.75,1.85) {$+$};

\node[edgesign] at (-0.3,-0.1) {$-$};

\node[edgesign] at (-0.55,0.40) {$-$};
\node[edgesign] at (-0.4,0.40) {$-$};

\node[edgesign] at (0.35,0.48) {$-$};
\node[edgesign] at (0.50,0.48) {$-$};

\node[edgesign] at (1.40,0.45) {$+$};
\node[edgesign] at (1.55,0.45) {$+$};

\node[edgesign] at (2.40,0.54) {$+$};
\node[edgesign] at (2.55,0.54) {$+$};

\node[edgesign] at (2.58,1) {$+$};

\node[edgesign] at (2.36,1.47) {$-$};
\node[edgesign] at (2.51,1.47) {$-$};

\node[edgesign] at (3.5,1.4) {$-$};
\node[edgesign] at (3.65,1.4) {$-$};

\node[edgesign] at (4.45,1.6) {$+$};
\node[edgesign] at (4.6,1.6) {$+$};

\node[enddot] at (A) {};
\node[enddot] at (B) {};

\end{scope}
\end{tikzpicture}
\caption{Signs along $\overline{L_t}$}
\label{fig:ex-extendedpresnakegraph-signed}
\end{figure}

The right endpoint of $\overline{L_t}$ lies in the interior of the upper rightmost edge of $\widetilde{\mathbb R^2}$, but this terminal contact is not an edge-crossing occurrence and contributes no sign. For $t=2/5$, $\overline{L_t}$ is chosen as a sufficiently small admissible perturbation of
\[
L_{t,\varepsilon}(u)=u(5,2)-(\varepsilon,0)
\]
The intersection of $L_{t,\varepsilon}$ with the unique crossed horizontal edge $y=1$ is $(5/2-\varepsilon,1)$, to the left of its midpoint $(5/2,1)$. Thus rule (3) assigns $+$ to this edge-crossing occurrence.
\end{exam}

\begin{lemm}\label{lem:parity-position-label}\index{Position label}
For a positive reduced fraction $t=p/q$, the position label $i_t$ is determined by the following congruence classes:
\[
\begin{array}{c|ccc}
(p,q)\pmod2&(0,1)&(1,1)&(1,0)\\ \hline
i_t&\sigma(1)&\sigma(2)&\sigma(3).
\end{array}
\]
\end{lemm}
\begin{proof}
The table holds for the initial Farey triple $(0/1,1/1,1/0)$. The vectors of numerators and denominators of two neighboring reduced fractions are distinct nonzero vectors in $\mathbb F_2^2$, so their sum is the remaining nonzero vector. Likewise, the new position label at each child in the GM tree is the label missing from the two inherited entries. Induction on the Farey tree therefore preserves the table.
\end{proof}

We now define the geometric objects used below. For a positive reduced fraction $t=p/q$, put $A=(0,0)$ and $B=(q,p)$. List the triangles traversed by the triangle-passage occurrences of $L_t=\gamma^L_{AB}$, in order from $A$ to $B$, as $\Delta_0,\ldots,\Delta_N$, and define
\[
\mathcal{PG}(t):=\bigcup_{j=0}^{N}\overline{\Delta_j}
\]
Also put
\[
M_t:=\frac{A+B}{2}=\left(\frac q2,\frac p2\right)
\]
and let $E_t$ denote the triangulation edge with midpoint $M_t$.

For a finite sign word $w=(\varepsilon_1,\ldots,\varepsilon_m)$, write
\[
w^\dagger:=(\bar\varepsilon_m,\ldots,\bar\varepsilon_1)
\]
for the word obtained by reversing the order and changing every sign, where $\bar{+}=-$ and $\bar{-}=+$.

\begin{lemm}\label{lem:half-turn-strongly-admissible}\index{Sign sequence}\index{Edge type}\index{Generalized strongly admissible sequence}
For a positive reduced fraction $t$, let $v_t$ be the finite sign word of $L_t$ with both endpoint signs chosen to be $+$. Then, for some finite sign word $P$ (possibly empty),
\[
v_t=+\,P\,+^{k_t}\,P^\dagger\,+
\]
Here $+^{k_t}$ is the empty word when $k_t=0$.
\end{lemm}
\begin{proof}
Write $t=p/q$, $A=(0,0)$, and $B=(q,p)$. The parameter $u\in[0,1]$ specifies a point on the segment: $uB=(uq,up)$ equals $A$ at $u=0$ and $B$ at $u=1$. Choose $0<\varepsilon<\varepsilon_t$, traverse the left admissible perturbation $uB-(\varepsilon,0)$ and the right translate $uB+(\varepsilon,0)$ from their initial to their terminal points, and list the signs assigned to their triangle passages and edge crossings. On the left, count the initial horizontal edge crossing but not the terminal one; on the right, count the terminal horizontal edge crossing but not the initial one. Denote the resulting finite sign words by $w^-$ and $w^+$, respectively.

The half-turn
\[
\rho_{A,B}(X)=A+B-X
\]
preserves the marked triangulation and each edge type, and sends the left translate to the right translate. Reparametrizing by $u\mapsto1-u$ reverses both the occurrence order and every sign. The initial horizontal edge counted on the left is sent to the terminal horizontal edge counted on the right. Thus $w^+=(w^-)^\dagger$, with no cyclic shift at either end of the finite words.

We next locate the possible differences between the two perturbations. All marked points belong to $\frac12\mathbb Z^2$. If $A+\lambda(B-A)\in\frac12\mathbb Z^2$, then $2\lambda p,2\lambda q\in\mathbb Z$. Since $\gcd(p,q)=1$, Lemma~\ref{lem:ax+by} gives integers $a,b$ with $ap+bq=1$. Hence $2\lambda=a(2\lambda p)+b(2\lambda q)\in\mathbb Z$. Thus the only possible marked point in the open segment $AB$ occurs at $\lambda=1/2$, namely $M_t$. Choose sufficiently small disjoint neighborhoods of the endpoints and of $M_t$. Outside them, the two curves can be continuously deformed into one another without meeting $\mathcal V$ or changing the crossed edges or their order. Their triangle-passage and edge-crossing occurrences, and the signs assigned to those occurrences, therefore agree outside these neighborhoods.

Remove the endpoint signs from the word of $L_t$ and split the remaining internal word at the central edge crossing. Let $P$ be the word before the center. The sequence of triangle and edge passages along $AB$ is reversed by the half-turn; away from the center, left and right push-offs receive the same signs. Reversing the orientation after the half-turn changes every sign, so the word after the center is $P^\dagger$.

The coordinates of $M_t$ show that $E_t$ is horizontal, diagonal, or vertical according as $(p,q)$ is congruent to $(0,1)$, $(1,1)$, or $(1,0)$ modulo $2$. By Lemma~\ref{lem:parity-position-label}, the edge type of $E_t$ is therefore $j=\sigma^{-1}(i_t)$, where $i_t$ is the position label of $t$. Rule (3) assigns $k_{\sigma(j)}=k_{i_t}=k_t$ signs to this crossing. Since $M_t$ lies to the right of the pure left push-off, these signs are all $+$, giving the central word $+^{k_t}$. Together with Lemma~\ref{lem:endpoint-rule-Lt}, this proves the asserted word representation.
\end{proof}

\begin{lemm}\label{rem:difference-(0,1)(1,infty)}
Let $0<t<\infty$ and $S(t)=(a_0,\dots,a_n)$. In terms of the entries, the symmetry of Lemma~\ref{lem:half-turn-strongly-admissible} gives the following properties.
\begin{itemize}
\item[(0)] $n$ is odd.
\item[(1)] $a_0=2+k_1+k_2+k_3$.
\item[(2)] If $t=\frac{1}{1}$, then $
S(t)=(2+k_1+k_2+k_3,2+k_{\sigma(2)})$.
\item[(3)] If $t\in(0,1)$, then $a_1=1$ and $a_n\neq1$. If $t\in(1,\infty)$, then $a_1\neq1$ and $a_n=1$.
\item[(4)] If $t=\frac12$, then $n=3$ and $a_3=a_2+1+k_t$. If $t=\frac21$, then $n=3$ and $a_1=a_2+1+k_t$.
\item[(5)] If $t\in(0,\frac12)\cup(\frac12,1)$, then $a_2+1=a_n$, and $a_{2+i}=a_{n-i}$ for $i=1,2,\dots,\frac{n-5}{2}$
(there are no such relations when $n=5$). Moreover,
\[
a_{\frac{n+3}{2}}=a_{\frac {n+1}{2}}+(-1)^{\frac{n+1}{2}}k_t
\]
holds.
\item[(6)] If $t\in(1,\frac21)\cup(\frac21,\infty)$, then
$a_1=a_{n-1}+1$, and $a_{1+i}=a_{n-i-1}$ for $i=1,2,\dots,\frac{n-5}{2}$ (there are no such relations when $n=5$). Moreover,
\[
a_{\frac{n+1}{2}}=a_{\frac {n-1}2}+(-1)^{\frac{n-1}{2}}k_t
\]
holds.
\item[(7)] Put $\sigma^\ast:=\sigma\circ(1\ 3)$, and denote the $(k_1,k_2,k_3,\sigma^\ast)$-generalized strongly admissible sequence associated with $1/t$ by $S^\ast(1/t)$. Then
\[
S^\ast(1/t)=(a_0,a_n,a_{n-1},\dots,a_1)
\]
and the corresponding two-sided periodic sequences are reversals of one another up to an index shift. This relation also holds for $t=0,\infty$.
\end{itemize}
\end{lemm}
\begin{proof}
The initial negative block counted in the proof of Lemma~\ref{lem:endpoint-rule-Lt} gives (1). Counting signs directly for $t=1$ gives (2).

Write $t=p/q\in(0,1)$ in lowest terms, and consider the word $P$ used in the proof of Lemma~\ref{lem:half-turn-strongly-admissible}. The first internal crossing is with $x+y=1$, whose midpoint lies to the left of the curve. The next triangle also cuts off the vertex $(1,0)$ on the right. Thus $P$ begins with at least $k_{\sigma(2)}+1$ negative signs. If $q=2$, then $p=1$ and $P=-^{k_{\sigma(2)}+1}$. If $q>2$, the triangle immediately after the first vertical edge $x=1$ lies before the center and cuts off $(1,1)$ on the left, so $P$ also contains positive signs.

Suppose $P$ has $d\geq1$ maximal runs of equal signs. The last sign of $P$ is opposite to the first sign of $P^\dagger$, so the central block $+^{k_t}$ merges with one of them. On adding the endpoint signs $+$, the initial $+$ forms a separate run, whereas the terminal $+$ merges with the last positive run of $P^\dagger$. Thus $v_t$ has $2d+1$ runs, and $n=2d+1$ in $S(t)$. This proves (0) and (3) for $0<t<1$.
By the preceding observation, $d=1$ occurs only for $t=1/2$. In that case $a_3=a_2+1+k_t$, proving the first half of (4). If $d\geq2$, the endpoint correspondence gives $a_n=a_2+1$, and the other noncentral pairs give $a_{2+i}=a_{n-i}$. The last sign of $P$ is negative for odd $d$ and positive for even $d$, so the difference between the two central run lengths is $(-1)^{d+1}k_t$. This is the central relation in (5).

Reflection $(x,y)\mapsto(y,x)$ changes the weight arrangement to $\sigma^*$ and sends a left push-off to a right push-off. Reflection changes every sign in the internal word. On the other hand, the relation between the left and right internal words established in the proof of Lemma~\ref{lem:half-turn-strongly-admissible} reverses their order and changes their signs. Combining these two operations shows that the internal word of the left push-off in the dual tree is the reversal of the original internal word. Since both endpoint signs are chosen to be $+$, the entire word $v_t$ is reversed as well. The added initial block has length $K-1$ in both cases, giving the reciprocal formula in (7). This formula and the case $0<t<1$ give (0) and (3) for $t>1$, the second half of (4), and (6). For $t=1$, assertion (0) follows from (2).

The sequence in (7) is a cyclic shift of the reversal of $(a_0,\dots,a_n)$, which proves the assertion about two-sided periodic sequences. At the endpoints, the definitions give $S^*(\infty)=S(0)$ and $S^*(0)=S(\infty)$. Both sequences have two entries, so the same relation holds.
\end{proof}

Let $0<t<1$ and $S(t)=(a_0,\ldots,a_n)$. By Lemma~\ref{lem:endpoint-rule-Lt}, list the signs assigned to the triangle-passage and edge-crossing occurrences of $L_t$ in their order from $A$ to $B$, and record their run lengths. Then
\[
S_-(t):=(a_1,\ldots,a_n)
\]
This sign word begins with a single $+$ followed by $a_2$ negative signs. Change only this initial $+$ to $-$, and denote the resulting run-length sequence by $S'_-(t)$. Since $a_1=1$ by Lemma~\ref{rem:difference-(0,1)(1,infty)}~(3),
\[
S'_-(t)=(a_2+1,a_3,\ldots,a_n)
\]
holds.

\begin{coro}\label{cor:local-reversal-rules}
The following two sign-reversal relations hold.
\begin{itemize}
\item[(1)] Let $p\geq2$ and $S'_{-}(1/p)=(b_1,\dots,b_m)$. In $\mathcal{PG}(1/p)$, reverse the signs assigned to the crossing occurrence of the edge $E_{1/p}$ with midpoint
\[
M_{1/p}=\left(\frac p2,\frac12\right)
\]
and also reverse the sign of the last triangle-passage occurrence before $L_{1/p}$ reaches its endpoint. The resulting run-length sequence is
\[
(b_m,b_{m-1},\dots,b_1-1,1)
\]
as displayed.
\item[(2)] Let $p\geq2$ and $S'_{-}(p/(p+1))=(a_1,\dots,a_\ell)$. In $\mathcal{PG}(p/(p+1))$, reverse the signs assigned to the crossing occurrence of the edge $E_{p/(p+1)}$ with midpoint
\[
M_{p/(p+1)}=\left(\frac{p+1}{2},\frac p2\right)
\]
The resulting run-length sequence is
\[
(a_\ell,a_{\ell-1},\dots,a_1)
\]
as displayed.
\end{itemize}
In either case, if the corresponding $k_t$ is zero, rule (3) assigns no sign to the crossing occurrence of $E_t$, so there is no sign to reverse at that occurrence.
\end{coro}
\begin{proof}
By the word representation in Lemma~\ref{lem:half-turn-strongly-admissible}, the word corresponding to $S'_-(t)$ for $0<t<1$ is
\[
v'_t=-\,P\,+^{k_t}\,P^\dagger\,+
\]
Reversing the signs at the central edge gives
\[
-\,P\,-^{k_t}\,P^\dagger\,+=(v'_t)^\dagger
\]
so the run-length sequence is reversed. This proves (2).
For (1), additionally change the last triangle sign from $+$ to $-$. Since $P$ begins with a negative sign, the initial negative run of $v'_t$ has length $b_1\geq2$. Thus removing one sign from the final positive run of $(v'_t)^\dagger$ leaves a nonempty run, followed by one new negative sign. The resulting run lengths are $(b_m,\ldots,b_1-1,1)$. When $k_t=0$, the center is an empty word, and the same identity holds.
\end{proof}

The main theorem of this section is as follows.

\begin{theorem}\label{continued-fraction-theorem2}\index{Generalized Cohn matrix}\index{Generalized strongly admissible sequence}
For every reduced fraction $t\in(0,\infty]$,
$C_t=F_{S(t)}$.
\end{theorem}

Before proving this theorem, we show that it immediately implies Theorem~\ref{thm:length-is-number}, whose proof was postponed at the end of Section~\ref{sec:GM-distance}.
\par\medskip
\noindent\textbf{Theorem~\ref{thm:length-is-number}} (restated)\textbf{.}\,
Fix $(k_1,k_2,k_3)\in\mathbb Z_{\geq0}^3$ and $\sigma\in\mathfrak S_3$. Let $p,q$ be relatively prime nonnegative integers, not both zero, and put $t=p/q$, with the convention $1/0=\infty$. For $A=(0,0)$ and $B=(q,p)$,
\[|\gamma_{AB}^R|=|\gamma_{AB}^L|=m_t\]
In particular, $d(A,B)=m_t$.
\par\medskip
\begin{proof}
If $t=0$ or $t=\infty$, then $A,B$ are neighboring lattice points on a single triangulation edge. Each of $\gamma_{AB}^R$ and $\gamma_{AB}^L$ consists of just one triangle passage, so its fence poset is empty and its GM length is $1$. The initial vertex of the GM tree gives $m_0=m_\infty=1$. Since GM lengths are positive integers, $d(A,B)=1$ also follows.

Now suppose $t\in(0,\infty)$ and write $S(t)=(a_0,\dots,a_n)$. Corollary~\ref{cor:matrix-combinatorial}, Theorem~\ref{continued-fraction-theorem2}, and Proposition~\ref{prop:C21=m} give
\[
N(a_1,\dots,a_n)=(F_{S(t)})_{21}=(C_t)_{21}=m_t
\]
By Lemma~\ref{lem:endpoint-rule-Lt}, this continuant is the GM length of $L_t=\gamma_{AB}^L$, so
\[
|\gamma_{AB}^L|=N(a_1,\dots,a_n)=m_t.
\]
Lemma~\ref{lem:same-length-lemma0} gives $|\gamma_{AB}^R|=|\gamma_{AB}^L|$. Finally, Theorem~\ref{thm:gm-distance} gives
\[
d(A,B)=|\gamma_{AB}^R|=|\gamma_{AB}^L|=m_t
\]
as required.
\end{proof}

To prove the theorem, we first consider $t\in(0,1)$. The sequences $S'_-(t)$ defined above satisfy the following proposition.

\begin{prop}\label{prop:presnake-relation}
Let \((r,t,s)\in\mathrm F\mathbb T\) with \(0<t<1\). Then \(0\leq r<t<s\leq1\), and the following hold.
\begin{itemize}
    \item[(1)] Suppose \(r=\frac01\) and \(s\neq\frac11\). If \(S'_{-}(s)=(b_1,\dots,b_m)\), then
    \[
    S'_{-}(t)=(b_m,b_{m-1},\dots,b_1-1,1,2+k_{\sigma(2)}+k_{\sigma(3)})
\]
    holds.
    \item[(2)] Suppose \(r\neq\frac01\) and \(s=\frac11\). If \(S'_{-}(r)=(a_1,\dots,a_\ell)\), then
    \[
    S'_{-}(t)=(2+k_{\sigma(2)},2+k_1+k_2+k_3,a_{\ell},\dots,a_1)
    \]
    holds.
    \item[(3)] Suppose \(r\neq\frac01\) and \(s\neq\frac11\). If \(S'_{-}(r)=(a_1,\dots,a_\ell)\) and \(S'_{-}(s)=(b_1,\dots,b_m)\), then
    \[
    S'_{-}(t)=(b_{m},\dots,b_1-1,1,2+k_1+k_2+k_3,a_{\ell},\dots,a_1)
    \]
    holds.
\end{itemize}
\end{prop}

\begin{proof}[Proof of Proposition~\ref{prop:presnake-relation} (1) and (2)]

We first prove (1). Under the assumptions \(r=\frac01\) and \(s\neq\frac11\), there is an integer \(p\in\ZZ_{>1}\) such that $
s=\frac{1}{p}$
Thus it suffices to prove the assertion for $
t=\frac{1}{p+1}$
Since \(p\geq2\), the final \(2+k_{\sigma(2)}+k_{\sigma(3)}\) signs in \(\mathcal{PG}(\frac{1}{p+1})\) are all \(+\), and the preceding sign is \(-\); see Figure~\ref{fig:presnake-1-8}. Figures~\ref{fig:presnake-1-8} and~\ref{fig:presnake-1-7} show $k_1=k_2=k_3=1$ and $p=7$.

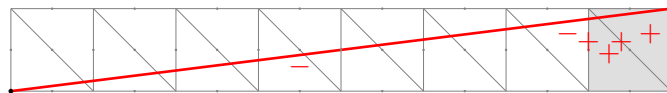
\begin{figure}[ht]
    \centering
\begin{tikzpicture}[x=0.01cm,y=0.01cm,scale=0.42]
  \definecolor{annred}{RGB}{255,0,0}

  \tikzset{
    gridline/.style={gray,line width=0.15pt,line cap=butt,line join=miter},
    gammaL/.style={draw=annred,line width=1pt,line cap=round,line join=round},
    enddot/.style={circle,fill=black,inner sep=0pt,minimum size=1.8pt},
    midpt/.style={circle,draw=gray,fill=gray,inner sep=0pt,minimum size=0.5pt},
    shadecell/.style={fill=gray!25,draw=none}
  }

  \begin{scope}[shift={(180,170)}, x={(260,0)}, y={(0,260)}]

    \path[shadecell] (7,0) -- (7,1) -- (8,0) -- cycle;
    \path[shadecell] (7,1) -- (8,1) -- (8,0) -- cycle;

    \foreach \i in {0,...,7}{
      \draw[gridline] (\i,1) -- ({\i+1},0);
    }

    \foreach \i in {0,...,8}{
      \draw[gridline] (\i,0) -- (\i,1);
    }
    \foreach \j in {0,1}{
      \draw[gridline] (0,\j) -- (8,\j);
    }

    \coordinate (A) at (0,0);
    \coordinate (B) at (8,1);

    \draw[gammaL] (A) -- (B);

    \foreach \i in {0,...,7}{
      \foreach \j in {0,1}{
        \node[midpt] at ({\i+0.5},\j) {};
      }
    }

    \foreach \i in {0,...,8}{
      \node[midpt] at (\i,0.5) {};
    }

    \foreach \i in {0,...,7}{
      \node[midpt] at ({\i+0.5},0.5) {};
    }

    \node[enddot] at (A) {};
    \node[enddot] at (B) {};
    \node[] at (3.5,0.3) {\textcolor{red}{$-$}};
    \node[] at (6.75,0.7) {\textcolor{red}{$-$}};
    \node[] at (7,0.6) {\textcolor{red}{$+$}};
    \node[] at (7.25,0.45) {\textcolor{red}{$+$}};
    \node[] at (7.4,0.6) {\textcolor{red}{$+$}};
    \node[] at (7.75,0.7) {\textcolor{red}{$+$}};
  \end{scope}
\end{tikzpicture}
    \caption{\(\mathcal{PG}(\frac{1}{p+1})\)}
    \label{fig:presnake-1-8}
\end{figure}

Remove the rightmost unit square, consisting of two right triangles, from \(\mathcal{PG}(\frac{1}{p+1})\). Denote the union of the remaining closed triangles by \(\mathcal{SPG}(\frac{1}{p+1})\). Compare the sign word obtained as $L_{1/(p+1)}$ traverses this region with the sign word of $L_{1/p}$ in \(\mathcal{PG}(\frac1p)\); compare the white part of Figure~\ref{fig:presnake-1-8} with Figure~\ref{fig:presnake-1-7}.

\begin{figure}[ht]
    \centering
    \begin{tikzpicture}[x=0.01cm,y=0.01cm,scale=0.42]
  \definecolor{annred}{RGB}{255,0,0}

  \tikzset{
    gridline/.style={gray,line width=0.15pt,line cap=butt,line join=miter},
    gammaL/.style={draw=annred,line width=1pt,line cap=round,line join=round},
    enddot/.style={circle,fill=black,inner sep=0pt,minimum size=1.8pt},
    midpt/.style={circle,draw=gray,fill=gray,inner sep=0pt,minimum size=0.5pt}
  }

  \begin{scope}[shift={(180,170)}, x={(260,0)}, y={(0,260)}]

    \foreach \i in {0,...,6}{
      \draw[gridline] (\i,1) -- ({\i+1},0);
    }

    \foreach \i in {0,...,7}{
      \draw[gridline] (\i,0) -- (\i,1);
    }
    \foreach \j in {0,1}{
      \draw[gridline] (0,\j) -- (7,\j);
    }

    \coordinate (A) at (0,0);
    \coordinate (B) at (7,1);

    \def\xleftshift{0}
    \draw[gammaL] ({-\xleftshift},0) -- ({7-\xleftshift},1);

    \foreach \i in {0,...,6}{
      \foreach \j in {0,1}{
        \node[midpt] at ({\i+0.5},\j) {};
      }
    }

    \foreach \i in {0,...,7}{
      \node[midpt] at (\i,0.5) {};
    }

    \foreach \i in {0,...,6}{
      \node[midpt] at ({\i+0.5},0.5) {};
    }

    \node[enddot] at (A) {};
    \node[enddot] at (B) {};
     \node[] at (3.5,0.3) {\textcolor{red}{$+$}};
    \node[] at (6.75,0.7) {\textcolor{red}{$+$}};
  \end{scope}
\end{tikzpicture}
    \caption{\(\mathcal{PG}(\frac{1}{p})\)}
    \label{fig:presnake-1-7}
\end{figure}
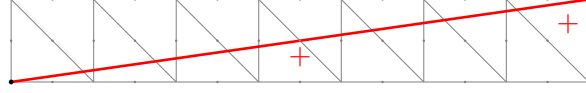
To make the comparison precise, put $v=(p,1)$ and $w=(p+1,1)$. The following intersection coordinates are calculated for the segments $\ell_{1/p}$ and $\ell_{1/(p+1)}$ before pushing them off. The latter meets $x=p$ at $P=(p,p/(p+1))$. The triangle between the segments $[0,v]$ and $[0,P]$ is
\[
T=\operatorname{conv}\{0,v,P\}
\]
and its only lattice points are $0,v$. Examining the heights $0,1/2,1$ also gives
\[
T\cap\tfrac12\mathbb Z^2=\{0,v/2,w/2,v\}
\]
Consequently, deforming one segment into the other while keeping the endpoint connections inside their first and last triangles preserves the triangle-passage sequence. The only edge midpoints whose avoidance sides can change are $v/2$ and $w/2$.

The point $w/2$ is the midpoint of the new segment, and the pure left push-offs of both the old and new segments pass to its left. Thus no sign changes there. In contrast, at the \emph{same geometric edge} $E_{1/p}$ with midpoint $v/2$, the old pure left push-off passes to the left of the midpoint, whereas the new segment passes to its right. Thus the $k_{1/p}$ assigned signs are reversed. If the weight is zero, both words are empty.
In the last triangle, the old word has terminal sign $+$, whereas the new subpath runs from the diagonal edge to $x=p$, cutting off $(p,0)$ on the right, and hence has sign $-$. All other triangle and edge signs agree. Here both words have initial sign $-$, as required by the convention for $S'_-$.

By Corollary~\ref{cor:local-reversal-rules}~(1), the run lengths of the sign word in \(\mathcal{SPG}(\frac{1}{p+1})\) are
\[
(b_m,b_{m-1},\dots,b_1-1,1)
\]
Count the vertical edge crossing at the cut point $P$ in the latter part. Inside the rightmost unit square, this crossing, the two triangle passages, and the intervening diagonal edge crossing all have positive signs. Thus the latter part is a run of length $2+k_{\sigma(2)}+k_{\sigma(3)}$, which does not merge with the final negative sign of the former part. This proves (1).

Next we prove (2). First, if $r=\frac12$, that is, $p=1$, counting the signs from the initial point according to the rules gives
\[
S(\tfrac23)=(K-1,1,1+k_{\sigma(2)},K-1,
2+k_{\sigma(2)}+k_{\sigma(3)},2+k_{\sigma(2)})
\]
Together with $S'_{-}(1/2)=(2+k_{\sigma(2)},2+k_{\sigma(2)}+k_{\sigma(3)})$, this verifies the formula in (2).

Now let $p\geq2$. Under the assumptions \(r\neq\frac01\) and \(s=\frac11\), we write, for an integer \(p\in\ZZ_{>1}\), $
r=\frac{p}{p+1}$
Thus it suffices to prove the assertion for $
t=\frac{p+1}{p+2}$
Since \(p\geq2\), the first \(4+k_{\sigma(1)}+2k_{\sigma(2)}+k_{\sigma(3)}\) signs in \(\mathcal{PG}(\frac{p+1}{p+2})\) consist of \(2+k_{\sigma(2)}\) negative signs followed by \(2+k_1+k_2+k_3\) positive signs; see Figure~\ref{fig:presnake-4-5}. Figures~\ref{fig:presnake-4-5} and~\ref{fig:presnake-3-4} show $k_1=k_2=k_3=1$ and $p=3$.

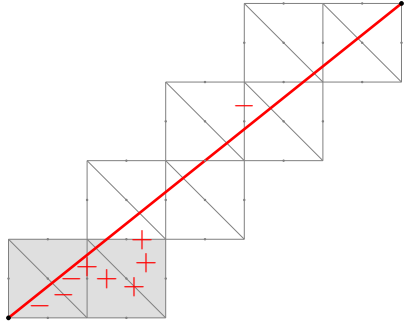
\begin{figure}[ht]
    \centering
   \begin{tikzpicture}[x=0.01cm,y=0.01cm,scale=0.4]
  \definecolor{annred}{RGB}{255,0,0}

  \tikzset{
    gridline/.style={gray,line width=0.15pt,line cap=butt,line join=miter},
    gammaL/.style={draw=annred,line width=1pt,line cap=round,line join=round},
    enddot/.style={circle,fill=black,inner sep=0pt,minimum size=1.8pt},
    midpt/.style={circle,draw=gray,fill=gray,inner sep=0pt,minimum size=0.5pt},
    shadecell/.style={fill=gray!25,draw=none}
  }

  \begin{scope}[shift={(180,170)}, x={(260,0)}, y={(0,260)}]

    \path[shadecell] (0,0) -- (0,1) -- (1,1) -- (1,0) --cycle;
    \path[shadecell] (1,0) -- (2,0) -- (2,1) -- (1,1) -- cycle;

    \foreach \xa/\ya/\xb/\yb in {
      0/0/1/0, 1/0/2/0,
      0/1/1/1, 1/1/2/1, 2/1/3/1,
      1/2/2/2, 2/2/3/2, 3/2/4/2,
      2/3/3/3, 3/3/4/3, 4/3/5/3,
      3/4/4/4, 4/4/5/4,
      0/0/0/1,
      1/0/1/1, 1/1/1/2,
      2/0/2/1, 2/1/2/2, 2/2/2/3,
      3/1/3/2, 3/2/3/3, 3/3/3/4,
      4/2/4/3, 4/3/4/4,
      5/3/5/4
    }{
      \draw[gridline] (\xa,\ya) -- (\xb,\yb);
    }

    \foreach \i/\j in {
      0/0, 1/0, 1/1, 2/1,
      2/2, 3/2, 3/3, 4/3
    }{
      \draw[gridline] (\i,{\j+1}) -- ({\i+1},\j);
    }

    \coordinate (A) at (0,0);
    \coordinate (B) at (5,4);

    \draw[gammaL] (A) -- (B);

    \foreach \xa/\ya/\xb/\yb in {
      0/0/1/0, 1/0/2/0,
      0/1/1/1, 1/1/2/1, 2/1/3/1,
      1/2/2/2, 2/2/3/2, 3/2/4/2,
      2/3/3/3, 3/3/4/3, 4/3/5/3,
      3/4/4/4, 4/4/5/4,
      0/0/0/1,
      1/0/1/1, 1/1/1/2,
      2/0/2/1, 2/1/2/2, 2/2/2/3,
      3/1/3/2, 3/2/3/3, 3/3/3/4,
      4/2/4/3, 4/3/4/4,
      5/3/5/4
    }{
      \node[midpt] at ({(\xa+\xb)/2},{(\ya+\yb)/2}) {};
    }

    \foreach \i/\j in {
      0/0, 1/0, 1/1, 2/1,
      2/2, 3/2, 3/3, 4/3
    }{
      \node[midpt] at ({\i+0.5},{\j+0.5}) {};
    }

    \node[enddot] at (A) {};
    \node[enddot] at (B) {};
    \node[] at (0.4,0.15) {\textcolor{red}{$-$}};
    \node[] at (0.7,0.3) {\textcolor{red}{$-$}};
    \node[] at (0.8,0.5) {\textcolor{red}{$-$}};
    \node[] at (1,0.65) {\textcolor{red}{$+$}};
    \node[] at (1.25,0.5) {\textcolor{red}{$+$}};
    \node[] at (1.6,0.4) {\textcolor{red}{$+$}};
    \node[] at (1.75,0.7) {\textcolor{red}{$+$}};
    \node[] at (1.7,1) {\textcolor{red}{$+$}};
    \node[] at (3,2.7) {\textcolor{red}{$-$}};
  \end{scope}
\end{tikzpicture}
    \caption{\(\mathcal{PG}(\frac{p+1}{p+2})\)}
    \label{fig:presnake-4-5}
\end{figure}

Remove the two leftmost unit squares, consisting of four right triangles in total, from \(\mathcal{PG}(\frac{p+1}{p+2})\). Denote the union of the remaining closed triangles by \(\mathcal{SPG}(\frac{p+1}{p+2})\). Compare the sign word obtained as $L_{(p+1)/(p+2)}$ traverses this region with the sign word of $L_{p/(p+1)}$ in \(\mathcal{PG}(\frac{p}{p+1})\); see Figures~\ref{fig:presnake-4-5} and~\ref{fig:presnake-3-4}.

\begin{figure}[ht]
    \centering
    \begin{tikzpicture}[x=0.01cm,y=0.01cm,scale=0.4]
  \definecolor{annred}{RGB}{255,0,0}

  \tikzset{
    gridline/.style={gray,line width=0.15pt,line cap=butt,line join=miter},
    gammaL/.style={draw=annred,line width=1pt,line cap=round,line join=round},
    enddot/.style={circle,fill=black,inner sep=0pt,minimum size=1.8pt},
    midpt/.style={circle,draw=gray,fill=gray,inner sep=0pt,minimum size=0.5pt}
  }

  \begin{scope}[shift={(180,170)}, x={(260,0)}, y={(0,260)}]

    \foreach \xa/\ya/\xb/\yb in {
      0/0/1/0, 1/0/2/0,
      0/1/1/1, 1/1/2/1, 2/1/3/1,
      1/2/2/2, 2/2/3/2, 3/2/4/2,
      2/3/3/3, 3/3/4/3,
      0/0/0/1,
      1/0/1/1, 1/1/1/2,
      2/0/2/1, 2/1/2/2, 2/2/2/3,
      3/1/3/2, 3/2/3/3,
      4/2/4/3
    }{
      \draw[gridline] (\xa,\ya) -- (\xb,\yb);
    }

    \foreach \i/\j in {
      0/0, 1/0, 1/1, 2/1, 2/2, 3/2
    }{
      \draw[gridline] (\i,{\j+1}) -- ({\i+1},\j);
    }

    \coordinate (A) at (0,0);
    \coordinate (B) at (4,3);

    \draw[gammaL] (A) -- (B);

    \foreach \xa/\ya/\xb/\yb in {
      0/0/1/0, 1/0/2/0,
      0/1/1/1, 1/1/2/1, 2/1/3/1,
      1/2/2/2, 2/2/3/2, 3/2/4/2,
      2/3/3/3, 3/3/4/3,
      0/0/0/1,
      1/0/1/1, 1/1/1/2,
      2/0/2/1, 2/1/2/2, 2/2/2/3,
      3/1/3/2, 3/2/3/3,
      4/2/4/3
    }{
      \node[midpt] at ({(\xa+\xb)/2},{(\ya+\yb)/2}) {};
    }

    \foreach \i/\j in {
      0/0, 1/0, 1/1, 2/1, 2/2, 3/2
    }{
      \node[midpt] at ({\i+0.5},{\j+0.5}) {};
    }

    \node[enddot] at (A) {};
    \node[enddot] at (B) {};
    \node[] at (2,1.7) {\textcolor{red}{$+$}};
  \end{scope}
\end{tikzpicture}
    \caption{\(\mathcal{PG}(\frac{p}{p+1})\)}
    \label{fig:presnake-3-4}
\end{figure}
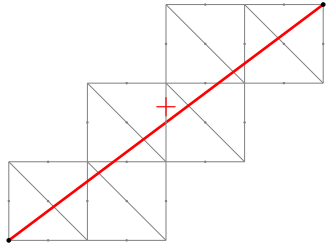

Put $u=(p+1,p)$, $v=(1,1)$, and $w=u+v$. The segment $[0,w]$ meets $y=1$ at $P=v+(1/(p+1),0)$. Comparing after translation by $v$, the region between the old segment $[0,u]$ and the new segment $[P-v,u]$ is
\[
T=\operatorname{conv}\{0,u,(1/(p+1),0)\}
\]
Since $\det(u,v)=1$, every half-lattice point has a unique expression $2z=\alpha u-\delta v$. The condition $z\in T$ is
\[
\delta\geq0,\qquad\alpha+\delta\leq2,\qquad p\alpha\geq\delta.
\]
Enumerating the integer pairs $\alpha,\delta$ gives
\[
T\cap\tfrac12\mathbb Z^2=\{0,u/2,u,(u-v)/2\}
\]
The only lattice points are $0,u$, so the triangle-passage sequences agree apart from the connections in the endpoint triangles.

After translation, the additional point is $v+(u-v)/2=w/2$. Both the old and new pure left push-offs pass to its left, so its signs do not change. Only the signs at the \emph{same edge} $v+E_{p/(p+1)}$, whose midpoint is $v+u/2$, change. The old pure left push-off passes to the left of this midpoint, whereas the new subpath passes to its right; hence the $k_{p/(p+1)}$ signs are reversed.
The first triangle of the new subpath is traversed from the horizontal edge to the diagonal edge, cutting off $v+(1,0)$ on the right. Its negative sign agrees with the initial sign chosen in $S'_-(p/(p+1))$. The last triangle has terminal sign $+$ in both words. Corollary~\ref{cor:local-reversal-rules} (2) therefore gives the run lengths $(a_\ell,\ldots,a_1)$ for the latter part.

Count the horizontal edge crossing at the cut point $P$ in the former part. In the four removed triangles, the initial sign, the first diagonal edge crossing, and the next triangle give $-^{2+k_{\sigma(2)}}$. The subsequent vertical edge crossing, two triangles, diagonal edge crossing, and final horizontal edge crossing give $+^{K-1}$. The latter part begins with a negative sign, so there is no merging, and the formula in (2) follows.
\end{proof}

\index{Upper Christoffel word}To prove Proposition~\ref{prop:presnake-relation} (3), we introduce upper Christoffel words. For a reduced fraction $a/b$ with $0\leq a\leq b$ and $b>0$, define a finite word in $X,Y$,
\[
\mathrm{uch}_{a/b}:=w_1\cdots w_b
\]
by
\[
w_i=
\begin{cases}
X & \text{if $\lceil ia/b\rceil-\lceil(i-1)a/b\rceil=1$},\\
Y & \text{if $\lceil ia/b\rceil-\lceil(i-1)a/b\rceil=0$}
\end{cases}
\qquad(1\leq i\leq b).
\]
Here $\lceil x\rceil$ is the least integer greater than or equal to $x$. This word is called the \emph{upper Christoffel word} associated with $a/b$.

\begin{exam}
The upper Christoffel word $\mathrm{uch}_{2/5}$ associated with $\frac25$ is $XYXYY$; see also Figure~\ref{fig:christoffel}.
\begin{figure}[ht]
    \centering
    \begin{tikzpicture}
    \draw[thick,gray] (0,0) grid (5,2);
    \draw[red,thick] (0,0) -- (5,2);
    \draw[line width=2pt] (0,0) -- (1,1)-- (2,1) -- (3,2) -- (5,2);
    \node at (0.5,-0.5) {$X$};
    \node at (1.5,-0.5) {$Y$};
    \node at (2.5,-0.5) {$X$};
    \node at (3.5,-0.5) {$Y$};
    \node at (4.5,-0.5) {$Y$};
    \end{tikzpicture}
    \caption{The upper Christoffel word $\mathrm{uch}_{2/5}$}
    \label{fig:christoffel}
\end{figure}
\end{exam}

\begin{theorem}\label{thm:christoffel}\index{Upper Christoffel word}
Let \((r,t,s)\in\mathrm F\TT\) with $0<t<1$. Then $0\leq r<t<s\leq1$, and
\[
\mathrm{uch}_t=\mathrm{uch}_s\cdot \mathrm{uch}_r
\]
where $\cdot$ denotes concatenation of the two finite words in the indicated order.
\end{theorem}
\begin{proof}
Write $r=\frac ab$ and $s=\frac cd$ in lowest terms.
Since $(r,t,s)\in\mathrm F\mathbb T$, we have $r<t<s$ and $
bc-ad=1$. Hence $
t=\frac{a+c}{b+d}$
Following the definition of upper Christoffel words, put
\[
\Delta_i(u):=\left\lceil iu\right\rceil-\left\lceil (i-1)u\right\rceil
\]
Then the $i$th letter of $\mathrm{uch}_u$ is $X$ if $\Delta_i(u)=1$, and $Y$ if $\Delta_i(u)=0$.

First we show $\Delta_i(t)=\Delta_i(s)$ for $1\leq i\leq d$.
\[
s-t=\frac{c}{d}-\frac{a+c}{b+d}
=\frac{bc-ad}{d(b+d)}
=\frac{1}{d(b+d)}
\]
For $1\leq i\leq d-1$, this gives
\[
0<is-it<\frac{1}{d}
\]
Since $\gcd(c,d)=1$, we have $is\notin\mathbb Z$,
and its fractional part is at least $1/d$. Therefore
\[
\lceil it\rceil=\lceil is\rceil \qquad (1\le i\le d-1)
\]
Moreover,
\[
dt=c-\frac{1}{b+d}
\]
gives $\lceil dt\rceil=c=\lceil ds\rceil$.
Thus $
\Delta_i(t)=\Delta_i(s)\qquad (1\le i\le d)$
follows.

Next we show $\Delta_{d+j}(t)=\Delta_j(r)$ for $1\leq j\leq b$.
\[
(d+j)t-c
=\frac{(d+j)(a+c)}{b+d}-c
=jr+\frac{j-b}{b(b+d)}
\]
Hence, for $1\leq j\leq b-1$,
\[
0<jr-\bigl((d+j)t-c\bigr)<\frac{1}{b}
\]
Since $\gcd(a,b)=1$, we have $jr\notin\mathbb Z$, and its fractional part is at least $1/b$. Thus
\[
\lceil (d+j)t\rceil-c=\lceil jr\rceil \qquad (1\le j\le b-1)
\]
For $j=b$,
\[
(d+b)t-c=a=br
\]
so $\lceil(d+b)t\rceil-c=a=\lceil br\rceil$. Hence
\[
\Delta_{d+j}(t)
=\lceil(d+j)t\rceil-\lceil(d+j-1)t\rceil
=\lceil jr\rceil-\lceil(j-1)r\rceil
=\Delta_j(r)
\]
follows.
The first equality $\Delta_i(t)=\Delta_i(s)$ shows that the first $d$ letters of $\mathrm{uch}_t$ form $\mathrm{uch}_s$. The second equality $\Delta_{d+j}(t)=\Delta_j(r)$ shows that the remaining $b$ letters form $\mathrm{uch}_r$. Therefore
\[
\mathrm{uch}_t=\mathrm{uch}_s\cdot\mathrm{uch}_r
\]
as required.
\end{proof}
The geometric decomposition corresponding to this concatenation can be described explicitly as follows.

\begin{coro}\label{lem:decomposition-lemma}
For \((r,t,s)\in\mathrm F\TT\) with $0<r<t<s<1$, the region \(\mathcal{PG}(t)\) decomposes, in order from the lower left to the upper right, into a translate of \(\mathcal{PG}(s)\), one unit square consisting of two right triangles, and a translate of \(\mathcal{PG}(r)\).
\end{coro}
\begin{proof}
Put $r=a/b$, $s=c/d$, $u=(b,a)$, and $v=(d,c)$. Then $\det(u,v)=bc-ad=1$.
The intersections of $[0,u+v]$ with $x=d$ and $y=c$ are
\[
P_s=v-\frac{(0,1)}{b+d},\qquad P_r=v+\frac{(1,0)}{a+c}
\]
Let $T_s$ be the triangle between $[0,v]$ and $[0,P_s]$, and $T_r$ the triangle between $[0,u]$ and $[(1,0)/(a+c),u]$.
For a lattice point $z=\alpha u+\beta v\in T_s$,
\[
\alpha\geq0,\quad\beta\geq\alpha,\quad b\alpha+d\beta\leq d
\]
so only $z=0,v$ are possible. Likewise, if $z=\alpha u-\delta v\in T_r$, then
\[
\delta\geq0,\quad\alpha+\delta\leq1,\quad a\alpha\geq c\delta
\]
so only $z=0,u$ are possible.
Thus moving these segments while connecting their endpoints inside the first and last triangles does not cross a lattice vertex, and preserves the order of crossed edges and triangles. The point $P_s$ lies on the vertical edge of the last triangle of $[0,v]$, and $P_r-v$ lies on the horizontal edge of the first triangle of $[0,u]$. Hence the triangle sequences in the first and last parts agree with those of $\mathcal{PG}(s)$ and $v+\mathcal{PG}(r)$, respectively.
Between $P_s$ and $P_r$, the segment traverses the two triangles of the unit square whose upper left vertex is $v$. The interior of this square is disjoint from the interiors of the preceding and following parts. This gives the required decomposition.
\end{proof}

\begin{lemm}\label{lem:substrip-contributions}
Let $r=a/b$ and $s=c/d$ be reduced fractions satisfying $0<r<s<1$ and $bc-ad=1$, and put
\[
t=r\oplus s=\frac{a+c}{b+d},\qquad u=(b,a),\qquad v=(d,c)
\]
Let $P_s$ and $P_r$ be the intersections of the segment $\ell_t$, before pushing it off, with $x=d$ and $y=c$, respectively. Let $\mathcal{SPG}(s)$ be the union of the closures of the triangles whose interiors are traversed by the open segment from the origin to $P_s$. Define $\mathcal{SPG}(r)$ similarly using the segment from $P_r$ to $u+v$. List signs in the order of passage along a sufficiently small pure left push-off $L_t$, with initial sign $-$ as in $S'_-(t)$ and terminal sign $+$. The actual cuts are made at the intersections of $L_t$ with $x=d,y=c$; even when the reference points $P_s,P_r$ are edge midpoints, approach them from the pure left push-off side. Count neither cut-edge crossing in either substrip, but assign both to the intervening square. If
\[
S'_{-}(s)=(b_1,\dots,b_m),\qquad S'_{-}(r)=(a_1,\dots,a_\ell)
\]
then the run-length sequences obtained by listing, in order along $L_t$, the signs assigned to the triangle-passage and edge-crossing occurrences in the respective unions are
\[
(b_m,\dots,b_1-1,1),\qquad(a_\ell,\dots,a_1)
\]
If $k_s=0$, the crossing occurrence of the edge with midpoint $v/2$ contributes no sign; if $k_r=0$, the crossing occurrence of the edge with midpoint $v+u/2$ contributes no sign.
\end{lemm}
\begin{proof}
The intersections of $\ell_t$ with $x=d$ and $y=c$ are, respectively,
\[
P_s=\frac d{b+d}(u+v)=v-\frac{(0,1)}{b+d},\qquad
P_r=\frac c{a+c}(u+v)=v+\frac{(1,0)}{a+c}
\]
Put $T_s:=\operatorname{conv}\{0,v,P_s\}$ and $T_r:=\operatorname{conv}\{0,u,(1,0)/(a+c)\}$. By the proof of the preceding corollary, moving the segments inside $T_s$ and $v+T_r$ while preserving connections inside the endpoint triangles gives the same triangle-passage sequences. To determine the changes of signs, we enumerate the half-lattice points in these regions.

Since $\det(u,v)=1$, the vectors $u,v$ form a basis of $\mathbb Z^2$. Writing $2z=\alpha u+\beta v$ and examining barycentric coordinates gives
\[
T_s\cap\tfrac12\mathbb Z^2
=\{0,v/2,v\}\quad\text{or}\quad
\{0,v/2,v,(u+v)/2\},
\]
\[
T_r\cap\tfrac12\mathbb Z^2
=\{0,u/2,u\}\quad\text{or}\quad
\{0,u/2,u,(u-v)/2\}.
\]
Indeed, $z\in T_s$ implies $\alpha\geq0$, $\beta\geq\alpha$, and $b\alpha+d\beta\leq2d$. For $z\in T_r$, putting $\delta:=-\beta$ gives $\delta\geq0$, $\alpha+\delta\leq2$, and $a\alpha\geq c\delta$. Enumerating the integer pairs $\alpha,\beta$ satisfying these inequalities gives the two displayed possibilities in each case. After translating $T_r$ by $v$, the possible additional point is $(u+v)/2$, the midpoint of $[0,u+v]$. For $w=u+v$,
\[
\det(w,-u)=\det(w,v)=1>0
\]
so the two curve portions being compared avoid $(u+v)/2$ on the same side. Thus, even when this point occurs, neither the order of triangle-passage and edge-crossing occurrences nor their assigned signs changes there.

Consequently, the only edge in the first part whose signs change is the same edge $E_s$ with midpoint $v/2$. The curve $L_s$ avoids the midpoint on the left, whereas the first part of $L_t$ passes on its right, so its $k_s$ signs are reversed. The last triangle of the first part is traversed from its diagonal edge to its vertical edge, cutting off $v-(0,1)$ on the right. Hence the old terminal sign $+$ becomes $-$. All other triangle signs are determined by the passage sequence and remain unchanged.
By the word representation in Lemma~\ref{lem:half-turn-strongly-admissible}, reversing the central edge signs in the word of $S'_-(s)$ gives the $\dagger$ of the original word. Changing the final $+$ to $-$ then gives run lengths $(b_m,\ldots,b_1-1,1)$. Since $b_1\geq2$, no run has length zero.

Likewise, in the latter part, the only edge whose signs change is the same edge $v+E_r$ with midpoint $v+u/2$, and its $k_r$ signs are reversed. The first triangle of this part is traversed from the horizontal edge to the diagonal edge, cutting off $v+(1,0)$ on the right, so its sign agrees with the initial sign $-$ of $S'_-(r)$. Both terminal signs are $+$. Thus the latter word is obtained from the word of $S'_-(r)$ by reversing only its central edge signs. The same word representation gives run lengths $(a_\ell,\ldots,a_1)$. A central edge of weight zero contributes an empty word, so the same calculations apply in that case as well.
\end{proof}

\begin{proof}[Proof of Proposition~\ref{prop:presnake-relation} (3)]
Figure~\ref{fig:decomposition} shows $r=\frac13$, $t=\frac25$, and $s=\frac12$.
Write $r=a/b$, $s=c/d$, $u=(b,a)$, and $v=(d,c)$, and use the cut points $P_s,P_r$ of Lemma~\ref{lem:substrip-contributions}.
By Corollary~\ref{lem:decomposition-lemma}, the region traversed by $L_t$ consists, in order, of $\mathcal{SPG}(s)$, a unit square, and $\mathcal{SPG}(r)$. Lemma~\ref{lem:substrip-contributions} gives the run lengths $(b_m,\ldots,b_1-1,1)$ and $(a_\ell,\ldots,a_1)$ for the first and last words.
The central part counts the vertical edge crossing corresponding to $P_s$, two triangle passages, the intervening diagonal edge crossing, and the horizontal edge crossing corresponding to $P_r$. Both triangles cut off $v$ on the left. The vertical intersection lies $1/(b+d)$ below $v$, and the horizontal intersection lies $1/(a+c)$ to its right, so their midpoints lie to the right of the curve. If $a+c=2$, the horizontal intersection is the midpoint, but the pure left push-off moves the intersection to its left, so the midpoint is still on the right. The diagonal intersection is displaced from $v$ by $(1,-1)/(a+b+c+d)$, so its midpoint is also on the right. Hence the central word is $+^{2+k_{\sigma(1)}+k_{\sigma(2)}+k_{\sigma(3)}}=+^{K-1}$. The last sign of the first part and the first sign of the last part are both negative, so no runs merge. This proves the formula in (3).

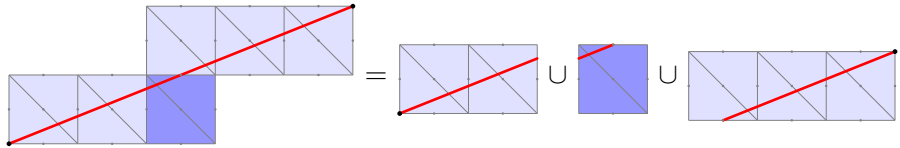
\begin{figure}[ht]
    \centering
    \begin{tikzpicture}[x=0.01cm,y=0.01cm,scale=0.35,baseline=14mm]
  \definecolor{annred}{RGB}{255,0,0}

  \tikzset{
    gridline/.style={gray,line width=0.15pt,line cap=butt,line join=miter},
    gammaL/.style={draw=annred,line width=1pt,line cap=round,line join=round},
    enddot/.style={circle,fill=black,inner sep=0pt,minimum size=1.8pt},
    midpt/.style={circle,draw=gray,fill=gray,inner sep=0pt,minimum size=0.5pt}
  }

  \begin{scope}[shift={(180,170)}, x={(260,0)}, y={(0,260)}]

    \foreach \i/\j in {0/0, 1/0, 2/1, 3/1, 4/1}{
      \fill[blue!12] (\i,\j) rectangle ({\i+1},{\j+1});
    }
    \fill[blue!42] (2,0) rectangle (3,1);

    \foreach \xa/\ya/\xb/\yb in {
      0/0/1/0, 1/0/2/0, 2/0/3/0,
      0/1/1/1, 1/1/2/1, 2/1/3/1, 3/1/4/1, 4/1/5/1,
      2/2/3/2, 3/2/4/2, 4/2/5/2,
      0/0/0/1,
      1/0/1/1,
      2/0/2/1, 2/1/2/2,
      3/0/3/1, 3/1/3/2,
      4/1/4/2,
      5/1/5/2
    }{
      \draw[gridline] (\xa,\ya) -- (\xb,\yb);
    }

    \foreach \i/\j in {
      0/0, 1/0, 2/0, 2/1, 3/1, 4/1
    }{
      \draw[gridline] (\i,{\j+1}) -- ({\i+1},\j);
    }

    \coordinate (A) at (0,0);
    \coordinate (B) at (5,2);

    \draw[gammaL] (A) -- (B);

    \foreach \xa/\ya/\xb/\yb in {
      0/0/1/0, 1/0/2/0, 2/0/3/0,
      0/1/1/1, 1/1/2/1, 2/1/3/1, 3/1/4/1, 4/1/5/1,
      2/2/3/2, 3/2/4/2, 4/2/5/2,
      0/0/0/1,
      1/0/1/1,
      2/0/2/1, 2/1/2/2,
      3/0/3/1, 3/1/3/2,
      4/1/4/2,
      5/1/5/2
    }{
      \node[midpt] at ({(\xa+\xb)/2},{(\ya+\yb)/2}) {};
    }

    \foreach \i/\j in {
      0/0, 1/0, 2/0, 2/1, 3/1, 4/1
    }{
      \node[midpt] at ({\i+0.5},{\j+0.5}) {};
    }

    \node[enddot] at (A) {};
    \node[enddot] at (B) {};

  \end{scope}
\end{tikzpicture}
$=$
\begin{tikzpicture}[x=0.01cm,y=0.01cm,scale=0.35,baseline=10mm]
  \definecolor{annred}{RGB}{255,0,0}

  \tikzset{
    gridline/.style={gray,line width=0.15pt,line cap=butt,line join=miter},
    gammaL/.style={draw=annred,line width=1pt,line cap=round,line join=round},
    enddot/.style={circle,fill=black,inner sep=0pt,minimum size=1.8pt},
    midpt/.style={circle,draw=gray,fill=gray,inner sep=0pt,minimum size=0.5pt}
  }

  \begin{scope}[shift={(180,170)}, x={(260,0)}, y={(0,260)}]

    \fill[blue!12] (0,0) rectangle (2,1);

    \foreach \i in {0,...,1}{
      \draw[gridline] (\i,1) -- ({\i+1},0);
    }

    \foreach \i in {0,...,2}{
      \draw[gridline] (\i,0) -- (\i,1);
    }
    \foreach \j in {0,1}{
      \draw[gridline] (0,\j) -- (2,\j);
    }

    \draw[gammaL] (0,0) -- (2,0.8);

    \foreach \i in {0,...,1}{
      \foreach \j in {0,1}{
        \node[midpt] at ({\i+0.5},\j) {};
      }
    }

    \foreach \i in {0,...,2}{
      \node[midpt] at (\i,0.5) {};
    }

    \foreach \i in {0,...,1}{
      \node[midpt] at ({\i+0.5},0.5) {};
    }

    \node[enddot] at (0,0) {};

  \end{scope}
\end{tikzpicture}
$\cup$
\begin{tikzpicture}[x=0.01cm,y=0.01cm,scale=0.35,baseline=10mm]
  \definecolor{annred}{RGB}{255,0,0}

  \tikzset{
    gridline/.style={gray,line width=0.15pt,line cap=butt,line join=miter},
    gammaL/.style={draw=annred,line width=1pt,line cap=round,line join=round},
    enddot/.style={circle,fill=black,inner sep=0pt,minimum size=1.8pt},
    midpt/.style={circle,draw=gray,fill=gray,inner sep=0pt,minimum size=0.5pt}
  }

  \begin{scope}[shift={(180,170)}, x={(260,0)}, y={(0,260)}]

    \fill[blue!42] (2,0) rectangle (3,1);

    \draw[gridline] (2,1) -- (3,0);

    \draw[gridline] (2,0) -- (2,1);
    \draw[gridline] (3,0) -- (3,1);
    \draw[gridline] (2,0) -- (3,0);
    \draw[gridline] (2,1) -- (3,1);

    \draw[gammaL] (2,0.8) -- (2.5,1);

    \node[midpt] at (2.5,0) {};
    \node[midpt] at (2.5,1) {};

    \node[midpt] at (2,0.5) {};
    \node[midpt] at (3,0.5) {};

    \node[midpt] at (2.5,0.5) {};

  \end{scope}
\end{tikzpicture}
$\cup$
\begin{tikzpicture}[x=0.01cm,y=0.01cm,scale=0.35,baseline=20mm]
  \definecolor{annred}{RGB}{255,0,0}

  \tikzset{
    gridline/.style={gray,line width=0.15pt,line cap=butt,line join=miter},
    gammaL/.style={draw=annred,line width=1pt,line cap=round,line join=round},
    enddot/.style={circle,fill=black,inner sep=0pt,minimum size=1.8pt},
    midpt/.style={circle,draw=gray,fill=gray,inner sep=0pt,minimum size=0.5pt}
  }

  \begin{scope}[shift={(180,170)}, x={(260,0)}, y={(0,260)}]

    \fill[blue!12] (2,1) rectangle (5,2);

    \foreach \i in {2,...,4}{
      \draw[gridline] (\i,2) -- ({\i+1},1);
    }

    \foreach \i in {2,...,5}{
      \draw[gridline] (\i,1) -- (\i,2);
    }
    \foreach \j in {1,2}{
      \draw[gridline] (2,\j) -- (5,\j);
    }

    \draw[gammaL] (2.5,1) -- (5,2);

    \foreach \i in {2,...,4}{
      \foreach \j in {1,2}{
        \node[midpt] at ({\i+0.5},\j) {};
      }
    }

    \foreach \i in {2,...,5}{
      \node[midpt] at (\i,1.5) {};
    }

    \foreach \i in {2,...,4}{
      \node[midpt] at ({\i+0.5},1.5) {};
    }

    \node[enddot] at (5,2) {};

  \end{scope}
\end{tikzpicture}
    \caption{Decomposition of $\mathcal{PG}(t)$}
    \label{fig:decomposition}
\end{figure}

The first and last parts are shaded light blue, and the central unit square is shaded dark blue. In the word decomposition above, the signs of the edges entering and leaving the central square are included in the central part.
\end{proof}

The following matrix lemma translates Proposition~\ref{prop:presnake-relation} directly into the GC matrix recurrence. Put
\[
H:=\begin{bmatrix}K&-1\\1&0\end{bmatrix},
\qquad
R:=\begin{bmatrix}1&1\\0&-1\end{bmatrix}.
\]
For $u\in(0,1)\cap\mathbb Q$, put
\[
S(u)=(K-1,1,c_2,\dots,c_n),\qquad
Q_u:=F_{S'_-(u)}
\]
At the boundary, define
\[
Q_0:=\begin{bmatrix}1&-k_{\sigma(1)}\\0&1\end{bmatrix},
\qquad
Q_1:=\begin{bmatrix}k_{\sigma(2)}+2&1\\k_{\sigma(2)}+1&1\end{bmatrix}
\]
Direct calculation gives $C_0=HQ_0$ and $C_1=HQ_1$.

\begin{lemm}\label{lem:matrix-recursion-strongly-admissible}
For every $u\in(0,1)\cap\mathbb Q$,
\[
F_{S(u)}=HQ_u
\]
Moreover, for $(r,t,s)\in\mathrm F\mathbb T$ with $t\in(0,1)$,
\[
Q_t=(Q_rHQ_s)^T
\]
where the boundary matrices above are used when $r=0$ or $s=1$.
\end{lemm}
\begin{proof}
The first identity follows immediately from
\[
F_{K-1}F_1F_c=HF_{c+1}
\]
For the second, note that
\[
F_{x-1}F_1=F_xR,\qquad RF_{K-1}=H^T
\]
If neither $r$ nor $s$ is a boundary point, Proposition~\ref{prop:presnake-relation} (3) gives
\[
Q_t=F_{b_m}\cdots F_{b_2}F_{b_1-1}F_1F_{K-1}F_{a_\ell}\cdots F_{a_1}
=Q_s^TRF_{K-1}Q_r^T=Q_s^TH^TQ_r^T=(Q_rHQ_s)^T.
\]

If $r=0$, put $L:=2+k_{\sigma(2)}+k_{\sigma(3)}
=K-1-k_{\sigma(1)}$. Proposition~\ref{prop:presnake-relation} (1), together with
\[
RF_L=H^TQ_0^T
\]
gives the same identity. If $s=1$, it follows from Proposition~\ref{prop:presnake-relation} (2) and
\[
F_{k_{\sigma(2)}+2}F_{K-1}=Q_1^TH^T
\]
For $t=1/2$, when both are boundary points, we calculate directly
\[
S(\tfrac12)=(K-1,1,1+k_{\sigma(2)},
2+k_{\sigma(2)}+k_{\sigma(3)})
\]
and multiply the corresponding matrices to verify $Q_{1/2}=(Q_0HQ_1)^T$.
\end{proof}

\begin{proof}[Proof of Theorem~\ref{continued-fraction-theorem2}]
We use the fact that a matrix in $SL(2,\mathbb R)$ with nonzero $(2,1)$ entry is uniquely determined by its bottom row and trace.

For $t=1$, the result follows directly from $S(1)=(K-1,k_{\sigma(2)}+2)$. Suppose $t\in(0,1)$, and let $(r,t,s)$ be the Farey triple with middle entry $t$. Using $Q_0,Q_1$ above in the boundary cases, write inductively
\[
C_r=HQ_r,\qquad C_s=HQ_s,\qquad
X:=Q_rHQ_s=\begin{bmatrix}x_{11}&x_{12}\\x_{21}&x_{22}\end{bmatrix}
\]
The GC matrix recurrence and Lemma~\ref{lem:matrix-recursion-strongly-admissible} give
\[
C_t=HX-D_t,\qquad F_{S(t)}=HX^T
\]
Since $(C_t)_{21}=x_{11}$, the trace formula for $C_t$ gives
\[
x_{12}-x_{21}=k_t
\]
Thus $C_t$ and $F_{S(t)}$ both have bottom row $(x_{11},x_{21})$, and
\[
\operatorname{tr}(F_{S(t)})=Kx_{11}-k_t=\operatorname{tr}(C_t)
\]
also holds.
Both determinants are $1$, and $x_{11}=m_t>0$, so $C_t=F_{S(t)}$. This proves the assertion for $t\in(0,1]$.

The case $t=\infty$ follows directly from $S(\infty)=(1+k_{\sigma(1)}+k_{\sigma(2)},1)$. For $1<t<\infty$, put $u=1/t$ and, in the dual tree, write
\[
S^*(u)=(b_0,b_1,\dots,b_m),\qquad
A:=F_{b_1}\cdots F_{b_m}
=\begin{bmatrix}p&q\\r&z\end{bmatrix}
\]
By the case already proved, $C_u^*=F_{b_0}A$, with $b_0=K-1$ and $b_1=1$. Corollary~\ref{cor:dual-remark} gives its bottom row as
\[
(p,q)=(m_t,u_u^*)
\]
and $k_u^*=k_t$. The trace formula gives
\[
r=p-q-k_t
\]
whose right-hand side equals $u_t$ by Proposition~\ref{prop:t-1/t-relation-gen}. On the other hand, Lemma~\ref{rem:difference-(0,1)(1,infty)}~(7) gives
\[
S(t)=(b_0,b_m,\dots,b_1)
\]
so $F_{S(t)}$ and $C_t$ both have bottom row
\[
(p,r)=(m_t,u_t)
\]
Moreover,
\[
\operatorname{tr}(F_{S(t)})=(K-1)p+q+r
=Kp-k_t=\operatorname{tr}(C_t).
\]
Again both determinants are $1$ and $p=m_t>0$, so $F_{S(t)}=C_t$. This completes all cases.
\end{proof}

\begin{rema}
The matrices $C_{\frac01}$ and $F_{S(\frac01)}$ do not agree. Indeed, $C_{\frac01}$ has negative entries.
\end{rema}

\chapter{Generalized Discrete Markov Spectra}\label{chap:generalized-discrete-markov-spectra}
By Chapter~\ref{chap:generalized-cohn-matrices}, we have defined generalized Markov numbers, characteristic numbers, GM distances, generalized Cohn matrices, and generalized strongly admissible sequences, and established their basic properties. In this chapter we use these constructions to define a discrete set of values naturally associated with generalized Markov numbers and show that every such value is realized as both a Lagrange constant and a Markov constant. This connects the arithmetic, combinatorial, geometric, and matrix-theoretic results of Part~II with the Lagrange and Markov spectra.

We first define the generalized discrete Markov spectrum and state the main theorem of this text. We then realize its elements as Lagrange constants of quadratic irrationals and as Markov constants of binary quadratic forms with rational coefficients. Thus generalized Markov numbers give explicit values in the Lagrange and Markov spectra. Specializing to $(k_1,k_2,k_3)=(0,0,0)$, we derive the classical Markov theorem from the preceding definitions and results. We next approximate irrational slopes by rational slopes and prove that the bi-infinite sequences obtained from irrational-slope lines have Markov value $3+k_1+k_2+k_3$. Finally, we discuss the relation between the $(0,0,0)$- and $(2,2,2)$-types and natural generalizations of Frobenius's uniqueness conjecture.

The discussion of generalized discrete Markov spectra and uniqueness in this chapter is based primarily on \cite{gyoda-generalized}. For the uniqueness conjecture and its generalizations, we also refer to \cite{frobenius,gyomatsu}. Our proof of Markov's theorem follows the basic strategy of the traditional accounts \cite{bombieri2,aig,reutenauer}, but introduces right and left mechanical words to make the relation with strongly admissible sequences explicit. Example~\ref{ex:noninjective-gm-tree} is due to Nakabayashi Shoma.

\section{Definitions and Main Theorems}\label{sec:generalized-discrete-definition}\index{Generalized discrete Markov spectrum}
For $(k_1,k_2,k_3)\in\mathbb Z_{\geq0}^3$ and $\sigma\in\mathfrak S_3$, define
\begin{align*}
\mathcal M_{k_1,k_2,k_3,\sigma}
&:=\left\{\frac{\sqrt{((3+k_1+k_2+k_3)m_t-k_t)^2-4}}{m_t}
\ \middle|\ t\in\mathbb Q_{\geq0}\cup\{\infty\}\right\},\\
\mathcal M_{k_1,k_2,k_3}
&:=\bigcup_{\sigma\in\mathfrak S_3}\mathcal M_{k_1,k_2,k_3,\sigma}.
\end{align*}
Here $m_t$ is the GM number with position label $i_t$, and $k_t:=k_{i_t}$. Throughout this section put $K:=3+k_1+k_2+k_3$. We call $\mathcal M_{k_1,k_2,k_3}$ the \emph{$(k_1,k_2,k_3)$-generalized discrete Markov spectrum}. The main theorem of this section, and of this text, is the following.

\begin{theorem}\label{thm:markov-value-gen}\index{Lagrange spectrum}\index{Generalized discrete Markov spectrum}\index{Lagrange constant}
Fix $(k_1,k_2,k_3)\in\mathbb Z_{\geq0}^3$ and $\sigma\in\mathfrak S_3$. We use the conventions $1/0=\infty$ and $1/\infty=0$.
For every reduced fraction $t\in[0,\infty]$, let $(m_t,i_t)$ be the corresponding GM number together with its position label, and let $S(t)$ be the corresponding generalized strongly admissible sequence. For a finite sequence $S$ of positive integers, write $\alpha_S:=[\overline S]$. Then
\[
\mathcal L(\alpha_{S(t)})
=\mathcal L(\alpha_{S^\ast(1/t)})
=\frac{\sqrt{((3+k_1+k_2+k_3)m_t-k_t)^2-4}}{m_t}.
\]
In particular, $\mathcal M_{k_1,k_2,k_3}\subset\mathcal L$.
\end{theorem}
Assuming this theorem, Theorem~\ref{thm:quadratic-markov-lagrange1770} immediately gives the following result.
\begin{theorem}\label{thm:markov-value-gen2}\index{Markov constant}
With the notation of Theorem~\ref{thm:markov-value-gen}, for a finite sequence $S$ of positive integers put
$Q_S=(x-\alpha_Sy)(x-\alpha'_Sy)$, where $\alpha'_S$ is the quadratic conjugate of $\alpha_S$. Then, for every reduced fraction $t\in[0,\infty]$,
\[
\mathcal M(Q_{S(t)})
=\mathcal M(Q_{S^\ast(1/t)})
=\frac{\sqrt{((3+k_1+k_2+k_3)m_t-k_t)^2-4}}{m_t}.
\]
\end{theorem}
We now prove Theorem~\ref{thm:markov-value-gen} using the preceding results.

\begin{proof}[Proof of Theorem~\ref{thm:markov-value-gen}]
By Theorem~\ref{thm:L=S} and the reversal invariance of $\mathcal S$, the equality
$\mathcal L(\alpha_{S(t)})=\mathcal L(\alpha_{S^\ast(1/t)})$
follows if the periodic blocks agree up to reversal and cyclic shift. This is precisely Lemma~\ref{rem:difference-(0,1)(1,infty)}(7).
It remains to establish the formula for $\mathcal L(\alpha_{S(t)})$.

First consider $t=0,\infty$. Put
\[
a_t:=\begin{cases}
1+k_{\sigma(2)}+k_{\sigma(3)}&(t=0),\\
1+k_{\sigma(1)}+k_{\sigma(2)}&(t=\infty).
\end{cases}
\]
Then $S(t)=(a_t,1)$, $m_t=1$, and $k_t$ is respectively $k_{\sigma(1)}$ or $k_{\sigma(3)}$. Hence
\[
Km_t-k_t=a_t+2,\qquad
F_{S(t)}=\begin{bmatrix}a_t+1&a_t\\1&1\end{bmatrix}.
\]
The $(2,1)$-entries of the two cyclic shifts are $1$ and $a_t$. Since $a_t\geq1$, Theorem~\ref{thm:lagrange1770-continued-fraction-matrix} gives
\[
\mathcal L(\alpha_{S(t)})
=\sqrt{(a_t+2)^2-4}
=\frac{\sqrt{(Km_t-k_t)^2-4}}{m_t}.
\]
We may therefore assume $t\in(0,\infty)$.

Write $S(t)=(a_0,\dots,a_n)$ and put
$S_i:=(a_i,\dots,a_n,a_0,\dots,a_{i-1})$ for $0\leq i\leq n$.
Theorem~\ref{thm:lagrange1770-continued-fraction-matrix} gives
\[
\mathcal L(\alpha_{S(t)})
=\max\left\{\frac{\sqrt{(\operatorname{tr}(F_{S_i}))^2-(-1)^{n+1}\cdot4}}
{(F_{S_i})_{21}}\ \middle|\ 0\leq i\leq n\right\}.
\]
By Lemma~\ref{rem:difference-(0,1)(1,infty)}(0), $(-1)^{n+1}=1$.
Since $F_{S_i}$ is a cyclic permutation of the factors in $F_{a_0}\cdots F_{a_n}$, the cyclic invariance of the trace,
\[
\operatorname{tr}(A_0\cdots A_n)
=\operatorname{tr}(A_i\cdots A_nA_0\cdots A_{i-1}),
\]
implies $\operatorname{tr}(F_{S_i})=\operatorname{tr}(F_{S(t)})$ for every $i$.
Thus the numerator is independent of $i$, and it suffices to minimize $(F_{S_i})_{21}$.
For $i=0$, Theorems~\ref{continued-fraction-theorem2} and~\ref{thm:Mt-description} give
\[
F_{S_0}=F_{S(t)}=C_t,\qquad
\operatorname{tr}(C_t)=Km_t-k_t,\qquad
(F_{S_0})_{21}=m_t.
\]
It remains to prove
\[
\min\{(F_{S_i})_{21}\mid 0\leq i\leq n\}=(F_{S_0})_{21}.
\]

Since $(F_{S_0})_{21}=N(a_1,\dots,a_n)$, this is a statement about continuants. Extend the indices periodically by $a_{j+n+1}=a_j$ and set
\[
w_k:=(a_{k+1},a_{k+2},\dots,a_{k+n}),\qquad
N_t:=\{N(w_k)\mid0\leq k\leq n\}.
\]
We must show that $N(w_0)$ is the least element of $N_t$.
Projecting $\widetilde{\mathbb R^2}$ to the triangulated once-punctured torus turns $\overline{L_t}$ into a loop; see Remark~\ref{rem:one-punctured-torus}.
For each $k$, let $\overline{L_t}(w_k)$ be the portion of this loop whose triangle-passage and edge-crossing signs, in occurrence order, have successive constant-sign run lengths $w_k$.

Reconnect the endpoints to lattice points according to Table~\ref{table1}, and denote the resulting arc by $\widetilde{L_t}(w_k)$.
When a correction spans two triangles, omit the terminal edge crossing and join directly to the opposite lattice vertex in the next triangle. The target vertex is the same as in the corresponding one-triangle case immediately above it in the table. Keep the intersections with all remaining edges, and the curve beyond the correction region, fixed. Figure~\ref{fig:modification} illustrates the construction.

\begin{table}[htbp]
\centering
\begingroup
\tikzset{
  ep edge/.style={draw=black!65,line width=.55pt},
  ep arc/.style={draw=annred,line width=.85pt},
  ep sign/.style={text=annred,fill=white,inner sep=.65pt,font=\small},
  ep point/.style={circle,fill=annred,inner sep=1.15pt}
}
\newcommand{\epUpper}{\draw[ep edge] (0,1)--(1,1)--(1,0)--cycle;}
\newcommand{\epLower}{\draw[ep edge] (0,0)--(0,1)--(1,0)--cycle;}
\newcommand{\epSquare}{\draw[ep edge] (0,0) rectangle (1,1);}
\newcommand{\epCell}[2][-.16]{%
  \begin{tikzpicture}[x=1.22cm,y=1.22cm,line cap=round,line join=round,
    baseline=(current bounding box.center)]
    \path[use as bounding box] (-1.13,#1) rectangle (1.13,1.16);
    #2
  \end{tikzpicture}%
}
\newcommand{\epFixedArc}[3]{%
  \draw[ep arc] #1--#2--#3;
  \node[ep point] at #1 {};
}
\newcommand{\epAfter}[1]{\epCell{\begin{scope}[xshift=-.61cm]
  \ifcase#1\relax
  \or %
    \epUpper
    \epFixedArc{(0,1)}{(1,{693/830})}{(1.17,.88)}
    \node[ep sign] at (.76,.51) {$-$};
  \or %
    \epUpper
    \epFixedArc{(0,1)}{(1,{1377/1850})}{(1.17,.85)}
    \node[ep sign] at (.76,.51) {$-$};
  \or %
    \epUpper
    \epFixedArc{(1,0)}{({1311/1775},1)}{(.78,1.14)}
    \node[ep sign] at (.43,.80) {$+$};
  \or %
    \epUpper
    \epFixedArc{(1,0)}{({383/690},1)}{(.62,1.14)}
    \node[ep sign] at (.43,.80) {$+$};
  \or %
    \epLower
    \epFixedArc{(0,0)}{({3901/6300},{2399/6300})}{(.94,.49)}
    \node[ep sign] at (.27,.54) {$+$};
  \or %
    \epLower
    \epFixedArc{(0,0)}{({1997/5450},{3453/5450})}{(.72,.75)}
    \node[ep sign] at (.58,.18) {$+$};
  \or %
    \epLower
    \epFixedArc{(0,0)}{({263/428},{165/428})}{(.77,.55)}
    \node[ep sign] at (.24,.48) {$-$};
  \or %
    \epLower
    \epFixedArc{(1,0)}{(0,{383/2000})}{(-.14,.16)}
    \node[ep sign] at (.25,.50) {$+$};
  \or %
    \epLower
    \epFixedArc{(1,0)}{(0,{293/475})}{(-.14,.58)}
    \node[ep sign] at (.30,.20) {$+$};
  \or %
    \epLower
    \epFixedArc{(0,1)}{({749/3950},0)}{(.16,-.13)}
    \node[ep sign] at (.48,.18) {$-$};
  \or %
    \epLower
    \epFixedArc{(0,1)}{({227/500},0)}{(.42,-.13)}
    \node[ep sign] at (.63,.12) {$-$};
  \or %
    \epUpper
    \epFixedArc{(1,1)}{({1188/2575},{1387/2575})}{(.24,.46)}
    \node[ep sign] at (.51,.82) {$-$};
  \or %
    \epUpper
    \epFixedArc{(1,1)}{({413/515},{102/515})}{(.28,.08)}
    \node[ep sign] at (.51,.82) {$-$};
  \or %
    \epUpper
    \epFixedArc{(1,1)}{({73/136},{63/136})}{(.46,.25)}
    \node[ep sign] at (.81,.43) {$+$};
  \or %
    \epLower
    \epFixedArc{(0,0)}{({1581/2200},{619/2200})}{(.83,.52)}
    \node[ep sign] at (.24,.48) {$-$};
  \or %
    \epUpper
    \epFixedArc{(1,1)}{({207/680},{473/680})}{(.18,.35)}
    \node[ep sign] at (.81,.43) {$+$};
  \fi
\end{scope}}}
\newcommand{\epBefore}[1]{%
  \ifnum#1>14\relax\def\epBottom{-1.16}\else\def\epBottom{-.16}\fi
  \epCell[\epBottom]{%
  \ifcase#1\relax
  \or %
    \begin{scope}[xshift=-.61cm]
      \epUpper
      \draw[ep arc] (.34,.66)--(1.17,.88);
      \node[ep point] at (.34,.66) {};
      \node[ep sign] at (.76,.47) {$-$};
    \end{scope}
  \or %
    \begin{scope}[xshift=-.61cm]
      \epSquare
      \draw[ep edge] (0,1)--(.22,.78) (.43,.57)--(1,0);
      \draw[ep arc] (.43,.39)--(1.17,.85);
      \node[ep point] at (.43,.39) {};
      \node[ep sign] at (.32,.68) {$-$};
      \node[ep sign] at (.68,.74) {$-$};
    \end{scope}
  \or %
    \begin{scope}[xshift=-.61cm]
      \epUpper
      \draw[ep arc] (.57,.43)--(.78,1.14);
      \node[ep point] at (.57,.43) {};
      \node[ep sign] at (.42,.80) {$+$};
    \end{scope}
  \or %
    \begin{scope}[xshift=-.61cm]
      \epSquare
      \draw[ep edge] (0,1)--(.42,.58) (.66,.34)--(1,0);
      \draw[ep arc] (.30,.45)--(.62,1.14);
      \node[ep point] at (.30,.45) {};
      \node[ep sign] at (.54,.46) {$+$};
      \node[ep sign] at (.78,.73) {$+$};
    \end{scope}
  \or %
    \begin{scope}[xshift=-.61cm]
      \epLower
      \draw[ep arc] (0,.17)--(.94,.49);
      \node[ep point] at (0,.17) {};
      \node[ep sign] at (.27,.47) {$+$};
    \end{scope}
  \or %
    \draw[ep edge] (-1,1)--(0,1)--(1,0)--(0,0)--cycle
      (0,1)--(0,.40) (0,.12)--(0,0);
    \draw[ep arc] (-.10,.48)--(.72,.75);
    \node[ep point] at (-.10,.48) {};
    \node[ep sign] at (0,.26) {$+$};
    \node[ep sign] at (.30,.28) {$+$};
  \or %
    \begin{scope}[xshift=-.61cm]
      \epLower
      \draw[ep arc] (.25,0)--(.77,.55);
      \node[ep point] at (.25,0) {};
      \node[ep sign] at (.25,.45) {$-$};
    \end{scope}
  \or %
    \begin{scope}[xshift=-.61cm]
      \epLower
      \draw[ep arc] (-.14,.16)--(.66,.34);
      \node[ep point] at (.66,.34) {};
      \node[ep sign] at (.24,.57) {$+$};
    \end{scope}
  \or %
    \begin{scope}[xshift=-.61cm]
      \epSquare
      \draw[ep edge] (0,1)--(.42,.58) (.66,.34)--(1,0);
      \draw[ep arc] (-.14,.58)--(.43,.73);
      \node[ep point] at (.43,.73) {};
      \node[ep sign] at (.54,.46) {$+$};
      \node[ep sign] at (.24,.25) {$+$};
    \end{scope}
  \or %
    \begin{scope}[xshift=-.61cm]
      \epLower
      \draw[ep arc] (.16,-.13)--(.34,.66);
      \node[ep point] at (.34,.66) {};
      \node[ep sign] at (.47,.17) {$-$};
    \end{scope}
  \or %
    \begin{scope}[xshift=-.61cm]
      \epSquare
      \draw[ep edge] (0,1)--(.22,.78) (.43,.57)--(1,0);
      \draw[ep arc] (.42,-.13)--(.59,.52);
      \node[ep point] at (.59,.52) {};
      \node[ep sign] at (.32,.68) {$-$};
      \node[ep sign] at (.17,.43) {$-$};
    \end{scope}
  \or %
    \begin{scope}[xshift=-.61cm]
      \epUpper
      \draw[ep arc] (.24,.46)--(1,.73);
      \node[ep point] at (1,.73) {};
      \node[ep sign] at (.55,.82) {$-$};
    \end{scope}
  \or %
    \draw[ep edge] (-1,1)--(0,1)--(1,0)--(0,0)--cycle
      (0,1)--(0,.64) (0,.40)--(0,0);
    \draw[ep arc] (-.72,.08)--(.12,.27);
    \node[ep point] at (.12,.27) {};
    \node[ep sign] at (0,.52) {$-$};
    \node[ep sign] at (-.51,.79) {$-$};
  \or %
    \begin{scope}[xshift=-.61cm]
      \epUpper
      \draw[ep arc] (.46,.25)--(.73,1);
      \node[ep point] at (.73,1) {};
      \node[ep sign] at (.82,.43) {$+$};
    \end{scope}
  \or %
    \begin{scope}[xshift=-.61cm]
      \draw[ep edge] (0,1)--(1,0)--(1,-1)--(0,0)--cycle
        (0,0)--(.18,0) (.42,0)--(1,0);
      \draw[ep arc] (.55,-.08)--(.83,.52);
      \node[ep point] at (.55,-.08) {};
      \node[ep sign] at (.30,0) {$-$};
      \node[ep sign] at (.23,.43) {$-$};
    \end{scope}
  \or %
    \begin{scope}[xshift=-.61cm]
      \draw[ep edge] (0,1)--(1,0)--(1,-1)--(0,0)--cycle
        (0,0)--(.52,0) (.76,0)--(1,0);
      \draw[ep arc] (.18,-.65)--(.45,.10);
      \node[ep point] at (.45,.10) {};
      \node[ep sign] at (.64,0) {$+$};
      \node[ep sign] at (.78,-.43) {$+$};
    \end{scope}
  \fi
  }%
}
\setlength{\tabcolsep}{4pt}
\renewcommand{\arraystretch}{1.12}
\begin{tabular}{|c|c||c|c|}
\multicolumn{1}{c}{\shortstack{Endpoints of\\$\overline{L_t}(w_k)$}}&
\multicolumn{1}{c}{\shortstack{Correction to\\$\widetilde{L_t}(w_k)$}}&
\multicolumn{1}{c}{\shortstack{Endpoints of\\$\overline{L_t}(w_k)$}}&
\multicolumn{1}{c}{\shortstack{Correction to\\$\widetilde{L_t}(w_k)$}}\\\hline
\epBefore{1}&\epAfter{1}&\epBefore{8}&\epAfter{8}\\\hline
\epBefore{2}&\epAfter{2}&\epBefore{9}&\epAfter{9}\\\hline
\epBefore{3}&\epAfter{3}&\epBefore{10}&\epAfter{10}\\\hline
\epBefore{4}&\epAfter{4}&\epBefore{11}&\epAfter{11}\\\hline
\epBefore{5}&\epAfter{5}&\epBefore{12}&\epAfter{12}\\\hline
\epBefore{6}&\epAfter{6}&\epBefore{13}&\epAfter{13}\\\hline
\epBefore{7}&\epAfter{7}&\epBefore{14}&\epAfter{14}\\\hline
\epBefore{15}&\epAfter{15}&\epBefore{16}&\epAfter{16}\\\hline
\end{tabular}
\endgroup
\par\smallskip
\caption{Endpoint corrections}
\label{table1}
\end{table}

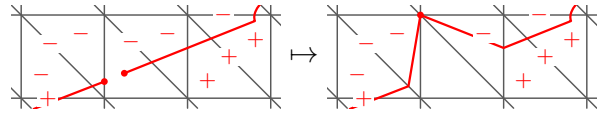
\begin{figure}[htbp]
\centering
\begingroup
\tikzset{
  ep example edge/.style={draw=black!65,line width=.55pt},
  ep example arc/.style={draw=annred,line width=.85pt},
  ep example sign/.style={text=annred,fill=white,inner sep=.5pt,
    font=\scriptsize},
  ep example point/.style={circle,fill=annred,inner sep=.9pt}
}
\newcommand{\epExampleGrid}{%
  \foreach \y in {0,1}
    \draw[ep example edge] (-.12,\y)--(3.12,\y);
  \foreach \x in {0,1,2,3}
    \draw[ep example edge] (\x,-.12)--(\x,1.12);
  \foreach \k in {0,1,2,3,4}
    \draw[ep example edge]
      ({\k-1.12},1.12)--({\k+.12},-.12);
}
\newcommand{\epExampleSigns}[1]{%
  \node[ep example sign] at (.24,.28) {$-$};
  \node[ep example sign] at (.35,.65) {$-$};
  \node[ep example sign] at (.72,.78) {$-$};
  \node[ep example sign] at (.32,0) {$+$};
  \ifnum#1=0\relax
    \node[ep example sign] at (1.32,.68) {$-$};
  \fi
  \node[ep example sign] at (1.76,.77) {$-$};
  \node[ep example sign] at (2.24,.21) {$+$};
  \node[ep example sign] at (2.50,.50) {$+$};
  \node[ep example sign] at (2.82,.70) {$+$};
  \node[ep example sign] at (2.43,1) {$-$};
}
\begin{tikzpicture}[x=1.1cm,y=1.1cm,line cap=round,line join=round]
  \path[use as bounding box] (-.12,-.12) rectangle (6.92,1.12);
  \begin{scope}[overlay]
    \clip (-.12,-.12) rectangle (3.12,1.12);
    \epExampleGrid
    \draw[ep example arc] (.15,-.13)--(1,.20);
    \draw[ep example arc] (1.24,.30)--(2.81,.928)
      .. controls (2.78,1.01) and (2.84,1.10) .. (2.90,1.16);
    \node[ep example point] at (1,.20) {};
    \node[ep example point] at (1.24,.30) {};
    \epExampleSigns{0}
  \end{scope}
  \node at (3.40,.50) {$\mapsto$};
  \begin{scope}[xshift=4.18cm,overlay]
    \clip (-.12,-.12) rectangle (3.12,1.12);
    \epExampleGrid
    \draw[ep example arc] (.15,-.13)--({101/118},{17/118})--(1,1);
    \draw[ep example arc] (1,1)--(2,.604)--(2.81,.928)
      .. controls (2.78,1.01) and (2.84,1.10) .. (2.90,1.16);
    \node[ep example point] at (1,1) {};
    \epExampleSigns{1}
  \end{scope}
\end{tikzpicture}
\endgroup
\caption{An example of endpoint correction}
\label{fig:modification}
\end{figure}

We check the displacement of the endpoints. For a positive-slope line, triangle signs on opposite sides of a horizontal or vertical edge are opposite. If the signs on opposite sides of a diagonal edge agree, the incoming and outgoing edges have the same common vertex. Moreover, an edge sign agrees with at least one adjacent triangle sign. Thus a maximal constant-sign block contains either one triangle passage or two adjacent triangle passages, and the incoming and outgoing edges of these passages meet at the same lattice vertex. The correction in the table joins both cuts bordering such a block to this vertex.

The two corresponding copies of the omitted run, in successive periods, differ by translation through $(q,p)$. The endpoints of the corrected arc for the intervening word $w_k$ therefore have the same displacement. After an integer translation, they are $A=(0,0)$ and $B=(q,p)$.

In the first and last retained triangles, the correction joins a lattice endpoint to the opposite edge without adding any intermediate edge crossings. The original straight line never visits the same triangle twice, so the correction creates neither a self-intersection nor consecutive crossings of the same triangulation edge. Hence $\widetilde{L_t}(w_k)\in\mathcal A_0(A,B)$.

Choose the original sign of each retained endpoint triangle in the endpoint rule. The corrected sign word is then obtained from the expanded word $w_k$ by deleting the edge-sign blocks at its ends. Its internal sign word is also a contiguous subword of the original internal word. Equation~\eqref{eq:gm-transfer} and $U_\pm\geq E_2$ therefore give
\[
|\widetilde{L_t}(w_k)|\leq N(w_k).
\]
If only one triangle is traversed, its GM length is $1$, and the same inequality holds.
We already know that $N(w_0)=(F_{S_0})_{21}=m_t$. Thus Theorems~\ref{thm:gm-distance} and~\ref{thm:length-is-number} imply
\[
N(w_0)=m_t=d(A,B)
\leq|\widetilde{L_t}(w_k)|\leq N(w_k).
\]
Consequently $N(w_0)$ is the least element of $N_t$.
Substitution into the formula for the Lagrange constant, together with the endpoint cases and the equality established at the beginning, proves
\[
\mathcal L(\alpha_{S(t)})
=\mathcal L(\alpha_{S^\ast(1/t)})
=\frac{\sqrt{((3+k_1+k_2+k_3)m_t-k_t)^2-4}}{m_t}
\]
for every reduced fraction $t\in[0,\infty]$. Since $\sigma$ was arbitrary, $\mathcal M_{k_1,k_2,k_3}\subset\mathcal L$ follows.
\end{proof}

The proof first restricts the possible values to
\[
\mathcal L(\alpha_{S(t)})
=\max\left\{\frac{\sqrt{(\operatorname{tr}(F_{S_i}))^2-(-1)^{n+1}\cdot4}}
{(F_{S_i})_{21}}\ \middle|\ 0\leq i\leq n\right\}.
\]
It then determines which cyclic shift has the least $(2,1)$-entry and identifies this entry. This step uses the minimality of the GM distance and the fact that its minimum is the corresponding GM number. The argument applies to general GM numbers.

\begin{exam}
Let $(k_1,k_2,k_3,\sigma)=(1,2,0,\textrm{id})$ and $t=\tfrac25$. Since $i_{\frac25}=1$,
\[
S\!\left(\tfrac25\right)=(5,1,3,3,1,5,4,1,3,4).
\]
The sequences $w$ that give $N(w)\in N_{\frac25}$ are
\begin{align*}
&(1,3,3,1,5,4,1,3,4),\ (3,3,1,5,4,1,3,4,5),\ (3,1,5,4,1,3,4,5,1),
\ (1,5,4,1,3,4,5,1,3),\\
&(5,4,1,3,4,5,1,3,3),\ (4,1,3,4,5,1,3,3,1),\ (1,3,4,5,1,3,3,1,5),\ (3,4,5,1,3,3,1,5,4),\\
&(4,5,1,3,3,1,5,4,1),\ (5,1,3,3,1,5,4,1,3).
\end{align*}
Their continuants are, respectively,
\[
8227,\ 32957,\ 12039,\ 12041,\ 32937,\ 8261,\ 9997,\ 31881,\ 12199,\ 11127.
\]
The least value is $8227$. The corresponding arcs are shown in Table~\ref{table2}. If an endpoint in an uncorrected diagram lies on an edge, only the signs displayed on that edge are included in $w$. In every case the corrected GM length is at most $N(w)$.
\[
F_{S(\frac25)}=
\begin{bmatrix}
N(5,1,3,3,1,5,4,1,3,4)&N(5,1,3,3,1,5,4,1,3)\\
N(1,3,3,1,5,4,1,3,4)&N(1,3,3,1,5,4,1,3)
\end{bmatrix}
=\begin{bmatrix}47431&11127\\8227&1930\end{bmatrix}.
\]
Hence $\alpha_{S(\frac25)}=(\sqrt{2436508317}+45501)/16454$ and
\[
Q_{S(\frac25)}=x^2-\frac{45501}{8227}xy-\frac{11127}{8227}y^2,
\]
and
\[
\mathcal L\!\left(\frac{\sqrt{2436508317}+45501}{16454}\right)
=\mathcal M\!\left(x^2-\frac{45501}{8227}xy-\frac{11127}{8227}y^2\right)
=\frac{\sqrt{2436508317}}{8227}.
\]
\begin{table}[htbp]
\centering
\begingroup
\pgfdeclarelayer{ec foreground}
\pgfsetlayers{main,ec foreground}
\tikzset{
  ec edge/.style={draw=black!65,line width=.34pt},
  ec arc/.style={draw=annred,line width=.65pt},
  ec sign/.style={text=annred,fill=white,inner sep=.12pt,
    font=\scriptsize,scale=.72,transform shape},
  ec point/.style={circle,fill=annred,inner sep=.95pt}
}
\newcommand{\ecTile}[5]{%
  \draw[ec edge] (#1,#2) rectangle ({#1+1},{#2+1});
  \draw[ec edge] (#1,{#2+1})--({#1+1},#2);
  \begin{pgfonlayer}{ec foreground}
  \if\relax\detokenize{#3}\relax\else
    \pgfmathsetmacro{\ecCross}{((#1+#2+1)-\ecIntercept)/1.4-#1}
    \ifdim\ecCross pt>.5pt
      \node[ec sign] at ({#1+.28},{#2+.72}) {$#3#3$};
    \else
      \node[ec sign] at ({#1+.70},{#2+.30}) {$#3#3$};
    \fi
  \fi
  \if\relax\detokenize{#4}\relax\else
    \node[ec sign] at ({#1+.80},{#2+.81}) {$#4$};
  \fi
  \if\relax\detokenize{#5}\relax\else
    \node[ec sign] at ({#1+.21},{#2+.25}) {$#5$};
  \fi
  \end{pgfonlayer}
}
\newcommand{\ecBeforeDiagonal}[5]{%
  \ifnum\ecPhase=0\relax
    \ecTile{#1}{#2}{#3}{#4}{#5}%
  \else
    \ecTile{#1}{#2}{}{#4}{#5}%
  \fi
}
\newcommand{\ecHorizontal}[3]{%
  \begin{pgfonlayer}{ec foreground}
  \def\ecGivenSign{#3}\def\ecPositiveSign{+}
  \ifx\ecGivenSign\ecPositiveSign
    \node[ec sign] at ({#1-.49},#2) {$#3$};
  \else
    \node[ec sign] at (#1,#2) {$#3$};
  \fi
  \end{pgfonlayer}
}
\newcommand{\ecBeforeHorizontal}[3]{%
  \ifnum\ecPhase=0\relax\ecHorizontal{#1}{#2}{#3}\fi
}
\newcommand{\ecUpperTriangle}[3]{%
  \draw[ec edge] (#1,{#2+1})--({#1+1},{#2+1})--({#1+1},#2)--cycle;
  \begin{pgfonlayer}{ec foreground}
  \node[ec sign] at ({#1+.79},{#2+.31}) {$#3$};
  \end{pgfonlayer}
}
\newcommand{\ecLowerTriangle}[3]{%
  \draw[ec edge] (#1,#2)--(#1,{#2+1})--({#1+1},#2)--cycle;
  \begin{pgfonlayer}{ec foreground}
  \node[ec sign] at ({#1+.21},{#2+.26}) {$#3$};
  \end{pgfonlayer}
}
\newcommand{\ecDiagram}[2]{%
  \begin{tikzpicture}[x=.65cm,y=.65cm,line cap=round,line join=round,
    baseline=(current bounding box.center)]
  \def\ecPhase{#2}%
  \def\ecBottom{0}\def\ecTop{3}
  \pgfmathsetmacro{\ecExitAngle}{atan(2/5)}
  \pgfmathsetmacro{\ecEntryAngle}{180+\ecExitAngle}
  \ifnum#1=1\relax\def\ecTop{2}\fi
  \ifnum#1=6\relax\def\ecTop{2}\fi
  \ifnum#1=10\relax\def\ecBottom{-1}\def\ecTop{2}\fi
  \path[use as bounding box] (-.08,{\ecBottom-.09}) rectangle (5.08,{\ecTop+.09});
  \ifcase#1\relax
  \or %
    \def\ecIntercept{0}
    \ecTile{0}{0}{-}{-}{+}\ecTile{1}{0}{+}{-}{+}\ecTile{2}{0}{+}{+}{+}
    \ecTile{2}{1}{-}{-}{-}\ecTile{3}{1}{-}{-}{+}\ecTile{4}{1}{+}{+}{+}
    \ecHorizontal{2.65}{1}{+}
    \ifnum#2=0\relax
      \draw[ec arc] (0,.20).. controls (.10,.20) and (.13,.10)..(.20,.08)
        --(2.20,.88)--(2.47,1)--(2.80,1.12)--(4.80,1.92)--(4.80,2);
      \node[ec point] at (0,.20) {};\node[ec point] at (4.80,2) {};
    \else
      \draw[ec arc] (0,0)--({5/7},{2/7})--(2.20,.88)--(2.47,1)
        --(2.80,1.12)--({30/7},{12/7})--(5,2);
      \node[ec point] at (0,0) {};\node[ec point] at (5,2) {};
    \fi
  \or %
    \def\ecIntercept{0}
    \ecBeforeDiagonal{0}{0}{-}{-}{}\ecTile{1}{0}{+}{-}{+}\ecTile{2}{0}{+}{+}{+}
    \ecTile{2}{1}{-}{-}{-}\ecTile{3}{1}{-}{-}{+}\ecTile{4}{1}{+}{+}{+}
    \ecTile{4}{2}{-}{-}{-}
    \ecHorizontal{2.65}{1}{+}\ecHorizontal{4.48}{2}{-}
    \ifnum#2=0\relax
      \draw[ec arc] ({5/7},{2/7})--(2.20,.88)--(2.47,1)--(2.80,1.12)
        --({5+.18*cos(\ecEntryAngle)},{2+.18*sin(\ecEntryAngle)})
        arc[start angle=\ecEntryAngle,end angle=90,radius=.18];
      \node[ec point] at ({5/7},{2/7}) {};\node[ec point] at (5,2.18) {};
    \else
      \draw[ec arc] (0,1)--(1,.40)--(2.20,.88)--(2.47,1)--(2.80,1.12)
        --({5+.18*cos(\ecEntryAngle)},{2+.18*sin(\ecEntryAngle)})
        arc[start angle=\ecEntryAngle,end angle=135,radius=.18]--(5,3);
      \node[ec point] at (0,1) {};\node[ec point] at (5,3) {};
    \fi
  \or %
    \def\ecIntercept{.4}
    \ecTile{0}{0}{+}{-}{+}\ecTile{1}{0}{+}{+}{+}\ecTile{1}{1}{-}{-}{-}
    \ecTile{2}{1}{-}{-}{+}\ecTile{3}{1}{+}{+}{+}\ecTile{3}{2}{-}{-}{-}
    \ecLowerTriangle{4}{2}{+}
    \ecHorizontal{1.65}{1}{+}\ecHorizontal{3.48}{2}{-}
    \ifnum#2=0\relax
      \draw[ec arc] (0,.40)--(1.20,.88)--(1.47,1)--(1.80,1.12)
        --({4+.18*cos(\ecEntryAngle)},{2+.18*sin(\ecEntryAngle)})
        arc[start angle=\ecEntryAngle,end angle=\ecExitAngle,radius=.18]
        --({33/7},{16/7});
      \node[ec point] at (0,.40) {};\node[ec point] at ({33/7},{16/7}) {};
    \else
      \draw[ec arc] (0,0)--({3/7},{4/7})--(1.20,.88)--(1.47,1)--(1.80,1.12)
        --({4+.18*cos(\ecEntryAngle)},{2+.18*sin(\ecEntryAngle)})
        arc[start angle=\ecEntryAngle,end angle=90,radius=.18]--(5,2);
      \node[ec point] at (0,0) {};\node[ec point] at (5,2) {};
    \fi
  \or %
    \def\ecIntercept{.4}
    \ecUpperTriangle{0}{0}{-}
    \ecTile{1}{0}{+}{+}{+}\ecTile{1}{1}{-}{-}{-}\ecTile{2}{1}{-}{-}{+}
    \ecTile{3}{1}{+}{+}{+}\ecTile{3}{2}{-}{-}{-}\ecTile{4}{2}{-}{-}{+}
    \ecHorizontal{1.65}{1}{+}\ecHorizontal{3.48}{2}{-}
    \ifnum#2=0\relax
      \draw[ec arc] ({3/7},{4/7})--(1.20,.88)--(1.47,1)--(1.80,1.12)
        --({4+.18*cos(\ecEntryAngle)},{2+.18*sin(\ecEntryAngle)})
        arc[start angle=\ecEntryAngle,end angle=\ecExitAngle,radius=.18]--(5,2.40);
      \node[ec point] at ({3/7},{4/7}) {};\node[ec point] at (5,2.40) {};
    \else
      \draw[ec arc] (0,1)--(1,.80)--(1.20,.88)--(1.47,1)--(1.80,1.12)
        --({4+.18*cos(\ecEntryAngle)},{2+.18*sin(\ecEntryAngle)})
        arc[start angle=\ecEntryAngle,end angle=\ecExitAngle,radius=.18]
        --({33/7},{16/7})--(5,3);
      \node[ec point] at (0,1) {};\node[ec point] at (5,3) {};
    \fi
  \or %
    \def\ecIntercept{.8}
    \ecTile{0}{0}{+}{+}{+}\ecTile{0}{1}{-}{-}{-}\ecTile{1}{1}{-}{-}{+}
    \ecTile{2}{1}{+}{+}{+}\ecTile{2}{2}{-}{-}{-}\ecTile{3}{2}{-}{-}{+}
    \ecBeforeDiagonal{4}{2}{+}{}{+}
    \ecHorizontal{.65}{1}{+}\ecHorizontal{2.48}{2}{-}
    \ifnum#2=0\relax
      \draw[ec arc] (0,.80)--(.20,.88)--(.47,1)--(.80,1.12)
        --({3+.18*cos(\ecEntryAngle)},{2+.18*sin(\ecEntryAngle)})
        arc[start angle=\ecEntryAngle,end angle=\ecExitAngle,radius=.18]--({31/7},{18/7});
      \node[ec point] at (0,.80) {};\node[ec point] at ({31/7},{18/7}) {};
    \else
      \draw[ec arc] (0,0)--({1/7},{6/7})--(.20,.88)--(.47,1)--(.80,1.12)
        --({3+.18*cos(\ecEntryAngle)},{2+.18*sin(\ecEntryAngle)})
        arc[start angle=\ecEntryAngle,end angle=\ecExitAngle,radius=.18]--(4,2.40)--(5,2);
      \node[ec point] at (0,0) {};\node[ec point] at (5,2) {};
    \fi
  \or %
    \def\ecIntercept{-.2}
    \ecTile{0}{0}{-}{-}{-}\ecTile{1}{0}{-}{-}{+}\ecTile{2}{0}{+}{+}{+}
    \ecTile{2}{1}{-}{-}{-}\ecTile{3}{1}{-}{-}{+}\ecTile{4}{1}{+}{-}{+}
    \ecHorizontal{2.48}{1}{-}
    \ifnum#2=0\relax
      \draw[ec arc] (.49,0)--({6/7},{1/7})--({3+.18*cos(\ecEntryAngle)},{1+.18*sin(\ecEntryAngle)})
        arc[start angle=\ecEntryAngle,end angle=\ecExitAngle,radius=.18]--(5,1.80);
      \node[ec point] at (.49,0) {};\node[ec point] at (5,1.80) {};
    \else
      \draw[ec arc] (0,0)--({6/7},{1/7})
        --({3+.18*cos(\ecEntryAngle)},{1+.18*sin(\ecEntryAngle)})
        arc[start angle=\ecEntryAngle,end angle=\ecExitAngle,radius=.18]
        --({31/7},{11/7})--(5,2);
      \node[ec point] at (0,0) {};\node[ec point] at (5,2) {};
    \fi
  \or %
    \def\ecIntercept{.2}
    \ecTile{0}{0}{-}{-}{+}\ecTile{1}{0}{+}{+}{+}\ecTile{1}{1}{-}{-}{-}
    \ecTile{2}{1}{-}{-}{+}\ecTile{3}{1}{+}{-}{+}\ecTile{4}{1}{+}{+}{+}
    \draw[ec edge] (4,2)--(4,3)--(5,2)--cycle;
    \ecHorizontal{1.48}{1}{-}\ecBeforeHorizontal{4.65}{2}{+}
    \ifnum#2=0\relax
      \draw[ec arc] (0,.20)--({2+.18*cos(\ecEntryAngle)},{1+.18*sin(\ecEntryAngle)})
        arc[start angle=\ecEntryAngle,end angle=\ecExitAngle,radius=.18]
        --(4.20,1.88)--(4.47,2);
      \node[ec point] at (0,.20) {};\node[ec point] at (4.47,2) {};
    \else
      \draw[ec arc] (0,0)--({4/7},{3/7})
        --({2+.18*cos(\ecEntryAngle)},{1+.18*sin(\ecEntryAngle)})
        arc[start angle=\ecEntryAngle,end angle=\ecExitAngle,radius=.18]
        --({29/7},{13/7})--(5,2);
      \node[ec point] at (0,0) {};\node[ec point] at (5,2) {};
    \fi
  \or %
    \def\ecIntercept{.2}
    \ecBeforeDiagonal{0}{0}{-}{-}{}\ecTile{1}{0}{+}{+}{+}\ecTile{1}{1}{-}{-}{-}
    \ecTile{2}{1}{-}{-}{+}\ecTile{3}{1}{+}{-}{+}\ecTile{4}{1}{+}{+}{+}
    \ecTile{4}{2}{-}{-}{-}
    \ecHorizontal{1.48}{1}{-}\ecHorizontal{4.65}{2}{+}
    \ifnum#2=0\relax
      \draw[ec arc] ({4/7},{3/7})--({2+.18*cos(\ecEntryAngle)},{1+.18*sin(\ecEntryAngle)})
        arc[start angle=\ecEntryAngle,end angle=\ecExitAngle,radius=.18]
        --(4.20,1.88)--(4.47,2)--(4.80,2.12)--(5,2.20);
      \node[ec point] at ({4/7},{3/7}) {};\node[ec point] at (5,2.20) {};
    \else
      \draw[ec arc] (0,1)--(1,.60)
        --({2+.18*cos(\ecEntryAngle)},{1+.18*sin(\ecEntryAngle)})
        arc[start angle=\ecEntryAngle,end angle=\ecExitAngle,radius=.18]
        --(4.20,1.88)--(4.47,2)--(4.80,2.12)--({34/7},{15/7})--(5,3);
      \node[ec point] at (0,1) {};\node[ec point] at (5,3) {};
    \fi
  \or %
    \def\ecIntercept{.6}
    \ecTile{0}{0}{+}{+}{+}\ecTile{0}{1}{-}{-}{-}\ecTile{1}{1}{-}{-}{+}
    \ecTile{2}{1}{+}{-}{+}\ecTile{3}{1}{+}{+}{+}\ecTile{3}{2}{-}{-}{-}
    \ecLowerTriangle{4}{2}{+}
    \ecHorizontal{.48}{1}{-}\ecHorizontal{3.65}{2}{+}
    \ifnum#2=0\relax
      \draw[ec arc] (0,.60)--({1+.18*cos(\ecEntryAngle)},{1+.18*sin(\ecEntryAngle)})
        arc[start angle=\ecEntryAngle,end angle=\ecExitAngle,radius=.18]
        --(3.20,1.88)--(3.47,2)--(3.80,2.12)--({32/7},{17/7});
      \node[ec point] at (0,.60) {};\node[ec point] at ({32/7},{17/7}) {};
    \else
      \draw[ec arc] (0,0)--({2/7},{5/7})
        --({1+.18*cos(\ecEntryAngle)},{1+.18*sin(\ecEntryAngle)})
        arc[start angle=\ecEntryAngle,end angle=\ecExitAngle,radius=.18]
        --(3.20,1.88)--(3.47,2)--(3.80,2.12)--(4,2.20)--(5,2);
      \node[ec point] at (0,0) {};\node[ec point] at (5,2) {};
    \fi
  \or %
    \def\ecIntercept{-.4}
    \draw[ec edge] (0,0)--(1,0)--(1,-1)--cycle;
    \ecTile{0}{0}{-}{-}{-}\ecTile{1}{0}{-}{-}{+}\ecTile{2}{0}{+}{-}{+}
    \ecTile{3}{0}{+}{+}{+}\ecTile{3}{1}{-}{-}{-}\ecTile{4}{1}{-}{-}{+}
    \ecBeforeHorizontal{.48}{0}{-}\ecHorizontal{3.65}{1}{+}
    \ifnum#2=0\relax
      \draw[ec arc] (.82,0)
        arc[start angle=180,end angle=\ecExitAngle,radius=.18]
        --(3.20,.88)--(3.47,1)--(3.80,1.12)--(5,1.60);
      \node[ec point] at (.82,0) {};
      \node[ec point] at (5,1.60) {};
    \else
      \draw[ec arc] (0,0)--({1+.18*cos(135)},{.18*sin(135)})
        arc[start angle=135,end angle=\ecExitAngle,radius=.18]
        --(3.20,.88)--(3.47,1)--(3.80,1.12)--({32/7},{10/7})--(5,2);
      \node[ec point] at (0,0) {};\node[ec point] at (5,2) {};
    \fi
  \fi
  \end{tikzpicture}%
}
\setlength{\tabcolsep}{4pt}
\renewcommand{\arraystretch}{1.02}
\small
\begin{tabular}{|>{\centering\arraybackslash}m{3.8cm}|>{\centering\arraybackslash}m{3.5cm}|>{\centering\arraybackslash}m{3.5cm}|>{\centering\arraybackslash}m{1.35cm}|}
$w$&$\overline{L_t}(w)$&$\widetilde{L_t}(w)$&$N(w)$\\\hline
$(1,3,3,1,5,4,1,3,4)$&\ecDiagram{1}{0}&\ecDiagram{1}{1}&$8227$\\\hline
$(3,3,1,5,4,1,3,4,5)$&\ecDiagram{2}{0}&\ecDiagram{2}{1}&$32957$\\\hline
$(3,1,5,4,1,3,4,5,1)$&\ecDiagram{3}{0}&\ecDiagram{3}{1}&$12039$\\\hline
$(1,5,4,1,3,4,5,1,3)$&\ecDiagram{4}{0}&\ecDiagram{4}{1}&$12041$\\\hline
$(5,4,1,3,4,5,1,3,3)$&\ecDiagram{5}{0}&\ecDiagram{5}{1}&$32937$\\\hline
$(4,1,3,4,5,1,3,3,1)$&\ecDiagram{6}{0}&\ecDiagram{6}{1}&$8261$\\\hline
$(1,3,4,5,1,3,3,1,5)$&\ecDiagram{7}{0}&\ecDiagram{7}{1}&$9997$\\\hline
$(3,4,5,1,3,3,1,5,4)$&\ecDiagram{8}{0}&\ecDiagram{8}{1}&$31881$\\\hline
$(4,5,1,3,3,1,5,4,1)$&\ecDiagram{9}{0}&\ecDiagram{9}{1}&$12199$\\\hline
$(5,1,3,3,1,5,4,1,3)$&\ecDiagram{10}{0}&\ecDiagram{10}{1}&$11127$\\\hline
\end{tabular}
\endgroup
\vspace{2mm}
\caption{Endpoint corrections and $N(w)$}\label{table2}
\end{table}

\end{exam}
\FloatBarrier
\begingroup
\setlength{\intextsep}{6pt}
\begin{exam}
We list several quadratic irrationals $\alpha$ and their Lagrange constants obtained from Theorem~\ref{thm:markov-value-gen}. Table~\ref{table000} treats $(k_1,k_2,k_3)=(0,0,0)$; when the three coefficients are equal, the values do not depend on $\sigma$. Tables~\ref{table001-1}, \ref{table001-2}, and~\ref{table001-3} treat $(0,0,1)$, and Tables~\ref{table011-1}, \ref{table011-2}, and~\ref{table011-3} treat $(0,1,1)$. Tables~\ref{table111} and~\ref{table222} treat $(1,1,1)$ and $(2,2,2)$, respectively. Table~\ref{table120-1} treats $(k_1,k_2,k_3,\sigma)=(1,2,0,\mathrm{id})$. In each table the rows are ordered by their GM numbers.
\begin{table}[H]
\centering
\resizebox{\linewidth}{!}{%
\begin{tabular}{|>{\centering\arraybackslash}m{1cm}|>{\centering\arraybackslash}m{5cm}|>{\centering\arraybackslash}m{4cm}|>{\centering\arraybackslash}m{1cm}|>{\centering\arraybackslash}m{3cm}|}
$t$&$S(t)$& $\alpha=[\overline{S(t)}]$ &$m_t$&$\mathcal L(\alpha)$\\
\hline
&&&&\\[-2mm]
$\dfrac{0}{1}$&$(1,1)$&$\dfrac{\sqrt{5}+1}{2}$&$1$&$\sqrt{5}$\\[4mm]
\hline
&&&&\\[-2mm]
$\dfrac{1}{1}$&$(2,2)$&${\sqrt{2}+1}$&$2$&$2\sqrt{2}$\\[4mm]
\hline
&&&&\\[-2mm]
$\dfrac{1}{2}$&$(2,1,1,2)$&$\dfrac{\sqrt{221}+11}{10}$&$5$&$\dfrac{\sqrt{221}}{5}$\\[4mm]
\hline
&&&&\\[-2mm]
$\dfrac{1}{3}$&$(2,1,1,1,1,2)$&$\dfrac{\sqrt{1517}+29}{26}$&$13$&$\dfrac{\sqrt{1517}}{13}$\\[4mm]
\hline
&&&&\\[-2mm]
$\dfrac{2}{3}$&$(2,1,1,2,2,2)$&$\dfrac{\sqrt{7565}+63}{58}$&$29$&$\dfrac{\sqrt{7565}}{29}$\\[4mm]
\hline
&&&&\\[-2mm]
$\dfrac{1}{4}$&$(2,1,1,1,1,1,1,2)$&$\dfrac{5\sqrt{26}+19}{17}$&$34$&$\dfrac{10\sqrt{26}}{17}$\\[4mm]
\hline
&&&&\\[-2mm]
$\dfrac{1}{5}$&$(2,1,1,1,1,1,1,1,1,2)$&$\dfrac{\sqrt{71285}+199}{178}$&$89$&$\dfrac{\sqrt{71285}}{89}$\\[4mm]
\hline
&&&&\\[-2mm]
$\dfrac{3}{4}$&$(2,1,1,2,2,2,2,2)$&$\dfrac{\sqrt{257045}+367}{338}$&$169$&$\dfrac{\sqrt{257045}}{169}$\\[4mm]
\hline
\end{tabular}}
\vspace{1mm}
\caption{$(k_1,k_2,k_3)=(0,0,0)$}\label{table000}
\end{table}
\begin{table}[H]
\centering
\resizebox{\linewidth}{!}{\begin{tabular}{|>{\centering\arraybackslash}m{1cm}|>{\centering\arraybackslash}m{5cm}|>{\centering\arraybackslash}m{4cm}|>{\centering\arraybackslash}m{1cm}|>{\centering\arraybackslash}m{3cm}|}
$t$&$S(t)$& $\alpha=[\overline{S(t)}]$ &$m_t$&$\mathcal L(\alpha)$\\
\hline
&&&&\\[-2mm]
$\dfrac{0}{1}$&$(2,1)$&$\sqrt{3}+1$&$1$&$2\sqrt{3}$\\[4mm]
\hline
&&&&\\[-2mm]
$\dfrac{1}{1}$&$(3,2)$&$\dfrac{\sqrt{15}+3}{2}$&$2$&$\sqrt{15}$\\[4mm]
\hline
&&&&\\[-2mm]
$\dfrac{1}{2}$&$(3,1,1,3)$&$\dfrac{5\sqrt{29}+23}{14}$&$7$&$\dfrac{5\sqrt{29}}{7}$\\[4mm]
\hline
&&&&\\[-2mm]
$\dfrac{1}{3}$&$(3,1,2,1,1,3)$&$\dfrac{7\sqrt{51}+43}{25}$&$25$&$\dfrac{14\sqrt{51}}{25}$\\[4mm]
\hline
&&&&\\[-2mm]
$\dfrac{2}{3}$&$(3,1,1,3,3,2)$&$\dfrac{\sqrt{11235}+83}{53}$&$53$&$\dfrac{2\sqrt{11235}}{53}$\\[4mm]
\hline
&&&&\\[-2mm]
$\dfrac{1}{4}$&$(3,1,2,1,1,2,1,3)$&$\dfrac{\sqrt{15293}+107}{62}$&$93$&$\dfrac{\sqrt{15293}}{31}$\\[4mm]
\hline
&&&&\\[-2mm]
$\dfrac{1}{5}$&$(3,1,2,1,2,1,1,2,1,3)$&$\dfrac{3\sqrt{53207}+599}{346}$&$346$&$\dfrac{3\sqrt{53207}}{173}$\\[4mm]
\hline
&&&&\\[-2mm]
$\dfrac{3}{4}$&$(3,1,1,3,2,3,3,2)$&$\dfrac{\sqrt{308765}+435}{278}$&$417$&$\dfrac{\sqrt{308765}}{139}$\\[4mm]
\hline
\end{tabular}}
\vspace{1mm}
\caption{ $(k_1,k_2,k_3,\sigma)=(0,0,1,\mathrm{id})$}\label{table001-1}
\end{table}
\newpage
\begin{table}[H]
\centering
\resizebox{\linewidth}{!}{%
\begin{tabular}{|>{\centering\arraybackslash}m{1cm}|>{\centering\arraybackslash}m{5cm}|>{\centering\arraybackslash}m{4cm}|>{\centering\arraybackslash}m{1cm}|>{\centering\arraybackslash}m{3cm}|}
$t$&$S(t)$& $\alpha=[\overline{S(t)}]$ &$m_t$&$\mathcal L(\alpha)$\\
\hline
&&&&\\[-2mm]
$\dfrac{0}{1}$&$(2,1)$&$\sqrt{3}+1$&$1$&$2\sqrt{3}$\\[4mm]
\hline
&&&&\\[-2mm]
$\dfrac{1}{1}$&$(3,3)$&$\dfrac{\sqrt{13}+3}{2}$&3&$\sqrt{13}$\\[4mm]
\hline
&&&&\\[-2mm]
$\dfrac{1}{2}$&$(3,1,2,3)$&$\dfrac{\sqrt{399}+17}{10}$&$10$&$\dfrac{\sqrt{399}}{5}$\\[4mm]
\hline
&&&&\\[-2mm]
$\dfrac{1}{3}$&$(3,1,2,2,1,3)$&$\dfrac{\sqrt{21605}+127}{74}$&$37$&$\dfrac{\sqrt{21605}}{37}$\\[4mm]
\hline
&&&&\\[-2mm]
$\dfrac{2}{3}$&$(3,1,2,3,3,3)$&$\dfrac{\sqrt{47523}+185}{109}$&$109$&$\dfrac{2\sqrt{47523}}{109}$\\[4mm]
\hline
&&&&\\[-2mm]
$\dfrac{1}{4}$&$(3,1,2,1,2,2,1,3)$&$\dfrac{5\sqrt{3003}+237}{137}$&$137$&$\dfrac{10\sqrt{3003}}{137}$\\[4mm]
\hline
&&&&\\[-2mm]
$\dfrac{1}{5}$&$(3,1,2,1,2,2,1,2,1,3)$&$\dfrac{\sqrt{4173845}+1769}{1022}$&$511$&$\dfrac{\sqrt{4173845}}{511}$\\[4mm]
\hline
&&&&\\[-2mm]
$\dfrac{3}{4}$&$(3,1,2,3,3,3,3,3)$&$\dfrac{\sqrt{5654883}+2018}{1189}$&$1189$&$\dfrac{2\sqrt{5654883}}{1189}$\\[4mm]
\hline
\end{tabular}}
\vspace{1mm}
\caption{ $(k_1,k_2,k_3,\sigma)=(0,0,1,(1\ 2 \ 3))$}\label{table001-2}
\end{table}
\begin{table}[H]
\centering
\resizebox{\linewidth}{!}{%
\begin{tabular}{|>{\centering\arraybackslash}m{1cm}|>{\centering\arraybackslash}m{5cm}|>{\centering\arraybackslash}m{4cm}|>{\centering\arraybackslash}m{1cm}|>{\centering\arraybackslash}m{3cm}|}
$t$&$S(t)$& $\alpha=[\overline{S(t)}]$ &$m_t$&$\mathcal L(\alpha)$\\
\hline
&&&&\\[-2mm]
$\dfrac{0}{1}$&$(1,1)$&$\dfrac{\sqrt{5}+1}{2}$&$1$&$\sqrt{5}$\\[4mm]
\hline
&&&&\\[-2mm]
$\dfrac{1}{1}$&$(3,2)$&$\dfrac{\sqrt{15}+3}{2}$&$2$&$\sqrt{15}$\\[4mm]
\hline
&&&&\\[-2mm]
$\dfrac{1}{2}$&$(3,1,1,2)$&$\dfrac{3\sqrt{11}+8}{5}$&$5$&$\dfrac{6\sqrt{11}}{5}$\\[4mm]
\hline
&&&&\\[-2mm]
$\dfrac{1}{3}$&$(3,1,1,1,1,2)$&$\dfrac{15\sqrt{3}+21}{13}$&$13$&$\dfrac{30\sqrt{3}}{13}$\\[4mm]
\hline
&&&&\\[-2mm]
$\dfrac{1}{4}$&$(3,1,1,1,1,1,1,2)$&$\dfrac{\sqrt{4623}+55}{34}$&$34$&$\dfrac{\sqrt{4623}}{17}$\\[4mm]
\hline
&&&&\\[-2mm]$\dfrac{2}{3}$&$(3,1,1,3,2,2)$&$\dfrac{\sqrt{2669}+41}{26}$&$39$&$\dfrac{\sqrt{2669}}{13}$\\[4mm]
\hline
&&&&\\[-2mm]
$\dfrac{1}{5}$&$(3,1,1,1,1,1,1,1,1,2)$&$\dfrac{\sqrt{31683}+144}{89}$&$89$&$\dfrac{2\sqrt{31683}}{89}$\\[4mm]
\hline
&&&&\\[-2mm]
$\dfrac{1}{6}$&$(3,1,1,1,1,1,1,1,1,1,1,2)$&$\dfrac{\sqrt{217155}+377}{233}$&$233$&$\dfrac{2\sqrt{217155}}{233}$\\[4mm]
\hline
\end{tabular}}
\vspace{1mm}
\caption{ $(k_1,k_2,k_3,\sigma)=(0,0,1,(1\ 3 \ 2))$}\label{table001-3}
\end{table}
\newpage
\begin{table}[H]
\centering
\resizebox{\linewidth}{!}{%
\begin{tabular}{|>{\centering\arraybackslash}m{1cm}|>{\centering\arraybackslash}m{5cm}|>{\centering\arraybackslash}m{4cm}|>{\centering\arraybackslash}m{1cm}|>{\centering\arraybackslash}m{3cm}|}
$t$&$S(t)$& $\alpha=[\overline{S(t)}]$ &$m_t$&$\mathcal L(\alpha)$\\
\hline
&&&&\\[-2mm]
$\dfrac{0}{1}$&$(3,1)$&$\dfrac{\sqrt{21}+3}{2}$&$1$&$\sqrt{21}$\\[4mm]
\hline
&&&&\\[-2mm]
$\dfrac{1}{1}$&$(4,3)$&$\dfrac{4\sqrt{3}+6}{3}$&$3$&$\dfrac{8\sqrt{3}}{3}$\\[4mm]
\hline
&&&&\\[-2mm]
$\dfrac{1}{2}$&$(4,1,2,4)$&$\dfrac{\sqrt{1023}+29}{13}$&$13$&$\dfrac{2\sqrt{1023}}{13}$\\[4mm]
\hline
&&&&\\[-2mm]
$\dfrac{1}{3}$&$(4,1,3,2,1,4)$&$\dfrac{3\sqrt{2567}+139}{61}$&$61$&$\dfrac{6\sqrt{2567}}{61}$\\[4mm]
\hline
&&&&\\[-2mm]
$\dfrac{2}{3}$&$(4,1,2,4,4,3)$&$\dfrac{\sqrt{49506}+195}{89}$&$178$&$\dfrac{2\sqrt{49506}}{89}$\\[4mm]
\hline
&&&&\\[-2mm]
$\dfrac{1}{4}$&$(4,1,3,1,2,3,1,4)$&$\dfrac{44\sqrt{273}+666}{291}$&$291$&$\dfrac{88\sqrt{273}}{291}$\\[4mm]
\hline
&&&&\\[-2mm]
$\dfrac{1}{5}$&$(4,1,3,1,3,2,1,3,1,4)$&$\dfrac{531\sqrt{43}+3191}{1393}$&$1393$&$\dfrac{1062\sqrt{43}}{1393}$\\[4mm]
\hline
&&&&\\[-2mm]
$\dfrac{3}{4}$&$(4,1,2,4,3,4,4,3)$&$\dfrac{2\sqrt{9600702}+5431}{2479}$&$2479$&$\dfrac{4\sqrt{9600702}}{2479}$\\[4mm]
\hline
\end{tabular}}
\vspace{1mm}
\caption{ $(k_1,k_2,k_3,\sigma)=(0,1,1,\mathrm{id})$}\label{table011-1}
\end{table}
\begin{table}[H]
\centering
\resizebox{\linewidth}{!}{\begin{tabular}{|>{\centering\arraybackslash}m{1cm}|>{\centering\arraybackslash}m{5cm}|>{\centering\arraybackslash}m{4cm}|>{\centering\arraybackslash}m{1cm}|>{\centering\arraybackslash}m{3cm}|}
$t$&$S(t)$& $\alpha=[\overline
{S(t)}]$ &$m_t$&$\mathcal L(\alpha)$\\
\hline
&&&&\\[-2mm]
$\dfrac{0}{1}$&$(2,1)$&$\sqrt{3}+1$&$1$&$2\sqrt{3}$\\[4mm]
\hline
&&&&\\[-2mm]
$\dfrac{1}{1}$&$(4,3)$&$\dfrac{4\sqrt{3}+6}{3}$&$3$&$\dfrac{8\sqrt{3}}{3}$\\[4mm]
\hline
&&&&\\[-2mm]
$\dfrac{1}{2}$&$(4,1,2,3)$&$\dfrac{2\sqrt{39}+11}{5}$&$10$&$\dfrac{4\sqrt{39}}{5}$\\[4mm]
\hline
&&&&\\[-2mm]
$\dfrac{1}{3}$&$(4,1,2,2,1,3)$&$\dfrac{\sqrt{8463}+82}{37}$&$37$&$\dfrac{2\sqrt{8463}}{37}$\\[4mm]
\hline
&&&&\\[-2mm]
$\dfrac{1}{4}$&$(4,1,2,1,2,2,1,3)$&$\dfrac{\sqrt{469221}+611}{274}$&$137$&$\dfrac{\sqrt{469221}}{137}$\\[4mm]
\hline
&&&&\\[-2mm]
$\dfrac{2}{3}$&$(4,1,2,4,3,3)$&$\dfrac{2\sqrt{30102}+305}{139}$&$139$&$\dfrac{4\sqrt{30102}}{139}$\\[4mm]
\hline
&&&&\\[-2mm]
$\dfrac{1}{5}$&$(4,1,2,1,2,2,1,2,1,3)$&$\dfrac{6\sqrt{45298}+1140}{511}$&$511$&$\dfrac{12\sqrt{45298}}{511}$\\[4mm]
\hline
&&&&\\[-2mm]
$\dfrac{2}{5}$&$(4,1,2,2,1,4,3,1,2,3)$&$\dfrac{22\sqrt{43662}+4050}{1839}$&$1839$&$\dfrac{44\sqrt{43662}}{1839}$\\[4mm]
\hline
\end{tabular}}
\vspace{1mm}
\caption{ $(k_1,k_2,k_3,\sigma)=(0,1,1,(1\ 2 \ 3))$}\label{table011-2}
\end{table}
\newpage
\begin{table}[H]
\centering
\resizebox{\linewidth}{!}{\begin{tabular}{|>{\centering\arraybackslash}m{1cm}|>{\centering\arraybackslash}m{5cm}|>{\centering\arraybackslash}m{4cm}|>{\centering\arraybackslash}m{1cm}|>{\centering\arraybackslash}m{3cm}|}
$t$&$S(t)$& $\alpha=[\overline{S(t)}]$ &$m_t$&$\mathcal L(\alpha)$\\
\hline
&&&&\\[-2mm]
$\dfrac{0}{1}$&$(2,1)$&$\sqrt{3}+1$&$1$&$2\sqrt{3}$\\[4mm]
\hline
&&&&\\[-2mm]
$\dfrac{1}{1}$&$(4,2)$&$\sqrt6+2$&$2$&$2\sqrt{6}$\\[4mm]
\hline
&&&&\\[-2mm]
$\dfrac{1}{2}$&$(4,1,1,3)$&$\dfrac{12\sqrt{2}+15}{7}$&$7$&$\dfrac{24\sqrt{2}}{7}$\\[4mm]
\hline
&&&&\\[-2mm]
$\dfrac{1}{3}$&$(4,1,2,1,1,3)$&$\dfrac{\sqrt{15621}+111}{50}$&$25$&$\dfrac{\sqrt{15621}}{25}$\\[4mm]
\hline
&&&&\\[-2mm]
$\dfrac{2}{3}$&$(4,1,1,4,3,2)$&$\dfrac{4\sqrt{1743}+138}{67}$&$67$&$\dfrac{8\sqrt{1743}}{67}$\\[4mm] 
\hline
&&&&\\[-2mm]
$\dfrac{1}{4}$&$(4,1,2,1,1,2,1,3)$&$\dfrac{\sqrt{53823}+207}{93}$&$93$&$\dfrac{2\sqrt{53823}}{93}$\\[4mm]
\hline
&&&&\\[-2mm]
$\dfrac{1}{5}$&$(4,1,2,1,2,1,1,2,1,3)$&$\dfrac{12\sqrt{1299}+386}{173}$&$346$&$\dfrac{24\sqrt{1299}}{173}$\\[4mm]
\hline
&&&&\\[-2mm]
$\dfrac{3}{4}$&$(4,1,1,4,2,3,4,2)$&$\dfrac{\sqrt{2729103}+1356}{661}$&$661$&$\dfrac{2\sqrt{2729103}}{661}$\\[4mm]
\hline
\end{tabular}}
\vspace{1mm}
\caption{ $(k_1,k_2,k_3,\sigma)=(0,1,1,(1\ 3 \ 2))$}\label{table011-3}
\end{table}
\begin{table}[H]
\centering
\resizebox{\linewidth}{!}{\begin{tabular}{|>{\centering\arraybackslash}m{1cm}|>{\centering\arraybackslash}m{5cm}|>{\centering\arraybackslash}m{4cm}|>{\centering\arraybackslash}m{1cm}|>{\centering\arraybackslash}m{3cm}|}
$t$&$S(t)$& $\alpha=[\overline{S(t)}]$ &$m_t$&$\mathcal L(\alpha)$\\
\hline
&&&&\\[-2mm]
$\dfrac{0}{1}$&$(3,1)$&$\dfrac{\sqrt{21}+3}{2}$&$1$&$\sqrt{21}$\\[4mm]
\hline
&&&&\\[-2mm]
$\dfrac{1}{1}$&$(5,3)$&$\dfrac{\sqrt{285}+15}{6}$&$3$&$\dfrac{\sqrt{285}}{3}$\\[4mm]
\hline
&&&&\\[-2mm]
$\dfrac{1}{2}$&$(5,1,2,4)$&$\dfrac{5\sqrt{237}+71}{26}$&$13$&$\dfrac{5\sqrt{237}}{13}$\\[4mm]
\hline
&&&&\\[-2mm]
$\dfrac{1}{3}$&$(5,1,3,2,1,4)$&$\dfrac{11\sqrt{1101}+339}{122}$&$61$&$\dfrac{11\sqrt{1101}}{61}$\\[4mm]
\hline
&&&&\\[-2mm]
$\dfrac{2}{3}$&$(5,1,2,5,4,3)$&$\dfrac{\sqrt{1692597}+1167}{434}$&$217$&$\dfrac{\sqrt{1692597}}{217}$\\[4mm]
\hline
&&&&\\[-2mm]
$\dfrac{1}{4}$&$(5,1,3,1,2,3,1,4)$&$\dfrac{\sqrt{3045021}+1623}{582}$&$291$&$\dfrac{\sqrt{3045021}}{291}$\\[4mm]
\hline
&&&&\\[-2mm]
$\dfrac{1}{5}$&$(5,1,3,1,3,2,1,3,1,4)$&$\dfrac{\sqrt{69839445}+7775}{2786}$&$1393$&$\dfrac{\sqrt{69839445}}{1393}$\\[4mm]
\hline
&&&&\\[-2mm]
$\dfrac{3}{4}$&$(5,1,2,5,3,4,5,3)$&$\dfrac{\sqrt{485629365}+19735}{7346}$&$3673$&$\dfrac{\sqrt{485629365}}{3673}$\\[4mm]
\hline
\end{tabular}}
\vspace{1mm}
\caption{ $(k_1,k_2,k_3)=(1,1,1)$}\label{table111}
\end{table}
\newpage
\begin{table}[H]
\centering
\resizebox{\linewidth}{!}{\begin{tabular}{|>{\centering\arraybackslash}m{1cm}|>{\centering\arraybackslash}m{5cm}|>{\centering\arraybackslash}m{4cm}|>{\centering\arraybackslash}m{1cm}|>{\centering\arraybackslash}m{3cm}|}
$t$&$S(t)$& $\alpha=[\overline{S(t)}]$ &$m_t$&$\mathcal L(\alpha)$\\
\hline
&&&&\\[-2mm]
$\dfrac{0}{1}$&$(5,1)$&$\dfrac{3\sqrt{5}+5}{2}$&$1$&$3\sqrt{5}$\\[4mm]
\hline
&&&&\\[-2mm]
$\dfrac{1}{1}$&$(8,4)$&$3\sqrt{2}+4$&$4$&$6\sqrt{2}$\\[5mm]
\hline
&&&&\\[-2mm]
$\dfrac{1}{2}$&$(8,1,3,6)$&$\dfrac{3\sqrt{221}+43}{10}$&$25$&$\dfrac{3\sqrt{221}}{5}$\\[4mm]
\hline
&&&&\\[-2mm]
$\dfrac{1}{3}$&$(8,1,5,3,1,6)$&$\dfrac{3\sqrt{1517}+113}{26}$&$169$&$\dfrac{3\sqrt{1517}}{13}$\\[4mm]
\hline
&&&&\\[-2mm]
$\dfrac{2}{3}$&$(8,1,3,8,6,4)$&$\dfrac{3\sqrt{7565}+247}{58}$&$841$&$\dfrac{3\sqrt{7565}}{29}$\\[4mm]
\hline
&&&&\\[-2mm]
$\dfrac{1}{4}$&$(8,1,5,1,3,5,1,6)$&$\dfrac{15\sqrt{26}+74}{17}$&$1156$&$\dfrac{30\sqrt{26}}{17}$\\[4mm]
\hline
&&&&\\[-2mm]
$\dfrac{1}{5}$&$(8,1,5,1,5,3,1,5,1,6)$&$\dfrac{3\sqrt{71285}+775}{178}$&$7921$&$\dfrac{3\sqrt{71285}}{89}$\\[4mm]
\hline
&&&&\\[-2mm]
$\dfrac{3}{4}$&$(8,1,3,8,4,6,8,4)$&$\dfrac{3\sqrt{257045}+1439}{338}$&$28561$&$\dfrac{3\sqrt{257045}}{169}$\\[4mm]
\hline
\end{tabular}}
\vspace{1mm}
\caption{ $(k_1,k_2,k_3)=(2,2,2)$}\label{table222}
\end{table}
\begin{table}[H]
\centering
\resizebox{\linewidth}{!}{\begin{tabular}{|>{\centering\arraybackslash}m{1cm}|>{\centering\arraybackslash}m{5cm}|>{\centering\arraybackslash}m{4cm}|>{\centering\arraybackslash}m{1cm}|>{\centering\arraybackslash}m{3cm}|}
$t$&$S(t)$& $\alpha=[\overline{S(t)}]$ &$m_t$&$\mathcal L(\alpha)$\\
\hline
&&&&\\[-2mm]
$\dfrac{0}{1}$&$(3,1)$&$\dfrac{\sqrt{21}+3}{2}$&$1$&$\sqrt{21}$\\[4mm]
\hline
&&&&\\[-2mm]
$\dfrac{1}{1}$&$(5,4)$&$\dfrac{\sqrt{30}+5}{2}$&$4$&$\sqrt{30}$\\[5mm]
\hline
&&&&\\[-2mm]
$\dfrac{1}{2}$&$(5,1,3,4)$&$\dfrac{10\sqrt{26}+47}{17}$&$17$&$\dfrac{20\sqrt{26}}{17}$\\[4mm]
\hline
&&&&\\[-2mm]
$\dfrac{1}{3}$&$(5,1,3,3,1,4)$&$\dfrac{\sqrt{723}+25}{9}$&$81$&$\dfrac{2\sqrt{723}}{9}$\\[4mm]
\hline
&&&&\\[-2mm]
$\dfrac{2}{3}$&$(5,1,3,5,4,4)$&$\dfrac{\sqrt{5004165}+2061}{746}$&$373$&$\dfrac{\sqrt{5004165}}{373}$\\[4mm]
\hline
&&&&\\[-2mm]
$\dfrac{1}{4}$&$(5,1,3,1,3,3,1,4)$&$\dfrac{\sqrt{1340963}+1077}{386}$&$386$&$\dfrac{\sqrt{1340963}}{193}$\\[4mm]
\hline
&&&&\\[-2mm]
$\dfrac{1}{5}$&$(5,1,3,1,3,3,1,3,1,4)$&$\dfrac{\sqrt{16635}+120}{43}$&$1849$&$\dfrac{2\sqrt{16635}}{43}$\\[4mm]
\hline
&&&&\\[-2mm]
$\dfrac{3}{4}$&$(5,1,3,5,4,4,5,4)$&$\dfrac{2\sqrt{150737006}+22602}{8185}$&$8185$&$\dfrac{4\sqrt{150737006}}{8185}$\\[4mm]
\hline
\end{tabular}}
\vspace{1mm}
\caption{ $(k_1,k_2,k_3,\sigma)=(1,2,0,\mathrm{id})$}\label{table120-1}
\end{table}
\end{exam}
\endgroup
\clearpage

\section{Mechanical Words}\label{sec:mechanical-words}
In this section we introduce the mechanical words needed for the proof of Markov's theorem in the next section. Since our goal is to relate them to strongly admissible sequences, we define them using slopes at least $1$, rather than the more usual convention in which the slope is at most $1$.

\begin{defi}\label{def:mechanical-word}\index{Mechanical word!slope}\index{Mechanical word!intercept}\index{Mechanical word!left}\index{Mechanical word}\index{Mechanical word!right}
Let $t\in[1,\infty]$ and $\theta\in\mathbb R$. First suppose that $1\le t<\infty$. Orient the line
$\ell_{t,\theta}:y=tx+\theta$ in the direction in which both coordinates increase, that is, from lower left to upper right. For each $n\in\mathbb Z$, let
$P_n=((n-\theta)/t,n)$ be the intersection of $\ell_{t,\theta}$ with the horizontal line $y=n$. Let
\[
    r_n=\left\lceil \frac{n-\theta}{t}\right\rceil
\]
be the $x$-coordinate of the nearest lattice point on or to the right of $P_n$, and let
\[
    l_n=\left\lfloor \frac{n-\theta}{t}\right\rfloor
\]
be the $x$-coordinate of the nearest lattice point on or to the left of $P_n$. If $P_n$ itself is a lattice point, that lattice point is regarded as belonging to both the right and the left side. Put
\[
    \varepsilon^R_n=r_{n+1}-r_n,
    \qquad
    \varepsilon^L_n=l_{n+1}-l_n.
\]
Since $t\ge1$, we have $\varepsilon^R_n,\varepsilon^L_n\in\{0,1\}$.

The bi-infinite word $\mathbf b^R=(b^R_n)_{n\in\mathbb Z}$ defined by
\[
    b^R_n=
    \begin{cases}
    X & \text{if }\varepsilon^R_n=1,\\
    Y & \text{if }\varepsilon^R_n=0
    \end{cases}
\]
is called the \emph{right mechanical word} of slope $t$ and intercept $\theta$.
Similarly, the bi-infinite word $\mathbf b^L=(b^L_n)_{n\in\mathbb Z}$ defined by
\[
    b^L_n=
    \begin{cases}
    X & \text{if }\varepsilon^L_n=1,\\
    Y & \text{if }\varepsilon^L_n=0
    \end{cases}
\]
is called the \emph{left mechanical word} of slope $t$ and intercept $\theta$. Right and left mechanical words are collectively called \emph{mechanical words}.

Finally, for $t=\infty$, both the right and the left mechanical word are defined to be $\cdots YYY\cdots$, independently of the intercept.
\end{defi}

\begin{exam}
For the line $\ell_{\frac52,\frac14}:y=\frac52x+\frac14$, the right mechanical word is the bi-infinite purely periodic word with period $XYXYY$. The left mechanical word is the same word; see the left panel of Figure~\ref{fig:mechanical}. In this example the left and right mechanical words coincide, but the situation changes when the line passes through a lattice point. For $\ell_{\frac52,0}:y=\frac52x$, the right mechanical word again has period $XYXYY$, while the left mechanical word has period $YYXYX$; see the middle panel of Figure~\ref{fig:mechanical}. These two words agree after a shift, but the example shows where the distinction between the two conventions comes from. If the slope is irrational and the line passes through a lattice point, then the left and right mechanical words differ only around the unique lattice point through which the line passes; see the right panel of Figure~\ref{fig:mechanical}.
All panels show finite portions of the words, with the letters ordered from bottom to top.
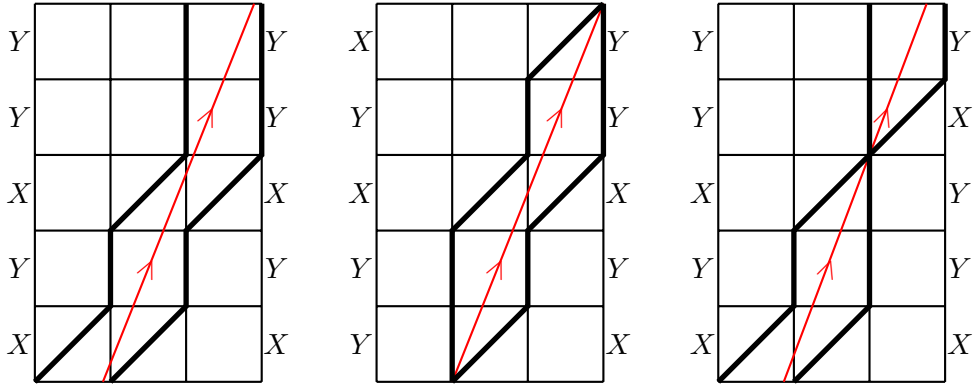
\begin{figure}[ht]
\centering
    \begin{tikzpicture}
    \draw[thick] (-1,0) grid (2,5);
    \draw[red,thick] (-0.1,0) -- (1.9,5);
    \draw[line width=2pt] (0,0) -- (1,1)-- (1,2) -- (2,3) -- (2,5);
    \draw[line width=2pt] (-1,0) -- (0,1)-- (0,2) -- (1,3) -- (1,5);
    \node at (2.2,0.5) {$X$};
    \node at (2.2,1.5) {$Y$};
    \node at (2.2,2.5) {$X$};
    \node at (2.2,3.5) {$Y$};
    \node at (2.2,4.5) {$Y$};
     \node at (-1.2,0.5) {$X$};
    \node at (-1.2,1.5) {$Y$};
    \node at (-1.2,2.5) {$X$};
    \node at (-1.2,3.5) {$Y$};
    \node at (-1.2,4.5) {$Y$};
    \node[red] at (0.5,1.5) {\rotatebox{65}{$>$}};
    \node[red] at (1.3,3.5) {\rotatebox{65}{$>$}};
\end{tikzpicture}
\quad 
\begin{tikzpicture}
    \draw[thick] (-1,0) grid (2,5);
    \draw[red,thick] (0,0) -- (2,5);
    \draw[line width=2pt] (0,0) -- (1,1)-- (1,2) -- (2,3) -- (2,5);
    \draw[line width=2pt] (0,0) -- (0,1)-- (0,2) -- (1,3) --(1,4)-- (2,5);
    \node at (2.2,0.5) {$X$};
    \node at (2.2,1.5) {$Y$};
    \node at (2.2,2.5) {$X$};
    \node at (2.2,3.5) {$Y$};
    \node at (2.2,4.5) {$Y$};
     \node at (-1.2,0.5) {$Y$};
    \node at (-1.2,1.5) {$Y$};
    \node at (-1.2,2.5) {$X$};
    \node at (-1.2,3.5) {$Y$};
    \node at (-1.2,4.5) {$X$};
    \node[red] at (0.6,1.5) {\rotatebox{65}{$>$}};
    \node[red] at (1.4,3.5) {\rotatebox{65}{$>$}};
\end{tikzpicture}
\quad 
\begin{tikzpicture}
    \draw[thick] (-1,0) grid (2,5);
    \draw[red,thick] ({1-3/sqrt(7)},0) -- ({1+2/sqrt(7)},5);
    \draw[line width=2pt] (0,0) -- (1,1)-- (1,2) -- (1,3) --(2,4)-- (2,5);
    \draw[line width=2pt] (-1,0) -- (0,1)-- (0,2) -- (1,3) --(1,4)-- (1,5);
    \node at (2.2,0.5) {$X$};
    \node at (2.2,1.5) {$Y$};
    \node at (2.2,2.5) {$Y$};
    \node at (2.2,3.5) {$X$};
    \node at (2.2,4.5) {$Y$};
     \node at (-1.2,0.5) {$X$};
    \node at (-1.2,1.5) {$Y$};
    \node at (-1.2,2.5) {$X$};
    \node at (-1.2,3.5) {$Y$};
    \node at (-1.2,4.5) {$Y$};
    \node[red] at (0.45,1.5) {\rotatebox{70}{$>$}};
    \node[red] at (1.2,3.5) {\rotatebox{70}{$>$}};
\end{tikzpicture}
    \caption{Examples of mechanical words}
    \label{fig:mechanical}
\end{figure}
\end{exam}

We record several elementary properties of mechanical words.
\begin{lemm}\label{rem:mechanical-slope-period}\index{Mechanical word!periodicity}
Let $t=\frac pq\in[1,\infty)$ be a reduced fraction. Then every mechanical word of slope $t$ has period $p$, and one period contains exactly $q$ occurrences of $X$. Thus, if $w$ is one period, then
\[
    \frac{|w|}{|w|_X}=\frac{p}{q}=t.
\]
\end{lemm}

\begin{proof}
It is enough to prove the assertion for right mechanical words; the proof for left mechanical words is identical. We have
\[
    r_n=\left\lceil \frac{n-\theta}{t}\right\rceil
    =
    \left\lceil \frac{q(n-\theta)}{p}\right\rceil.
\]
Since $t=\frac pq$,
\[
    r_{n+p}
    =
    \left\lceil \frac{q(n+p-\theta)}{p}\right\rceil
    =
    \left\lceil \frac{q(n-\theta)}{p}+q\right\rceil
    =
    r_n+q.
\]
Therefore
\[
    \varepsilon^R_{n+p}
    =
    r_{n+p+1}-r_{n+p}
    =
    r_{n+1}-r_n
    =
    \varepsilon^R_n,
\]
so $p$ is a period. Moreover, the number of $X$'s in one period is
\[
    \sum_{n=0}^{p-1}\varepsilon^R_n
    =
    r_p-r_0
    =
    q.
\]
Hence $|w|/|w|_X=p/q=t$. The same argument, with the floor function in place of the ceiling function, proves the assertion for left mechanical words.
\end{proof}

Next we relate shifts of mechanical words to changes of intercept.
\begin{lemm}\label{lem:shift-changes-intercept}\index{Mechanical word!intercept}\index{Shift of a word}
Let $k\in\mathbb Z$, and define the shift of a word by
\[
    (T^k\mathbf b)_n=b_{n+k}.
\]
Then, for every $t\in[1,\infty]$ and every intercept $\theta$,
\[
    T^k\mathbf b^R(t,\theta)=\mathbf b^R(t,\theta-k),
    \qquad
    T^k\mathbf b^L(t,\theta)=\mathbf b^L(t,\theta-k).
\]
In particular, a shift of a mechanical word is again a mechanical word of the same slope.
\end{lemm}

\begin{proof}
The case $t=\infty$ is clear because the word is $\cdots YYY\cdots$. Assume $1\le t<\infty$.
For right mechanical words, put
\[
    r_n(\theta)=\left\lceil\frac{n-\theta}{t}\right\rceil.
\]
Then
\[
    r_n(\theta-k)
    =
    \left\lceil\frac{n-(\theta-k)}{t}\right\rceil
    =
    \left\lceil\frac{n+k-\theta}{t}\right\rceil
    =
    r_{n+k}(\theta).
\]
Therefore
\[
    r_{n+1}(\theta-k)-r_n(\theta-k)
    =
    r_{n+k+1}(\theta)-r_{n+k}(\theta),
\]
and hence $T^k\mathbf b^R(t,\theta)=\mathbf b^R(t,\theta-k)$.
The proof for left mechanical words is identical, using
$l_n(\theta)=\lfloor(n-\theta)/t\rfloor$.
\end{proof}

For rational slopes, left mechanical words can also be represented as right mechanical words.
\begin{lemm}\label{lem:rational-left-is-right}\index{Mechanical word!left}\index{Mechanical word!right}
Assume that $t\in([1,\infty)\cap\mathbb Q)\cup\{\infty\}$. Then every mechanical word of slope $t$ can be written as a right mechanical word of the same slope $t$.
\end{lemm}

\begin{proof}
There is nothing to prove if the word is already a right mechanical word. Let $\mathbf b=\mathbf b^L(t,\theta)$ be a left mechanical word of slope $t$.
If $t=\infty$, then $\mathbf b=\cdots YYY\cdots$, which is also a right mechanical word. Assume $1\le t<\infty$, and write $t=p/q$ in lowest terms, with $p,q\in\mathbb Z_{>0}$ and $p\ge q$. The left mechanical word is determined by the differences of
\[
    l_n=\left\lfloor \frac{n-\theta}{t}\right\rfloor
    =
    \left\lfloor \frac{q(n-\theta)}{p}\right\rfloor.
\]
Put $x_n=(n-\theta)/t$. Since $t=p/q$, we have $x_{n+p}=x_n+q$. Hence the fractional parts $\{x_n\}$ are periodic in $n$ with period $p$, and the set
\[
    \{\{x_n\}\mid n\in\mathbb Z\}
\]
is finite. Choose $\delta>0$ sufficiently small so that
\[
    \{x_n\}+\delta<1
\]
for all $n$ with $\{x_n\}\ne0$. Such a $\delta$ exists because the set of fractional parts is finite; for example one may take
\[
    0<\delta<
    \min\{\,1-\{x_n\}\mid n=0,\dots,p-1,\ \{x_n\}\ne0\,\},
\]
with arbitrary $0<\delta<1$ if the set on the right is empty. Then, for all $n\in\mathbb Z$,
\[
    \left\lceil x_n+\delta\right\rceil
    =
    \left\lfloor x_n\right\rfloor+1.
\]
Indeed, this is immediate if $x_n$ is an integer, and otherwise follows from $\{x_n\}+\delta<1$. Put $\theta'=\theta-t\delta$. Then
\[
    \frac{n-\theta'}{t}=\frac{n-\theta}{t}+\delta=x_n+\delta.
\]
For the sequence
\[
    r_n=\left\lceil \frac{n-\theta'}{t}\right\rceil
\]
defining the right mechanical word $\mathbf b^R(t,\theta')$, we have $r_n=l_n+1$ for every $n$. Taking differences gives $r_{n+1}-r_n=l_{n+1}-l_n$. Thus $\mathbf b^L(t,\theta)$ and $\mathbf b^R(t,\theta')$ give the same letter at every position, i.e.
$\mathbf b^L(t,\theta)=\mathbf b^R(t,\theta')$.
\end{proof}

\begin{lemm}\label{lem:intercept-rational}\index{Mechanical word!intercept}\index{Shift of a word}
Let $t\in([1,\infty)\cap\mathbb Q)\cup\{\infty\}$. Then the mechanical word of slope $t$ is uniquely determined up to shift.
\end{lemm}

\begin{proof}
For $t=\infty$ the only word is $\cdots YYY\cdots$, so the assertion is clear. Write $t=p/q$ in lowest terms, with $p,q\in\mathbb Z_{>0}$ and $p\ge q$. By Lemma~\ref{lem:rational-left-is-right}, it suffices to consider right mechanical words.
The right mechanical word $\mathbf b^R(t,\theta)$ is determined by the differences of
\[
    r_n=\left\lceil \frac{n-\theta}{t}\right\rceil.
\]
Put $\alpha=q/p$ and $\beta=-q\theta/p$. Then $r_n=\lceil n\alpha+\beta\rceil$, and
\[
    b_n=X
    \Longleftrightarrow
    \lceil (n+1)\alpha+\beta\rceil-\lceil n\alpha+\beta\rceil=1.
\]
The right-hand side is unchanged when $\beta$ is increased by $1$, so we regard $\beta$ as a point of $\mathbb R/\mathbb Z$. Put
\[
    u_\beta(n)
    =
    \lceil (n+1)\alpha+\beta\rceil-
    \lceil n\alpha+\beta\rceil.
\]
Then
\[
    u_{\beta+\alpha}(n)
    =
    \lceil (n+2)\alpha+\beta\rceil-
    \lceil (n+1)\alpha+\beta\rceil
    =
    u_\beta(n+1).
\]
Thus increasing $\beta$ by $\alpha$ corresponds to shifting the word by one position.

Since $\alpha=q/p$ and $\gcd(p,q)=1$, the points
\[
    0,\alpha,2\alpha,\ldots,(p-1)\alpha
\]
in $\mathbb R/\mathbb Z$ are just a permutation of
\[
    0,\frac1p,\frac2p,\ldots,\frac{p-1}{p}.
\]
A right mechanical word is constant on each half-open interval $((j-1)/p,j/p]$, $j=1,\ldots,p$, on $\mathbb R/\mathbb Z$. The map $\beta\mapsto\beta+\alpha$ cyclically permutes these intervals. Therefore all words obtained by changing the intercept agree up to shift.
\end{proof}

\begin{prop}\label{prop:rational-slopes-distinct}\index{Mechanical word!slope}\index{Reversal of a word}
Let $t,t'\in([1,\infty)\cap\mathbb Q)\cup\{\infty\}$ with $t\ne t'$. Let $\mathbf b$ be a mechanical word of slope $t$, and let $\mathbf b'$ be a mechanical word of slope $t'$. Then $\mathbf b$ and $\mathbf b'$ are not shift-equivalent, and they are not reversals of each other up to shift.
\end{prop}
\begin{proof}
For a periodic word $\mathbf c$, let $D_X(\mathbf c)$ denote the proportion of $X$'s in one period. This is independent of the choice of period, and is unchanged by shifts and reversal.

The mechanical word of slope $t=\infty$ is $\cdots YYY\cdots$, so $D_X(\mathbf b)=0$ in this case. If $t=p/q\in[1,\infty)\cap\mathbb Q$ is written in lowest terms, Lemma~\ref{rem:mechanical-slope-period} gives $D_X(\mathbf b)=q/p=1/t$. Thus, in general, $D_X(\mathbf b)=1/t$, with the convention $1/\infty=0$.

If $\mathbf b$ and $\mathbf b'$ were shift-equivalent, then $D_X(\mathbf b)=D_X(\mathbf b')$, and hence $1/t=1/t'$, so $t=t'$, contradicting the assumption. Thus they are not shift-equivalent.

Similarly, if $\mathbf b$ and $\mathbf b'$ were reversals of each other up to shift, that is, if $\mathbf b=T^k((\mathbf b')^{\ast})$ for some $k\in\mathbb Z$, then reversal and shift would again preserve $D_X$. Hence $D_X(\mathbf b)=D_X(\mathbf b')$, forcing $t=t'$, again a contradiction.
\end{proof}

Finally, in this section, we relate rational-slope mechanical words to strongly admissible sequences in the case $(k_1,k_2,k_3)=(0,0,0)$. This relation is the key point in the proof of Markov's theorem.

\index{Substitution of words}Define the substitution
\[
\iota(X)=(2,2),\qquad \iota(Y)=(1,1),
\]
and extend it by concatenation to finite, one-sided infinite, and bi-infinite words. For a bi-infinite word $\mathbf w=(w_i)_{i\in\mathbb Z}$, the block $\iota(w_i)$ occupies positions $2i,2i+1$. We write $\mathcal S(\mathbf w):=\mathcal S(\iota(\mathbf w))$; thus the Markov value of an $X,Y$-word always means the value of its expanded integer sequence.

\begin{prop}\label{prop:mechanical-sign}\index{Mechanical word!periodicity}\index{Mechanical word}
Let $\mathbf w$ be a mechanical word of rational slope $t$, allowing $t=1/0$. Then, in the case $(k_1,k_2,k_3)=(0,0,0)$, there exists $j\in\mathbb Z$ such that
\[
\iota(\mathbf w)=T^j(\dots,S(t),S(t),S(t),\dots),
\]
where $S(t)$ is the corresponding generalized strongly admissible sequence.
\end{prop}
\begin{proof}
The case $t=1/0$ follows by direct inspection. Assume $1\le t<\infty$, and write $t=p/q$ in lowest terms. Let $\mathbf w_0$ be the right mechanical word defined by the line obtained by extending $\overline{L_t}$ periodically by translation through $(q,p)$. On an interval contributing $X$, this line passes through four triangles of $\widetilde{\mathbb R^2}$ in the $x$-direction, with signs $-,-,+,+$. On an interval contributing $Y$, it passes through two triangles in the $y$-direction, with signs $-,+$; see Figure~\ref{fig:mechanical2}. Thus replacing $X$ by $(2,2)$ and $Y$ by $(1,1)$ in $\mathbf w_0$ gives the periodic repetition of $S(t)$, up to a choice of index origin.

By Lemmas~\ref{lem:rational-left-is-right} and~\ref{lem:intercept-rational}, any mechanical word $\mathbf w$ of slope $t$ is $T^h\mathbf w_0$ for some $h\in\mathbb Z$. Each letter is replaced by two entries, so $\iota(\mathbf w)=T^{2h}\iota(\mathbf w_0)$. Absorbing the choice of origin into $j$ proves the assertion.
\begin{figure}[ht]
\centering
    \begin{tikzpicture}
    \draw[gray] (-1,0) grid (1,1);
    \draw[red,thick] (-0.5,0) -- (0.4,1);
    \draw[line width=2pt] (0,0) -- (1,1);
    \draw[gray] (-1,1) -- (0,0);
    \draw[gray] (0,1) -- (1,0);
    \node at (1.2,0.5) {$X$};
    \node[red] at (-0.75,0.25) {$-$};
    \node[red] at (-0.3,0.75) {$-$};
    \node[red] at (0.25,0.5) {$+$};
    \node[red] at (0.5,0.75) {$+$};
    \node[red] at (0,0.55) {\rotatebox{55}{$>$}};
\end{tikzpicture}\quad 
\begin{tikzpicture}
    \draw[gray] (0,0) grid (1,1);
    \draw[red,thick] (0.2,0) -- (0.7,1);
    \draw[line width=2pt] (1,0) -- (1,1);
    \draw[gray] (0,1) -- (1,0);
    \node at (1.2,0.5) {$Y$};
    \node[red] at (0.2,0.35) {$-$};
    \node[red] at (0.75,0.5) {$+$};
    \node[red] at (0.47,0.55) {\rotatebox{65}{$>$}};
\end{tikzpicture}
    \caption{The relation between the sign rule and right mechanical words}
    \label{fig:mechanical2}
\end{figure}
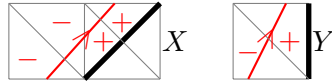
\end{proof}

\section{Markov's Theorem}\label{sec:markov-theorem}
In this section we consider the case $(k_1,k_2,k_3)=(0,0,0)$. In this case the $(1,1)$-, $(2,1)$-, and $(3,1)$-entries at each vertex of $\mathrm{M}\mathbb T(0,0,0,\sigma)$ do not depend on $\sigma\in\mathfrak S_3$, so the fraction labeling of GM numbers is the same for every choice of $\sigma$. We therefore omit $\sigma$. Applying Theorem~\ref{thm:markov-value-gen} in this case gives $\mathcal M_{0,0,0}\subset \mathcal L\cap(0,3)\subset \mathcal M\cap(0,3)$. The reverse inclusion is the classical theorem known as Markov's theorem, which we now prove.
\begin{theorem}[Markov's Theorem]\label{thm:markov}\index{Lagrange spectrum!values below three}\index{Lagrange spectrum}\index{Markov spectrum}\index{Markov's theorem}
$\mathcal M_{0,0,0}=\mathcal L\cap(0,3)=\mathcal M\cap(0,3)$.
\end{theorem}
It suffices to prove $\mathcal M_{0,0,0}\supset\mathcal M\cap(0,3)$. By Corollary~\ref{cor:M=S}, an element of $\mathcal M\cap(0,3)$ is represented by a bi-infinite sequence $\mathbf b$ of positive integers with $\mathcal S(\mathbf b)\in(0,3)$. It remains to classify the sequences with $\mathcal S(\mathbf b)<3$. We first study the broader condition $\mathcal S(\mathbf b)\leq3$; the equality case will be needed in the substitution argument below. From the definition
\[
\mathcal S(\mathbf b):=\sup_{h\in\mathbb Z}\ell_h(\mathbf b)
=\sup_{h\in\mathbb Z}([b_h;b_{h+1},\dots]+[0;b_{h-1},b_{h-2},\dots]).
\]
the following observation is immediate.
\begin{prop}
If $\mathcal S(\mathbf b)\leq3$, then $b_i=1$ or $b_i=2$ for every $i\in\mathbb Z$.
\end{prop}
In the remainder of the proof, we work under the assumption that each $b_i$ is either $1$ or $2$.
We introduce the notation
\[
\ell(\dots,b_{h-2},b_{h-1}\mid b_h,b_{h+1},\dots)
:=[b_h;b_{h+1},\dots]+[0;b_{h-1},b_{h-2},\dots].
\]
This lets us write computations of \(\ell_h(\mathbf b)\) without specifying the index \(h\) explicitly.
\begin{lemm}\label{lem:reverse-invariant}\index{Reversal of a word}
For one-sided infinite sequences $u,v$ and $c\in\mathbb Z_{\geq1}$,
$\ell(u^{\ast}\mid c,v)=\ell(v^{\ast}\mid c,u)$. Here $u^{\ast}$ denotes the reversal of $u$. In particular, $\mathcal S(\mathbf b)\leq3$ if and only if $\mathcal S(\mathbf b^{\ast})\leq3$.
\end{lemm}
\begin{proof}
This follows immediately from the computation
    \[
\ell(u^{\ast} \mid c,v)=[c;v]+[0;u]
=\left(c+\frac{1}{[u]}\right)+\frac{1}{[v]}
=\left(c+\frac{1}{[v]}\right)+\frac{1}{[u]}
=[c;u]+[0;v]=\ell(v^{\ast} \mid c,u).
\]
\end{proof}
We also have the following characterization.
\begin{prop}\label{prop:characterization-sequence}
The condition $\mathcal S(\mathbf b)\leq3$ is equivalent to the following two conditions on $\mathbf b$.
\begin{itemize}
    \item[(i)] The sequence $\mathbf b$ contains neither the consecutive block $(1,2,1)$ nor the consecutive block $(2,1,2)$.
    \item[(ii)] If $\mathbf b$ or $\mathbf b^{\ast}$ can be written as $(u^{\ast},1,1,2,2,v)$ using one-sided infinite sequences $u,v$, then $[v]\leq[u]$.
\end{itemize}
\end{prop}
\begin{proof}
We first show that $\mathcal S(\mathbf b)\leq3$ implies (i) and (ii). Since the continued-fraction values $[u]$ and $[v]$ are greater than $1$, we obtain
\[
\ell(u^{\ast},1\mid2,1,v)=[2;1,v]+[0;1,u]>[2;1,1]+[0;1,1]=\frac{5}{2}+\frac{1}{2}=3.
\]
Thus $(1,2,1)$ cannot occur in $\mathbf b$. The displayed inequality follows from \eqref{eq:pk/qk-pk-1/qk-1}, applied to the positive real values $[u]$ and $[v]$. If $\mathbf b$ contains $(2,1,2)$, the preceding result lets us extend it to $(2,2,1,2)$. Then
\[
\ell(u^{\ast},2\mid2,1,2,v)=[2;1,2,v]+[0;2,u]>[2;1,2]+[0;2,1]=\frac{8}{3}+\frac{1}{3}=3.
\]
This proves (i). If $[v]>[u]$, then
\[
\ell(u^{\ast},1,1\mid2,2,v)=[2;2,v]+[0;1,1,u]>[2;2,v]+[0;1,1,v]=\frac{5[v]+2}{2[v]+1}+\frac{[v]+1}{2[v]+1}=3.
\]
so (ii) also holds.

Conversely, assume (i) and (ii). We distinguish the three possible cuts: (1) $\ell(u^{\ast}\mid1,v)$, (2) $\ell(u^{\ast},1\mid2,v)$, and (3) $\ell(u^{\ast},2\mid2,v)$. In case (1), we have
\[
\ell(u^{\ast}\mid1,v)=[1;v]+[0;u]<[1;1]+[0;1]=3.
\]
Hence the required inequality holds.

In case (2), since neither $(1,2,1)$ nor $(2,1,2)$ occurs, we may write $u=(1,u')$ and $v=(2,v')$. Hence
\[
\ell(u^{\ast},1\mid2,v)=\ell((u')^{\ast},1,1\mid2,2,v')=[2;2,v']+[0;1,1,u']\leq[2;2,u']+[0;1,1,u']=3.
\]
Thus the inequality holds.
 
In case (3), we further divide into (3-1) $v=(1,v')$ and (3-2) $v=(2,v')$. In case (3-1), the absence of $(2,1,2)$ implies $v'=(1,v'')$. By Lemma~\ref{lem:reverse-invariant},
\[
\ell(u^{\ast},2\mid2,1,v')=\ell(u^{\ast},2\mid2,1,1,v'')=\ell((v'')^{\ast},1,1\mid2,2,u)\leq 3.
\]
where the final inequality is obtained by applying the same argument as in case (2) to $\mathbf b^{\ast}$. In case (3-2),
\[
\ell(u^{\ast},2\mid2,2,v')=[2;2,v']+[0;2,u]<[2;2]+[0;2]=3.
\]
and the proof is complete.
\end{proof}

\begin{rema}
Condition (ii) will be used below in the following form. If $x,y\in\{1,2\}$ and $\mathbf b=(\dots,x,w^{\ast},1,1,2,2,w,y,\dots)$, then
\[
    \text{$x=1,y=2$}\Rightarrow\text{$|w|$ is odd},\quad \text{$x=2,y=1$}\Rightarrow\text{$|w|$ is even}
\]
This is just a reformulation of condition (ii).
\end{rema}
When we write $1^n$ or $2^m$ as part of a sequence, it means that $1$ appears $n$ times consecutively or that $2$ appears $m$ times consecutively.
\begin{theorem}\label{thm:degen-constant-regular}\index{Constant type}\index{Degenerate type}\index{Regular type}
If $\mathcal S(\mathbf b)\le3$, then $\mathbf b$ is one of the following types.
\begin{itemize}
    \item[(1)] $(\overline{1}^{\ast},2,2,\overline1)$ or $(\overline{2}^{\ast},1,1,\overline2)$ (degenerate type)
    \item[(2)] $(\overline1^{\ast},\overline{1})$ or $(\overline2^{\ast},\overline{2})$ (constant type)
    \item[(3)] $(\dots,1^{m_{i-1}},2^{n_{i-1}},1^{m_i},2^{n_i},1^{m_{i+1}},2^{n_{i+1}},\dots)$ (regular type)
\end{itemize}
where $m_j$ and $n_j$ are even for every $j$.
\end{theorem}
\begin{proof}
First consider a finite maximal block $1^m$ and suppose that $m$ is odd. Proposition~\ref{prop:characterization-sequence} (i) gives $m\ge3$ and excludes a following block consisting of a single $2$. If the following block has length at least $3$, including an infinite block, then $(1,1,1,2,2,2)$ occurs, contradicting condition (ii). Hence the following block is $2,2$. Write
\[
\mathbf b=(u^{\ast},1,1,2,2,v),\qquad
u=(1^{m-2},2,\ldots).
\]
Let $p$ be the length of the initial block of $1$'s in $v$, allowing $p=\infty$. Condition (ii) gives $[v]\le[u]$. Comparing the first differing entries in the continued fractions shows that $p$ is finite and odd, and $p\le m-2$: if $p>m-2$, including $p=\infty$, or if $p<m-2$ is even, then $[v]>[u]$. Repeating the argument gives a strictly decreasing sequence of positive odd block lengths, which is impossible. Thus every finite maximal block of $1$'s has even length.

Now suppose a finite maximal block $2^n$ has odd length. Condition (i) gives $n\ge3$, and condition (ii) implies that the preceding block of $1$'s has length $2$. Write
\[
\mathbf b=(u^{\ast},1,1,2,2,v),\qquad
v=(2^{n-2},1,\ldots).
\]
If $q$ is the length of the initial block of $2$'s in $u$, allowing $q=\infty$, the same first-difference comparison in $[v]\le[u]$ gives a finite odd $q\le n-2$. Iterating to the left again gives a strictly decreasing sequence of positive odd integers. Hence every finite maximal block of $2$'s also has even length.

Suppose next that $\mathbf b$ has a constant left tail of $1$'s but is not constant. The first nonconstant block must have length $2$, so
\[
\mathbf b=(\overline1^{\ast},2,2,v).
\]
If $v\ne\overline1$, write $v=(1^p,2,\ldots)$ with finite $p\ge0$. Condition (ii) forces $p$ to be odd, contradicting the preceding result. Hence this is the degenerate sequence $(\overline1^{\ast},2,2,\overline1)$.

For a constant left tail of $2$'s, the first block of $1$'s likewise has length $2$, giving $(\overline2^{\ast},1,1,v)$. If $v\ne\overline2$, write $v=(2^p,1,\ldots)$. Applying condition (ii) to the reversed sequence at this boundary gives
\[
[\overline2]\le[2^p,1,\ldots],
\]
which again forces $p$ to be odd. Thus the only nonconstant possibility is $(\overline2^{\ast},1,1,\overline2)$. Constant right tails are handled by reversal and Lemma~\ref{lem:reverse-invariant}. If neither tail is constant, all maximal blocks are finite, and the sequence has the stated regular form with even block lengths.
\end{proof}
We next compute the values for the degenerate and constant types.
\begin{prop}\label{prop:degen-constant}\index{Constant type}\index{Degenerate type}
If $\mathbf b$ is of degenerate type, then $\mathcal S(\overline{1}^{\ast},2,2,\overline1)=\mathcal S(\overline{2}^{\ast},1,1,\overline2)=3$. If $\mathbf b$ is of constant type, then $\mathcal S(\overline1^{\ast},\overline{1})=\sqrt{5}$ and $\mathcal S(\overline2^{\ast},\overline{2})=2\sqrt{2}$.
\end{prop}
\begin{proof}
The values $\mathcal S(\overline1^{\ast},\overline1)$ and $\mathcal S(\overline2^{\ast},\overline2)$ are obtained directly from the definition, since every cut gives the same decomposition.

We show that $\mathcal S(\overline{1}^{\ast},2,2,\overline1)=3$.
\[\ell(\overline{1}^{\ast},1\mid2,2,\overline1)=[2;2,\overline{1}]+[0;1,1,\overline1]=\frac{5[\overline1]+2}{2[\overline1]+1}+\frac{[\overline1]+1}{2[\overline1]+1}=3\]
By this computation and Lemma~\ref{lem:reverse-invariant},
\[\ell(\overline{1}^{\ast},1,2\mid2,\overline1)=3\]
Moreover,
\[
\ell(u^{\ast},1\mid1,v)=[1;v]+[0;1,u]<2+1=3\]
Thus only cuts between $1$ and $2$ need be considered. The remaining pattern is $\ell(\overline{1}^{\ast},2,2\mid1,\overline1)$, which equals $\ell(\overline{1}^{\ast},1\mid1,2,2,\overline1)$ by Lemma~\ref{lem:reverse-invariant} and is less than $3$.

Next we prove $\mathcal S(\overline{2}^{\ast},1,1,\overline2)=3$.
\[\ell(\overline{2}^{\ast},1,1\mid 2,\overline2)=[2;2,\overline{2}]+[0;1,1,\overline2]=\frac{5[\overline2]+2}{2[\overline2]+1}+\frac{[\overline2]+1}{2[\overline2]+1}=3\]
By this computation and Lemma~\ref{lem:reverse-invariant},
\[\ell(\overline{2}^{\ast},2\mid2,1,1,\overline2)=3\]
Moreover,
\[
\ell(u^{\ast},2\mid2,2,v)=[2;2,v]+[0;2,u]<\frac{5}{2}+\frac{1}{2}=3\]
Together with the preceding equality, this shows that every cut between two $2$'s gives a value at most $3$. The inequality $\ell(\overline2^{\ast},1\mid1,\overline2)<3$ was proved in the first part. The remaining pattern is $\ell(\overline2^{\ast},2\mid1,1,\overline2)$, which equals $\ell(\overline2^{\ast},1\mid1,\overline2)$ by Lemma~\ref{lem:reverse-invariant} and is therefore less than $3$.
\end{proof}

By Theorem~\ref{thm:degen-constant-regular}, if $\mathcal S(\mathbf b)\le3$, then $\mathbf b$ can be written using the blocks $(1,1)$ and $(2,2)$. Equivalently, up to a shift, $\mathbf b=\iota(\mathbf w)$ for a bi-infinite word $\mathbf w$ in $X,Y$. Once a block decomposition and an index origin have been chosen, we also use $\mathbf b$ for this associated $X,Y$-word. Its Markov value is understood through $\iota$ as above, and continued fractions of words in $X,Y$ likewise mean those of their expansions under $\iota$. With this convention, the preceding classification becomes
\begin{itemize}
    \item[(1)] $\overline{X}^{\ast}Y\overline X$ or $\overline{Y}^{\ast}X\overline Y$ (degenerate type)
    \item[(2)] $\overline X^{\ast}\overline X$ or $\overline Y^{\ast}\overline Y$ (constant type)
    \item[(3)] $\cdots X^{k_{i-1}}Y^{\ell_{i-1}}X^{k_i}Y^{\ell_i}X^{k_{i+1}}Y^{\ell_{i+1}}\cdots$ (regular type)
\end{itemize}
Proposition~\ref{prop:characterization-sequence} can now be restated as follows.

\begin{prop}\label{prop:characterization2}
For a bi-infinite sequence $\mathbf b$, the following are equivalent.
\begin{itemize}
\item[(1)] $\mathcal S(\mathbf b)\leq3$.
\item[(2)] The sequence $\mathbf b$ can be written as a sequence in the letters $X,Y$, and the following two conditions hold.
\begin{itemize}
    \item[(i)] If $\mathbf b$ has a representation $\mathbf b=u^{\ast}YXv$, then either $u=v$, or there exist words $w,u',v'$ in $X,Y$ such that $\mathbf b=(u')^{\ast}Xw^{\ast}YXwYv'$.
    \item[(ii)] If $\mathbf b$ has a representation $\mathbf b=u^{\ast}XYv$, then either $u=v$, or there exist words $w,u',v'$ in $X,Y$ such that $\mathbf b=(u')^{\ast}Yw^{\ast}XYwXv'$.
\end{itemize}
\end{itemize}
\end{prop}
\begin{proof}
By Theorem~\ref{thm:degen-constant-regular}, condition (1) implies that $\mathbf b$ can be written in the letters $X,Y$. At a central block $YX$, Proposition~\ref{prop:characterization-sequence} (ii) gives $[v]\le[u]$. If $u\ne v$, let $w$ be their common initial word. Each letter expands into two entries, so the first differing continued-fraction entries have the even index $2|w|$, starting at index $0$. Thus $[v]<[u]$ holds exactly when the next letter of $u$ is $X$ and that of $v$ is $Y$. This is condition (i). Reversing the sequence gives condition (ii) at a central block $XY$.

Conversely, writing $\mathbf b$ in $X,Y$ already implies condition (i) of Proposition~\ref{prop:characterization-sequence}. The same comparison at the first differing entries shows that the two conditions in (2) imply its condition (ii). Hence $\mathcal S(\mathbf b)\le3$.
\end{proof}
The following consequence will be useful.
\begin{prop}\label{prop:regular-ki=1}\index{Regular type}
Let $\mathbf b=\cdots X^{k_{i-1}}Y^{\ell_{i-1}}X^{k_i}Y^{\ell_i}X^{k_{i+1}}Y^{\ell_{i+1}}\cdots$ be a regular bi-infinite sequence with $\mathcal S(\mathbf b)\le3$. Then either $k_i=1$ for every $i\in\mathbb Z$, or $\ell_i=1$ for every $i\in\mathbb Z$.
\end{prop}
\begin{proof}
Suppose neither alternative holds. Then a finite part of $\mathbf b$ has one of the forms
\[
Y^m(XY)^kX^n\quad\text{or}\quad X^m(YX)^kY^n,
\qquad m,n\ge2.
\]
Choose such a part with $k\ge0$ minimal. The case $k=0$ contains $YYXX$ or $XXYY$ and contradicts Proposition~\ref{prop:characterization2}. Hence $k\ge1$.

Consider the first form and the last central block $YX$ in it. The outward words are
\[
u=X(YX)^{k-1}Y^m\cdots,
\qquad v=X^{n-1}\cdots.
\]
Condition (i), in lexicographic order with $Y<X$, requires $v\preceq u$. Comparing the initial letters gives $n=2$. Starting after this block $X^2$, another $XX$ before the first $YY$ would make $v\succ u$. If no $YY$ occurs before the alternating word $X(YX)^k$ has been completed, the comparison likewise gives $v\succ u$ at the position where $u$ has $YY$. Thus
\[
v=X(YX)^{k'}YY\cdots
\quad\text{for some }0\le k'<k.
\]
The original word therefore contains $X^2(YX)^{k'}Y^2$, contrary to the minimality of $k$. The other form is treated by interchanging $X,Y$ and using condition (ii).
\end{proof}
\index{X-type@$X$-type}\index{Y-type@$Y$-type}In the preceding proposition, we call the first case, where $k_i=1$ for every $i$, the \emph{$X$-type}, and the second case, where $\ell_i=1$ for every $i$, the \emph{$Y$-type}.
\begin{defi}\label{def:characteristic-word-sequences}\index{Characteristic sequence of a word}
Let $\mathbf b=\cdots X^{k_{i-1}}Y^{\ell_{i-1}}X^{k_i}Y^{\ell_i}X^{k_{i+1}}Y^{\ell_{i+1}}\cdots$ be a regular bi-infinite sequence; here we do not impose $\mathcal S(\mathbf b)\le3$. Define
\begin{align*}
    C^X(\mathbf b):=&(\dots,(0)^{k_{i-1}-1},\ell_{i-1},(0)^{k_{i}-1},\ell_i,(0)^{k_{i+1}-1},\ell_{i+1},\dots),\\
    C^Y(\mathbf b):=&(\dots,k_{i-1},(0)^{\ell_{i-1}-1},k_{i},(0)^{\ell_i-1},k_{i+1},(0)^{\ell_{i+1}-1},\dots)\\
\end{align*}
We call these sequences the \emph{$X$-characteristic sequence} and the \emph{$Y$-characteristic sequence} of $\mathbf b$, respectively.
\end{defi}
If $\mathbf b$ is of $X$-type, then the $X$-characteristic sequence has all entries at least $1$, whereas the $Y$-characteristic sequence consists only of $0$'s and $1$'s. For $Y$-type the roles are reversed.
In terms of characteristic sequences, Proposition~\ref{prop:characterization2} can be reformulated as follows.
\begin{prop}\label{prop:characterization3}\index{Characteristic sequence of a word}\index{Regular type}
For a regular bi-infinite sequence $\mathbf b$, the following three conditions are equivalent.
\begin{itemize}
\item[(1)] $\mathcal S(\mathbf b)\leq3$.
\item[(2)]
The $X$-characteristic sequence $C^X(\mathbf b)=(x_i)_{i\in\mathbb Z}$ is defined and, for every $i\in\mathbb Z$,
 \begin{align*}
     (x_{i}-1,x_{i+1},x_{i+2},\dots)&\preceq(x_{i-1},x_{i-2},x_{i-3},\dots),\\
     (x_{i}-1,x_{i-1},x_{i-2},\dots)&\preceq(x_{i+1},x_{i+2},x_{i+3},\dots)
\end{align*}
hold.
 \item[(3)]
The $Y$-characteristic sequence $C^Y(\mathbf b)=(y_i)_{i\in\mathbb Z}$ is defined and, for every $i\in\mathbb Z$,
 \begin{align*}
     (y_{i}-1,y_{i+1},y_{i+2},\dots)&\preceq(y_{i-1},y_{i-2},y_{i-3},\dots),\\
      (y_{i}-1,y_{i-1},y_{i-2},\dots)&\preceq(y_{i+1},y_{i+2},y_{i+3},\dots)
 \end{align*}
hold.
\end{itemize}
Here $\preceq$ denotes lexicographic order.
\end{prop}
\begin{proof}
We prove the equivalence of (1) and (2). Index the occurrences of $X$ so that $x_i$ is the number of $Y$'s between the $i$th and $(i+1)$st occurrences of $X$; locally the word is
\[
\cdots XY^{x_{i-1}}XY^{x_i}XY^{x_{i+1}}X\cdots.
\]
If $x_i=0$, the two required inequalities hold automatically because their left-hand sides start with $-1$ and their right-hand sides with nonnegative integers. There is then no $XY$ or $YX$ boundary in this empty block to check.

If $x_i\ge1$, at the first $XY$ boundary in $XY^{x_i}$ the outward words are
\[
u=Y^{x_{i-1}}XY^{x_{i-2}}X\cdots,
\qquad
v=Y^{x_i-1}XY^{x_{i+1}}X\cdots.
\]
At the first difference of two such words, the word with the shorter initial $Y$-block has the letter $X$, and the other has $Y$. Therefore Proposition~\ref{prop:characterization2} (ii) is exactly
\[
(x_i-1,x_{i+1},x_{i+2},\ldots)
\preceq(x_{i-1},x_{i-2},x_{i-3},\ldots).
\]
At the last $YX$ boundary of this block, the outward words are
\[
u=Y^{x_i-1}XY^{x_{i-1}}X\cdots,
\qquad
v=Y^{x_{i+1}}XY^{x_{i+2}}X\cdots.
\]
Condition (i) is exactly the second displayed inequality in (2). This accounts for every $XY$ and $YX$ boundary, proving the equivalence. Interchanging $X,Y$ proves the equivalence of (1) and (3).
\end{proof}

\index{Substitution of words}We now temporarily forget that $X$ and $Y$ stand for $(2,2)$ and $(1,1)$, and consider the free group on the letters $X,Y$, denoted by $\mathfrak F(X,Y)$. Define automorphisms $\lambda,\rho\in\operatorname{Aut}\mathfrak F(X,Y)$ by
\[
\lambda\colon
\begin{cases}
X\mapsto X,\\
Y\mapsto XY,
\end{cases}
\qquad
\rho\colon
\begin{cases}
X\mapsto XY,\\
Y\mapsto Y.
\end{cases}
\]
Their inverses are
\[
\lambda^{-1}\colon
\begin{cases}
X\mapsto X,\\
Y\mapsto X^{-1}Y,
\end{cases}
\qquad
\rho^{-1}\colon
\begin{cases}
X\mapsto XY^{-1},\\
Y\mapsto Y.
\end{cases}
\]
Let $\tau\in\{\lambda,\rho\}$ and let $\mathbf b=(b_i)_{i\in\mathbb Z}$ be a bi-infinite word. Concatenate the words $\tau(b_i)$ in order,
\[
\cdots\tau(b_{-2})\tau(b_{-1})\,\big|\,\tau(b_0)\tau(b_1)\tau(b_2)\cdots,
\]
and reindex so that the first letter of $\tau(b_0)$ has index $0$. This defines $\tau(\mathbf b)$. A different choice of index origin changes the result only by a shift, and all the conditions considered here---the value of $\mathcal S$, regularity, and the characteristic-sequence conditions---are shift-invariant.

We apply the inverse substitutions to positive bi-infinite words only in the following two cases. If $\mathbf b$ is of $X$-type, then it has a unique expression
\[
\mathbf b=\cdots XY^{\ell_{i-1}}XY^{\ell_i}XY^{\ell_{i+1}}\cdots
\qquad(\ell_i\ge1).
\]
Reading the initial $XY$ of each block as $\rho(X)$ and every remaining $Y$ as $\rho(Y)$ gives a unique parsing and defines
\[
\rho^{-1}(\mathbf b)
=\cdots XY^{\ell_{i-1}-1}XY^{\ell_i-1}XY^{\ell_{i+1}-1}\cdots.
\]
Similarly, if $\mathbf b$ is of $Y$-type, write uniquely
\[
\mathbf b=\cdots X^{k_{i-1}}YX^{k_i}YX^{k_{i+1}}Y\cdots
\qquad(k_i\ge1).
\]
Reading the terminal $XY$ in each block as $\lambda(Y)$ and each preceding $X$ as $\lambda(X)$ defines
\[
\lambda^{-1}(\mathbf b)
=\cdots X^{k_{i-1}-1}YX^{k_i-1}YX^{k_{i+1}-1}Y\cdots.
\]
An exponent $0$ means that the corresponding word is empty. Whenever $\rho^{-1}$ or $\lambda^{-1}$ is applied to a bi-infinite word below, it is understood in this sense.

\begin{lemm}\label{lem:lambda-rho}\index{X-type@$X$-type}\index{Y-type@$Y$-type}
Let $\mathbf b$ be a regular bi-infinite word satisfying $\mathcal S(\mathbf b)\le3$. Then
\[
\mathcal S(\lambda(\mathbf b))\le3,
\qquad
\mathcal S(\rho(\mathbf b))\le3.
\]
Moreover, if $\mathbf b$ is of $X$-type, then $\mathcal S(\rho^{-1}(\mathbf b))\le3$, and if $\mathbf b$ is of $Y$-type, then $\mathcal S(\lambda^{-1}(\mathbf b))\le3$.
\end{lemm}

\begin{proof}
The characteristic sequences satisfy
\[
C^X(\rho(\mathbf b))
=C^X(\mathbf b)+(\overline1^{\ast},\overline1),
\qquad
C^Y(\lambda(\mathbf b))
=C^Y(\mathbf b)+(\overline1^{\ast},\overline1).
\]
Proposition~\ref{prop:characterization3} therefore gives
$\mathcal S(\lambda(\mathbf b))\le3$ and
$\mathcal S(\rho(\mathbf b))\le3$.

We next consider the inverse substitutions. It is enough to treat the case in which $\mathbf b$ is of $X$-type; the $Y$-type case follows by interchanging $X$ with $Y$ and $\rho$ with $\lambda$. Write
\[
\mathbf b=\cdots XY^{x_{i-1}}XY^{x_i}XY^{x_{i+1}}\cdots
\qquad(x_i\ge1).
\]
Then $C^X(\mathbf b)=(x_i)_{i\in\mathbb Z}$ and
\[
\rho^{-1}(\mathbf b)
=\cdots XY^{x_{i-1}-1}XY^{x_i-1}XY^{x_{i+1}-1}\cdots.
\]
Put $d_i:=x_i-1\ge0$. Subtracting $1$ from every component on both sides of the two lexicographic inequalities in Proposition~\ref{prop:characterization3} gives
\begin{align*}
(d_i-1,d_{i+1},d_{i+2},\ldots)
&\preceq(d_{i-1},d_{i-2},d_{i-3},\ldots),\\
(d_i-1,d_{i-1},d_{i-2},\ldots)
&\preceq(d_{i+1},d_{i+2},d_{i+3},\ldots).
\end{align*}

If indices with $d_i>0$ occur infinitely often in both directions, then $\rho^{-1}(\mathbf b)$ is regular and
$C^X(\rho^{-1}(\mathbf b))=(d_i)_{i\in\mathbb Z}$. The displayed inequalities and Proposition~\ref{prop:characterization3} imply
$\mathcal S(\rho^{-1}(\mathbf b))\le3$.

If $d_i=0$ for every $i$, then
$\rho^{-1}(\mathbf b)=\overline X^{\ast}\overline X$ is of constant type. Proposition~\ref{prop:degen-constant} gives
\[
\mathcal S(\rho^{-1}(\mathbf b))=2\sqrt2<3.
\]

It remains to consider the case in which the set of indices with $d_i>0$ is nonempty and is bounded below or above. Suppose first that it is bounded below, and let $i$ be its least element. Then $x_{i-1}=x_{i-2}=\cdots=1$, so the first inequality of Proposition~\ref{prop:characterization3} gives
\[
(d_i,x_{i+1},x_{i+2},\ldots)\preceq(1,1,1,\ldots).
\]
Every component on the left is a positive integer and $d_i\ge1$. Hence $d_i=1$ and $x_{i+1}=x_{i+2}=\cdots=1$. Thus $d_i=1$ and every other $d_j$ is $0$. If the set of indices with $d_i>0$ is bounded above, applying the second inequality at its greatest element gives the same conclusion. In either case,
\[
\rho^{-1}(\mathbf b)=\overline X^{\ast}Y\overline X
\]
is of degenerate type, and Proposition~\ref{prop:degen-constant} gives
$\mathcal S(\rho^{-1}(\mathbf b))=3$.

These cases exhaust all possibilities. Therefore an $X$-type word satisfies
$\mathcal S(\rho^{-1}(\mathbf b))\le3$. The corresponding assertion for a $Y$-type word follows in the same way.
\end{proof}

We now use these conditions to describe the inequality $\mathcal S(\mathbf b)\le3$ in terms of mechanical words.
\begin{lemm}\label{lem:mechanical-substitution-right-left}\index{Mechanical word}\index{Mechanical word!slope}
Let $\mathbf c$ be a mechanical word of slope $t\in[1,\infty]$. Then $\lambda(\mathbf c)$ and $\rho(\mathbf c)$ are also mechanical words. More precisely, their slopes are respectively
\[
    L(t)=2-\frac1t,
    \qquad
    R(t)=t+1,
\]
with the conventions $L(\infty)=2$ and $R(\infty)=\infty$.
\end{lemm}
\begin{proof}
The case $t=\infty$ is immediate. Indeed, $\mathbf c=\cdots YYY\cdots$. Under $\lambda$ this becomes the alternating word $\cdots XYXYXY\cdots$, which is mechanical of slope $2$, while under $\rho$ it remains $\cdots YYY\cdots$, which is mechanical of slope $\infty$.

The case $t=1$ is also immediate. Then $\mathbf c=\cdots XXX\cdots$. Under $\lambda$ it remains $\cdots XXX\cdots$, which is mechanical of slope $1$, and under $\rho$ it becomes the alternating word $\cdots XYXYXY\cdots$, which is mechanical of slope $2$.

Assume $1<t<\infty$. We first treat the case where $\mathbf c$ is a right mechanical word. Write $\mathbf c=\mathbf b^R(t,\theta)$ and put
\[
    r_n=\left\lceil \frac{n-\theta}{t}\right\rceil.
\]
Then
\[
    c_n=X
    \Longleftrightarrow r_{n+1}-r_n=1,
    \qquad
    c_n=Y
    \Longleftrightarrow r_{n+1}-r_n=0.
\]

First consider $\lambda$. Let $P_n$ be the position at which $\lambda(c_n)$ begins, normalized by $P_0=0$. Since $\lambda(X)=X$ and $\lambda(Y)=XY$, we have $P_n=2n-(r_n-r_0)$. For each $n$, write
\[
    r_n=\frac{n-\theta}{t}+\delta_n,
    \qquad 0\le\delta_n<1.
\]
Put $t_\lambda=2-1/t$ and $\theta_\lambda=r_0+\theta/t$. We show that $\lambda(\mathbf c)$ is the right mechanical word of slope $t_\lambda$ and intercept $\theta_\lambda$. Define
\[
    R_m:=\left\lceil \frac{m-\theta_\lambda}{t_\lambda}\right\rceil.
\]
Then
\[
\begin{aligned}
    P_n=2n-(r_n-r_0) =\left(2-\frac1t\right)n+r_0+\frac{\theta}{t}-\delta_n 
    =t_\lambda n+\theta_\lambda-\delta_n.
\end{aligned}
\]
Thus
\[
    \frac{P_n-\theta_\lambda}{t_\lambda}
    =n-\frac{\delta_n}{t_\lambda},
\]
and since $0\le\delta_n/t_\lambda<1$, we obtain $R_{P_n}=n$.

If $c_n=X$, then $r_{n+1}-r_n=1$, so $P_{n+1}=P_n+1$. Hence $R_{P_n+1}-R_{P_n}=R_{P_{n+1}}-R_{P_n}=1$. This corresponds to $\lambda(c_n)=X$.

If $c_n=Y$, then $r_{n+1}-r_n=0$, so $P_{n+1}=P_n+2$. Moreover
\[
    \frac{P_n+1-\theta_\lambda}{t_\lambda}
    =
    n+\frac{1-\delta_n}{t_\lambda},
\]
and $0<(1-\delta_n)/t_\lambda<1$, so $R_{P_n+1}=n+1$. Since we already know $R_{P_k}=k$ for every $k$, we also have
$R_{P_n+2}=R_{P_{n+1}}=n+1$. Therefore
\[
    R_{P_n+1}-R_{P_n}=1,
    \qquad
    R_{P_n+2}-R_{P_n+1}=0.
\]
This corresponds to $\lambda(c_n)=XY$.

We have shown that
\[
    \lambda(\mathbf c)
    =
    \mathbf b^R\left(2-\frac1t,\ r_0+\frac{\theta}{t}\right).
\]
In particular, $\lambda(\mathbf c)$ is a mechanical word of slope $2-1/t$.

Next consider $\rho$. Let $Q_n$ be the position at which $\rho(c_n)$ begins, normalized by $Q_0=0$. Since $\rho(X)=XY$ and $\rho(Y)=Y$, we have $Q_n=n+(r_n-r_0)$. Put $t_\rho=t+1$ and $\theta_\rho=\theta-r_0$. We show that $\rho(\mathbf c)$ is the right mechanical word of slope $t_\rho$ and intercept $\theta_\rho$. Define
\[
    S_m:=\left\lceil \frac{m-\theta_\rho}{t_\rho}\right\rceil.
\]
As above, write $r_n=(n-\theta)/t+\delta_n$. Then
\[
\begin{aligned}
    \frac{Q_n-\theta_\rho}{t_\rho}
    =
    \frac{n+r_n-r_0-\theta+r_0}{t+1}
    =
    r_n-\frac{t\delta_n}{t+1}.
\end{aligned}
\]
Since $0\le t\delta_n/(t+1)<1$, we have $S_{Q_n}=r_n$.

If $c_n=Y$, then $r_{n+1}-r_n=0$, so $Q_{n+1}=Q_n+1$. Hence $S_{Q_n+1}-S_{Q_n}=S_{Q_{n+1}}-S_{Q_n}=r_{n+1}-r_n=0$. This corresponds to $\rho(c_n)=Y$.

If $c_n=X$, then $r_{n+1}-r_n=1$, so $Q_{n+1}=Q_n+2$. The condition $r_{n+1}-r_n=1$ is equivalent to $\delta_n<1/t$. In this case
\[
    \frac{Q_n+1-\theta_\rho}{t_\rho}
    =
    r_n+\frac{1-t\delta_n}{t+1},
\]
and $0<(1-t\delta_n)/(t+1)<1$, so $S_{Q_n+1}=r_n+1$. Moreover $S_{Q_n+2}=S_{Q_{n+1}}=r_{n+1}=r_n+1$. Therefore
\[
    S_{Q_n+1}-S_{Q_n}=1,
    \qquad
    S_{Q_n+2}-S_{Q_n+1}=0.
\]
This corresponds to $\rho(c_n)=XY$.

Thus
\[
    \rho(\mathbf c)=\mathbf b^R(t+1,\theta-r_0),
\]
and in particular $\rho(\mathbf c)$ is a mechanical word of slope $t+1$.

It remains to consider the case where $\mathbf c$ is a left mechanical word. Write $\mathbf c=\mathbf b^L(t,\theta)$ and put
\[
    l_n=\left\lfloor \frac{n-\theta}{t}\right\rfloor.
\]
The beginning positions of $\lambda(c_n)$ and $\rho(c_n)$ are respectively
\[
    P_n=2n-(l_n-l_0),
    \qquad
    Q_n=n+(l_n-l_0).
\]
Repeating the same calculation with floor functions gives
\[
    \lambda(\mathbf c)
    =
    \mathbf b^L\left(
        2-\frac1t,\,
        l_0+\frac{\theta}{t}-\left(1-\frac1t\right)
    \right),\quad 
    \rho(\mathbf c)
    =
    \mathbf b^L(t+1,\theta-l_0-1).
\]
Thus the assertion also holds for left mechanical words.
\end{proof}

\begin{lemm}\label{lem:infinite-desubstitution-gives-S-equals-3}\index{Substitution of words}
Let $\mathbf b$ be a regular bi-infinite word satisfying $\mathcal S(\mathbf b)\le3$. Put $\mathbf b^{(0)}=\mathbf b$, and suppose that for every $N\ge0$ there exist $\sigma_N\in\{\lambda,\rho\}$ and a bi-infinite word $\mathbf b^{(N+1)}$ such that
\[
\mathbf b^{(N)}=\sigma_N(\mathbf b^{(N+1)}).
\]
If none of the words $\mathbf b^{(N)}$ is of constant type, then $\mathcal S(\mathbf b)=3$.
\end{lemm}

\begin{proof}
We first note that, as long as $\mathbf b^{(N)}$ is regular, the word $\mathbf b^{(N+1)}$ also satisfies $\mathcal S(\mathbf b^{(N+1)})\le3$. A regular word $\mathbf b^{(N)}$ is of $X$-type or of $Y$-type, and the corresponding positive inverse substitution is respectively $\rho^{-1}$ or $\lambda^{-1}$. The two inverse substitutions are simultaneously available only for the alternating word $\cdots XYXY\cdots$. Both inverse images are then constant, contrary to the hypothesis. Hence the inverse image determined by $\sigma_N$ is unique, and Lemma~\ref{lem:lambda-rho} gives $\mathcal S(\mathbf b^{(N+1)})\le3$. By Theorem~\ref{thm:degen-constant-regular} and the hypothesis, $\mathbf b^{(N+1)}$ is either regular or degenerate.

We next record how symmetric subwords propagate under the substitutions. Suppose that a bi-infinite word contains $w^{\ast}YXw$. A direct calculation shows that its image under $\lambda$ contains
\[
(\lambda(w)X)^{\ast}YX\lambda(w)X,
\]
whereas its image under $\rho$ contains
\[
(Y\rho(w))^{\ast}YXY\rho(w).
\]
Indeed,
\[
\lambda(w^{\ast})X=(\lambda(w)X)^{\ast},
\qquad
Y\rho(w^{\ast})=(Y\rho(w))^{\ast},
\]
and both identities follow immediately by induction on $|w|$. Thus, after one substitution, the word on the right of the central $YX$ is replaced by
\[
w\longmapsto\lambda(w)X
\qquad\text{or}\qquad
w\longmapsto Y\rho(w),
\]
and its length increases by at least $1$.

We divide the argument according to whether a degenerate word occurs.

Suppose first that $\mathbf b^{(N)}$ is degenerate for some $N$. In the notation of Proposition~\ref{prop:degen-constant}, it is one of
\[
\overline X^{\,*}Y\overline X,
\qquad
\overline Y^{\,*}X\overline Y.
\]
For every $M\ge0$, the first contains $(X^M)^{\ast}YXX^M$ and the second contains $(Y^M)^{\ast}YXY^M$. Apply $\sigma_{N-1},\ldots,\sigma_0$ successively. By the propagation rule above, the original word $\mathbf b=\mathbf b^{(0)}$ contains a subword $u_M^{\ast}YXu_M$, and $|u_M|\to\infty$ as $M\to\infty$.

Now suppose that no $\mathbf b^{(N)}$ is degenerate. Then every $\mathbf b^{(N)}$ is regular and consequently contains the subword $YX$. For any $N\ge1$, regard such an occurrence in $\mathbf b^{(N)}$ as $\emptyset^{\ast}YX\emptyset$ and apply $\sigma_{N-1},\ldots,\sigma_0$. The propagation rule produces in $\mathbf b$ a subword $w_N^{\ast}YXw_N$. Since the length increases by at least $1$ at every stage, $|w_N|\ge N$.

Thus in either case $\mathbf b$ contains subwords of the form $w^{\ast}YXw$ with $|w|$ arbitrarily large. Returning to $X=(2,2)$ and $Y=(1,1)$, we obtain arbitrarily long subwords of the form
\[
w^{\ast},1,1,2,2,w.
\]
For suitable one-sided infinite sequences $\alpha$ and $\beta$, the word can therefore be cut, for arbitrarily long $w$, in the form
\[
\mathbf b=(\alpha^{\ast},w^{\ast},1,1,2,2,w,\beta).
\]
At this cut,
\[
\mathcal S(\mathbf b)
\ge [2;2,w,\beta]+[0;1,1,w,\alpha].
\]
As $|w|\to\infty$, finite-window approximation for continued fractions shows that the right-hand side tends to $3$. Hence $\mathcal S(\mathbf b)\ge3$. The reverse inequality is part of the hypothesis, so $\mathcal S(\mathbf b)=3$.
\end{proof}

\begin{prop}\label{prop:mechanical-characterization}\index{Markov's theorem}\index{Mechanical word}
Let $\mathbf b$ be a bi-infinite sequence of positive integers with $\mathcal S(\mathbf b)<3$. Then, up to a shift, $\mathbf b=\iota(\mathbf w)$ for a mechanical word $\mathbf w$ of slope $t\in([1,\infty)\cap\mathbb Q)\cup\{\infty\}$.
\end{prop}
\begin{proof}
By Theorem~\ref{thm:degen-constant-regular}, after shifting the integer sequence $\mathbf b$ we may write $\mathbf b=\iota(\mathbf w)$ for a word $\mathbf w$ in $X,Y$. Then $\mathcal S(\mathbf w)=\mathcal S(\mathbf b)<3$. Proposition~\ref{prop:degen-constant} excludes the degenerate type, so $\mathbf w$ is either constant or regular.
If it is constant, then it is
\[
    \cdots XXX\cdots
    \quad\text{or}\quad
    \cdots YYY\cdots,
\]
which are mechanical words of slopes $1$ and $\infty$, respectively.

Assume that $\mathbf w$ is regular. Then at least one of $\rho^{-1}$ and $\lambda^{-1}$ is defined on $\mathbf w$. Choose an available inverse substitution and denote its image by $\mathbf w^{(1)}$. If $\mathbf w^{(1)}$ is regular, repeat the same operation and define $\mathbf w^{(2)}$, continuing until a word $\mathbf w^{(n)}$ is no longer regular. If this process does not reach a constant word in finitely many steps, then there are two possibilities: either it reaches a degenerate word in finitely many steps, or it remains regular forever. These are the only possibilities because every word for which the operation is defined still satisfies $\mathcal S\le3$.
In the first case, if a finite number of desubstitutions reaches a degenerate word, then the degenerate word contains symmetric subwords of the form $v^{\ast}YXv$ with $|v|$ arbitrarily large. The same argument as in the proof of Lemma~\ref{lem:infinite-desubstitution-gives-S-equals-3} then gives $\mathcal S(\mathbf w)=3$, contradicting $\mathcal S(\mathbf w)<3$.

In the second case, Lemma~\ref{lem:infinite-desubstitution-gives-S-equals-3} again gives $\mathcal S(\mathbf w)=3$, a contradiction. Hence $\mathbf w$ is obtained from a constant mechanical word by applying finitely many of $\lambda$ and $\rho$. By Lemma~\ref{lem:mechanical-substitution-right-left}, $\mathbf w$ is a mechanical word.
Moreover, its slope is obtained from $1$ or $\infty$ by applying finitely many times the transformations
\[
    t\mapsto 2-\frac1t,
    \qquad
    t\mapsto t+1.
\]
Therefore the slope belongs to $([1,\infty)\cap\mathbb Q)\cup\{\infty\}$, as required.
\end{proof}

We now prove Markov's theorem.
\begin{proof}[Proof of Theorem~\ref{thm:markov}]
By Proposition~\ref{prop:mechanical-characterization}, if $\mathcal S(\mathbf b)<3$, then $\mathbf b$ is a shift of $\iota(\mathbf w)$ for a mechanical word $\mathbf w$ of slope $t\in([1,\infty)\cap\mathbb Q)\cup\{\infty\}$. By Proposition~\ref{prop:mechanical-sign}, after choosing the index origin we have
\[
\mathbf b=(\dots,S(t),S(t),S(t),\dots).
\]
Theorems~\ref{thm:L=S} and~\ref{thm:markov-value-gen} then give
\[
\mathcal S(\mathbf b)=\frac{\sqrt{(3m_t)^2-4}}{m_t}.
\]
This proves $\mathcal M\cap(0,3)\subset\mathcal M_{0,0,0}$, and hence Markov's theorem.
\end{proof}

The following theorem gives an important restriction on irrational numbers whose Lagrange constants belong to $\mathcal M_{0,0,0}$. Just as binary quadratic forms, or equivalently bi-infinite sequences, with Markov value below $3$ are highly constrained, so are the corresponding irrational numbers.

\begin{theorem}\label{lem:lagrange-below-three-quadratic}\index{Lagrange spectrum!values below three}\index{Markov's theorem}
Let $\alpha$ be irrational and suppose that $\mathcal L(\alpha)<3$. Then there exists $t\in([1,\infty)\cap\mathbb Q)\cup\{\infty\}$ such that
\[
\alpha\sim[\overline{S(t)}].
\]
In particular, $\alpha$ is a quadratic irrational.
\end{theorem}

\begin{proof}
Write
\[
\alpha=[a_0;a_1,a_2,\ldots],
\]
and put
\[
\alpha_n:=[a_n;a_{n+1},a_{n+2},\ldots]
\qquad(n\ge0),
\qquad
\beta_n:=[a_n;a_{n-1},\ldots,a_1]
\qquad(n\ge1).
\]
For $n\ge1$, define
\[
\lambda_n(\alpha):=\alpha_{n+1}+\frac1{\beta_n},
\]
and call $((\alpha_{n+1},\beta_n))_{n\ge1}$ the pair sequence of $\alpha$. By Theorem~\ref{thm:characterization-lagrange1770},
\[
\mathcal L(\alpha)=\limsup_{n\to\infty}\lambda_n(\alpha).
\]
Since $\lambda_n(\alpha)>\alpha_{n+1}>a_{n+1}$, the inequality $\mathcal L(\alpha)<3$ implies that $a_n\in\{1,2\}$ for all sufficiently large $n$. Replacing $\alpha$, if necessary, by one of its complete quotients, we may assume that $a_n\in\{1,2\}$ for every $n\ge1$. This replacement does not change the equivalence class of $\alpha$ by Theorem~\ref{thm:characterization-equivalent}, nor its Lagrange constant by Proposition~\ref{prop:equivalent-lagrange1770}.

The proof of Theorem~\ref{thm:LsubsetS} supplies a bi-infinite sequence $\mathbf b$ determined by a pair associated with an accumulation point realizing $\mathcal L(\alpha)$ and satisfying
\[
\mathcal S(\mathbf b)=\mathcal L(\alpha)<3.
\]
By Propositions~\ref{prop:mechanical-characterization} and~\ref{prop:mechanical-sign}, there exists $t\in([1,\infty)\cap\mathbb Q)\cup\{\infty\}$ such that $\mathbf b$ is a shift of
\[
(\ldots,S(t),S(t),S(t),\ldots).
\]
Let $Y$ be a suitable cyclic shift of $S(t)$, so that we may write
\[
\mathbf b=\cdots YYY\cdots.
\]
By Theorem~\ref{thm:limit-point}, shifting $\mathbf b$ again gives a sequence associated with an accumulation point of the pair sequence of $\alpha$. We may therefore place the cut at the beginning of a copy of $Y$.

Let $(\theta_0,\eta_0)$ be the accumulation pair corresponding to this cut, and choose a subsequence $(n_i)_i$ such that
\[
(\alpha_{n_i+1},\beta_{n_i})\longrightarrow(\theta_0,\eta_0).
\]
The continued-fraction expansion of $\theta_0$ is $YYY\cdots$. The argument in the proof of Proposition~\ref{prop:theta-eta} shows that, for every $k\ge1$ and all sufficiently large $i$, the first $k|Y|$ partial quotients of $\alpha_{n_i+1}$ agree with $Y^k$. Hence $Y^k$ occurs as a finite block in the continued-fraction expansion of $\alpha$ for every $k$.

Suppose, toward a contradiction, that the continued-fraction expansion of $\alpha$ is not eventually $YYY\cdots$. For each $k\ge1$, choose an occurrence of $Y^k$ and extend it to the right by copies of $Y$ as far as possible. By assumption this extension stops after finitely many copies. Let $q_k$ be the position immediately following the last complete copy of $Y$. Then at least $k$ copies of $Y$ occur immediately before $q_k$, while
\[
(a_{q_k},a_{q_k+1},\ldots,a_{q_k+|Y|-1})\ne Y.
\]
There are only finitely many words of length $|Y|$ in the alphabet $\{1,2\}$. Passing to a subsequence in $k$, we may suppose that the word on the left is a fixed word $W\ne Y$.

The sequence of pairs $(\alpha_{q_k},\beta_{q_k-1})$ is bounded. Passing to a further subsequence, assume that it converges to a pair $(\theta,\eta)$, and let
\[
\mathbf c=(\ldots,c_{-2},c_{-1},c_0,c_1,c_2,\ldots)
\]
be the bi-infinite sequence determined by this pair. Put $m:=|Y|$. Since at least $k$ copies of $Y$ occur immediately before $q_k$, for every $r\ge1$ and every $k\ge r$ we have
\[
(a_{q_k-rm},a_{q_k-rm+1},\ldots,a_{q_k-1})=Y^r.
\]
On the other hand,
\[
\beta_{q_k-1}
=[a_{q_k-1};a_{q_k-2},\ldots,a_1]
\longrightarrow
\eta=[c_{-1};c_{-2},\ldots].
\]
By the argument in the proof of Proposition~\ref{prop:theta-eta}, for each fixed $r$ the first $rm$ partial quotients on the two sides agree for all sufficiently large $k$. Since the entries to the left of $q_k$ occur in reverse order in $\beta_{q_k-1}$, this means
\[
(c_{-rm},c_{-rm+1},\ldots,c_{-1})=Y^r.
\]
As $r$ is arbitrary, the entire left-hand side of $\mathbf c$ is $\cdots YYY$.

By the choice of the subsequence,
\[
(a_{q_k},a_{q_k+1},\ldots,a_{q_k+m-1})=W,
\]
and
\[
\alpha_{q_k}=[a_{q_k};a_{q_k+1},\ldots]
\longrightarrow
\theta=[c_0;c_1,\ldots].
\]
The same argument from Proposition~\ref{prop:theta-eta} gives
\[
(c_0,c_1,\ldots,c_{m-1})=W.
\]
Thus
\[
\mathbf c=\cdots YYYV
\]
for a right-infinite word $V$ whose first $m$ entries form $W$. In particular, $V\ne YYY\cdots$ because $W\ne Y$.

The word $\mathbf c$ is also determined by an accumulation pair of the pair sequence of $\alpha$. By Theorem~\ref{thm:limit-point}, every $\ell_h(\mathbf c)$ is an accumulation point of $(\lambda_n(\alpha))_{n\ge1}$. Hence
\[
\mathcal S(\mathbf c)\le\mathcal L(\alpha)<3.
\]
Proposition~\ref{prop:mechanical-characterization} therefore implies that $\mathbf c$ is periodic.

Let $p$ be a period of $\mathbf c$, and put $m=|Y|$. Both $\mathbf c$ and $\cdots YYY\cdots$ have period $pm$. Any position can therefore be shifted into their common left tail by a multiple of $pm$, so the two sequences agree everywhere. This contradicts $V\ne YYY\cdots$.

It follows that the continued-fraction expansion of $\alpha$ is eventually $YYY\cdots$. By Theorem~\ref{thm:characterization-equivalent},
\[
\alpha\sim[\overline Y]\sim[\overline{S(t)}].
\]
Finally, Theorem~\ref{thm:characterization-quadratic} shows that $\alpha$ is a quadratic irrational.
\end{proof}

\section{Lagrange and Markov Constants from Lines of Irrational Slope}
\label{sec:irrational-slope-boundary}
In the preceding section, for $(k_1,k_2,k_3)=(0,0,0)$, we proved that the bi-infinite sequences with Markov value below $3$ arise from mechanical words of rational slope. We now fix general data $(k_1,k_2,k_3,\sigma)$ and replace rational-slope lines by lines of irrational slope. The purpose of this section is to determine the boundary value obtained in this way.

Fix $(k_1,k_2,k_3)\in\mathbb Z_{\geq0}^3$ and $\sigma\in\mathfrak S_3$, and put
\[
K:=3+k_1+k_2+k_3.
\]
\index{Regular line}We consider only positive slopes, in accordance with the definition of the generalized strongly admissible sequence $S(t)$ for reduced fractions $t\in[0,\infty]$ in Chapter~\ref{chap:generalized-cohn-matrices}. A positive-slope line is oriented in the direction of increasing $x$-coordinate. Reversing the orientation merely reverses the resulting sequence and therefore does not change its $\mathcal S$-value. We call a line \emph{regular} if it avoids the marked-point set $\mathcal V$ of $\widetilde{\mathbb R^2}$.

\begin{defi}\label{def:line-sign-sequence}\index{Regular line}
Let $l$ be an oriented regular line of positive slope. Apply the triangle-crossing and edge-crossing rules to every triangle-passage occurrence and every edge-crossing occurrence of $l$, and list the resulting signs in occurrence order. Decompose this bi-infinite sign word into maximal consecutive blocks of equal signs and record their lengths. The resulting bi-infinite sequence of positive integers is denoted by
\[
\mathbf b(l)=(b_n)_{n\in\mathbb Z}.
\]
The index origin is arbitrary, so $\mathbf b(l)$ is defined only up to shift. For rational slope the sequence is periodic; for irrational slope it is generally aperiodic.
\end{defi}

For an irrational slope, the sequence is genuinely indexed by all of $\mathbb Z$. Indeed, the line meets the locally finite triangulation in a discrete sequence of passages unbounded in both directions. By Theorem~\ref{thm:density}, applied to the slope and its negative, every forward and backward tail of the line is dense modulo $\mathbb Z^2$ in the torus. On a short transversal, the two possible triangle signs occur on nonempty open subintervals, and hence each sign occurs infinitely often in both directions. Edge crossings insert blocks of at most $\max\{k_1,k_2,k_3\}$ signs and cannot eliminate all of these sign changes. Thus every run is finite, and the changes of sign are unbounded in both directions.

Since $\mathcal S(\mathbf b)$ is shift-invariant, the arbitrary choice of index origin causes no ambiguity.

\index{Extended sign block}When a finite integer block
\[
W=(b_r,b_{r+1},\ldots,b_s)
\]
of $\mathbf b(l)$ is treated geometrically, we keep not only the finite sign word whose maximal constant-sign runs have these lengths, but also one sign immediately before it and one sign immediately after it. We call the resulting word the \emph{extended sign block} associated with $W$ and denote it by $\widehat W$. Thus $\widehat W$ begins with the sign opposite to the first run of length $b_r$ and ends with the sign opposite to the last run of length $b_s$. Preserving these two boundary signs prevents the first and last runs from merging with adjacent runs and therefore preserves the integer block $W$ exactly.

Our main result in this section is the following.

\begin{theorem}\label{thm:irrational-slope-boundary-value}\index{Irrational-slope boundary value}
Fix $(k_1,k_2,k_3)\in\mathbb Z_{\ge0}^3$ and $\sigma\in\mathfrak S_3$, and put $K=3+k_1+k_2+k_3$. Let $l$ be a regular line of positive irrational slope. Then
\[
\mathcal S(\mathbf b(l))=K.
\]
\end{theorem}

We prepare four lemmas. The first says that the quantity $\ell_n$ defining $\mathcal S$ can be approximated to arbitrary accuracy from a finite window.

\begin{lemm}\label{lem:finite-window-continuity-S}
For every $\varepsilon>0$ there exists $N\ge1$ such that, whenever two bi-infinite sequences of positive integers $\mathbf a=(a_n)_{n\in\mathbb Z}$ and $\mathbf c=(c_n)_{n\in\mathbb Z}$ satisfy
\[
a_i=c_i\qquad(-N\le i\le N),
\]
one has
\[
|\ell_0(\mathbf a)-\ell_0(\mathbf c)|<\varepsilon.
\]
\end{lemm}

\begin{proof}
Choose $N\ge1$ so that
\[
\frac{4}{N(N+1)}<\varepsilon.
\]
We first compare the continued fractions in the positive direction. Put
\[
u_+:=[a_0;a_1,\ldots,a_N]=[c_0;c_1,\ldots,c_N],
\]
and let $Q_N^+$ and $Q_{N-1}^+$ be the denominators of this convergent and the preceding convergent. Set
\[
\xi_{\mathbf a}^+:=[a_{N+1};a_{N+2},\ldots],
\qquad
\xi_{\mathbf c}^+:=[c_{N+1};c_{N+2},\ldots].
\]
Both tails are greater than $1$, and
\[
[a_0;a_1,\ldots]=[a_0;a_1,\ldots,a_N,\xi_{\mathbf a}^+],
\qquad
[c_0;c_1,\ldots]=[c_0;c_1,\ldots,c_N,\xi_{\mathbf c}^+].
\]
Applying Lemma~\ref{lem:finite-infinite-comparison} to each infinite continued fraction and comparing it with $u_+$, then using the triangle inequality and Corollary~\ref{cor:qk>k}, gives
\[
\left|[a_0;a_1,\ldots]-[c_0;c_1,\ldots]\right|
<
\frac{2}{Q_N^+(Q_N^++Q_{N-1}^+)}
\le
\frac{2}{N(N+1)}.
\]

Similarly, put
\[
u_-:=[0;a_{-1},\ldots,a_{-N}]=[0;c_{-1},\ldots,c_{-N}],
\]
let $Q_N^-$ and $Q_{N-1}^-$ be the corresponding two denominators, and define
\[
\xi_{\mathbf a}^-:=[a_{-N-1};a_{-N-2},\ldots],
\qquad
\xi_{\mathbf c}^-:=[c_{-N-1};c_{-N-2},\ldots].
\]
The same argument gives
\[
\left|[0;a_{-1},a_{-2},\ldots]-[0;c_{-1},c_{-2},\ldots]\right|
<
\frac{2}{Q_N^-(Q_N^-+Q_{N-1}^-)}
\le
\frac{2}{N(N+1)}.
\]
Adding the two bounds yields
\[
|\ell_0(\mathbf a)-\ell_0(\mathbf c)|
<\frac{4}{N(N+1)}<\varepsilon.
\]
\end{proof}

The next lemma says that every finite block arising from an irrational-slope line also occurs in a generalized strongly admissible sequence of rational slope.

\begin{lemm}\label{lem:irrational-block-rational-approx}
Let $l$ be a regular line of positive irrational slope. For every finite block $W$ occurring in $\mathbf b(l)$, there exist reduced positive fractions $t=p/q$ with arbitrarily large denominator and regular lines $l_t$ of slope $t$ such that $W$ occurs in $\mathbf b(l_t)$. Moreover, $\mathbf b(l_t)$ is a shift of the periodic sequence
\[
{}^\infty S(t)^\infty=(\ldots,S(t),S(t),S(t),\ldots).
\]
In particular, $W$ occurs in ${}^\infty S(t)^\infty$.
\end{lemm}

\begin{proof}
Fix a finite integer block $W$ and its extended sign block $\widehat W$. The block $\widehat W$ is determined by a finite segment of $l$: one must know the order in which that segment meets the triangles and edges of $\widetilde{\mathbb R^2}$, the local configuration at each incidence, the type of every crossed edge, and the side on which the line passes the relevant marked point. The set $\mathcal V$ is locally finite, and $l$ is regular. Hence all of this finite incidence data, and therefore $\widehat W$, remains unchanged under sufficiently small changes of the slope and intercept.

Write the slope and intercept of $l$ as $\tau$ and $\theta$. Choose a reduced positive fraction $t=p/q$ sufficiently close to $\tau$, with $q$ arbitrarily large. Let $l_t^0$ be the line through $(-\varepsilon,0)$ and $(q-\varepsilon,p)$ for sufficiently small $\varepsilon>0$. The exceptional translates that meet $\mathcal V$ form a discrete set, so $l_t^0$ may be chosen regular. Its sign pattern is invariant under translation by $(q,p)$. Take the half-open fundamental segment from $(-\varepsilon,0)$ to $(q-\varepsilon,p)$, with the initial--terminal cut used to define $\overline{L_t}$. For a sufficiently small left translate, grouping the signs on this segment into runs gives exactly $S(t)$. Moreover, Lemma~\ref{rem:difference-(0,1)(1,infty)} (0) says that $S(t)$ has an even number of entries. The last sign in one period is therefore opposite to the first sign in the next, so no two runs merge at the seam. Consequently,
\[
\mathbf b(l_t^0)={}^\infty S(t)^\infty
\]
up to a shift of the indices.

Translating $l_t^0$ by an integer vector $(m,n)$ changes its intercept by $n-tm$. Since $t=p/q$ is reduced,
\[
\{n-tm\mid m,n\in\mathbb Z\}=\frac1q\mathbb Z.
\]
Thus, by taking $q$ sufficiently large, an integer translate $l_t$ of $l_t^0$ can be chosen with intercept arbitrarily close to $\theta$. Integer translation preserves $\mathcal V$ and both sign rules, so $l_t$ is regular and $\mathbf b(l_t)$ is still a shift of ${}^\infty S(t)^\infty$. Choosing the slope and translated intercept inside the stability neighborhood from the first paragraph makes $l_t$ produce the same extended sign block $\widehat W$. The two boundary signs then ensure that the integer block $W$ itself is preserved exactly.
\end{proof}

The third lemma shows that the collection of finite blocks arising from a regular line of fixed irrational slope is independent of its intercept.

\begin{lemm}\label{lem:same-irrational-slope-same-language}\index{Irrational rotation!density}
Let $l$ and $l'$ be two regular lines of the same positive irrational slope. Then the finite blocks occurring in $\mathbf b(l)$ are exactly the finite blocks occurring in $\mathbf b(l')$.
\end{lemm}

\begin{proof}
Let the common slope be $\tau\notin\mathbb Q$. Translating $l$ by an integer vector $(m,n)$ changes its intercept by $n-\tau m$. Since $\tau$ is irrational, Theorem~\ref{thm:density}, applied to $-\tau$, implies that
\[
\{n-\tau m\mid m,n\in\mathbb Z\}
\]
is dense in $\mathbb R$. Indeed, the orbit $\{-\tau m\bmod1:m\in\mathbb Z\}$ is dense in $\mathbb R/\mathbb Z$, and the integer $n$ can then be chosen to approximate any prescribed real value.

Fix a finite block $W$ of $\mathbf b(l)$ and retain its extended sign block $\widehat W$. Because $l$ is regular, the finite segment producing $\widehat W$ is stable under sufficiently small parallel translations. By the density just noted, an integer translate of $l$ can be made arbitrarily close to $l'$ on this finite region. Integer translations preserve the sign rules, and the two boundary signs prevent the first and last runs from merging with neighboring runs. Hence $W$ occurs in $\mathbf b(l')$. Reversing the roles of $l$ and $l'$ proves the converse inclusion.
\end{proof}

Finally, we record the values of the periodic sequences obtained from rational slopes.

\begin{lemm}\label{lem:rational-values-approach-boundary}\index{Generalized discrete Markov spectrum!accumulation point}\index{Irrational-slope boundary value}
For every reduced positive fraction $t\in(0,\infty)\cap\mathbb Q$,
\[
\mathcal S({}^\infty S(t)^\infty)<K.
\]
Moreover, if $(t_j)_{j\ge0}$ is a sequence of distinct reduced positive fractions converging to an irrational number $\tau$, then
\[
\lim_{j\to\infty}\mathcal S({}^\infty S(t_j)^\infty)=K.
\]
\end{lemm}

\begin{proof}
Let $(m_t,i_t)$ be the $(k_1,k_2,k_3,\sigma)$-GM number and component position corresponding to $t$, and put $k_t=k_{i_t}$. By Theorems~\ref{thm:markov-value-gen} and~\ref{thm:L=S},
\[
\mathcal S({}^\infty S(t)^\infty)
=
\frac{\sqrt{(Km_t-k_t)^2-4}}{m_t},
\]
which is strictly less than $K$.

Now write $t_j=p_j/q_j$ in lowest terms. Since the $t_j$ are distinct and converge to the irrational number $\tau$, we have $p_j+q_j\to\infty$. Indeed, a subsequence on which $p_j+q_j$ remained bounded could contain only finitely many reduced fractions, contradicting either distinctness or irrationality of the limit.

Write
\[
S(t_j)=(a_0^{(j)},\ldots,a_{n_j}^{(j)}).
\]
One period of $\overline{L_{t_j}}$ has at least $p_j+q_j-1$ triangle-passage occurrences, and every such occurrence contributes one sign before adjacent equal signs are grouped. Therefore
\[
\sum_{r=0}^{n_j}a_r^{(j)}
\ge p_j+q_j-1
\longrightarrow\infty.
\]
On the other hand, for every $r$,
\[
a_r^{(j)}
<\ell_r({}^\infty S(t_j)^\infty)
\le\mathcal S({}^\infty S(t_j)^\infty)
<K.
\]
Since $K$ and the $a_r^{(j)}$ are integers, $a_r^{(j)}\le K-1$. The preceding sum can therefore diverge only if $n_j\to\infty$.

By Theorem~\ref{continued-fraction-theorem2} and Corollary~\ref{cor:matrix-combinatorial},
\[
m_{t_j}=(F_{S(t_j)})_{21}
=N(a_1^{(j)},\ldots,a_{n_j}^{(j)})
=q_{n_j},
\]
where $q_{n_j}$ is the denominator of the finite continued fraction $[a_0^{(j)};\ldots,a_{n_j}^{(j)}]$. Corollary~\ref{cor:qk>k} gives $q_{n_j}\ge n_j$, and hence $m_{t_j}\to\infty$.

Finally, $k_{t_j}\in\{k_1,k_2,k_3\}$, so the sequence $(k_{t_j})$ is bounded. Therefore
\[
\frac{\sqrt{(Km_{t_j}-k_{t_j})^2-4}}{m_{t_j}}
=
\sqrt{\left(K-\frac{k_{t_j}}{m_{t_j}}\right)^2-
      \frac4{m_{t_j}^2}}
\longrightarrow K.
\]
\end{proof}

We now prove the theorem.

\begin{proof}[Proof of Theorem~\ref{thm:irrational-slope-boundary-value}]
We first prove
\[
\mathcal S(\mathbf b(l))\le K.
\]
Fix $r\in\mathbb Z$. It is enough to prove $\ell_r(\mathbf b(l))\le K$.

Let $\varepsilon>0$. By Lemma~\ref{lem:finite-window-continuity-S}, if $N$ is sufficiently large, the value of $\ell_r$ is determined to within $\varepsilon$ by the finite central block
\[
b_{r-N},\ldots,b_r,\ldots,b_{r+N}.
\]
By Lemma~\ref{lem:irrational-block-rational-approx}, this block also occurs in the periodic sequence ${}^\infty S(t)^\infty$ for some rational slope $t$. Hence, at a suitable position $j$,
\[
\ell_r(\mathbf b(l))
\le \ell_j({}^\infty S(t)^\infty)+\varepsilon.
\]
Lemma~\ref{lem:rational-values-approach-boundary} then gives
\[
\ell_r(\mathbf b(l))
\le\mathcal S({}^\infty S(t)^\infty)+\varepsilon
<K+\varepsilon.
\]
Since $\varepsilon$ is arbitrary, $\ell_r(\mathbf b(l))\le K$. As $r$ is arbitrary,
\[
\mathcal S(\mathbf b(l))\le K.
\]

We next prove the reverse inequality. Let $\tau$ be the slope of $l$, and choose distinct reduced fractions $t_j\in(0,\infty)\cap\mathbb Q$ with $t_j\to\tau$. By Lemma~\ref{lem:rational-values-approach-boundary},
\[
\mathcal S({}^\infty S(t_j)^\infty)\longrightarrow K.
\]
For each $j$, choose a position at which the periodic sequence attains its $\mathcal S$-value, and shift the sequence so that this position is $0$. Denote the shifted sequence by $\mathbf c^{(j)}$. Since a periodic sequence has only finitely many candidate positions modulo its period, such a position exists, and
\[
\ell_0(\mathbf c^{(j)})
=
\mathcal S({}^\infty S(t_j)^\infty).
\]

We first obtain a coordinatewise convergent subsequence. For every $j$ and $n\in\mathbb Z$,
\[
c_n^{(j)}
<\ell_n(\mathbf c^{(j)})
\le\mathcal S(\mathbf c^{(j)})
<K.
\]
Because $K$ is an integer and the entries are positive integers,
\[
c_n^{(j)}\in\{1,2,\ldots,K-1\}.
\]
Order the integers as $0,1,-1,2,-2,\ldots$. By repeatedly passing to subsequences and then taking the diagonal subsequence, we may assume that, for every fixed $n\in\mathbb Z$, the coordinate $c_n^{(j)}$ is eventually constant. Put
\[
c_n:=\lim_{j\to\infty}c_n^{(j)},
\qquad
\mathbf c=(c_n)_{n\in\mathbb Z}.
\]

We claim that
\[
\ell_0(\mathbf c^{(j)})\longrightarrow\ell_0(\mathbf c).
\]
Given $\varepsilon>0$, choose $N$ from Lemma~\ref{lem:finite-window-continuity-S}. For all sufficiently large $j$,
\[
c_n^{(j)}=c_n
\qquad(-N\le n\le N).
\]
The lemma therefore gives
\[
|\ell_0(\mathbf c^{(j)})-\ell_0(\mathbf c)|<\varepsilon
\]
for all sufficiently large $j$. Since $\ell_0(\mathbf c^{(j)})\to K$, we obtain
\[
\ell_0(\mathbf c)=K.
\]

We now show that every finite block of $\mathbf c$ also occurs in $\mathbf b(l)$. It is enough to consider blocks centered at $0$. Fix $N\ge0$ and set
\[
W=(c_{-N},c_{-N+1},\ldots,c_N).
\]
For all sufficiently large $j$,
\[
(c_{-N}^{(j)},c_{-N+1}^{(j)},\ldots,c_N^{(j)})=W.
\]
Thus $W$ occurs, for infinitely many $j$, in the periodic sequence arising from a regular rational-slope line of slope $t_j$. For each such occurrence, retain the corresponding extended sign block. Once the run lengths in $W$ are fixed, the extended sign block is determined by the sign of its first run; hence there are at most two possibilities. After passing to a subsequence, we may assume that all extended sign blocks are the same word $\widehat W$.

Every triangle passage contributes one sign. Therefore the number of triangle passages, and hence the number of intervening edge crossings, is bounded in terms of the fixed word $\widehat W$. Up to integer translation, only finitely many ordered local crossing configurations can occur. Passing to a further subsequence, we may assume that the complete local data are the same for every occurrence: the order of all passages, the types of the crossed edges, and the side on which the segment passes every relevant marked point are fixed.

Put $L=|\widehat W|$. A segment producing $\widehat W$ has at most $L$ triangle-passage occurrences, and the number of edge crossings is bounded in terms of $L$ as well. There are only finitely many triangle types, each of diameter at most some constant $D$. Take the smallest closed segment containing the passages that produce this occurrence of $\widehat W$, and extend it by Euclidean length $1$ at each end; call the resulting segment $\gamma_j$. Its Euclidean length satisfies
\[
2\le\ell_{\mathrm{Euc}}(\gamma_j)\le C(L)D+2.
\]
Integer translations preserve the triangulation and both sign rules. Translate $\gamma_j$ so that its midpoint lies in $[0,1]^2$. Passing to a subsequence, its midpoint and length converge. Since $t_j\to\tau$, the segments converge to a nondegenerate segment $\gamma$ of slope $\tau$ through a point $P_\infty\in[0,1]^2$.

Let the complete limiting line containing $\gamma$ be
\[
l_{\theta_0}:y=\tau x+\theta_0.
\]
Since $\tau$ is irrational, this line passes through at most one point of $\mathcal V$: the difference of two distinct marked points lies in $(\frac12\mathbb Z)^2$, so the slope of the line joining them is rational. We distinguish whether this possible exceptional point lies on the finite segment relevant to $\widehat W$.

If it does not, the incidence data and the signs defining $\widehat W$ are locally constant under sufficiently small parallel translations. Hence there is an open interval $I$ containing $\theta_0$ such that every line
\[
l_\theta:y=\tau x+\theta,
\qquad \theta\in I,
\]
realizes the same extended sign block $\widehat W$ on the relevant finite segment.

Suppose instead that the limiting segment passes through a marked point $z_0\in\mathcal V$. None of the approximating segments $\gamma_j$ passes through $z_0$. After taking a further subsequence, all $\gamma_j$ pass on the same side of $z_0$. Translating the limiting line slightly toward that side preserves all other incidences and reproduces the fixed local configuration at $z_0$. Thus, in this case, there is a one-sided open interval
\[
I=(\theta_0,\theta_0+\delta)
\quad\text{or}\quad
I=(\theta_0-\delta,\theta_0)
\]
such that every $l_\theta$ with $\theta\in I$ realizes the same extended sign block $\widehat W$.

In either case, the intercepts for which a line of slope $\tau$ meets a marked point form the countable set
\[
E_\tau:=\{b-\tau a\mid(a,b)\in\mathcal V\}.
\]
Choose $\theta\in I\setminus E_\tau$. Then $l_\theta$ is a regular line of slope $\tau$ and contains the integer block $W$ exactly. By Lemma~\ref{lem:same-irrational-slope-same-language}, the collection of finite blocks produced by regular lines of slope $\tau$ is independent of the intercept. Hence $W$ occurs in the original sequence $\mathbf b(l)$.

Finally, let $\varepsilon>0$. Choose $N$ in Lemma~\ref{lem:finite-window-continuity-S} sufficiently large and take the central block
\[
W=(c_{-N},\ldots,c_0,\ldots,c_N).
\]
This block occurs in $\mathbf b(l)$, say centered at a position $q$. Therefore
\[
\ell_q(\mathbf b(l))
>\ell_0(\mathbf c)-\varepsilon
=K-\varepsilon.
\]
Thus $\mathcal S(\mathbf b(l))\ge K-\varepsilon$. Since $\varepsilon$ is arbitrary,
\[
\mathcal S(\mathbf b(l))\ge K.
\]
Together with the opposite inequality, this proves
\[
\mathcal S(\mathbf b(l))=K=3+k_1+k_2+k_3.
\]
\end{proof}

\begin{rema}
For $(k_1,k_2,k_3)=(0,0,0)$, the theorem says that substituting $X\mapsto(2,2)$ and $Y\mapsto(1,1)$ into an irrational-slope mechanical word produces a bi-infinite sequence with boundary value $3$. Thus the theorem may be viewed as a sign-rule formulation of the fact that the accumulation point of the discrete values arising from rational-slope generalized strongly admissible sequences is $3+k_1+k_2+k_3$.
\end{rema}

We conclude the section by showing that the same value is realized as a Lagrange constant.

\begin{coro}\label{cor:irrational-slope-lagrange-value}\index{Irrational-slope boundary value}
Fix $k_1,k_2,k_3\in\mathbb Z_{\ge0}$ and $\sigma\in\mathfrak S_3$. Let $l$ be a regular line of positive irrational slope $\tau$, and let
\[
\mathbf b(l)=(b_n)_{n\in\mathbb Z}
\]
be the bi-infinite sequence obtained from the $(k_1,k_2,k_3,\sigma)$-sign rules. For every $r\in\mathbb Z$, put
\[
\alpha_r:=[0;b_r,b_{r+1},b_{r+2},\ldots].
\quad\text{Then}\quad
\mathcal L(\alpha_r)=K.
\]
\end{coro}

\begin{proof}
It is enough to prove the case $r=0$, since the other cases differ only by a shift of indices. Put
\[
\alpha:=[0;b_0,b_1,b_2,\ldots].
\]
By Theorem~\ref{thm:characterization-lagrange1770},
\[
\mathcal L(\alpha)
=
\limsup_{n\to\infty}
\left(
[b_n;b_{n+1},b_{n+2},\ldots]
+[0;b_{n-1},b_{n-2},\ldots,b_0]
\right).
\]
For the bi-infinite sequence $\mathbf b(l)$,
\[
\ell_n(\mathbf b(l))
=
[b_n;b_{n+1},b_{n+2},\ldots]
+[0;b_{n-1},b_{n-2},\ldots].
\]
We compare the finite and infinite backward continued fractions directly. For $n\ge1$, set
\[
x_n:=[0;b_{n-1},b_{n-2},\ldots,b_0],
\qquad
\xi:=[b_{-1};b_{-2},b_{-3},\ldots]>1,
\]
so that
\[
y_n:=[0;b_{n-1},b_{n-2},\ldots]
=[0;b_{n-1},b_{n-2},\ldots,b_0,\xi].
\]
Let $q_n$ be the denominator of $x_n$ and $q_{n-1}$ the denominator of its preceding convergent. Lemma~\ref{lem:finite-infinite-comparison} gives
\[
|x_n-y_n|
=
\frac1{q_n(\xi q_n+q_{n-1})}
<\frac1{q_n(q_n+q_{n-1})}
<\frac1{q_n^2}.
\]
By Corollary~\ref{cor:qk>k}, $q_n\ge n$, so $|x_n-y_n|\to0$. Hence
\[
\mathcal L(\alpha)
=
\limsup_{n\to\infty}\ell_n(\mathbf b(l)).
\]
Theorem~\ref{thm:irrational-slope-boundary-value} gives $\mathcal S(\mathbf b(l))=K$, so $\ell_n(\mathbf b(l))\le K$ for all $n\in\mathbb Z$. Therefore
\[
\mathcal L(\alpha)
=
\limsup_{n\to\infty}\ell_n(\mathbf b(l))
\le K.
\]

For the reverse inequality, let $\varepsilon>0$. Since $\mathcal S(\mathbf b(l))=K$, there exists $r_0\in\mathbb Z$ such that
\[
\ell_{r_0}(\mathbf b(l))>K-\varepsilon.
\]
By Lemma~\ref{lem:finite-window-continuity-S}, choose $N\ge1$ so large that any bi-infinite sequence with the same central block
\[
W:=(b_{r_0-N},b_{r_0-N+1},\ldots,b_{r_0+N})
\]
has, at the corresponding central position, an $\ell$-value differing from $\ell_{r_0}(\mathbf b(l))$ by less than $\varepsilon$. Thus every occurrence of $W$ as a central block in $\mathbf b(l)$ gives an $\ell$-value greater than $K-2\varepsilon$.

We show that $W$ occurs infinitely often in the positive direction. Retain its extended sign block $\widehat W$ and write
\[
l:y=\tau x+\theta.
\]
Extend the finite portion of $l$ producing $\widehat W$ slightly, and choose an integer $m_0$ and $M\ge1$ so that its $x$-coordinates lie in $[m_0-M,m_0+M]$. Because $l$ is regular, all incidences relevant to this segment have positive distance from the boundaries at which a sign rule can change. Hence there is an open neighborhood $I$ of
\[
\xi_0:=\tau m_0+\theta\pmod 1
\]
in $\mathbb R/\mathbb Z$ with the following property: whenever $m\in\mathbb Z$ satisfies
\[
\tau m+\theta\pmod 1\in I,
\]
the part of $l$ with $x\in[m-M,m+M]$ has, up to an integer translation, the same crossing order, local configurations, and signs as the reference segment. It therefore contains the same extended sign block $\widehat W$, and hence the integer block $W$ exactly.

By Theorem~\ref{thm:density}, the forward rotation orbit
\[
\{\tau m+\theta\pmod1\mid m\in\mathbb Z_{\ge0}\}
\]
and each of its tails are dense in $\mathbb R/\mathbb Z$. We may consequently choose integers
\[
m_1<m_2<m_3<\cdots,
\qquad
m_{i+1}>m_i+2M,
\]
such that $\tau m_i+\theta\pmod1\in I$ for every $i$. The corresponding portions of $l$ are pairwise disjoint and ordered in the positive $x$-direction, so they yield positions
\[
q_1<q_2<q_3<\cdots
\]
at which $W$ occurs as a central block. For every $i$,
\[
\ell_{q_i}(\mathbf b(l))>K-2\varepsilon.
\]
It follows that
\[
\limsup_{n\to\infty}\ell_n(\mathbf b(l))
\ge K-2\varepsilon.
\]
Since $\varepsilon>0$ is arbitrary,
\[
\limsup_{n\to\infty}\ell_n(\mathbf b(l))\ge K.
\]
Combining the two inequalities gives
\[
\mathcal L(\alpha)
=
\limsup_{n\to\infty}\ell_n(\mathbf b(l))
=K.
\]
\end{proof}

The same boundary value is also realized as a Markov constant of explicit real indefinite binary quadratic forms.

\begin{coro}\label{cor:irrational-slope-markov-value}\index{Markov spectrum}
Fix $k_1,k_2,k_3\in\mathbb Z_{\ge0}$ and $\sigma\in\mathfrak S_3$, and put
\[
K=3+k_1+k_2+k_3.
\]
Let $l$ be a regular line of positive irrational slope, and let
\[
\mathbf b(l)=(b_n)_{n\in\mathbb Z}
\]
be the bi-infinite sequence obtained from the $(k_1,k_2,k_3,\sigma)$-sign rules. For every $r\in\mathbb Z$, set
\[
\theta_r=[b_r;b_{r+1},b_{r+2},\ldots],
\qquad
\beta_r=-[0;b_{r-1},b_{r-2},\ldots],
\]
and define
\[
Q_r(x,y):=(x-\theta_r y)(x-\beta_r y).
\]
Then $Q_r$ is a real indefinite binary quadratic form and
\[
\mathcal M(Q_r)=K.
\]
In particular, $K\in\mathcal M$.
\end{coro}

\begin{proof}
Since $\theta_r>1$ and $-1<\beta_r<0$, the form $Q_r$ is canonical reduced. Its two roots are $\theta_r$ and $\beta_r=-1/\eta_r$, where
\[
\eta_r=[b_{r-1};b_{r-2},b_{r-3},\ldots]>1.
\]
Thus the bi-infinite sequence associated with $Q_r$ in Theorem~\ref{thm:markov-infinity-sequence} is precisely $\mathbf b(l)$, up to a shift of the index origin. Therefore
\[
\mathcal M(Q_r)=\mathcal S(\mathbf b(l)).
\]
By Theorem~\ref{thm:irrational-slope-boundary-value}, the right-hand side is $K$. Hence
\[
\mathcal M(Q_r)=K,
\]
as claimed.
\end{proof}

\section{\texorpdfstring{The Relation between the $(0,0,0)$- and $(2,2,2)$-Types}{The Relation between the (0,0,0)- and (2,2,2)-Types}}
Among GM numbers, the $(0,0,0)$-type, consisting of the usual Markov numbers, and the $(2,2,2)$-type have a special relation. To describe it, consider the equation for $k_1=k_2=k_3=2$:
\begin{align}\label{Diophantine5}
    x^2+y^2+z^2+2yz+2zx+2xy=9xyz.
\end{align}

\begin{theorem}\label{thm:0-2-relation-number}\index{Generalized Markov number!squaring correspondence}\index{Generalized Markov number}
If a triple $(a,b,c)$ of positive integers satisfies the Markov equation, then $(a^2,b^2,c^2)$ satisfies the $(2,2,2)$-GM equation. Conversely, if a triple $(A,B,C)$ of positive integers satisfies the $(2,2,2)$-GM equation, then $A,B,C$ are perfect squares and $(\sqrt A,\sqrt B,\sqrt C)$ is a positive integer solution of the Markov equation.
\end{theorem}
\begin{proof}
If $(a,b,c)$ satisfies the Markov equation, then
\[
a^2+b^2+c^2=3abc.
\]
Squaring gives
\[
a^4+b^4+c^4+2a^2b^2+2b^2c^2+2c^2a^2
=9a^2b^2c^2,
\]
which is exactly \eqref{Diophantine5} with $(x,y,z)=(a^2,b^2,c^2)$.

For the converse, use the solution tree in Theorem~\ref{Diophantinetheorem}. Its root $(1,1,1)$ is the componentwise square of the Markov triple $(1,1,1)$. Suppose a vertex is $(a^2,b^2,c^2)$ for a positive Markov triple $(a,b,c)$. The Vieta jump in its first component gives
\[
\left(\frac{b^4+2b^2c^2+c^4}{a^2},b^2,c^2\right)
=\left(\left(\frac{b^2+c^2}{a}\right)^2,b^2,c^2\right).
\]
The Markov equation gives
\[
\frac{b^2+c^2}{a}=3bc-a,
\]
so this is the componentwise square of the Markov triple $(3bc-a,b,c)$. The other two jumps are identical after a cyclic permutation. Thus every vertex of the solution tree is the componentwise square of a positive integer triple. By Theorem~\ref{Diophantinetheorem}, every positive integer solution $(A,B,C)$ of \eqref{Diophantine5} occurs in this tree, and hence $A=a^2$, $B=b^2$, $C=c^2$. Taking the positive square root in
\[
(A+B+C)^2=9ABC
\]
gives $a^2+b^2+c^2=3abc$, so $(a,b,c)$ is a positive integer solution of the Markov equation.
\end{proof}
This relation also gives a simple relation between the discrete Markov spectrum and the $(2,2,2)$-generalized discrete Markov spectrum.
\begin{theorem}\label{thm:0-2-relation}\index{Generalized discrete Markov spectrum}
If $r\in\mathcal M_{0,0,0}$, then $3r\in\mathcal M_{2,2,2}$. Conversely, if $R\in\mathcal M_{2,2,2}$, then $R/3\in\mathcal M_{0,0,0}$.
\end{theorem}
\begin{proof}
For $r\in\mathcal M_{0,0,0}$, there is a Markov number $m$ such that
\[
r=\frac{\sqrt{9m^2-4}}m.
\]
Therefore
\[
3r=\frac{3\sqrt{9m^2-4}}m
=\frac{\sqrt{81m^4-36m^2}}{m^2}
=\frac{\sqrt{(9m^2-2)^2-4}}{m^2}.
\]
By Theorem~\ref{thm:0-2-relation-number}, $m^2$ is a $(2,2,2)$-GM number. Hence $3r\in\mathcal M_{2,2,2}$.

Conversely, if $R\in\mathcal M_{2,2,2}$, there is a $(2,2,2)$-GM number $M$ such that
\[
R=\frac{\sqrt{(9M-2)^2-4}}M.
\]
By Theorem~\ref{thm:0-2-relation-number}, $M=m^2$ for some Markov number $m$. Thus
\[
\frac R3
=\frac13\frac{\sqrt{(9m^2-2)^2-4}}{m^2}
=\frac{\sqrt{9m^2-4}}m
\in\mathcal M_{0,0,0}.
\]
\end{proof}

\section{Frobenius's Uniqueness Conjecture and Its Generalizations}
We discuss natural generalizations of Frobenius's uniqueness conjecture. The original conjecture is as follows.
\begin{conj}\label{uniqueness-conjecture}\index{Frobenius's uniqueness conjecture}
For every Markov number $c$, there is a unique Markov triple $(a,b,c)$ with $a\le b\le c$.
\end{conj}

\begin{lemm}\label{lem:sorted-triple-reciprocal-label}\index{Fraction label!reciprocal labels}
Let $k_1=k_2=k_3=k$, and fix $\sigma\in\mathfrak S_3$. For each positive reduced fraction $t$, let
\[
((x,h),(m_t,i_t),(z,j))
\]
be the vertex corresponding to the unique Farey vertex whose middle entry is $t$, and let $T(t)$ be the triple obtained by arranging $(x,m_t,z)$ in nondecreasing order. Put
\[
T(0)=T(\infty)=(1,1,1).
\]
Then, for $s,t\in\mathbb Q_{\geq0}\cup\{\infty\}$,
\[
T(t)=T(s)\quad\Longleftrightarrow\quad
s=t\ \text{or}\ s=1/t,
\]
where $0^{-1}=\infty$ and $\infty^{-1}=0$. Declare two labels equivalent when they are equal or reciprocal. The map $T$ induces a bijection from the resulting equivalence classes to the nondecreasing positive integer solutions of the $(k,k,k)$-GM equation, and the largest component of $T(t)$ is $m_t$. Consequently, both
\[
\mathbb Q\cap[0,1],\qquad
(\mathbb Q\cap[1,\infty))\cup\{\infty\}
\]
are complete sets of representatives for these equivalence classes.
\end{lemm}
\begin{proof}
The existence and uniqueness of the Farey vertex follow from Proposition~\ref{prop:property-farey}(2), and Proposition~\ref{prop:maximal-middle-component} shows that its second component $m_t$ is the unique largest component.

Because the three parameters are equal, forgetting the position labels in $\mathrm M\mathbb T(k,k,k,\sigma)$ gives the same ordered numerical binary tree for every $\sigma$. Its root is $(1,k+2,1)$, and its two child rules are
\[
(a,b,c)\longmapsto
\left(a,\frac{a^2+kab+b^2}{c},b\right),\qquad
(a,b,c)\longmapsto
\left(b,\frac{b^2+kbc+c^2}{a},c\right).
\]
We show that no ordered numerical triple appears twice in this tree. At a nonroot vertex $(a,b,c)$, we have $b>\max\{a,c\}$. If $a=c$, Lemma~\ref{singular} forces the vertex to be the root $(1,k+2,1)$, so $a\ne c$ at a nonroot vertex. Put
\[
d:=\frac{a^2+kac+c^2}{b}.
\]
This Vieta jump in the largest component is a positive integer, and Corollary~\ref{nonsingular-induction-general} gives $d<\max\{a,c\}$. Hence the parent is uniquely determined as
\[
\begin{cases}
(a,c,d)&(a<c),\\
(d,a,c)&(a>c).
\end{cases}
\]
Each application of this rule strictly decreases the second component from $b$ to $\max\{a,c\}$, so iteration reaches the root. The ordered numerical triple therefore uniquely determines its position in the tree.

The parent rule commutes with interchanging the first and third components. Thus the positions of $(a,b,c)$ and $(c,b,a)$ are obtained from one another by interchanging left and right at every step from the root. By Corollary~\ref{cor:dual-remark}, their fraction labels are reciprocal. More precisely, since the parameters are equal, forgetting position labels in $\mathrm M\mathbb T(k,k,k,\sigma)$ and $\mathrm M\mathbb T(k,k,k,\sigma^*)$ gives the same ordered numerical tree. Hence the two positions have labels $t$ and $1/t$ within $\mathrm M\mathbb T(k,k,k,\sigma)$ itself, and $T(t)=T(1/t)$.

Now take a nondecreasing positive integer solution $(u,v,w)$. The triple $(1,1,1)$ corresponds to the class $\{0,\infty\}$, and the other singular solution $(1,1,k+2)$ corresponds to the root $(1,k+2,1)$ and the fixed label $1$. By Lemma~\ref{singular}, a nonsingular solution has three distinct components. There are therefore exactly two orders with the largest component in the middle: $(u,w,v)$ and $(v,w,u)$. Iterating the uniquely determined parent rule places each order at exactly one position in the tree. Interchanging the outer components shows that their labels are reciprocal. There is no third label with the same nondecreasing triple. This proves both the stated equivalence and surjectivity, and hence the bijection. Each equivalence class plainly has exactly one representative in each of the two displayed sets.
\end{proof}

Taking $k=0$ in Lemma~\ref{lem:sorted-triple-reciprocal-label}, we can restate the original conjecture using the Markov number tree $\mathrm M\mathbb T(0,0,0,\sigma)$ and its fraction labels.
\begin{conj}\label{conj:injectivity-conjecture}\index{Frobenius's uniqueness conjecture}
Fix $\sigma\in\mathfrak S_3$ and consider $\mathrm M\mathbb T(0,0,0,\sigma)$. For reduced fractions $t,s\in[1,\infty]$, let $m_t,m_s$ be the Markov numbers with labels $t,s$. Then $m_t=m_s$ implies $t=s$.
\end{conj}
The following conjecture explains the spectral significance of this elementary formulation.
\begin{conj}\label{uniqueness-conjecture2}
For every $L\in\mathcal M_{0,0,0}$, if $L=\mathcal L(\alpha)=\mathcal L(\beta)$, then $\alpha$ and $\beta$ are $GL(2,\mathbb Z)$-equivalent.
\end{conj}
In other words, an irrational number giving a prescribed Lagrange-spectrum value below $3$ should be unique up to $GL(2,\mathbb Z)$-equivalence. By Markov's theorem, the two conjectures are equivalent.
\begin{prop}\label{prop:conjecture-equivalent}\index{Frobenius's uniqueness conjecture}
Conjecture~\ref{uniqueness-conjecture}, or equivalently Conjecture~\ref{conj:injectivity-conjecture} by Lemma~\ref{lem:sorted-triple-reciprocal-label}, is equivalent to Conjecture~\ref{uniqueness-conjecture2}.
\end{prop}
\begin{proof}
Assume Conjecture~\ref{conj:injectivity-conjecture}. Let $L\in\mathcal M_{0,0,0}$, and let $\alpha,\beta$ be irrational numbers with $L=\mathcal L(\alpha)=\mathcal L(\beta)$. Since $L<3$, Theorem~\ref{lem:lagrange-below-three-quadratic} gives $t,s\in([1,\infty)\cap\mathbb Q)\cup\{\infty\}$ such that
\[
\alpha\sim\alpha_t:=[\overline{S(t)}],\qquad
\beta\sim\alpha_s:=[\overline{S(s)}].
\]
If $\alpha\not\sim\beta$, then $t\ne s$: otherwise $\alpha\sim\alpha_t=\alpha_s\sim\beta$. Conjecture~\ref{conj:injectivity-conjecture} therefore gives $m_t\ne m_s$. On the other hand, Proposition~\ref{prop:equivalent-lagrange1770} and Theorem~\ref{thm:markov-value-gen} imply
\begin{equation}\label{eq:la}
\mathcal L(\alpha)
=\mathcal L(\alpha_t)
=\frac{\sqrt{9m_t^2-4}}{m_t}
=\sqrt{9-\frac4{m_t^2}},\quad \mathcal L(\beta)
=\mathcal L(\alpha_s)
=\frac{\sqrt{9m_s^2-4}}{m_s}
=\sqrt{9-\frac4{m_s^2}}.
\end{equation}
The function $m\mapsto\sqrt{9-4/m^2}$ is strictly increasing for $m\ge1$. Thus $m_t\ne m_s$ gives $\mathcal L(\alpha)\ne\mathcal L(\beta)$, a contradiction. Hence $\alpha\sim\beta$, proving Conjecture~\ref{uniqueness-conjecture2}.

Conversely, assume Conjecture~\ref{uniqueness-conjecture2}, and let $t,s\in([1,\infty)\cap\mathbb Q)\cup\{\infty\}$ with $t\ne s$. Put
\[
\alpha_t:=[\overline{S(t)}],\qquad
\alpha_s:=[\overline{S(s)}].
\]
By Propositions~\ref{prop:rational-slopes-distinct} and~\ref{prop:mechanical-sign}, the proportions of entries equal to $2$ in the periodic sequences ${}^\infty S(t)^\infty$ and ${}^\infty S(s)^\infty$ are $1/t$ and $1/s$, respectively, with $1/\infty=0$. Indeed, replacing $X$ by $(2,2)$ and $Y$ by $(1,1)$ preserves the proportion of $X$'s as the proportion of entries equal to $2$. These proportions are invariant under shifts and reversal. Since $t\ne s$, the sequences cannot agree under either operation. Theorem~\ref{thm:characterization-equivalent} therefore gives $\alpha_t\not\sim\alpha_s$. If $m_t=m_s$, Theorem~\ref{thm:markov-value-gen} would give
\[
\mathcal L(\alpha_t)
=\frac{\sqrt{9m_t^2-4}}{m_t}
=\frac{\sqrt{9m_s^2-4}}{m_s}
=\mathcal L(\alpha_s).
\]
This common value belongs to $\mathcal M_{0,0,0}$, so Conjecture~\ref{uniqueness-conjecture2} would imply $\alpha_t\sim\alpha_s$, again a contradiction. Thus $m_t\ne m_s$, proving Conjecture~\ref{conj:injectivity-conjecture}.
\end{proof}

As a generalization of Conjecture~\ref{uniqueness-conjecture}, consider the following problem.
\begin{prob}\label{uniqueness-question}
For every $(k_1,k_2,k_3)$-GM number $c$, is there a unique $(k_1,k_2,k_3)$-GM triple $(a,b,c)$ satisfying $a\le b\le c$?
\end{prob}
Similarly, Conjecture~\ref{uniqueness-conjecture2} suggests the following problem.
\begin{prob}\label{uniqueness-question2}
For every $L\in\mathcal M_{k_1,k_2,k_3}$, if $L=\mathcal L(\alpha)=\mathcal L(\beta)$, are $\alpha$ and $\beta$ $GL(2,\mathbb Z)$-equivalent?
\end{prob}
For general parameters, the spectral value depends on $k_t$ as well as $m_t$, so the preceding equivalence proof does not apply directly.

Counterexamples to Problem~\ref{uniqueness-question} occur when $k_1,k_2,k_3$ are pairwise distinct. Indeed, $(1,81,17)$ and $(7,81,2)$ are both $(1,2,0)$-GM triples, and are even positive integer solutions of the $(1,2,0)$-GM equation in the displayed orders. Their nondecreasing rearrangements $(1,17,81)$ and $(2,7,81)$ are two distinct GM triples with the same largest component $81$.

Problem~\ref{uniqueness-question2} also has counterexamples. Theorem~\ref{thm:markov-value-gen} provides a natural pair with equal Lagrange constants: for $S(t)=(a_0,\ldots,a_n)$, put
\[
\alpha=[\overline{S(t)}],\qquad
\beta=[\overline{S^*(1/t)}].
\]
Lemma~\ref{rem:difference-(0,1)(1,infty)}(7) gives $S^*(1/t)=(a_0,a_n,\ldots,a_1)$, and Theorem~\ref{thm:markov-value-gen} gives $\mathcal L(\alpha)=\mathcal L(\beta)$. However, $GL(2,\mathbb Z)$-equivalence of purely periodic continued fractions is determined by cyclic shifts of their forward periods. Thus $\alpha\not\sim\beta$ whenever $S(t)$ is not cyclically equivalent to its reversal.

For example, take $(k_1,k_2,k_3)=(1,2,0)$, $\sigma=\mathrm{id}$, and $t=1/3$. Then
\[
S(t)=(5,1,3,3,1,4),\qquad
S^*\left(\frac1t\right)=(5,4,1,3,3,1).
\]
Put
\[
\alpha=[\overline{5,1,3,3,1,4}],\qquad
\beta_0=[\overline{5,4,1,3,3,1}].
\]
Theorem~\ref{thm:markov-value-gen} gives
\[
\mathcal L(\alpha)=\mathcal L(\beta_0)=\frac{2\sqrt{723}}9.
\]
Each period contains exactly one $5$, but its next entry is $1$ in $S(t)$ and $4$ in $S^*(1/t)$. Hence the periods do not agree under any cyclic shift. By Theorem~\ref{thm:characterization-equivalent}, $\alpha$ and $\beta_0$ are not $GL(2,\mathbb Z)$-equivalent, giving a counterexample to Problem~\ref{uniqueness-question2}.

There are also counterexamples that do not arise by reversing a period. The counterexample to Problem~\ref{uniqueness-question} gives one. Put
\[
L=\frac{\sqrt{((3+0+1+2)81-2)^2-4}}{81}
=\frac{2\sqrt{723}}9.
\]
Two distinct quadratic irrationals $\alpha,\beta$ satisfy $L=\mathcal L(\alpha)=\mathcal L(\beta)$, one corresponding to the $(1,2,0)$-GM triple $(1,81,17)$ and the other to $(7,81,2)$.
The first triple corresponds to a vertex of $\mathrm M\mathbb T(1,2,0,\mathrm{id})$, where $81$ has fraction label $1/3$, and gives
\[
\alpha=[\overline{5,1,3,3,1,4}]=\frac{\sqrt{723}+25}9.
\]
The second corresponds to a vertex of $\mathrm M\mathbb T(1,2,0,(1\ 2\ 3))$, where the label is $2/3$, and gives
\[
\beta=[\overline{5,1,1,5,3,2}]=\frac{\sqrt{723}+23}9.
\]
These periods do not agree, and neither is a cyclic shift or reversal of the other. The equality therefore comes from GM numbers in different trees, rather than from the same periodic continued fraction.

Because these examples use different GM trees, they do not disprove injectivity of $t\mapsto m_t$ in a fixed tree. The following example does.
\begin{exam}\label{ex:noninjective-gm-tree}
For $(k_1,k_2,k_3)=(0,0,6)$ and $\sigma=(1\ 3\ 2)$, the tree $\mathrm M\mathbb T(0,0,6,\sigma)$ satisfies
\[
m_{\frac15}=m_{\frac23}=89.
\]
Its root is $((1,3),(2,1),(1,2))$. Repeatedly taking the left child changes the middle fraction label through $1,1/2,1/3,1/4,1/5$. The leftmost GM pair remains $(1,3)$, and the weights at the updated positions are $k_1=k_2=0$. The generation rule therefore gives
\[
m_1=2,\quad m_{\frac12}=5,\quad
m_{\frac13}=\frac{1+5^2}{2}=13,\quad
m_{\frac14}=\frac{1+13^2}{5}=34,\quad
m_{\frac15}=\frac{1+34^2}{13}=89.
\]
On the other hand, taking the right child of the root's left child $((1,3),(5,2),(2,1))$ gives middle label $2/3$ and
\[
m_{\frac23}=\frac{5^2+6\cdot5\cdot2+2^2}{1}=89.
\]
Thus distinct fraction labels in the same GM tree can give the same GM number.
\end{exam}

This example belongs to an infinite family. Define the Fibonacci numbers by $F_0=0$, $F_1=1$, and $F_{n+1}=F_n+F_{n-1}$ for $n\ge1$. For every integer $r\ge0$, put
\[
k_r=\frac{F_{60r+11}-29}{10}\in\mathbb Z_{\geq0}.
\]
Then, in $\mathrm M\mathbb T(0,0,k_r,(1\ 3\ 2))$,
\[
m_{\frac1{30r+5}}=m_{\frac23}=F_{60r+11}.
\]
To verify these assertions, first consider $\mathrm M\mathbb T(0,0,k,(1\ 3\ 2))$ for any $k\in\mathbb Z_{\geq0}$. Its leftmost branch has $x_n:=m_{1/n}$ with
\[
x_1=2,\qquad x_2=5,\qquad
x_{n+1}=\frac{1+x_n^2}{x_{n-1}}\quad(n\geq2),
\]
because the updated positions have weight zero. Put $g_n:=F_{2n+1}$ for $n\geq0$. The Fibonacci recurrence gives $g_0=1$, $g_1=2$, $g_2=5$, and $g_{n+1}=3g_n-g_{n-1}$ for $n\geq1$. Consequently,
\[
g_{n-1}g_{n+1}-g_n^2
=3g_ng_{n+1}-g_n^2-g_{n+1}^2
=g_ng_{n+2}-g_{n+1}^2
\]
is independent of $n$ and equals $g_0g_2-g_1^2=1$. Thus $g_n$ satisfies the same initial conditions and recurrence as $x_n$, proving $m_{1/n}=F_{2n+1}$. The right child of the root's left child also gives
\[
m_{2/3}=5^2+5\cdot2k+2^2=29+10k.
\]
Finally, direct iteration of the Fibonacci recurrence modulo $10$ gives $(F_{60},F_{61})\equiv(0,1)\pmod{10}$. Applying the recurrence to these two consecutive entries yields $F_{n+60}\equiv F_n\pmod{10}$ for every $n\geq0$. Since $F_{11}=89$ and $F_{60r+11}\geq F_{11}$, the number $k_r$ is an integer at least $6$. Substituting $n=30r+5$ and $k=k_r$ in the formulas above proves the claimed equality. The case $r=0$ is Example~\ref{ex:noninjective-gm-tree}.

In these counterexamples the three parameters are not all equal. The following conjecture, restricted to equal parameters, remains open.
\begin{conj}\label{uniqueness-conjecture-gen}\index{Generalized uniqueness conjecture}
Fix $k\in\mathbb Z_{\geq0}$ and $\sigma\in\mathfrak S_3$, and consider $\mathrm M\mathbb T(k,k,k,\sigma)$. For reduced fractions $t,s\in[0,1]$, let $m_t,m_s$ be the corresponding GM numbers. Then $m_t=m_s$ implies $t=s$.
\end{conj}
The interval $[0,1]$ selects one representative from each reciprocal pair of labels $t,1/t$. More precisely, Corollary~\ref{cor:dual-remark} sends an occurrence with label $t$ in the $\sigma$-tree to an occurrence with label $1/t$ in the $\sigma^*=\sigma\circ(1\ 3)$-tree. When $k_1=k_2=k_3=k$, the numerical recurrence is independent of $\sigma$. By Lemma~\ref{lem:sorted-triple-reciprocal-label}, the nondecreasing rearrangement of each GM triple corresponds to exactly one label in $[0,1]$, and its largest component is $m_t$. Thus Conjecture~\ref{uniqueness-conjecture-gen} is equivalent to an affirmative answer to Problem~\ref{uniqueness-question} for $(k,k,k)$.

\chapter{Further Topics}
\label{chap:more-topics}

This chapter collects several directions related to the Lagrange spectrum, the Markov spectrum, generalized Markov numbers, generalized Cohn matrices, and the generalized discrete Markov spectrum. These topics may at first appear rather specialized, but they touch many areas: Diophantine approximation, continued fractions, combinatorics on words, hyperbolic geometry, cluster algebras from surfaces, arithmetic geometry, and dynamical systems. The list is not meant to be exhaustive, and it deliberately overlaps with some of the historical discussion in Chapter~\ref{chap:background-organization}. Some of the topics lie outside the author's own area of expertise, and some recent works cited here are still preprints; the purpose is only to provide entry points for further reading.

\begin{enumerate}
\item
\index{Hall's ray}\index{Freiman's constant}\textbf{Hall's ray, Freiman's constant, and the transition region}

Between the discrete part below $3$, which is governed by Markov's theorem, and the region where Hall's ray begins, lies one of the most complicated parts of the spectrum. Hall proved that $[6,\infty)$ is contained in the Lagrange spectrum, and Freiman determined the initial point $c_F$ of the largest half-line contained in it \cite{hall1947,freiman1975}. Thus a natural next problem is to understand the transition region $[3,c_F)$. The description by bi-infinite continued fractions studied in Chapters~\ref{chap:lagrange-spectrum} and~\ref{chap:markov-spectrum} remains one of the basic tools for this purpose. Standard references include the monograph of Cusick--Flahive and the more recent dynamical and fractal account of Lima--Matheus--Moreira--Romana \cite{cusick-flahive,classical-dynamical}.

\item
\index{Markov spectrum!difference from the Lagrange spectrum}\textbf{Where do $\mathcal L$ and $\mathcal M$ begin to differ?}

This text used the inclusion $\mathcal L\subset\mathcal M$ and the equality below $3$. A natural question is how large $\mathcal M\setminus\mathcal L$ is and how close to $3$ it begins to appear. Moreira proved that, for every half-line $(-\infty,t)$, the intersections of the Lagrange and Markov spectra with this half-line have the same Hausdorff dimension \cite{moreira18}. On the other hand, Erazo--Lima--Matheus--Moreira--Vieira proved that $\inf(\mathcal M\setminus\mathcal L)=3$, showing that the two spectra already differ immediately above $3$ \cite{elmmv}. Thus they are extremely close from the point of view of dimension, but separate at once as sets.

\item
\index{Dynamical Lagrange spectrum}\index{Dynamical Markov spectrum}\textbf{Dynamical Lagrange and Markov spectra}

Classically, $\mathcal L$ and $\mathcal M$ are expressed using the shift $\sigma$ on bi-infinite continued-fraction sequences and the real-valued function $\ell_0$ that sums the two continued fractions at the origin. More generally, for a set $X$, a bijection $\phi\colon X\to X$, and a real-valued function $f\colon X\to\mathbb R$, define
\[
L_{\phi,f}:=\left\{\limsup_{n\to\infty}f(\phi^n(x))\ \middle|\ x\in X\right\},\qquad M_{\phi,f}:=\left\{\sup_{n\in\mathbb Z}f(\phi^n(x))\ \middle|\ x\in X\right\}.
\]
Extended real values are allowed when needed. The first definition uses only forward iterates and does not require $\phi$ to be bijective. \index{Hausdorff dimension}These are dynamical Lagrange and Markov spectra. Cerqueira--Matheus--Moreira start with a smooth area-preserving diffeomorphism of a compact surface and a horseshoe of Hausdorff dimension less than $1$. For generic small area-preserving perturbations and generic smooth real-valued functions, they prove that the dimensions of $L_{\phi,f}\cap(-\infty,t)$ and $M_{\phi,f}\cap(-\infty,t)$ on the corresponding horseshoe agree and depend continuously on $t$ \cite{cerqueira-matheus-moreira2018}. Cerqueira--Moreira--Roma\~na study related questions for geodesic flows on negatively curved surfaces \cite{cerqueira-moreira-romana2022}.

\item
\index{Translation surface}\textbf{Lagrange spectra of translation surfaces}

The classical Lagrange spectrum also has a geometric interpretation in terms of how deeply a geodesic on the modular surface enters a cusp. Starting from this interpretation, one can define analogous spectra in Teichm\"uller dynamics. Hubert--Marchese--Ulcigrai introduced Lagrange spectra for closed $SL(2,\mathbb R)$-invariant loci in the moduli space of translation surfaces \cite{HubertMarcheseUlcigrai2015}. Artigiani--Marchese--Ulcigrai then proved that the Lagrange spectrum of a Veech surface has a Hall ray \cite{artigiani-marchese-ulcigrai2016}. Thus there are meaningful analogues of the Lagrange spectrum beyond the modular surface.

\item
\index{Lagrange spectrum!multiplicative analogue}\textbf{A multiplicative analogue of the Lagrange spectrum}

The classical Lagrange spectrum is related to approximation properties of the arithmetic progression $n\alpha$ modulo $1$. In contrast, Akiyama--Kaneko introduced a multiplicative analogue using fractional parts of the geometric progression $\alpha\beta^n$ \cite{akiyama-kaneko2021}. In particular, when $\beta$ is a Pisot number they prove closedness results for the spectrum; they also describe differences between the case where $\beta$ is an integer and the case where $\beta$ is a quadratic unit, including the existence of intervals and the structure of the first accumulation points and isolated points below them \cite{akiyama-kaneko2021,akiyama-kaneko-corrigendum2022}. More recent work of Akiyama--Kamae--Kaneko extends formulas relating this multiplicative spectrum to symbolic dynamics to broader polynomial and recurrence-theoretic settings \cite{akiyama-kamae-kaneko2025}.

\item
\index{Asymmetric approximation}\index{Inhomogeneous approximation}\textbf{Asymmetric and inhomogeneous approximation spectra}

This text focused on the usual Lagrange and Markov constants, but other spectra arise when one treats left and right approximations asymmetrically or adds inhomogeneous terms. Tornheim's asymmetric approximation is a classical example \cite{tornheim}. For inhomogeneous minima of binary quadratic forms, the series of papers by Barnes and Swinnerton-Dyer is a standard classical reference \cite{BarnesSwinnertonDyer1952InhomogeneousMinimaI,BarnesSwinnertonDyer1952InhomogeneousMinimaII,BarnesSwinnertonDyer1954InhomogeneousMinimaIII}. The minimization problems for binary quadratic forms and the continued-fraction descriptions developed in this text give a useful basis for understanding such variants. It is natural to ask whether the discrete values arising from generalized Markov numbers also appear in spectra other than the standard $\mathcal L$ and $\mathcal M$.

\item
\index{Frobenius's uniqueness conjecture!prime and prime-power cases}\index{Frobenius's uniqueness conjecture}\textbf{Frobenius's uniqueness conjecture and partial results}

Frobenius's uniqueness conjecture is one of the best-known open problems about Markov numbers. Aigner's book presents the conjecture together with Markov's theorem, fraction labels, and perfect matchings as a coherent story \cite{aig}. The conjecture remains open, but it is known when the largest component is prime by work of Button, and for prime powers by results of Schmutz, Lang--Tan, and Zhang \cite{button,lang-tan,zhang}. For generalized Markov numbers, unrestricted injectivity fails by Example~\ref{ex:noninjective-gm-tree}. Conjecture~\ref{uniqueness-conjecture-gen} for equal parameters remains open. Gyoda--Maruyama prove it when the largest component is prime, and also for prime powers $p^m$ with $m\geq2$ under additional conditions on the common coefficient $k$ and the prime $p$ \cite[Theorems 1.6 and 1.7]{gyo-maru}.

\item
\index{Aigner-type conjecture}\index{Fraction label}\textbf{Order by fraction labels and Aigner-type conjectures}

Frobenius's fraction labels do not merely enumerate Markov numbers; they also provide coordinates for studying their order. Lee--Li--Rabideau--Schiffler gave precise inequalities that determine the order of Markov numbers from slopes and lattice data \cite{llrs}. McShane gave a new proof of a related conjecture using convexity of length functions in hyperbolic geometry \cite{mcshane}. For generalized Markov numbers, analogous and generalized results in the case $k=k_1=k_2=k_3$ appear in work of Banaian and Banaian--Huang \cite{banaian,banaian-huang}.

\item
\textbf{Christoffel words and Sturmian words}\index{Sturmian word}

The correspondence between Markov numbers and reduced fractions is closely related to Christoffel words, which encode lattice segments, and to their non-periodic analogues, Sturmian words. Reutenauer made the correspondence between Christoffel words and Markov triples explicit \cite{reut2009} and developed this subject systematically from a combinatorial viewpoint in his monograph \cite{reutenauer}, whose second edition appeared in 2026 \cite{reutenauer2026}. Cohn's matrix construction sends words to products in $SL(2,\mathbb Z)$ and thereby links traces, continued fractions, and Markov numbers \cite{cohn,cohn1971}. The generalized Cohn matrices in this text can be viewed as a generalization of this classical passage from words to matrices. A careful study of reversal, cyclic shift, and the classification of primitive words also clarifies the meaning of the matrix descriptions in Chapter~\ref{chap:generalized-cohn-matrices}.

\item
\index{Simple closed geodesic}\textbf{Simple closed geodesics on the once-punctured torus}

Markov numbers are closely related to lengths of simple closed geodesics on the once-punctured torus. Cohn described Markov forms using geodesics on this torus \cite{cohn1971}, and McShane--Rivin studied lengths of simple geodesics as a norm on homology \cite{mcshane-rivin-norm,mcshane-rivin-simple-curves}. In this viewpoint, fraction labels correspond to slopes of primitive lattice vectors in the plane before points differing by integer vectors are identified to form the torus, and Markov numbers record the associated geodesic lengths. Recent work of Fisac translates the simple length spectrum into combinatorics of cyclic shift classes of integer sequences and gives a new formulation of the uniqueness conjecture \cite{Fisac2025}.

\item
\index{Markoff map}\textbf{Markoff maps, Bowditch space, and McShane identities}

If Markov triples are allowed to take complex values and are regarded as functions on the trivalent tree, one obtains the theory of Markoff maps. Bowditch related Markoff triples to quasifuchsian representations of the once-punctured torus group and derived Bowditch conditions and variants of McShane identities \cite{Bowditch1998}. \index{Character variety}In this direction, the Markov equation is not merely an integer equation; it is a trace identity on the character variety of $SL(2,\mathbb C)$-representations of the free group $F_2$. Although this text mainly treats integer-valued triples, the same tree structure and mutation operations also occur on complex character varieties.

\item
\index{Lambda-length@$\lambda$-length}\textbf{Decorated Teichm\"uller space and $\lambda$-lengths}

\index{Ptolemy relation}In Penner's decorated Teichm\"uller theory, arcs on a surface are assigned positive real numbers called $\lambda$-lengths, and diagonal exchange in a quadrilateral is governed by the Ptolemy relation \cite{penner1987}. This is one geometric origin of the modern principle that flips of triangulations correspond to mutations in cluster algebras. In higher Teichm\"uller theory, Fock--Goncharov introduced positive structures and cluster coordinates on moduli spaces of local systems \cite{fg06,fg09}. The generalized Markov equations in this text are therefore connected not only to formal algebraic modifications, but also to positivity, Ptolemy-type relations, and the geometry of mutation.

\item
\index{Cluster algebra!from surfaces}\textbf{Cluster algebras from surfaces and snake graph calculus}\index{Cluster algebra}\index{Snake graph}

Fomin--Shapiro--Thurston constructed the correspondence between tagged triangulations of bordered surfaces and seeds of cluster algebras \cite{fst}. \index{Lambda-length@$\lambda$-length}Fomin--Thurston related this to the geometry of $\lambda$-lengths and interpreted cluster variables from surfaces as normalized $\lambda$-lengths \cite{ft}. \index{Perfect matching}Musiker--Schiffler--Williams expressed cluster variables from surfaces by perfect matchings of snake graphs \cite{msw}, and Canakci--Schiffler developed the relation between snake graph calculus and continued fractions \cite{snake-graph-calculus2013,canakci-schiffler2018}. The fence posets and GM distances in Chapter~\ref{chap:fence-posets-gm-distance} can be understood naturally by comparing them with such perfect-matching formulas.

\item
\index{Generalized Cohn matrix!cluster structure}\index{Generalized Cohn matrix}\textbf{Further structure of generalized Cohn matrices}

The generalized Cohn matrices treated in this text realize generalized Markov numbers as matrix entries, but they also carry richer structure. In the classical case, Veselov identifies Cohn-matrix indices with Springborn's Markov fractions and describes the associated continued fractions by concatenation on the Conway topograph \cite{veselov2026}, connecting the index in Definition~\ref{def:gc-index} with the arithmetic of rational approximations. \index{Markov--monodromy matrix}In work of Gyoda--Maruyama--Sato, both generalized Cohn matrices and the parallel family called Markov--monodromy matrices are introduced as families of matrices in $SL(2,\mathbb Z)$, and they recover the tree of positive integer solutions of the generalized Markov equation \cite{gyoda-maruyama-sato}. In work of Banaian--Gyoda, these matrices are lifted to matrices with Laurent-polynomial entries, giving cluster structures on generalized Cohn and Markov--monodromy matrices \cite{bana-gyo}. Thus generalized Cohn matrices naturally connect the theory developed here with cluster algebras and the combinatorics of surfaces.

\item
\textbf{The place of the generalized discrete Markov spectrum}

The main theorem places the discrete values built from generalized Markov numbers in the Lagrange spectrum. For fixed coefficients, put $K=3+k_1+k_2+k_3$. The values have the form
\[
\lambda_i(m)=\sqrt{\left(K-\frac{k_i}{m}\right)^2-\frac4{m^2}}
\qquad(i\in\{1,2,3\},\ m\in\mathbb Z_{>0}).
\]
\index{Generalized discrete Markov spectrum!accumulation point}Accumulation of distinct values requires $m\to\infty$, so $K$ is the only possible accumulation point. Conversely, Lemma~\ref{lem:rational-values-approach-boundary} provides a sequence converging to $K$. Thus the fixed-coefficient value set has exactly the accumulation point $K$. Describing intersections of sets from different coefficient triples and the distribution obtained by varying the coefficients remains a further problem. The triple $(0,0,0)$ recovers the classical discrete part below $3$, whereas general coefficients also produce values above $3$.

\item
\index{q-deformation@$q$-deformation}\textbf{$q$-deformations, mirror deformations, and weighted perfect matchings}

Polynomial and Laurent-polynomial deformations of Markov numbers have been studied actively in recent years. Morier-Genoud--Ovsienko introduced $q$-rationals and $q$-continued fractions and related them to the Farey tree and triangulations \cite{morier-genoud-ovsienko2020}. Kantarci Oguz gave a combinatorial model for $q$-deformed Markov numbers using directed posets and rank matrices \cite{oguz2025}. Evans--Jouteur--Morier-Genoud--Ovsienko described $q$-Markov numbers by $q$-deformed Cohn matrices and weighted perfect matchings of snake graphs \cite{evans-jouteur-morier-genoud-ovsienko2025}. Bittmann--Jouteur--Kantarci Oguz--Molander--Yildirim introduced mirror Markov numbers and connected deformed Markov equations, mutations, and orbifold geometry \cite{BittmannJouteurKantarciOguzMolanderYildirim2026}. Generalized Markov numbers may eventually fit into similar weighted or deformed frameworks.

\item
\index{Frieze pattern}\textbf{Frieze patterns and Markov numbers}

\index{Ptolemy relation}Conway--Coxeter frieze patterns are closely related to cluster algebras of type $A$, triangulations, and Ptolemy relations. \index{Perfect matching}Propp explained the combinatorics of frieze patterns and Markov numbers through a model using perfect matchings, giving an intuitive explanation of positivity and the Laurent phenomenon \cite{propp}. Morier-Genoud's survey on frieze patterns is also a useful entry point from classical friezes to modern cluster algebras \cite{morier-genoud2015}. Although the generalized Cohn matrices and fence posets of this text are not frieze patterns themselves, they share the same underlying features: Ptolemy-type relations and perfect matchings, or equivalently order ideals.

\item
\index{Markov--Hurwitz equation}\textbf{Markov--Hurwitz equations and higher-dimensional analogues}

The classical Markov equation has three variables, but higher-dimensional analogues such as
\[
 x_1^2+x_2^2+\cdots+x_n^2=a x_1x_2\cdots x_n+k
\]
have also been studied. For $n\geq4$, Gamburd--Magee--Ronan obtained asymptotic formulas for the number of integer points outside an exceptional set of solution families, assuming that the remaining set is infinite \cite{gamburd-magee-ronan2019}. In higher dimensions the Vieta-jumping graph is no longer a simple trivalent tree, and questions about orbits of integer points, growth, and geometry of numbers become central. It is natural to ask whether generalized Markov equations can also be extended by increasing the number of variables, and whether any connection with spectra survives.

\item
\index{Markov equation!over finite fields}\textbf{Markov equations over finite fields}\index{Markov equation}

One may also study the Markov equation over finite fields $\mathbb F_p$. Then the Vieta involutions generate a graph on a finite set of solutions. Bourgain--Gamburd--Sarnak studied the action of Vieta involutions on congruence solutions of the Markov surface and gave applications to strong approximation and sieve theory \cite{bourgain-gamburd-sarnak2016}. Chen proved that, except for finitely many primes $p$, the group generated by Vieta involutions and coordinate permutations acts transitively on
\[
\left\{(x,y,z)\in\mathbb F_p^3\setminus\{(0,0,0)\}\ \middle|\ x^2+y^2+z^2=3xyz\right\}
\]
\cite{Chen2024Nonabelian}. The origin must be excluded because it is fixed by the group. For generalized Markov equations it is natural to ask what connected components the congruence-solution graphs have and how strong approximation depends on $(k_1,k_2,k_3)$; results in this direction appear in \cite{courcy-litman-mizuno,KingsburyNeuschotz2026}.

\item
\index{Markov-type K3 surface}\textbf{Markov-type K3 surfaces and arithmetic dynamics}

Markov-type equations also appear in the dynamics of K3 surfaces and character varieties. Fuchs--Litman--Silverman--Tran studied orbits of automorphism groups on Markov-type K3 surfaces, including orbit decompositions over finite fields and arithmetic-dynamical properties \cite{fuchs-litman-silverman-tran2022}. In the classical Markov surface, Vieta involutions generate integer points; on K3 surfaces analogous involutions produce more complicated dynamics on more elaborate geometric structures. This viewpoint moves Markov-type equations from trees of integer solutions to actions of automorphism groups on algebraic varieties, and gives an important reference point for considering the algebro-geometric meaning of generalized Markov equations.

\item
\index{Symplectic geometry}\textbf{Symplectic geometry and $\mathbb{CP}^2$}

Markov triples also occur in exceptional bundles on $\mathbb{CP}^2$, weighted projective planes, and Lagrangian cell complexes. Classically, Rudakov used Markov numbers in the classification of exceptional bundles on $\mathbb{CP}^2$ \cite{Rudakov1989MarkovNumbers}. More recently, Evans--Smith studied the relation between Markov numbers and Lagrangian cell complexes in $\mathbb{CP}^2$, showing that Markov numbers arise naturally in symplectic geometry \cite{EvansSmith2018}. In this direction, the Markov equation appears away from Diophantine approximation, in contexts closer to surface degenerations, mirror symmetry, and Floer theory. It remains an open problem to identify what geometric objects are classified by the generalized Markov equations of this text, or what kind of mirror-side deformation they represent.

\item
\index{Hirzebruch--Jung continued fraction}\textbf{Toric geometry and Hirzebruch--Jung continued fractions}

\index{Cyclic quotient singularity}Although this text mainly used regular continued fractions to study the Lagrange and Markov spectra, continued fractions also arise naturally in toric geometry. For a positive integer $r$ and integers $b_1,\dots,b_r\geq2$, the Hirzebruch--Jung continued fraction
\[
  [b_1,\ldots,b_r]_-
  =
  b_1-\frac{1}{\displaystyle b_2-\frac{1}{\displaystyle
  \ddots-\frac{1}{b_r}}}
\]
describes the minimal resolution of a two-dimensional cyclic quotient singularity $\frac{1}{m}(1,q)$, where $0<q<m$ and $\gcd(m,q)=1$. Namely, if $\frac{m}{q}=[b_1,\ldots,b_r]_-$, then the exceptional curves form a chain whose self-intersection numbers are $-b_1,\ldots,-b_r$. For references, see \cite{fulton1993}, \cite{cox-little-schenck2011}, and \cite{popescu-pampu2007}.

\index{Wahl singularity}From this viewpoint Markov numbers are related to degenerations of algebraic surfaces and Wahl singularities. While the preceding item approached the same circle of ideas from the side of symplectic geometry, this is a birational-geometric viewpoint.

Urzua--Zuniga studied the birational-geometric structure of Markov numbers using the Hirzebruch--Jung continued fractions of Wahl singularities associated with Markov triples \cite{urzua-zuniga2023}.

Similar correspondences with cyclic quotient singularities also occur for generalized Markov numbers. For instance, $k$-Wahl chains are Hirzebruch--Jung continued fractions obtained inductively from $[k+2]$ and have been studied as a class including cyclic quotient singularities arising from $k$-generalized Markov triples \cite{gyoda-maruyama-sato,sato2026}.

\item
\index{Markov triple!growth law}\index{Markov triple}\textbf{Growth laws}

Zagier studied the growth of Markov triples. Precisely, if
\[
M_Z(X):=\#\left\{(a,b,c)\in\mathbb Z_{>0}^3\ \middle|\ a\leq b\leq c\leq X,\quad a^2+b^2+c^2=3abc\right\},
\]
then $M_Z(X)\sim C(\log X)^2$ for a positive constant $C$ \cite{zagier1982}. This counts largest components with their multiplicities as triples. Identifying it with the number of distinct Markov numbers requires the uniqueness conjecture. Interpreted as counting simple closed geodesics on the once-punctured torus, this belongs to the same broad circle of ideas as Mirzakhani's theorem on the growth of simple closed geodesics \cite{mirzakhani2008}. For generalized Markov numbers, one may simultaneously count depth in the tree, denominators of fraction labels, values of the numbers, and the associated spectral values; this may reveal growth laws different from the classical case. In computational experiments it is important to specify clearly which parameter is being counted and whether repeated numerical values are counted with multiplicity.

\item
\index{Cluster algebra!mutation invariants}\index{Generalized Markov equation}\index{Cluster algebra}\textbf{Markov-type equations as mutation invariants of cluster algebras}

Cluster mutations often preserve polynomial invariants or positive integer solutions of Diophantine equations. Chen--Li classified sign-equivalence in mutation classes and gave applications to Markov-type equations \cite{chen-li}. Recent preprints by Chen--Li, Bao--Li, and Chen--Jia study Markov-type equations from the viewpoints of mutation invariants, cluster symmetries, and tropicalization \cite{chen-li2,bao-li,chen-jia}. The generalized Markov equation in this text is another example of a mutation-preserved quantity read as an equation for positive integer solutions. Classifying which cluster-algebraic invariants give rise to good Diophantine equations is a natural way to extend the theory of generalized Markov numbers.
\end{enumerate}

These topics show that the generalized discrete Markov spectrum studied in this text is not an isolated construction. It is related to several streams running from classical Diophantine approximation to cluster algebras, hyperbolic geometry, and arithmetic geometry. The purpose of this chapter is to indicate several paths through which readers can move further in these directions.

\appendix

\chapter{Proofs of Standard Facts Used in the Text}\label{chap:standard-facts}
\section{The Bolzano--Weierstrass Theorem}
\begin{theorem}\label{thm:bolzano}\index{Bolzano--Weierstrass theorem}
Every bounded real sequence $(x_n)_{n\geq1}$ has a convergent subsequence.
\end{theorem}

\begin{proof}
Since $(x_n)$ is bounded, there exist real numbers $a_1,b_1$ such that
\[
a_1 \le x_n \le b_1 \qquad (n\in\mathbb{N})
\]
for all $n$. Put $I_1=[a_1,b_1]$.

Next, bisect $I_1$ at its midpoint. Then at least one of the two half-intervals contains infinitely many terms of the sequence $(x_n)$. Indeed, if each of the two half-intervals contained only finitely many terms, then the total number of terms contained in $I_1$ would be finite, contradicting the fact that all terms of $(x_n)$ lie in $I_1$.

Choose one of the half-intervals of $I_1$ that contains infinitely many terms of $(x_n)$, and denote it by $I_2$. In the same way, once $I_k=[a_k,b_k]$ has been defined, bisect it and define $I_{k+1}=[a_{k+1},b_{k+1}]$ to be one of the two halves that contains infinitely many terms of $(x_n)$. In this way we obtain a sequence of closed intervals
\[
I_1 \supset I_2 \supset I_3 \supset \cdots
\]
such that, for each $k$,
\[
b_k-a_k=\frac{b_1-a_1}{2^{k-1}}.
\]
Moreover, each $I_k$ contains infinitely many terms of $(x_n)$.

We now choose a subsequence from these intervals. First choose one term $x_{n_1}$ belonging to $I_1$. Since $I_2$ contains infinitely many terms, we may choose a term belonging to $I_2$ whose index is larger than $n_1$; call it $x_{n_2}$. Continuing in the same way, since $I_k$ contains infinitely many terms, we may choose a term belonging to $I_k$ whose index is larger than $n_{k-1}$; call it $x_{n_k}$. Then
$n_1<n_2<\cdots$, and for each $k$ we have $x_{n_k}\in I_k$. Thus $(x_{n_k})$ is a subsequence of $(x_n)$.

It remains to show that this subsequence converges. Since the closed intervals are nested,
\[
a_1 \le a_2 \le a_3 \le \cdots, \qquad
b_1 \ge b_2 \ge b_3 \ge \cdots.
\]
The sequence $(a_k)$ is bounded above, so by completeness of the real numbers the supremum
\[
x:=\sup\{a_k \mid k\in\mathbb{N}\}
\]
exists.

For each fixed $k$, since $I_{j}\subset I_k$ for all $j\ge k$, and also $a_j\le b_k$ for $j<k$ by the nesting and monotonicity above, $b_k$ is an upper bound of the set $\{a_j\}$. Hence $x\le b_k$. On the other hand, by definition we have $a_k\le x$. Therefore
\[
x\in [a_k,b_k]=I_k \qquad (k\in\mathbb{N}).
\]

Consequently, for each $k$, both $x_{n_k}$ and $x$ belong to $I_k$. Hence
\[
|x_{n_k}-x|\le b_k-a_k=\frac{b_1-a_1}{2^{k-1}}.
\]
The right-hand side tends to $0$ as $k\to\infty$, and therefore $x_{n_k}$ converges to $x$. Thus $(x_n)$ has a convergent subsequence.
\end{proof}

\section{The Cayley--Hamilton Theorem}
\begin{theorem}\label{thm:CH}\index{Cayley--Hamilton theorem}
Let $L$ be a field, $n\geq1$, and $A\in M_n(L)$. For the characteristic polynomial
\[
\chi_A(t)=\det(tE_n-A)=t^n+\sum_{j=0}^{n-1}c_jt^j,
\]
one has
\[
\chi_A(A)=A^n+\sum_{j=0}^{n-1}c_jA^j=0.
\]
Here $E_n$ is the identity matrix and $A^0=E_n$.
\end{theorem}

\begin{proof}
Empty sums are understood to be zero. The adjugate identity gives
\[
\operatorname{adj}(tE_n-A)(tE_n-A)=\det(tE_n-A)E_n=\chi_A(t)E_n.
\]
Each entry of the adjugate has degree at most $n-1$, so write
\[
\operatorname{adj}(tE_n-A)=\sum_{j=0}^{n-1}B_jt^j,\qquad B_j\in M_n(L).
\]
Substitution yields
\[
B_{n-1}t^n+\sum_{j=1}^{n-1}(B_{j-1}-B_jA)t^j-B_0A
=\left(t^n+\sum_{j=0}^{n-1}c_jt^j\right)E_n.
\]
Comparing coefficients gives
\[
B_{n-1}=E_n,\quad B_{j-1}=B_jA+c_jE_n\ (1\leq j\leq n-1),\quad -B_0A=c_0E_n.
\]
Successive substitution gives
\[
B_j=A^{n-1-j}+\sum_{k=j+1}^{n-1}c_kA^{k-1-j}\qquad(0\leq j\leq n-1).
\]
In particular, $B_0=A^{n-1}+\sum_{k=1}^{n-1}c_kA^{k-1}$, whence
\[
\chi_A(A)=A^n+\sum_{k=1}^{n-1}c_kA^k+c_0E_n=B_0A+c_0E_n=0.
\]
\end{proof}

\begin{coro}\label{cor:CH}\index{Cayley--Hamilton theorem}
Let $A=\begin{bsmallmatrix}a&b\\ c&d\end{bsmallmatrix}\in M_2(L)$. Then
\[
A^2-\operatorname{tr}(A)A+\det(A)E_2
=A^2-(a+d)A+(ad-bc)E_2=0.
\]
\end{coro}

\begin{proof}
We have
\[
\chi_A(t)=\det(tE_2-A)
=\det\begin{bmatrix}t-a&-b\\ -c&t-d\end{bmatrix}
=(t-a)(t-d)-bc
=t^2-(a+d)t+(ad-bc).
\]
Therefore, by the theorem,
\[
\chi_A(A)=A^2-(a+d)A+(ad-bc)E_2=0.
\]
\end{proof}

\section{Density of Irrational Rotations}

\begin{theorem}\label{thm:density}\index{Irrational rotation!density}
Let $\tau\in\mathbb R$. The set
\[
\{n+\tau m+\mathbb Z\mid m,n\in\mathbb Z\}\subset\mathbb R/\mathbb Z
\]
is dense in $\mathbb R/\mathbb Z$ if and only if $\tau\notin\mathbb Q$.
If $\tau$ is irrational, then for every $\theta\in\mathbb R$ and $N\in\mathbb Z_{\ge0}$, the set
\[
\{\theta+m\tau+\mathbb Z\mid m\in\mathbb Z,\ m\ge N\}
\]
is dense in $\mathbb R/\mathbb Z$. Thus every forward orbit and each of its tails are dense.
\end{theorem}

\begin{proof}
Since $n+\mathbb Z=0$ for integers $n$, the first displayed set equals $\{m\tau+\mathbb Z\mid m\in\mathbb Z\}$. If $\tau=p/q$ is rational with $q\geq1$, this set has at most $q$ points and is not dense.

Suppose $\tau$ is irrational, and take a nonempty open interval $I$ in the circle, of length $\ell>0$. Choose $H\geq1$ with $1/H<\ell$. Partition $[0,1)$ into $H$ half-open intervals of length $1/H$. Two of the $H+1$ fractional parts of $0,\tau,\dots,H\tau$ lie in the same interval. Thus some $1\leq q\leq H$ satisfies
\[
0<\delta:=\|q\tau\|_{\mathbb Z}<\frac1H,
\qquad \|x\|_{\mathbb Z}:=\min_{r\in\mathbb Z}|x-r|.
\]
The positivity follows from irrationality, and $q\tau+\mathbb Z$ equals either $\delta+\mathbb Z$ or $-\delta+\mathbb Z$. Put $K=\lfloor1/\delta\rfloor$. The circular gaps between the points $0,\delta,\dots,K\delta$ are all at most $\delta$. Reflection preserves these gaps, so the same holds for the points $kq\tau+\mathbb Z$ with $0\leq k\leq K$. Since $\ell>\delta$, the interval $I$ contains one of these points. All indices $kq$ are nonnegative, so the forward orbit $\{m\tau+\mathbb Z\mid m\in\mathbb Z_{\ge0}\}$ is dense. The two-sided orbit contains this forward orbit and is therefore dense as well.

Finally, for every $\theta\in\mathbb R$ and $N\in\mathbb Z_{\ge0}$, the tail $\{\theta+m\tau+\mathbb Z\mid m\in\mathbb Z,\ m\ge N\}$ is the translate of the forward orbit by $\theta+N\tau+\mathbb Z$. Translation preserves density, which proves the remaining assertion.
\end{proof}
\backmatter
\bibliographystyle{alpha-fullauthors}
\bibliography{myrefs}

\begin{thebibliography}{ELMMV24}

\bibitem[Aig13]{aig}
M.~Aigner.
\newblock {\em {{Markov}'s theorem and 100 years of the uniqueness conjecture:
  {A} mathematical journey from irrational numbers to perfect matchings}}.
\newblock Springer, Cham, 2013.

\bibitem[AK21]{akiyama-kaneko2021}
S.~Akiyama and H.~Kaneko.
\newblock {Multiplicative analogue of {Markoff--Lagrange} spectrum and {Pisot}
  numbers}.
\newblock {\em Adv. Math.}, 380:107547, 2021.

\bibitem[AK22]{akiyama-kaneko-corrigendum2022}
S.~Akiyama and H.~Kaneko.
\newblock {Corrigendum to ``Multiplicative analogue of {Markoff--Lagrange}
  spectrum and {Pisot} numbers'' [{Adv. Math.} 380 (2021) 107547]}.
\newblock {\em Adv. Math.}, 394:107996, 2022.

\bibitem[AKK25]{akiyama-kamae-kaneko2025}
S.~Akiyama, T.~Kamae, and H.~Kaneko.
\newblock {Exponential Diophantine approximation and symbolic dynamics}.
\newblock {\em Math. Z.}, 311:70, 2025.

\bibitem[AMU16]{artigiani-marchese-ulcigrai2016}
M.~Artigiani, L.~Marchese, and C.~Ulcigrai.
\newblock {The {Lagrange} spectrum of a {Veech} surface has a {Hall} ray}.
\newblock {\em Groups Geom. Dyn.}, 10(4):1287--1337, 2016.

\bibitem[Ban26]{banaian}
E.~Banaian.
\newblock {Orderings on $k$-Markov numbers}.
\newblock {\em Ramanujan J.}, 71, 2026.
\newblock Article 12; preprint arXiv:2512.04026v2 [math.NT].

\bibitem[BG26]{bana-gyo}
E.~Banaian and Y.~Gyoda.
\newblock {Cluster algebraic interpretation of generalized Markov numbers and
  their matrixizations}, 2026.
\newblock arXiv:2507.06900v3 [math.CO].

\bibitem[BGS16]{bourgain-gamburd-sarnak2016}
J.~Bourgain, A.~Gamburd, and P.~Sarnak.
\newblock {Markoff triples and strong approximation}.
\newblock {\em C. R. Math. Acad. Sci. Paris}, 354(2):131--135, 2016.

\bibitem[BH26]{banaian-huang}
E.~Banaian and M.~Huang.
\newblock {Orderings of generalized $k$-Markov numbers}, 2026.
\newblock arXiv:2604.17445 [math.NT].

\bibitem[BJKMY26]{BittmannJouteurKantarciOguzMolanderYildirim2026}
L.~Bittmann, P.~Jouteur, E.~{Kantarc{\i} O{\u{g}}uz}, M.~Molander, and
  E.~Y{\i}ld{\i}r{\i}m.
\newblock {A mirror deformation of Markov numbers}, 2026.
\newblock arXiv:2602.14802 [math.CO].

\bibitem[BL26]{bao-li}
L.~Bao and F.~Li.
\newblock {The approach of cluster symmetry to Diophantine equations}, 2026.
\newblock arXiv:2508.02005v3 [math.NT].

\bibitem[Bom07]{bombieri2}
E.~Bombieri.
\newblock {Continued fractions and the {Markoff} tree}.
\newblock {\em Expo. Math.}, 25(3):187--213, 2007.

\bibitem[Bow98]{Bowditch1998}
B.~H. Bowditch.
\newblock {Markoff triples and quasifuchsian groups}.
\newblock {\em Proc. Lond. Math. Soc.}, 77(3):697--736, 1998.

\bibitem[BSD52a]{BarnesSwinnertonDyer1952InhomogeneousMinimaI}
E.~S. Barnes and H.~P.~F. Swinnerton-Dyer.
\newblock {The inhomogeneous minima of binary quadratic forms ({I})}.
\newblock {\em Acta Math.}, 87:259--323, 1952.

\bibitem[BSD52b]{BarnesSwinnertonDyer1952InhomogeneousMinimaII}
E.~S. Barnes and H.~P.~F. Swinnerton-Dyer.
\newblock {The inhomogeneous minima of binary quadratic forms ({II})}.
\newblock {\em Acta Math.}, 88:279--316, 1952.

\bibitem[BSD54]{BarnesSwinnertonDyer1954InhomogeneousMinimaIII}
E.~S. Barnes and H.~P.~F. Swinnerton-Dyer.
\newblock {The inhomogeneous minima of binary quadratic forms ({III})}.
\newblock {\em Acta Math.}, 92:199--234, 1954.

\bibitem[But98]{button}
J.~O. Button.
\newblock {The uniqueness of the prime {Markoff} numbers}.
\newblock {\em J. Lond. Math. Soc.}, 58(1):9--17, 1998.

\bibitem[CF89]{cusick-flahive}
T.~W. Cusick and M.~E. Flahive.
\newblock {\em {The Markoff and Lagrange spectra}}, volume~30 of {\em
  {Mathematical Surveys and Monographs}}.
\newblock American Mathematical Society, Providence, RI, 1989.

\bibitem[Che24]{Chen2024Nonabelian}
W.~Y. Chen.
\newblock {Nonabelian level structures, {Nielsen} equivalence, and {Markoff}
  triples}.
\newblock {\em Ann. of Math.}, 199(1):301--443, 2024.

\bibitem[CJ25]{chen-jia}
Z.~Chen and Z.~Jia.
\newblock {Tropicalization and cluster asymptotic phenomenon of generalized
  Markov equations}, 2025.
\newblock arXiv:2511.03428v2 [math.NT].

\bibitem[CL25a]{chen-li2}
Z.~Chen and Z.~Li.
\newblock {A cluster theory approach from mutation invariants to Diophantine
  equations}, 2025.
\newblock arXiv:2501.09435 [math.NT].

\bibitem[CL25b]{chen-li}
Z.~Chen and Z.~Li.
\newblock {Sign-equivalence in cluster algebras: Classification and
  applications to Markov-type equations}.
\newblock {\em J. Pure Appl. Algebra}, 229(10):108058, 2025.

\bibitem[CLS11]{cox-little-schenck2011}
D.~A. Cox, J.~B. Little, and H.~K. Schenck.
\newblock {\em {Toric varieties}}, volume 124 of {\em Graduate Studies in
  Mathematics}.
\newblock American Mathematical Society, Providence, RI, 2011.

\bibitem[CMM18]{cerqueira-matheus-moreira2018}
A.~Cerqueira, C.~Matheus, and C.~G. Moreira.
\newblock {Continuity of {Hausdorff} dimension across generic dynamical
  {Lagrange} and {Markov} spectra}.
\newblock {\em J. Mod. Dyn.}, 12:151--174, 2018.

\bibitem[CMR22]{cerqueira-moreira-romana2022}
A.~Cerqueira, C.~G. Moreira, and S.~Roma{\~n}a.
\newblock {Continuity of {Hausdorff} dimension across generic dynamical
  {Lagrange} and {Markov} spectra {II}}.
\newblock {\em Ergodic Theory Dynam. Systems}, 42(6):1898--1907, 2022.

\bibitem[Coh55]{cohn}
H.~Cohn.
\newblock {Approach to Markoff's minimal forms through modular functions}.
\newblock {\em Ann. of Math.}, 61(1):1--12, 1955.

\bibitem[Coh71]{cohn1971}
H.~Cohn.
\newblock {Representation of {Markoff}'s binary quadratic forms by geodesics on
  a perforated torus}.
\newblock {\em Acta Arith.}, 18(1):125--136, 1971.

\bibitem[{\c{C}}S13]{snake-graph-calculus2013}
{\.I}.~{\c{C}}anak{\c{c}}{\i} and R.~Schiffler.
\newblock {Snake graph calculus and cluster algebras from surfaces}.
\newblock {\em J. Algebra}, 382:240--281, 2013.

\bibitem[CS14]{chsh14}
L.~O. Chekhov and M.~Shapiro.
\newblock {{Teichm{\"u}ller} spaces of {Riemann} surfaces with orbifold points
  of arbitrary order and cluster variables}.
\newblock {\em Int. Math. Res. Not. IMRN}, 2014(10):2746--2772, 2014.

\bibitem[{\c{C}}S18]{canakci-schiffler2018}
{\.I}.~{\c{C}}anak{\c{c}}{\i} and R.~Schiffler.
\newblock {Cluster algebras and continued fractions}.
\newblock {\em Compos. Math.}, 154(3):565--593, 2018.

\bibitem[dCILM26]{courcy-litman-mizuno}
M.~de~Courcy-Ireland, M.~Litman, and Y.~Mizuno.
\newblock {Divisibility by $p$ for Markoff-like surfaces}, 2026.
\newblock arXiv:2509.02187v3 [math.NT].

\bibitem[Dir42]{dirichlet1842}
P.~G.~L. Dirichlet.
\newblock {Verallgemeinerung eines Satzes aus der Lehre von den
  Kettenbr{\"u}chen nebst einigen Anwendungen auf die Theorie der Zahlen}.
\newblock {\em Ber. K. Preuss. Akad. Wiss. Berlin}, pages 93--95, 1842.

\bibitem[EJMGO25]{evans-jouteur-morier-genoud-ovsienko2025}
S.~Evans, P.~Jouteur, S.~Morier-Genoud, and V.~Ovsienko.
\newblock {On $q$-deformed {Markov} numbers. {Cohn} matrices and perfect
  matchings with weighted edges}, 2025.
\newblock arXiv:2507.19080 [math.CO].

\bibitem[ELMMV24]{elmmv}
H.~Erazo, D.~Lima, C.~Matheus, C.~G. Moreira, and S.~Vieira.
\newblock {$\inf(M\setminus L)=3$}, 2024.
\newblock arXiv:2411.06933 [math.NT].

\bibitem[ES18]{EvansSmith2018}
J.~D. Evans and I.~Smith.
\newblock {Markov numbers and Lagrangian cell complexes in the complex
  projective plane}.
\newblock {\em Geom. Topol.}, 22(2):1143--1180, 2018.

\bibitem[Eul44]{euler1737}
L.~Euler.
\newblock {De fractionibus continuis dissertatio}.
\newblock {\em Comment. Acad. Sci. Petropolitanae}, 9:98--137, 1744.
\newblock Written in 1737.

\bibitem[FG06]{fg06}
V.~V. Fock and A.~B. Goncharov.
\newblock {Moduli spaces of local systems and higher {T}eichm{\"{u}}ller
  theory}.
\newblock {\em Publ. Math. Inst. Hautes \'{E}tudes Sci.}, 103:1--211, 2006.

\bibitem[FG09]{fg09}
V.~V. Fock and A.~B. Goncharov.
\newblock {Cluster ensembles, quantization and the dilogarithm}.
\newblock {\em Ann. Sci. \'Ecole Norm. Sup. (4)}, 42(6):865--930, 2009.

\bibitem[Fis25]{Fisac2025}
D.~Fisac.
\newblock {Markov's conjecture on integral necklaces}.
\newblock {\em Bull. Lond. Math. Soc.}, 57(12):4122--4131, 2025.

\bibitem[FLST24]{fuchs-litman-silverman-tran2022}
E.~Fuchs, M.~Litman, J.~H. Silverman, and A.~Tran.
\newblock {Orbits on {K3} surfaces of {Markoff} type}.
\newblock {\em Exp. Math.}, 33(4):663--700, 2024.

\bibitem[Fre68]{freiman}
G.~A. Freiman.
\newblock {Noncoincidence of the {Markov} and {Lagrange} spectra}.
\newblock {\em Mat. Zametki}, 3(2):195--200, 1968.
\newblock In Russian.

\bibitem[Fre75]{freiman1975}
G.~A. Freiman.
\newblock {\em {{Diophantine} approximation and the geometry of numbers ({The
  Markov} problem)}}.
\newblock Kalinin State University, Kalinin, 1975.
\newblock In Russian.

\bibitem[Fro13]{frobenius}
G.~Frobenius.
\newblock {\"Uber die Markoffschen Zahlen}.
\newblock {\em Sitzungsber. Kgl. Preuss. Akad. Wiss.}, pages 458--487, 1913.

\bibitem[FST08]{fst}
S.~Fomin, M.~Shapiro, and D.~Thurston.
\newblock {Cluster algebras and triangulated surfaces. Part {I}: Cluster
  complexes}.
\newblock {\em Acta Math.}, 201:83--146, 2008.

\bibitem[FT18]{ft}
S.~Fomin and D.~Thurston.
\newblock {Cluster algebras and triangulated surfaces. {Part II}: Lambda
  lengths}.
\newblock {\em Mem. Amer. Math. Soc.}, 255(1223), 2018.

\bibitem[Ful93]{fulton1993}
W.~Fulton.
\newblock {\em {Introduction to toric varieties}}, volume 131 of {\em Annals of
  Mathematics Studies}.
\newblock Princeton University Press, Princeton, NJ, 1993.

\bibitem[FZ02]{fzi}
S.~Fomin and A.~Zelevinsky.
\newblock {Cluster algebras {I}: {F}oundations}.
\newblock {\em J. Amer. Math. Soc.}, 15:497--529, 2002.

\bibitem[FZ07]{fziv}
S.~Fomin and A.~Zelevinsky.
\newblock {Cluster algebras {IV}: {C}oefficients}.
\newblock {\em Compos. Math.}, 143:112--164, 2007.

\bibitem[GM23]{gyomatsu}
Y.~Gyoda and K.~Matsushita.
\newblock {Generalization of Markov Diophantine equation via generalized
  cluster algebra}.
\newblock {\em Electron. J. Combin.}, 30(4):P4.10, 2023.

\bibitem[GM26]{gyo-maru}
Y.~Gyoda and S.~Maruyama.
\newblock {Uniqueness theorem of generalized Markov numbers that are prime
  powers}.
\newblock {\em Integers}, 26:A91, 2026.

\bibitem[GMR19]{gamburd-magee-ronan2019}
A.~Gamburd, M.~Magee, and R.~Ronan.
\newblock {An asymptotic formula for integer points on {Markoff--Hurwitz}
  varieties}.
\newblock {\em Ann. of Math.}, 190(3):751--809, 2019.

\bibitem[GMS25]{gyoda-maruyama-sato}
Y.~Gyoda, S.~Maruyama, and Y.~Sato.
\newblock {SL(2,Z)-matrixizations of generalized Markov numbers}, 2025.
\newblock arXiv:2407.08203v3 [math.NT].

\bibitem[Gyo26]{gyoda-generalized}
Y.~Gyoda.
\newblock {Generalized discrete Markov spectra}, 2026.
\newblock arXiv:2512.04547v5 [math.NT].

\bibitem[Hal47]{hall1947}
M.~Hall, Jr.
\newblock {On the sum and product of continued fractions}.
\newblock {\em Ann. of Math.}, 48(4):966--993, 1947.

\bibitem[HMU15]{HubertMarcheseUlcigrai2015}
P.~Hubert, L.~Marchese, and C.~Ulcigrai.
\newblock {Lagrange spectra in Teichm{\"u}ller dynamics via renormalization}.
\newblock {\em Geom. Funct. Anal.}, 25(1):180--255, 2015.

\bibitem[Hur91]{hurwitz}
A.~Hurwitz.
\newblock {Ueber die angen\"aherte Darstellung der Irrationalzahlen durch
  rationale Br\"uche}.
\newblock {\em Math. Ann.}, 39:279--284, 1891.

\bibitem[{Kan}25]{oguz2025}
E.~{Kantarc{\i} O{\u{g}}uz}.
\newblock {Oriented posets, rank matrices and $q$-deformed {Markov} numbers}.
\newblock {\em Discrete Math.}, 348(2):114256, 2025.

\bibitem[Kid22]{kida}
M.~Kida.
\newblock {\em {Renbunsu [Continued Fractions]}}.
\newblock Daigaku Sugaku Spotlight Series 9. Kindai Kagaku Sha, 2022.
\newblock In Japanese.

\bibitem[KN26]{KingsburyNeuschotz2026}
N.~Kingsbury-Neuschotz.
\newblock {Strong approximation for the relative character variety of the
  four-times punctured sphere}, 2026.
\newblock arXiv:2603.04096v3 [math.NT].

\bibitem[KZ73]{korkin-zolotarev1873}
A.~Korkine and G.~Zolotareff.
\newblock {Sur les formes quadratiques}.
\newblock {\em Math. Ann.}, 6:366--389, 1873.

\bibitem[Lag70]{lagrange1770}
J.-L. Lagrange.
\newblock {Additions au m{\'e}moire sur la r{\'e}solution des {\'e}quations
  num{\'e}riques}.
\newblock {\em M\'em. Acad. Roy. Sci. Belles-Lettres Berlin}, 24:111--180,
  1770.
\newblock Volume for 1768, published in 1770.

\bibitem[Lio44]{Liouville}
J.~Liouville.
\newblock {Remarques sur des classes tr{\`e}s-{\'e}tendues de quantit{\'e}s
  dont la valeur n'est ni rationnelle ni m{\^e}me r{\'e}ductible {\`a} des
  irrationnelles alg{\'e}briques}.
\newblock {\em C. R. Acad. Sci. Paris}, 18:883--885, 1844.

\bibitem[LLRS23]{llrs}
K.~Lee, L.~Li, M.~Rabideau, and R.~Schiffler.
\newblock {On the ordering of the Markov numbers}.
\newblock {\em Adv. Appl. Math.}, 143:102453, 2023.

\bibitem[LMMR20]{classical-dynamical}
D.~Lima, C.~Matheus, C.~G. Moreira, and S.~Roma{\~n}a.
\newblock {\em {Classical and dynamical Markov and Lagrange spectra: Dynamical,
  fractal and arithmetic aspects}}.
\newblock World Scientific, 2020.

\bibitem[LT07]{lang-tan}
M.~L. Lang and S.~P. Tan.
\newblock {A simple proof of the Markoff conjecture for prime powers}.
\newblock {\em Geom. Dedicata}, 129:15--22, 2007.

\bibitem[Mar79]{mar1}
A.~Markoff.
\newblock {Sur les formes quadratiques binaires ind\'efinies}.
\newblock {\em Math. Ann.}, 15:381--406, 1879.

\bibitem[Mar80]{mar2}
A.~Markoff.
\newblock {Sur les formes quadratiques binaires ind{\'e}finies ({Second
  m{\'e}moire})}.
\newblock {\em Math. Ann.}, 17:379--399, 1880.

\bibitem[McS21]{mcshane}
G.~McShane.
\newblock {Convexity and Aigner's conjectures}, 2021.
\newblock arXiv:2101.03316 [math.NT].

\bibitem[MG15]{morier-genoud2015}
S.~Morier-Genoud.
\newblock {Coxeter's frieze patterns at the crossroads of algebra, geometry and
  combinatorics}.
\newblock {\em Bull. Lond. Math. Soc.}, 47(6):895--938, 2015.

\bibitem[MGO20]{morier-genoud-ovsienko2020}
S.~Morier-Genoud and V.~Ovsienko.
\newblock {$q$-deformed rationals and $q$-continued fractions}.
\newblock {\em Forum Math. Sigma}, 8:e13, 2020.

\bibitem[Mir08]{mirzakhani2008}
M.~Mirzakhani.
\newblock {Growth of the number of simple closed geodesics on hyperbolic
  surfaces}.
\newblock {\em Ann. of Math.}, 168(1):97--125, 2008.

\bibitem[Mor18]{moreira18}
C.~G. Moreira.
\newblock {Geometric properties of the Markov and Lagrange spectra}.
\newblock {\em Ann. of Math.}, 188(1):145--170, 2018.

\bibitem[MR95a]{mcshane-rivin-norm}
G.~McShane and I.~Rivin.
\newblock {A norm on homology of surfaces and counting simple geodesics}.
\newblock {\em Int. Math. Res. Not. IMRN}, 1995(2):61--69, 1995.

\bibitem[MR95b]{mcshane-rivin-simple-curves}
G.~McShane and I.~Rivin.
\newblock {Simple curves on hyperbolic tori}.
\newblock {\em C. R. Acad. Sci. Paris Ser. I Math.}, 320(12):1523--1528, 1995.

\bibitem[MSW11]{msw}
G.~Musiker, R.~Schiffler, and L.~Williams.
\newblock {Positivity for cluster algebras from surfaces}.
\newblock {\em Adv. Math.}, 227:2241--2308, 2011.

\bibitem[MSW13]{msw2}
G.~Musiker, R.~Schiffler, and L.~Williams.
\newblock {Bases for cluster algebras from surfaces}.
\newblock {\em Compos. Math.}, 149(2):217--263, 2013.

\bibitem[Pen87]{penner1987}
R.~C. Penner.
\newblock {The decorated {Teichm{\"u}ller} space of punctured surfaces}.
\newblock {\em Comm. Math. Phys.}, 113:299--339, 1987.

\bibitem[Per21a]{perron1}
O.~Perron.
\newblock {\"Uber die Approximation irrationaler Zahlen durch rationale. I}.
\newblock {\em S.-B. Heidelberg Akad. Wiss.}, 1921.
\newblock Abh. 4, 17 pp.

\bibitem[Per21b]{perron2}
O.~Perron.
\newblock {\"Uber die Approximation irrationaler Zahlen durch rationale. II}.
\newblock {\em S.-B. Heidelberg Akad. Wiss.}, 1921.
\newblock Abh. 8, 12 pp.

\bibitem[PP07]{popescu-pampu2007}
P.~Popescu-Pampu.
\newblock {The geometry of continued fractions and the topology of surface
  singularities}.
\newblock In J.-P. Brasselet and T.~Suwa, editors, {\em Singularities in
  Geometry and Topology 2004}, volume~46 of {\em Advanced Studies in Pure
  Mathematics}, pages 119--195. Mathematical Society of Japan, Tokyo, 2007.

\bibitem[Pro20]{propp}
J.~Propp.
\newblock {The combinatorics of frieze patterns and Markoff numbers}.
\newblock {\em Integers}, 20:A12, 2020.

\bibitem[Reu09]{reut2009}
C.~Reutenauer.
\newblock {Christoffel words and Markoff triples}.
\newblock {\em Integers}, 9:A26, 2009.

\bibitem[Reu19]{reutenauer}
C.~Reutenauer.
\newblock {\em {From Christoffel words to Markoff numbers}}.
\newblock Oxford University Press, 2019.

\bibitem[Reu26]{reutenauer2026}
C.~Reutenauer.
\newblock {\em {From Christoffel words to Markoff numbers}}.
\newblock Oxford University Press, second edition, 2026.

\bibitem[Rot55]{roth1955}
K.~F. Roth.
\newblock {Rational approximations to algebraic numbers}.
\newblock {\em Mathematika}, 2(1):1--20, 1955.

\bibitem[RS20]{RabideauSchiffler2020}
M.~Rabideau and R.~Schiffler.
\newblock {Continued fractions and orderings on the {Markov} numbers}.
\newblock {\em Adv. Math.}, 370:107231, 2020.

\bibitem[Rud89]{Rudakov1989MarkovNumbers}
A.~N. Rudakov.
\newblock {The {Markov} numbers and exceptional bundles on {$\mathbb{P}^2$}}.
\newblock {\em Math. USSR-Izv.}, 32(1):99--112, 1989.
\newblock English translation of the Russian original in Izv. Akad. Nauk SSSR
  Ser. Mat. 52 (1988), no. 1, 100--112.

\bibitem[Sat26]{sato2026}
Y.~Sato.
\newblock {{$k$}-{Wahl} chains and cyclic quotient singularities}, 2026.
\newblock arXiv:2603.27126 [math.AG].

\bibitem[Tor55]{tornheim}
L.~Tornheim.
\newblock {Asymmetric minima of quadratic forms and asymmetric Diophantine
  approximation}.
\newblock {\em Duke Math. J.}, 22:287--294, 1955.

\bibitem[UZ25]{urzua-zuniga2023}
G.~Urz{\'u}a and J.~P. Z{\'u}{\~n}iga.
\newblock {The birational geometry of {Markov} numbers}.
\newblock {\em Mosc. Math. J.}, 25(2):197--248, 2025.

\bibitem[Ves26]{veselov2026}
A.~P. Veselov.
\newblock {Markov fractions and Cohn matrices}, 2026.
\newblock arXiv:2604.17401v2 [math.NT].

\bibitem[Zag82]{zagier1982}
D.~Zagier.
\newblock {On the number of {Markoff} numbers below a given bound}.
\newblock {\em Math. Comp.}, 39(160):709--723, 1982.

\bibitem[Zha07]{zhang}
Y.~Zhang.
\newblock {An elementary proof of uniqueness of Markoff numbers which are prime
  powers}, 2007.
\newblock arXiv:math/0606283v2 [math.NT].

\end{thebibliography}

\chapter*{List of Symbols}
This list collects the main notation, grouped by subject. References point to
definitions or to statements explaining the notation; page numbers refer to
this volume. GM and GC abbreviate generalized Markov and generalized Cohn,
respectively.

\begingroup
\small
\setlength{\tabcolsep}{5pt}
\renewcommand{\arraystretch}{1.08}
\setlength{\LTleft}{0pt}
\setlength{\LTright}{0pt plus 1fill}
\begin{longtable}{@{}>{\raggedright\arraybackslash}p{.27\textwidth}%
  >{\raggedright\arraybackslash}p{.55\textwidth}%
  >{\raggedright\arraybackslash}p{\dimexpr.18\textwidth-4\tabcolsep\relax}@{}}
\toprule
\textbf{Symbol} & \textbf{Meaning} & \textbf{Reference} \\
\midrule
\endfirsthead
\toprule
\textbf{Symbol} & \textbf{Meaning} & \textbf{Reference} \\
\midrule
\endhead
\midrule
\multicolumn{3}{r@{}}{\footnotesize Continued on the next page}\\
\endfoot
\bottomrule
\endlastfoot
\multicolumn{3}{@{}l}{\textbf{Continued fractions and classical spectra}}\\*[3pt]
$[a_0;a_1,\dots,a_n]$ &
Finite continued fraction with partial quotients $a_0,\dots,a_n$; regular expansions use the terminal convention stated in the definition. &
Def.~\ref{def:regular-continued-fraction}\newline p.~\pageref{def:regular-continued-fraction} \\
\addlinespace[3pt]
$[a_0;a_1,a_2,\dots]$ &
Infinite regular continued fraction, defined as the limit of its finite truncations. &
Def.~\ref{def:infinite-regular-continued-fraction}\newline p.~\pageref{def:infinite-regular-continued-fraction} \\
\addlinespace[3pt]
$p_n,\ q_n$ &
Numerator and denominator of the $n$th convergent, computed by the standard continued-fraction recurrences. &
Prop.~\ref{prop:recursion}\newline p.~\pageref{prop:recursion} \\
\addlinespace[3pt]
$\alpha_n,\ \beta_n$ &
Forward complete quotient $[a_n;a_{n+1},\dots]$ and reversed finite continued fraction $[a_n;a_{n-1},\dots,a_1]$ (for $n\geq1$). &
Thm.~\ref{thm:characterization-lagrange1770}\newline p.~\pageref{thm:characterization-lagrange1770} \\
\addlinespace[3pt]
$[\overline{a_0,\dots,a_{r-1}}]$ &
Purely periodic continued fraction obtained by repeating the displayed finite block indefinitely. &
Def.~\ref{def:periodic-continued-fraction}\newline p.~\pageref{def:periodic-continued-fraction} \\
\addlinespace[3pt]
$GL(2,\mathbb Z)$ &
Unimodular group of integer $2\times2$ matrices with determinant $1$ or $-1$. &
Def.~\ref{def:unimodular-group}\newline p.~\pageref{def:unimodular-group} \\
\addlinespace[3pt]
$\alpha\sim\beta$ &
Unimodular equivalence of irrational numbers under the fractional linear action of $GL(2,\mathbb Z)$. &
Def.~\ref{def:unimodular-equivalence-irrational}\newline p.~\pageref{def:unimodular-equivalence-irrational} \\
\addlinespace[3pt]
$\alpha'$ &
Quadratic conjugate of the quadratic irrational $\alpha$, obtained by changing the sign of its square root. &
Def.~\ref{def:quadratic-irrational}\newline p.~\pageref{def:quadratic-irrational} \\
\addlinespace[3pt]
$I_2,\ R_2$ &
Sets of quadratic irrationals and reduced quadratic irrationals; $I_2(d)$ and $R_2(d)$ specify the discriminant. &
Sect.~\ref{sec:quadratic-irrational}\newline p.~\pageref{sec:quadratic-irrational} \\
\addlinespace[3pt]
$\mathcal L(\alpha)$ &
Lagrange constant of an irrational number $\alpha$, measuring the quality of infinitely many rational approximations. &
Def.~\ref{def:lagrange1770-spectrum}\newline p.~\pageref{def:lagrange1770-spectrum} \\
\addlinespace[3pt]
$\mathcal L$ &
Lagrange spectrum: the set of all Lagrange constants of irrational numbers, including the possible value $\infty$. &
Def.~\ref{def:lagrange1770-spectrum}\newline p.~\pageref{def:lagrange1770-spectrum} \\
\addlinespace[3pt]
$\ell_n(\mathbf a)$ &
Sum of the forward continued fraction at index $n$ and the reciprocal continued fraction extending leftward. &
Cor.~\ref{cor:perron-formula}\newline p.~\pageref{cor:perron-formula} \\
\addlinespace[3pt]
$\mathscr A$ &
Set of all bi-infinite sequences of positive integers used in the symbolic descriptions of the spectra. &
Thm.~\ref{thm:LsubsetS}\newline p.~\pageref{thm:LsubsetS} \\
\addlinespace[3pt]
$\mathcal S(\mathbf a)$ &
Supremum $\sup_{n\in\mathbb Z}\ell_n(\mathbf a)$ associated with a bi-infinite sequence of positive integers; it may be infinite. &
Thm.~\ref{thm:LsubsetS}\newline p.~\pageref{thm:LsubsetS} \\
\addlinespace[3pt]
$\mathcal S$ &
Set of all values $\mathcal S(\mathbf a)$ for $\mathbf a\in\mathscr A$; it coincides with the Markov spectrum. &
Cor.~\ref{cor:M=S}\newline p.~\pageref{cor:M=S} \\
\addlinespace[3pt]
$F_{(a_0,\dots,a_r)}$ &
Continued-fraction matrix: the ordered product of the matrices $\left[\begin{smallmatrix}a_j&1\\1&0\end{smallmatrix}\right]$ for the displayed sequence. &
Thm.~\ref{thm:lagrange1770-continued-fraction-matrix}\newline p.~\pageref{thm:lagrange1770-continued-fraction-matrix} \\
\addlinespace[3pt]
$D(Q)$ &
Discriminant $b^2-4ac$ of the binary quadratic form $Q(x,y)=ax^2+bxy+cy^2$. &
Def.~\ref{def:markov-spectrum}\newline p.~\pageref{def:markov-spectrum} \\
\addlinespace[3pt]
$\mathcal M(Q)$ &
Markov constant of an indefinite binary quadratic form $Q$ without nonzero lattice zeros; infinite when its absolute-value infimum is zero. &
Def.~\ref{def:markov-spectrum}\newline p.~\pageref{def:markov-spectrum} \\
\addlinespace[3pt]
$\mathcal M$ &
Markov spectrum: the set of all Markov constants of the admissible real binary quadratic forms, including $\infty$. &
Def.~\ref{def:markov-spectrum}\newline p.~\pageref{def:markov-spectrum} \\
\addlinespace[3pt]
$\mathcal Q,\ \mathcal R$ &
Set of indefinite real binary quadratic forms without nonzero lattice zeros, and its subset of canonical reduced forms, respectively. &
Thm.~\ref{thm:equivalent-reduced}\newline p.~\pageref{thm:equivalent-reduced} \\
\addlinespace[3pt]
$Q\sim R,\ O_Q$ &
Unimodular equivalence of quadratic forms, and the orbit of $Q$ under integral changes of variables with determinant $\pm1$. &
Def.~\ref{def:unimodular-equivalence-forms}\newline p.~\pageref{def:unimodular-equivalence-forms} \\
\addlinespace[3pt]
\addlinespace[4pt]
\multicolumn{3}{@{}l}{\textbf{Generalized Markov numbers and geometry}}\\*[3pt]
$k_1,k_2,k_3,\ K$ &
Nonnegative integer parameters of the GM equation; $K=3+k_1+k_2+k_3$ abbreviates their sum plus three. &
Eq.~\eqref{Diophantine}\newline p.~\pageref{Diophantine} \\
\addlinespace[3pt]
$\sigma,\ \sigma^\ast$ &
Permutation specifying the GM tree and edge weights; $\sigma^\ast=\sigma\circ(1\ 3)$ exchanges the horizontal and vertical assignments. &
Rem.~\ref{rem:mirror}\newline p.~\pageref{rem:mirror} \\
\addlinespace[3pt]
$\mathrm M\mathbb T$ &
GM tree $\mathrm M\mathbb T(k_1,k_2,k_3,\sigma)$: the ordered binary tree of triples of GM pairs, with each newly produced number in the middle. &
Def.~\ref{def:gm-tree}\newline p.~\pageref{def:gm-tree} \\
\addlinespace[3pt]
$\mathrm F\mathbb T$ &
Farey tree rooted at $(0/1,1/1,1/0)$; its ordered-tree correspondence supplies the fraction labels. &
Sect.~\ref{sec:farey-label}\newline p.~\pageref{sec:farey-label} \\
\addlinespace[3pt]
$r\oplus s$ &
Mediant of fractions: $(a/b)\oplus(c/d)=(a+c)/(b+d)$; used to generate children in the Farey tree. &
Sect.~\ref{sec:farey-label}\newline p.~\pageref{sec:farey-label} \\
\addlinespace[3pt]
$(m_t,i_t)$ &
GM number and position label attached to fraction $t$; $i_t$ records the original equation coordinate. &
Def.~\ref{def:gm-fraction-label}\newline p.~\pageref{def:gm-fraction-label} \\
\addlinespace[3pt]
$k_t$ &
Parameter $k_{i_t}$ attached to the position label of $t$, rather than directly to a geometric edge type. &
Thm.~\ref{thm:Mt-description}\newline p.~\pageref{thm:Mt-description} \\
\addlinespace[3pt]
$u_t$ &
Interior characteristic number: $0<u_t<m_t$ and $m_ru_t\equiv m_s\pmod{m_t}$ for the Farey vertex $(r,t,s)$. &
Def.~\ref{def:characteristic-number}\newline p.~\pageref{def:characteristic-number} \\
\addlinespace[3pt]
$u_0,\ u_\infty$ &
Auxiliary endpoint values $u_0=-k_{\sigma(1)}$ and $u_\infty=1$; these are outside the domain of characteristic numbers. &
Def.~\ref{def:characteristic-number}\newline p.~\pageref{def:characteristic-number} \\
\addlinespace[3pt]
$P_S$ &
Fence poset associated with a finite positive integer sequence $S$, whose entries determine runs of cover orientations. &
Thm.~\ref{thm:continued-fraction-order-ideal}\newline p.~\pageref{thm:continued-fraction-order-ideal} \\
\addlinespace[3pt]
$\mathcal J(P),\ N(P)$ &
Set and number of order ideals of $P$; $N(a_0,\ldots,a_n)$ abbreviates $N(P_{(a_0,\ldots,a_n)})$, with $N()=1$. &
Def.~\ref{def:order-ideal}\newline p.~\pageref{def:order-ideal} \\
\addlinespace[3pt]
$\widetilde{\mathbb R}^{\,2},\ \mathcal V$ &
Triangulated plane and its marked points: lattice points together with edge midpoints selected by the positive parameters. &
Sect.~\ref{sec:GM-distance}\newline p.~\pageref{sec:GM-distance} \\
\addlinespace[3pt]
$\varepsilon(\gamma)$ &
Finite sign sequence assigned to the triangle-passage and edge-crossing occurrences of an oriented generalized arc. &
Def.~\ref{def:sign-assignment}\newline p.~\pageref{def:sign-assignment} \\
\addlinespace[3pt]
$S(\gamma)$ &
Positive integer sequence of maximal constant-sign run lengths in $\varepsilon(\gamma)$, using the chosen endpoint signs. &
Prop.~\ref{prop:length}\newline p.~\pageref{prop:length} \\
\addlinespace[3pt]
$P_\gamma,\ |\gamma|$ &
Fence poset determined by the interior signs of $\gamma$, and its number of order ideals, the GM length. &
Prop.~\ref{prop:length}\newline p.~\pageref{prop:length} \\
\addlinespace[3pt]
$d(A,B)$ &
Infimum of GM lengths of simple generalized arcs joining distinct lattice points; this quantity is not an ordinary metric. &
Def.~\ref{def:gm-distance}\newline p.~\pageref{def:gm-distance} \\
\addlinespace[3pt]
$\gamma^L_{AB},\ \gamma^R_{AB}$ &
Pure left and right push-offs of the segment oriented from $A$ to $B$, with the endpoints fixed. &
Sect.~\ref{sec:GM-distance}\newline p.~\pageref{sec:GM-distance} \\
\addlinespace[3pt]
$U_+,\ U_-$ &
Transfer matrices encoding allowable adjacent membership states of order ideals; their products compute GM lengths. &
Eq.~\eqref{eq:gm-transfer}\newline p.~\pageref{eq:gm-transfer} \\
\addlinespace[3pt]
$C_t$ &
Generalized Cohn matrix attached to fraction $t$, with bottom row $(m_t,u_t)$ and trace $Km_t-k_t$. &
Thm.~\ref{thm:Mt-description}\newline p.~\pageref{thm:Mt-description} \\
\addlinespace[3pt]
$I_t$ &
Index of the GC matrix $C_t$: $I_t=(C_t)_{22}/(C_t)_{21}=u_t/m_t$, strictly increasing with the fraction label. &
Def.~\ref{def:gc-index}\newline p.~\pageref{def:gc-index} \\
\addlinespace[3pt]
$L_t,\ \overline{L_t}$ &
Pure left push-off toward $(q,p)$ for $t=p/q>0$, and its admissible translated-segment perturbation with specified endpoint-crossing conventions. &
Lem.~\ref{lem:admissible-perturbation-Lt}\newline p.~\pageref{lem:admissible-perturbation-Lt} \\
\addlinespace[3pt]
$S(t)$ &
Generalized strongly admissible sequence attached to fraction $t$; interior labels use sign runs of $\overline{L_t}$, with endpoints defined separately. &
Sect.~\ref{sec:description-cohn}\newline p.~\pageref{sec:description-cohn} \\
\addlinespace[3pt]
$w^\dagger$ &
Sign word obtained from $w$ by reversing its order and changing every $+$ to $-$ and conversely. &
Lem.~\ref{lem:half-turn-strongly-admissible}\newline p.~\pageref{lem:half-turn-strongly-admissible} \\
\addlinespace[3pt]
\addlinespace[4pt]
\multicolumn{3}{@{}l}{\textbf{Words and generalized spectra}}\\*[3pt]
$\mathcal M_{k_1,k_2,k_3,\sigma}$ &
Spectral values obtained from the GM numbers and position labels in the tree with fixed parameters and permutation. &
Sect.~\ref{sec:generalized-discrete-definition}\newline p.~\pageref{sec:generalized-discrete-definition} \\
\addlinespace[3pt]
$\mathcal M_{k_1,k_2,k_3}$ &
Generalized discrete Markov spectrum: the union of $\mathcal M_{k_1,k_2,k_3,\sigma}$ over all $\sigma\in\mathfrak S_3$. &
Sect.~\ref{sec:generalized-discrete-definition}\newline p.~\pageref{sec:generalized-discrete-definition} \\
\addlinespace[3pt]
$\alpha_S$ &
The quadratic irrational $[\overline S]$ defined by the purely periodic continued fraction with finite positive-integer period $S$. &
Thm.~\ref{thm:markov-value-gen}\newline p.~\pageref{thm:markov-value-gen} \\
\addlinespace[3pt]
$Q_S$ &
The binary quadratic form $(x-\alpha_Sy)(x-\alpha'_Sy)$, where $\alpha'_S$ is the quadratic conjugate of $\alpha_S$. &
Thm.~\ref{thm:markov-value-gen2}\newline p.~\pageref{thm:markov-value-gen2} \\
\addlinespace[3pt]
$\mathbf b^R(t,\theta)$,\newline$\mathbf b^L(t,\theta)$ &
Right and left mechanical words of slope $t\in[1,\infty]$ and intercept $\theta$, defined using ceiling and floor functions. &
Def.~\ref{def:mechanical-word}\newline p.~\pageref{def:mechanical-word} \\
\addlinespace[3pt]
$\iota$ &
The substitution $X\mapsto(2,2)$ and $Y\mapsto(1,1)$; word values satisfy $\mathcal S(\mathbf w)=\mathcal S(\iota(\mathbf w))$. &
Sect.~\ref{sec:mechanical-words}\newline p.~\pageref{sec:mechanical-words} \\
\addlinespace[3pt]
$|w|,\ |w|_X$ &
Length of the finite word $w$, and the number of its letters equal to $X$, respectively. &
Lem.~\ref{rem:mechanical-slope-period}\newline p.~\pageref{rem:mechanical-slope-period} \\
\addlinespace[3pt]
$C^X(\mathbf b),\ C^Y(\mathbf b)$ &
Sequences counting intervening $Y$ letters between consecutive $X$ letters, and intervening $X$ letters between consecutive $Y$ letters. &
Def.~\ref{def:characteristic-word-sequences}\newline p.~\pageref{def:characteristic-word-sequences} \\
\addlinespace[3pt]
$\lambda,\ \rho$ &
Word substitutions given by $\lambda(X)=X$, $\lambda(Y)=XY$, $\rho(X)=XY$, and $\rho(Y)=Y$. &
Sect.~\ref{sec:markov-theorem}\newline p.~\pageref{sec:markov-theorem} \\
\addlinespace[3pt]
$\mathbf b(l)$ &
Bi-infinite sequence of maximal constant-sign run lengths produced by the sign rules along the oriented regular line $l$. &
Def.~\ref{def:line-sign-sequence}\newline p.~\pageref{def:line-sign-sequence} \\
\addlinespace[3pt]
$\widehat W$ &
Sign block representing the integer block $W$, extended by one sign at each end to preserve its run boundaries. &
Sect.~\ref{sec:irrational-slope-boundary}\newline p.~\pageref{sec:irrational-slope-boundary} \\
\addlinespace[3pt]
\end{longtable}
\endgroup

\renewcommand{\indexintro}{Selected terms, with page references to definitions and principal discussions.}
\printindex

\end{document}